\documentclass[openany]{amsbook}
\usepackage[dvips]{graphicx}
\usepackage{amscd}
\usepackage{amsmath}
\usepackage{amsfonts}
\usepackage{amssymb}
\newcommand{\frogbracedef}{\newcommand{\frogbrace}[2]{\protect\underset{##1}{\protect\underbrace{##2,\dots,##2}}}}
\newcommand{\frogelbowdef}{\newcommand{\frogelbow}{\ }}
\newcommand{\firkelbowdef}{\newcommand{\firkelbow}{\protect\frogelbow}}
\newcommand{\LOSitemdef}{\newcommand{\LOSitem}[3]{\noindent\parbox[t]{0.4\textwidth}{\raggedright ##1}\kern0.05\textwidth\parbox[t]{0.3\textwidth}{\raggedright ##2}\kern0.05\textwidth\parbox[t]{0.2\textwidth}{\raggedleft ##3}\medskip}}
\newcommand{\olatheoremformat}{}
\renewcommand{\olatheoremformat}{\hskip-\parindent\upshape\bfseries}
\newtheorem{theorem}{\olatheoremformat Theorem}[chapter]
\newtheorem{corollary}[theorem]{\olatheoremformat Corollary}
\newtheorem{lemma}[theorem]{\olatheoremformat Lemma}
\newtheorem{proposition}[theorem]{\olatheoremformat Proposition}
\theoremstyle{definition}
\newtheorem{definition}[theorem]{\olatheoremformat Definition}
\newtheorem{remark}[theorem]{\olatheoremformat Remark}
\newtheorem{scholium}[theorem]{\olatheoremformat Scholium}
\newtheorem{observation}[theorem]{\olatheoremformat Observation}
\newtheorem{example}[theorem]{\olatheoremformat Example}
\newtheorem{case}{\olatheoremformat Case}

\newtheorem{claim}[theorem]{\olatheoremformat Claim}
\theoremstyle{remark}
\newtheorem*{acknowledgements}{Acknowledgements}

\newlength{\displayboxwidth}

\renewcommand{\theenumi}{\roman{enumi}}

\newcounter{figurelink}
\newcounter{tablelink}
\newcounter{glimmer}
\renewcommand{\subset}{\subseteq}
\renewcommand{\supset}{\supseteq}
\def\hooklongrightarrow{\lhook\joinrel\longrightarrow} 
\newbox\frogdown
\newlength\frogdrop
\def\hookdownarrow{\setlength{\unitlength}{0.4pt}\setbox\frogdown=\hbox to 0pt{\hss $\displaystyle \downarrow $\hss }\setlength{\frogdrop}{2.5\ht\frogdown}\raisebox{0pt}[5\unitlength][\frogdrop]
{\begin{picture}(0,5)(10,0)
\put(5,0){\oval(10,10)[t]}
\end{picture}\lower\ht\frogdown\box\frogdown}}
\def\openone
{\mathchoice
{\hbox{\upshape \small1\kern-3.3pt\normalsize1}}
{\hbox{\upshape \small1\kern-3.3pt\normalsize1}}
{\hbox{\upshape \tiny1\kern-2.3pt\SMALL1}}
{\hbox{\upshape \Tiny1\kern-2pt\tiny1}}}
\newbox\ipbox
\newcommand{\ip}[2]{\left\langle #1\mathrel{\mathchoice
{\setbox\ipbox=\hbox{$\displaystyle \left\langle\mathstrut #1#2\right\rangle$}
\vrule height\ht\ipbox width0.25pt depth\dp\ipbox}
{\setbox\ipbox=\hbox{$\textstyle \left\langle\mathstrut #1#2\right\rangle$}
\vrule height\ht\ipbox width0.25pt depth\dp\ipbox}
{\setbox\ipbox=\hbox{$\scriptstyle \left\langle\mathstrut #1#2\right\rangle$}
\vrule height\ht\ipbox width0.25pt depth\dp\ipbox}
{\setbox\ipbox=\hbox{$\scriptscriptstyle \left\langle\mathstrut #1#2\right\rangle$}
\vrule height\ht\ipbox width0.25pt depth\dp\ipbox}
} #2\right\rangle}
\newcommand{\diracb}[1]{\left\langle #1\mathrel{\mathchoice
{\setbox\ipbox=\hbox{$\displaystyle \left\langle\mathstrut #1\right.$}
\vrule height\ht\ipbox width0.25pt depth\dp\ipbox}
{\setbox\ipbox=\hbox{$\textstyle \left\langle\mathstrut #1\right.$}
\vrule height\ht\ipbox width0.25pt depth\dp\ipbox}
{\setbox\ipbox=\hbox{$\scriptstyle \left\langle\mathstrut #1\right.$}
\vrule height\ht\ipbox width0.25pt depth\dp\ipbox}
{\setbox\ipbox=\hbox{$\scriptscriptstyle \left\langle\mathstrut #1\right.$}
\vrule height\ht\ipbox width0.25pt depth\dp\ipbox}
}\right. }
\newcommand{\dirack}[1]{\left. \mathrel{\mathchoice
{\setbox\ipbox=\hbox{$\displaystyle \left.\mathstrut #1\right\rangle$}
\vrule height\ht\ipbox width0.25pt depth\dp\ipbox}
{\setbox\ipbox=\hbox{$\textstyle \left.\mathstrut #1\right\rangle$}
\vrule height\ht\ipbox width0.25pt depth\dp\ipbox}
{\setbox\ipbox=\hbox{$\scriptstyle \left.\mathstrut #1\right\rangle$}
\vrule height\ht\ipbox width0.25pt depth\dp\ipbox}
{\setbox\ipbox=\hbox{$\scriptscriptstyle \left.\mathstrut #1\right\rangle$}
\vrule height\ht\ipbox width0.25pt depth\dp\ipbox}
} #1\right\rangle}
\long\def\TeXButton#1#2{#2}%

\newcommand{\frogkeywords}{\renewcommand{\thefootnote}{}\footnotetext{\textit{Key words and phrases. } $C^*$-algebras, Fourier basis, irreducible representations, Hilbert space, wavelets, radix-representations, lattices, iterated function systems.}}

\numberwithin{equation}{chapter}
\newlength{\asubwidth}
\newlength{\qedskip}
\newlength{\qedadjust}
\frogbracedef
\newcommand{\frogleg}{\\}
\frogelbowdef
\firkelbowdef
\LOSitemdef
\input cyracc.def

\makeatletter
\AtBeginDocument{  \renewcommand\mainmatter{\cleardoublepage\pagenumbering{arabic}}  \renewcommand\backmatter{\if@openright\cleardoublepage\else\clearpage\fi}  \renewcommand\chapter{\if@openright\cleardoublepage\else\clearpage\fi
    \thispagestyle{plain}\global\@topnum\z@
    \@afterindenttrue \secdef\@chapter\@schapter}  }
\def\ps@headings{\ps@empty
  \def\@evenhead{\normalfont\scriptsize
      \rlap{\thepage}\hfil \leftmark{}{}\hfil}
  \def\@oddhead{\normalfont\scriptsize \hfil
      \rightmark{}{}\hfil \llap{\thepage}}
  \let\@mkboth\markboth
  \def\partmark{\@secmark\markright\partrunhead\partname}
  \def\chaptermark{
    \@secmark\markright\chapterrunhead{}}
  \let\sectionmark\@gobble
}
\renewenvironment{proof}[1][\proofname]{\par \normalfont
  \topsep6\p@\@plus6\p@ \trivlist \itemindent\normalparindent
  \item[\hskip\labelsep\scshape
      \hskip-\itemindent\itshape     #1\@addpunct{.}]\ignorespaces
}{  \qed\endtrivlist
}
\makeatother
\begin{document}
\frontmatter
\title[REPRESENTATION THEORY AND NUMERICAL AF-INVARIANTS]{{\huge Representation Theory and}\\{\huge Numerical AF-invariants\bigskip}\\The representations and centralizers\\of certain states on $\mathcal{O}_{d}$}
\author{Ola~Bratteli}
\address{Mathematics Institute\\
University of Oslo\\
PB 1053 Blindern\\
N-0316 Oslo\\
Norway}
\email{bratteli@math.uio.no}
\urladdr{http://www.math.uio.no/\symbol{126}bratteli/}
\author{Palle~E.~T.~ Jorgensen}
\address{Department of Mathematics\\
The University of Iowa\\
Iowa City, IA 52242-1419\\
U.S.A.}
\email{jorgen@math.uiowa.edu}
\urladdr{http://www.math.uiowa.edu/\symbol{126}jorgen/}
\author{Vasyl' Ostrovs'ky\u{\i}}
\address{Institute of Mathematics\\
National Academy of Sciences of Ukraine\\
Tereshchenkivs'ka str.~3\\
252601 Kiev 4\\
Ukraine}
\email{vo@imath.kiev.ua}
\thanks{Work supported by the Norwegian Research Council, the University of Oslo, the
U.S. National Science Foundation, and the CRDF research grant UM1-311.}

\begin{abstract}
Let $\mathcal{O}_{d}$ be the Cuntz algebra on generators $S_{1},\dots,S_{d}$,
$2\leq d<\infty$, and let $\mathcal{D}_{d}\subset\mathcal{O}_{d}$ be the
abelian subalgebra generated by monomials $S_{\alpha_{{}}}^{{}}S_{\alpha_{{}}%
}^{\ast}=S_{\alpha_{1}}^{{}}\cdots S_{\alpha_{k}}^{{}}S_{\alpha_{k}}^{\ast
}\cdots S_{\alpha_{1}}^{\ast}$ where $\alpha=\left(  \alpha_{1}\dots\alpha
_{k}\right)  $ ranges over all multi-indices formed from $\left\{
1,\dots,d\right\}  $. In any representation of $\mathcal{O}_{d}$,
$\mathcal{D}_{d}$ may be simultaneously diagonalized. Using $S_{i}^{{}}\left(
S_{\alpha}^{{}}S_{\alpha}^{\ast}\right)  =\left(  S_{i\alpha}^{{}}S_{i\alpha
}^{\ast}\right)  S_{i}^{{}}$, we show that the operators $S_{i}$ from a
general representation of $\mathcal{O}_{d}$ may be expressed directly in terms
of the spectral representation of $\mathcal{D}_{d}$. We use this in describing
a class of type $\mathrm{III}$ representations of $\mathcal{O}_{d}$ and
corresponding endomorphisms, and the heart of the memoir is a description of an
associated family of AF-algebras arising as the fixed-point algebras of the
associated modular automorphism groups. Chapters \ref{Rem}--\ref{App.Fur} are
devoted to finding effective methods to decide isomorphism and non-isomorphism
in this class of AF-algebras.
\end{abstract}\maketitle
\tableofcontents


\markboth{\shorttitle}{\shorttitle}

\chapter*{Preface}

The present memoir consists of two parts. The first part encompasses Chapters
\ref{Gen}--\ref{Sub}, and is concerned with the description of a class of
representations of the Cuntz algebra $\mathcal{O}_{d}$, starting out with a
very general description of such representations. The second part encompasses
Chapters \ref{Sub}--\ref{App.Fur} and is a description of a class of
AF-algebras with constant incidence matrices $J$ of the special form
(\ref{eqBrunt.1}). The two parts are thus connected by Chapter \ref{Sub},
where it is explained how these AF-algebras arise as the fixed-point algebras
of modular automorphism groups associated to certain states on $\mathcal{O}%
_{d}$. Readers who are not interested in representation theory can therefore
read the memoir from Chapter \ref{Rem}. Since the special examples we study can
be understood very concretely, we hope that the memoir may serve as an
invitation for graduate students who want to study isomorphism and invariants
in more general settings.
\bigskip

\begin{center}
\settowidth{\displayboxwidth}
{Oslo, Iowa City and Kiev, June 1999, January 2003, and May 2003}
Oslo, Iowa City and Kiev, June 1999, January 2003, and May 2003
\linebreak
\makebox[\displayboxwidth]{Ola Bratteli\hfill Palle E. T.
Jorgensen\hfill Vasyl' Ostrovs'ky\u\i}
\end{center}

\chapter*{Introduction}

\renewcommand{\thesubsection}{\arabic{subsection}}
The origin of this monograph was the
desire to understand certain concrete operator
relations arising as realizations
of filtering processes in signal theory.
Our research did, however, lead us naturally in
the direction of
analyzing certain noncommutative
dynamical systems and their fixed
points and states. In this introduction, we give an
overview of the contents of the monograph
in three stages: First (\ref{IntGen})~a very general
discussion, then (\ref{IntWay})~a discussion with
more specific definitions and details,
and finally (\ref{IntDet})~a detailed technical account
of what the paper actually contains.

\subsection{\label{IntGen}General discussion and motivation}

This monograph is centered around the issue of
distinguishing a particular family of AF-algebras,
those which arise as the
centralizers of certain states on Cuntz algebras---or
equivalently, as
the fixed-point algebras
under certain one-parameter subgroups of the gauge action. (The
Cuntz algebras are the range of a functor from Hilbert
space into $C^*$-algebras. The term AF-algebra\label{LOSAFalgebra_0}
is short for approximately
finite-dimensional $C^*$-algebra. Both the Cuntz algebras and the AF-algebras play a
role in several areas of mathematics, e.g., $K$-theory and dynamical systems, and in
applications, for example to statistical mechanics
\cite{BrRoII}.) While AF-algebras generally are classified 
up to stable
isomorphism by the equivalence classes of their
Bratteli diagrams, or by the isomorphism classes of their ordered 
dimension groups,
none
of the three items in this triple of equivalence 
classes, of AF-algebras,
Bratteli diagrams, or ordered dimension groups, respectively,
is
especially amenable to
computation. In this memoir, we therefore try to approach the subject via
classes of concrete examples, which do in fact admit algorithms for 
distinguishing
isomorphism classes, in particular by the
computation of numerical invariants which distinguish these classes. The
examples are chosen so they illustrate in a concrete manner the main issues of
computation in each of the three incarnations, AF-algebras, Bratteli diagrams,
or dimension groups.
The terms in this paragraph will be defined
in Section \ref{IntWay} of this introduction.

The first part of our memoir is in a sense divorced from the central theme,
that of numerical invariants. But we feel that it helps the reader see how the
main issue fits into a wider context, and hopefully it will help the reader
make connections to the areas of mathematics and its applications 
which are touched
upon in a more indirect manner. The ubiquitous Cuntz algebras and their
representations have a surprising number of applications in a wide range of
areas, also outside the subject of operator algebras, such as wavelets and
dynamical systems, to mention just two. One of the recent areas of application
of dimension groups and Perron--Frobenius structures is to notions of
equivalence for symbolic dynamics, and more specifically subshifts.
Section 7.5 in the book
\cite{LiMa95}
serves as an excellent background reference. In particular, it describes how
substitution matrices and dimension groups\label{LOSdimensiongroup_0}
are derived from symbolic systems
in dynamics, and it gives the basics on some of the numerical invariants used
there. The matrices that we encounter in the context of AF-$C^*$-algebras
are
a subclass of those used for the general dynamics problem, and our notion
of equivalence is weaker. As we show in \cite{BJKR00}, our
notion of $C^*$-equivalence for two dimension groups corresponds to a weaker
concept of equivalence than the ones of strong shift equivalence and shift
equivalence. The latter are nicely surveyed in
\cite{Wag99},
which also recounts the history of the Williams problem, and its
resolution. 

Let us give a short rundown of this problem.
In the rest of this paragraph, let the
term `matrix' mean `matrix over the nonnegative integers'.
Then two square matrices $J$, $K$
are \emph{elementary shift equivalent} if there exist
matrices $A$, $B$ such that $J=BA$ and
$K=AB$. We say that $J$, $K$ are \emph{shift equivalent of lag} $k$ if
$AJ=KA$, $BK=JB$, $BA=J^k$, and $AB=K^k$ hold for some matrices $A$, $B$, and
they are \emph{shift equivalent} if they are shift equivalent
of some lag $k$ in $\mathbb{N}$. Thus shift equivalence
of lag $1$ is the same as elementary shift equivalence
\cite[Proposition 7.3.2]{LiMa95}.
The matrices $J$, $K$ are \emph{strongly shift equivalent} if they can
be connected by a finite chain of elementary shift equivalent 
matrices. (Note that elementary shift equivalence is not an 
equivalence relation; it is not transitive.) Strong shift equivalence
trivially implies shift equivalence \cite[Theorem 7.3.3]{LiMa95},
but the converse is the longstanding Williams conjecture, and it is not
true \cite{KiRo92}, even when $J$ and $K$ are irreducible
\cite{KiRo99}.
The conjecture was stated in \cite{Wil73}.
The properly weaker notion of $C^*$-equivalence or weak equivalence
of two square matrices $A$, $B$ is defined by the existence of
a sequence  $A_1$, $B_1$, $A_2$, $B_2$, $\dots$ of matrices
and sequences $\left( n_k\right) $, $\left( m_k\right) $
of natural numbers such that
\begin{align*}
J^{n_k}&=B_kA_k\\
K^{m_k}&=A_{k+1}B_k
\end{align*}
for $k=1,2,\dots$.

One notes that the numerical invariants that are known in the
subject have a history based on computation, and that applies to
the issue of $C^*$-equivalence, also called weak equivalence, as well. And the
invariants for weak equivalence are different from those used in \cite{Wag99}
and the papers cited there. The latter are based on sign-gyration conditions,
cohomology, and the $K$-theory group $K_2$: see \cite[Section 1]{Wag99}.
A part of the paper \cite{Wag99}
further describes a class of one-dimensional systems for which equivalence up
to homeomorphism is decided by weak equivalence ($C^*$-equivalence) for the
corresponding dimension groups.

Chapter \ref{Gen} is a discussion of the spectral decomposition of representations
of the Cuntz algebras relative to canonical MASAs (maximal abelian
subalgebras). It serves as a pointer to how the particular AF-algebras in
Chapter \ref{Sub},
which form the basis for our analysis, arise. But the dimension groups of the
AF-algebras which result as fixed-point subalgebras from our
operator-theoretic construct
turn out also to come up in other areas of mathematics, and
readers primarily interested in dynamics may wish to consult
\cite{LiMa95}, or
\cite{BoHa93}
for this alternative use of the same special
dimension groups. These problems have a flavor of algebraic number theory; see
also
\cite{BJKR01}.
Chapter \ref{Fre} specialises the general setup of
Chapter \ref{Gen} to a certain class of
representations of the Cuntz algebras, and
the chapter
is motivated partly by the
representations arising from wavelet theory. These chapters describe basic
theory of Cuntz-relation representations
which has also been considered in \cite[Section 2.5.3]{OsSa99}.
Chapter \ref{KMS} is devoted to so-called
KMS states (short for
Kubo--Martin--Schwinger states) from statistical physics. More specifically,
we look at KMS states of the dynamic defined by
one-parameter groups of quasi-free automorphisms on a Cuntz $C^*$-algebra. One
theorem gives a characterization of when there is, or is not, such a
KMS state. Chapter \ref{Sub} determines when,
in this situation,
the fixed-point algebra
is AF, and gives an outline of how to go about computing the corresponding
Bratteli diagrams. Since a main aim
of the memoir
is to make an otherwise technical subject
concrete, many examples are offered to help the reader appreciate the abstract
theorems, and to motivate the numerical invariants that are introduced later.
The rest of the memoir is devoted to the analysis of
the isomorphism classes of these AF-algebras.
Thus the initial chapters are mainly devoted to setting up the class
of AF-algebras
to be considered. But we feel that they are of independent interest, and have
relevance outside the particular issue of AF-algebra classification.
Part \ref{PartRep}
is
representation-theoretic, while Part \ref{PartInv} describes numerical invariants
for the AF-algebra
classification context, and associated Bratteli diagrams.
There are discussions of special cases of the general problem which are
related to:
(i)~the Williams problem (the issue of understanding the two main
ways of classifying dimension groups in dynamics---see \cite{Wag99} for more on
this connection); (ii)~Krieger's theorem and conjugacy of the corresponding
actions on the Cuntz algebras (see \cite{LiMa95}); (iii)~invariants related to
the Perron--Frobenius eigenvalue; and (iv)~more algebraic invariants with
examples, or situations where the invariants become complete.

The
one-parameter groups of automorphisms that we consider
act upon the Cuntz algebras, and they are called gauge actions; they
scale the generators with a unitary gauge. Our analysis is focused on the
corresponding fixed-point algebras: their classification, and their
significance in, for example, symbolic dynamics.
But there are analogous and interesting questions for more general $C^*$-algebras
which arise in dynamical systems, for instance the Cuntz--Krieger algebras $\mathcal{O}_A$.
These are algebras, naturally generalizing the Cuntz algebras, defined on
generators $S_i$ and relations, the relations
being given by a $0$-$1$ matrix $A$. The matrix $A$ in turn determines, in a standard
fashion, see \cite{LiMa95}, a state space $S_A$ for a dynamical system
called a subshift. While we will not treat them here, our methods suggest
directions for future work on other gauge actions, as illustrated by recent
papers of Exel, Laca, and Vershik, for example
\cite{ExVe02}
and
\cite{ExLa99,ExLa00}.
A gauge action on an algebra
scales the generators by a unitary phase, and a recent paper
\cite{Exe01}
defines the notions of gauge action and KMS states
in the setting of Cuntz--Krieger algebras. While our present analysis for the
Cuntz algebras is relatively managable for computations, involving only
Perron--Frobenius theory for finite matrices, the analogous analysis of Exel if
carried through would instead be based on an infinite-dimensional version of
the Perron--Frobenius theorem, due to Ruelle
(see \cite{Rue89} and \cite{Rue02}).
In Ruelle's theorem, the
finite matrix from the classical Perron--Frobenius theorem is replaced by an 
operator, now called the Ruelle transfer operator
(see \cite{Jor01a}, \cite[Section 3.1]{BrJo02b}),
and
the Perron--Frobenius left eigenvector from the classical case, by a measure on 
the state space $S_A$. The Ruelle transfer operator acts on functions on the 
compact space $S_A$, and the right Perron--Frobenius eigenvector is 
a continuous function on $S_A$.  Following the model from the simpler case of 
the Cuntz algebra, and using the Ruelle--Perron--Frobenius theorem, Exel in 
\cite{Exe01} computes the KMS states for the
Cuntz--Krieger algebras, and it would be interesting to analyze the issues of
fixed-point subalgebras and state centralizers in this context. Relevant for
our present analysis is the fact that unital endomorphisms of
$\mathcal{B}(\mathcal{H})$
are known
to correspond to representations of the Cuntz algebras
$\mathcal{O}_d$,
and the integer $d$
is the Powers index. The corresponding index in the setting of
\cite{Exe01} is not
necessarily integer-valued, but some methods in our memoir point to
generalizations and future directions of research in this and other algebraic
settings of symbolic dynamics.

The use of certain representations of $\mathcal{O}_d$ in
wavelet theory is outlined in \cite[Exercises 1--11, 1--12, and 2--25]{BrJo02b}.
The starting observations there are the facts
(i)~that the system of frequency-subband filters
with $d$ subbands is known to satisfy
the relations which define representations of
$\mathcal{O}_d$, and
(ii)~that the multiresolution approach
to wavelets is based on subband filters.

Yet other uses of the representations
of the Cuntz algebras $\mathcal{O}_d$ to dynamics include
the work of J.~Ball and V.~Vinnikov
\cite{BaVi02a,BaVi02b}
on Lax--Phillips scattering theory for multivariable
linear systems. These two
papers develop a functional model
for representations of $\mathcal{O}_d$
in a setting which is more special
than that of Theorem \ref{ThmGen.1} below.

The original AF-algebra\label{LOSAFalgebra_00}
in \cite{Bra72} had Pascal's triangle as its 
diagram and arose as the fixed-point
algebra of the standard infinite-product-type action of the circle on 
the infinite tensor product
of a sequence of full $2 \times 2$ matrix algebras.
Infinite-tensor-product-type actions of more general
compact groups have been considered by Anthony Wassermann, Adrian 
Ocneanu and others,
and the fixed-point algebras are related to the AF-algebras occuring 
in Vaughan Jones´s
subfactor theory,
which again has important applications in knot theory, 
low-dimensional manifolds, conformal
field theory, etc.; see \cite{EvKa98,GHJ89}. But this 
is far outside the scope of
the present monograph.

As noted in \cite{MaRo01}, the list of
applications of Bratteli diagrams includes such areas of applied
mathematics as the fast Fourier transform, and more general
Fourier transforms on finite or infinite nonabelian groups. These transforms
in turn are widely used in computational mathematics, in
error-correction codes, and even in P.~Shor's fast algorithm for
integer factorization on a quantum computer \cite{Sho97,Sho99}.
As reviewed in \cite{MaRo01}, the effectiveness of the
Bratteli diagrams for development of algorithms pays off
especially well when a sequence of operations is involved, with
a succession of steps, perhaps with a scaling similarity,
and with embeddings from one step into the next.
In the case of matrix algebras, the embeddings are in the
category of rings, while for groups, the embeddings are
group homomorphisms. The number of embeddings
may be infinite, as it is in the applications in the
present memoir, or it may be finite, as in the case
of the fast-Fourier-transform algorithms of Cooley and Tukey,
see \cite{MaRo01}.

The discrete Fourier transform for the cyclic group
$\mathbb{Z}\diagup n\mathbb{Z}$
involves the unitary matrix
$\frac{1}{\sqrt{n}}\left( e^{i2\pi jk/n}\right) _{j,k=0}^{n-1}$
and is of
complexity $O(n^2)$. If $n$ is a composite number, such as
$n=2^k$, then there are choices of Bratteli diagrams built
on $\mathbb{Z}\diagup n\mathbb{Z}$ which turn the
discrete Fourier transform into a
fast Fourier transform of complexity $O(n\log n)$. For each
of the Bratteli diagrams formed from a sequence of subgroups
of $\mathbb{Z}\diagup n\mathbb{Z}$, the corresponding
fast Fourier transform involves a summation
over the paths of the Bratteli diagram.

The subject of this memoir is developing into a long-term programme by
the authors, also in collaboration
with Kim and Roush
\cite{BJKR00,BJKR01,BJKR02}. The results
are often
somewhat technical.
The paper
\cite{Wag99}
provides
a beautiful survey of the Williams
problem, its history, and its resolution. The subject of equivalences, strong,
shift equivalence, or weak equivalence, is situated neither directly in
operator algebras, nor perhaps precisely in dynamical systems or wavelet
analysis, but rather intermingles these, and more, areas of mathematics. 

\addvspace{\bigskipamount}
Let us now add more specific definitions and go into more details:

\subsection{\label{IntWay}Three ways of measuring isomorphism}

During the sixties and seventies it was established that there is a one-to-one
canonical correspondence between the following three sets
\cite{Bra72,Dix67,Eff81,EHS80,Ell76,Gli60}:

\begin{enumerate}
\item \label{PrefSet1}the isomorphism classes of AF-algebras,\medskip

\item \label{PrefSet2}the isomorphism classes of certain ordered abelian
groups, called dimension groups,
\end{enumerate}

\noindent and finally

\begin{enumerate}
\setcounter{enumi}{2}

\item \label{PrefSet3}the equivalence classes of certain combinatorial
objects, called Bratteli diagrams.
\end{enumerate}

\noindent In more recent times, this has led to an undercurrent of
misunderstanding that AF-algebras, which are complex objects, are classified
by dimension groups, which are easy objects, and that this is the end of the
story. However, as anyone who has worked with these matters knows, although
for special subsets it may be easier to work with one of the three sets
mentioned above rather than another, in general the computation of isomorphism
classes in any of the three categories is equally difficult. Although
dimension groups are easy objects, their isomorphism classes in general are
not! One may even be tempted to flip the coin around and say that dimension
groups are classified by AF-algebras. If one thinks about isomorphism classes,
this is logically true, but the only completely general method to decide
isomorphism classes in all the cases is to resort to the computation of the
equivalence relation for the associated Bratteli diagrams. This problem is not
only hard in general, it is even undecidable: There is no general recursive
algorithm to decide if two effective presentations of Bratteli diagrams yield
equivalent diagrams \cite{MuPa98}. In this memoir, we will encounter this
problem in a very special situation, and try to resolve it in a modest way by
introducing various numerical invariants which are easily computable from the
diagram. In the situation that the AF-embeddings are given by a constant
primitive nonsingular matrix, the classification problem
has, after the writing of this monograph, been proved in general to
be decidable
\cite{BJKR00,BJKR01,BJKR02}.

Recall that an \emph{AF-algebra}\label{LOSAFalgebra_1} is a separable
$C^{\ast}$-algebra with the property that for any $\varepsilon>0$, any finite
subset of the algebra can be approximated with elements of some
finite-dimensional $\ast$-subalgebra with the precision given by $\varepsilon
$. An AF-algebra is \emph{stable} if it is isomorphic to its tensor product
with the compact operators on a separable Hilbert space. A \emph{dimension
group}\label{LOSdimensiongroup_1} is a countable abelian group with an order
satisfying the Riesz interpolation property and which is unperforated. The
\emph{Bratteli diagram}\label{LOSBrattelidiagrams_1} is described in
\cite{Bra72}, \cite{Dav96}, and \cite{EvKa98}, and the equivalence relation is
also described in detail in \cite{Bra99} and in Remark \ref{RemRemJun.6}. (All
these concepts will be treated in some detail in Chapter \ref{Rem}, where it
is also explained that the stability assumption on the AF-algebra can be
removed by putting more structure on the group and the diagram.)

Recently there has been a fruitful interaction between the theory of dynamical
systems, analytic number theory, and $C^{\ast}$-algebras. In \cite{BoCo95},
the authors show how $\beta$-KMS states may be used in understanding the
Riemann zeta function, and vice versa. In \cite{UeWa00}, \cite{Cho87},
\cite{EvGo94}, \cite{GHJ89}, \cite{EvKi97}, \cite{HPS92}, \cite{KMW98}, and
\cite{KiKu97}, certain dynamical systems are used to generate new simple
$C^{\ast}$-algebras from the Cuntz algebras, and to better understand the
corresponding isomorphism classes of $C^{\ast}$-algebras. The results in
Chapter \ref{Sub} should be contrasted with results of Izumi \cite{Izu93} and
Watatani \cite{Wat90} which deal with crossed product constructions built from
the Cuntz algebras $\mathcal{O}_{d}$.\label{LOSOd_1} Here we study the
AF-subalgebras of $\mathcal{O}_{d}$ formed from the one-parameter automorphism
groups of Chapter \ref{KMS}.

It follows from the definition of the $\mathcal{O}_{d}$-relations that they
are well adapted to $d$-multiresolutions of the kind used in wavelets and
fractal analysis. The number $d$ represents the scaling factor of the wavelet.
This viewpoint was exploited in recent papers \cite{BJP96}, \cite{BrJo97a},
\cite{BrJo97b}, and \cite{DaLa98}. While the representations for these
applications are type $\mathrm{I}$, the focus in the present memoir is type
$\mathrm{III}$ representations of $\mathcal{O}_{d}$, and a family of
associated AF-$C^{\ast}$-algebras $\mathfrak{A}_{L}$\label{LOSAL_1}
($\subset\mathcal{O}_{d}$ for some $d$). These representations arise from a
modified version of the technique which we used in generating wavelet
representations. This starting point in fact yields a general decomposition
result for representations of $\mathcal{O}_{d}$ which seems to be of
independent interest. To understand better the resulting decomposition
structure, we will establish that the centralizers of these states are simple
AF-algebras,\label{LOSAFalgebra_2} and that the Bratteli diagrams have
stationary incidence matrices $J$ of a special form given in (\ref{eqCyc.2}).
(In general the centralizer of a state $\omega$ on a $C^*$-algebra
$\mathfrak{A}$ is defined as the set of $x\in\mathfrak{A}$ such that
$\omega\left( xy\right) =\omega\left( yx\right) $ for all $y\in\mathfrak{A}$.
If $\omega$ is a KMS state with
respect to a dynamical one-parameter group $\sigma$, as in (\ref{eqKMS.3}),
then this is equivalent to
saying that $x$ is a fixed point for the dynamic $\sigma$.)
Clearly the rank of the corresponding dimension group is an invariant, but it
appeared at first sight that different matrices $J$ and $J^{\prime}$ would
yield non-isomorphic AF-algebras $\mathfrak{A}_{J}$ and $\mathfrak
{A}_{J^{\prime}}$. This turns out not to be the case, however, and the bulk of
the memoir concerns numerical AF-invariants. It is not easy to get invariants
that discriminate the most natural cases of algebras $\mathfrak
{A}_{J}$ that arise from this seemingly easy family of AF-algebras. There is a
connection to subshifts in dynamical systems, but if subshifts are
constructed, from $J$ and $J^{\prime}$, say, we note that strong shift
equivalence in the sense of Krieger is equivalent to $J=J^{\prime}$, while
isomorphism of the AF-algebras $\mathfrak{A}_{J}$ and $\mathfrak{A}%
_{J^{\prime}}$ turns out to be a much more delicate problem; see also
\cite{BJKR00}. While we do not have a complete set of numerical AF-invariants
for our algebras $\mathfrak{A}_{J}$, we do find interesting subfamilies of
$\mathfrak{A}_{J}$'s which do in fact admit a concrete isomorphism labeling,
and Part \ref{PartInv} of the memoir concentrates on these cases. In contrast,
we mention that
\cite{BJKR01} and \cite{BJKR02} do
have general criteria for $C^{\ast}%
$-isomorphism of the algebras $\mathfrak{A}_{J}$, but those conditions are
rather abstract in comparison with the explicit and numerical invariants which
are the focus of the present memoir.
A complete algorithm for $C^*$-isomorphism is written out
in \cite[pp.\ 1651--1653 and corrigendum]{BJKR01}.

\addvspace{\bigskipamount}
Let us now go into
even
more
nitty-gritty technical
details:

\subsection{\label{IntDet}Detailed outline}

The Cuntz algebra
$\mathcal{O}_{d}$\label{CuntzalgebraIntrobis} with generators $s_{i}%
$\label{LOSsi} and relations $s_{i}^{\ast}s_{j}^{{}}=\delta_{ij}1$,
$\sum_{i=1}^{d}s_{i}^{{}}s_{i}^{\ast}=1$ contains a natural abelian subalgebra
$\mathcal{D}_{d}\cong C\left(  \prod_{1}^{\infty}\mathbb{Z}_{d}\right)
$,\label{LOSDd_1} $\mathbb{Z}_{d}=\left\{  1,\dots,d\right\}  $\label{LOSZd}
(see \cite{Cun77}). We relate general representations of $\mathcal{O}_{d}$
with the spectral resolution of the restriction of the representation to
$\mathcal{D}_{d}$. From this, we read off cocycle formulations of the factor
property, irreducibility, and of equivalence of representations. We then
specialize to the representations associated with the GNS construction from
states $\omega=\omega_{p}$\label{LOSomega} on $\mathcal{O}_{d}$ indexed by
$p_{i}>0$, $\sum_{i=1}^{d}p_{i}=1$, given by
\[
\omega\left(  s_{\alpha_{1}}^{{}}\cdots s_{\alpha_{k}}^{{}}s_{\gamma_{l}%
}^{\ast}\cdots s_{\gamma_{1}}^{\ast}\right)  =\delta_{kl_{{}}}\delta
_{\alpha_{1}\gamma_{1}}\cdots\delta_{\alpha_{k}\gamma_{k}}p_{\alpha_{1}}\cdots
p_{\alpha_{k}}%
\]
where $\alpha$ and $\gamma$ are multi-indices formed from $\mathbb{Z}_{d}$.
The cyclic space $\mathcal{H}_{\omega_{p}}$,\label{LOSHomegap} for $\left(
p_{i}\right)  $ fixed, is shown in Chapter \ref{Fre} to have a bundle
structure over the set of all $\mathbb{Z}_{d\,}$-multi-indices with fiber
$\ell^{2}\left(  \mathbb{S}_{d}\right)  $ where $\mathbb{S}_{d}$ is the free
semigroup on $d$ generators.

Let $L=\left(  L_{1},\dots,L_{d}\right)  \in\mathbb{R}^{d}$,\label{LOSL} and
consider the one-parameter group $\sigma^{L}$ of $\ast$-au\-to\-mor\-phisms of
$\mathcal{O}_{d}$ defined by%
\[
\sigma_{t}^{L}\left(  s_{j}\right)  =\exp\left(  itL_{j}\right)
s_{j}.\label{LOSsigmatL_1}
\]
It can be shown (\cite{Eva80}, or Proposition \ref{ProKMS.1}) that $\sigma
^{L}$ admits a $\sigma^{L}$-KMS state, at some value $\beta$, if and only if
all $L_{j}$'s are nonzero and have the same sign. This value $\beta$ is then
unique and is defined as the solution of%
\[
\sum_{k=1}^{d}e^{-\beta L_{k}}=1,
\]
and the $\left(  \sigma^{L},\beta\right)  $-KMS
state\label{LOSsigmaLbetaKMSstate_1} is then also unique, and is the state
defined in the previous paragraph with $p_{k}=e^{-\beta L_{k}}$%
.\label{LOSpj_1} Note that the group $\sigma^{L}$ is periodic if and only if
any pair $L_{j}$, $L_{k}$ is rationally dependent. In that case, let
$\mathfrak{A}_{L}$ be the fixed-point subalgebra in $\mathcal{O}_{d}$ under
the action $\sigma^{L}$. We show in Chapter \ref{Sub} that $\mathfrak{A}_{L}$
is an AF-algebra (\cite{Bra72}, \cite{Eff81}) if and only if all $L_{k}$'s
have the same sign, and furthermore the $\mathfrak{A}_{L}$'s are then simple
with a unique trace state (namely the restriction of the state in the previous
paragraph to $\mathfrak{A}_{L}$). We compute the Bratteli diagrams of the
$\mathfrak{A}_{L}$ in this case, and show (using a result from \cite{Con77})
that the endomorphism $\rho\left(  a\right)  :=\sum_{i=1}^{d}s_{i}^{{}}%
as_{i}^{\ast}$ restricts to a shift (in the sense of Powers \cite{Pow88}) on
each of the algebras $\mathfrak{A}_{L}$, i.e., $\bigcap_{k=1}^{\infty}\rho
^{k}\left(  \mathfrak{A}_{L}\right)  =\mathbb{C}1$.

While the dimension group $D\left(  \mathfrak{A}_{L}\right)  $, described in
\cite{Bra72}, \cite{Ell76}, \cite{Eff81}, and \cite{Dav96} in principle is a
complete AF-invariant, we have mentioned that its structure is not immediately
transparent. For the present AF-algebras $\mathfrak{A}_{L}$, the
classification is facilitated by the display of specific numerical invariants,
derived from $D\left(  \mathfrak{A}_{L}\right)  $, but at the same time
computable directly in terms of the given data $\left(  L_{1},\dots
,L_{d}\right)  $. These invariants are described in Chapters \ref{Brunt}%
--\ref{ClmN} where their connection to $\operatorname*{Ext}$ is partially explained.

Let us give a short road map to the various invariants introduced and where it
proved that they are invariants (sometimes in restricted settings):
$K_{0}\left(  \mathfrak{A}_{L}\right)  $\label{LOSK0AL_1} in
(\ref{eqRemNewBis.6}), (\ref{eqRemNewBis.19}); $\tau\left(  K_{0}\left(
\mathfrak{A}_{L}\right)  \right)  $ in (\ref{eqRemNewBis.22}); $D\left(
\mathfrak{A}_{L}\right)  =\left(  K_{0}\left(  \mathfrak{A}_{L}\right)
,K_{0}\left(  \mathfrak{A}_{L}\right)  _{+},\left[  \openone\right]  \right)
$\label{LOSDAL_1} in (\ref{eqRemNewBis.30}); $\ker\tau$\label{LOSkertau_1} in
(\ref{eqRemNewBis.31}); $\mathbb{Q}\left[  \lambda\right]  $ together with the
prime factors of $\lambda$ before (\ref{eqBrunt.1}); $N$, $D$,
$\operatorname*{Prim}\left(  m_{N}\right)  $,\label{LOSPrim_1}
$\operatorname*{Prim}\left(  Q_{N-D}\right)  $, $\operatorname*{Prim}\left(
R_{D}\right)  $ in Theorem \ref{CorCyc.9}; $K_{0}\left(  \mathfrak{A}%
_{L}\right)  \otimes_{\mathbb{Z}}\mathbb{Z}_{n}$ and $\ker\tau\otimes
_{\mathbb{Z}}\mathbb{Z}_{n}$ in Chapter \ref{Ker}, $M$ in (\ref{eqKer.p});
$\operatorname*{rank}\mathcal{L}_{k}$\label{LOSLk_1} in Corollary
\ref{CorCyc.plus2}; class in $\operatorname*{Ext}\left(  \tau\left(
K_{0}\left(  \mathfrak{A}_{L}\right)  \right)  ,\ker\tau\right)  $ in Chapter
\ref{Ext}; $D_{\lambda}\left(  K_{0}\left(  \mathfrak{A}_{L}\right)  \right)
$ in (\ref{eq10.55})--(\ref{eq10.56}); $I\left(  J\right)  $\label{LOSIJ_1} in
(\ref{eqClmN.L}) and Corollary \ref{CorClmN.E1}. In general it is very hard to
find complete invariants apart from $D\left(  \mathfrak{A}_{L}\right)  $, even
for special subclasses; but if the Perron--Frobenius eigenvalue $\lambda
$\label{LOSPerronFrobeniuseigenvalue_1} of $J$ is rational (and thus integral)
and $N=2$ and $\operatorname*{Prim}\left(  \lambda\right)
=\operatorname*{Prim}\left(  m_{2}/\lambda\right)  $, then
$\operatorname*{Prim}\left(  \lambda\right)  $ is a complete invariant by
Proposition \ref{ProApp.Cla2}. The same is true if the condition
$\operatorname*{Prim}\left(  \lambda\right)  =\operatorname*{Prim}\left(
m_{2}/\lambda\right)  $ is replaced by $\operatorname*{Prim}\left(
\lambda+\frac{m_{2\mathstrut}}{\lambda}\right)  \subseteq\operatorname*{Prim}%
\left(  \lambda\right)  $ by Proposition \ref{ProApp.ClaNew3}.

In Chapter \ref{APP.EXA} we give a complete classification of the class
$\lambda=2$, $N=2,3,4$. This class contains $28$ specimens, and it turns out
that all of them are non-isomorphic except for a subset consisting of the
three specimens in Figure \ref{BratDiagsN4g2s1}.

The most striking classification result for a restricted, but infinite, class
of examples in this memoir is that if $\lambda=m_{N}$, then $\left(
N,\operatorname*{Prim}\lambda,I\left(  J\right)  \right)  $ is a
\emph{complete} invariant. This is proved in Theorem \ref{ThmClmNJun.18}.

In addition to these formal invariants there are very efficient methods to
decide non-isomorphism when $\lambda$ is rational based on a quantity
$\tau\left(  v\right)  =%
\ip{\alpha}{v}%
$\label{LOStauv_1} defined in (\ref{eqajla})--(\ref{eqjvlv}); see Theorem
\ref{theorem10.10}, Remark \ref{RemSesquilabialNew.11}, Corollaries
\ref{corollary10.11}--\ref{corollary10.12}, Scholium \ref{scholium10.21}. In
fact $I\left(  J\right)  =\lambda^{N-1}%
\ip{\alpha}{v}%
=\lambda^{N-1}$ times the inner product of the right and left
Perron--Frobenius eigenvectors $\alpha$%
,\label{LOSPerronFrobeniuslefteigenvector_1} $v$%
\label{LOSPerronFrobeniusrighteigenvector_1} of $J$, normalized so that
$\alpha_{1}=v_{1}=1$.

In
two new papers
\cite{BJKR01,BJKR02}
written jointly
with K.H. Kim and F. Roush,
we give a completely general
algorithm to decide $C^*$-equivalence of two square
matrices $A$ and $B$.
In its most general form, however, this general
algorithm is not merely based on the computation
of a complete set of invariants, and the
methods of the present monograph are still
most useful for the particular examples
arising from the noncommutative
dynamical systems we are considering.

\clearpage\thispagestyle{empty}\mainmatter

\part{\label{PartRep}Representation Theory}

\ \markboth{\shorttitle}{\shorttitle} \cleardoublepage

\chapter{\label{Gen}General representations of $\mathcal{O}_{d}$ on a
separable Hilbert space}

Representations of the Cuntz algebras $\mathcal{O}_{d}$\label{LOSOd_2} play a
role in several recent papers; see, e.g., \cite{EvGo94}, \cite{BJP96},
\cite{Cho87}, \cite{Cho93}, \cite{Cho95}, \cite{DoRo87}, \cite{Wat90},
\cite{Izu93}. Since $\mathcal{O}_{d}$ is purely infinite, there are few
results that cover all representations. The following result does just that,
and serves as a ``noncommutative spectral resolution''. We will use the
convention that $S_{i}$ denotes the representative of the Cuntz algebra
generator $s_{i}$ in any given representation.

Let $\Omega=\prod_{1}^{\infty}\mathbb{Z}_{d}$,\label{LOSOmega_1} and let
$\sigma$ be the right shift on $\Omega$:\label{LOSsigma}%
\begin{equation}
\sigma\left(  x_{1},x_{2},\dots\right)  =\left(  x_{2},x_{3},\dots\right)
\label{eqGen.1}%
\end{equation}
Define sections $\sigma_{i}$ of $\sigma$ by\label{LOSsigmai}%
\begin{equation}
\sigma_{i}\left(  x_{1},x_{2},\dots\right)  =\left(  i,x_{1},x_{2}%
,\dots\right)  , \label{eqGen.8}%
\end{equation}
for $i=1,\dots,d$. Note that $\sigma$ is a $d$-to-$1$ map and that the
sections are injective. The sets\label{LOSOmegai}%
\begin{equation}
\Omega_{i}=\sigma_{i}\left(  \Omega\right)  \label{eqGen.23}%
\end{equation}
form a partition of the Cantor set $\Omega$ into clopen sets. The $\sigma_{i}$
are right inverses of $\sigma$:%
\begin{equation}
\sigma\sigma_{i}=\operatorname*{id} \label{eqGen.24}%
\end{equation}
for $i=1,\dots,d$. If $\mu$ is a probability measure on $\Omega$, we say that
$\mu$ is $\sigma$-quasi-invariant if%
\begin{equation}
\mu\left(  E\right)  =0\Longrightarrow\mu\left(  \sigma^{-1}\left(  E\right)
\right)  =0 \label{eqGenJun.Y}%
\end{equation}
for all Borel sets $E\subseteq\Omega$, and we say that $\mu$ is $\sigma_{i}%
$-quasi-invariant if%
\begin{equation}
\mu\left(  E\right)  =0\Longrightarrow\mu\left(  \sigma_{i}^{-1}\left(
E\right)  \right)  =0 \label{eqGenJun.Z}%
\end{equation}
for all Borel sets $E\subseteq\Omega$, where
\begin{equation}%
\begin{aligned}
\sigma_{i}^{-1}\left( E\right) &=\left\{ x\mid\sigma_{i}\left( x\right
) =\left( i,x_{1},x_{2},\dots\right) \in E\right\} , \\
\sigma^{-1}\left( E\right) &=\left\{ x\mid\sigma\left( x\right) =\left
( x_{2},x_{3},\dots\right) \in E\right\} .
\end{aligned}
\label{eqGenJun.uther}%
\end{equation}
The set of $d$ conditions (\ref{eqGenJun.Z}) is equivalent to the set of $d$
conditions
\begin{equation}
\mu\left(  \sigma_{i}\left(  F\right)  \right)  \Longrightarrow\mu\left(
F\right)  =0 \label{eqGenJun.onemore}%
\end{equation}
for $i=1,\dots,d$. The condition (\ref{eqGenJun.Y}) is implied by, but does
not imply, the condition%
\begin{equation}
\mu\left(  \sigma\left(  F\right)  \right)  =0\Longrightarrow\mu\left(
F\right)  =0. \label{eqGenJun.onemorebis}%
\end{equation}
(Conditions like (\ref{eqGenJun.onemore}) and (\ref{eqGenJun.onemorebis}) make
sense in this setting since $\sigma_{i}$ and $\sigma$ are local
homeomorphisms, and thus map measurable sets into measurable sets.) Note that%
\begin{equation}
\sigma^{-1}\left(  E\right)  =\bigcup_{i}\sigma_{i}\left(  E\right)
\label{eqGenJun.okseskaft}%
\end{equation}
for all sets $E$ by the chain: $x\in\bigcup_{i}\sigma_{i}\left(  E\right)
\iff x=iy$ for some $i=1,\dots,d$, $y\in E\iff\sigma\left(  x\right)  \in
E\iff x\in\sigma^{-1}\left(  E\right)  $. Note also that if $\mu$ is both
$\sigma$- and $\sigma_{i}$-quasi-invariant, we have the connection%
\[
\frac{d\mu\left(  \sigma_{i}\left(  y\right)  \right)  }{d\mu\left(  y\right)
}=\frac{d\mu\left(  \sigma_{i}\left(  y\right)  \right)  }{d\mu\left(
\sigma\sigma_{i}\left(  y\right)  \right)  }=\left(  \frac{d\mu\circ\sigma
}{d\mu}\right)  ^{-1}\left(  \sigma_{i}\left(  y\right)  \right)
\]
between the Radon--Nikodym derivatives.\medskip

Note also that the two quasi-invariance conditions (\ref{eqGenJun.Y}) and
(\ref{eqGenJun.Z}) together imply the $d$ equivalences%
\begin{equation}
\mu\left(  E\right)  =0\iff\mu\left(  \sigma_{i}E\right)  =0,\qquad
i=1,\dots,d, \label{eqGenJun.9}%
\end{equation}
for all Borel sets $E\subseteq\Omega$, and that (\ref{eqGenJun.9}) implies the
$\sigma$-invariance (\ref{eqGenJun.Y}). (When referring to (\ref{eqGenJun.Z})
(or (\ref{eqGenJun.9})) in the following we mean ``(\ref{eqGenJun.Z}) (or
(\ref{eqGenJun.9})) for all $i=1,\dots,d\,$''). Let us prove this.\medskip

\underline{(\ref{eqGenJun.Y}) \& (\ref{eqGenJun.Z}) $\Rightarrow$
(\ref{eqGenJun.9})}: If (\ref{eqGenJun.Y}) and (\ref{eqGenJun.Z}) hold and
$\mu\left(  E\right)  =0$, it follows from (\ref{eqGenJun.Y}) that $\mu\left(
\sigma^{-1}\left(  E\right)  \right)  =0$ and hence from
(\ref{eqGenJun.okseskaft}) that $\mu\left(  \sigma_{i}\left(  E\right)
\right)  =0$ for all $i$. Conversely, if $\mu\left(  \sigma_{i}\left(
E\right)  \right)  =0$ for some $i$, then since $E=\sigma_{i}^{-1}\sigma
_{i}^{{}}E$ it follows from (\ref{eqGenJun.Z}) that $\mu\left(  E\right)  =0$.\medskip

\underline{(\ref{eqGenJun.9}) $\Rightarrow$ (\ref{eqGenJun.Y})}: Assume that
(\ref{eqGenJun.9}) holds and that $\mu\left(  E\right)  =0$. Then $\mu\left(
\sigma_{i}\left(  E\right)  \right)  =0$ for all $i$ by (\ref{eqGenJun.9}) and
hence $\mu\left(  \sigma^{-1}\left(  E\right)  \right)  =0$ by
(\ref{eqGenJun.okseskaft}).\medskip

The condition (\ref{eqGenJun.9}) does not, however, imply $\sigma_{i}%
$-quasi-invariance (\ref{eqGenJun.Z}), by the following example: $d=2$,
$\mu=\delta$-measure on $\left(  1,1,1,\dots\right)  $. Then (\ref{eqGenJun.9}%
) holds for all $E$, but (\ref{eqGenJun.Z}) fails for $E=\left\{  \left(
2,1,1,1,\dots\right)  \right\}  $ and $i=2$. In this case $\mu$ is $\sigma
$-quasi-invariant and $\sigma_{1}$-quasi-invariant, but not $\sigma_{2}%
$-quasi-invariant, so (\ref{eqGenJun.Y}) does not imply (\ref{eqGenJun.Z}).
More interestingly, the converse implication is always valid:

\begin{proposition}
\label{ProGenJun.1} If $\mu$ is a probability measure on $\Omega$ and $\mu$ is
$\sigma_{i}$-quasi-invariant for $i=1,\dots,d$, then $\mu$ is $\sigma$-quasi-invariant.
\end{proposition}

\begin{proof}
Put $\rho_{i}\left(  y\right)  =\frac{d\mu\left(  \sigma_{i}\left(  y\right)
\right)  }{d\mu\left(  y\right)  }$. Since the maps $\sigma_{i}$ are injective
and have disjoint ranges, there is actually one function $G$ such that%
\begin{equation}
G\left(  \sigma_{i}\left(  y\right)  \right)  =\rho_{i}\left(  y\right)  .
\label{eqGenJun.10}%
\end{equation}
One now proves as in (\ref{eqGen.25}) below (the tacit assumption there that
$\mu$ is $\sigma$-quasi-invariant is not needed for this) that%
\begin{equation}
\int_{\Omega}g\left(  x\right)  \,d\mu\left(  x\right)  =\int_{\Omega}R\left(
g\right)  \left(  x\right)  \,d\mu\left(  x\right)  , \label{eqGenJun.11}%
\end{equation}
where%
\begin{equation}
R\left(  g\right)  \left(  x\right)  =\sum_{\substack{y\\\sigma\left(
y\right)  =x}}G\left(  y\right)  g\left(  y\right)  . \label{eqGenJun.12}%
\end{equation}
Note that the Ruelle operator $R$ has the property%
\begin{equation}
R\left(  f\circ\sigma\right)  =R\left(  1\right)  \cdot f, \label{eqGenJun.13}%
\end{equation}
by the computation%
\[
R\left(  f\circ\sigma\right)  \left(  x\right)  =\sum_{\substack{y\\\sigma
\left(  y\right)  =x}}G\left(  y\right)  f\left(  \sigma\left(  y\right)
\right)  =f\left(  x\right)  \sum_{\substack{y\\\sigma\left(  y\right)
=x}}G\left(  y\right)  =f\left(  x\right)  R\left(  1\right)  \left(
x\right)  .
\]
We observe that $R\left(  1\right)  $ is a positive function by
(\ref{eqGenJun.12}) and it is $\mu$-integrable, as can be seen by using
(\ref{eqGenJun.11}) on $g=1$. If $f$ is a positive bounded function on
$\Omega$ we have from (\ref{eqGenJun.13}):%
\begin{equation}
\mu\left(  f\circ\sigma\right)  =\mu\left(  R\left(  f\circ\sigma\right)
\right)  =\mu\left(  R\left(  1\right)  \cdot f\right)  . \label{eqGenJun.14}%
\end{equation}
Putting $f$ equal to characteristic functions, the $\sigma$-quasi-invariance
(\ref{eqGenJun.Y}) of $\mu$ is immediate.
\end{proof}

As a final note on invariance, observe that (\ref{eqGenJun.9}) implies the
condition%
\[
\mu\left(  E\right)  =0\iff\mu\left(  \sigma\left(  E\right)  \right)  =0,
\]
by the following reasoning. Assume (\ref{eqGenJun.9}) throughout. If
$\mu\left(  \sigma\left(  E\right)  \right)  =0$, then $\mu\left(  \sigma
_{i}\sigma\left(  E\right)  \right)  =0$ for all $i$, but as $E\subseteq
\bigcup_{i}\sigma_{i}\sigma\left(  E\right)  $ it follows that $\mu\left(
E\right)  =0$. Conversely, if $\mu\left(  E\right)  =0$, write $E$ as
$E=\bigcup_{i}\sigma_{i}\left(  E_{i}\right)  $, and then $\mu\left(
\sigma_{i}\left(  E_{i}\right)  \right)  \leq\mu\left(  E\right)  =0$, and
hence, as $\sigma\left(  E\right)  =\bigcup_{i}\sigma_{i}\left(  E_{i}\right)
$, we have $\mu\left(  \sigma\left(  E\right)  \right)  =0$.\medskip

We now come to the main result in this chapter.

\begin{theorem}
\label{ThmGen.1}For any nondegenerate representation $s_{i}\mapsto S_{i}$ of
$\mathcal{O}_{d}$ on a separable Hilbert space $\mathcal{H}$, there exists a
probability measure $\mu$ on $\Omega$ which is $\sigma_{i}$-quasi-invariant
for $i=1,\dots,d$ \textup{(}and thus $\sigma$-quasi-invariant by Proposition
\textup{\ref{ProGenJun.1}),} and a measurable direct integral
decomposition\label{LOSH_1}%
\begin{equation}
\mathcal{H}=\int_{\Omega}^{\oplus}\mathcal{H}\left(  x\right)  \,d\mu\left(
x\right)  \label{eqGen.2}%
\end{equation}
of $\mathcal{H}$ into Hilbert spaces $\mathcal{H}\left(  x\right)  $ such that
the spaces $\mathcal{H}\left(  x\right)  $ are constant \textup{(}have fixed
dimension\/\textup{)} over $\sigma$-orbits in $\Omega$, and there exists a
measurable field $\Omega\ni x\mapsto U\left(  x\right)  $ of unitary operators
such that if%
\begin{equation}
\xi=\int_{\Omega}^{\oplus}\xi\left(  x\right)  \,d\mu\left(  x\right)
\label{eqGen.3}%
\end{equation}
is a vector in $\mathcal{H}$, then\label{LOSSi_1}\label{LOSSistar}%
\begin{align}
S_{i}\xi &  =\int_{\Omega}^{\oplus}\left(  S_{i}\xi\right)  \left(  x\right)
\,d\mu\left(  x\right)  ,\label{eqGen.4}\\
S_{i}^{\ast}\xi &  =\int_{\Omega}^{\oplus}\left(  S_{i}^{\ast}\xi\right)
\left(  x\right)  \,d\mu\left(  x\right)  , \label{eqGen.5}%
\end{align}
where%
\begin{align}
\left(  S_{i}^{{}}\xi\right)  \left(  x\right)   &  =\chi_{i}^{{}}\left(
x\right)  \rho\left(  x\right)  U\left(  x\right)  \xi\left(  \sigma\left(
x\right)  \right)  ,\label{eqGen.6}\\
\left(  S_{i}^{\ast}\xi\right)  \left(  x\right)   &  =\rho\left(  \sigma
_{i}\left(  x\right)  \right)  ^{-1}U\left(  \sigma_{i}\left(  x\right)
\right)  ^{\ast}\xi\left(  \sigma_{i}\left(  x\right)  \right)  .
\label{eqGen.7}%
\end{align}
Here%
\begin{equation}
\rho\left(  x\right)  =\left(  \frac{d\mu\left(  \sigma\left(  x\right)
\right)  }{d\mu\left(  x\right)  }\right)  ^{\frac{1}{2}}, \label{eqGen.9}%
\end{equation}
so that%
\begin{equation}
\rho\left(  \sigma_{i}\left(  x\right)  \right)  ^{-1}=\left(  \frac
{d\mu\left(  \sigma_{i}\left(  x\right)  \right)  }{d\mu\left(  x\right)
}\right)  ^{1/2} \label{eqGenJun.X}%
\end{equation}
and%
\begin{equation}
\chi_{i}^{{}}\left(  x\right)  =%
\begin{cases}1&\text{if }x_{1}=i,\\0&\text
{otherwise.}\end{cases}%
\label{eqGen.10}%
\end{equation}
Conversely, if $\mu$, $\mathcal{H}$, $x\mapsto\mathcal{H}\left(  x\right)  $
and $x\mapsto U\left(  x\right)  $ satisfy all the conditions in the initial
part of the theorem, the formulae \textup{(\ref{eqGen.4})--(\ref{eqGenJun.X})}
define a nondegenerate representation of $\mathcal{O}_{d}$\label{LOSOd_3} on
$\mathcal{H}$.
\end{theorem}

\begin{remark}
\label{RemGenJun.2}At the outset the formula \textup{(\ref{eqGen.6})} does not
make sense, since $U\left(  x\right)  $ is an operator on $\mathcal{H}\left(
x\right)  $, while $\xi\left(  \sigma\left(  x\right)  \right)  \in
\mathcal{H}\left(  \sigma\left(  x\right)  \right)  $. Here we have actually
made an identification of the Hilbert spaces $\mathcal{H}\left(  x\right)  $
over each orbit of $\sigma$. The Hilbert spaces over each $\sigma$-invariant
set have constant dimension $\mu$-almost everywhere by the argument after
\textup{(\ref{eqGenJun.angstrom})} below. Hence if we define%
\[
\Omega_{\left(  n\right)  }=\left\{  x\in\Omega\mid\dim\left(  \mathcal{H}%
\left(  x\right)  \right)  =n\right\}
\]
for $n=1,2,\dots,\aleph_{0}$, then the sets $\Omega_{\left(  n\right)  }$ are
$\mu$-measurable and $\sigma$-, as well as $\sigma_{i}$-, invariant up to sets
of measure zero. If $\mathcal{H}_{n}$ is the Hilbert space of dimension $n$
for $n=1,\dots,\aleph_{0}$, then we may identify $\mathcal{H}\left(  x\right)
$ with $\mathcal{H}_{n}$ for $x\in\Omega_{\left(  n\right)  }$, and we have
the decomposition%
\[
\mathcal{H=}\int_{\Omega}^{\oplus}\Omega\left(  x\right)  \,d\mu\left(
x\right)  =\bigoplus_{n=1}^{\aleph_{0}}\int_{\Omega_{\left(  n\right)  }%
}^{\oplus}\mathcal{H}\left(  x\right)  \,d\mu\left(  x\right)  =\bigoplus
_{n=1}^{\aleph_{0}}\mathcal{H}_{n}\otimes L^{2}\left(  \Omega_{\left(
n\right)  },d\mu\right)  .
\]
Hence we may view $U\left(  x\right)  $ as a unitary operator on
$\mathcal{H}_{n}$ for all $x\in\Omega_{\left(  n\right)  }$. Since the
$\Omega_{\left(  n\right)  }$ are $\sigma$- and $\sigma_{i}$-invariant, the
formula \textup{(\ref{eqGen.6})} is meaningful, and expressions like the one
on the third line of \textup{(\ref{eqGen.29})} make sense since $\xi\left(
x\right)  $ and $\eta\left(  \sigma\left(  x\right)  \right)  $ lie in the
same Hilbert space. The direct sum above is a decomposition of the
representation of $\mathcal{O}_{d}$. See also Remark \textup{\ref{RemGen.3},}
and see the book \cite[Section 2.5.3]{OsSa99} for more details.
\end{remark}

Before proving this theorem, let us consider the intertwiner space between two
representations $s_{i}\mapsto S_{i}$ and $s_{i}\mapsto T_{i}$. Recall that an
operator $T$ intertwines these representations if and only if it intertwines
the operators $S_{i}$, $T_{i}$:%
\begin{equation}
TS_{i}=T_{i}T,\qquad i=1,\dots,d. \label{eqGen.11}%
\end{equation}
``Only if'' is obvious. As for ``if'', note that if $T$ satisfies
(\ref{eqGen.11}), then
\begin{align}
T_{i}^{\ast}T  &  =\sum_{j=1}^{d}T_{i}^{\ast}TS_{j}^{{}}S_{j}^{\ast
}\label{eqGen.12}\\
&  =\sum_{j=1}^{d}T_{i}^{\ast}T_{j}^{{}}TS_{j}^{\ast}\nonumber\\
&  =TS_{i}^{\ast}.\nonumber
\end{align}

\begin{theorem}
\label{ThmGen.2}Let $S_{i}$, $\tilde{S}_{i}$ be representations of
$\mathcal{O}_{d}$ on separable Hilbert spaces\label{LOSH_2}%
\begin{align}
\mathcal{H}  &  =\int_{\Omega}^{\oplus}\mathcal{H}\left(  x\right)
\,d\mu\left(  x\right) \label{eqGen.13}\\%
\intertext{and}%
\widetilde{\mathcal{H}}  &  =\int_{\Omega}^{\oplus}\widetilde{\mathcal{H}%
}\left(  x\right)  \,d\tilde{\mu}\left(  x\right)  \label{eqGen.14}%
\end{align}
as defined in Theorem \textup{\ref{ThmGen.1}.} Partition $\Omega$ into three
$\sigma$-invariant Borel sets%
\begin{equation}
\Omega=\Omega_{0}\cup\Omega_{1}\cup\Omega_{2} \label{eqGen.15}%
\end{equation}
such that $\mu$ and $\tilde{\mu}$ are equivalent on $\Omega_{0}$, $\tilde{\mu
}\left(  \Omega_{1}\right)  =0$, and $\mu\left(  \Omega_{2}\right)  =0$. Then
an operator $T\in\mathcal{B}\left(  \mathcal{H},\smash{\widetilde{\mathcal{H}%
}}\right)  \vphantom{\widetilde{\mathcal{H}}}$ is an intertwiner between the
two representations, i.e.,%
\begin{equation}
TS_{i}=\tilde{S}_{i}T, \label{eqGen.16}%
\end{equation}
if and only if $T$ has a measurable decomposition%
\begin{equation}
T=\int_{\Omega_{0}}^{\oplus}T\left(  x\right)  \,d\mu\left(  x\right)
\label{eqGen.17}%
\end{equation}
where $T\left(  x\right)  \in\mathcal{B}\left(  \mathcal{H}\left(  x\right)
,\smash{\widetilde{\mathcal{H}}}\left(  x\right)  \right)  \vphantom
{\widetilde{\mathcal{H}}}$ and%
\begin{equation}
T\left(  x\right)  U\left(  x\right)  =\tilde{U}\left(  x\right)  T\left(
\sigma\left(  x\right)  \right)  \label{eqGen.18}%
\end{equation}
for almost all $x\in\Omega_{0}$.
\end{theorem}

\begin{remark}
\label{RemGen.3}In particular, if $S_{i}=T_{i}$, the commutant of the
representation consists of all decomposable operators%
\begin{equation}
T=\int_{\Omega}^{\oplus}T\left(  x\right)  \,d\mu\left(  x\right)
\label{eqGen.19}%
\end{equation}
such that%
\begin{equation}
T\left(  x\right)  U\left(  x\right)  =U\left(  x\right)  T\left(
\sigma\left(  x\right)  \right)  \label{eqGen.20}%
\end{equation}
for almost all $x\in\Omega$, and the center of the representation consists of
all decomposable operators%
\begin{equation}
T=\int_{\Omega}^{\oplus}\lambda\left(  x\right)  \openone_{\mathcal{H}\left(
x\right)  }\,d\mu\left(  x\right)  \label{eqGen.21}%
\end{equation}
where the scalar function $\lambda\in L^{\infty}\left(  \Omega,d\mu\right)  $
is $\sigma$-invariant. Thus the representation is a factor representation if
and only if the right shift on $L^{\infty}\left(  \Omega,d\mu\right)  $ is
ergodic. If in addition $\dim\left(  \mathcal{H}\left(  x\right)  \right)  =1$
for almost all $x$, then the representation is irreducible since
\textup{(\ref{eqGen.20})} then only has the trivial solutions $T\left(
x\right)  =$ const.

Note that if the right shift on $L^{2}\left(  \Omega,d\mu\right)  $ is
ergodic, then $\dim\left(  \mathcal{H}\left(  x\right)  \right)  $ is constant
for almost all $x$, and if $\mathcal{H}_{0}$ is a Hilbert space of that
dimension, then we have an isomorphism\label{LOSH_3}%
\[
\int_{\Omega}^{\oplus}\mathcal{H}\left(  x\right)  \,d\mu\left(  x\right)
\cong L^{2}\left(  \Omega,d\mu\right)  \otimes\mathcal{H}_{0}%
\]
and $U$ may be viewed as a measurable function from $\Omega$ into the unitary
group $\mathcal{U}\left(  \mathcal{H}_{0}\right)  $ on $\mathcal{H}_{0}$. The
element $T$ is then a function from $\Omega$ into $\mathcal{B}\left(
\mathcal{H}_{0}\right)  $ and \textup{(\ref{eqGen.20})} takes the form%
\[
T\left(  x\right)  =U\left(  x\right)  T\left(  \sigma\left(  x\right)
\right)  U^{\ast}\left(  x\right)  .
\]
Thus the commutant of the representation is canonically isomorphic to the
fixed-point algebra in%
\[
L^{\infty}\left(  \Omega,d\mu\right)  \otimes\mathcal{B}\left(  \mathcal{H}%
_{0}\right)  =L^{\infty}\left(  \Omega,\mathcal{B}\left(  \mathcal{H}%
_{0}\right)  \right)
\]
for the endomorphism%
\[
T\longmapsto U\left(  T\circ\sigma\right)  U^{\ast}.
\]
\end{remark}

Cocycle equivalence of functions with values in groups $G$ of unitaries have
been studied recently in ergodic theory; see, e.g., \cite{Par97,PaPo97}.
Equation (\ref{eqGen.18}) above in that setup is the assertion that $U$ and
$\tilde{U}$ (taking values in the corresponding $G$) are cohomologous.

\begin{proof}
[Proof of Theorem \textup{\ref{ThmGen.1}}]We will first verify that the
relations (\ref{eqGen.2})--(\ref{eqGen.10}) define a representation of
$\mathcal{O}_{d}$, and verify that its restriction to the abelian
subalgebra\label{LOSDd_2}
\begin{equation}
\mathcal{D}_{d}=C^{\ast}\left(  s_{\alpha}^{{}}s_{\alpha}^{\ast}\biggm|%
\alpha\in\vphantom{\coprod}\smash{\coprod_{1}^{\infty}}\mathbb{Z}_{d}\right)
\vphantom{\coprod_{1}^{\infty}} \label{eqGen.22}%
\end{equation}
is the spectral representation. If $g\in L^{1}\left(  \Omega,d\mu\right)  $,
we have%
\begin{align}
\int_{\Omega}g\left(  x\right)  \,d\mu\left(  x\right)   &  =\sum_{i}%
\int_{\Omega_{i}}g\left(  x\right)  \,d\mu\left(  x\right) \label{eqGen.25}\\
&  =\sum_{i}\int_{\Omega}g\left(  \sigma_{i}\left(  y\right)  \right)
\frac{d\mu\left(  \sigma_{i}\left(  y\right)  \right)  }{d\mu\left(  y\right)
}\,d\mu\left(  y\right) \nonumber\\
&  =\sum_{i}\int_{\Omega}g\left(  \sigma_{i}\left(  y\right)  \right)
\frac{d\mu\left(  \sigma_{i}\left(  y\right)  \right)  }{d\mu\left(
\sigma\sigma_{i}\left(  y\right)  \right)  }\,d\mu\left(  y\right) \nonumber\\
&  =\sum_{i}\int_{\Omega}g\left(  \sigma_{i}\left(  y\right)  \right)
\rho\left(  \sigma_{i}\left(  y\right)  \right)  ^{-2}\,d\mu\left(  y\right)
\nonumber\\
&  =\int_{\Omega}\vphantom{\sum_{\substack{x\\\sigma\left( x\right) =y}}%
}\left(  \;\vphantom{\sum_{}}\smash{\sum_{\substack{x\\\sigma\left
( x\right) =y}}}g\left(  x\right)  G\left(  x\right)  \right)  \,d\mu\left(
y\right)  ,\nonumber
\end{align}
where $G\left(  x\right)  =\rho\left(  x\right)  ^{-2}$. (If it happens that
$\sum_{x\colon\sigma\left(  x\right)  =y}G\left(  x\right)  =1$, the relation
(\ref{eqGen.25}) says that $\mu$ is $\sigma$-invariant, and $\mu$ is then what
is called a $G$-measure in \cite{Kea72}.) Applying (\ref{eqGen.25}) to
$g\left(  x\right)  =f\left(  x,\sigma\left(  x\right)  \right)  $ we obtain%
\begin{equation}
\int_{\Omega}f\left(  x,\sigma\left(  x\right)  \right)  \,d\mu\left(
x\right)  =\sum_{i}\int_{\Omega}f\left(  \sigma_{i}\left(  y\right)
,y\right)  \rho\left(  \sigma_{i}\left(  y\right)  \right)  ^{-2}\,d\mu\left(
y\right)  . \label{eqGen.26}%
\end{equation}
Defining $S_{i}$ by (\ref{eqGen.4}) and (\ref{eqGen.6}), we see immediately
from the $\chi_{i}^{{}}\left(  x\right)  $ term that the ranges of $S_{i}$ are
mutually orthogonal, and if $\xi\in\mathcal{H}$, then from (\ref{eqGen.26}):%
\begin{align}
\left\|  S_{i}\xi\right\|  _{2}^{2}  &  =\int_{\Omega}\chi_{i}^{{}}\left(
x\right)  \rho\left(  x\right)  ^{2}\left\|  \xi\left(  \sigma\left(
x\right)  \right)  \right\|  ^{2}\,d\mu\left(  x\right) \label{eqGen.27}\\
&  =\int_{\Omega_{i}}\rho\left(  x\right)  ^{2}\left\|  \xi\left(
\sigma\left(  x\right)  \right)  \right\|  ^{2}\,d\mu\left(  x\right)
\nonumber\\
&  =\int_{\Omega}\rho\left(  \sigma_{i}\left(  y\right)  \right)  ^{2}\left\|
\xi\left(  y\right)  \right\|  ^{2}\rho\left(  \sigma_{i}\left(  y\right)
\right)  ^{-2}\,d\mu\left(  y\right) \nonumber\\
&  =\left\|  \xi\right\|  ^{2}\nonumber
\end{align}
so each $S_{i}$ is an isometry, and hence%
\begin{equation}
S_{i}^{\ast}S_{j}^{{}}=\delta_{ij}\openone. \label{eqGen.28}%
\end{equation}
Furthermore%
\begin{align}%
\ip{S_{i}^{\ast}\xi}{\eta}%
&  =%
\ip{\xi}{S_{i}^{}\eta}%
\label{eqGen.29}\\
&  =\int_{\Omega}%
\ip{\xi\left( x\right) }{\left( S_{i}\eta\right
) \left( x\right) }%
\,d\mu\left(  x\right) \nonumber\\
&  =\int_{\Omega}\chi_{i}^{{}}\left(  x\right)  \rho\left(  x\right)
\ip{U\left( x\right) ^{\ast}\xi\left( x\right) }{\eta\left( \sigma
\left( x\right) \right) }%
\,d\mu\left(  x\right) \nonumber\\
&  =\int_{\Omega_{i}}\rho\left(  x\right)
\ip{U\left( x\right) ^{\ast}\xi\left( x\right) }{\eta\left( \sigma
\left( x\right) \right) }%
\,d\mu\left(  x\right) \nonumber\\
&  =\int_{\Omega}\rho\left(  \sigma_{i}\left(  y\right)  \right)
\ip{U\left( \sigma_{i}\left( y\right) \right) ^{\ast
}\xi\left( \sigma_{i}\left( y\right) \right) }{\eta\left( y\right) }%
\rho\left(  \sigma_{i}\left(  y\right)  \right)  ^{-2}\,d\mu\left(  y\right)
\nonumber\\
&  =\int_{\Omega}\rho\left(  \sigma_{i}\left(  y\right)  \right)  ^{-1}%
\ip{U\left( \sigma_{i}\left( y\right) \right) ^{\ast
}\xi\left( \sigma_{i}\left( y\right) \right) }{\eta\left( y\right) }%
\,d\mu\left(  y\right)  ,\nonumber
\end{align}
and the expression (\ref{eqGen.7}) for $S_{i}^{\ast}$ follows.

If $\alpha=\left(  \alpha_{1}\alpha_{2}\dots\alpha_{n}\right)  $%
\label{LOSalpha} with $\alpha_{k}\in\mathbb{Z}_{d}$, define\label{LOSsalpha_1}%
\label{LOSSalpha_1}%
\begin{equation}
s_{\alpha_{{}}}=s_{\alpha_{1}}s_{\alpha_{2}}\cdots s_{\alpha_{n}},\qquad
S_{\alpha_{{}}}=S_{\alpha_{1}}S_{\alpha_{2}}\cdots S_{\alpha_{n}}.
\label{eqGen.30}%
\end{equation}
One verifies from (\ref{eqGen.6}) and (\ref{eqGen.7}) that%
\begin{multline}
S_{\alpha}\xi\left(  x\right)  =\chi_{\alpha_{1}}^{{}}\left(  x\right)
\chi_{\alpha_{2}}^{{}}\left(  \sigma\left(  x\right)  \right)  \cdots
\chi_{\alpha_{n}}^{{}}\left(  \sigma^{n-1}\left(  x\right)  \right)
\label{eqGen.31}\\
\times\rho\left(  x\right)  \rho\left(  \sigma\left(  x\right)  \right)
\cdots\rho\left(  \sigma^{n-1}\left(  x\right)  \right) \\
\times U\left(  x\right)  U\left(  \sigma\left(  x\right)  \right)  \cdots
U\left(  \sigma^{n-1}\left(  x\right)  \right)  \xi\left(  \sigma^{n}\left(
x\right)  \right)
\end{multline}
and%
\begin{multline}
S_{\alpha}^{\ast}\xi\left(  x\right)  =\rho\left(  \sigma_{\alpha_{n}}\left(
x\right)  \right)  ^{-1}\rho\left(  \sigma_{\alpha_{n-1}}\sigma_{\alpha_{n}%
}\left(  x\right)  \right)  ^{-1}\cdots\rho\left(  \sigma_{\alpha_{1}}%
\cdots\sigma_{\alpha_{n}}\left(  x\right)  \right)  ^{-1}\label{eqGen.32}\\
\times U\left(  \sigma_{\alpha_{n}}\left(  x\right)  \right)  ^{\ast}U\left(
\sigma_{\alpha_{n-1}}\sigma_{\alpha_{n}}\left(  x\right)  \right)  ^{\ast
}\cdots U\left(  \sigma_{\alpha_{1}}\cdots\sigma_{\alpha_{n}}\left(  x\right)
\right)  ^{\ast}\xi\left(  \sigma_{\alpha_{1}}\cdots\sigma_{\alpha_{1}}\left(
x\right)  \right)  .
\end{multline}
Combining (\ref{eqGen.31})--(\ref{eqGen.32}) with the relations%
\begin{equation}
\begin{gathered} \sigma_{\alpha_{n}}\sigma^{n}\left( x\right) =\sigma
^{n-1}\left( x\right) ,\\ {}\mathrel{\mkern4.5mu\vdots\mkern4.5mu}{}%
\\ \sigma_{\alpha_{1}}\cdots\sigma_{\alpha_{n}}\sigma^{n}\left( x\right
) =x, \end{gathered} \label{eqGen.33}%
\end{equation}
which are valid if $x=\left(  \alpha_{1},\dots,\alpha_{n},x_{n+1}%
,\dots\right)  $, we obtain%
\begin{equation}
S_{\alpha}^{{}}S_{\alpha}^{\ast}\xi\left(  x\right)  =\chi_{\alpha}^{{}%
}\left(  x\right)  \xi\left(  x\right)  , \label{eqGen.34}%
\end{equation}
where%
\begin{equation}
\chi_{\alpha}^{{}}\left(  x_{1},x_{2},x_{3},\dots\right)  =\delta_{\alpha
_{1}x_{1}}\delta_{\alpha_{2}x_{2}}\cdots\delta_{\alpha_{n}x_{n}}.
\label{eqGen.35}%
\end{equation}
This proves firstly that
\begin{equation}
\sum_{i=1}^{d}S_{i}^{{}}S_{i}^{\ast}=\openone; \label{eqGen.36}%
\end{equation}
and (\ref{eqGen.28}) and (\ref{eqGen.36}) show that $s_{i}\mapsto S_{i}$ is
indeed a representation of the Cuntz relations. Secondly, (\ref{eqGen.34})
shows that $\mathcal{D}_{d}$ maps onto the algebra of operators on
$\mathcal{H}$ of the form%
\begin{equation}
\int_{\Omega}\lambda\left(  x\right)  \openone_{\mathcal{H}\left(  x\right)
}\,d\mu\left(  x\right)  \label{eqGen.37}%
\end{equation}
where $\lambda$ ranges over all continuous functions on the Cantor set
$\Omega$. Thus the restriction of the representation $s_{i}\mapsto S_{i}$ to
$\mathcal{D}_{d}$ is indeed the spectral representation.

To show the main part of Theorem \ref{ThmGen.1}, i.e., the existence of the
objects $\mathcal{H}\left(  x\right)  $, $d\mu\left(  x\right)  $, $U\left(
x\right)  $, one does indeed start with a spectral measure $\mu$ for the
restriction of the representation to $\mathcal{D}_{d}$. The spectrum of
$\mathcal{D}_{d}$ is $\Omega$, so this gives the decomposition (\ref{eqGen.2}%
), and the action of $\mathcal{D}_{d}$ on $\mathcal{H}$ is given by
(\ref{eqGen.34}). If $f\in C\left(  \Omega\right)  =\mathcal{D}_{d}$, and
$M_{f}$ is the representative of $f$ in $\mathcal{H}$:\label{LOSMf}%
\begin{equation}
M_{f}=\int_{\Omega}^{\oplus}f\left(  x\right)  \openone_{\mathcal{H}\left(
x\right)  }\,d\mu\left(  x\right)  , \label{eqGen.38}%
\end{equation}
then%
\begin{equation}
M_{f\circ\sigma}=\sum_{i=1}^{d}S_{i}^{{}}M_{f}^{{}}S_{i}^{\ast}
\label{eqGen.39}%
\end{equation}
and the quasi-invariance of $\mu$ under $\sigma$ follows. Thus one may define
$\rho\left(  x\right)  $ by (\ref{eqGen.9}). Similarly, if $f\in C\left(
\Omega\right)  =\mathcal{D}_{d}$ has support in $\sigma_{i}\left(
\Omega\right)  =i\Omega=\Omega_{i}$, one verifies that%
\begin{equation}
M_{f\circ\sigma_{i}}=S_{i}^{\ast}M_{f}^{{}}S_{i}^{{}}.
\label{eqGenJun.angstrom}%
\end{equation}
Thus the two representations of $C\left(  \Omega_{i}\right)  $ given by
$f\mapsto M_{f}$ on $M_{\chi_{i}^{{}}}\mathcal{H}$ and $f\mapsto
M_{f\circ\sigma_{i}}$ on $\mathcal{H}$ are unitarily equivalent. In
particular, this means that $\dim\left(  \mathcal{H}\left(  x\right)  \right)
=\dim\left(  \mathcal{H}\left(  \sigma_{i}\left(  x\right)  \right)  \right)
$ for $\mu$-almost all $x$, so the constancy of $\dim\left(  \mathcal{H}%
\left(  x\right)  \right)  $ almost everywhere over the orbits of $\sigma
_{1},\dots,\sigma_{d}$ follows. But (\ref{eqGen.24}) then implies that
$\dim\left(  \mathcal{H}\left(  x\right)  \right)  $ is constant on $\sigma
$-orbits (actually the two forms of constancy are equivalent). Also it follows
from the unitary equivalence (\ref{eqGenJun.angstrom}) that $\mu$ is
quasi-invariant under $\sigma_{i}$ and that $\rho\left(  \sigma_{i}\left(
x\right)  \right)  ^{-1}=\left(  \frac{d\mu\left(  \sigma_{i}\left(  x\right)
\right)  }{d\mu\left(  x\right)  }\right)  ^{1/2}$ exists. See \cite{Nel69} or
\cite[Section 2.5.3]{OsSa99}
for details on spectral multiplicity theory. Now, one may define
a representation $s_{i}\mapsto T_{i}$ of $\mathcal{O}_{d}$ on $\mathcal{H}$ by%
\begin{equation}
\left(  T_{i}\xi\right)  \left(  x\right)  =\chi_{i}^{{}}\left(  x\right)
\rho\left(  x\right)  \xi\left(  \sigma\left(  x\right)  \right)  .
\label{eqGen.40}%
\end{equation}
One checks that this is indeed a representation of $\mathcal{O}_{d}$ by the
first part of the proof, and by the proof of (\ref{eqGen.34}) it follows that%
\begin{equation}
T_{\alpha}^{{}}T_{\alpha}^{\ast}=S_{\alpha}^{{}}S_{\alpha}^{\ast}
\label{eqGen.41}%
\end{equation}
for all multi-indices $\alpha$. Define an operator $U$ by\label{LOSU_1}%
\begin{equation}
U=\sum_{i=1}^{d}S_{i}^{{}}T_{i}^{\ast}. \label{eqGen.42}%
\end{equation}
Using the Cuntz relations in a standard manner, one checks that $U$ is a
unitary operator, and
\begin{equation}
S_{i}=UT_{i} \label{eqGen.43}%
\end{equation}
for $i=1,\dots,d$. Putting
\begin{equation}
i\alpha:=\left(  i\,\alpha_{1}\alpha_{2}\dots\alpha_{n}\right)  ,
\label{eqGen.44}%
\end{equation}
we have by (\ref{eqGen.41})%
\begin{align}
T_{i\alpha}^{{}}T_{i\alpha}^{\ast}  &  =S_{i\alpha}^{{}}S_{i\alpha}^{\ast
}\label{eqGen.45}\\
&  =S_{i}^{{}}S_{\alpha}^{{}}S_{\alpha}^{\ast}S_{i}^{\ast}\nonumber\\
&  =UT_{i}^{{}}T_{\alpha}^{{}}T_{\alpha}^{\ast}T_{i}^{\ast}U_{{}}^{\ast
}\nonumber\\
&  =UT_{i\alpha}^{{}}T_{i\alpha}^{\ast}U_{{}}^{\ast},\nonumber
\end{align}
and hence $U$ commutes with the representatives on $\mathcal{H}$ of the
algebra $\mathcal{D}_{d}$. Hence $U$ has a decomposition\label{LOSU_2}%
\begin{equation}
U=\int_{\Omega}^{\oplus}U\left(  x\right)  \,d\mu\left(  x\right)
\label{eqGen.46}%
\end{equation}
where $\Omega\ni x\mapsto U\left(  x\right)  $ is a measurable field on
unitaries. It now follows from (\ref{eqGen.40}) and (\ref{eqGen.43}) that
$S_{i}$ has the form (\ref{eqGen.6}). This ends the proof of Theorem
\ref{ThmGen.1}.
\end{proof}

\begin{proof}
[Proof of Theorem \textup{\ref{ThmGen.2}}]Adopt the assumptions in Theorem
\ref{ThmGen.2} and let $T$ be an intertwiner between the two representations.
In particular this means that $T$ intertwines the two spectral representations
of $\mathcal{D}_{d}$ on $\mathcal{H}$ and $\widetilde{\mathcal{H}}$,
respectively, i.e.,%
\begin{equation}
TS_{\alpha}^{{}}S_{\alpha}^{\ast}=\tilde{S}_{\alpha}^{{}}\tilde{S}_{\alpha
}^{\ast}T \label{eqGen.47}%
\end{equation}
for all multi-indices $\alpha$. But this is equivalent to $\Omega$ having the
decomposition (\ref{eqGen.15}) and $T$ having the measurable decomposition%
\begin{equation}
T=\int_{\Omega_{0}}^{\oplus}T\left(  x\right)  \,d\mu\left(  x\right)
\label{eqGen.48}%
\end{equation}
where $T\left(  x\right)  \in\mathcal{B}\left(  \mathcal{H}\left(  x\right)
,\smash{\widetilde{\mathcal{H}}}\left(  x\right)  \right)  \vphantom
{\widetilde{\mathcal{H}}}$. We now compute, using (\ref{eqGen.6}),%
\begin{align}
TS_{i}\xi\left(  x\right)   &  =T\left(  x\right)  \left(  S_{i}\xi\right)
\left(  x\right) \label{eqGen.49}\\
&  =\chi_{i}^{{}}\left(  x\right)  \rho\left(  x\right)  T\left(  x\right)
U\left(  x\right)  \xi\left(  \sigma\left(  x\right)  \right) \nonumber
\end{align}
and%
\begin{align}
\tilde{S}_{i}T\xi\left(  x\right)   &  =\chi_{i}^{{}}\left(  x\right)
\rho\left(  x\right)  \tilde{U}\left(  x\right)  \left(  T\xi\right)  \left(
\sigma\left(  x\right)  \right) \label{eqGen.50}\\
&  =\chi_{i}^{{}}\left(  x\right)  \rho\left(  x\right)  \tilde{U}\left(
x\right)  T\left(  \sigma\left(  x\right)  \right)  \xi\left(  \sigma\left(
x\right)  \right)  .\nonumber
\end{align}
Using the intertwining property (\ref{eqGen.16}) we thus deduce that%
\begin{equation}
T\left(  x\right)  U\left(  x\right)  =\tilde{U}\left(  x\right)  T\left(
\sigma\left(  x\right)  \right)  . \label{eqGen.51}%
\end{equation}
Conversely, if $T$ satisfies (\ref{eqGen.51}), the intertwining follows from
(\ref{eqGen.49}) and (\ref{eqGen.50}). This ends the proof of Theorem
\ref{ThmGen.2}.
\end{proof}

\chapter{\label{Fre}The free group on $d$ generators}

In this chapter we will construct certain representations of $\mathcal{O}_{d}%
$\label{LOSOd_4} in the Hilbert spaces $\mathcal{H}$\label{LOSH_4} where the
decomposition in Theorem \ref{ThmGen.1} takes the form%
\begin{equation}
\int_{\Omega}^{\oplus}\mathcal{H}\left(  x\right)  \,d\mu\left(  x\right)
\cong L^{2}\left(  \Omega,d\mu\right)  \otimes\mathcal{H}_{0}. \label{eq2.5A}%
\end{equation}
We will equip $\Omega=\prod_{1}^{\infty}\mathbb{Z}_{d}$ with the product
measure $\mu=\mu_{p}$ defined from a choice of weights $\left(  p_{i}\right)
_{i=1}^{d}$, with $p_{i}>0$, and $\sum_{i}p_{i}=1$. Then the representation
(\ref{eqGen.6})--(\ref{eqGen.7}) takes the form%
\begin{align}
\left(  S_{i}\xi\right)  \left(  x\right)   &  =\delta_{i}\left(
x_{1}\right)  \frac{1}{\sqrt{p_{i}}}U\left(  x\right)  \xi\left(
\sigma\left(  x\right)  \right)  ,\label{eqFre.1}\\
\left(  S_{i}^{\ast}\xi\right)  \left(  x\right)   &  =\sqrt{p_{i}}U\left(
\sigma_{i}\left(  x\right)  \right)  ^{\ast}\xi\left(  \sigma_{i}\left(
x\right)  \right)  . \label{eqFre.2}%
\end{align}
The simplest case of this is when $\dim\mathcal{H}_{0}=1$ and $U\left(
x\right)  \equiv1$. Then the corresponding operators $S_{i}$ of (\ref{eqFre.1}%
) act on scalar functions in $L^{2}\left(  \Omega,\mu\right)  $. The constant
function $\openone$ in $L^{2}\left(  \Omega,\mu\right)  $ satisfies
$S_{i}^{\ast}\openone=\sqrt{p_{i}}\openone$, and the state $\omega_{\openone
}=%
\ip{\openone}{\,\cdot\,\openone}%
$ on $\mathcal{B}\left(  L^{2}\left(  \Omega,\mu\right)  \right)  $ satisfies%
\begin{equation}
\omega_{\openone}\left(  \rho\left(  A\right)  \right)  =\omega_{\openone
}\left(  A\right)  ,\qquad A\in\mathcal{B}\left(  L^{2}\left(  \Omega
,\mu\right)  \right)  \label{eqFre.3}%
\end{equation}
where%
\[
\rho\left(  A\right)  =\sum S_{i}^{{}}AS_{i}^{\ast}.
\]
This is the representation defined by the Cuntz states \cite[Theorem
4.1]{BJP96}.

It is well known, see, e.g., \cite{BJP96}, that there is a correspondence
between representations of $\mathcal{O}_{d}$ (for some $d$ including
$d=\infty$) and endomorphisms of $\mathcal{B}\left(  \mathcal{H}\right)  $. An
endomorphism $\rho$ of $\mathcal{B}\left(  \mathcal{H}\right)  $ has a finite
Powers index $d$ if the commutant of $\rho\left(  \mathcal{B}\left(
\mathcal{H}\right)  \right)  $ is isomorphic to $M_{d}\left(  \mathbb{C}%
\right)  $\label{LOSMd}, and then the corresponding representation is of
$\mathcal{O}_{d}$. Two representations $\pi$, $\tilde{\pi}$ of $\mathcal{O}%
_{d}$ define the same endomorphism if and only if there exists a $g$ in the
group $U\left(  d\right)  $ of complex unitary $d\times d$ matrices such that
$\tilde{\pi}=\pi\circ\alpha_{g}$ where $\alpha_{g}$ is the canonical $U\left(
d\right)  $-action on $\mathcal{O}_{d}$ rotating the generators.

There is precisely one conjugacy class of endomorphisms of $\mathcal{B}\left(
\mathcal{H}\right)  $ with an invariant vector state $\omega$, i.e.,%
\begin{equation}
\omega\circ\rho=\omega,\label{eqFre.4}%
\end{equation}
see (\ref{eqFre.3}) and \cite{Pow67,Pow88} or \cite[Theorem 4.2]{BJP96}. We
showed in \cite{BrJo97b} and \cite{BrJo97a} that the theory of wavelets gives
examples of endomorphisms in different conjugacy classes. In this memoir, we
will also look at endomorphisms of von Neumann algebras not of type \textup{I.}

Scale-two wavelet representations are constructed from measurable functions on
$\mathbb{T}$ subject to $\left|  m\left(  z\right)  \right|  ^{2}+\left|
m\left(  -z\right)  \right|  ^{2}=2$. If%
\[%
\begin{cases}
m_{1}\left(  z\right)  =m\left(  z\right)  \\
m_{2}\left(  z\right)  =z\bar{m}\left(  -z\right)
\end{cases}%
\]
then\label{LOSSi_2}%
\begin{equation}
\left(  S_{j}\xi\right)  \left(  z\right)  =m_{j}\left(  z\right)  \xi\left(
z^{2}\right)  ,\qquad j=1,2 \label{eqFre.5}%
\end{equation}
define a representation of $\mathcal{O}_{2}$. Giving the wavelet
representations (of $\mathcal{O}_{2}$) in the form of Theorem \ref{ThmGen.1}
amounts to representing the commuting operators (in fact
projections)\label{LOSSalpha_2}
\[
S_{j_{1}}^{{}}\cdots S_{j_{k}}^{{}}S_{j_{k}}^{\ast}\cdots S_{j_{1}}^{\ast
}=S_{\alpha}^{{}}S_{\alpha}^{\ast}\qquad\left(  \alpha=\left(  j_{1}\dots
j_{k}\right)  \right)
\]
as multiplication operators on some $L^{2}\left(  \Omega,\mathcal{H}%
_{0}\right)  $. Such a representation will involve a $2$-adic completion, but
will perhaps not be explicit enough for practical applications: In the
representation, the operator%
\begin{multline}
\left(  S_{\alpha}^{{}}S_{\alpha}^{\ast}\xi\right)  \left(  z\right)
=\frac{1}{2^{k}}m_{j_{1}}\left(  z\right)  \cdots m_{j_{k}}\left(  z^{2^{k-1}%
}\right) \label{eqFre.6}\\
\times\sum_{w^{2^{k}}=1}\bar{m}_{j_{1}}\left(  wz\right)  \bar{m}_{j_{2}%
}\left(  w^{2}z^{2}\right)  \cdots\bar{m}_{j_{k}}\left(  w^{2^{k-1}}%
z^{2^{k-1}}\right)  \xi\left(  wz\right)
\end{multline}
must be multiplication by a characteristic function of a set $E_{\alpha}$ in
the $2$-adic completion.

We postpone the details to a later paper.

Returning to the computation of the measurable field $\Omega\ni x\mapsto
U\left(  x\right)  $ of unitary operators in Theorem \ref{ThmGen.1}, we do the
calculation for the $\left(  p_{1},\dots,p_{d}\right)  $-product measure on
$\Omega=\prod_{1}^{\infty}\mathbb{Z}_{d}$,\label{LOSOmega_2} and with the
resulting representation $s_{i}\mapsto S_{i}$ of type \textup{III.} (More
details on $\beta$-KMS states\label{LOSsigmaLbetaKMSstate_2} and the
$\mathbb{T}^{d}\subset U\left(  d\right)  $ action on $\mathcal{O}_{d}$ are
included in Chapter \ref{KMS} below.) We show there that if\label{LOSpj_2}
\begin{equation}
p_{j}=e^{-\beta L_{j}},\qquad j=1,\dots,d, \label{eqFre.7}%
\end{equation}
and $L=\left(  L_{1},\dots,L_{d}\right)  \in\mathbb{R}^{d}$, $L_{j}>0$, then
the state $\omega$ on $\mathcal{O}_{d}$ given by
\begin{equation}
\omega\left(  s_{i_{1}}^{{}}\cdots s_{i_{k}}^{{}}s_{j_{l}}^{\ast}\cdots
s_{j_{1}}^{\ast}\right)  =\delta_{k_{{}}l_{{}}}\delta_{i_{1}j_{1}}\cdots
\delta_{i_{k}j_{k}}p_{i_{1}}p_{i_{2}}\cdots p_{i_{k}} \label{eqFre.8}%
\end{equation}
is a (unique) $\beta$-KMS state for the one-parameter subgroup of
$\mathbb{T}^{d}$ defined by $L$, i.e., $t\mapsto\left(  e^{itL_{1}}%
,e^{itL_{2}},\dots,e^{itL_{d}}\right)  $. (For $\omega$ to be a state, $\beta$
must be chosen such that $\sum_{j}p_{j}=1$, and then $\omega_{\beta}$ is the
gauge-invariant extension to $\mathcal{O}_{d}$ of the product state defined on
$\operatorname*{UHF}\nolimits_{d}\cong\bigotimes_{1}^{\infty}M_{d}%
$\label{LOSUHFd} as $\bigotimes_{1}^{\infty}\varphi$, where $\varphi$ is the
state on $M_{d}$ defined by the density matrix $\operatorname*{diag}\left(
p_{1},\dots,p_{d}\right)  $.) Let $s_{i}\mapsto T_{i}$ be the representation
of $\mathcal{O}_{d}$ which is induced from $\omega$ via the GNS construction.
Let $\mathbb{F}_{d}$\label{LOSFd} be the free group on $d$ generators
$g_{1},\dots,g_{d}$, and let $\mathbb{F}_{d}\ni g\mapsto\lambda\left(
g\right)  $ be the trace representation of $\mathbb{F}_{d}$. Recall the trace
$\operatorname*{tr}$ on $C_{r}^{\ast}\left(  \mathbb{F}_{d}\right)  $ is given
by
\[
\operatorname*{tr}\left(  \lambda\left(  g\right)  \right)  =%
\begin{cases}
1 & \text{if }g=e,\\
0 & \text{if }g\neq e.
\end{cases}%
\]
The Hilbert space $\ell^{2}\left(  \mathbb{F}_{d}\right)  $ has as orthonormal
basis the functions%
\[
\left\{  \xi_{g}\mid g\in\mathbb{F}_{d}\right\}  \text{\qquad where }\xi
_{g}\left(  x\right)  =%
\begin{cases}
1 & \text{if }x=g,\\
0 & \text{if }x\neq g,
\end{cases}%
\]
and
\begin{equation}
\operatorname*{tr}\left(  A\right)  =%
\ip{\xi_{e}}{\lambda\left( A\right) \xi_{e}}%
,\qquad A\in C_{r}^{\ast}\left(  \mathbb{F}_{d}\right)  . \label{eqFre.9}%
\end{equation}
Let $\mathcal{H}_{0}=\ell^{2}\left(  \mathbb{F}_{d}\right)  $. For
multi-indices $\alpha$, set $p^{\alpha}=p_{\alpha_{1}}\cdots p_{\alpha_{k}}$.

\begin{proposition}
\label{ProFre.1}The state defined by \textup{(\ref{eqFre.8}),} i.e.,%
\begin{equation}
\omega\left(  s_{\alpha}^{{}}s_{\gamma}^{\ast}\right)  =\delta_{\alpha\gamma
}p^{\alpha} \label{eqFre.10}%
\end{equation}
is the vector state defined by $\openone\otimes\xi_{e}$ in the representation
on $L_{\mu}^{2}\left(  \prod_{1}^{\infty}\mathbb{Z}_{d},\ell^{2}\left(
\mathbb{F}_{d}\right)  \right)  =L_{\mu}^{2}\left(  \prod_{1}^{\infty
}\mathbb{Z}_{d}\right)  \otimes\ell^{2}\left(  \mathbb{F}_{d}\right)  $ by%
\begin{align}
\left(  T_{j}\xi\right)  \left(  x\right)   &  =e^{\frac{1}{2}\beta L_{j}}%
\chi_{j}^{{}}\left(  x\right)  \lambda\left(  g_{j}\right)  \xi\left(
\sigma\left(  x\right)  \right)  ,\label{eqFre.11}\\
\left(  T_{j}^{\ast}\xi\right)  \left(  x\right)   &  =e^{-\frac{1}{2}\beta
L_{j}}\lambda\left(  g_{j}^{-1}\right)  \xi\left(  \sigma_{i}\left(  x\right)
\right)  , \label{eqFre.12}%
\end{align}
so, in particular, the operators $U\left(  x\right)  $ from Theorem
\textup{\ref{ThmGen.1}} are independent of the product measure $\mu$ when the
representation is realized in $L_{\mu}^{2}\left(  \prod_{1}^{\infty}%
\mathbb{Z}_{d},\ell^{2}\left(  \mathbb{F}_{d}\right)  \right)  $, i.e., on
vector-valued functions on the group $\prod_{1}^{\infty}\mathbb{Z}_{d}$ with
values in $\ell^{2}\left(  \mathbb{F}_{d}\right)  $ and with $\mu$ equal to
the product measure on $\prod_{1}^{\infty}\mathbb{Z}_{d}$ relative to
$p_{j}=e^{-\beta L_{j}}$, $j=1,\dots,d$.
\end{proposition}

\begin{proof}
Note that the representation $T_{j}$ in (\ref{eqFre.11})--(\ref{eqFre.12}) is
of the form $T_{j}=S_{j}\otimes\lambda\left(  g_{j}\right)  $ where $S_{j}$ is
the representation in Theorem \ref{ThmGen.1} corresponding to the
scalar-valued case with $\mu$ product measure and $U\equiv1$. We then use

\begin{lemma}
\label{LemFre.2}Let $\left(  S_{j}\right)  $ be a representation of
$\mathcal{O}_{d}$ in a Hilbert space $\mathcal{L}$ and let $\mathcal{H}_{0}$
be a second Hilbert space. If $\left(  A_{j}\right)  _{j=1}^{d}$ are operators
in $\mathcal{H}_{0}$, then $T_{j}:=S_{j}\otimes A_{j}$ define a representation
of $\mathcal{O}_{d}$ in $\mathcal{L}\otimes\mathcal{H}_{0}$ if and only if the
$A_{i}$'s are unitary.
\end{lemma}

\begin{proof}
We have%
\begin{align*}
T_{j}^{\ast}T_{k}^{{}}  &  =S_{j}^{\ast}S_{k}^{{}}\otimes A_{j}^{\ast}%
A_{k}^{{}}\\
&  =\delta_{jk}\openone_{\mathcal{L}}\otimes A_{j}^{\ast}A_{k}^{{}}.
\end{align*}
Hence $T_{j}^{\ast}T_{k}^{{}}=\delta_{jk}\openone_{\mathcal{L}\otimes
\mathcal{H}_{0}}$ holds if and only if each $A_{j}$ is isometric.

We have%
\[
\sum_{j}T_{j}^{{}}T_{j}^{\ast}=\sum_{j}S_{j}^{{}}S_{j}^{\ast}\otimes A_{j}%
^{{}}A_{j}^{\ast}.
\]
But the projections $S_{j}^{{}}S_{j}^{\ast}$ are mutually orthogonal. So
$\sum_{j}T_{j}^{{}}T_{j}^{\ast}=\openone_{\mathcal{L}\otimes\mathcal{H}_{0}}$
if and only if each $A_{j}$ is coisometric. The result follows.
\end{proof}

Now we apply the lemma to $A_{j}=\lambda\left(  g_{j}\right)  $,
$\mathcal{H}_{0}=\ell^{2}\left(  \mathbb{F}_{d}\right)  $, and it remains to
check that the vector state
\begin{equation}
\Omega_{0}:=\openone\otimes\xi_{e}\in L^{2}\left(  \vphantom{\prod}%
\smash{\prod_{1}^{\infty}}\mathbb{Z}_{d},\mu\right)  \vphantom{\prod
_{1}^{\infty}}\otimes\ell^{2}\left(  \mathbb{F}_{d}\right)  \label{eqFre.13}%
\end{equation}
yields the state $\omega$ in (\ref{eqFre.10}). Let $g_{\alpha_{{}}}%
=g_{\alpha_{1}}g_{\alpha_{2}}\cdots g_{\alpha_{k}}$ for multi-indices
$\alpha=\left(  \alpha_{1}\dots\alpha_{k}\right)  $, $\alpha_{i}\in
\mathbb{Z}_{d}$, $1\leq i\leq k$. Then $T_{\alpha}=S_{\alpha}\otimes
\lambda\left(  g_{\alpha}\right)  $, and%
\begin{align*}%
\ip{\Omega_{0}}{T_{\alpha}^{{}}T_{\gamma}^{\ast}\Omega_{0}}%
&  =%
\ip{\xi_{e}}{\lambda\left( g_{\alpha}^{{}}g_{\gamma
}^{-1}\right) \xi_{e}}%
\ip{\openone}{S_{\alpha}^{{}}S_{\gamma}^{\ast}\openone}%
\\
&  =\delta_{\alpha\gamma}p_{\alpha}%
\end{align*}
where $p^{\alpha}=p_{\alpha_{1}}p_{\alpha_{2}}\cdots p_{\alpha_{k}}$, and
where we used Theorem \ref{ThmGen.1} for the scalar-valued representation
$S_{j}$ in $L^{2}\left(  \prod_{1}^{\infty}\mathbb{Z}_{d}\right)  $ and the
observations from above on the trace of $\mathbb{F}_{d}$. The term%
\[%
\ip{\xi_{e}}{\lambda\left( g_{\alpha}^{{}}g_{\gamma
}^{-1}\right) \xi_{e}}%
=\operatorname*{tr}\left(  g_{\alpha}^{{}}g_{\gamma}^{-1}\right)
\]
is nonzero (and therefore $=1$) if and only if $g_{\alpha}=g_{\gamma}$, i.e.,%
\[
g_{\alpha_{1}}g_{\alpha_{2}}\cdots g_{\alpha_{k}}=g_{\gamma_{1}}g_{\gamma_{2}%
}\cdots g_{\gamma_{l}}.
\]
Since we are in the free group, this happens precisely when $k=l$ and
$g_{\alpha_{1}}=g_{\gamma_{1}},\dots,g_{\alpha_{k}}=g_{\gamma_{k}}$.
\end{proof}

We now turn to the characterization of the cyclic subspace generated by the
representation $\pi_{\omega}$ from Proposition \ref{ProFre.1} when the state
$\omega$ is given as in (\ref{eqFre.8}), (\ref{eqFre.10}).
Let\label{LOSHOmega0}%
\begin{equation}
\mathcal{H}_{\Omega_{0}}:=\left[  \pi_{\omega}\left(  \mathcal{O}_{d}\right)
\Omega_{0}\right]  \label{eqFre.14}%
\end{equation}
where $\Omega_{0}=\openone\otimes\xi_{e}$\label{LOSOmega0_1} and $\pi_{\omega
}\left(  s_{i}\right)  =S_{i}\otimes\lambda\left(  g_{i}\right)  $.

Let $g_{1},\dots,g_{d}$ be the generators of $\mathbb{F}_{d}$, and let
$\mathbb{S}_{d}\subset\mathbb{F}_{d}$\label{LOSSd} be the corresponding free
semigroup, i.e., $\mathbb{S}_{d}$ consists of elements $g_{\alpha_{{}}%
}=g_{\alpha_{1}}g_{\alpha_{2}}\cdots g_{\alpha_{k}}$ (containing no inverses
of any $g_{i}$, $i=1,\dots,d$) indexed by $\alpha=\left(  \alpha_{1}%
\dots\alpha_{k}\right)  $, $\alpha_{i}\in\mathbb{Z}_{d}=\left\{
1,\dots,d\right\}  $, $1\leq i\leq k$, with $k$ depending on $\alpha$. Let
$\mathbb{S}_{d}^{-1}=\left\{  s^{-1}\mid s\in\mathbb{S}_{d}\right\}  $, and
let $\mathcal{H}_{-}:=\left[  \lambda\left(  \mathbb{S}_{d}^{-1}\right)
\xi_{e}\right]  \subset\ell^{2}\left(  \mathbb{F}_{d}\right)  $%
\label{LOSHminus_1} be the closed linear span in $\ell^{2}\left(
\mathbb{F}_{d}\right)  $ of the vectors $\left\{  \lambda\left(
s^{-1}\right)  \xi_{e}\mid s\in\mathbb{S}_{d}\right\}  $. For a multi-index
$\alpha$, and $v\in\mathcal{H}_{-}$, define $\xi_{v}^{\left(  \alpha\right)
}\colon\prod_{1}^{\infty}\mathbb{Z}_{d}\rightarrow\ell^{2}\left(
\mathbb{F}_{d}\right)  $ by%
\begin{align}
\xi_{v}^{\left(  \alpha\right)  }\left(  x\right)   &  =\chi_{\alpha}^{{}%
}\left(  x\right)  \lambda\left(  g_{\alpha}\right)  v\label{eqFre.15}\\
(  &  =\chi_{\alpha}^{{}}\otimes\lambda\left(  g_{\alpha}\right)  v),\qquad
x\in\prod_{1}^{\infty}\mathbb{Z}_{d},\;x=\left(  x_{1},x_{2},\dots\right)
,\nonumber
\end{align}
where%
\[
\chi_{\alpha}^{{}}\left(  x\right)  =\delta_{\alpha_{1}x_{1}}\delta
_{\alpha_{2}x_{2}}\cdots\delta_{\alpha_{k}x_{k}}%
\]
and we use the convention%
\begin{equation}
\xi_{v}^{\left(  \text{\textup{\o}}\right)  }\left(  x\right)  =v.
\label{eqFreNew.17}%
\end{equation}

\begin{lemma}
\label{LemFre.3}The cyclic subspace $\mathcal{H}_{\Omega_{0}}=\left[
\pi_{\omega}\left(  \mathcal{O}_{d}\right)  \Omega_{0}\right]  \subset L_{\mu
}^{2}\left(  \prod_{1}^{\infty}\mathbb{Z}_{d},\ell^{2}\left(  \mathbb{F}%
_{d}\right)  \right)  $\label{LOSHOmega_1} generated by $\Omega_{0}%
=\openone\otimes\xi_{e}$ in the representation defined by $T_{i}=S_{i}%
\otimes\lambda\left(  g_{i}\right)  $, where $\left(  S_{i}\xi\right)  \left(
x\right)  =e^{\frac{1}{2}\beta L_{i}}\chi_{i}^{{}}\left(  x\right)  \xi\left(
\sigma\left(  x\right)  \right)  $, $\xi\in L^{2}\left(  \prod_{1}^{\infty
}\mathbb{Z}_{d},\mu\right)  $, is the closure in $L_{\mu}^{2}\left(  \prod
_{1}^{\infty}\mathbb{Z}_{d},\ell^{2}\left(  \mathbb{F}_{d}\right)  \right)  $
of the linear span of the functions $\xi_{v}^{\left(  \alpha\right)  }$ in
\textup{(\ref{eqFre.15}).}
\end{lemma}

\begin{proof}
From (\ref{eqFre.11})--(\ref{eqFre.12}), we have
\begin{equation}
T_{\alpha}^{{}}T_{\gamma}^{\ast}\left(  \openone\otimes\xi_{e}\right)  \left(
x\right)  =e^{\left(  \beta/2\right)  \left(  L\left(  \alpha\right)
-L\left(  \gamma\right)  \right)  }\chi_{\alpha}^{{}}\left(  x\right)
\lambda\left(  g_{\alpha}^{{}}g_{\gamma}^{-1}\right)  \xi_{e},
\label{eqFreNew.18}%
\end{equation}
where $\alpha=\left(  \alpha_{1}\dots\alpha_{k}\right)  $, $\gamma=\left(
\gamma_{1}\dots\gamma_{l}\right)  $ are multi-indices, and
\begin{equation}
L\left(  \alpha\right)  =\sum_{i}L_{\alpha_{i}}=\sum_{j}\#_{j}\left(
\alpha\right)  L_{j},\qquad\#_{j}\left(  \alpha\right)  =\#\left\{  \alpha
_{i}\mid\alpha_{i}=j\right\}  . \label{eqFreNew.19}%
\end{equation}
Since $v_{\gamma}:=\lambda\left(  g_{\gamma}^{-1}\right)  \xi_{e}%
\in\mathcal{H}_{-}$, the result follows.
\end{proof}

\begin{remark}
\label{RemFre.4}The cyclic subspace $\mathcal{H}_{\Omega_{0}}$%
\label{LOSHOmega_2} is a proper subspace in $L_{\mu}^{2}\left(  \prod
_{1}^{\infty}\mathbb{Z}_{d},\ell^{2}\left(  \mathbb{F}_{d}\right)  \right)  $.
If $i\neq j$, define $\xi\left(  x\right)  =\delta_{i}\left(  x_{1}\right)
\lambda\left(  g_{j}\right)  \xi_{e}$, $x\in\prod_{1}^{\infty}\mathbb{Z}_{d}$.
Then $\xi$ is orthogonal to $\mathcal{H}_{\Omega_{0}}$. For this, we need only
show that $\xi$ is orthogonal to the vectors $\xi_{v}^{\left(  \alpha\right)
}$ in (\ref{eqFre.15}). We have
\begin{align*}%
\ip{\xi}{\xi_{v}^{\left( \alpha\right) }}%
&  =\sum_{r=1}^{d}p_{r}\delta_{i}\left(  r\right)  \delta_{\alpha_{1}}\left(
r\right)  \int_{\prod_{1}^{\infty}\mathbb{Z}_{d}}\chi_{\alpha_{2}\dots
\alpha_{k}}^{{}}\left(  x\right)  \,d\mu\left(  x\right)
\ip{\lambda\left( g_{j}\right) \xi_{e}}{\lambda
\left( g_{\alpha}\right) v}%
\\
&  =p_{\alpha_{1}}\delta_{i\alpha_{1}}p_{\alpha_{2}}\cdots p_{\alpha_{k}}%
\ip{\xi_{e}}{\lambda\left( g_{j}^{-1}g_{\alpha}^{{}}\right) v}%
.
\end{align*}
Since $v\in\mathcal{H}_{-}$, it is enough to show that%
\[
\delta_{i\alpha_{1}}%
\ip{\xi_{e}}{\lambda\left( g_{j}^{-1}g_{\alpha}^{{}}\right) \xi_{s^{-1}}}%
\]
vanishes for $s\in\mathbb{S}_{d}$. The second factor is $\operatorname*{tr}%
\left(  g_{j}^{-1}g_{\alpha}^{{}}s_{{}}^{-1}\right)  $, and this is
nonvanishing only if $g_{j}s=g_{\alpha}$. But the first factor is
$\delta_{i\alpha_{1}}$, so we must have $\alpha_{1}=i$ for the product to be
$\neq0$. Hence $g_{j}s=g_{i}g_{\alpha_{2}}\cdots g_{\alpha_{k}}$ must hold at
a place where the product is $\neq0$. But this is impossible in the free group
$\mathbb{F}_{d}$.
\end{remark}

The vectors $\xi_{v}^{\left(  \alpha\right)  }$ in (\ref{eqFre.15}) are
indexed by the multi-indices $\alpha=\left(  \alpha_{1}\dots\alpha_{k}\right)
$ and vectors $v\in\mathcal{H}_{-}$. Using these we get the following explicit
formula for the operators $T_{i}=\pi_{\omega}\left(  s_{i}\right)  $.

\begin{proposition}
\label{ProFre.5}The generators $T_{i}$, $1\leq i\leq d$, for the cyclic
$\mathcal{O}_{d\,}$-representation of the state defined by
\textup{(\ref{eqFre.10})} act as follows on the vectors $\xi_{v}^{\left(
\alpha\right)  }\left(  x\right)  =\chi_{\alpha}^{{}}\left(  x\right)
\lambda\left(  g_{\alpha}\right)  v$ defined by \textup{(\ref{eqFre.15})} and
\textup{(\ref{eqFreNew.17})} \textup{(}$x\in\prod_{1}^{\infty}\mathbb{Z}_{d}$,
$v\in\mathcal{H}_{-}$\textup{):}%
\begin{align}
T_{i}\left(  \xi_{v}^{\left(  \alpha\right)  }\right)   &  =e^{\left(
\beta/2\right)  L_{i}}\xi_{v}^{\left(  i\alpha\right)  },\label{eqFreNew.20}\\
T_{i}^{\ast}\left(  \xi_{v}^{\left(  \alpha\right)  }\right)   &  =e^{-\left(
\beta/2\right)  L_{i}}\delta_{i\alpha_{1}}\xi_{v}^{\left(  \alpha_{2}%
\alpha_{3}\dots\alpha_{k}\right)  },\label{eqFreNew.21}\\
&  \text{\qquad if }\alpha\neq\text{\textup{\o}},\nonumber\\
T_{i}^{\ast}\left(  \xi_{v}^{\left(  \text{\textup{\o}}\right)  }\right)   &
=e^{-\left(  \beta/2\right)  L_{i}}\xi_{\lambda\left(  g_{i}^{-1}\right)
v}^{\left(  \text{\textup{\o}}\right)  },\label{eqFreNew.22}\\
T_{\gamma}^{{}}T_{\gamma}^{\ast}\left(  \xi_{v}^{\left(  \alpha\right)
}\right)   &  =\chi_{\gamma}^{{}}\left(  \alpha\right)  \xi_{v}^{\left(
\alpha\right)  }\label{eqFreNew.23}\\
&  \text{\qquad if }\left|  \gamma\right|  \leq\left|  \alpha\right|
\text{\textup{(}i.e., }l\leq k,\;\gamma=\left(  \gamma_{1}\dots\gamma
_{l}\right)  ,\;\alpha=\left(  \alpha_{1}\dots\alpha_{k}\right)
\text{\textup{),}}\nonumber\\%
\intertext{and}%
\left(  T_{\gamma}^{{}}T_{\gamma}^{\ast}\right)  \xi_{v}^{\left(
\text{\textup{\o}}\right)  }  &  =\chi_{\gamma}^{{}}\otimes
v\label{eqFreNew.24}\\
&  \text{\qquad with }\gamma=\left(  \gamma_{1}\dots\gamma_{l}\right)
.\nonumber
\end{align}
Hence, if $l>k$,
\begin{equation}
T_{\eta}^{{}}T_{\gamma}^{\ast}\xi_{v}^{\left(  \alpha\right)  }=e^{\left(
\beta/2\right)  \left(  L\left(  \eta\right)  -L\left(  \gamma\right)
\right)  }\delta_{\gamma_{1}\alpha_{1}}\cdots\delta_{\gamma_{k}\alpha_{k}}%
\xi_{\lambda\left(  g_{\gamma_{l}}^{-1}\cdots g_{\gamma_{k+1}}^{-1}\right)
v}^{\left(  \eta\right)  }. \label{eqFreNew.25}%
\end{equation}
\end{proposition}

\begin{proof}
We compute the action of $T_{i}^{{}}$ and $T_{i}^{\ast}$ directly from the
formulas given in Proposition \ref{ProFre.1}. We have%
\begin{align*}
T_{i}^{{}}\xi_{v}^{\left(  \alpha\right)  }\left(  x\right)   &  =e^{\frac
{1}{2}\beta L_{i}}\chi_{i}^{{}}\left(  x\right)  \chi_{\alpha}^{{}}\left(
\sigma\left(  x\right)  \right)  \lambda\left(  g_{i}\right)  \lambda\left(
g_{\alpha}\right)  v\\
&  =e^{\frac{1}{2}\beta L_{i}}\delta_{i_{{}}}\left(  x_{1}\right)
\delta_{\alpha_{1}}\left(  x_{2}\right)  \cdots\delta_{\alpha_{k}}\left(
x_{k+1}\right)  \lambda\left(  g_{i}g_{\alpha}\right)  v\\
&  =e^{\frac{1}{2}\beta L_{i}}\chi_{\left(  i\alpha\right)  }^{{}}\left(
x\right)  \lambda\left(  g_{\left(  i\alpha\right)  }\right)  v\\
&  =e^{\left(  \beta/2\right)  L_{i}}\xi_{v}^{\left(  i\alpha\right)  }\left(
x\right)  ,
\end{align*}
proving (\ref{eqFreNew.20}), and%
\begin{align*}
T_{i}^{\ast}\xi_{v}^{\left(  \alpha\right)  }\left(  x\right)   &
=e^{-\left(  \beta/2\right)  L_{i}}\chi_{\alpha}^{{}}\left(  \sigma_{i}^{{}%
}\left(  x\right)  \right)  \lambda\left(  g_{i}^{-1}\right)  \lambda\left(
g_{\alpha}^{{}}\right)  v\\
&  =e^{-\left(  \beta/2\right)  L_{i}}\delta_{i\alpha_{1}}\chi_{\left(
\alpha_{2}\alpha_{3}\dots\alpha_{k}\right)  }^{{}}\left(  x\right)
\lambda\left(  g_{i}^{-1}g_{\alpha}^{{}}\right)  v\\
&  =e^{-\left(  \beta/2\right)  L_{i}}\delta_{i\alpha_{1}}\chi_{\left(
\alpha_{2}\dots\alpha_{k}\right)  }^{{}}\left(  x\right)  \lambda\left(
g_{i}^{-1}g_{i}^{{}}g_{\alpha_{2}}^{{}}\cdots g_{\alpha_{k}}^{{}}\right)  v\\
&  =e^{-\left(  \beta/2\right)  L_{i}}\delta_{i\alpha_{1}}\chi_{\left(
\alpha_{2}\dots\alpha_{k}\right)  }^{{}}\left(  x\right)  \lambda\left(
g_{\left(  \alpha_{2}\dots\alpha_{k}\right)  }\right)  v\\
&  =e^{-\left(  \beta/2\right)  L_{i}}\delta_{i\alpha_{1}}\xi_{v}^{\left(
\alpha_{2}\dots\alpha_{k}\right)  }\left(  x\right)  ,
\end{align*}
proving (\ref{eqFreNew.21}). The stated formulas for $T_{i}^{\ast}\xi
_{v}^{\left(  \text{\textup{\o}}\right)  }$ and $T_{\gamma}^{{}}T_{\gamma
}^{\ast}\xi_{v}^{\left(  \alpha\right)  }$, $v\in\mathcal{H}_{-}$, result from
the following covariance principle: $\chi_{i}^{{}}\otimes v=\xi_{\lambda
\left(  g_{i}^{-1}\right)  v}^{\left(  i\right)  }$, and, more generally,
$\chi_{\alpha}^{{}}\otimes v=\xi_{\lambda\left(  g_{\alpha}^{-1}\right)
v}^{\left(  \alpha\right)  }$. The formula (\ref{eqFreNew.24}) is a special
case of (\ref{eqGen.34}).
\end{proof}

\begin{corollary}
\label{CorFre.6}Let $v$ be an arbitrary vector in $\mathcal{H}_{-}$. As
$\gamma=\left(  \gamma_{1}\dots\gamma_{l}\right)  $ ranges over all
$\mathbb{Z}_{d\,}$-multi-indices, the closed linear span of $\left\{  \left(
T_{\gamma}^{{}}T_{\gamma}^{\ast}\right)  \xi_{v}^{\left(  \text{\textup{\o}%
}\right)  }\right\}  _{\gamma}$ in $L_{\mu}^{2}\left(  \prod_{1}^{\infty
}\mathbb{Z}_{d},\ell^{2}\left(  \mathbb{F}_{d}\right)  \right)  $ is $L_{\mu
}^{2}\left(  \prod_{1}^{\infty}\mathbb{Z}_{d}\right)  \otimes v$ where $\mu$
is still the $\left(  p_{i}\right)  $-product measure on $\prod_{1}^{\infty
}\mathbb{Z}_{d}$.
\end{corollary}

\begin{proof}
This is immediate from (\ref{eqFreNew.24}).
\end{proof}

\begin{remark}
\label{RemFre.7}Let $\mathcal{H}_{\Omega_{0}}$\label{LOSOmega0_2} be the
cyclic subspace of the representation of $\mathcal{O}_{d}$ induced from the
state $\omega_{\left(  p\right)  }\left(  s_{\alpha}^{{}}s_{\gamma}^{\ast
}\right)  =p_{\alpha}\delta_{\alpha\gamma}$. Then $L^{2}\left(  \prod
_{1}^{\infty}\mathbb{Z}_{d}\right)  \otimes\mathcal{H}_{-}$ is a proper
subspace in $\mathcal{H}_{\Omega_{0}}$. For example, the vector $\chi_{i}^{{}%
}\otimes\xi_{g_{i}}$ is in $\mathcal{H}_{\Omega_{0}}\ominus\left(
L^{2}\otimes\mathcal{H}_{-}\right)  $\label{LOSHOmega_3}.
\end{remark}

\begin{proof}
We check that $%
\ip{\chi_{i}^{{}}\otimes\xi_{g_{i}}}{f\otimes v}%
=0$ for all $f\in L^{2}\left(  \prod_{1}^{\infty}\mathbb{Z}_{d}\right)  $, and
$v\in\mathcal{H}_{-}$. We may assume that $v=\xi_{s^{-1}}$ for $s\in
\mathbb{S}_{d}$ ($=$ the free semigroup on the generators $\left\{
g_{j}\right\}  _{j=1}^{d}$). The inner product is%
\[
p_{i}\int f\left(  i,x\right)  \,d\mu_{\left(  p\right)  }\left(  x\right)
\ip{\xi_{g_{i}}}{\xi_{s^{-1}}}%
,
\]
and $%
\ip{\xi_{g_{i}}}{\xi_{s^{-1}}}%
=\operatorname*{tr}\left(  g_{i}^{-1}s_{{}}^{-1}\right)  =0$, since there is
no solution $s\in\mathbb{S}_{d}$ to the equation $sg_{i}=e$.
\end{proof}

Summarizing Remarks \ref{RemFre.4} and \ref{RemFre.7} we have%
\[
L_{\mu}^{2}\otimes\mathcal{H}_{-}^{{}}\subsetneqq\mathcal{H}_{\Omega_{0}}^{{}%
}\subsetneqq L_{\mu}^{2}\otimes\ell_{{}}^{2}\left(  \mathbb{F}_{d}\right)  .
\]

\begin{definition}
\label{DefFreNew.8}Let $\mathcal{H}$ and $\mathcal{H}_{-}$ be Hilbert spaces,
and let $\mathcal{M}$ be a set. We say that $\mathcal{H}$ is \emph{fibered}
over $\mathcal{M}$ with fibers isomorphic to $\mathcal{H}_{-}$ if there are
isometries $i_{\alpha}$, indexed by $\alpha\in\mathcal{M}$, $i_{\alpha}%
\colon\mathcal{H}_{-}\rightarrow\mathcal{H}$, such that $\mathcal{H}$ is the
closed linear span of $\left\{  i_{\alpha}\left(  \mathcal{H}_{-}\right)
\mid\alpha\in\mathcal{M}\right\}  $.
\end{definition}

\noindent\textbf{The Fibered Space.} Let $\omega_{\left(  p\right)  }$ be the
state from above, and let $\mathcal{H}_{\Omega_{0}}$ be the cyclic space of
the $\mathcal{O}_{d}$-representation. Let $\mathcal{M}\left(  \mathbb{Z}%
_{d}\right)  $ be the set of all multi-indices formed from $\mathbb{Z}_{d}$.
Then $\mathcal{H}_{\Omega_{0}}$ is fibered \textup{(}as a Hilbert
space\textup{)} over $\mathcal{M}\left(  \mathbb{Z}_{d}\right)  $, the fiber
over each $\alpha$ in $\mathcal{M}\left(  \mathbb{Z}_{d}\right)  $ is a copy
of $\mathcal{H}_{-}$.

To prove this, let $\alpha\in\mathcal{M}\left(  \mathbb{Z}_{d}\right)  $, and
define $\mathcal{H}\left(  \alpha\right)  \cong\mathcal{H}_{-}$ by%
\[
\mathcal{H}\left(  \alpha\right)  =\left\{  \chi_{\alpha}^{{}}\otimes
\lambda\left(  g_{\alpha}\right)  v\mid v\in\mathcal{H}_{-}\right\}  .
\]
The isomorphism $\mathcal{H}\left(  \alpha\right)  \cong\mathcal{H}_{-}$ and
the isometries $i_{\alpha}$ are then made explicit by using the scale given by
the following identity:%
\[
\left\|  \chi_{\alpha}^{{}}\otimes\lambda\left(  g_{\alpha}\right)  v\right\|
_{\mathcal{H}_{\Omega_{0}}}^{2}=p^{\alpha}\left\|  v\right\|  _{\mathcal{H}%
_{-}}^{2},
\]
where $p^{\alpha}=p_{\alpha_{1}}p_{\alpha_{2}}\cdots p_{\alpha_{k}}>0$. The
convention for the empty index \textup{\o} in $\mathcal{M}\left(
\mathbb{Z}_{d}\right)  $ is that the fiber $\mathcal{H}\left(
\text{\textup{\o}}\right)  $ over \textup{\o} is\label{LOSHminus_2}%
\[
\mathcal{H}\left(  \text{\textup{\o}}\right)  =\left\{  \openone\otimes v\mid
v\in\mathcal{H}_{-}\right\}
\]
where $\openone=\chi_{\text{\textup{\o}}}^{{}}$ denotes the constant function
``one'' in $L^{2}\left(  \prod_{1}^{\infty}\mathbb{Z}_{d}\right)  $. The
action of $T_{i}^{{}}$, $T_{i}^{\ast}$ on the fibers is given by Proposition
\ref{ProFre.5}. In particular it follows from (\ref{eqGen.34}) that the action
of $L^{\infty}\left(  \prod_{1}^{\infty}\mathbb{Z}_{d}\right)  $ is given by%
\[
\pi_{\omega}\left(  f\right)  \xi_{v}^{\left(  \alpha\right)  }=\left(
f\chi_{\alpha}^{{}}\right)  \otimes\lambda\left(  g_{\alpha}\right)  v,\qquad
v\in\mathcal{H}_{-}.
\]

\chapter[$\beta$-KMS states]{\label{KMS}$\beta$-KMS states for one-parameter
subgroups%
\protect\frogleg
of the action of $\mathbb{T}^{d}$ on $\mathcal{O}_{d}$}

Consider the action of $\mathbb{T}^{d}$\label{LOSTd} on $\mathcal{O}_{d}%
$\label{LOSOd_5} given by%
\begin{equation}
\sigma\left(  z_{1},\dots,z_{d}\right)  \left(  s_{i}\right)  =z_{i}s_{i}.
\label{eqKMS.1}%
\end{equation}
If $L=\left(  L_{1},\dots,L_{d}\right)  \in\mathbb{R}^{d}$, consider the
one-parameter group\label{LOSsigmatL_2}%
\begin{equation}
\sigma_{t}^{L}\left(  x\right)  =\sigma\left(  e^{itL_{1}},\dots,e^{itL_{d}%
}\right)  \left(  x\right)  . \label{eqKMS.2}%
\end{equation}
In general, if $\mathfrak{A}$ is a $C^{\ast}$-algebra and $t\mapsto\sigma_{t}$
is a one-parameter group of $\ast$-au\-to\-mor\-phisms of $\mathfrak{A}$, and
$\beta\in\mathbb{R}$, recall that a state $\omega$ over $\mathfrak{A}$ is
defined to be a $\sigma$-KMS state at value $\beta$, or a $\left(
\sigma,\beta\right)  $-KMS state\label{LOSsigmaLbetaKMSstate_3} if%
\begin{equation}
\omega\left(  x\sigma_{i\beta}\left(  y\right)  \right)  =\omega\left(
yx\right)  \label{eqKMS.3}%
\end{equation}
for all $x$, $y$ in a norm-dense $\sigma$-invariant $\ast$-algebra of $\sigma
$-analytic elements of $\mathfrak{A}$ (see \cite[Section 5.3.1]{BrRoII} for
several alternative formulations of this condition). It is well known that if
$L=\left(  1,1,\dots,1\right)  $, so that $\sigma$ is the so-called gauge
group, the group $\sigma^{L}$ has a KMS state at value $\beta$ if and only if
$\beta=\log d$, and this state is unique and is given by%
\begin{equation}
\omega\left(  s_{\alpha}^{{}}s_{\gamma}^{\ast}\right)  =\delta_{\alpha\gamma
}d^{-\left|  \alpha\right|  }; \label{eqKMS.4}%
\end{equation}
see \cite[Example 5.3.27]{BrRoII}, \cite{BEH80}, or \cite{OlPe78}. We first
note that the latter result has an easy extension to more general
one-parameter subgroups.

\begin{proposition}
\label{ProKMS.1}The one-parameter group $\sigma^{L}$ admits a KMS state at
some value $\beta$ if and only if $L_{1},L_{2},\dots,L_{d}$ all are nonzero
and have the same sign. This value $\beta$ is then unique and is given as the
real solution of
\begin{equation}
\sum_{k=1}^{d}e^{-\beta L_{k}}=1. \label{eqKMS.5}%
\end{equation}
The $\sigma^{L}$-KMS state $\omega$ at value $\beta$ is then also unique, and
is given by\label{LOSomegasalphasgamma}%
\begin{equation}
\omega\left(  s_{\alpha}^{{}}s_{\gamma}^{\ast}\right)  =\delta_{\alpha\gamma
}e^{-\beta\sum_{k=1}^{\left|  \alpha\right|  }L_{\alpha_{k}}}. \label{eqKMS.6}%
\end{equation}
\end{proposition}

\begin{proof}
If $\omega$ is a KMS state at value $\beta$, then%
\begin{equation}
\omega\left(  s_{\alpha}^{{}}s_{\gamma}^{\ast}\right)  =\omega\left(
s_{\gamma}^{\ast}\sigma_{i\beta}\left(  s_{\alpha}^{{}}\right)  \right)
=e^{-\beta\sum_{k=1}^{\left|  \alpha\right|  }L_{\alpha_{k}}}\omega\left(
s_{\gamma}^{\ast}s_{\alpha}^{{}}\right)  . \label{eqKMS.7}%
\end{equation}
If $\alpha=\left(  k\right)  $, $\gamma=\left(  j\right)  $, this says%
\begin{equation}
\omega\left(  s_{k}^{{}}s_{j}^{\ast}\right)  =\delta_{kj}e^{-\beta L_{k}}.
\label{eqKMS.8}%
\end{equation}
But%
\begin{equation}
1=\sum_{k=1}^{d}\omega\left(  s_{k}^{{}}s_{k}^{\ast}\right)  =\sum_{k=1}%
^{d}e^{-\beta L_{k}} \label{eqKMS.9}%
\end{equation}
and hence $\beta$ is a solution of (\ref{eqKMS.5}). But this equation has
solutions $\beta$ if and only if all $L_{k}$ are nonzero, and all have the
same sign; and, in that case, the solution $\beta$ is unique. For
definiteness, assume that all $L_{k}$ are positive, and then the solution
$\beta$ of (\ref{eqKMS.5}) is also positive. Because of the Cuntz relations,
the element $s_{\gamma}^{\ast}s_{\alpha}^{{}}$ is either $\openone$ (if
$\gamma=\alpha$), $0$, or of the form $s_{\delta}^{{}}$ or $s_{\delta}^{\ast}$
for some $\delta$. But from (\ref{eqKMS.7}), we have%
\begin{equation}
\omega\left(  s_{\delta}\right)  =e^{-\beta\sum_{k=1}^{\left|  \delta\right|
}L_{k}}\omega\left(  s_{\delta}\right)  \label{eqKMS.10}%
\end{equation}
and thus $\omega\left(  s_{\delta}\right)  =0$ for all nonempty strings
$\delta$. Hence it follows from (\ref{eqKMS.7}) again that%
\begin{equation}
\omega\left(  s_{\alpha}^{{}}s_{\gamma}^{\ast}\right)  =e^{-\beta\sum
_{k=1}^{\left|  \alpha\right|  }L_{\alpha_{k}}}\delta_{\alpha\gamma},
\label{eqKMS.11}%
\end{equation}
which is (\ref{eqKMS.6}). But this expression does indeed define a state by
Proposition \ref{ProFre.1}. The case that all $L_{k}$ are negative is treated
similarly, so this proves Proposition \ref{ProKMS.1}.
\end{proof}

The KMS states and the one-parameter subgroups of automorphisms were also used
in recent papers \cite{EvKi97,KiKu96,KiKu97} where crossed products
$\mathcal{O}_{d}\rtimes_{\sigma}\mathbb{R}$ were studied. The states
(\ref{eqKMS.6}) seem to have first appeared in \cite{Cun80}, \cite{Eva80} and
\cite{Eva82}.

The result in Proposition \ref{ProKMS.1} is also related to results in
\cite{Lac98}, where KMS states for one-parameter subgroups of the dual actions
of actions of lattice semigroups of endomorphisms scaling tracial states of a
given $C^{\ast}$-algebra are considered. It turns out that the KMS states have
non-trivial symmetries apart from invariance under the one-parameter
semigroup, and in particular an ``explanation'' is given of the fact that our
states given by (\ref{eqKMS.6}) have the $\delta_{\alpha\gamma}$ term which
forces them to live on the maximal abelian subalgebra $\mathcal{D}_{d}$ which
is the closure of the linear span of the monomials $s_{\alpha}^{{}}s_{\alpha
}^{\ast}$, i.e., the fixed-point algebra of the canonical coaction of
$\mathbb{F}_{d}$.

Let us comment a bit on the representations defined by the state $\omega$ in
(\ref{eqKMS.8}). For definiteness, assume that $L_{1},\dots,L_{d}$ are all
strictly positive. Let $\mathfrak{A}_{L}$ be the fixed-point algebra of the
modular automorphism group $\sigma^{\left(  L\right)  }$. We prove in
Proposition \ref{ProSub.1} and Remark \ref{RemSubNew.11}, below, that
$\mathfrak{A}_{L}$ is an AF-algebra. We consider the following algebras:

\begin{enumerate}
\item [$\mathcal{O}_{d}$]$=$\label{LOSOd_6} closed linear span of all
$s_{\alpha}^{{}}s_{\gamma}^{\ast}$,

\item[$\mathcal{O}_{d}^{\mathbb{T}}$] $=\operatorname*{UHF}_{d}=$%
\label{LOSOdT} closed linear span of all $s_{\alpha}^{{}}s_{\gamma}^{\ast}$
with $\left|  \alpha\right|  =\left|  \gamma\right|  $\newline $=$ fixed-point
algebra of the gauge action of $\mathbb{T}$,

\item[$\mathfrak{A}_{L}$] $=$\label{LOSAL_2} closed linear span of all
$s_{\alpha}^{{}}s_{\gamma}^{\ast}$ with $L\left(  \alpha\right)  =L\left(
\gamma\right)  $\newline $=$ fixed-point algebra of the action $\sigma
^{\left(  L\right)  }$, where $L\left(  \alpha\right)  $ is defined by
(\ref{eqSub.2}),

\item[$\mathcal{O}_{d}^{\mathbb{T}^{d}}$] $=\operatorname*{GICAR}_{d}%
=$\label{LOSOdTd} closed linear span of all $s_{\alpha}^{{}}s_{\gamma}^{\ast}$
where $\left|  \alpha\right|  =\left|  \gamma\right|  $ and $\gamma$ is a
permutation of $\alpha$\newline $=$ fixed-point algebra of the gauge action of
$\mathbb{T}^{d}$,

\item[$\mathcal{D}_{d}$] $=$\label{LOSDd_3} closed linear span of all
$s_{\alpha}^{{}}s_{\alpha}^{\ast}$ (see (\ref{eqGen.22}))\newline $=$
fixed-point algebra of the coaction of $\mathbb{F}_{d}$.
\end{enumerate}

We have the following inclusions:%
\[
\begin{picture}(295,120)(-201.5,-60)
\put(-173,-42){\includegraphics
[bb=0 77 240 161,height=84bp,width=240bp]{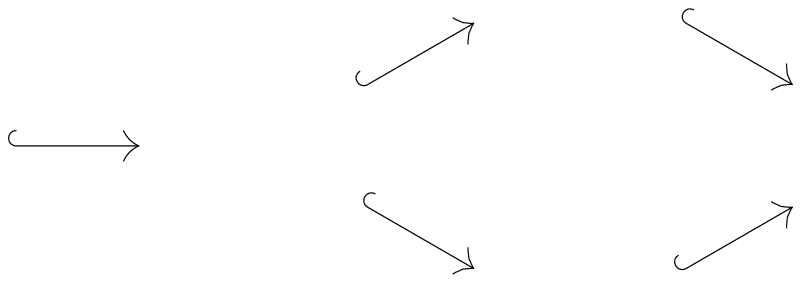}}
\put(-201.5,0){\makebox(295,0){\hbox to 0bp
{\hss$\displaystyle\mathcal{D}_{d}$\hss}\hskip108bp\hbox to 0bp
{\hss$\displaystyle\mathcal{O}_{d}^{\mathbb{T}^{d}}$\hss}\hskip187bp\hbox
to 0bp
{\hss$\displaystyle\mathcal{O}_{d},$\hss}}}
\put(0,54){\makebox(0,0){$\displaystyle\mathcal{O}_{d}^{\mathbb{T}}$}}
\put(0,-54){\makebox(0,0){$\displaystyle\mathfrak{A}_{L}$}}
\end{picture}%
\]
where $\mathfrak{A}_{L}=\mathcal{O}_{d}^{\mathbb{T}^{d}}$ if and only if
$L_{1},L_{2},\dots,L_{d}$ are rationally independent and $\mathfrak{A}%
_{L}=\mathcal{O}_{d}^{\mathbb{T}}$ if and only if $L_{1}=L_{2}=\dots=L_{d}$.
In general $\mathfrak{A}_{L}$ has a skew position relative to $\mathcal{O}%
_{d}^{\mathbb{T}}=\operatorname*{UHF}\nolimits_{d}$.

We will here only analyze the representations given by the state
(\ref{eqKMS.6}) in the case that $L=\left(  L_{1},\dots,L_{d}\right)  $ is in
a class extensively studied in the remainder of the memoir: Each pair
$L_{i},L_{j}$ is rationally dependent. We have to refer to results in Chapters
\ref{Sub} and \ref{Rem}. By a renormalization (see remarks after
(\ref{eqSub.1})) we may assume that the $L_{i}$'s are positive (nonzero)
integers, and that $\gcd\left(  L_{1},\dots,L_{d}\right)  =1$. Then the
associated one-parameter group $\sigma^{\left(  L\right)  }$ is periodic with
period $2\pi$, so we may view $\sigma^{\left(  L\right)  }$ as a
representation of $\mathbb{T}$ in the automorphism group of $\mathcal{O}_{d}$.
Since $\gcd\left(  L_{1},\dots,L_{d}\right)  =1$, it follows, from the
Euclidean algorithm and (\ref{eqSub.5}), that the spectral
subspaces\label{LOSOdsigmaLn}%
\begin{equation}
\mathcal{O}_{d}^{\sigma}\left(  n\right)  =\left\{  x\in\mathcal{O}_{d}%
\mid\sigma_{t}^{\left(  L\right)  }\left(  x\right)  =e^{int}x\right\}
\label{eqKMSNew.13}%
\end{equation}
are nonzero for all $n\in\mathbb{\hat{T}=Z}$ \cite{BrRoI}. But we argue in the
beginning of Chapter \ref{Rem} that the fixed-point algebra\label{LOSAL_3}%
\begin{equation}
\mathfrak{A}_{L}=\mathcal{O}_{d}^{\sigma}\left(  0\right)  \label{eqKMSNew.14}%
\end{equation}
is a simple unital AF-algebra with a unique trace state $\tau=\omega
|_{\mathfrak{A}_{L}}$. Since all the spaces $\mathcal{O}_{d}^{\sigma}\left(
n\right)  \mathcal{O}_{d}^{\sigma}\left(  n\right)  ^{\ast}$ are ideals in
$\mathfrak{A}_{L}$, it follows further that%
\begin{equation}
\mathcal{O}_{d}^{\sigma}\left(  n\right)  \mathcal{O}_{d}^{\sigma}\left(
n\right)  ^{\ast}=\mathfrak{A}_{L} \label{eqKMSNew.15}%
\end{equation}
for all $n\in\mathbb{N}$. If $\left(  \hat{\sigma},\hat{\omega}\right)  $
denotes the pair of extensions of $\left(  \sigma,\omega\right)  $ to the weak
closure $\mathcal{O}_{d}^{\prime\prime}$ of $\mathcal{O}_{d}$ in the
representation defined by $\omega$, it follows from (\ref{eqKMSNew.15}) that
the $\Gamma$-spectrum of the extension is%
\begin{equation}
\Gamma\left(  \hat{\sigma}\right)  =\mathbb{\hat{T}=Z.} \label{eqKMSNew.16}%
\end{equation}
Also, since $\omega$ is a $\sigma^{\left(  L\right)  }$-KMS state at value
$\beta$, where $\beta$ is defined by (\ref{eqKMS.5}), it follows that%
\begin{equation}
t\longmapsto\hat{\sigma}\left(  -t\beta\right)  \label{eqKMSNew.17}%
\end{equation}
is the modular automorphism group of $\hat{\omega}$; see \cite[Definition
5.3.1 and Theorem 5.3.10]{BrRoII}.

Now, since $\tau=\omega|_{\mathfrak{A}_{L}}$ is the unique trace state on
$\mathfrak{A}_{L}$, it defines a type $\mathrm{II}_{1}$ factor representation
of $\mathfrak{A}_{L}$. Using (\ref{eqKMSNew.15}) in the form%
\begin{equation}
\mathcal{O}_{d}^{\sigma}\left(  n\right)  \mathfrak{A}_{L}\mathcal{O}%
_{d}^{\sigma}\left(  n\right)  ^{\ast}=\mathfrak{A}_{L}, \label{eqKMSNew.18}%
\end{equation}
it follows that the representation of $\mathcal{O}_{d}$ defined by $\omega$,
restricted to $\mathfrak{A}_{L}$, is quasiequivalent to the trace
representation, and in particular $\mathfrak{A}_{L}^{\prime\prime}$ is
isomorphic to the hyperfinite $\mathrm{II}_{1}$ factor,%
\begin{equation}
\mathfrak{A}_{L}^{\prime\prime}\cong\mathcal{R}. \label{eqKMSNew.19}%
\end{equation}
Using the definition (\ref{eqSub.2}), we see that%
\begin{equation}
s_{\alpha}^{{}}s_{\gamma}^{\ast}\in\mathcal{O}_{d}^{\sigma}\left(  n\right)
\iff L\left(  \alpha\right)  -L\left(  \gamma\right)  =n\text{\qquad for all
multi-indices }\alpha,\gamma. \label{eqKMSNew.20}%
\end{equation}
Thus, using (\ref{eqKMS.6}), we see that, if $y\in\mathfrak{A}_{L}$ and
$x\in\mathcal{O}_{d}^{\sigma}\left(  n\right)  $ with $x^{\ast}xy=y=yx^{\ast
}x$, then%
\begin{equation}
\omega\left(  xyx^{\ast}\right)  =e^{-\beta n}\omega\left(  y\right)  ,
\label{eqKMSNew.21}%
\end{equation}
as follows from (\ref{eqKMSNew.17}), (\ref{eqKMSNew.19}), and
(\ref{eqKMSNew.21}). To see this, consider for example $x=s_{\alpha}^{{}%
}s_{\gamma}^{\ast}$ with $L\left(  \alpha\right)  -L\left(  \gamma\right)
=n$. Let $y$ be the initial projection of the partial isometry $x$, i.e.,%
\[
y=x^{\ast}x=s_{\gamma}^{{}}s_{\alpha}^{\ast}s_{\alpha}^{{}}s_{\gamma}^{\ast
}=s_{\gamma}^{{}}s_{\gamma}^{\ast}\in\mathfrak{A}_{L}.
\]
Then%
\[
\omega\left(  y\right)  =e^{-\beta L\left(  \gamma\right)  }%
\]
from (\ref{eqKMS.6}). But%
\begin{align*}
xyx^{\ast}  &  =s_{\alpha}^{{}}s_{\gamma}^{\ast}s_{\gamma}^{{}}s_{\gamma
}^{\ast}s_{\gamma}^{{}}s_{\alpha}^{\ast}\\
&  =s_{\alpha}^{{}}s_{\alpha}^{\ast},
\end{align*}
and so%
\begin{align*}
\omega\left(  xyx^{\ast}\right)   &  =e^{-\beta L\left(  \alpha\right)  }\\
&  =e^{-\beta\left(  L\left(  \alpha\right)  -L\left(  \gamma\right)  \right)
}e^{-\beta L\left(  \gamma\right)  }\\
&  =e^{-\beta n}\omega\left(  y\right)  .
\end{align*}
An elaboration of this computation proves (\ref{eqKMSNew.21}).

It now follows from (\ref{eqKMSNew.17}), (\ref{eqKMSNew.19}),
(\ref{eqKMSNew.21}), and \cite{Con73} or \cite[Proposition 29.1]{Str81} that
$\mathcal{O}_{d}^{\prime\prime}$ is the hyperfinite $\mathrm{III}_{e^{-\beta}%
}$-factor. The factor $\mathcal{O}_{d}^{\prime\prime}$ can be written as the
crossed product of $\mathfrak{A}_{L}^{\prime\prime}\otimes\mathcal{B}\left(
\mathcal{H}\right)  $ ($=$ the hyperfinite $\mathrm{II}_{\infty}$-factor) by
an automorphism scaling the trace by $e^{-\beta}$, something which is
reflected in (\ref{eqKMSNew.21}). This automorphism is described in the end of
Chapter \ref{Rem}, and should not be confused with a stabilized version of the
endomorphism $\lambda=\sum_{i}s_{i}^{{}}\cdot s_{i}^{\ast}$, except when
$L_{1}=L_{2}=\dots=L_{d}=1$.

We defer a detailed analysis of the case when the $L_{i}$'s are not pairwise
rationally dependent to a later paper. Although $\mathfrak{A}_{L}$ is still an
AF-algebra\label{LOSAFalgebra_3}, it is no longer simple, and it does not have
a unique trace state. For example if $d=2$ and $L_{1}$, $L_{2}$ are rationally
independent, then $\mathfrak{A}_{L}$ is the GICAR algebra which is a
primitive, non-simple $C^{\ast}$-algebra, and the extremal boundary of the
compact convex set of trace states is homeomorphic to the unit interval
$\left[  0,1\right]  $; see \cite{Bra72}, \cite{Ren80}, or \cite[Examples
III.5.5 and IV.3.7]{Dav96}. Hence the analysis of the algebras $\mathfrak
{A}_{L}$ will be radically different for general $L$ than in the remaining
chapters of the present memoir.

\chapter{\label{Sub}Subalgebras of $\mathcal{O}_{d}$}

In this chapter we will study the fixed-point subalgebras of $\mathcal{O}_{d}$
under the one-parameter groups $\sigma=\sigma^{\left(  L\right)  }$ of
automorphisms defined by\label{LOSsigmatL_3}%
\begin{equation}
\sigma_{t}^{\left(  L\right)  }\left(  s_{j}\right)  =e^{itL_{j}}s_{j},\qquad
j=1,\dots,d, \label{eqSub.1}%
\end{equation}
where we will assume that all the $L_{j}$ have the same sign and any pair
$\left(  L_{j},L_{k}\right)  $ is rationally dependent. By a renormalization
of the variable $t$ we may, and will, assume that all the $L_{k}$ are positive
integers and that the greatest common divisor of $L_{1},\dots,L_{d}$ is $1$.
The group $\sigma_{t}^{\left(  L\right)  }$ is then periodic with period
$2\pi$. If $\alpha=\left(  \alpha_{1}\dots\alpha_{k}\right)  $ is a
multi-index with $\alpha_{m}\in\mathbb{Z}_{d}$, we define the weight
function\label{LOSLalpha_1}%
\begin{equation}
L\left(  \alpha\right)  =\sum_{j=1}^{d}\#_{j}\left(  \alpha\right)  L_{j}%
=\sum_{m=1}^{k}L_{\alpha_{m}} \label{eqSub.2}%
\end{equation}
where\label{LOShashjalpha}
\begin{equation}
\#_{j}\left(  \alpha\right)  =\#\left\{  \alpha_{i}\mid\alpha_{i}=j\right\}
\label{eqSub.3}%
\end{equation}
and using the standard multi-index notation\label{LOSsalpha_2}%
\begin{equation}
s_{\alpha}=s_{\alpha_{1}}s_{\alpha_{2}}\cdots s_{\alpha_{k}} \label{eqSub.4}%
\end{equation}
we have%
\begin{equation}
\sigma_{t}\left(  s_{\alpha}^{{}}s_{\gamma}^{\ast}\right)  =e^{it\left(
L\left(  \alpha\right)  -L\left(  \gamma\right)  \right)  }s_{\alpha}^{{}%
}s_{\gamma}^{\ast}. \label{eqSub.5}%
\end{equation}
Since these elements span $\mathcal{O}_{d}$, it follows that the eigenspace
$\mathcal{O}_{d}^{\sigma}\left(  n\right)  $ in $\mathcal{O}_{d}$ is the
closed linear span of the $s_{\alpha}^{{}}s_{\gamma}^{\ast}$ with $L\left(
\alpha\right)  -L\left(  \gamma\right)  =n$. In particular, the fixed-point
algebra $\mathfrak{A}_{L}=\mathcal{O}_{d}^{\sigma}=\mathcal{O}_{d}^{\sigma
}\left(  0\right)  $ is the closure of the linear span of $s_{\alpha}^{{}%
}s_{\gamma}^{\ast}$ with $L\left(  \alpha\right)  =L\left(  \gamma\right)  $.

The first result on $\mathfrak{A}_{L}$ is that it is an AF-algebra, i.e., it
is the closure of the union of an increasing sequence of finite-dimensional subalgebras:

\begin{proposition}
\label{ProSub.1}Let $L_{1},\dots,L_{d}$ be integers and consider the periodic
one-pa\-ram\-e\-ter group $\sigma$ of $\ast$-au\-to\-mor\-phisms of
$\mathcal{O}_{d}$ defined by%
\begin{equation}
\sigma_{t}\left(  S_{j}\right)  =e^{itL_{j}}S_{j}. \label{eqSub.6}%
\end{equation}
Then the following conditions are equivalent.

\begin{enumerate}
\item \label{ProSub.1(1)}The fixed-point algebra $\mathfrak{A}_{L}$ is an AF-algebra.

\item \label{ProSub.1(2)}All the $L_{i}$ have the same sign \textup{(}in
particular none are zero\textup{).}

\item \label{ProSub.1(3)}There is a $\beta\in\mathbb{R}$ such that
$\mathcal{O}_{d}$ admits a $\left(  \sigma,\beta\right)  $-KMS state.
\end{enumerate}

\noindent Furthermore, if these conditions are not fulfilled, $\mathfrak
{A}_{L}$ contains an isometry which is not unitary.
\end{proposition}

\begin{proof}
(\ref{ProSub.1(2)})$\Leftrightarrow$(\ref{ProSub.1(3)}) was established in
Proposition \ref{ProKMS.1}.

(\ref{ProSub.1(1)}) $\Rightarrow$ (\ref{ProSub.1(2)}): Assume that
(\ref{ProSub.1(2)}) does not hold. Then there exist $i,j\in\left\{
1,\dots,d\right\}  $ with $L_{i}>0$, $L_{j}<0$. Put%
\begin{equation}
s=s_{i}^{-L_{j}}s_{j}^{L_{i}}. \label{eqSub.7}%
\end{equation}
Then $s$ is an isometry in $\mathfrak{A}_{L}$ which is not unitary. Hence
$\mathfrak{A}_{L}$ cannot be an AF-algebra. This also establishes the final
remark in Proposition \ref{ProSub.1}.

(\ref{ProSub.1(2)}) $\Rightarrow$ (\ref{ProSub.1(1)}): We may assume that all
$L_{i}$ are positive. We have noted that $\mathcal{O}_{d}^{\sigma}%
=\mathcal{O}_{d}^{\sigma}\left(  0\right)  $ is the closure of the linear span
of $s_{\alpha}^{{}}s_{\gamma}^{\ast}$ with $L\left(  \alpha\right)  =L\left(
\gamma\right)  $. If $L\left(  \alpha\right)  =L\left(  \gamma\right)  $, we
define\label{LOSgrade}%
\begin{equation}
\operatorname*{grade}\left(  s_{\alpha}^{{}}s_{\gamma}^{\ast}\right)
=L\left(  \alpha\right)  , \label{eqSub.8}%
\end{equation}
and we set $\operatorname*{grade}\left(  \openone\right)
=\operatorname*{grade}\left(  0\right)  =0$. Now, if $s_{\alpha}^{{}}%
s_{\gamma}^{\ast}$, $s_{\delta}^{{}}s_{\varepsilon}^{\ast}$ are in
$\mathcal{O}_{d}^{\sigma}$ then either the product $s_{\alpha}^{{}}s_{\gamma
}^{\ast}s_{\delta}^{{}}s_{\varepsilon}^{\ast}$ is zero, or we have
$\gamma=\delta\gamma^{\prime}$ and the product is $s_{\alpha}^{{}%
}s_{\varepsilon\gamma^{\prime}}^{\ast}$, or we have $\delta=\gamma
\delta^{\prime}$ and the product is $s_{\alpha\delta^{\prime}}^{{}%
}s_{\varepsilon}^{\ast}$. In the latter two cases $\operatorname*{grade}%
\left(  s_{\alpha}^{{}}s_{\gamma}^{\ast}s_{\delta}^{{}}s_{\varepsilon}^{\ast
}\right)  =\max\left(  \operatorname*{grade}\left(  s_{\alpha}^{{}}s_{\gamma
}^{\ast}\right)  ,\operatorname*{grade}\left(  s_{\gamma}^{{}}s_{\delta}%
^{\ast}\right)  \right)  $, and, in the former case, $\operatorname*{grade}%
\left(  s_{\alpha}^{{}}s_{\gamma}^{\ast}s_{\delta}^{{}}s_{\varepsilon}^{\ast
}\right)  =0$. Thus in general,%
\begin{equation}
\operatorname*{grade}\left(  s_{\alpha}^{{}}s_{\gamma}^{\ast}s_{\varepsilon
}^{{}}s_{\delta}^{\ast}\right)  \leq\max\left(  \operatorname*{grade}\left(
s_{\alpha}^{{}}s_{\gamma}^{\ast}\right)  ,\operatorname*{grade}\left(
s_{\varepsilon}^{{}}s_{\delta}^{\ast}\right)  \right)  . \label{eqSub.9}%
\end{equation}
Thus if we define%
\begin{equation}
\mathfrak{A}_{n}=%
\operatorname*{lin\,span}%
\left\{  s_{\alpha}^{{}}s_{\gamma}^{\ast}\mid L\left(  \alpha\right)
=L\left(  \gamma\right)  \leq n\right\}  , \label{eqSub.10}%
\end{equation}
then $\mathfrak{A}_{n}$ is a $\ast$-algebra, and $\mathfrak{A}_{n}$ is
finite-dimensional since $L_{i}>0$ for $i=1,\dots,d$. Since any $s_{\alpha
}^{{}}s_{\gamma}^{\ast}\in\mathcal{O}_{d}^{\sigma}$ has a grade, it follows
that $\bigcup_{n}\mathfrak{A}_{n}$ is dense in $\mathcal{O}_{d}^{\sigma
}=\mathfrak{A}_{L}$. Thus $\mathfrak{A}_{L}$ is an AF-algebra, and Proposition
\ref{ProSub.1} is proved.
\end{proof}

We refer to \cite{Bra72} and Remark \ref{RemRemJun.6} for AF-algebras and
Bratteli diagrams,\label{LOSBrattelidiagrams_2} to \cite{Bla86} for
$K$-theory, and to
\cite{RLL00},
\cite{Dav96} and \cite{GoHa87} for good recent treatments
of both. In order to analyze the AF-algebra $\mathfrak
{A}_{L}$ further, it turns out to be convenient to define subalgebras
$\mathfrak{A}_{n}$ in a more sophisticated way than above, and this is the
object of the following. Note that, except for simple cases (like $d=2$), the
finite-dimensional subalgebras introduced below are larger than $\mathfrak
{A}_{n}$. The main structure theorem on $\mathfrak{A}_{L}$ is the following.

\begin{theorem}
\label{ThmSub.2}Let $L_{1}\leq L_{2}\leq\dots\leq L_{d}$ be positive integers
such that the greatest common divisor of $L_{1},\dots,L_{d}$ is $1$. It
follows that $\mathfrak{A}_{L}$ is a simple AF-algebra with a unique trace
state defined as follows: Let $\beta$ be the positive real number such that%
\begin{equation}
\sum_{i=1}^{d}e^{-\beta L_{i}}=1, \label{eqSub.11}%
\end{equation}
and put\label{LOSpj_3}%
\begin{equation}
p_{i}=e^{-\beta L_{i}}. \label{eqSub.12}%
\end{equation}
Then the unique trace state is the restriction to $\mathfrak{A}_{L}$ of the
state $\omega$ defined on $\mathcal{O}_{d}$ by%
\begin{equation}
\omega\left(  s_{\alpha}^{{}}s_{\gamma}^{\ast}\right)  =\delta_{\alpha\gamma
}p^{\alpha} \label{eqSub.13}%
\end{equation}
where\label{LOSpalpha}%
\begin{equation}
p^{\alpha}=p_{\alpha_{1}}p_{\alpha_{2}}\cdots p_{\alpha_{k}} \label{eqSub.14}%
\end{equation}
for $\alpha=\left(  \alpha_{1}\dots\alpha_{k}\right)  $.
\end{theorem}

\begin{remark}
\label{RemSub.3}During the proof of Theorem \textup{\ref{ThmSub.2},} we will
show that the AF-algebra $\mathfrak{A}_{L}$ has a Bratteli diagram which
stabilizes after a finite number of steps to having constant incidence
matrices. This diagram may be described explicitly as follows: The nodes are
indexed by $\left(  n,m\right)  $, where $n$ indexes the rows, $n=0,1,2,\dots
$, and $m$ indexes the nodes in the row, $m=0,1,\dots,L_{d}-1$. Some of the
nodes in the first rows may correspond to the algebra $0$: for example,
$\left(  0,m\right)  $, $m=0,1,\dots,L_{d}-1$, correspond to the algebras
$M_{1}=\mathbb{C},0,0,0,\dots,0$. The embedding from one row to the next is
given as follows: There are lines from $\left(  n-1,0\right)  $ to $\left(
n,m\right)  $ if and only if $m=L_{k}-1$ for some $k$, and the number of lines
between these nodes is equal to the number of such $k$'s. There is a single
line from $\left(  n-1,m\right)  $ to $\left(  n,m-1\right)  $ for
$m=1,\dots,L_{d}-1$. Finally, to obtain the actual Bratteli diagram, one
should throw away all nodes corresponding to the algebra $0$ as well as the
edges emanating from such nodes. The assumption that the greatest common
divisor of $L_{1},\dots,L_{d}$ is $1$ will imply that there are just finitely
many such nodes. It will be clear from the proof how this pattern appears. We
give some examples in the figures below.

\begin{figure}[ptb]
\begin{picture}(19,226)(0,-12)
\put(0,0){\includegraphics
[bb=110 32 129 240,clip,height=208bp,width=19bp]{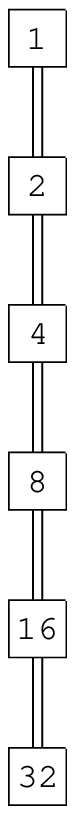}}
\put(10.5,208){\makebox(0,0)[b]{$m=0$}}
\put(10.5,0){\makebox(0,0)[t]{$\vdots$}}
\put(-6,190){\makebox(0,0)[r]{$n=0$}}
\put(-6,147.5){\makebox(0,0)[r]{$n=1$}}
\put(-6,105){\makebox(0,0)[r]{$n=2$}}
\put(-6,62.5){\makebox(0,0)[r]{$\vdots\mkern18mu$}}
\end{picture}
\caption{$d=2$; $L_{1}=1$, $L_{2}=1$; $\beta=\ln2$.}%
\label{BratDiag11}%
\end{figure}

\begin{figure}[ptb]
\begin{picture}(288,532)(0,-12)
\put(0,40){\includegraphics
[bb=232 48 296 528,clip,height=480bp,width=64bp]{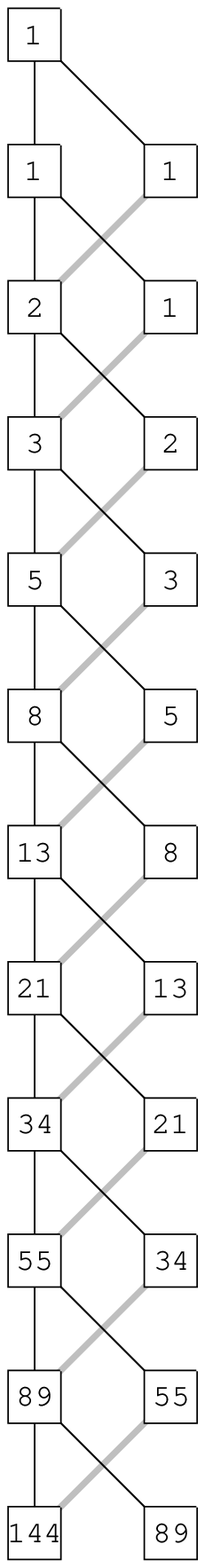}}
\put(10,514){\makebox(0,0)[b]{$m=0$}}
\put(55,514){\makebox(0,0)[b]{$m=1$}}
\put(10,40){\makebox(0,0)[t]{$\vdots$}}
\put(55,40){\makebox(0,0)[t]{$\vdots$}}
\put(-6,495){\makebox(0,0)[r]{$n=0$}}
\put(-6,450.5){\makebox(0,0)[r]{$n=1$}}
\put(-6,406){\makebox(0,0)[r]{$n=2$}}
\put(-6,361.5){\makebox(0,0)[r]{$\vdots\mkern18mu$}}
\put(124,0){\includegraphics
[bb=238 120 402 640,clip,height=520bp,width=164bp]{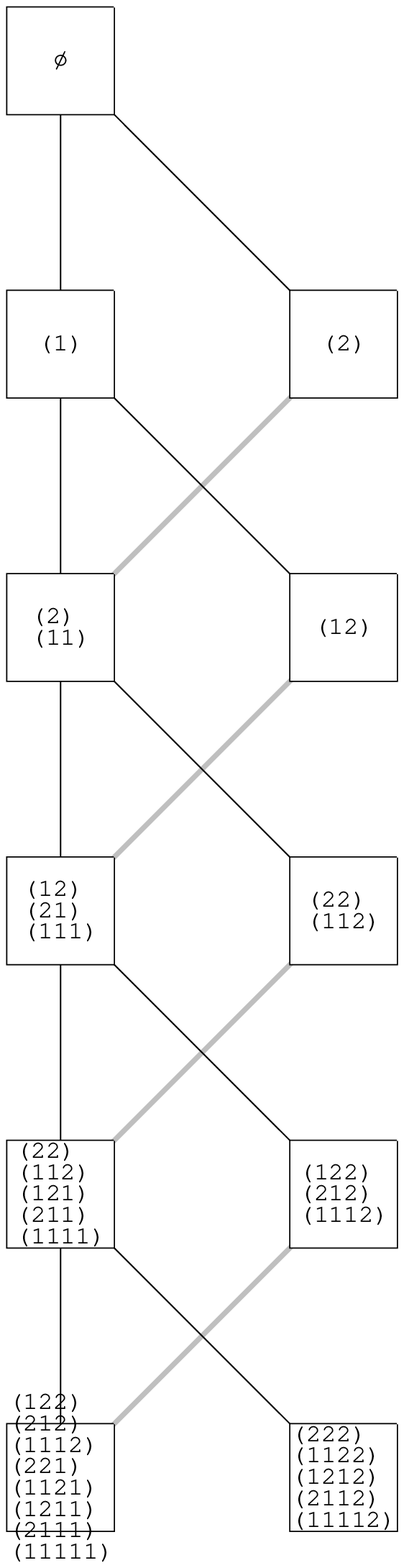}}
\put(149.5,514){\makebox(0,0)[b]{$m=0$}}
\put(262.5,514){\makebox(0,0)[b]{$m=1$}}
\put(149.5,0){\makebox(0,0)[t]{$\vdots$}}
\put(262.5,0){\makebox(0,0)[t]{$\vdots$}}
\put(294,464){\makebox(0,0)[l]{$n=0$}}
\put(294,351.5){\makebox(0,0)[l]{$n=1$}}
\put(294,239){\makebox(0,0)[l]{$n=2$}}
\put(294,126.5){\makebox(0,0)[l]{$\mkern18mu\vdots$}}
\end{picture}
\caption{$d=2$; $L_{1}=1$, $L_{2}=2$; $\beta=-\ln\left(  \left(  \sqrt
{5}-1\right)  /2\right)  $. Then the Bratteli diagram is given by the
Fibonacci sequence. Detail on the right shows the multi-indices for each node
in the top five rows.}%
\label{BratDiag12}%
\end{figure}

\begin{figure}[ptb]
\begin{picture}(220,528)(0,-12)
\put(0,0){\includegraphics
[bb=223 75 498 720,clip,height=516bp,width=220bp]{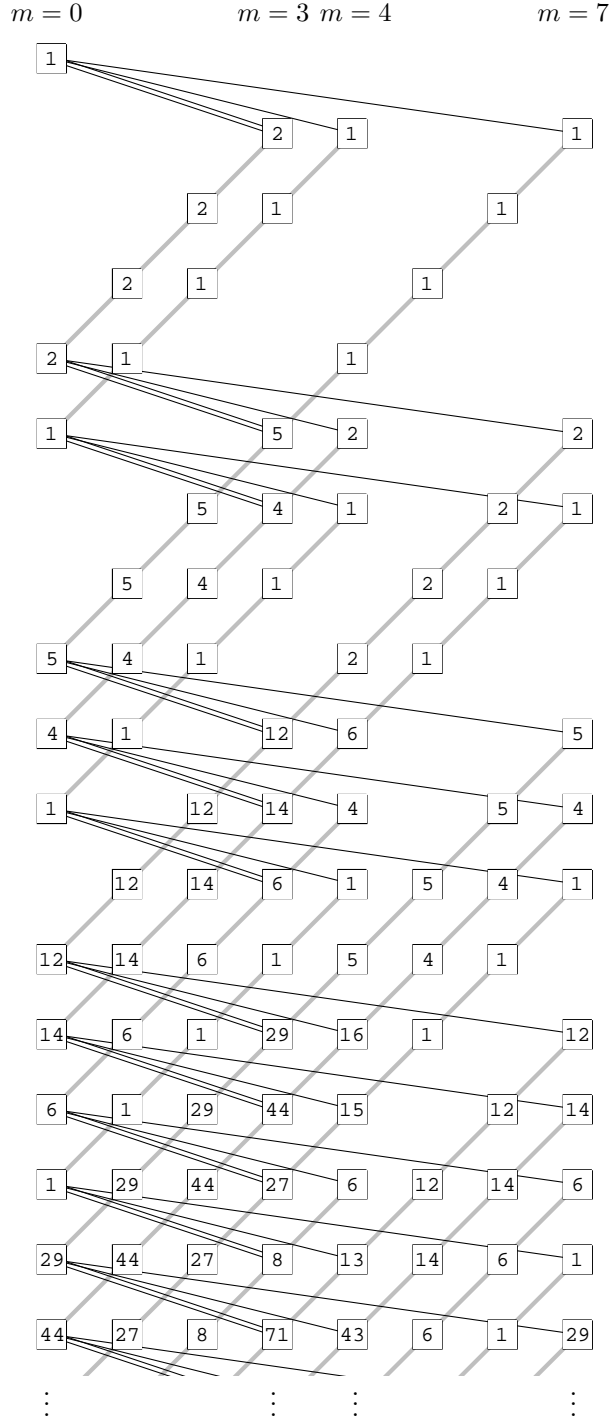}}
\put(8.6,511){\makebox(0,0)[b]{$m=0$}}
\put(8.6,0){\makebox(0,0)[t]{$\vdots$}}
\put(94,511){\makebox(0,0)[b]{$m=3$}}
\put(94,0){\makebox(0,0)[t]{$\vdots$}}
\put(125,511){\makebox(0,0)[b]{$m=4$}}
\put(125,0){\makebox(0,0)[t]{$\vdots$}}
\put(207,511){\makebox(0,0)[b]{$m=7$}}
\put(207,0){\makebox(0,0)[t]{$\vdots$}}
\end{picture}
\caption{$d=4$; $L=\{4,4,5,8\}$; first matrix column
$=(0\;0\;0\;2\;1\;0\;0\;1)^{\mathrm{t}}$; $\beta=-\ln x$ where $x=(-2+\sqrt[3]%
{100+12\sqrt{69}}+\sqrt[3]{100-12\sqrt{69}})/6\approx0.7549$ solves
$2x^{4}+x^{5}+x^{8}=1$. (Actually $x$ solves $x^{2}+x^{3}=1$.) See the $n=5$
case in Example \ref{ExaRemNewBis.3}.}%
\label{BratDiag4458}%
\end{figure}

\begin{figure}[ptb]
\begin{picture}(288,566)(0,-12)
\put(0,0){\includegraphics
[bb=249 36 341 590,clip,height=554bp,width=92bp]{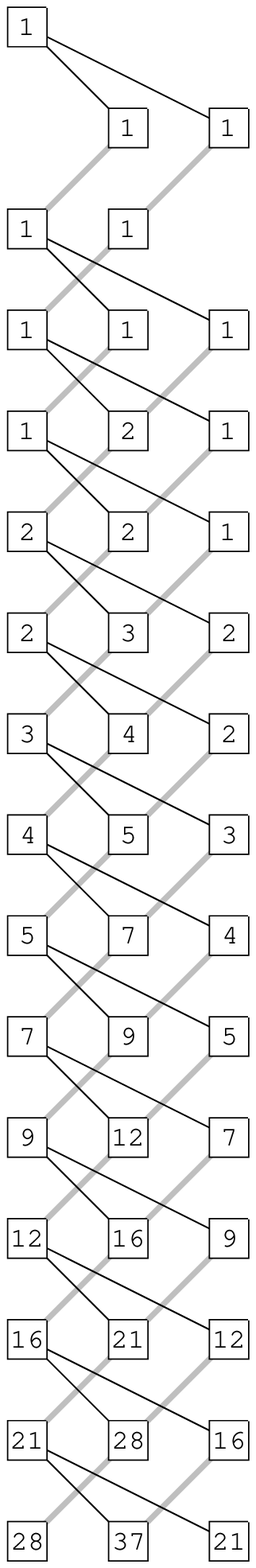}}
\put(10,548){\makebox(0,0)[b]{$m=0$}}
\put(46,548){\makebox(0,0)[b]{$m=1$}}
\put(82,548){\makebox(0,0)[b]{$m=2$}}
\put(10,0){\makebox(0,0)[t]{$\vdots$}}
\put(46,0){\makebox(0,0)[t]{$\vdots$}}
\put(82,0){\makebox(0,0)[t]{$\vdots$}}
\put(-6,530.5){\makebox(0,0)[r]{$n=0$}}
\put(-6,493.5){\makebox(0,0)[r]{$n=1$}}
\put(-6,456.5){\makebox(0,0)[r]{$n=2$}}
\put(-6,419.5){\makebox(0,0)[r]{$\vdots\mkern18mu$}}
\put(128,249){\includegraphics
[bb=100 60 260 360,clip,height=300bp,width=160bp]{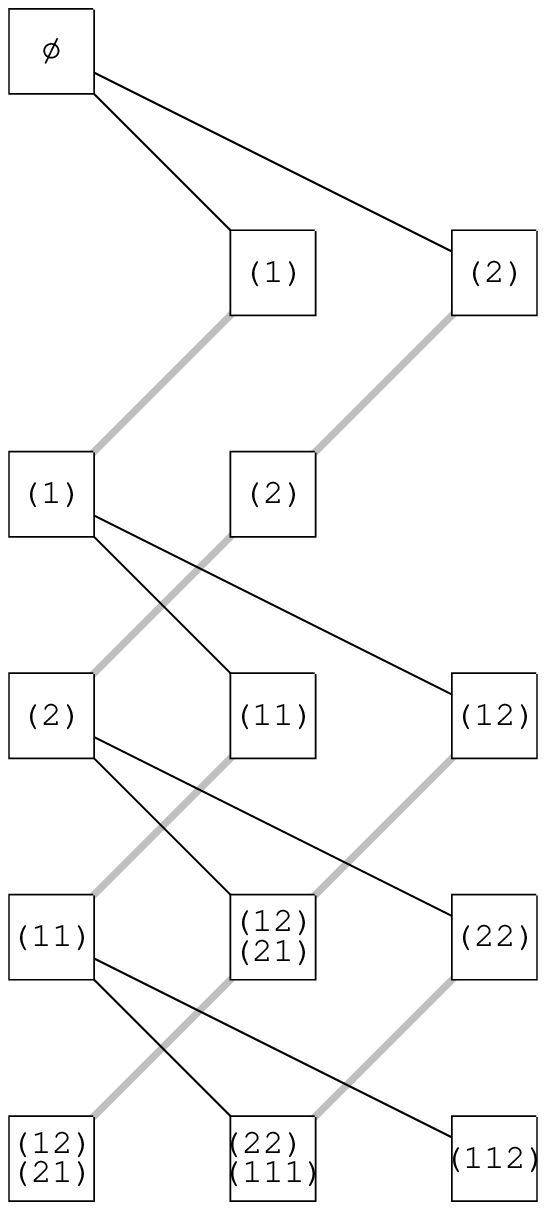}}
\put(144,548){\makebox(0,0)[b]{$m=0$}}
\put(208,548){\makebox(0,0)[b]{$m=1$}}
\put(272,548){\makebox(0,0)[b]{$m=2$}}
\put(144,249){\makebox(0,0)[t]{$\vdots$}}
\put(208,249){\makebox(0,0)[t]{$\vdots$}}
\put(272,249){\makebox(0,0)[t]{$\vdots$}}
\put(294,520.33){\makebox(0,0)[l]{$n=0$}}
\put(294,456.33){\makebox(0,0)[l]{$n=1$}}
\put(294,392.33){\makebox(0,0)[l]{$n=2$}}
\put(294,328.33){\makebox(0,0)[l]{$\mkern18mu\vdots$}}
\end{picture}
\caption{$d=2$; $L=\left\{  2,3\right\}  $; first matrix column
$=(0\;1\;1)^{\mathrm{t}}$; $\beta=-\ln x$ where $x>0$ solves $x^{2}+x^{3}=1$.
Detail on the right shows the multi-indices for each node in the top five
rows. See the proof of Lemma \ref{LemSub.6}.}%
\label{BratDiag23}%
\end{figure}

\begin{figure}[ptb]
\begin{picture}(168,566)(0,-12)
\put(0,0){\includegraphics
[bb=211 36 379 590,clip,height=554bp,width=168bp]{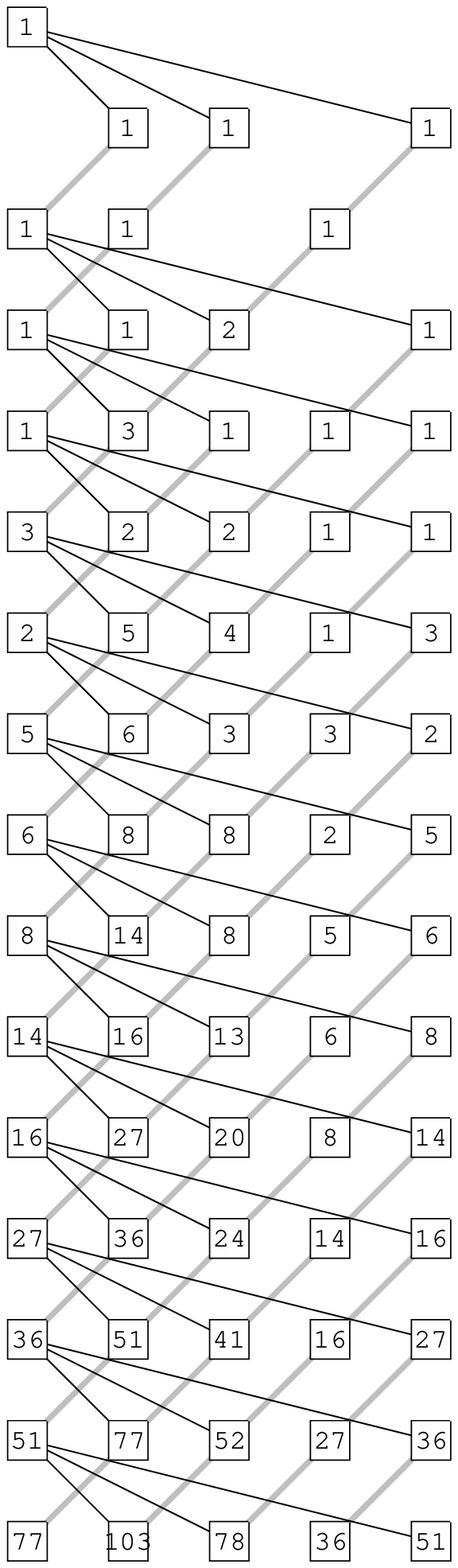}}
\put(10.5,548){\makebox(0,0)[b]{$m=0$}}
\put(47,548){\makebox(0,0)[b]{$m=1$}}
\put(84,548){\makebox(0,0)[b]{$m=2$}}
\put(157,548){\makebox(0,0)[b]{$m=4$}}
\put(10.5,0){\makebox(0,0)[t]{$\vdots$}}
\put(47,0){\makebox(0,0)[t]{$\vdots$}}
\put(84,0){\makebox(0,0)[t]{$\vdots$}}
\put(157,0){\makebox(0,0)[t]{$\vdots$}}
\put(-6,530.5){\makebox(0,0)[r]{$n=0$}}
\put(-6,493.5){\makebox(0,0)[r]{$n=1$}}
\put(-6,456.5){\makebox(0,0)[r]{$n=2$}}
\put(-6,419.5){\makebox(0,0)[r]{$\vdots\mkern18mu$}}
\end{picture}
\caption{$d=3$; $L=\left\{  2,3,5\right\}  $; first matrix column
$=(0\;1\;1\;0\;1)^{\mathrm{t}}$; $\beta=-\ln x$ where $x>0$ solves
$x^{2}+x^{3}+x^{5}=1$.}%
\label{BratDiag235}%
\end{figure}

We will show that the unique positive eigenvalue of the constant incidence
matrix (the Frobenius eigenvalue) is $e^{\beta}$.
\end{remark}

Before proving Theorem \ref{ThmSub.2} and Remark \ref{RemSub.3}, we look at
some examples.

Figure \ref{BratDiag11} is the CAR-algebra of type $2^{\infty}$; see
\cite{Gli60}, \cite{Bra72}, \cite{Eff81}, \cite{BrRoII} and \cite{Ren80}.

Figure \ref{BratDiag12} is the AF-algebra with same dimension group as the
rotation algebra $\mathfrak{A}_{\theta}$ for $\theta=\frac{\sqrt{5}-1}{2}=$
the golden ratio. Pimsner and Voiculescu \cite{PiVo80} embedded $\mathfrak
{A}_{\theta}$ into this AF-algebra.

Figure \ref{BratDiag4458} illustrates that the Bratteli diagram is more
``slow'' in stabilizing when the $L_{i}$-numbers are dispersed. Figures
\ref{BratDiag23} and \ref{BratDiag235} illustrate how the multi-indices build
up as the sizes of the matrix algebras increase going down the Bratteli
diagram. The significance of the choices of $L_{i}$-numbers will become more
clear in Chapter \ref{Rem} below where we study isomorphism invariants for the
AF-algebras $\mathfrak{A}_{L}$ in general.

Figure \ref{BratDiag23} represents $\left(
\begin{smallmatrix}
0 & 1 & 1
\end{smallmatrix}
\right)  $, the first of two AF-algebras which share Perron--Frobenius
eigenvalue $\lambda=e^{\beta}$\label{LOSPerronFrobeniuseigenvalue_2} where
$a=e^{-\beta}\approx0.7549$ is the real root of $x^{2}+x^{3}=1$. The other
one, $\left(
\begin{smallmatrix}
1 & 0 & 0 & 0 & 1
\end{smallmatrix}
\right)  $, is obtained from $x+x^{5}=1$, which has the same positive root
$a$. (See Remark \ref{RemSub.9} and Chapter \ref{Rem} for more details on the
Perron--Frobenius eigenvalue.) Yet these two AF-algebras are non-isomorphic,
since their dimension groups have rank $3$ and $5$, respectively. (See Theorem
\ref{CorCyc.9}.) They correspond to the pair of lattice points $\left(
2,3\right)  $, $\left(  1,5\right)  $ that is illustrated in Figure
\ref{Lfromx}.

\begin{figure}[ptb]
\begin{picture}(324,529)(-12,10)
\put(0,0){\includegraphics
[bb=147 0 579 732,clip,height=549bp,width=324bp]{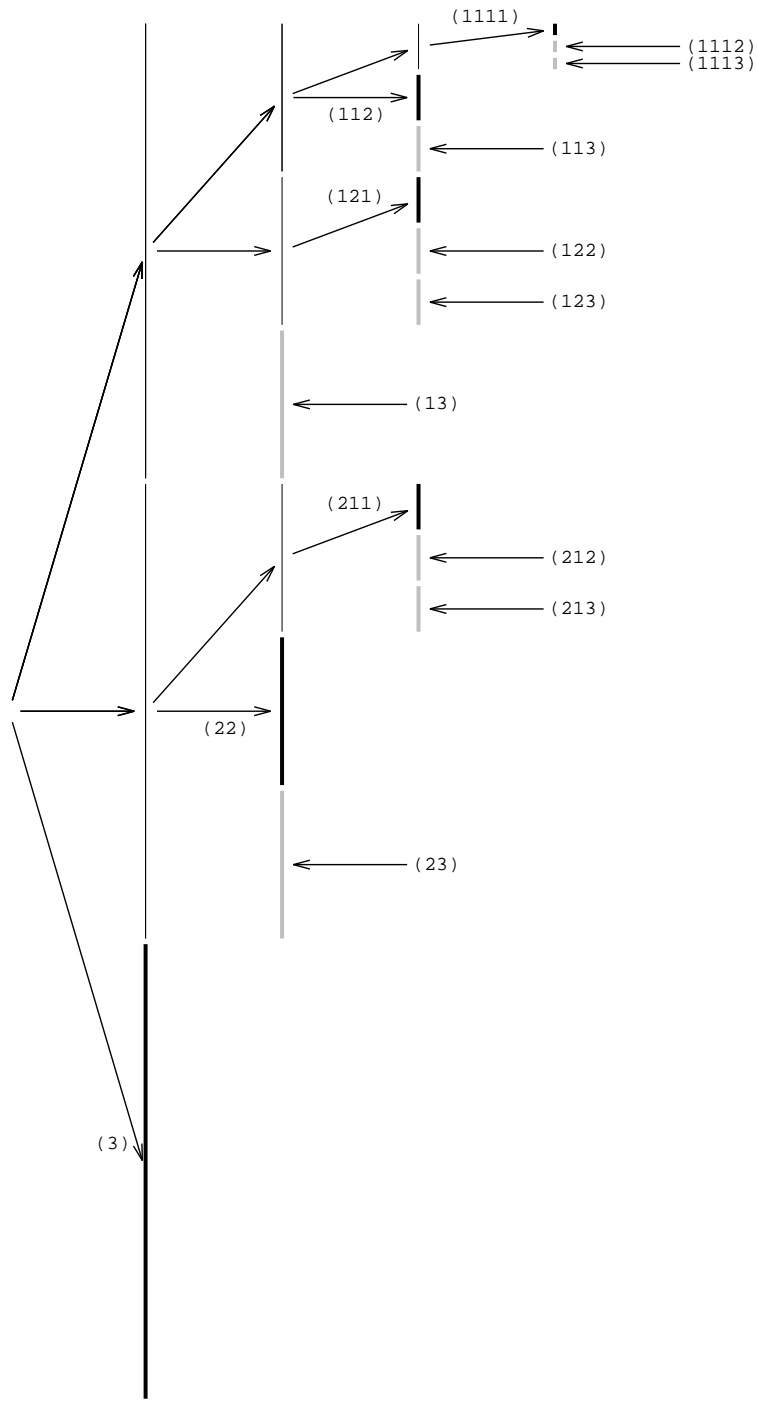}}
\end{picture}
\caption{Illustration of procedure in proof of Lemma \ref{LemSub.6}, with
$d=3$, $L=\{1,2,4\}$. Compare with Figure \ref{BratDiag124} and Example
\ref{ExaSub.7}.}%
\label{DustDiag124}%
\end{figure}

Figure \ref{DustDiag124} illustrates the procedure in the proof of Lemma
\ref{LemSub.6}, below.

\begin{figure}[ptb]
\begin{picture}(360,537)(0,0)
\put(0,0){\includegraphics
[bb=215 46 401 616,clip,height=513bp,width=167.4bp]{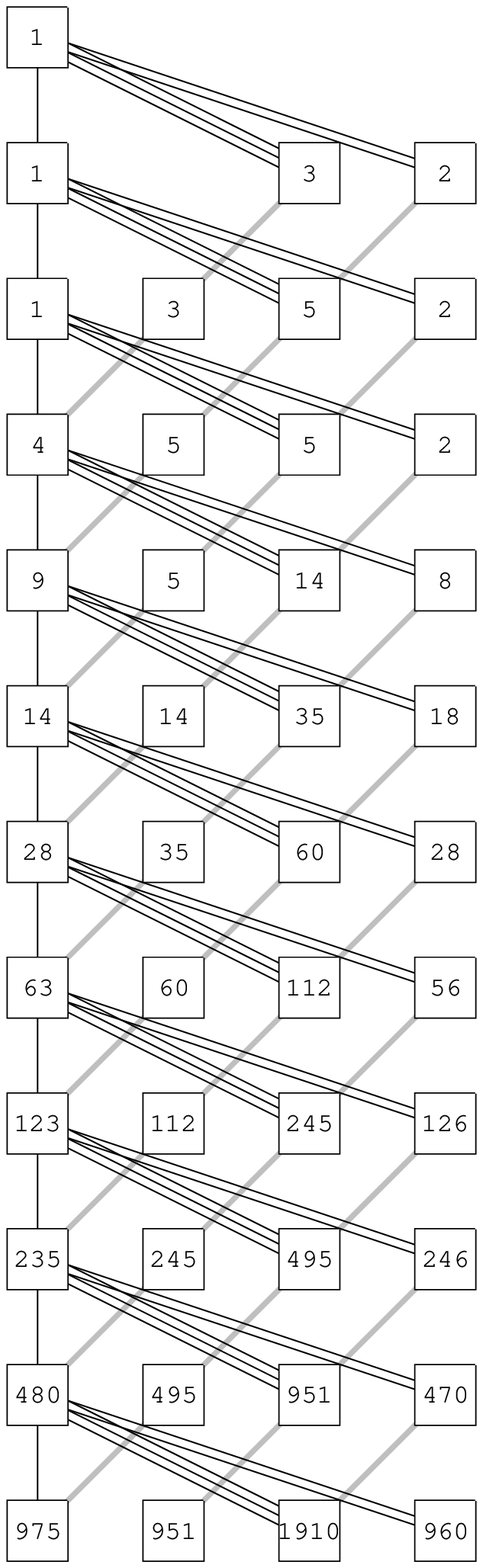}}
\put(192.6,0){\includegraphics
[bb=215 46 401 616,clip,height=513bp,width=167.4bp]{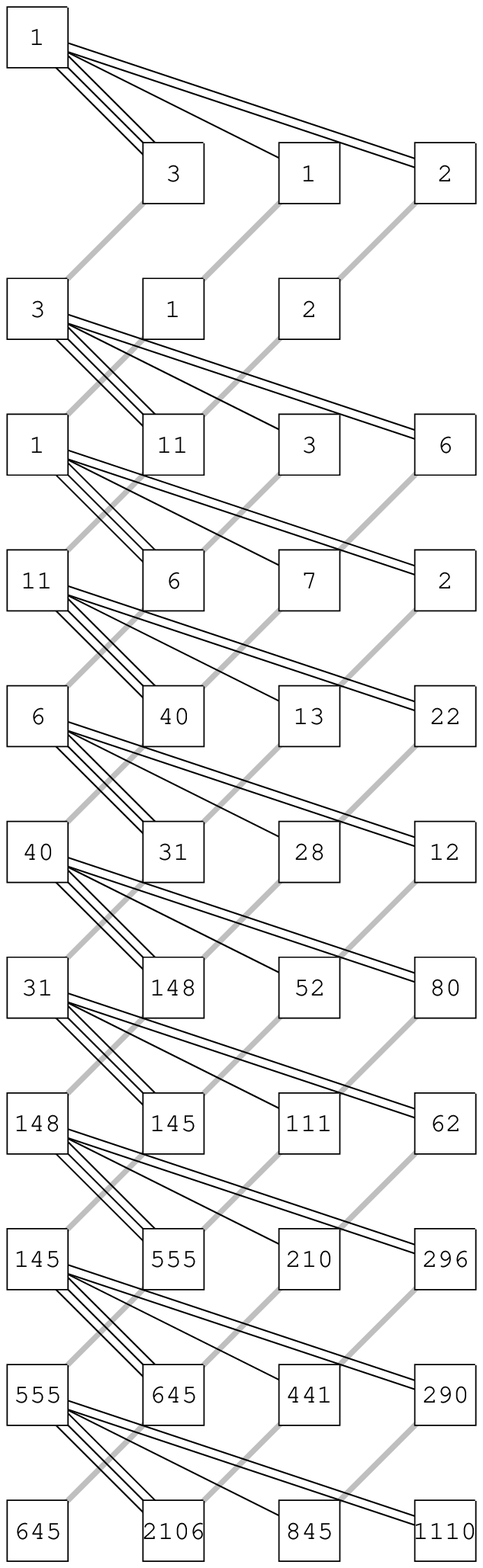}}
\end{picture}
\caption{$d=6$; $L=\{1,3,3,3,4,4\}$ (left), $L=\left\{  2,2,2,3,4,4\right\}  $
(right). These define non-isomorphic algebras (see Chapter \ref{APP.EXA}).}%
\label{BratDiagsd6}%
\end{figure}

Let $\tau$ be the additive real character defined on the dimension group
$K_{0}\left(  \mathfrak{A}_{L}\right)  $\label{K0frakAL} by the trace state,
\cite{Eff81}. Figure \ref{BratDiagsd6} represents two examples with the same
$\ker\left(  \tau\right)  $ ($\cong\mathbb{Z}^{3}$), the same $\tau\left(
K_{0}\right)  $ ($\cong\mathbb{Z}\left[  \frac{1}{2}\right]  $) but still
non-isomorphic AF-algebras, as they represent different elements of
$\operatorname*{Ext}\left(  \mathbb{Z}\left[  \frac{1}{2}\right]
,\mathbb{Z}^{3}\right)  $. (Details in Chapter \ref{Ext}.)

We will prove Theorem \ref{ThmSub.2} and Remark \ref{RemSub.3} via a series of
lemmas. First a definition:

\begin{definition}
\label{DefSub.4}A set $\left\{  e_{\alpha\gamma}\right\}  _{\alpha,\gamma\in
I}$ of elements of a $C^{\ast}$-algebra $\mathfrak{A}$, doubly indexed by a
finite set $I$, is said to be a \emph{system of matrix units}%
\label{LOSmatrixunits} if

\begin{enumerate}
\item \label{DefSub.4(1)}$e_{\alpha\gamma}e_{\xi\eta}=\delta_{\gamma\xi
}e_{\alpha\eta}$,

\item \label{DefSub.4(2)}$e_{\gamma\alpha}^{\ast}=e_{\alpha\gamma}^{{}}$\textup{.}
\end{enumerate}
\end{definition}

In that case, matrices $\left(  A_{\alpha\gamma}\right)  _{\alpha,\gamma\in
I}$ over $\mathbb{C}$ may be represented in $\mathfrak{A}$ as follows:
$\left(  A_{\alpha\gamma}\right)  \mapsto\sum_{\alpha} \sum_{\gamma}%
A_{\alpha\gamma}e_{\alpha\gamma}$. Note that we do \emph{not} assume that the
projection $\sum_{\alpha} e_{\alpha\alpha}$ is the identity of $\mathfrak{A}$.

\begin{lemma}
\label{LemSub.5}Let $L_{1},L_{2},\dots,L_{d}$ be positive integers and let
$\sigma=\sigma^{L}$ be the associated automorphism group
\textup{(\ref{eqSub.5}).} Let $\mathcal{B}_{I}=\left\{  s_{\alpha}^{{}%
}s_{\gamma}^{\ast}\right\}  _{\alpha,\gamma\in I}$ be a finite set of elements
of $\mathfrak{A}_{L}=\mathcal{O}_{d}^{\sigma}$. The doubly indexed set
$\mathcal{B}_{I}$ is then a set of matrix units if and only if there is an
$n\in\mathbb{N}$ such that $L\left(  \alpha\right)  =n$ for all $\alpha\in I$.
\end{lemma}

\begin{proof}
Consider arbitrary multi-indices $\alpha$, $\gamma$, $\xi$, and $\eta$ built
from $\mathbb{Z}_{d}$. The product
\begin{equation}
\left(  s_{\alpha}^{{}}s_{\gamma}^{\ast}\right)  \left(  s_{\xi}^{{}}s_{\eta
}^{\ast}\right)  \label{eqSub.9bis}%
\end{equation}
is nonzero only if either $\gamma$ is of the form $\gamma=\left(  \xi
\gamma^{\prime}\right)  $, or $\xi$ is of the form $\xi=\left(  \gamma
\xi^{\prime}\right)  $. If each of the factors in (\ref{eqSub.9bis}) is in
$\mathcal{O}_{d}^{\sigma}$, then $L\left(  \alpha\right)  =L\left(
\gamma\right)  $ and $L\left(  \xi\right)  =L\left(  \eta\right)  $. Recall
that the grade of the first factor is $L\left(  \alpha\right)  $, and that of
the second is $L\left(  \xi\right)  $. If the two factors have different
grades, and if the product is nonzero, then $\gamma^{\prime}\neq{}$\textup{\o}
or $\xi^{\prime}\neq{}$\textup{\o}. In the first case, the product is
$s_{\alpha}^{{}}s_{\left(  \eta\gamma^{\prime}\right)  }^{\ast}$, and in the
second it is $s_{\left(  \alpha\xi^{\prime}\right)  }^{{}}s_{\eta}^{\ast}$. In
either case, if $\gamma^{\prime}\neq{}$\textup{\o} or $\xi^{\prime}\neq{}%
$\textup{\o}, the product of the two elements from $\mathcal{B}_{I}$ will be
nonzero with $\gamma\neq\xi$, see (\ref{eqSub.9bis}), so $\mathcal{B}$ will
not then be a set of matrix units, i.e., condition (\ref{DefSub.4(1)}) of
Definition \ref{DefSub.4} will not hold. This proves the ``only if'' part of
Lemma \ref{LemSub.5}.

Conversely, if there exists an $n$ such that $L\left(  \alpha\right)  =n$ for
all $\alpha\in I$, then the case $\gamma=\left(  \xi\gamma^{\prime}\right)  $
with $\gamma^{\prime}\neq{}$\textup{\o} is excluded since $L\left(  \xi
\gamma^{\prime}\right)  =L\left(  \xi\right)  +L\left(  \gamma^{\prime
}\right)  $. For if $L\left(  \gamma\right)  =L\left(  \xi\right)  =n$, then
$L\left(  \gamma^{\prime}\right)  =0$, and $\gamma^{\prime}={}$\textup{\o}.
The same argument also excludes $\xi=\left(  \gamma\xi^{\prime}\right)  $ with
$\xi^{\prime}\neq{}$\textup{\o}. It follows that condition (\ref{DefSub.4(1)})
of Definition \ref{DefSub.4} will be satisfied for $e_{\alpha\gamma}^{{}%
}=s_{\alpha}^{{}}s_{\gamma}^{\ast}$ with $I$ as an index set.
\end{proof}

\begin{lemma}
\label{LemSub.6}Let $d\in\mathbb{N}$ and let $L_{1},\dots,L_{d}$ be positive
integers. Define $L\left(  \alpha\right)  =\sum_{j}\#_{j}\left(
\alpha\right)  L_{j}$\label{LOSLalpha_2} on all finite multi-indices $\alpha$
from $\mathbb{Z}_{d}$ as in \textup{(\ref{eqSub.2}).}
Define\label{LOSLinversen}%
\begin{equation}
L^{-1}\left(  n\right)  =\left\{  \alpha\mid L\left(  \alpha\right)
=n\right\}  \label{eqSub.15}%
\end{equation}
and put\label{LOSEn}%
\begin{equation}
E_{n}=\left\{  \gamma\mid\gamma=\left(  \alpha i\right)  \text{ where
}L\left(  \alpha\right)  <n\text{ and }L\left(  \alpha\right)  +L_{i}%
>n\right\}  . \label{eqSub.16}%
\end{equation}
Then%
\begin{equation}
\sum_{\alpha\in L^{-1}\left(  n\right)  }s_{\alpha}^{{}}s_{\alpha}^{\ast}%
+\sum_{\gamma\in E_{n}}s_{\gamma}^{{}}s_{\gamma}^{\ast}=1, \label{eqSub.17}%
\end{equation}
i.e., the projections in the family $\left\{  s_{\alpha}^{{}}s_{\alpha}^{\ast
}\mid\alpha\in L^{-1}\left(  n\right)  \right\}  \cup\left\{  s_{\gamma}^{{}%
}s_{\gamma}^{\ast}\mid\gamma\in E_{n}\right\}  $ are mutually orthogonal with
sum $1$.
\end{lemma}

\setcounter{case}{0}

\begin{proof}
Let us use the shorthand notation%
\[
\left(  \alpha\right)  =e_{\alpha\alpha}^{{}}\qquad\left(  =s_{\alpha}^{{}%
}s_{\alpha}^{\ast}\right)  .
\]
It follows from the computations in the proof of Lemma \ref{LemSub.5} that,
given two projections $\left(  \alpha\right)  $, $\left(  \gamma\right)  $,
then $\left(  \alpha\right)  $, $\left(  \gamma\right)  $ are either mutually
orthogonal, or one is contained in the other; and the latter case occurs in,
and only in, the following two cases:

\begin{case}
\label{LemSub.6Case(1)}$\alpha=\gamma\alpha^{\prime}$. Then $\left(
\alpha\right)  \leq\left(  \gamma\right)  $. Or,
\end{case}

\begin{case}
\label{LemSub.6Case(2)}$\gamma=\alpha\gamma^{\prime}$. Then $\left(
\gamma\right)  \leq\left(  \alpha\right)  $ (with strict inequalities if and
only if $\alpha^{\prime}\neq{}$\textup{\o}, $\gamma^{\prime}\neq{}$\textup{\o
}, respectively).
\end{case}

Using this, it follows easily from case-by-case considerations that the
projections in the family\label{LOSAn}%
\[
A_{n}:=\left\{  \left(  \alpha\right)  \mid\alpha\in L^{-1}\left(  n\right)
\cup E_{n}\right\}
\]
are mutually orthogonal. For example, the projections $\left(  \alpha\right)
$, $\alpha\in L^{-1}\left(  n\right)  $ are mutually orthogonal by Lemma
\ref{LemSub.5}, and if $\alpha\in L^{-1}\left(  n\right)  $ and $\gamma
=\left(  \delta i\right)  \in E_{n}$ with $L\left(  \delta\right)  <n$,
$L\left(  \delta\right)  +i>n$, then both Case \ref{LemSub.6Case(1)} and
\ref{LemSub.6Case(2)} are excluded, so $\left(  \alpha\right)  \left(
\gamma\right)  =0$; and similarly, if $\alpha=\left(  \varepsilon j\right)  $
and $\gamma=\left(  \delta i\right)  $ are in $E_{n}$, then $\left(
\alpha\right)  \left(  \gamma\right)  \neq0$ implies $\alpha=\gamma$. It
remains to show that the projections in these two families add up to $1$. If
not, there would exist a multi-index $\left(  \delta\right)  $ such that
$\left(  \delta\right)  $ is orthogonal to all projections in the two
families. If then $L\left(  \delta\right)  <n$, we could find a $\delta
^{\prime}$ such that $\delta\delta^{\prime}\in L^{-1}\left(  n\right)  $ or
$\delta\delta^{\prime}\in E_{n}$, but since $\left(  \delta\right)  \left(
\delta\delta^{\prime}\right)  =\left(  \delta\delta^{\prime}\right)  \neq0$,
this is impossible. If $L\left(  \delta\right)  =n$, then $\delta\in
L^{-1}\left(  n\right)  $, which is impossible. If $L\left(  \delta\right)
>n$, write $\delta=\left(  \delta_{1}\delta_{2}\dots\delta_{k}\right)  $. If
there exists an $m<k$ such that $\sum_{i=1}^{m}L_{\delta_{i}}=n$, then
$\left(  \delta\right)  \leq\left(  \left(  \delta_{1}\dots\delta_{m}\right)
\right)  $, which is impossible; and, if not, there is an $m$ with $\sum
_{i=1}^{m}L_{\delta_{i}}<n$ and $\sum_{i=1}^{m+1}L_{\delta_{i}}>n$. But then
$\left(  \delta_{1}\dots\delta_{m+1}\right)  \in E_{n}$, and $\left(
\delta\right)  \leq\left(  \left(  \delta_{1}\dots\delta_{m+1}\right)
\right)  $, so this is equally impossible. Thus the projections in the two
families add up to $1$, and Lemma \ref{LemSub.6} is proved.
\end{proof}

\begin{figure}[ptb]
\begin{picture}(360,563)(0,-12)
\put(0,0){\includegraphics
[bb=257 252 751 1008,clip,height=551bp,width=360bp]{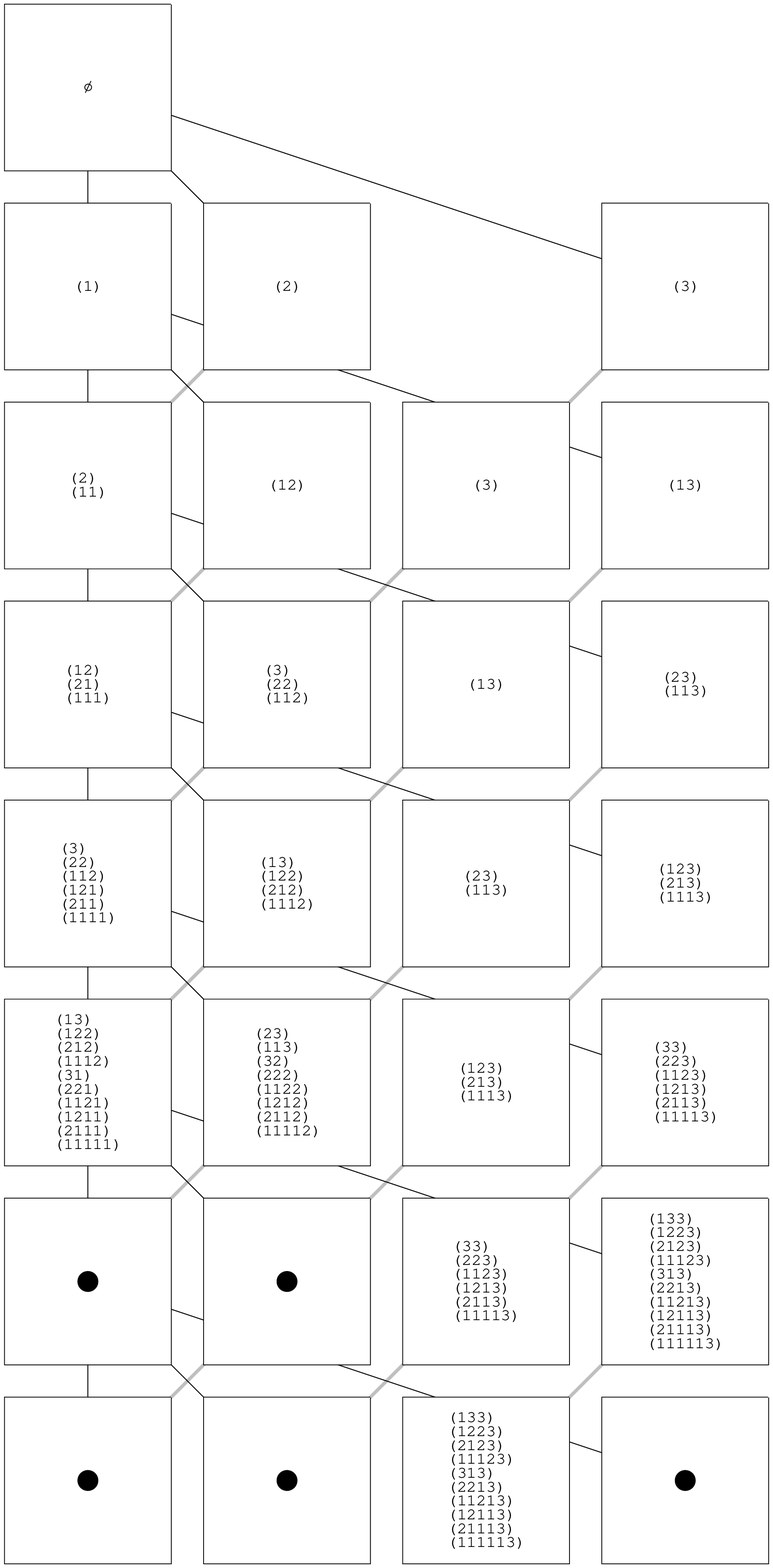}}
\end{picture}
\caption{$d=3$; $L=\left\{  1,2,4\right\}  $; first matrix column
$=(1\;1\;0\;1)^{\mathrm{t}}$. Compare with Figure \ref{DustDiag124} and
Example \ref{ExaSub.7}.}%
\label{BratDiag124}%
\end{figure}

\begin{figure}[ptb]
\begin{picture}(339,545)
\put(0,533){\makebox(0,12)[l]{\textbf{Level 1:}}}
\put(0,273){\includegraphics
[bb=82 0 281 363,height=272.25bp,width=149.25bp]{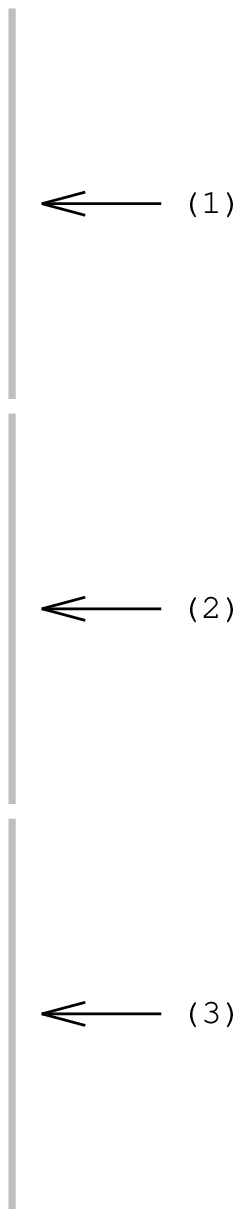}}
\put(150,533){\makebox(0,12)[l]{\textbf{Level 3:}}}
\put(149,273){\includegraphics
[bb=55 0 308 363,height=272.25bp,width=189.75bp]{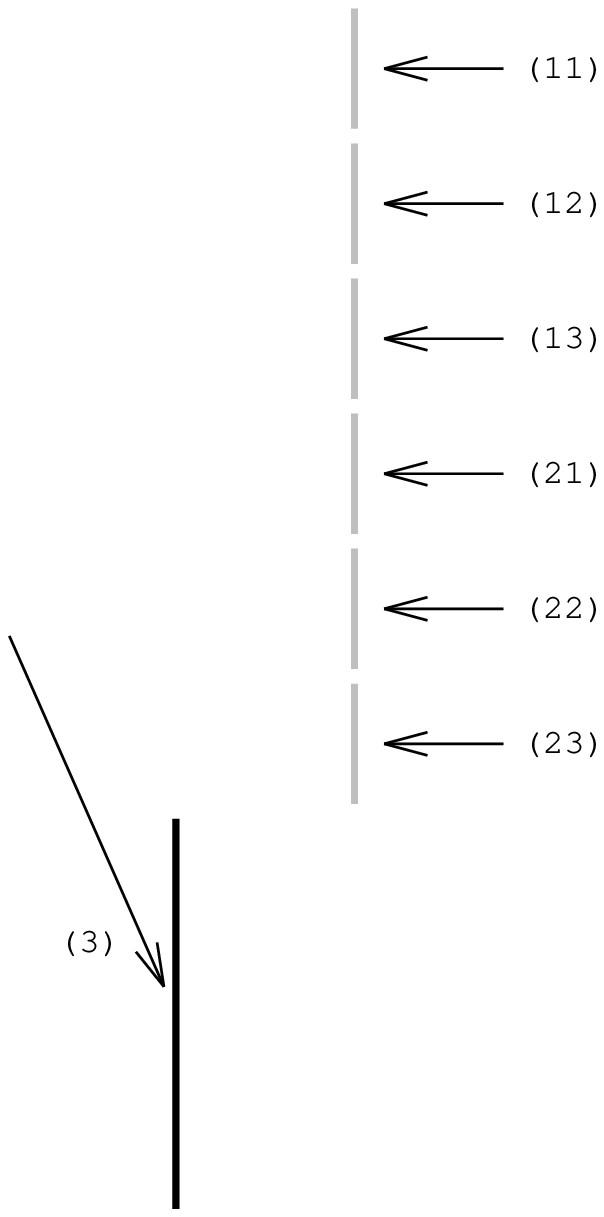}}
\put(0,260){\makebox(0,12)[l]{\textbf{Level 2:}}}
\put(0,0){\includegraphics
[bb=82 0 281 363,height=272.25bp,width=149.25bp]{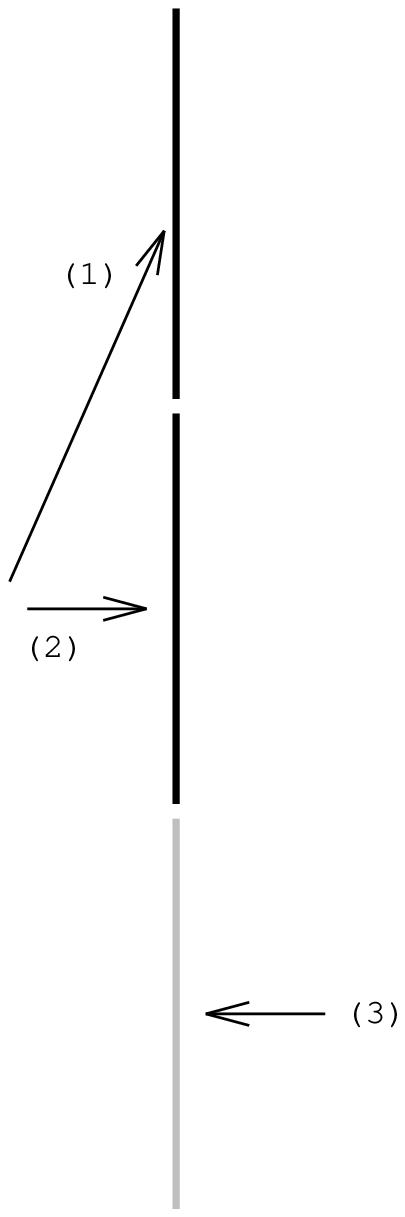}}
\put(150,260){\makebox(0,12)[l]{\textbf{Level 4:}}}
\put(149,0){\includegraphics
[bb=55 0 308 363,height=272.25bp,width=189.75bp]{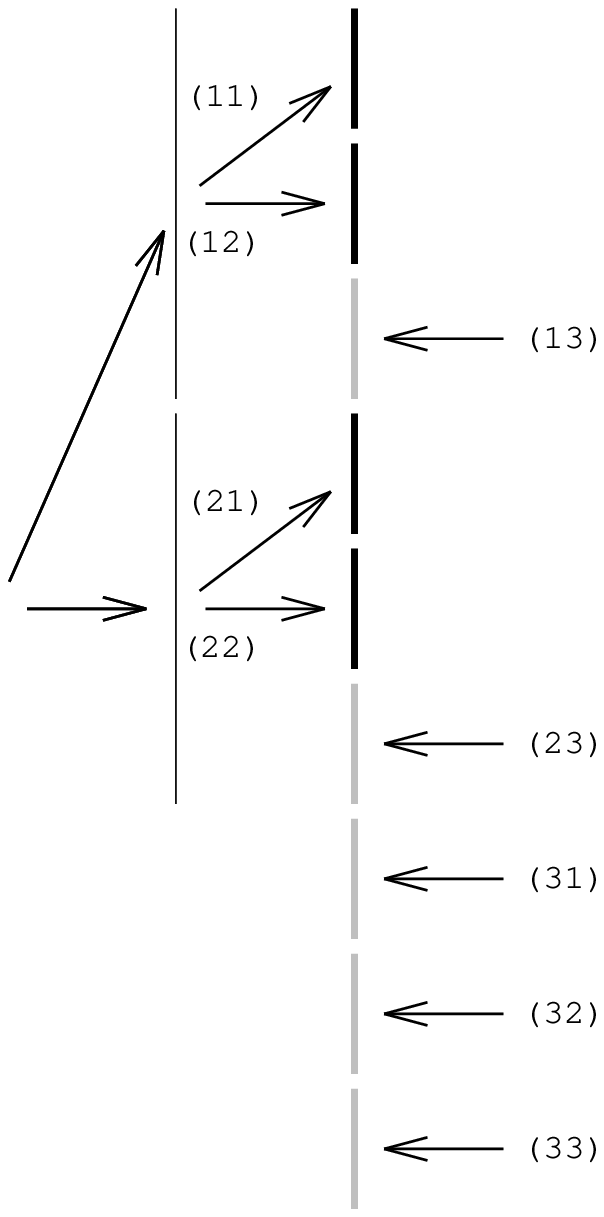}}
\end{picture}
\caption{$L=\{2,2,3\}$; levels 1--4. Compare Figures \ref{DustDiag223a},
\ref{DustDiag223b}, and \ref{DustDiag223c} with Figure \ref{BratDiag223}.}%
\label{DustDiag223a}%
\end{figure}

\begin{figure}[ptb]
\begin{picture}(312,545)
\put(0,533){\makebox(0,12)[l]{\textbf{Level 5:}}}
\put(0,273){\includegraphics
[bb=32 0 448 363,height=272.25bp,width=312bp]{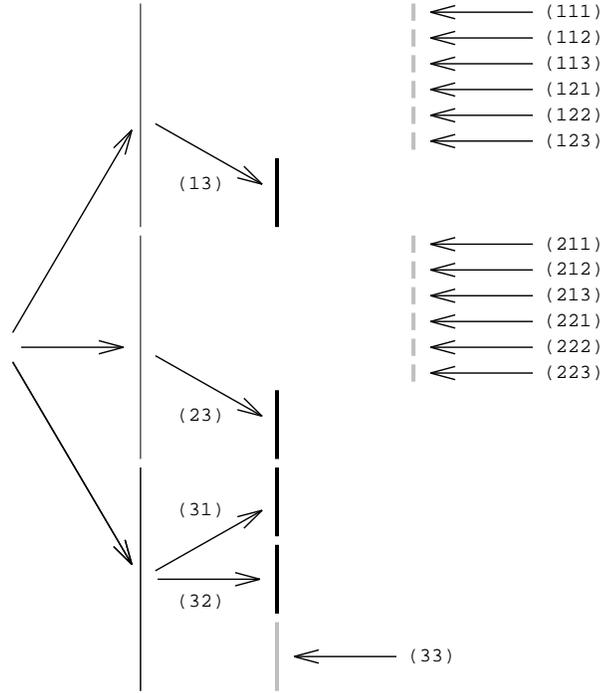}}
\put(0,260){\makebox(0,12)[l]{\textbf{Level 6:}}}
\put(0,0){\includegraphics
[bb=32 0 448 363,height=272.25bp,width=312bp]{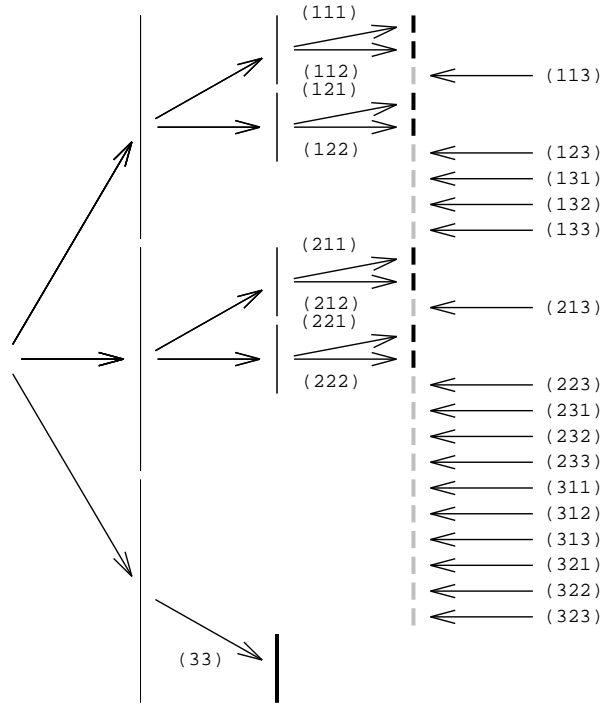}}
\end{picture}
\caption{$L=\{2,2,3\}$; levels 5--6.}%
\label{DustDiag223b}%
\end{figure}

\begin{figure}[ptb]
\begin{picture}(360,545)(9,0)
\put(0,533){\makebox(0,12)[l]{\textbf{Level 7:}}}
\put(0,0){\includegraphics
[bb=117 0 609 726,height=544.5bp,width=369bp]{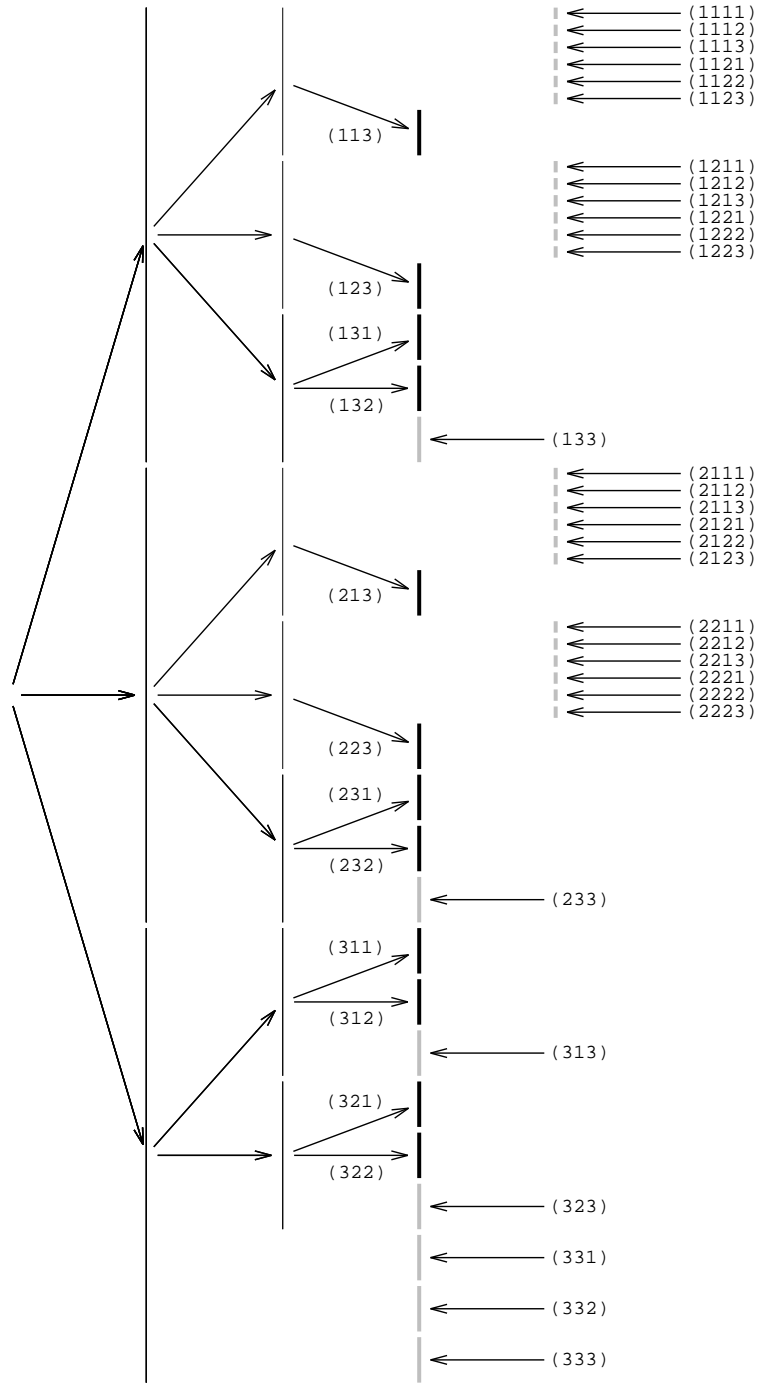}}
\end{picture}
\caption{$L=\{2,2,3\}$; level 7.}%
\label{DustDiag223c}%
\end{figure}

\begin{figure}[ptb]
\begin{picture}(201.3,554.4)
\put(0,0){\includegraphics
[bb=321 0 687 1008,height=554.4bp,width=201.3bp]{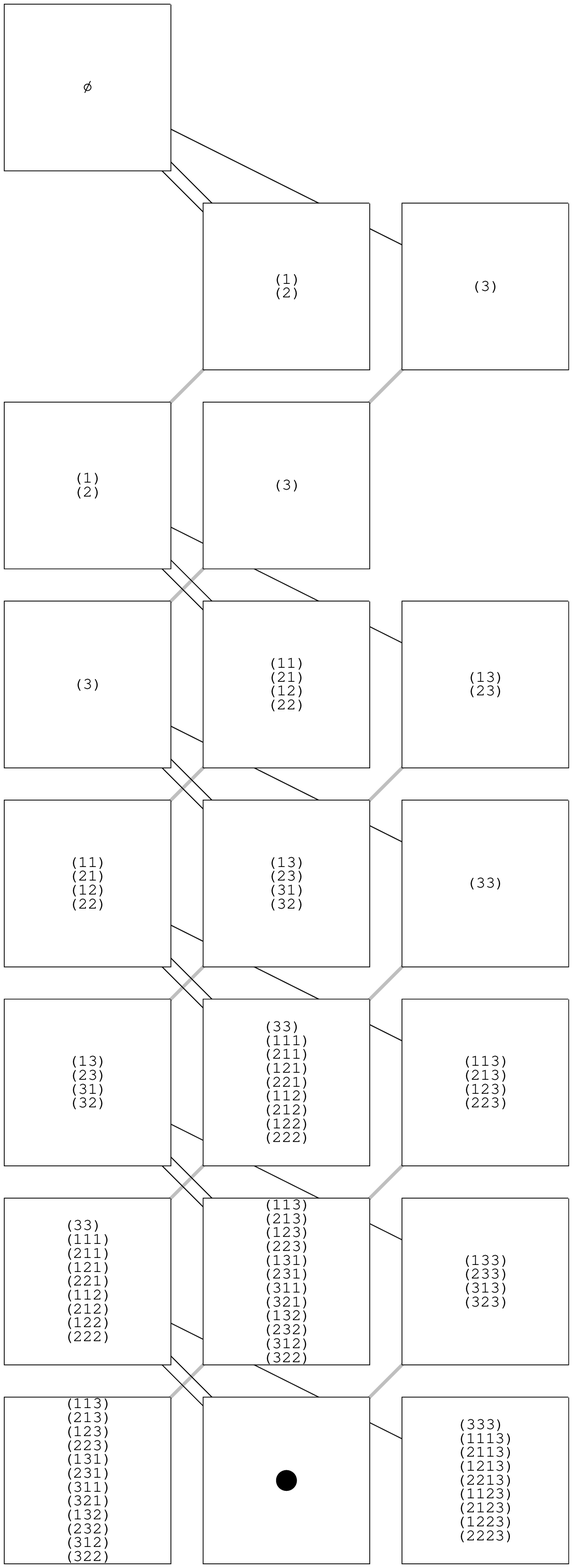}}
\end{picture}
\caption{$d=3$; $L=\left\{  2,2,3\right\}  $; first matrix column
$=(0\;2\;1)^{\mathrm{t}}$. Compare with Figures \ref{DustDiag223a},
\ref{DustDiag223b}, and \ref{DustDiag223c}.}%
\label{BratDiag223}%
\end{figure}

\begin{example}
\label{ExaSub.7}The procedure in the proof of Lemma \textup{\ref{LemSub.6}}
may be illustrated graphically as follows: Let $\alpha=\left(  \alpha_{1}%
\dots\alpha_{p}\right)  \in L^{-1}\left(  n\right)  $, and set
\[
E_{n}\left(  \alpha\right)  =\left\{  \gamma\mid\exists q<p\text{ such that
}\gamma=\left(  \alpha_{1}\dots\alpha_{q}\gamma_{q+1}\right)  \text{ and
}L\left(  \gamma\right)  >n\right\}  .
\]
For the example $d=3$, $L_{1}=1$, $L_{2}=2$, $L_{3}=4$, we have%
\begin{align*}
E_{4}\left(  \left(  1111\right)  \right)   &  =\left\{  \left(  1112\right)
,\left(  1113\right)  ,\left(  113\right)  ,\left(  13\right)  \right\}  ;\\
E_{4}\left(  \left(  121\right)  \right)  \diagdown E_{4}\left(  \left(
1111\right)  \right)   &  =\left\{  \left(  122\right)  ,\left(  123\right)
\right\}  ;\\
E_{4}\left(  \left(  22\right)  \right)   &  \text{\qquad contains a new
element }\left(  23\right)  ;\\
E_{4}\left(  \left(  211\right)  \right)   &  \text{\qquad contains the rest,
i.e., }\left(  212\right)  ,\left(  213\right)  .
\end{align*}
This is illustrated in Figure \ref{DustDiag124}. Elements from $L^{-1}\left(
4\right)  $ have arrows coming from the left ending at dark bars, while
elements from $E_{4}$ have arrows coming from the right ending at light bars.
The points in $L^{-1}\left(  4\right)  \cup E_{4}$ together represent Cuntz
algebra generators. The ordinary diagram\label{LOSBrattelidiagrams_3} for this
$\mathfrak{A}_{L}$ is illustrated in Figure \textup{\ref{BratDiag124}.}
Corresponding diagrams for $L_{1}=L_{2}=2$, $L_{3}=3$ are shown in Figures
\textup{\ref{DustDiag223a}, \ref{DustDiag223b}, \ref{DustDiag223c},} and
\textup{\ref{BratDiag223}. }
\end{example}

\begin{proof}
[Proof of Theorem \textup{\ref{ThmSub.2}} and Remark \textup{\ref{RemSub.3}}%
]Referring to Lemma \ref{LemSub.6}, define
\begin{equation}
E_{n}\left(  0\right)  =L^{-1}\left(  n\right)  \label{eqSub.18}%
\end{equation}
and\label{LOSEnm}%
\begin{equation}
E_{n}\left(  m\right)  =\left\{  \gamma\in E_{n}\mid L\left(  \gamma\right)
=n+m\right\}  \label{eqSub.19}%
\end{equation}
for $m=1,2,\dots,L_{d}-1$; for greater $m$'s, $E_{n}\left(  m\right)  $
becomes the empty set. $E_{n}\left(  m\right)  $ may also be the empty set for
some $m\in\left\{  0,\dots,L_{d}-1\right\}  $, but we will prove in a moment
that if the greatest common divisor of $L_{1},\dots,L_{d}$ is $1$, this only
happens for finitely many pairs $\left(  n,m\right)  $. Now, define
$\mathfrak{A}^{\left(  n,m\right)  }$ as the linear span of elements
$e_{\alpha\gamma}^{{}}=s_{\alpha}^{{}}s_{\gamma}^{\ast}$ with $\alpha
,\gamma\in E_{n}\left(  m\right)  $, $m=0,\dots,L_{d}-1$, with the convention
that $\mathfrak{A}^{\left(  n,m\right)  }$ is empty if $E_{n}\left(  m\right)
={}$\textup{\o} and $\mathfrak{A}^{\left(  0,0\right)  }=\mathbb{C}\openone$,
$\mathfrak{A}^{\left(  0,m\right)  }=0$ for $m=1,\dots,L_{d}-1$. It follows
from Lemma \ref{LemSub.5} that each $\mathfrak{A}^{\left(  n,m\right)  }$ is a
full $\#\left(  E_{n}\left(  m\right)  \right)  \times\#\left(  E_{n}\left(
m\right)  \right)  $ matrix algebra, and that the units of $\mathfrak
{A}^{\left(  n,m\right)  }$ are orthogonal and add up to $\openone$ as $m$
runs over $0,1,\dots,L_{d}-1$ for fixed $n$. Put\label{LOSAnm}
\begin{equation}
\mathfrak{A}_{n}=\bigoplus_{k=0}^{L_{d}-1}\mathfrak{A}^{\left(  n,k\right)  }.
\label{eqSub.20}%
\end{equation}
If $L\left(  \gamma\right)  =n$, then%
\begin{equation}
\left(  \gamma\right)  =\sum_{i=1}^{d}\left(  \gamma i\right)
\label{eqSub.21}%
\end{equation}
and%
\begin{equation}
\gamma i\in E_{n+1}\left(  L_{i}-1\right)  ,\qquad i=1,\dots,d;
\label{eqSub.22}%
\end{equation}
and hence $\mathfrak{A}^{\left(  n,0\right)  }$ is partially embedded in
$\mathfrak{A}^{\left(  n+1,m\right)  }$ with multiplicity equal to the number
of $k$'s such that $L_{k}-1=m$. We also have%
\begin{equation}
E_{n+1}\left(  m\right)  \subset E_{n}\left(  m+1\right)  \label{eqSub.23}%
\end{equation}
for $m=0,1,\dots,L_{k}-2$, and thus $\mathfrak{A}^{\left(  n,m+1\right)  }$ is
embedded into $\mathfrak{A}^{\left(  n+1,m\right)  }$ with multiplicity $1$
for $m=0,1,\dots,L_{k}-2$. It follows that $\mathfrak{A}_{n}$ is indeed an
increasing sequence of finite-dimensional subalgebras, and in particular
$\mathfrak{A}_{n}$ contains all monomials $s_{\alpha}^{{}}s_{\gamma}^{\ast}%
\in\mathfrak{A}_{L}$ of grade $\leq n$. Thus $\bigcup_{n}\mathfrak{A}_{n}$ is
dense in $\mathfrak{A}_{L}$, reestablishing that $\mathfrak{A}_{L}$ is an
AF-algebra, and the description of the embedding $\mathfrak{A}_{n}%
\hookrightarrow\mathfrak
{A}_{n+1}$ proves Remark \ref{RemSub.3}. The remaining statements in Theorem
\ref{ThmSub.2} will be proved after Lemma \ref{LemSub.8bis},
below.\renewcommand{\qed}{}
\end{proof}

By Proposition \ref{ProKMS.1}, the state defined on $\mathcal{O}_{d}$ by
(\ref{eqSub.13}) is a $\left(  \sigma,\beta\right)  $-KMS state. Thus the
restriction to $\mathfrak{A}_{L}=\mathcal{O}_{d}^{\sigma}$ is a trace state.
Now the embeddings $\mathfrak{A}_{n}\hookrightarrow\mathfrak{A}_{n+1}$ are
given by a constant embedding matrix $J:$ if, for example, $d=4$, $L_{1}=1$,
$L_{2}=L_{3}=3$, $L_{4}=4$, then%
\[
J=%
\begin{pmatrix}
1 & 1 & 0 & 0\\
0 & 0 & 1 & 0\\
2 & 0 & 0 & 1\\
1 & 0 & 0 & 0
\end{pmatrix}
.
\]
In general $J$ has the property that $J^{n}$ has strictly positive matrix
elements for some positive $n$. This is in fact equivalent to the property
that the numbers $L_{1},\dots,L_{d}$ have greatest common divisor $1$, which
may be seen as follows:

\begin{lemma}
\label{LemSub.8bis}Let $\mathcal{P}$ be the semigroup generated by
$L_{1},\dots,L_{d}$:\label{LOSP}%
\begin{equation}
\mathcal{P}=\left\{  \sum_{k=1}^{d}n_{k}L_{k}\biggm|n_{k}\in\mathbb{N}%
\cup\left\{  0\right\}  \right\}  . \label{eqSub.24}%
\end{equation}
Then $\mathbb{N\,}\diagdown\mathcal{P}$ is finite.
\end{lemma}

\begin{proof}
Since $L_{1},\dots,L_{d}$ have greatest common divisor $1$, there are
$n_{k}\in\mathbb{Z}$ such that%
\[
\sum_{k=1}^{d}n_{k}L_{k}=1,
\]
and hence there are $x_{1},x_{2}\in\mathcal{P}$ such that%
\[
x_{1}=x_{2}+1\text{.}%
\]
Now, if \emph{ad absurdum} $\mathbb{N\,}\diagdown\mathcal{P}$ is infinite we
may find arbitrarily large $y\in\mathbb{N\,}\diagdown\mathcal{P}$, but then
$y-x_{1}$, $y-x_{2}$ are not contained in $\mathcal{P}$; thus $y-x_{1}-x_{1}$,
$y-x_{1}-x_{2}$, $y-x_{2}-x_{2}$ are not in $\mathcal{P}$, etc., and thus we
can find arbitrarily long sequences of the form $\left(  z,z+1,z+2,\dots
,z+k\right)  $ not in $\mathcal{P}$. But as $\mathcal{P}$ contains
$\mathbb{N}L_{1}$, this is impossible. Thus $\mathbb{N\,}\diagdown\mathcal{P}$
is finite.
\end{proof}

\begin{proof}
[End of proof of Theorem \textup{\ref{ThmSub.2}}]Since any node in the
Bratteli diagram is connected to a node of the form $\left(  n,0\right)  $
further down, and $\left(  n,0\right)  $ is connected to all nodes $\left(
n+m,0\right)  $ where $m\in\mathcal{P}$, it follows that all nodes in a row
will be connected to all nodes in some row further down, which means that
$J^{n}$ has strictly positive matrix elements for some $n\in\mathbb{N}$.
Therefore $\mathfrak{A}_{L}$ is simple \cite{Bra72}, and $\mathfrak{A}$ has a
unique trace state \cite[Theorem 6.1]{Eff81}, \cite{Tor91}. This ends the
proof of Theorem \ref{ThmSub.2}.
\end{proof}

\begin{remark}
\label{RemSub.8}The semigroup $\mathcal{P}$ defined by \textup{(\ref{eqSub.24}%
)} can be read off the diagram of $\mathfrak{A}_{L}$ as follows:
$n\in\mathcal{P}$ if and only if the node $\left(  n,0\right)  $ actually
occurs in the diagram, i.e.,, if and only if $L^{-1}\left(  n\right)  \neq{}%
$\textup{\o}. To decide which $\left(  n,m\right)  $ actually occurs, start
with the vector $\left(  0,0\right)  =1$, $\left(  0,m\right)  =0$,
$m=1,\dots,L_{d}-1$, and apply the incidence matrix $J$. For example, in the
example illustrated in Figure \textup{\ref{BratDiag4458},} with $L=\left\{
4,4,5,8\right\}  $, we have $\mathcal{P}=\left\{  4,5,8,9,10,12,13,14,\dots
\right\}  $, while in the right-hand example in Figure
\textup{\ref{BratDiagsn4},} we have $\mathcal{P}=\left\{
3,6,7,9,10,12,13,14,\dots\right\}  $ \textup{(}both $\mathcal{P}$'s continuing
with no further gaps in the sequence\/\textup{).}
\end{remark}

\begin{figure}[ptb]
\begin{picture}(360,192)(-180,-96)
\put(-86,0){\makebox(0,0)[r]{$J={}$}}
\put(0,0){\makebox(0,0){
$
\begin{pmatrix}
\vphantom{m_1^1\vdots}m_1 & 1 & 0 & 0 & 0 & 0 & \cdots& 0 & 0 & 0\\
\vphantom{m_1^1\vdots}0 & 0 & 1 & 0 & 0 & 0 & \cdots& 0 & 0 & 0\\
\vphantom{m_1^1\vdots}0 & 0 & 0 & 1 & 0 & 0 & \cdots& 0 & 0 & 0\\
\vphantom{m_1^1\vdots}m_2 & 0 & 0 & 0 & 1 & 0 & \cdots& 0 & 0 & 0\\
\vphantom{m_1^1\vdots}0 & 0 & 0 & 0 & 0 & 1 & \cdots& 0 & 0 & 0\\
\vphantom{m_1^1\vdots}0 & 0 & 0 & 0 & 0 & 0 & \ddots& 0 & 0 & 0\\
\vphantom{m_1^1\vdots}\vdots& \vdots& \vdots& \vdots& \vdots&
& \ddots& \ddots& \vdots& \vdots\\
\vphantom{m_1^1\vdots}0 & 0 & 0 & 0 & 0 & 0 & & 0 & 1 & 0\\
\vphantom{m_1^1\vdots}0 & 0 & 0 & 0 & 0 & 0 & \cdots& 0 & 0 & 1\\
\vphantom{m_1^1\vdots}m_k & 0 & 0 & 0 & 0 & 0 & \cdots& 0 & 0 & 0
\end{pmatrix}
$}}
\put(90,81){\makebox(0,0)[l]{$\scriptstyle\leftarrow\text{ place }M_1$}}
\put(90,25.33){\makebox(0,0)[l]{$\scriptstyle\leftarrow\text{ place }M_2$}}
\put(90,-86){\makebox(0,0)[l]{$\scriptstyle\leftarrow\text{ place }M_k$}}
\end{picture}
\caption{Incidence matrix.}%
\label{IncMat2}%
\end{figure}

\begin{remark}
\label{RemSub.9}The result on the unique trace state cited at the end of the
proof above is actually related to the classical Perron--Frobenius
theorem\label{LOSPerronFrobeniustheorem} \cite{Gan59,Per07,Fro08}. If
$v^{\left(  n\right)  }$ is the value of the trace state on the minimal
projections in $\mathfrak{A}_{m}^{\left(  n\right)  }=\mathfrak{A}_{{}%
}^{\left(  n,m\right)  }$, and $v^{\left(  n\right)  }=\left(  v_{0}^{\left(
n\right)  },\dots,v_{L_{d}-1}^{\left(  n\right)  }\right)  $, then%
\begin{equation}
v^{\left(  n-1\right)  }=v^{\left(  n\right)  }J, \label{eqSub.25}%
\end{equation}
provided $n$ is so large that the Bratteli
diagram\label{LOSBrattelidiagrams_4} has stabilized, i.e., $\mathfrak{A}%
_{m}^{\left(  n\right)  }\neq\left\{  0\right\}  $ for $m=0,1,\dots,L_{d}-1$.
Since the components of $v^{\left(  n\right)  }$ have to be nonnegative, the
only solutions of \textup{(\ref{eqSub.25})} are such that each $v^{\left(
n\right)  }$ \textup{(}for large $n$\textup{)} is a multiple of the
Perron--Frobenius eigenvector $v$\label{LOSPerronFrobeniusrighteigenvector_2}
of $J$, i.e.,%
\begin{equation}
vJ=\lambda_{0}v. \label{eqSubNew.27}%
\end{equation}
Recall that the irreducibility of $J$ \textup{(}some power of $J$ has only
positive matrix units\textup{)} implies that $J$ has a simple positive
eigenvalue $\lambda_{0}$ such that $\lambda_{0}>\left|  \mu\right|  $ for all
other eigenvalues $\mu$, and the corresponding one-dimensional eigenspace is
spanned by a vector $v$ with positive components. All $v^{\left(  n\right)  }%
$'s for large $n$ are multiples of $v$. \textup{(}This is because of the
uniqueness of the normalized positive solution of \textup{(\ref{eqSubNew.27}%
),} together with the fact that $v^{\left(  n\right)  }={}$const.$\cdot
e^{-\beta n}v$ is indeed a solution. Note that the Perron--Frobenius
eigenvalue\label{LOSPerronFrobeniuseigenvalue_3} of $J$ is $e^{\beta}%
$.\textup{)} Thus $v$ may be computed explicitly in the examples by choosing
$n$ so large that the diagram has stabilized, and using
\textup{(\ref{eqSub.13})} and \textup{(\ref{eqSub.19})} to evaluate the trace
on the minimal projections in $\mathfrak{A}_{m}^{\left(  n\right)  }$. The
result is surprisingly simple; see \textup{(\ref{eqRemNewBis.17})} in the next chapter.

Let us give the details of the graphic description of the embedding of
$\mathfrak{A}_{n}$ into $\mathfrak{A}_{n+1}$. Suppose that the integral
weights of Theorem \textup{\ref{ThmSub.2}} are $1\leq L_{1}\leq L_{2}\leq
\dots\leq L_{d}$ with possible multiple occupancy. Let
\begin{equation}%
\begin{aligned}
L_{i}  &  =M_{1}\text{\qquad for }1\leq i\leq i_{1},\\
L_{i}  &  =M_{2}\text{\qquad for }i_{1}<i\leq i_{2},\\
{}  &  \mathrel{\mkern4.5mu\vdots\mkern4.5mu}\phantom{M_{2}}\qquad\vdots\\
L_{i}  &  =M_{k}\text{\qquad for }i_{k-1}<i\leq i_{k}=d.
\end{aligned}%
\label{eqLM}%
\end{equation}
Let $m_{j}=i_{j}-i_{j-1}$ be the multiplicities. Then, after stabilization,
the partial embedding of $\mathfrak{A}_{0}^{\left(  n-1\right)  }$ into the
factors $\mathfrak{A}_{m}^{\left(  n\right)  }$, $m=0,\dots,L_{d}-1$, are
given by the diagram in \textup{(\ref{eqBratBlank3232})} below \textup{(}%
illustrated in the case $L_{1}=1$\textup{):}%
\begin{equation}
\begin{minipage}{288\unitlength} \begin{picture}(288,144)(-19,-48) \put
(-19,0){\includegraphics[bb=0 98 288 190,height=92bp,width=288bp]{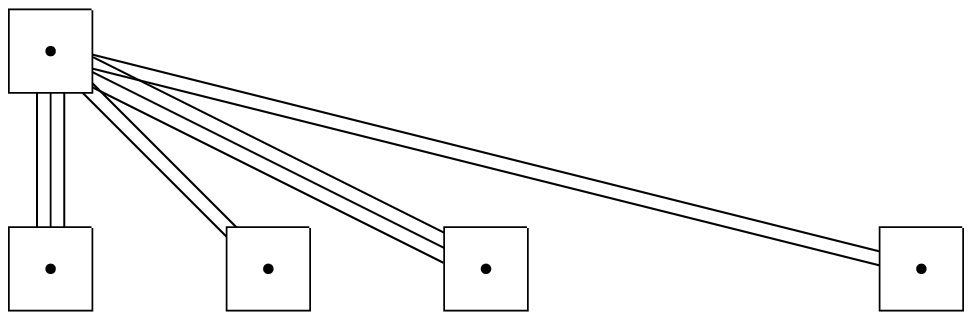}%
}\put(0,-21){\makebox(0,0)[b]{$M_{1}$}} \put(62.5,-21){\makebox(0,0)[b]{$M_{2}%
$}} \put(125,-21){\makebox(0,0)[b]{$M_{3}$}} \put(250,-21){\makebox
(0,0)[b]{$M_{k}$}} \put(187.5,15){\makebox(0,0){$\dots$}} \put(0,-42){\makebox
(0,0)[b]{$m_{1}$ lines}} \put(62.5,-42){\makebox(0,0)[b]{$m_{2}$ lines}}%
\put(125,-42){\makebox(0,0)[b]{$m_{3}$ lines}} \put(250,-42){\makebox
(0,0)[b]{$m_{k}$ lines.}} \end{picture}\end{minipage} \label{eqBratBlank3232}%
\end{equation}
Then $J^{-1}$ is given by the matrix%
\begin{equation}
\setlength{\unitlength}{1bp}\makebox[0bp]{\hss$\vcenter{\begin{picture}%
(360,212)(-180,-108) \put(-102,0){\makebox(0,0)[r]{$J^{-1}={}$}}%
\put(0,0){\makebox(0,0){ $ \begin{pmatrix} \vphantom{\frac{m_{j}}{m_{k}}%
\vdots} 0 & 0 & 0 & 0 & 0 & \cdots& \cdots& 0 & 0 & \frac{1}{m_{k}%
}\\ \vphantom{\frac{m_{j}}{m_{k}}\vdots} 1 & 0 & 0 & 0 & 0 & \cdots
& \cdots& 0 & 0 & -\frac{m_{1}}{m_{k}}\\ \vphantom{\frac{m_{j}}{m_{k}}\vdots
}0 & 1 & 0 & 0 & 0 & \cdots& \cdots& 0 & 0 & 0\\ \vphantom{\frac{m_{j}}{m_{k}%
}\vdots} 0 & 0 & 1 & 0 & 0 & \cdots& \cdots& 0 & 0 & 0\\ \vphantom{\frac
{m_{j}}{m_{k}}\vdots} 0 & 0 & 0 & 1 & 0 & \cdots& \cdots& 0 & 0 & -\frac
{m_{2}}{m_{k}}\\ \vphantom{\frac{m_{j}}{m_{k}}\vdots} \vdots& \vdots
& \vdots& \vdots& \ddots& \ddots& & & \vdots& \vdots\\ \vphantom{\frac{m_{j}%
}{m_{k}}\vdots} \vdots& \vdots& \vdots& \vdots& & \ddots& \ddots
& & \vdots& -\frac{m_{j}}{m_{k}}\\ \vphantom{\frac{m_{j}}{m_{k}}\vdots}%
\vdots& \vdots& \vdots& \vdots& & & \ddots& \ddots& \vdots& \vdots
\\ \vphantom{\frac{m_{j}}{m_{k}}\vdots} 0 & 0 & 0 & 0 & 0 & \cdots
& \cdots& 1 & 0 & 0\\ \vphantom{\frac{m_{j}}{m_{k}}\vdots}%
0 & 0 & 0 & 0 & 0 & \cdots& \cdots& 0 & 1 & 0 \end{pmatrix} $}} \put
(108,-32.5){\makebox(0,0)[l]{$\scriptstyle\leftarrow\text{ place } M_{j}+1$,}%
}\put(180,-43){\makebox(0,0)[r]{$\scriptstyle j=1,\dots,k-1$}} \end{picture}%
}$\hss} \label{eqSubNew.28}%
\end{equation}
The characteristic polynomial for the corresponding inverse $J^{-1}$ is
proportional to
\begin{equation}
p_{m}\left(  x\right)  =m_{k}x^{M_{k}}+m_{k-1}x^{M_{k-1}}+\dots+m_{1}x^{M_{1}%
}-1 \label{eqSubNew.29}%
\end{equation}
Since $\sum_{i=1}^{k}m_{i}e^{-\beta M_{i}}=1$, we see that $x=e^{-\beta}$ is
the unique positive root for this polynomial. Thus $e^{\beta}$ is the
Perron--Frobenius eigenvalue for $J$.
\end{remark}

\begin{remark}
\label{RemSubApr.11}Note that the implications \textup{(\ref{ProSub.1(1)})}
$\Leftarrow$ \textup{(\ref{ProSub.1(2)})} $\Leftrightarrow$
\textup{(\ref{ProSub.1(3)})} in Proposition \textup{\ref{ProSub.1}} remain
true even if $L_{1},\dots,L_{d}$ are not integers, by essentially the same
proof. This is because the action $\sigma^{\left(  L\right)  }$ defined by
\textup{(\ref{eqSub.1})} is almost periodic in all cases, and hence
$\mathfrak{A}_{L}$ is the closure of the linear span of $s_{\alpha}^{{}%
}s_{\gamma}^{\ast}$ with $L\left(  \alpha\right)  =L\left(  \gamma\right)  $
even in the general case, using \textup{(\ref{eqSub.5})} and the definition%
\[
L\left(  \alpha\right)  =\sum_{m=1}^{k}L_{\alpha_{m}}.
\]
It is no longer true that \textup{(\ref{ProSub.1(1)})} $\Rightarrow$
\textup{(\ref{ProSub.1(2)}).} Take for example $d=2$ and $L_{1}$, $L_{2}$
rationally independent irrational numbers of opposite sign. Then $\mathfrak
{A}_{L}$ is the GICAR algebra \cite{Dav96}.
\end{remark}

\clearpage\thispagestyle{empty}\setcounter{figurelink}{\value{figure}}

\part{\label{PartInv}Numerical AF-Invariants}

\ \markboth{\shorttitle}{\shorttitle} \cleardoublepage

\chapter{\label{Rem}\label{LOSdimensiongroup_2}The dimension group of
$\mathfrak{A}_{L}$}

\setcounter{figure}{\value{figurelink}} In this chapter and the following ones
we will construct isomorphism invariants for $\mathfrak{A}_{L}$ and try to
classify the $\mathfrak{A}_{L}$. It is known that there exists a complete
isomorphism invariant for AF-algebras $\mathfrak{A}$, namely the dimension
group. In the case that $\mathfrak{A}$ has a unit this is the triple $\left(
K_{0}\left(  \mathfrak{A}\right)  ,K_{0}\left(  \mathfrak{A}\right)
_{+},\left[  \openone\right]  \right)  $\label{LOSDAL_2} where $K_{0}\left(
\mathfrak
{A}\right)  $\label{LOSK0AL_2} is an abelian group, $K_{0}\left(  \mathfrak
{A}\right)  _{+}$\label{LOSK0ALplus_1} are the positive elements of
$K_{0}\left(  \mathfrak{A}\right)  $ relative to an order making $K_{0}\left(
\mathfrak
{A}\right)  $ into a Riesz ordered group without perforation, and $\left[
\openone\right]  $ is the class of the identity in $K_{0}\left(  \mathfrak
{A}\right)  $ (if $\mathfrak{A}$ is nonunital, replace $\left[  \openone
\right]  $ by the hereditary subset $\{\left[  p\right]  \mid p$ projection in
$\mathfrak{A}\}$ of $K_{0}\left(  \mathfrak{A}\right)  _{+}$). See
\cite{Eff81} for details on this and the following statements. (Connections to
ergodic theory are also described in \cite{Ver81a}, \cite{Ver81b}.) Let us now
specialize to the case that $\mathfrak{A}$ is given by a constant $N\times N$
incidence matrix $J$ (with nonnegative integer entries) which is primitive,
i.e., $J^{n}$ has only positive entries for some $n\in\mathbb{N}$. Then
$\mathfrak{A}$ is simple with a unique trace state $\tau$. In the case that
$K_{0}\left(  \mathfrak{A}\right)  \cong\mathbb{Z}^{N}$, this class of
AF-algebras (or rather dimension groups) has been characterized intrinsically
in \cite[Theorems 3.3 and 4.1]{Han81}. In general when $J$ is an $n\times
n=L_{d}\times L_{d}$ matrix with nonnegative entries, the dimension group is
the inductive limit%
\begin{equation}
\mathbb{Z}^{N}\overset{J}{\longrightarrow}\mathbb{Z}^{N}\overset
{J}{\longrightarrow}\mathbb{Z}^{N}\overset{J}{\longrightarrow}\cdots
\label{eqRemNewBis.1}%
\end{equation}
with order generated by the order defined by%
\begin{equation}
\left(  m_{1},\dots,m_{N}\right)  \geq0\Longleftrightarrow m_{i}%
\geq0\text{\qquad on }\mathbb{Z}^{N}. \label{eqRemNew.0}%
\end{equation}
This group can be computed explicitly as a subgroup of $\mathbb{R}^{N}$ as
follows when $\det\left(  J\right)  \neq0$ (as it is in our case): Put%
\begin{equation}
G_{m}=J^{-m}\left(  \mathbb{Z}^{N}\right)  ,\qquad m=0,1,\dots,
\label{eqRemNewBis.3}%
\end{equation}
and equip $G_{m}$ with the order%
\begin{equation}
G_{m}^{+}=J^{-m}\left(  \left(  \mathbb{Z}^{N}\right)  ^{+}\right)  .
\label{eqRemNewBis.4}%
\end{equation}
Then\label{LOSG0_1}
\begin{equation}
G_{0}\subset G_{1}\subset G_{2}\subset\cdots, \label{eqRemNewBis.5}%
\end{equation}
and%
\begin{equation}
K_{0}\left(  \mathfrak{A}_{L}\right)  =\bigcup_{m}G_{m}, \label{eqRemNewBis.6}%
\end{equation}
a subgroup of $\mathbb{R}^{N}$ (containing $\mathbb{Z}^{N}$), with order
defined by%
\begin{equation}
g\geq0\text{\qquad if }g\geq0\text{ in some }G_{m}. \label{eqRemNewBis.7}%
\end{equation}
The action of the trace state $\tau$ on $K_{0}\left(  \mathfrak{A}_{L}\right)
$ may be computed as follows: If $\lambda$ is the Frobenius eigenvalue of $J$,
and $\alpha=\left(  \alpha_{1},\dots,\alpha_{N}\right)  $ is a corresponding
eigenvector\label{LOSPerronFrobeniuslefteigenvector_2} in the sense%
\begin{equation}
\alpha J=\lambda\alpha\label{eqRemNew.1}%
\end{equation}
(i.e., $J^{t}\alpha^{t}=\lambda\alpha^{t}$, see \cite[pp.\ 33--37]{Eff81}),
then if $\alpha$ is suitably normalized (by multiplying with a positive
factor), the trace applied to something at the $m$'th stage of
\begin{equation}
\underset{1\mathstrut}{\mathbb{Z}^{N}}\longrightarrow\underset{2\mathstrut
}{\mathbb{Z}^{N}}\longrightarrow\underset{3\mathstrut}{\mathbb{Z}^{N}%
}\longrightarrow\dots\longrightarrow\overset{%
\begin{array}
[c]{c}%
g\\
\makebox[0pt]{\hss$\displaystyle\cap$\hss}\makebox[0pt] {\hss\rule
[-0.15pt]{0.275pt}{6pt}\hss}%
\end{array}
}{\underset{m\mathstrut}{\mathbb{Z}^{N}}}\longrightarrow\cdots
\label{eqRemNewBis.9}%
\end{equation}
is%
\begin{equation}
\tau\left(  g\right)  =\lambda^{-m+1}\ip{\alpha}{g}, \label{eqRem.1}%
\end{equation}
where $\ip{\,\cdot\, }{\,\cdot\,}$ here denotes the usual inner product in
$\mathbb{R}^{N}$, i.e., $\ip{\alpha}{g}=\sum_{i=1}^{N}\alpha_{i}g_{i}$. Taking
$\alpha$ as the Frobenius eigenvector in (\ref{eqRem.1}) makes the ansatz well
defined: if $g\in G_{m}\subset G_{m+1}$, then%
\begin{equation}
\lambda^{-m+1}\ip{\alpha}{g}=\lambda^{-\left(  m+1\right)  +1}\ip{\alpha}{Jg}.
\label{eqRemNewBis.11}
\end{equation}
Thus $\tau$ is an additive character on $K_{0}\left(  \mathfrak{A}_{L}\right)
$, and up to normalization the unique positive such. If we identify $\left[
\openone\right]  _{0}$ with $\left(  1,0,0,\dots\right)  $ in the first
$\mathbb{Z}^{N}$, the normalization of $\alpha$ is $\alpha_{1}=1$.

Elements of the kernel of the additive real-valued character $\tau$ on
$K_{0}\left(  \mathfrak{A}\right)  $ are called infinitesimal elements. Thus
$K_{0}\left(  \mathfrak{A}\right)  $ is an extension of $\tau\left(
K_{0}\left(  \mathfrak{A}\right)  \right)  $ by the kernel of $\tau$. But in
general it is not the trivial extension, i.e.,%
\[
K_{0}\left(  \mathfrak{A}_{L}\right)  \ncong\tau\left(  K_{0}\left(
\mathfrak{A}_{L}\right)  \right)  \oplus\left(  \text{kernel of }\tau\right)
,
\]
which complicates classification; see Chapter \ref{Ext}.

Suppose we calculate the groups $\tau\left(  K_{0}\left(  \mathfrak{A}%
_{L}\right)  \right)  $ and $\ker\left(  \tau_{L}\right)  $ for a specific
pair, given by $L$ and $L^{\prime}$, say. Then if one of the two groups
$\tau\left(  K_{0}\left(  \mathfrak{A}_{L}\right)  \right)  $ or $\ker\left(
\tau_{L}\right)  $ is different for $L$ and for $L^{\prime}$, the AF-algebras
$\mathfrak{A}_{L}$ and $\mathfrak{A}_{L^{\prime}}$ are non-isomorphic. We
show, however, in Chapter \ref{Ext} that the AF-algebras can be non-isomorphic
even if the two groups agree for $L$ and $L^{\prime}$.

It can then be shown that the range of the trace on projections is
\linebreak $\tau\left(  K_{0}\left(  \mathfrak{A}_{L}\right)  \right)
\cap\left[  0,1\right]  $.

When $K_{0}\left(  \mathfrak{A}_{L}\right)  $ is given concretely in
$\mathbb{R}^{N}$ as above, the trace can be computed as\label{LOStaug}%
\begin{equation}
\tau\left(  g\right)  =\ip{\alpha}{g}, \label{eqRemNewBis.12}%
\end{equation}
where $g\in m$'th term $\mathbb{Z}^{N}$ is identified with its image
$J^{-m+1}g$ in $\mathbb{R}^{N}$; and the positive cone in $K_{0}\left(
\mathfrak{A}_{L}\right)  \subset\mathbb{R}^{N}$ identifies with those $g$ such
that $\tau\left(  g\right)  >0$, or $g=0$.

Let us now specialize to the case that $J=J_{L}$ has the special form we are
interested in . So assume that $1\leq L_{1}\leq L_{2}\leq\dots\leq L_{d}$,
that the greatest common divisor of $L_{1},\dots,L_{d}$ is $1$, and put
\begin{equation}
\left\{  L_{1},\dots,L_{d}\right\}  =\left\{  \underset{m_{1}^{\mathstrut}%
}{M_{1}},\underset{\dots\vphantom{m_{1}^{\mathstrut}}}{\dots},\underset
{m_{k}^{\mathstrut}}{M_{k}}\right\}  , \label{eqRemNewBis.13}%
\end{equation}
where $m_{i}$ is the multiplicity of $M_{i}$. Put%
\begin{equation}
N=M_{k}=L_{d}. \label{eqRemNewBis.14}%
\end{equation}
Then the incidence matrix $J$ is%
\begin{equation}
J=\mkern27mu\left(
\begin{array}
[c]{ccccccccccccccc}%
\llap{${\scriptstyle1}\mkern33mu$}0 & 1 & 0 & \cdots & 0 & 0 & 0 & \cdots &
0 & 0 & 0 & \cdots & \cdots & 0 & 0\\
0 & 0 & 1 &  & 0 & 0 & 0 & \cdots & 0 & 0 & 0 & \cdots & \cdots & 0 & 0\\
\vdots &  & \ddots & \ddots &  & \vdots & \vdots &  & \vdots & \vdots & \vdots
&  &  & \vdots & \vdots\\
0 & 0 & 0 & \ddots & 1 & 0 & 0 & \cdots & 0 & 0 & 0 & \cdots & \cdots & 0 &
0\\
\llap{${\scriptstyle M_{1}}\mkern27mu$}m_{1} & 0 & 0 &  & 0 & 1 & 0 & \cdots &
0 & 0 & 0 & \cdots & \cdots & 0 & 0\\
0 & 0 & 0 & \cdots & 0 & 0 & 1 &  & 0 & 0 & 0 & \cdots & \cdots & 0 & 0\\
\vdots & \vdots & \vdots &  & \vdots &  & \ddots & \ddots &  & \vdots & \vdots
&  &  & \vdots & \vdots\\
0 & 0 & 0 & \cdots & 0 & 0 & 0 & \ddots & 1 & 0 & 0 & \cdots & \cdots & 0 &
0\\
\llap{${\scriptstyle M_{2}}\mkern27mu$}m_{2} & 0 & 0 & \cdots & 0 & 0 & 0 &  &
0 & 1 & 0 & \cdots & \cdots & 0 & 0\\
0 & 0 & 0 & \cdots & 0 & 0 & 0 & \cdots & 0 & 0 & 1 &  & \cdots & 0 & 0\\
\vdots & \vdots & \vdots &  & \vdots & \vdots & \vdots &  & \vdots &  & \ddots
& \ddots &  & \vdots & \vdots\\
\vdots & \vdots & \vdots &  & \vdots & \vdots & \vdots &  & \vdots & \vdots &
& \ddots & \ddots &  & \vdots\\
0 & 0 & 0 & \cdots & 0 & 0 & 0 & \cdots & 0 & 0 & 0 &  & \ddots & 1 & 0\\
0 & 0 & 0 & \cdots & 0 & 0 & 0 & \cdots & 0 & 0 & 0 & \cdots &  & 0 & 1\\
\llap{${\scriptstyle M_{k}}\mkern27mu$}m_{k} & 0 & 0 & \cdots & 0 & 0 & 0 &
\cdots & 0 & 0 & 0 & \cdots & \cdots & 0 & 0
\end{array}
\right)  . \label{eqRemNewBis.15}%
\end{equation}
Let $x=e^{-\beta}$ be the unique solution in $\left(  0,1\right)  $ of%
\begin{equation}
1-\sum_{i}m_{i}x^{M_{i}}=0. \label{eqRemNewBis.16}%
\end{equation}
If\label{LOSalphaeigenvector_1}%
\begin{equation}
\alpha=\left(  1,e^{-\beta},e^{-2\beta},\dots,e^{-\left(  N-1\right)  \beta
}\right)  , \label{eqRemNewBis.17}%
\end{equation}
then $\alpha$ is the left Frobenius eigenvector%
\begin{equation}
\alpha J=e^{\beta}\alpha. \label{eqRemNewBis.18}%
\end{equation}
As explained before, we have the identification\label{LOSK0AL_3}%
\begin{equation}
K_{0}\left(  \mathfrak{A}_{L}\right)  =\bigcup_{n=0}^{\infty}J^{-n}%
\mathbb{Z}^{N}\qquad(\subset\mathbb{R}^{N}) \label{eqRemNewBis.19}%
\end{equation}
with the trace functional%
\begin{equation}
\tau\left(  y\right)  =\ip{\alpha}{y},\qquad y\in\bigcup_{n=0}^{\infty}%
J^{-n}\mathbb{Z}^{N}. \label{eqRemNewBis.20}%
\end{equation}
Using%
\begin{equation}
\ip{\alpha}{J^{-n}k}=e^{-n\beta}\ip{\alpha}{k}=e^{-n\beta}\sum_{i=1}^{N}%
k_{i}e^{-\beta\left(  i-1\right)  } \label{eqRemNewBis.21}%
\end{equation}
for $k\in\mathbb{Z}^{N}$, $n\in\mathbb{N}$, together with the fact that the
range of the trace is a subgroup of the additive group $\mathbb{R}$, it is
clear that the range of the trace is $\mathbb{Z}\left[  e^{-\beta}\right]  $,
i.e.,%
\begin{equation}
\tau\left(  K_{0}\left(  \mathfrak{A}_{L}\right)  \right)  =\mathbb{Z}\left[
e^{-\beta}\right]  , \label{eqRemNewBis.22}%
\end{equation}
and, furthermore, from \cite{Eff81},%
\begin{equation}
\tau\left(  \left\{  p\mid p\text{ projection in }\mathfrak{A}_{L}\right\}
\right)  =\mathbb{Z}\left[  e^{-\beta}\right]  \cap\left[  0,1\right]  .
\label{eqRemNewBis.23}%
\end{equation}
Now, if $m=\left(  m_{1},\dots,m_{N}\right)  $ is an element of the $k$'th
group%
\begin{equation}
\mathbb{Z}^{N}\overset{J}{\longrightarrow}\mathbb{Z}^{N}\overset
{J}{\longrightarrow}\mathbb{Z}^{N}\overset{J}{\longrightarrow}\cdots
\label{eqRemNewBis.24}%
\end{equation}
and $m$ is an infinitesimal element then $\ip{\alpha}{m^{t}}=0$, i.e.,%
\begin{equation}
\sum_{i=1}^{N}m_{i}\left(  e^{-\beta}\right)  ^{i-1}=0 \label{eqRemNewBis.25}%
\end{equation}
(where we include zero terms!). This sum $\sum_{i=1}^{N}m_{i}x^{i-1}$ is a
multiple of the minimal polynomial $p_{\beta}\left(  x\right)  $ having
$e^{-\beta}$ as a root. If this minimal polynomial happens to be $1-\sum
_{i}m_{i}x^{M_{i}}$, which has degree $N$, then there are no nontrivial
infinitesimal elements, and%
\begin{equation}
K_{0}\left(  G\right)  \cong\mathbb{Z}\left[  e^{-\beta}\right]  .
\label{eqRemNew.pound}%
\end{equation}
If $p_{\beta}$ has degree $\deg\left(  p_{\beta}\right)  <N$,\label{LOSdeg_1}
it follows that\label{LOSmdeg_1}%
\begin{equation}
\sum_{i=1}^{N}m_{k}x^{k-1}=p\left(  x\right)  \cdot\left(  \text{arbitrary
polynomial of degree }\leq\left(  N-1\right)  \text{-}\deg p_{\beta}\right)  .
\label{eqRemNewBis.27}%
\end{equation}
It follows that the group of infinitesimal elements of the $m$'th group
$\mathbb{Z}^{N}$ is isomorphic to%
\begin{equation}
\mathbb{Z}^{N\text{-}\deg p_{\beta}}, \label{eqRemNewBis.28}%
\end{equation}
and as $J$ maps these groups into each other, we obtain the infinitesimal
elements as an inductive limit%
\begin{equation}
\mathbb{Z}^{N\text{-}\deg p_{\beta}}\overset{J_{0}}{\longrightarrow}%
\mathbb{Z}^{N\text{-}\deg p_{\beta}}\overset{J_{0}}{\longrightarrow}\cdots,
\label{eqRem.star}%
\end{equation}
where $J_{0}$ is a restriction of $J$, so $J_{0}$ is an injective matrix with
integer entries, but the entries are no longer necessarily positive, as we see
in the examples. See Chapter \ref{Pediferient} for more details on $J_{0}$.

In conclusion, the complete invariant\label{LOSDAL_3}%
\begin{equation}
\left(  K_{0}\left(  \mathfrak{A}_{L}\right)  ,K_{0}\left(  \mathfrak{A}%
_{L}\right)  _{+},\left[  \openone\right]  \right)  \label{eqRemNewBis.30}%
\end{equation}
of the algebras $\mathfrak{A}_{L}$ defines an extension%
\begin{equation}
0\longrightarrow\ker\left(  \tau\right)  \overset{i}{\longrightarrow}%
K_{0}\left(  \mathfrak{A}_{L}\right)  \overset{\tau}{\longrightarrow
}\mathbb{Z}\left[  e^{-\beta}\right]  \longrightarrow0 \label{eqRemNewBis.31}%
\end{equation}
together with an element $\left[  \openone\right]  $ of $K_{0}\left(
\mathfrak{A}_{L}\right)  $ such that\label{LOStau}%
\begin{equation}
\tau\left(  \left[  \openone\right]  \right)  =1. \label{eqRemNewBis.32}%
\end{equation}
See Chapter \ref{Ext} for more details on these extensions. Concretely,
$K_{0}\left(  \mathfrak{A}_{L}\right)  $ is the subgroup (\ref{eqRemNewBis.19}%
) of $\mathbb{R}^{N}$, $\tau$ is given by (\ref{eqRemNewBis.20}) and
(\ref{eqRemNewBis.17}) and\label{LOSoneonesquare}%
\begin{equation}
\left[  \openone\right]  =\left(  1,0,0,\dots,0\right)  \label{eqRemNewBis.33}%
\end{equation}
and\label{LOSK0ALplus_2}
\begin{equation}
K_{0}\left(  \mathfrak{A}_{L}\right)  _{+}=\left\{  0\right\}  \cup\left\{
v\in K_{0}\left(  \mathfrak{A}_{L}\right)  \mid%
\ip{\alpha}{v}%
>0\right\}  . \label{eqRemNewBis.34}%
\end{equation}

Note in passing that if $G$ is any countable abelian group which is an
extension\label{LOSG0_2}%
\[
0\longrightarrow G_{0}\hooklongrightarrow G\overset{\tau}{\longrightarrow
}\mathbb{Z}\left[  a\right]  \longrightarrow0
\]
where $G_{0}$ is a torsion-free abelian group and $a$ is a real number, and
$\mathbb{Z}\left[  a\right]  $ is equipped with the order coming from
$\mathbb{Z}\left[  a\right]  \subseteq\mathbb{R}$, and if $G$ is equipped with
the order $g>0$ if and only if $\tau\left(  g\right)  >0$, then $G$ is
unperforated and has the Riesz interpolation property, so $G$ is the dimension
group of an AF-algebra by Effros--Handelman--Shen's theorem \cite{EHS80},
\cite{Dav96}.

Another way of describing $\left(  K_{0}\left(  \mathfrak{A}_{L}\right)
,K_{0}\left(  \mathfrak{A}_{L}\right)  _{+},\left[  \openone\right]  \right)
$ which will be quite useful in the sequel is the following: Let $p_{L}\left(
x\right)  $ be $\left|  \det J\right|  $ times the characteristic polynomial
of $J^{-1}$, see (\ref{eqSubNew.29}), (\ref{eqRemNewBis.16}),
(\ref{eqRemNewBis.41}), and let $a=e^{-\beta}$ be the positive real root of
this polynomial (i.e., $1/a$ is the Perron--Frobenius eigenvalue of $J$).
Then\label{LOSK0AL_4}%
\begin{equation}
K_{0}\left(  \mathfrak{A}_{L}\right)  \cong\mathbb{Z}\left[  x\right]
\diagup\left(  p_{L}\right)  \label{eqRemFeb.X}%
\end{equation}
as additive groups, and the order on $K_{0}\left(  \mathfrak{A}_{L}\right)  $
is given by that $p+\mathbb{Z}\left[  x\right]  p_{L}\left(  x\right)  >0$ if
and only if%
\begin{equation}
p\left(  a\right)  >0 \label{eqRemFeb.Y}%
\end{equation}
(this condition is well defined since $p_{L}\left(  a\right)  =0$). The
element $\left[  \openone\right]  $ corresponds to $1+\mathbb{Z}\left[
x\right]  p_{L}\left(  x\right)  $ by this isomorphism. Application of
$J^{-1}$ on $K_{0}\left(  \mathfrak{A}_{L}\right)  $ (which is well defined by
(\ref{eqRemNewBis.19})) corresponds to multiplication by $x$, i.e.,
\begin{equation}
J^{-1}\left(  p\left(  x\right)  +\mathbb{Z}\left[  x\right]  p_{L}\left(
x\right)  \right)  =xp\left(  x\right)  +\mathbb{Z}\left[  x\right]
p_{L}\left(  x\right)  \label{eqRemFeb.Z}%
\end{equation}
where the left-hand polynomial is identified with its representative in
$K_{0}\left(  \mathfrak{A}_{L}\right)  $ given as in (\ref{eqRemFeb.uther}),
below. The isomorphism between the concrete realization of $K_{0}\left(
\mathfrak
{A}_{L}\right)  $ in (\ref{eqRemNewBis.19}) and $\mathbb{Z}\left[  x\right]
\diagup\left(  p_{L}\right)  $ is thus given by%
\begin{equation}
\left(  a_{0},\dots,a_{N-1}\right)  \longmapsto a_{0}+a_{1}x+\dots
+a_{N-1}x^{N-1}\mod{p_{L}\left( x\right) }, \label{eqRemFeb.uther}%
\end{equation}
and using this and (\ref{eqRemNewBis.21}) the statements above follow
immediately. Note also that in this picture\label{LOSkertau_2}%
\begin{equation}
\ker\tau=\mathbb{Z}\left[  x\right]  p_{a}\left(  x\right)  \diagup\left(
p_{L}\left(  x\right)  \right)  , \label{eqRemFeb.okseskaft}%
\end{equation}
where $p_{a}\in\mathbb{Z}\left[  x\right]  $ is the minimal polynomial of $a$,
which is a factor of $p_{L}$. Factorizing\label{LOSpLx_1}%
\begin{equation}
p_{L}\left(  x\right)  =p_{0}\left(  x\right)  p_{a}\left(  x\right)
\label{eqRemFeb.angstrom}%
\end{equation}
we thus have%
\begin{equation}
\ker\tau\cong\mathbb{Z}\left[  x\right]  \diagup\left(  p_{0}\left(  x\right)
\right)  . \label{eqRemFeb.a}%
\end{equation}
This viewpoint will be important in Chapter \ref{Pediferient} and later chapters.

One connection between the cone (\ref{eqRemNew.0}) and that of
(\ref{eqRemNewBis.34}) can be made by the use of \cite[Lemma 2]{FrSo92}, which
shows that a given element $g$ of $K_{0}\left(  \mathfrak{A}_{L}\right)
=\bigcup_{k\geq0}J_{L}^{-k}\mathbb{Z}_{{}}^{N}$ satisfies $\tau\left(
g\right)  >0$ if and only if there are $k\in\left\{  0,1,2,\dots\right\}  $,
$n=\left(  n_{1},\dots,n_{N}\right)  \in\mathbb{Z}^{N}$, such that $n_{i}>0$,
and $v\in\ker\left(  \tau\right)  $ such that%
\[
g=v+J_{L}^{-k}n.
\]

In applications, this ``concrete'' realization of $\left(  K_{0}\left(
\mathfrak{A}_{L}\right)  ,K_{0}\left(  \mathfrak{A}_{L}\right)  _{+},\left[
\openone\right]  \right)  $\label{LOSDAL_4} is often nevertheless not concrete
enough to decide isomorphism and non-isomorphism of the algebras $\mathfrak
{A}_{L}$, but there is a simple sufficient condition for isomorphism, namely equality:

\begin{corollary}
\label{CorRemNewBis.1}Let $1\leq L_{1}\leq\dots\leq L_{d}$ and $1\leq
L_{1}^{\prime}\leq\dots\leq L_{d^{\prime}}^{\prime}$ be two sets of integers,
each with greatest common divisor $1$. Assume%
\begin{align}
&  \text{the unique solutions }x,y\in\left(  0,1\right)  \text{ of the
equations }\label{eqRemNewBis.35}\\
&  1-\sum_{i}^{d}x^{L_{i}}=0\text{ and }1-\sum_{i}^{d^{\prime}}y^{L_{i}%
^{\prime}}=0\text{ are the same, i.e., }x=y;\nonumber
\end{align}
and%
\begin{equation}
L_{d}=L_{d^{\prime}}^{\prime}\text{\textup{(}}=N\text{\textup{)}\qquad
and\qquad}\bigcup_{n=0}^{\infty}J_{L}^{-n}\left(  \mathbb{Z}^{N}\right)
=\bigcup_{n=0}^{\infty}J_{L^{\prime}}^{-n}\left(  \mathbb{Z}^{N}\right)  .
\label{eqRemNewBis.36}%
\end{equation}
It follows that $\mathfrak{A}_{L}$ and $\mathfrak{A}_{L^{\prime}}$ are
isomorphic $C^{\ast}$-algebras.
\end{corollary}

\begin{proof}
By condition (\ref{eqRemNewBis.35}) and (\ref{eqRemNewBis.16}%
)--(\ref{eqRemNewBis.18}) the Perron--Frobenius eigenvalue $e^{\beta}$ and the
normalized left Perron--Frobenius eigenvector $\alpha$%
\label{LOSPerronFrobeniuslefteigenvector_3} are the same for $J_{L}$ and
$J_{L^{\prime}}$. But (\ref{eqRemNewBis.36}) states that $K_{0}\left(
\mathfrak{A}_{L}\right)  $ and $K_{0}\left(  \mathfrak
{A}_{L^{\prime}}\right)  $ are the same subgroup of $\mathbb{Q}^{N}$, and by
(\ref{eqRemNewBis.34}) the positive cones are the same. By
(\ref{eqRemNewBis.33}), $\left[  \openone\right]  $ is represented by the same
element of the two cases, and thus the complete invariants
(\ref{eqRemNewBis.30}) are the same. Thus $\mathfrak{A}_{L}$ and $\mathfrak
{A}_{L^{\prime}}$ are isomorphic $C^{\ast}$-algebras.
\end{proof}

Still we will see in the examples that the computation of $\bigcup
_{n=0}^{\infty}J_{L}^{-n}\left(  \mathbb{Z}^{N}\right)  $ is not so simple in
general. But there is one simple special case, namely when $m_{k}=1$ in
(\ref{eqRemNewBis.15}), i.e., $\left|  \det\left(  J_{L}\right)  \right|  =1$.
Then $J_{L}^{-1}$ is a matrix with integer entries, so $J_{L}\colon
\mathbb{Z}^{N}\rightarrow\mathbb{Z}^{N}$ is bijective and hence%
\begin{equation}
K_{0}\left(  \mathfrak{A}_{L}\right)  =\mathbb{Z}^{N} \label{eqRemNewBis.37}%
\end{equation}
by (\ref{eqRemNewBis.19}). It follows immediately from Corollary
\ref{CorRemNewBis.1} that

\begin{corollary}
\label{CorRemNewBis.2}Let $1\leq L_{1}\leq\dots\leq L_{d}$ and $1\leq
L_{1}^{\prime}\leq\dots\leq L_{d^{\prime}}^{\prime}$ be two sets of integers,
each with greatest common divisor $1$. Assume%
\begin{align}
&  \text{the unique solutions }x,y\in\left(  0,1\right)  \text{ of the
equations }\label{eqRemNewBis.38}\\
&  1-\sum_{i}^{d}x^{L_{i}}=0\text{ and }1-\sum_{i}^{d^{\prime}}y^{L_{i}%
^{\prime}}=0\text{ are the same, i.e., }x=y;\nonumber
\end{align}
and%
\begin{multline}
L_{d}=L_{d^{\prime}}^{\prime}\text{ and }L_{d-1}<L_{d}\text{ and }%
L_{d^{\prime}-1}^{\prime}<L_{d^{\prime}}^{\prime}\text{\textup{(}i.e., the
matrices }J_{L}\\
\text{ and }J_{L^{\prime}}\text{ have the same rank, and the lower left matrix
element is }1\text{\textup{).}} \label{eqRemNewBis.39}%
\end{multline}
It follows that $\mathfrak{A}_{L}$ and $\mathfrak{A}_{L^{\prime}}$ are
isomorphic $C^{\ast}$-algebras.
\end{corollary}

\begin{proof}
In this case $\left|  \det J_{L}\right|  =\left|  \det J_{L^{\prime}}\right|
=1$ so $K_{0}\left(  \mathfrak{A}_{L}\right)  =K_{0}\left(  \mathfrak
{A}_{L^{\prime}}\right)  =\mathbb{Z}^{N}$ and the result follows from
Corollary \ref{CorRemNewBis.1}.
\end{proof}

\begin{figure}[ptb]
\begin{picture}(360,564)
\put(0,213){\includegraphics
[bb=186 90 486 672,clip,width=180bp,height=349.2bp]{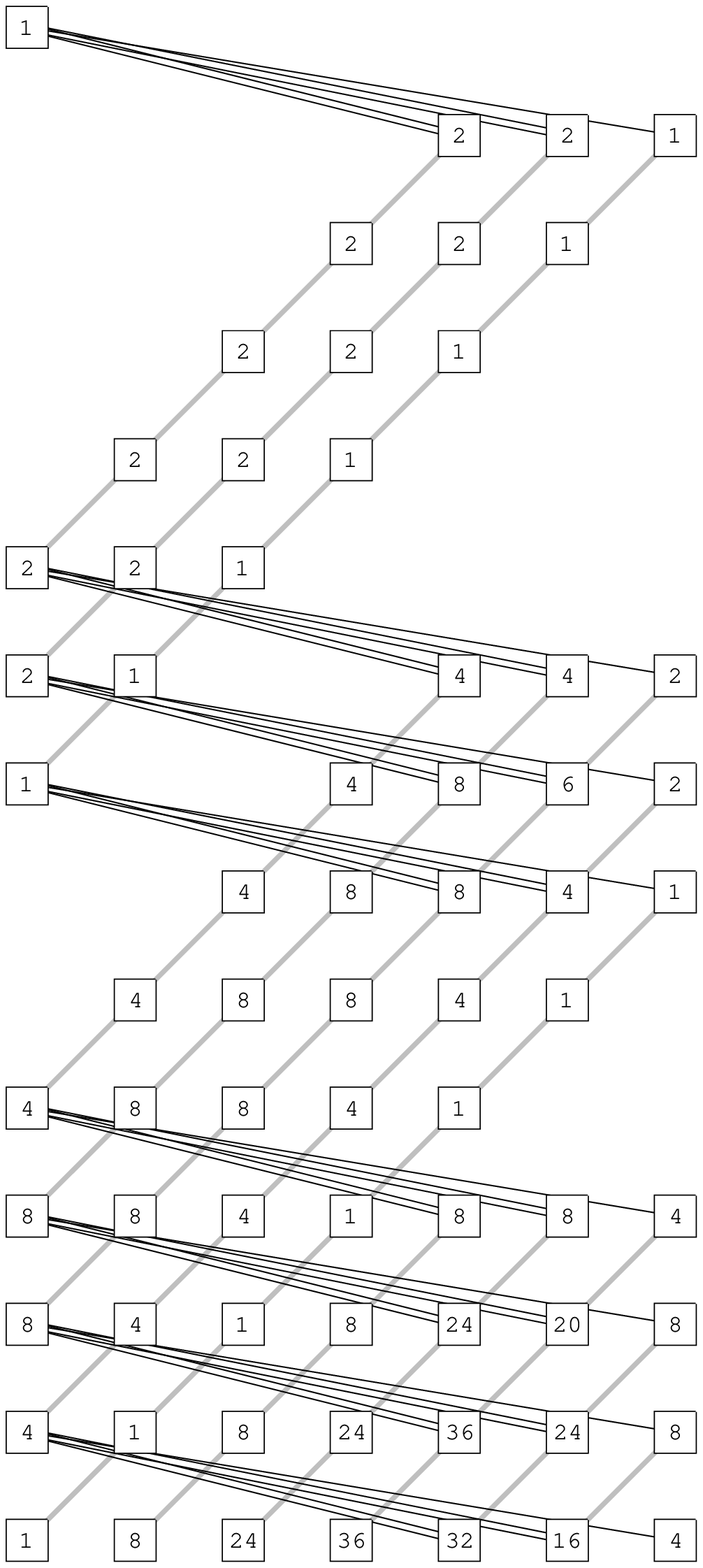}}
\put(0,0){\includegraphics
[bb=186 314 486 672,clip,width=180bp,height=214.8bp]{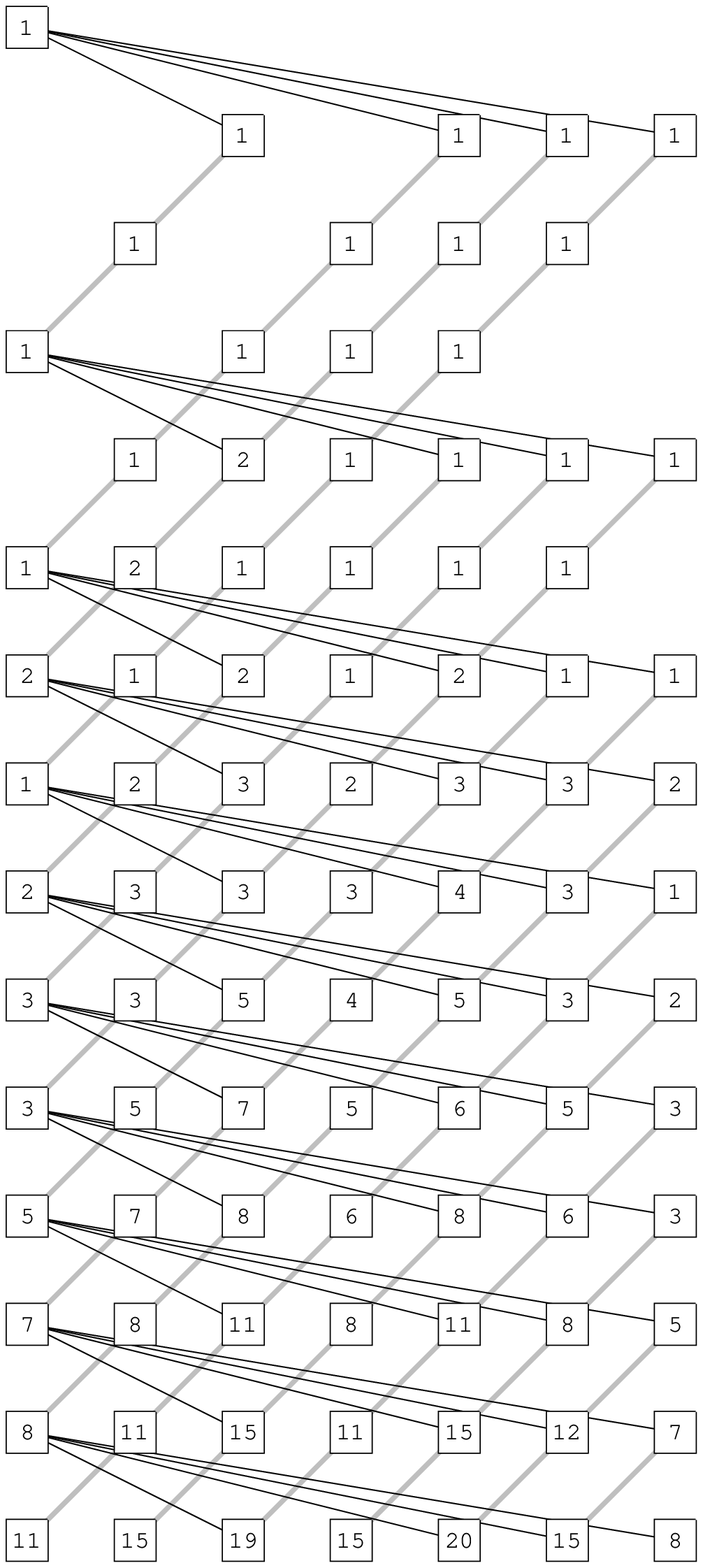}}
\put(180.8,49.2){\includegraphics
[bb=241 36 479 720,clip,width=179.2bp,height=515bp]{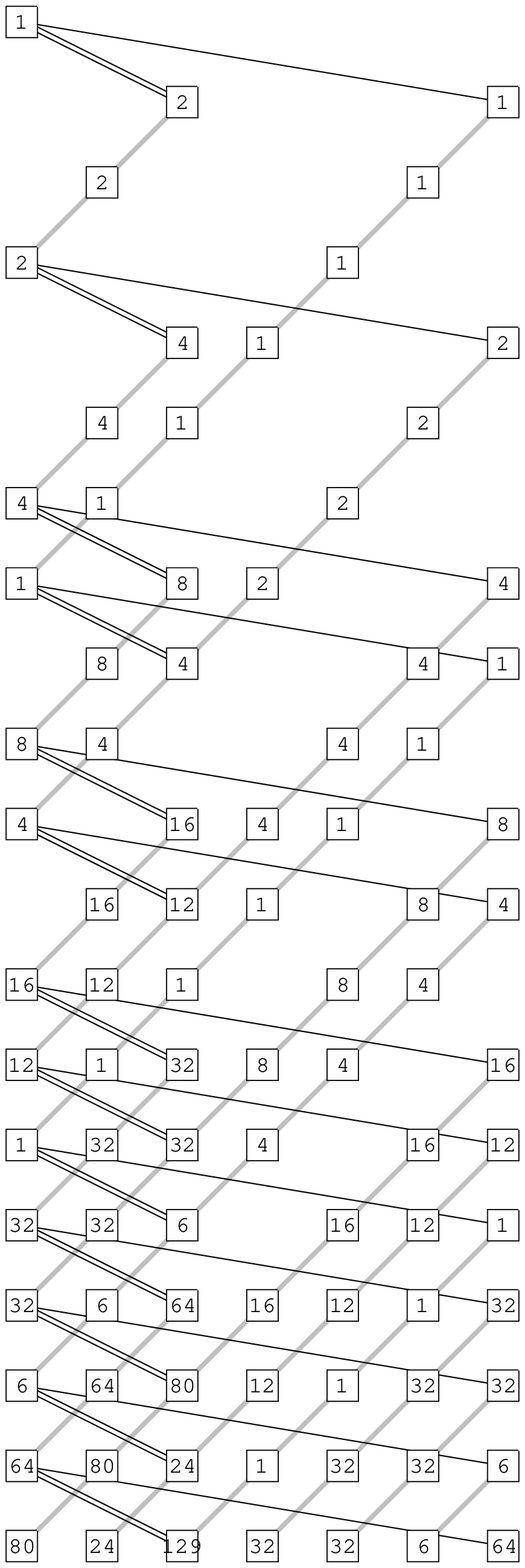}}
\end{picture}
\caption{$L=\{5,5,6,6,7\}$ (top left), $\{3,5,6,7\}$ (bottom left), and
$\{3,3,7\}$ (right), illustrating the $n=4$ case in Example
\ref{ExaRemNewBis.3}. These represent isomorphic algebras.}%
\label{BratDiagsn4}%
\end{figure}

In general we will see in the examples that the algebras $\mathfrak{A}_{L}$
for different $L$'s are ``almost never'' isomorphic. However, Corollary
\ref{CorRemNewBis.2} may be used to make some isomorphic tuples:

\begin{example}
\label{ExaRemNewBis.3}It is convenient from here and henceforth to write $J$
in the form\label{LOSJ_1}
\begin{equation}
J=%
\begin{pmatrix}
\vphantom{\vdots}m_{1} & 1 & 0 & \cdots & 0 & 0\\
\vphantom{\vdots}m_{2} & 0 & 1 & \cdots & 0 & 0\\
\vdots &  & \ddots & \ddots & \vdots & \vdots\\
\vphantom{\vdots}m_{N-2} & 0 &  & \ddots & 1 & 0\\
\vphantom{\vdots}m_{N-1} & 0 & 0 &  & 0 & 1\\
\vphantom{\vdots}m_{N} & 0 & 0 & \cdots & 0 & 0
\end{pmatrix}
, \label{eqRemNewBis.40}%
\end{equation}
instead of \textup{(\ref{eqRemNewBis.15}),} and then equation
\textup{(\ref{eqRemNewBis.16})} becomes\label{LOSpLx_2}%
\begin{equation}
p_{L}\left(  x\right)  =\sum_{j=1}^{N}m_{j}x^{j}-1=0. \label{eqRemNewBis.41}%
\end{equation}
As noted in \textup{(\ref{eqSubNew.28})--(\ref{eqSubNew.29})} this equation is
$m_{N}$ times the characteristic equation of
\begin{equation}
J^{-1}=%
\begin{pmatrix}
\vphantom{\frac{m_{j}}{m_{k}}\vdots}0 & 0 & \cdots & 0 & 0 & \frac{1}{m_{N}}\\
\vphantom{\frac{m_{j}}{m_{k}}\vdots}1 & 0 &  & 0 & 0 & -\frac{m_{1}}{m_{N}}\\
\vphantom{\frac{m_{j}}{m_{k}}\vdots}0 & 1 & \ddots &  & 0 & -\frac{m_{2}%
}{m_{N}}\\
\vphantom{\frac{m_{j}}{m_{k}}\vdots}\vdots & \vdots & \ddots & \ddots &  &
\vdots\\
\vphantom{\frac{m_{j}}{m_{k}}\vdots}0 & 0 & \cdots & 1 & 0 & -\frac{m_{N-2}%
}{m_{N}}\\
\vphantom{\frac{m_{j}}{m_{k}}\vdots}0 & 0 & \cdots & 0 & 1 & -\frac{m_{N-1}%
}{m_{N}}%
\end{pmatrix}
\label{eqThetaNyet}%
\end{equation}
The condition in Corollary \textup{\ref{CorRemNewBis.2}} is that $m_{N}=1$,
i.e., the polynomial \textup{(\ref{eqRemNewBis.41}) }should be monic. Now it
follows from Corollary \textup{\ref{CorRemNewBis.2}} that two \emph{monic}
polynomials of the form \textup{(\ref{eqRemNewBis.41})} give rise to
isomorphic algebras if they have the same degree $N$ and the root in $\left(
0,1\right)  $ is the same for the two polynomials \textup{(}under the overall
condition $\gcd\left(  \left\{  i\mid m_{i}\neq0\right\}  \right)
=1$\textup{).} \textup{(}This is no longer true if the polynomials are not
monic; see, e.g., the examples in Chapters \textup{\ref{APP.EXA}} and
\textup{\ref{ClmN}.)} To generate polynomials of the form
\textup{(\ref{eqRemNewBis.41})} with the same root, one may start with a fixed
polynomial of the required form, e.g.,%
\[
p_{0}\left(  x\right)  =x^{3}+x^{2}-1,
\]
and then multiply $p_{0}\left(  x\right)  $ with a polynomial%
\[
q\left(  x\right)  =x^{n}+k_{n-1}x^{n-1}+k_{n-2}x^{n-2}+\dots+k_{1}x+1.
\]
Choose the coefficients $k_{1},\dots,k_{n-1}$ as integers such that $m_{j}%
\geq0$ for all $j$ in%
\[
p_{0}\left(  x\right)  q\left(  x\right)  =x^{n+3}+\sum_{j=1}^{n+2}m_{j}%
x^{j}-1.
\]
This procedure, applied to $n=2,3,4,5$, gives the following values for the
possible first column%
\[%
\begin{pmatrix}
m_{1}\\
\vdots\\
m_{n+3}%
\end{pmatrix}
=%
\begin{pmatrix}
m_{1}\\
\vdots\\
m_{N}%
\end{pmatrix}
=%
\begin{pmatrix}
m_{1}\\
\vdots\\
m_{N-1}\\
1
\end{pmatrix}
\]
of the incidence matrix $J$:

\begin{description}
\item [$n=2$]Two isomorphic algebras:%
\[%
\begin{pmatrix}
1\\
0\\
0\\
0\\
1
\end{pmatrix}
,%
\begin{pmatrix}
0\\
0\\
1\\
1\\
1
\end{pmatrix}
.
\]

\item[$n=3$] Two isomorphic algebras:%
\[%
\begin{pmatrix}
0\\
1\\
0\\
0\\
1\\
1
\end{pmatrix}
,%
\begin{pmatrix}
0\\
0\\
0\\
1\\
2\\
1
\end{pmatrix}
.
\]

\item[$n=4$] Three isomorphic algebras, which are subalgebras of
$\mathcal{O}_{5}$, $\mathcal{O}_{4}$, $\mathcal{O}_{3}$, respectively:%
\[%
\begin{pmatrix}
0\\
0\\
0\\
0\\
2\\
2\\
1
\end{pmatrix}
,%
\begin{pmatrix}
0\\
0\\
1\\
0\\
1\\
1\\
1
\end{pmatrix}
,%
\begin{pmatrix}
0\\
0\\
2\\
0\\
0\\
0\\
1
\end{pmatrix}
.
\]
See Figure \textup{\ref{BratDiagsn4}.}

\item[$n=5$] There are $6+1$ possibilities to begin with,%
\[%
\begin{pmatrix}
0\\
0\\
0\\
0\\
1\\
2\\
2\\
1
\end{pmatrix}
,%
\begin{pmatrix}
0\\
0\\
0\\
1\\
1\\
1\\
1\\
1
\end{pmatrix}
,%
\begin{pmatrix}
0\\
1\\
0\\
0\\
0\\
1\\
1\\
1
\end{pmatrix}
,%
\begin{pmatrix}
0\\
0\\
1\\
1\\
0\\
0\\
1\\
1
\end{pmatrix}
,%
\begin{pmatrix}
1\\
0\\
0\\
0\\
0\\
0\\
1\\
1
\end{pmatrix}
,%
\begin{pmatrix}
0\\
0\\
0\\
2\\
1\\
0\\
0\\
1
\end{pmatrix}
,%
\begin{pmatrix}
0\\
1\\
0\\
1\\
0\\
0\\
0\\
1
\end{pmatrix}
,
\]
but in the last example $\gcd\left(  L\right)  =2$, so this falls outside our
scope. The remaining $6$ vectors give rise to isomorphic subalgebras of
$\mathcal{O}_{d}$ with $d=6,5,4,4,3,4$, respectively. Note that this shows
that $d$ is \emph{not} an invariant. The next-to-last example is illustrated
in Figure \textup{\ref{BratDiag4458}.}
\end{description}
\end{example}

\begin{remark}
\label{RemRemNewBis.4}The isomorphism of the algebras $\mathfrak{A}_{L}$ and
$\mathfrak{A}_{L^{\prime}}$ established in Example
\textup{\ref{ExaRemNewBis.3}} for various pairs $L$, $L^{\prime}$ was arrived
at in a quite roundabout way. In general it follows from \cite[Theorem
2.7]{Bra72} that $\mathfrak{A}_{L}$ and $\mathfrak{A}_{L^{\prime}}$ are stably
isomorphic if and only if there exist natural numbers $k_{1},k_{2},k_{3}%
,\dots$, $l_{1},l_{2},l_{3},\dots$, and matrices $A_{1},A_{2},\dots$,
$B_{1},B_{2},\dots$ with nonnegative integer matrix elements such that the
following diagram\label{LOSBrattelidiagrams_5} commutes:%
\begin{equation}
\setlength{\unitlength}{78bp}\begin{minipage}{1.67\unitlength} \begin{picture}%
(1.67,5)(-0.4,-4.895) \put(-0.093,-4.3565){\includegraphics
[bb=133 0 215 348,height=348bp,width=82bp]{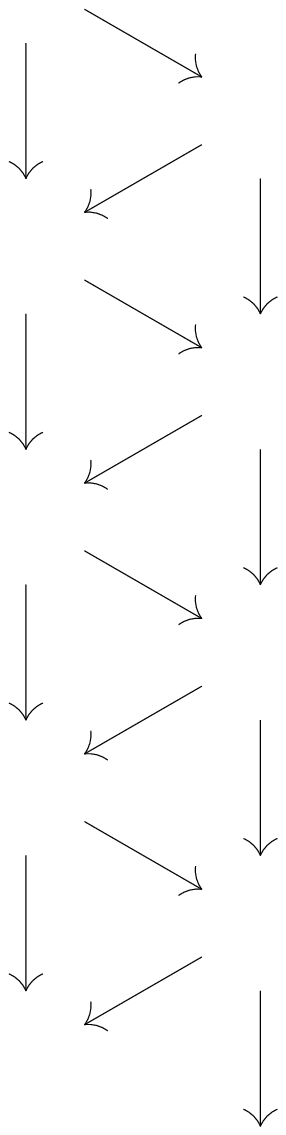}} \put(0,0){\makebox
(0,0){$\bullet$}} \put(0.866,-0.5){\makebox(0,0){$\bullet$}} \put
(0,-1){\makebox(0,0){$\bullet$}} \put(0.866,-1.5){\makebox(0,0){$\bullet$}%
}\put(0,-2){\makebox(0,0){$\bullet$}} \put(0.866,-2.5){\makebox(0,0){$\bullet
$}} \put(0,-3){\makebox(0,0){$\bullet$}} \put(0.866,-3.5){\makebox
(0,0){$\bullet$}} \put(0,-4){\makebox(0,0)[t]{$\vdots$}} \put
(0.866,-4.5){\makebox(0,0)[t]{$\vdots$}} \put(0.433,-0.225){\makebox
(0,0)[bl]{$A_{1}$}} \put(0.433,-1.225){\makebox(0,0)[bl]{$A_{2}$}}%
\put(0.433,-2.225){\makebox(0,0)[bl]{$A_{3}$}} \put(0.433,-3.225){\makebox
(0,0)[bl]{$A_{4}$}} \put(0.433,-0.725){\makebox(0,0)[br]{$B_{1}$}}%
\put(0.433,-1.725){\makebox(0,0)[br]{$B_{2}$}} \put(0.433,-2.725){\makebox
(0,0)[br]{$B_{3}$}} \put(-0.025,-0.5){\makebox(0,0)[r]{$J_{L}^{k_{1}}$}}%
\put(-0.025,-1.5){\makebox(0,0)[r]{$J_{L}^{k_{2}}$}} \put
(-0.025,-2.5){\makebox(0,0)[r]{$J_{L}^{k_{3}}$}} \put(0.891,-1){\makebox
(0,0)[l]{$J_{L^{\prime}}^{l_{1}}$}} \put(0.891,-2){\makebox
(0,0)[l]{$J_{L^{\prime}}^{l_{2}}$}} \put(0.891,-3){\makebox
(0,0)[l]{$J_{L^{\prime}}^{l_{3}}$}} \end{picture} \end{minipage}
\label{eqbeforeRemNewTer.42}%
\end{equation}
This means that%
\begin{equation}%
\begin{aligned}
J_{L}^{k_{i}} &= B_{i}A_{i},  \\
J_{L^{\prime}}^{l_{i}} &= A_{i+1}B_{i}
\end{aligned}
\label{eqRemNewTer.42}%
\end{equation}
for $i=1,2,\dots$. There are examples showing that the sequences $A$, $B$,
$k$, $l$ cannot always be taken to be constant when they exist \cite{BJKR00}.
In our case, when the $J_{L}$'s are nonsingular, the existence of constant
sequences would entail that $J_{L}$ and $J_{L^{\prime}}$ have the same
dimension, and $J_{L}^{k}$ be conjugate to $J_{L^{\prime}}^{l}$. Note in this
connection that $J_{L}$ is conjugate to $J_{L^{\prime}}$ if and only if
$L=L^{\prime}$, because the characteristic polynomial of $J_{L}$ completely
determines $L=\left(  L_{1},\dots,L_{d}\right)  $, as we have seen.

In the covariant version of this isomorphism problem, it is known from a
theorem by Krieger that the sequences can be taken to be constant. Let
$G\left(  L\right)  $ be the dimension group associated to $L$, and $\left(
\sigma_{L}\right)  _{\ast}$ the automorphism of $G\left(  L\right)  $
determined by $J_{L}$. Let now $\mathfrak{B}_{L}=\mathfrak{A}_{L}%
\otimes\mathcal{K}\left(  \ell^{2}\right)  $ be the stable AF-algebra
associated to $G\left(  L\right)  $, and $\sigma_{L}$ an automorphism of
$\mathfrak{B}_{L}$ such that the corresponding automorphism of $G\left(
L\right)  $ is $\left(  \sigma_{L}\right)  _{\ast}$. Then Krieger's theorem
\cite{Kri80} says that $\left(  G\left(  L\right)  ,\left(  \sigma_{L}\right)
_{\ast}\right)  $\label{LOSGLsigmaLstar} is isomorphic to $\left(  G\left(
L^{\prime}\right)  ,\left(  \sigma_{L^{\prime}}\right)  _{\ast}\right)  $ if
and only if there is a $k\in\mathbb{N}$ and nonnegative rectangular matrices
$A$, $B$ such that%
\begin{equation}%
\begin{aligned}
AJ_{L} &= J_{L^{\prime}}A,  \\
BJ_{L^{\prime}} &= J_{L}B,  \\
AB &= J_{L}^{k},  \\
BA &= J_{L^{\prime}}^{k}.
\end{aligned}
\label{eqRemNewTer.43}%
\end{equation}
If also $N>1$, it was proved recently in \cite{BrKi00b} that this is also
equivalent to outer conjugacy of $\sigma_{L}$ and $\sigma_{L^{\prime}}$. All
these results were proved in the more general setting of constant incidence
matrices. In the $J_{L}$ case, the conditions simply mean $L=L^{\prime}$.
\textup{In fact, the third condition in (\ref{eqRemNewTer.43})} implies that
both $A$ and $B$ are nonsingular. Hence, the first condition reads
$J_{L^{\prime}}=AJ_{L}A^{-1}$, and we conclude that $J_{L}$ and $J_{L^{\prime
}}$ have the same characteristic polynomial. Since the coefficients in the
characteristic polynomial of $J_{L}$ are the numbers in the first column of
$J_{L}$, it follows that $J_{L}=J_{L^{\prime}}$ as claimed. \textup{(}See also
\textup{(\ref{eqExt.lambdam})--(\ref{eqCycNew.1})} for more details.\textup{)}

Note that the Williams conjecture discussed at the end of Chapter 6 in
\cite{Eff81} has been settled in the negative in \cite{KiRo99}.
\end{remark}

Let us end this chapter by mentioning another corollary of results in
\cite{BrKi00b}, which classifies the actions $\sigma^{\left(  L\right)  }$ of
$\mathbb{T}$ on $\mathcal{O}_{d}$ defined by (\ref{eqSub.1}):

\begin{corollary}
\label{CorRemFeb.5}Let $1\leq L_{1}\leq\dots\leq L_{d}$ and $1\leq
L_{1}^{\prime}\leq\dots\leq L_{d^{\prime}}^{\prime}$ be two sets of integers,
each with greatest common divisor $1$. The following conditions are equivalent.

\begin{enumerate}
\item \label{CorRemFeb.5(1)}The automorphism $\sigma_{L}$ of $\mathfrak{A}%
_{L}\otimes\mathcal{K}\left(  \ell^{2}\right)  $ defined prior to
\textup{(\ref{eqRemNewTer.43})} is outer conjugate to $\sigma_{L^{\prime}}$.

\item \label{CorRemFeb.5(2)}$\left(  G\left(  L\right)  ,\left(  \sigma
_{L}\right)  _{\ast}\right)  $ is isomorphic to $\left(  G\left(  L^{\prime
}\right)  ,\left(  \sigma_{L^{\prime}}\right)  _{\ast}\right)  $.

\item \label{CorRemFeb.5(3)}The action $\sigma^{\left(  L\right)  }$ of
$\mathbb{T}$ on $\mathcal{O}_{d}$ defined by \textup{(\ref{eqSub.1})} is outer
conjugate to the action $\sigma^{\left(  L^{\prime}\right)  }$.

\item \label{CorRemFeb.5(4)}$\sigma^{\left(  L\right)  }$ and $\sigma^{\left(
L^{\prime}\right)  }$ are conjugate actions of $\mathbb{T}$.

\item \label{CorRemFeb.5(5)}$L=L^{\prime}$.
\end{enumerate}
\end{corollary}

\begin{proof}
We already noted above that (\ref{CorRemFeb.5(1)}) $\Leftrightarrow$
(\ref{CorRemFeb.5(2)}) is \cite[Corollary 1.5]{BrKi00b}. But
(\ref{CorRemFeb.5(2)}) $\Leftrightarrow$ (\ref{CorRemFeb.5(4)}) follows from
\cite[Corollary 4.1]{BrKi00b}. The implication (\ref{CorRemFeb.5(3)})
$\Rightarrow$ (\ref{CorRemFeb.5(2)}) follows by noting that the stabilization
of the dual actions of $\sigma^{\left(  L\right)  }$, $\sigma^{\left(
L^{\prime}\right)  }$ is outer conjugate to $\sigma_{L}$, $\sigma_{L^{\prime}%
}$ by Takai duality. The only remaining nontrivial implication is
(\ref{CorRemFeb.5(2)}) $\Rightarrow$ (\ref{CorRemFeb.5(5)}); as noted in
\cite{BJKR00}, the relations (\ref{eqRemNewTer.43}) imply that $J_{L}$ and
$J_{L^{\prime}}$ are similar, and thus have the same characteristic
polynomial. But by (\ref{eqExt.10}), the characteristic polynomial determines
$J_{L}$ and thus $L$ uniquely. Thus (\ref{CorRemFeb.5(2)}) $\Rightarrow$
(\ref{CorRemFeb.5(5)}).
\end{proof}

\begin{remark}
\label{RemRemJun.6}The equivalence relation of Bratteli diagrams referred to
in the second paragraph in the Introduction can be described as follows: The
diagram itself can be described as a sequence of incidence matrices%
\begin{equation}
J_{1},J_{2},J_{3},J_{4},\dots. \label{eqRemJun.53}%
\end{equation}
These are \textup{(}not necessarily square\/\textup{)} matrices with integer
nonnegative matrix units such that the number of columns in $J_{n+1}$ is equal
to the number of rows in $J_{n}$. One way of obtaining an equivalent diagram
is then to remove rows from the diagram and connect the remaining vertices by
edges with multiplicity given by the number of ways one can go from the upper
vertex to the lower along the original diagram. In terms of incidence
matrices, one picks an increasing sequence $1\leq n_{1}<n_{2}<n_{3}$ of
integers, and replaces the sequence \textup{(\ref{eqRemJun.53})} by%
\begin{equation}
J_{n_{2}-1}J_{n_{2}-2}\cdots J_{n_{1}},\qquad J_{n_{3}-1}J_{n_{3}-2}\cdots
J_{n_{2}},\qquad\dots. \label{eqRemJun.54}%
\end{equation}
The equivalence relation is then simply the equivalence relation on sequences
of incidence matrices generated by this relation. One has to apply the
relation or its inverse four times to go from one diagram to another. Roughly,
start from%
\[
\mathfrak{A}_{1}\longrightarrow\mathfrak{A}_{2}\longrightarrow\mathfrak{A}%
_{3}\longrightarrow\cdots
\]
by removing rows to obtain%
\[
\mathfrak{A}_{n_{1}}\longrightarrow\mathfrak{A}_{n_{2}}\longrightarrow
\mathfrak{A}_{n_{3}}\longrightarrow\cdots,
\]
then insert new rows to obtain%
\[
\mathfrak{A}_{n_{1}}\longrightarrow\mathfrak{B}_{m_{1}}\longrightarrow
\mathfrak{A}_{n_{2}}\longrightarrow\mathfrak{B}_{m_{2}}\longrightarrow\cdots,
\]
then remove rows to obtain%
\[
\mathfrak{B}_{m_{1}}\longrightarrow\mathfrak{B}_{m_{2}}\longrightarrow
\mathfrak{B}_{m_{3}}\longrightarrow\cdots,
\]
and finally insert rows to obtain%
\[
\mathfrak{B}_{1}\longrightarrow\mathfrak{B}_{2}\longrightarrow\mathfrak{B}%
_{3}\longrightarrow\cdots.
\]
One example from \cite{Bra72}, where the first and last steps are unnecessary,
is%
\[
\setlength{\unitlength}{30pt}\begin{picture}(12,11)
\ttfamily
\setcounter{glimmer}{1}
\multiput(0.75,8.75)(0,-2){5}{\framebox(0.5,0.5){\arabic{glimmer}%
}\addtocounter{glimmer}{\value{glimmer}}}
\multiput(1,8.75)(0,-2){4}{\line(0,-1){1.5}}
\put(1,0.75){\line(0,-1){0.5}}
\multiput(1.25,8.75)(0,-2){4}{\line(1,-1){1.5}}
\put(1.25,0.75){\line(1,-1){0.5}}
\setcounter{glimmer}{1}
\multiput(2.75,8.75)(0,-2){5}{\framebox(0.5,0.5){\arabic{glimmer}%
}\addtocounter{glimmer}{\value{glimmer}}}
\multiput(3,8.75)(0,-2){4}{\line(0,-1){1.5}}
\put(3,0.75){\line(0,-1){0.5}}
\multiput(2.75,8.75)(0,-2){4}{\line(-1,-1){1.5}}
\put(2.75,0.75){\line(-1,-1){0.5}}
\multiput(3.33333,9)(0.16667,0){15}{\makebox(0,0){$\cdot$}}
\put(3.5,5){\vector(1,0){2}}
\put(4.5,5.25){\makebox(0,0){\rmfamily Insert rows}}
\multiput(3.33333,1)(0.16667,0){15}{\makebox(0,0){$\cdot$}}
\setcounter{glimmer}{1}
\multiput(5.75,8.75)(0,-2){5}{\framebox(0.5,0.5){\arabic{glimmer}%
}\addtocounter{glimmer}{\value{glimmer}}}
\multiput(6.25,8.75)(0,-2){5}{\line(1,-1){0.5}}
\setcounter{glimmer}{1}
\multiput(6.75,9.75)(0,-2){5}{\framebox(0.5,0.5){\arabic{glimmer}%
}\addtocounter{glimmer}{\value{glimmer}}}
\multiput(6.75,9.75)(0,-2){5}{\line(-1,-1){0.5}}
\multiput(7.25,9.75)(0,-2){5}{\line(1,-1){0.5}}
\setcounter{glimmer}{1}
\multiput(7.75,8.75)(0,-2){5}{\framebox(0.5,0.5){\arabic{glimmer}%
}\addtocounter{glimmer}{\value{glimmer}}}
\multiput(7.75,8.75)(0,-2){5}{\line(-1,-1){0.5}}
\multiput(7.33333,10)(0.16667,0){21}{\makebox(0,0){$\cdot$}}
\put(8.5,5){\vector(1,0){2}}
\put(9.5,5.25){\makebox(0,0){\rmfamily Remove rows}}
\multiput(7.33333,2)(0.16667,0){21}{\makebox(0,0){$\cdot$}}
\setcounter{glimmer}{1}
\multiput(10.75,9.75)(0,-2){5}{\framebox(0.5,0.5){\arabic{glimmer}%
}\addtocounter{glimmer}{\value{glimmer}}}
\multiput(10.95,9.75)(0,-2){5}{\line(0,-1){1.5}}
\multiput(11.05,9.75)(0,-2){5}{\line(0,-1){1.5}}
\end{picture}
\]
Here the algebra is the UHF algebra $\bigotimes^{\infty}M_{2}$ of Glimm type
$2^{\infty}$, also illustrated in Figure \textup{\ref{BratDiag11}.} The
algebra to the left is the fixed-point subalgebra of $\bigotimes^{\infty}%
M_{2}$ under the infinite-product action $\bigotimes^{\infty}%
\operatorname*{Ad}\left(
\begin{smallmatrix}
1 & 0\\
0 & -1
\end{smallmatrix}
\right)  $ of $\mathbb{Z}_{2}$, and the figure shows that this fixed-point
algebra is isomorphic to the full algebra \cite{Sto70}. To show directly that
the pairs or triples of diagrams shown in Figures \textup{\ref{BratDiagsn4},
\ref{BratDiagsn2}, \ref{BratDiagsn2n4}, \ref{BratDiagsN4g2s1},} and
\textup{\ref{BratDiagN4samlam}} can be transformed into each other by this
method is presumably a much harder task, as it is to show directly that the
pair in Figure \textup{\ref{BratDiagsd6}} cannot be transformed into each other.
\end{remark}

\setcounter{figurelink}{\value{figure}}

\chapter{\label{Brunt}Invariants related to the Perron--Frobenius eigenvalue}

Let $J$, $K$ be two nonsingular $N\times N$ matrices with nonnegative matrix
elements which are primitive, i.e., for sufficiently large $n\in\mathbb{N}$,
$J^{n}$ and $K^{n}$ have only strictly positive matrix elements. Let
$\lambda_{1}$, $\lambda_{2}$ be the Perron--Frobenius eigenvalues of $J$, $K$.
Then $\lambda_{1}$, $\lambda_{2}$ are algebraic numbers, and $\mathbb{Q}%
\left[  \lambda_{1}\right]  $ and $\mathbb{Q}\left[  \lambda_{2}\right]  $ are
fields which are finite extensions of $\mathbb{Q}$. If $\lambda_{1}$ and
$\lambda_{2}$ are rational, they are integers since they satisfy a monic
equation. If in addition $N=1$, then the stable $C^{\ast}$-algebras associated
with the corresponding dimension groups characterized in (\ref{eqRemNewBis.1}%
)--(\ref{eqRemNewBis.34}) are $M_{\lambda_{i}^{\infty}}\otimes\mathcal{K}%
\left(  \mathcal{H}\right)  $, where $M_{\lambda^{\infty}}$ is the UHF algebra
of Glimm type $\lambda^{\infty}$ and $\mathcal{K}\left(  \mathcal{H}\right)  $
is the compact operators on a separable Hilbert space $\mathcal{H}$. It
follows from Glimm's theorem \cite{Gli60} that these algebras are isomorphic
if and only if $\lambda_{1}$ and $\lambda_{2}$ contain the same prime factors.
In particular, if $J=\left(  6\right)  $ and $K=\left(  12\right)  $ (as
$1\times1$ matrices), the associated $C^{\ast}$-algebras are isomorphic. See
also \cite[Example 9]{BJKR00}. This was partly generalized in
\cite[Theorem 10]{BJKR00}, where it was proved that if $J$, $K$ are
nonsingular primitive $N\times N$ matrices and the stable $C^{\ast}$-algebras
they define are isomorphic, then $\mathbb{Q}\left[  \lambda_{1}\right]
=\mathbb{Q}\left[  \lambda_{2}\right]  $ and $\lambda_{1}$, $\lambda_{2}$ are
products of the same primes over this field (i.e., primes in the subring
generated by the algebraic integers in the field). The example mentioned above
shows that $\lambda$ itself is not an invariant, and the purpose of this
chapter is to show that $\lambda$ itself is not an invariant in more
interesting examples of matrices of type (\ref{eqRemNewBis.40}),\label{LOSJ_2}%
\begin{equation}
J=%
\begin{pmatrix}
\vphantom{\vdots}m_{1} & 1 & 0 & \cdots & 0 & 0\\
\vphantom{\vdots}m_{2} & 0 & 1 & \cdots & 0 & 0\\
\vdots &  & \ddots & \ddots & \vdots & \vdots\\
\vphantom{\vdots}m_{N-2} & 0 &  & \ddots & 1 & 0\\
\vphantom{\vdots}m_{N-1} & 0 & 0 &  & 0 & 1\\
\vphantom{\vdots}m_{N} & 0 & 0 & \cdots & 0 & 0
\end{pmatrix}
, \label{eqBrunt.1}%
\end{equation}
where the $m_{i}$ are nonnegative integers, $m_{N}\neq0$ and $\gcd\left\{
i\mid m_{i}\neq0\right\}  =1$. The characteristic polynomial of $J$
is\label{LOSdett1J}
\begin{equation}
\det\left(  t\openone-J\right)  =t^{N}-m_{1}t^{N-1}-m_{2}t^{N-2}-\dots
-m_{N-1}t-m_{N} \label{eqBrunt.2}%
\end{equation}
and the Perron--Frobenius eigenvalue\label{LOSPerronFrobeniuseigenvalue_4}
$\lambda$ is the unique positive solution of this equation.

More examples of this kind where the $J$'s are $2\times2$ matrices can be
constructed by a machine developed in Chapter \ref{App.Cla}; see in particular
Example \ref{ExaApp.ClaNew.4} and remarks prior to Proposition
\ref{ProApp.ClaNew3}.

The example we shall give here is a modification of another example in
\cite[Example 9]{BJKR00}. For $a=2,3,4,\dots$, consider the monic polynomial%
\begin{equation}
p_{a}\left(  t\right)  =\left(  t-a^{2}\right)  \left(  t^{2}+at+a^{2}\right)
=t^{3}+\left(  -a^{2}+a\right)  t^{2}+\left(  -a^{3}+a^{2}\right)  t-a^{4}.
\label{eqBrunt.3}%
\end{equation}
The last three coefficients are negative for $a=2,3,\dots$, so this is the
characteristic polynomial of%
\begin{equation}
J=%
\begin{pmatrix}
a^{2}-a & 1 & 0\\
a^{3}-a^{2} & 0 & 1\\
a^{4} & 0 & 0
\end{pmatrix}
. \label{eqBrunt.4}%
\end{equation}
The spectrum of $J_{a}$ consists of the roots%
\begin{equation}
\operatorname*{sp}\left(  J_{a}\right)  =\left\{  a^{2},\left(  -\frac{1}%
{2}+\frac{i}{2}\sqrt{3}\right)  a,\left(  -\frac{1}{2}-\frac{i}{2}\sqrt
{3}\right)  a\right\}  \label{eqBrunt.5}%
\end{equation}
and hence we observe%
\begin{equation}
\operatorname*{sp}\left(  J_{a^{2}}\right)  =\left\{  \lambda^{2}\mid
\lambda\in\operatorname*{sp}\left(  J_{a}\right)  \right\}  .
\label{eqBrunt.6}%
\end{equation}
Thus $J_{a^{2}}$ and $J_{a}^{2}$ are conjugate over $\mathbb{Q}\left[
\sqrt{3}\right]  $, and hence over $\mathbb{Q}$. Now put%
\begin{equation}
K=J_{2}=%
\begin{pmatrix}
2 & 1 & 0\\
4 & 0 & 1\\
16 & 0 & 0
\end{pmatrix}
,\qquad J=J_{4}=%
\begin{pmatrix}
12 & 1 & 0\\
48 & 0 & 1\\
256 & 0 & 0
\end{pmatrix}
. \label{eqBrunt.7}%
\end{equation}
Then we compute that
\begin{equation}
JT=TK^{2} \label{eqBrunt.8}%
\end{equation}
for%
\begin{equation}
T=%
\begin{pmatrix}
1 & 0 & 0\\
-4 & 2 & 1\\
0 & 16 & -4
\end{pmatrix}
. \label{eqBrunt.9}%
\end{equation}
Let%
\begin{equation}
S=T^{-1}=%
\begin{pmatrix}
\vphantom{\vdots}1 & 0 & 0\\
\vphantom{\vdots}\frac{2}{3} & \frac{1}{6} & \frac{1}{24}\\
\vphantom{\vdots}\frac{8}{3} & \frac{2}{3} & -\frac{1}{12}%
\end{pmatrix}
. \label{eqBrunt.10}%
\end{equation}
It follows from (\ref{eqBrunt.8}) that%
\begin{equation}
SJ=K^{2}S. \label{eqBrunt.11}%
\end{equation}
For a given $n\in\mathbb{N}$, put%
\begin{equation}
A=K^{2n}S=SJ^{n},\qquad B=T. \label{eqBrunt.12}%
\end{equation}
It follows from (\ref{eqBrunt.8}), (\ref{eqBrunt.11}), and $ST=TS=\openone$
that%
\begin{equation}
J^{n}=TSJ^{n}=BA,\qquad K^{2n}=K^{2n}ST=AB. \label{eqBrunt.13}%
\end{equation}
This is a version of (\ref{eqRemNewTer.42}) except that $A$, $B$ are not
necessarily matrices with positive integer matrix elements, only rational
elements. To remedy this we now replace $J$, $K$ by scaled versions%
\begin{equation}
K_{d}=%
\begin{pmatrix}
2d & 1 & 0\\
4d^{2} & 0 & 1\\
16d^{3} & 0 & 0
\end{pmatrix}
,\qquad J_{d}=%
\begin{pmatrix}
12d^{2} & 1 & 0\\
48d^{4} & 0 & 1\\
256d^{6} & 0 & 0
\end{pmatrix}
, \label{eqBrunt.14}%
\end{equation}
where $d$ is an integer. One now checks that the eigenvalues of both
$K_{d}^{2}$ and $J_{d}$ are%
\begin{equation}
16d^{2},\qquad4\left(  -\frac{1}{2}\pm\frac{i}{2}\sqrt{3}\right)  d^{2},
\label{eqBrunt.15}%
\end{equation}
and then%
\begin{equation}
J_{d}T_{d}=T_{d}^{{}}K_{d}^{2},\qquad S_{d}J_{d}=K_{d}^{2}S_{d}^{{}},
\label{eqBrunt.16}%
\end{equation}
with%
\begin{equation}
T_{d}=%
\begin{pmatrix}
1 & 0 & 0\\
-4d^{2} & 2d & 1\\
0 & 16d^{3} & -4d^{2}%
\end{pmatrix}
,\qquad S_{d}=T_{d}^{-1}=%
\begin{pmatrix}
\vphantom{\vdots}1 & 0 & 0\\
\vphantom{\vdots}\frac{2}{3}d & \frac{1}{6d} & \frac{1}{24d^{3\mathstrut}}\\
\vphantom{\vdots}\frac{8}{3}d^{2} & \frac{2}{3} & -\frac{1}{12d^{2\mathstrut}}%
\end{pmatrix}
. \label{eqBrunt.17}%
\end{equation}
With this change, we note that $K_{d}^{n}$ is a multiple of an arbitrary large
power of $d$ with an integer matrix provided $n$ is large enough. Taking $n=4$
we compute%
\begin{equation}
K_{d}^{4}S_{d}^{{}}=%
\begin{pmatrix}
\vphantom{\vdots}192d^{4} & 12d^{2} & 1\\
\vphantom{\vdots}\frac{1184}{3}d^{5} & \frac{80}{3}d^{3} & \frac{2}{3}d\\
\vphantom{\vdots}\frac{2432}{3}d^{6} & \frac{128}{3}d^{4} & \frac{8}{3}d^{2}%
\end{pmatrix}
. \label{eqBrunt.18}%
\end{equation}
Choosing $d=3$ we see that%
\begin{equation}
A=K_{3}^{4}S_{3}^{{}}=S_{3}^{{}}J_{3}^{2} \label{eqBrunt.19}%
\end{equation}
is a positive integer matrix. Similarly%
\begin{equation}
T_{d}^{{}}K_{d}^{4}=%
\begin{pmatrix}
144d^{4} & 40d^{3} & 8d^{2}\\
640d^{6} & 96d^{5} & 48d^{4}\\
2048d^{8} & 512d^{7} & 256d^{6}%
\end{pmatrix}
\label{eqBrunt.20}%
\end{equation}
is a positive integer matrix whatever integer value $d$ has, and we put%
\begin{equation}
B=T_{3}^{{}}K_{3}^{4}=J_{3}^{2}T_{3}^{{}}. \label{eqBrunt.21}%
\end{equation}
Now, redefining%
\begin{equation}
K:=K_{3}=%
\begin{pmatrix}
6 & 1 & 0\\
36 & 0 & 1\\
432 & 0 & 0
\end{pmatrix}
,\qquad J:=J_{3}=%
\begin{pmatrix}
108 & 1 & 0\\
3888 & 0 & 1\\
186624 & 0 & 0
\end{pmatrix}
, \label{eqBrunt.22}%
\end{equation}
it follows from (\ref{eqBrunt.16}), (\ref{eqBrunt.17}), (\ref{eqBrunt.19}),
and (\ref{eqBrunt.21}) that
\begin{equation}%
\begin{aligned}
AJ &= K^{2}A,  \\
JB &= BK^{2},  \\
J^{4} &= BA ,  \\
K^{8} &= AB.
\end{aligned}
\label{eqBrunt.23}%
\end{equation}
Thus, $J$ and $K^{2}$ are shift equivalent in the sense of
(\ref{eqRemNewTer.43}), and in particular, $J$ and $K$ define isomorphic
AF-algebras by (\ref{eqRemNewTer.42}). But the Perron--Frobenius eigenvalues
of these matrices are $16d^{2}=12^{2}=144$ and $4d=12$, respectively. Hence
this eigenvalue in itself is not an isomorphism invariant.

\chapter{\label{Pediferient}The invariants $N$, $D$, $\operatorname*{Prim}%
\left(  m_{N}\right)  $, $\operatorname*{Prim}\left(  R_{D}\right)  $,
$\operatorname*{Prim}\left(  Q_{N-D}\right)  $}

In this chapter, we establish a triangular representation $J_{L}=\left(
\begin{tabular}
[c]{c|c}%
$J_{0_{\mathstrut}}$ & $Q$\\\hline
$0$ & $J_{D}$%
\end{tabular}
\right)  $\label{LOSJL_1} of a matrix $J_{L}$ in the standard form
(\ref{eqCyc.2}) such that the submatrices $J_{0}$ and $J_{D}$ are again in the
same standard form (with the exception that the integers corresponding to
$m_{1},\dots,m_{N}$ are no longer necessarily positive), and $\ker\left(
\tau\right)  $ is obtained from $J_{0}$ the same way $K_{0}\left(
\mathfrak{A}_{L}\right)  $ is obtained from $J_{L}$. We then use this for the
derivation of numerical $C^{\ast}$-isomorphism invariants.

Proposition \ref{ProPediferient.1}, Corollary \ref{CorCyc.4half}, and Theorem
\ref{ThmCyc.7} below account for the terms $\mathbb{Z}\left[  \frac{1}%
{k}\right]  $ (where $k\in\mathbb{Z}$, $k\geq2$) in $K_{0}\left(  \mathfrak
{A}_{L}\right)  $ and in $\ker\left(  \tau_{L}\right)  $ when they are
present, as they are in many examples; see, e.g., Examples \ref{Exa1} and
\ref{Exa2}. The convention regarding $L=\left(  L_{1},\dots,L_{d}\right)  $ is
as in Theorem \ref{ThmSub.2}. We assume $1\leq L_{1}\leq L_{2}\leq\dots\leq
L_{d}$, and we count the values of the $L_{i}$'s with multiplicity according
to:%
\begin{equation}
m_{j}:=\#\left\{  L_{i}\mid L_{i}=j\right\}  \label{eqCyc.1}%
\end{equation}
for $j=1,\dots,N$ where $N:=L_{d}$. Then the matrix $J=J_{L}$ takes the
form\label{LOSJ_3}%
\begin{equation}
J_{L}=%
\begin{pmatrix}
\vphantom{\vdots}m_{1} & 1 & 0 & \cdots & 0 & 0\\
\vphantom{\vdots}m_{2} & 0 & 1 & \cdots & 0 & 0\\
\vphantom{\vdots}m_{3} & 0 & 0 & \ddots & 0 & 0\\
\vdots & \vdots &  & \ddots & \ddots & \vdots\\
\vphantom{\vdots}m_{N-1} & 0 & 0 &  & 0 & 1\\
\vphantom{\vdots}m_{N} & 0 & 0 & \cdots & 0 & 0
\end{pmatrix}
. \label{eqCyc.2}%
\end{equation}
We always assume $\gcd\left\{  i\mid m_{i}\neq0\right\}  =1$. With this
convention, we have $m_{N}\geq1$. Let $a:=e^{-\beta}$ where $\beta$ is the
unique solution to%
\begin{equation}
\sum_{i}e^{-\beta L_{i}}=\sum_{j}m_{j}e^{-\beta j}=1. \label{eqCyc.3}%
\end{equation}
As explained in (\ref{eqRemNewBis.15})--(\ref{eqRemNewBis.18}), $\lambda
:=e^{\beta}$ is the Perron--Frobenius
eigenvalue\label{LOSPerronFrobeniuseigenvalue_5} for $J_{L}$.

The results in this chapter are somewhat technical. The matrix $J$ is given a
representation which admits a triangular form $\left(
\begin{tabular}
[c]{c|c}%
$J_{0_{\mathstrut}}$ & $Q$\\\hline
$0$ & $R$%
\end{tabular}
\right)  $ where $J_{0}$ and $R$ have the same type (\ref{eqCyc.2}) as $J$,
and $Q$ is of rank one (see Theorem \ref{ThmCyc.7}). Hence it is easy to read
off the determinants of $J$ and $J_{0}$. We use this to show that the prime
factors of these determinants are $C^{\ast}$-isomorphism invariants (Theorem
\ref{CorCyc.9}).

Each lattice $\mathbb{Z}^{N}$ is (linearly) isomorphic to the space
$\mathcal{V}_{N}$ of polynomials $f\in\mathbb{Z}\left[  x\right]  $ of degree
$\leq N-1$. This means that matrix multiplication by $J_{L}$ in $\mathbb{Z}%
^{N}$ is equivalent to an operation on the polynomials $\mathbb{Z}\left[
x\right]  $ of degree $\leq N-1$. This operation can be described by the
following explicit representation.

\begin{lemma}
\label{LemCyc.3}Define\label{LOSVN_1}
\begin{equation}
\mathcal{V}_{N}:=\left\{  f\left(  x\right)  \in\mathbb{Z}\left[  x\right]
\mid\deg f\leq N-1\right\}  . \label{eqCyc.11}%
\end{equation}
Let\label{LOSfm_1}%
\begin{equation}
f_{m}\left(  x\right)  :=m_{1}+m_{2}x+\dots+m_{N}x^{N-1}. \label{eqCyc.12}%
\end{equation}
Then matrix multiplication by $J$ in $\mathbb{Z}^{N}$ induces the following
operation $\tilde{J}$ on $\mathcal{V}_{N}$:%
\begin{equation}
\left(  \smash{\tilde{J}}f\right)  \left(  x\right)  =f\left(  0\right)
f_{m}\left(  x\right)  +\frac{f\left(  x\right)  -f\left(  0\right)  }%
{x},\qquad f\in\mathcal{V}_{N}. \label{eqCyc.13}%
\end{equation}
\end{lemma}

\begin{proof}
For $k=\left(  k_{1},\dots,k_{N}\right)  \in\mathbb{Z}^{N}$, let $f_{k}\left(
x\right)  =k_{1}+k_{2}x+\dots+k_{N}x^{N-1}$. Then%
\begin{align*}
\left(  \smash{\tilde{J}}f_{k}\right)  \left(  x\right)   &  =f_{Jk}\left(
x\right) \\
&  =m_{1}k_{1}+k_{2}+\left(  m_{2}k_{1}+k_{3}\right)  x+\dots\\
&  \qquad+\left(  m_{N-1}k_{1}+k_{N}\right)  x^{N-2}+m_{N}k_{1}x^{N-1}\\
&  =k_{1}\sum_{i=1}^{N}m_{i}x^{i-1}+\sum_{j=2}^{N}k_{j}x^{j-2}\\
&  =f_{k}\left(  0\right)  f_{m}\left(  x\right)  +\frac{f_{k}\left(
x\right)  -f_{k}\left(  0\right)  }{x},
\end{align*}
which proves the lemma.
\end{proof}

The construction of $K_{0}\left(  \mathfrak{A}_{L}\right)  $ and $\ker\left(
\tau_{L}\right)  $ involves the Frobenius eigenvector $\alpha=\left(
\alpha_{1},\alpha_{2},\dots,\alpha_{N}\right)  $ which
solves\label{LOSPerronFrobeniuslefteigenvector_4}%
\begin{equation}
\alpha J=\lambda\alpha\label{eqCyc.14}%
\end{equation}
where $\lambda=e^{\beta}$ is the Frobenius eigenvalue. (See (\ref{eqRemNew.1}).)

\begin{lemma}
\label{CorCyc.4}Let $a:=\lambda^{-1}=e^{-\beta}$.\label{LOSa} When normalized
with $\alpha_{1}=1$, the eigenvector $\alpha$ from \textup{(\ref{eqCyc.14})}
is\label{LOSalphaeigenvector_2}%
\begin{equation}
\alpha=\left(  1,a,a^{2},\dots,a^{N-1}\right)  . \label{eqCyc.15}%
\end{equation}
\end{lemma}

\begin{proof}
This was verified in (\ref{eqRemNewBis.17}).
\end{proof}

\begin{lemma}
\label{LemCyc.5}Let $\alpha=\left(  1,a,\dots,a^{N-1}\right)  $ be the
Frobenius eigenvector, and let%
\[
p_{a}\left(  x\right)  \in\mathbb{Z}\left[  x\right]
\]
be the minimal polynomial of $a=e^{-\beta}$. With the
identification\label{LOSVN_2}
\[
\mathbb{Z}^{N}\cong\mathcal{V}_{N}=\left\{  f\in\mathbb{Z}\left[  x\right]
\mid\deg f\leq N-1\right\}  ,
\]
as in \textup{(\ref{eqCyc.11}),} the following two conditions are equivalent
for $k=\left(  k_{1},\dots,k_{N}\right)  \in\mathbb{Z}^{N}$:

\begin{enumerate}
\item \label{LemCyc.5(1)}$k\in\left\{  \alpha\right\}  ^{\perp}$.

\item \label{LemCyc.5(2)}$p_{a}\left(  x\right)  |f_{k}\left(  x\right)  $,
where $f_{k}\left(  x\right)  =\sum_{i=1}^{N}k_{i}x^{i-1}$.\label{LOSfm_2}
\end{enumerate}
\end{lemma}

\begin{proof}
We have\label{LOSkalpha}
\begin{equation}%
\ip{k}{\alpha}%
=\sum_{i=1}^{N}k_{i}a^{i-1}=f_{k}\left(  a\right)  , \label{eqCyc.17}%
\end{equation}
showing that $f_{k}\left(  a\right)  =0$ if and only if $k\in\left\{
\alpha\right\}  ^{\perp}$. But $f_{k}\left(  x\right)  $ is divisible by
$p_{a}\left(  x\right)  $ if and only if $a$ is a root.
\end{proof}

\begin{corollary}
\label{CorCyc.6}If $D:=$ degree of $p_{a}\leq N-1$, then the subgroup
$\left\{  \alpha\right\}  ^{\perp}\cap\mathbb{Z}^{N}$ may be represented in
the form%
\[
\left\{  q\left(  x\right)  p_{a}\left(  x\right)  \mid q\left(  x\right)
\in\mathcal{V}_{N-D}\right\}  .
\]
If $D=N$, then $\left\{  \alpha\right\}  ^{\perp}\cap\mathbb{Z}^{N}=\left\{
0\right\}  $. In any case, $J$ leaves $\left\{  \alpha\right\}  ^{\perp}%
\cap\mathbb{Z}^{N}$ invariant, and if $D\leq N-1$, $J$ induces an operator
$q\mapsto J_{0}\left(  q\right)  $ on $\mathcal{V}_{N-D}$ by%
\begin{equation}
J\left(  qp_{a}\right)  \left(  x\right)  =\left(  J_{0}q\right)  \left(
x\right)  p_{a}\left(  x\right)  ,\qquad q\in\mathcal{V}_{N-D}.
\label{eqCyc.18}%
\end{equation}
\end{corollary}

\begin{proof}
The representation $q\left(  x\right)  p_{a}\left(  x\right)  $ is unique
since $p_{a}\left(  x\right)  $ is irreducible. To see that $\left\{
\alpha\right\}  ^{\perp}\cap\mathbb{Z}^{N}$ is invariant under $J$, use
(\ref{eqCyc.14}) directly, or substitute $x=a$ into (\ref{eqCyc.13}) as
follows: If $f\in\mathbb{Z}\left[  x\right]  $ satisfies $f\left(  a\right)
=0$, then
\[
\left(  \smash{\tilde{J}}f\right)  \left(  a\right)  =f\left(  0\right)
f_{m}\left(  a\right)  +\frac{f\left(  a\right)  -f\left(  0\right)  }%
{a}=f\left(  0\right)  a^{-1}-f\left(  0\right)  a^{-1}=0,
\]
where we used the identity $f_{m}\left(  a\right)  =a^{-1}$ which in turn is
equivalent to (\ref{eqCyc.3}).
\end{proof}

We need one more prelude to the main theorem of this chapter. As in Lemma
\ref{LemCyc.5}, let $p_{a}\in\mathbb{Z}\left[  x\right]  $ be the minimal
polynomial of $a=e^{-\beta}$, and let $p_{\lambda}$ be the minimal polynomial
of the Perron--Frobenius eigenvalue $\lambda=1/a=e^{\beta}$. It is clear that
these polynomials have the same degree $D$, and up to a sign\label{LOSplambda}%
\label{LOSpa}%
\begin{equation}
p_{\lambda}\left(  x\right)  =x^{D}p_{a}\left(  \frac{1}{x}\right)  ,\qquad
p_{a}\left(  x\right)  =x^{D}p_{\lambda}\left(  \frac{1}{x}\right)  .
\label{eqpp}%
\end{equation}
Since $\lambda$ is a root of the monic polynomial (\ref{eqExt.10}) in
$\mathbb{Z}\left[  x\right]  $, it follows that $p_{\lambda}$ is a monic
polynomial, and hence the constant term in $p_{a}\left(  x\right)  $ is $\pm
1$, i.e.,%
\begin{equation}
p_{a}\left(  0\right)  \in\left\{  \pm1\right\}  . \label{eqppm1}%
\end{equation}
(This also follows from (\ref{eqSubNew.29}), or (\ref{eqRemNewBis.41}).) We
will often fix the normalization of $p_{a}$ such that $p_{a}\left(  0\right)
=1$.

\begin{theorem}
\label{ThmCyc.7}Let $J$ be a matrix of the form \textup{(\ref{eqCyc.2})} with
the $m_{i}$ positive integers, $m_{N}\neq0$, $\gcd\left\{  i\mid m_{i}%
\neq0\right\}  =1$. Normalize the minimal polynomial $p_{a}\left(  x\right)  $
by $p_{a}\left(  0\right)  =1$. Decompose the polynomial $f_{m}\left(
x\right)  =m_{1}+m_{2}x+\dots+m_{N}x^{N-1}$, given in \textup{(\ref{eqCyc.12}%
),} by the Euclidean algorithm, yielding\label{LOSfm_3}%
\begin{equation}
f_{m}\left(  x\right)  =q_{m}\left(  x\right)  p_{a}\left(  x\right)
+r_{m}\left(  x\right)  , \label{eqCyc.20}%
\end{equation}
where $q_{m}\left(  x\right)  =\sum_{k=1}^{N-D}Q_{k}x^{k-1}$,\label{LOSqm}
$r_{m}\left(  x\right)  =\sum_{k=1}^{D}R_{k}x^{k-1}$.\label{LOSrm_1} It
follows that, in the basis%
\begin{equation}
\left\{  p_{a}\left(  x\right)  ,xp_{a}\left(  x\right)  ,\dots,x^{N-D-1}%
p_{a}\left(  x\right)  ,1,x,\dots,x^{D-1}\right\}  \label{eqpbasis}%
\end{equation}
for $\mathbb{Z}^{N}\cong\mathcal{V}_{N}$, the operator $J$ is given by%
\begin{equation}
J=\left(
\begin{tabular}
[c]{cccccc|cccccc}%
\vphantom{$\vdots$}$Q_{1}$ & $1$ & $0$ & $\cdots$ & $0$ & $0$ & $Q_{1}$ & $0$%
& $0$ & $\cdots$ & $0$ & $0$\\
\vphantom{$\vdots$}$Q_{2}$ & $0$ & $1$ & $\cdots$ & $0$ & $0$ & $Q_{2}$ & $0$%
& $0$ & $\cdots$ & $0$ & $0$\\
\vphantom{$\vdots$}$Q_{3}$ & $0$ & $0$ & $\ddots$ & $0$ & $0$ & $Q_{3}$ & $0$%
& $0$ & $\cdots$ & $0$ & $0$\\
$\vdots$ & $\vdots$ &  & $\ddots$ & $\ddots$ & $\vdots$ & $\vdots$ & $\vdots$%
& $\vdots$ & $\ddots$ & $\vdots$ & $\vdots$\\
\vphantom{$\vdots$}$Q_{N-D-1}$ & $0$ & $0$ &  & $0$ & $1$ & $Q_{N-D-1}$ & $0$%
& $0$ & $\cdots$ & $0$ & $0$\\
\vphantom{$\vdots$}$Q_{N-D_{\mathstrut}}$ & $0$ & $0$ & $\cdots$ & $0$ & $0$ &
$Q_{N-D}$ & $0$ & $0$ & $\cdots$ & $0$ & $0$\\\hline
\vphantom{$\vdots$}$0$ &  & \rlap{$\mkern2mu\cdots$} &  &  & $0$ & $R_{1}$ &
$1$ & $0$ & $\cdots$ & $0$ & $0$\\
\vphantom{$\vdots$} &  &  &  &  &  & $R_{2}$ & $0$ & $1$ & $\cdots$ & $0$ &
$0$\\
\vphantom{$\vdots$}\raisebox{-10pt}[0pt][0pt]{$\vdots$} &  & \raisebox
{-12pt}[0pt][0pt]{\rlap{\kern3pt{\Huge{$0$}}}} &  &  & \raisebox
{-10pt}[0pt][0pt]{$\vdots$} & $R_{3}$ & $0$ & $0$ & $\ddots$ & $0$ & $0$\\
\vphantom{$\vdots$} &  &  &  &  &  & $\vdots$ & $\vdots$ &  & $\ddots$ &
$\ddots$ & $\vdots$\\
\vphantom{$\vdots$} &  &  &  &  &  & $R_{D-1}$ & $0$ & $0$ &  & $0$ & $1$\\
\vphantom{$\vdots$}$0$ &  & \rlap{$\mkern2mu\cdots$} &  &  & $0$ & $R_{D}$ &
$0$ & $0$ & $\cdots$ & $0$ & $0$%
\end{tabular}
\right)  . \label{eqCyc.22}%
\end{equation}
In the extreme case $D=N$, the upper left-hand matrix in
\textup{(\ref{eqCyc.22})} disappears, and the lower right-hand matrix is just
\textup{(\ref{eqCyc.2}).} If $D=N-1$, the upper left-hand matrix is $\left(
Q_{1}\right)  $, and if $D=1$, the lower right-hand matrix is $\left(
R_{1}\right)  =\left(  R_{D}\right)  $. In general, the coefficients
$R_{1},\dots,R_{D}$ can be computed from the formula\label{LOSrm_2}
\begin{equation}
r_{m}\left(  x\right)  =\frac{1-p_{a}\left(  x\right)  }{x}. \label{eqrp}%
\end{equation}
\textup{(}Without the normalization $p_{a}\left(  0\right)  =1$, the upper
left-hand matrix elements $Q_{i}$ must be replaced by $p_{a}\left(  0\right)
Q_{i}$ where $p_{a}\left(  0\right)  \in\left\{  \pm1\right\}  $.\textup{)}
\end{theorem}

\begin{proof}
We leave the modifications needed to cope with the extremal cases $D=N,N-1,1$
to the reader, and consider the generic situation $1<D<N-1$. We use formula
(\ref{eqCyc.13}) in calculating $J$ in the basis defined from Lemma
\ref{LemCyc.5} and Corollary \ref{CorCyc.6}. Define\label{LOSej}%
\begin{equation}
e_{j}:=x^{j}p_{a}\left(  x\right)  ,\qquad j=0,\dots,N-D-1. \label{eqej}%
\end{equation}
Then $\left\{  e_{j}\right\}  $ is a basis for $\left\{  \alpha\right\}
^{\perp}\cap\mathbb{Z}$ by Lemma \ref{LemCyc.5}. Furthermore,%
\begin{align*}
J\left(  e_{0}\right)   &  =J\left(  p_{a}\right)  =p_{a}\left(  0\right)
f_{m}\left(  x\right)  +\frac{p_{a}\left(  x\right)  -p_{a}\left(  0\right)
}{x}\\
&  =p_{a}\left(  0\right)  q_{m}\left(  x\right)  p_{a}\left(  x\right)
+\underset{\text{remainder}}{\underbrace{\left(  p_{a}\left(  0\right)
r_{m}\left(  x\right)  +\frac{p_{a}\left(  x\right)  -p_{a}\left(  0\right)
}{x}\right)  }}.
\end{align*}
Since $\deg\left(  \left(  p_{a}\left(  x\right)  -p_{a}\left(  0\right)
\right)  /x\right)  <D$ it follows from Corollary \ref{CorCyc.6} that the
remainder is zero (this accounts for (\ref{eqrp})), and%
\[
J\left(  e_{0}\right)  =p_{a}\left(  0\right)  q_{m}\left(  x\right)
p_{a}\left(  x\right)  =p_{a}\left(  0\right)  \sum_{j=1}^{N-D}Q_{j}e_{j-1},
\]
which accounts for the upper left column in (\ref{eqCyc.22}) via formula
(\ref{eqCyc.20}). Since for $j>0$%
\[
J\left(  x^{j}p_{a}\right)  =0+\frac{x^{j}p_{a}-0}{x}=x^{j-1}p_{a},
\]
the rest of the left half of the matrix (\ref{eqCyc.22}) is accounted for.

For the rest of the entries in the formula (\ref{eqCyc.22}) for $J$, pick the
monomials $1,x,\dots,x^{D-1}$ as a basis for the remainder terms in the
Euclidean representation of $\mathbb{Z}^{N}\cong\mathcal{V}_{N}$. Using again
(\ref{eqCyc.13}), we get%
\[
J\left(  1\right)  =f_{m}\left(  x\right)  =q_{m}\left(  x\right)
p_{a}\left(  x\right)  +r_{m}\left(  x\right)  ,
\]
which accounts for the $\left(  N-D+1\right)  $'st column in (\ref{eqCyc.22}).

For $j$ such that $0<j<D$ we have, using (\ref{eqCyc.13}):%
\[
J\left(  x^{j}\right)  =0+\frac{x^{j}-0}{x}=x^{j-1},
\]
and that accounts for the remaining columns in (\ref{eqCyc.22}).
\end{proof}

\begin{corollary}
\label{CorCyc.8}Assume $1\leq D\leq N-1$. Then the relationship between the
determinants of $J$ and the restriction $J_{0}$ of $J$ to $\left\{
\alpha\right\}  ^{\perp}\cap\mathbb{Z}^{N}$ is given by%
\begin{equation}
\det\left(  J\right)  =\left(  -1\right)  ^{D-1}R_{D}\det\left(  J_{0}\right)
=\left(  -1\right)  ^{N-1}m_{N}, \label{eqdetdet1}%
\end{equation}
and%
\begin{equation}
\det\left(  J_{0}\right)  =\left(  -1\right)  ^{N-D-1}p_{a}\left(  0\right)
Q_{N-D}, \label{eqdetdet2}%
\end{equation}
and therefore%
\begin{equation}
p_{a}\left(  0\right)  Q_{N-D}=\left(  -1\right)  ^{D}\frac{m_{N}}{R_{D}},
\label{eqdetdettf}%
\end{equation}
which implies that $Q_{N-D}\neq0$.
\end{corollary}

\begin{proof}
Use the standard rules for computing determinants on (\ref{eqCyc.22}), and use
(\ref{eqCyc.20}).
\end{proof}

Note that the number $m_{N}$ is \emph{not} an isomorphism invariant. See, for
example, \textup{(\ref{eqBrunt.22})--(\ref{eqBrunt.23}),} or let us consider
the following example from \cite[Proposition 5 and following remark]{BJKR00}.
If%
\begin{equation}
J=%
\begin{pmatrix}
4 & 1\\
32 & 0
\end{pmatrix}
,\qquad J^{\prime}=%
\begin{pmatrix}
6 & 1\\
16 & 0
\end{pmatrix}
, \label{eqjjp}%
\end{equation}
then in both cases the dimension group $G_{J}$ \textup{(}resp.\ $G_{J^{\prime
}}$\textup{)} is $\mathbb{Z}\left[  \frac{1}{2}\right]  \oplus\mathbb{Z}%
\left[  \frac{1}{2}\right]  $ with order given by $\left(  x,y\right)  >0$
$\iff8x+y>0$. Furthermore, $a=\frac{1}{8}$, so the minimal polynomial is
$p_{a}\left(  x\right)  =8x-1$ in both cases. Clearly $m_{2}=32$ for $J$ and
$m_{2}^{\prime}=16$ for $J^{\prime}$, so $m_{N}=m_{2}$ is not an invariant.
But, as%
\[
4+32x=4p_{a}\left(  x\right)  +8,\qquad6+16x=2p_{a}\left(  x\right)  +8
\]
we have $R_{1}=8$ for both $J$ and $J^{\prime}$, so this does not \emph{a
priori} rule out that $R_{D}$ is an invariant. This is, however, ruled out by
\textup{(\ref{eqBrunt.22}),} where $R_{D}$ has the value $144$, $12$ for the
two matrices respectively. We will in fact prove in Theorem
\textup{\ref{CorCyc.9}} that the sets of prime factors of $m_{N}$, $R_{D}$,
respectively, are invariants. See (\ref{eqExtNew.pound}) and Figure
\ref{BratDiagsn2} for more on (\ref{eqjjp}).

\begin{proposition}
\label{ThmCyc.plus1}Let $\left(  J,\mathbb{Z}^{N}\right)  $, $D=\deg\left(
p_{a}\right)  $, and the trace $\tau\left(  \,\cdot\,\right)  =%
\ip{\,\cdot\,}{\alpha}%
$ be as described in Theorem \textup{\ref{ThmCyc.7}} and
\textup{(\ref{eqRemNewBis.20}).} Let $J_{0}$ denote the restriction of $J$ to
$\left\{  \alpha\right\}  ^{\perp}\cap\mathbb{Z}^{N}\cong\mathbb{Z}^{N-D}$.
Then\label{LOSkertau_3}%
\begin{equation}
\ker\left(  \tau\right)  =\bigcup_{n\geq0}J_{0}^{-n}\left(  \mathbb{Z}%
^{N-D}\right)  , \label{eqCyc.plus4}%
\end{equation}
where the equality refers to the identification \textup{(\ref{eqej}).}
\end{proposition}

\begin{proof}
This proposition is essentially also true in the more general situation where
$J$ is a primitive nonsingular matrix. Using the standard basis for
$\mathbb{Z}^{N}$, we saw in (\ref{eqRemNewBis.3})-(\ref{eqRemNewBis.6}) that%
\[
K_{0}\left(  \mathfrak{A}_{L}\right)  =\bigcup_{m=1}^{\infty}J^{-m}\left(
\mathbb{Z}^{N}\right)  .
\]
But $g=J^{-m}\left(  n\right)  $ is in $\ker\tau$ if and only if (using
(\ref{eqRemNew.1}) and (\ref{eqRem.1})):%
\[
0=\tau\left(  g\right)  =\tau\left(  J^{-m}\left(  n\right)  \right)  =%
\ip{\alpha}{J^{-m}n}%
=\lambda^{-m}%
\ip{\alpha}{n}%
,
\]
i.e., if and only if $n\in\mathbb{Z}^{N}\cap\left\{  \alpha\right\}  ^{\perp}%
$, that is,%
\[
J^{-m}\left(  n\right)  \cap\ker\tau=J^{-m}\left(  \mathbb{Z}^{N}\cap\left\{
\alpha\right\}  ^{\perp}\right)  .
\]
Using the basis (\ref{eqpbasis}) in Theorem \ref{ThmCyc.7}, this is
(\ref{eqdetdet1}). More specifically, we saw in (\ref{eqCyc.22}) of Theorem
\ref{ThmCyc.7} that $J$ takes the block form $\left(
\begin{tabular}
[c]{c|c}%
$J_{0_{\mathstrut}}$ & $Q$\\\hline
$0$ & $J_{R}$%
\end{tabular}
\right)  $\label{LOSJL_2} relative to the decomposition
\begin{equation}
\mathbb{Z}^{N}\cong L_{0}\oplus\mathbb{Z}^{D},\qquad L_{0}=\mathbb{Z}^{N-D}.
\label{eqCyc.plus5}%
\end{equation}
The submatrices $J_{0}$ and $J_{R}$ are both invertible in dimensions $N-D$
and $D$, respectively. Moreover (\ref{eqCyc.22}) shows that each of the
submatrices $J_{0}$ and $J_{R}$ has a form which is similar to that of $J$
itself. The $\left(  N-D\right)  \times D$ matrix $Q$ was also computed in
(\ref{eqCyc.22}). For $J^{-1}$, we therefore have the
formula\label{LOSJLinverse}%
\begin{equation}
J^{-1}=\left(
\begin{tabular}
[c]{c|c}%
$J_{0_{\mathstrut}}^{-1}$ & $-J_{0_{\mathstrut}}^{-1}QJ_{R_{\mathstrut}}^{-1}%
$\\\hline
$0$ & $J_{R}^{-1^{\mathstrut}}$%
\end{tabular}
\right)  \label{eqCyc.plus6}%
\end{equation}
and, similarly,%
\begin{equation}
J^{-n}=\left(
\begin{tabular}
[c]{c|c}%
$J_{0}^{-n}$ & $%
\begin{array}
[c]{cc}%
\ast & \ast\\
\ast & \ast
\end{array}
$\\\hline
$0$ & $J_{R}^{-n^{\mathstrut}}$%
\end{tabular}
\right)  . \label{eqCyc.plus7}%
\end{equation}
\end{proof}

\begin{theorem}
\label{CorCyc.9}The following numbers and sets of primes are isomorphism
invariants for the AF-algebras $\mathfrak{A}_{L}$, where the members of $L$
satisfy the hypothesis in Theorem \textup{\ref{ThmSub.2}:}

\begin{enumerate}
\item \label{CorCyc.9(1)}$N$, i.e., $L_{d}$,

\item \label{CorCyc.9(2)}the set of prime factors of $m_{N}$,

\item \label{CorCyc.9(3)}\hspace*{-6pt}, resp.\ \textup{(\ref{CorCyc.9(3)}%
)}$^{\prime}$, the set of prime factors of $Q_{N-D}$, resp.\ $R_{D}$, the
coefficient in the highest-order term in $q_{m}\left(  x\right)  $,
resp.\ $r_{m}\left(  x\right)  $, where
\begin{equation}
m_{1}+m_{2}x+\dots+m_{N}x^{N-1}=q_{m}\left(  x\right)  p_{a}\left(  x\right)
+r_{m}\left(  x\right)  , \label{eqPediferient.34}%
\end{equation}
and

\item \label{CorCyc.9(4)}$D=\deg\left(  p_{a}\right)  $.
\end{enumerate}

\noindent The invariants can be read off from the following commutative
diagram:%
\begin{equation}%
\begin{array}
[c]{ccccccccc}%
0 &  & 0 &  & 0 &  &  &  & 0\\
\downarrow &  & \downarrow &  & \downarrow &  &  &  & \downarrow\\
\mathbb{Z}^{N-D} & \overset{J_{0}}{\longrightarrow} & \mathbb{Z}^{N-D} &
\overset{J_{0}}{\longrightarrow} & \mathbb{Z}^{N-D} & \overset{J_{0}%
}{\longrightarrow} & \cdots & \underset{\text{\textup{ind}}}{\longrightarrow}%
& \ker\left(  \tau\right) \\
\downarrow &  & \downarrow &  & \downarrow &  &  &  & \downarrow\\
\mathbb{Z}^{N} & \overset{J}{\longrightarrow} & \mathbb{Z}^{N} & \overset
{J}{\longrightarrow} & \mathbb{Z}^{N} & \overset{J}{\longrightarrow} & \cdots
& \underset{\text{\textup{ind}}}{\longrightarrow} & K_{0}\left(  \mathfrak
{A}\right) \\
\downarrow &  & \downarrow &  & \downarrow &  &  &  & \downarrow\\
\mathbb{Z}^{D} & \overset{J_{R}}{\longrightarrow} & \mathbb{Z}^{D} &
\overset{J_{R}}{\longrightarrow} & \mathbb{Z}^{D} & \overset{J_{R}%
}{\longrightarrow} & \cdots & \underset{\text{\textup{ind}}}{\longrightarrow}%
& \mathbb{Z}\left[  \lambda^{-1}\right] \\
\downarrow &  & \downarrow &  & \downarrow &  &  &  & \downarrow\\
0 &  & 0 &  & 0 &  &  &  & 0
\end{array}
\label{eqPediferient.35}%
\end{equation}
where the vertical sequences of maps are short exact sequences, and the
horizontal maps are injective, and where $J$ has the form
\textup{(\ref{eqCyc.22}),} $J=\left(
\begin{tabular}
[c]{c|c}%
$J_{0_{\mathstrut}}$ & $Q$\\\hline
$0$ & $J_{R}$%
\end{tabular}
\right)  $.\label{LOSJL_3} This picture is also valid when $J$ is a general
nonsingular primitive $N\times N$ matrix, except that $J_{0}$, $Q$, $J_{R}$ do
not then have the special form in Theorem \textup{\ref{ThmCyc.7}.}
Nevertheless, $N$, $D$, $\operatorname*{Prim}\left(  \det J\right)  $,
$\operatorname*{Prim}\left(  \det J_{0}\right)  $, $\operatorname*{Prim}%
\left(  \det J_{R}\right)  $ are still invariants for stable $C^{\ast}%
$-isomorphism, where $\operatorname*{Prim}\left(  n\right)  $\label{LOSPrim_2}
denotes the set of prime factors of $n$ for any $n\in\mathbb{Z}\diagdown
\left\{  0\right\}  $.
\end{theorem}

\begin{remark}
\label{RemPediferient.Prim}The $\operatorname*{Prim}$-invariants are
independent in the following sense: In Chapter \textup{\ref{APP.EXA},} we give
examples $J$, $J^{\prime}$ for the same fixed values of $N$ and $D$
where\label{LOSPrimQND}%
\[
\operatorname*{Prim}\left(  Q_{N-D}^{{}}\right)  =\operatorname*{Prim}\left(
Q_{N-D}^{\prime}\right)
\]
but\label{LOSPrimRD}%
\[
\operatorname*{Prim}\left(  R_{D}^{{}}\right)  \neq\operatorname*{Prim}\left(
R_{D}^{\prime}\right)  ;
\]
and also examples with
\[%
\begin{cases}
\operatorname{Prim}\left( m_{N}^{{}}\right) =\operatorname{Prim}\left
( m_{N}^{\prime}\right)   \\
\operatorname{Prim}\left( R_{D}^{{}}\right) =\operatorname{Prim}\left
( R_{D}^{\prime}\right)
\end{cases}
\]
but%
\[
\operatorname*{Prim}\left(  Q_{N-D}^{{}}\right)  \neq\operatorname*{Prim}%
\left(  Q_{N-D}^{\prime}\right)  .
\]
\end{remark}

\begin{proof}
[Proof of Theorem \textup{\ref{CorCyc.9}}](\ref{CorCyc.9(1)})~We have already
commented that $N=L_{d}$ is the rank of the group $K_{0}\left(  \mathfrak
{A}_{L}\right)  $, so $N$ is an isomorphism invariant.

(\ref{CorCyc.9(2)})~If $n\in\mathbb{N}$, let again $\operatorname*{Prim}%
\left(  n\right)  $ denote the set of prime factors of $n$, with the
convention $\operatorname*{Prim}\left(  1\right)  =\varnothing$. If
$\mathfrak{A}_{L}$ and $\mathfrak{A}_{L^{\prime}}$ are isomorphic, it follows
from (\ref{eqRemNewTer.42}) by taking the determinant on both sides
that\label{LOSPrimmN}%
\begin{equation}%
\begin{aligned}
\operatorname*{Prim}\left( m_{N}\right)
& =\operatorname*{Prim}\left( \left| \det\left( B_{i}\right) \right| \right)
\cup\operatorname*{Prim}\left( \left| \det\left( A_{i}\right) \right
| \right) ,\\
\operatorname*{Prim}\left( m_{N^{\prime}}^{\prime}\right)
& =\operatorname*{Prim}\left( \left| \det\left( A_{i+1}\right) \right
| \right)
\cup\operatorname*{Prim}\left( \left| \det\left( B_{i}\right) \right
| \right) ,\\
\operatorname*{Prim}\left( m_{N}\right)
& =\operatorname*{Prim}\left( \left| \det\left( B_{i+1}\right) \right
| \right)
\cup\operatorname*{Prim}\left( \left| \det\left( A_{i+1}\right) \right
| \right) ,\\
\operatorname*{Prim}\left( m_{N^{\prime}}^{\prime}\right)
& =\operatorname*{Prim}\left( \left| \det A_{i+2}\right| \right)
\cup\operatorname*{Prim}\left( \left| \det B_{i+1}\right| \right) ,
\end{aligned}
\label{eqPediferient.36}%
\end{equation}
where we used Corollary \ref{CorCyc.8}. Hence
\[
\operatorname*{Prim}\left(  m_{N^{\prime}}^{\prime}\right)  \subset
\operatorname*{Prim}\left(  m_{N}\right)  \subset\operatorname*{Prim}\left(
m_{N^{\prime}}^{\prime}\right)
\]
so%
\begin{equation}
\operatorname*{Prim}\left(  m_{N^{\prime}}^{\prime}\right)
=\operatorname*{Prim}\left(  m_{N}\right)  . \label{eqPediferient.37}%
\end{equation}
Thus in particular $\operatorname*{Prim}\left(  m_{N}\right)  $ is an
isomorphism invariant, as claimed.

As the exact sequence\label{LOSkertau_4}
\begin{equation}
0\longrightarrow\ker\tau\hooklongrightarrow K_{0}\left(  \mathfrak{A}%
_{L}\right)  \longrightarrow\tau\left(  K_{0}\left(  \mathfrak{A}_{L}\right)
\right)  \longrightarrow0 \label{eqPediferient.38}%
\end{equation}
is uniquely determined by the dimension group $\left(  K_{0}\left(
\mathfrak{A}_{L}\right)  ,K_{0}\left(  \mathfrak{A}_{L}\right)  _{+}\right)
$, the group $\ker\tau$ is an isomorphism invariant. But if $J_{0}$ denotes
the restriction of $J$ to $\ker\tau=\mathbb{Z}^{N}\cap\left\{  \alpha\right\}
^{\perp}$, then $J_{0}$ identifies with the upper left-hand part of the matrix
(\ref{eqCyc.22}). But $Q_{N-D}\neq0$ by Corollary \ref{CorCyc.8} and hence
$\det J_{0}$ $\neq0$ by (\ref{eqdetdet2}). It follows from Proposition
\ref{ThmCyc.plus1} that $N-D=\operatorname*{rank}\left(  \ker\tau\right)  $ is
an isomorphism invariant and thus $D$ is so. Thus (\ref{CorCyc.9(4)}) is
proved. Furthermore, if $J^{\prime}$ is another nonsingular primitive
incidence matrix defining the same dimension group as $J$, it follows from
Proposition \ref{ThmCyc.plus1} that%
\begin{equation}
\bigcup_{n\geq0}\left(  J_{0}^{\prime}\right)  ^{-n}\left(  \mathbb{Z}%
^{N-D}\right)  \cong\bigcup_{n\geq0}J_{0}^{-n}\left(  \mathbb{Z}^{N-D}\right)
\label{eqPediferient.39}%
\end{equation}
and thus $J_{0}$ and $J_{0}^{\prime}$ are related as $J_{L}$ and
$J_{L^{\prime}}$ in (\ref{eqRemNewTer.42}), except that the $B_{i}$, $A_{i}$
now are just (necessarily nonsingular) integer matrices, without any
positivity. (See an elaboration of this in the following paragraph.) But
positivity did not play any role in the first part of the present proof, and
hence%
\[
\operatorname*{Prim}\left(  \left|  \det J_{0}\right|  \right)
=\operatorname*{Prim}\left(  \left|  \det J_{0}^{\prime}\right|  \right)  .
\]
But $\left|  \det J_{0}\right|  =\left|  Q_{N-D}\right|  $ and $\left|  \det
J_{0}^{\prime}\right|  =\left|  Q_{N-D}^{\prime}\right|  $, so
$\operatorname*{Prim}\left(  Q_{N-D}\right)  $ is an isomorphism invariant,
which shows (\ref{CorCyc.9(3)}).

By \cite{Eff81}, the groups $K_{0}\left(  \mathfrak{A}\left(  J_{L}\right)
\right)  $ and $K_{0}\left(  \mathfrak{A}\left(  J_{L^{\prime}}\right)
\right)  $ order isomorphic. Let $\theta$ be the corresponding order
isomorphism. It follows from (\ref{eqRemNewBis.34}) that $\theta$ restricts to
an isomorphism of $\ker\left(  \tau\right)  $ onto $\ker\left(  \tau^{\prime
}\right)  $. We have shown in Proposition \ref{ThmCyc.plus1} that $\ker\left(
\tau\right)  $ is constructed from $J_{0}$ the same way $K_{0}\left(
\mathfrak{A}\left(  J_{L}\right)  \right)  $ is gotten from $J_{L}$ as an
inductive limit. Now apply (\ref{eqCyc.plus4}) to both $\ker\left(
\tau\right)  $ and $\ker\left(  \tau^{\prime}\right)  $. Then the argument
from (\ref{eqRemNewTer.42}) yields%
\begin{align*}
J_{0}^{k_{i}}  &  =C_{i}E_{i},\\
\left(  J_{0}^{\prime}\right)  ^{l_{i}}  &  =E_{i+1}C_{i},
\end{align*}
where $k_{1},k_{2},\dots$, $l_{1},l_{2},\dots$ are natural numbers, and the
matrices $C_{1},C_{2},\dots$ and $E_{1},E_{2},\dots$ are $\left(  N-D\right)
\times\left(  N-D\right)  $ over $\mathbb{Z}$.

The argument which yields $\ker\left(  \tau\right)  $ as the inductive limit
$\bigcup_{n\geq0}J_{0}^{-n}\mathbb{Z}_{{}}^{N-D}$ in (\ref{eqCyc.plus4}) also
yields the following associated isomorphism:%
\begin{equation}
K_{0}\left(  \mathfrak{A}_{L}\right)  \diagup\ker\left(  \tau\right)
\cong\bigcup_{n\geq0}J_{R}^{-n}\mathbb{Z}_{{}}^{D}. \label{eqPalle.pound}%
\end{equation}
This follows by general category theory from the commutativity of the diagram
(\ref{eqPediferient.35}) and exactness of the vertical short exact sequences
of this diagram. Let us elaborate on this: Use induction, and
(\ref{eqCyc.plus4}) for $\ker\left(  \tau\right)  $, starting with the obvious
isomorphism%
\[
\mathbb{Z}^{N}\diagup\mathbb{Z}^{N-D}\cong\mathbb{Z}^{D}%
\]
given by
\[%
\begin{pmatrix}
\mathbb{Z}^{N-D}\\
\mathbb{Z}^{D}%
\end{pmatrix}
\ni%
\begin{pmatrix}
k\\
l
\end{pmatrix}
\longmapsto l\in\mathbb{Z}^{D},
\]
and arriving at%
\begin{equation}
J^{-n}%
\begin{pmatrix}
k\\
l
\end{pmatrix}
\overset{\rho_{n}}{\longmapsto}J_{R}^{-n}l. \label{eqPalle.pound1}%
\end{equation}
Since%
\begin{equation}
J^{-n}%
\begin{pmatrix}
k\\
l
\end{pmatrix}
=%
\begin{pmatrix}
J_{0}^{-n}k+Q_{n}^{\ast}l\\
J_{R}^{-n}l
\end{pmatrix}
\label{eqPalle.pound2}%
\end{equation}
for a suitable matrix $Q_{n}^{\ast}$ by (\ref{eqCyc.plus7}), we get
\begin{equation}
J^{-n}\mathbb{Z}^{N}\diagup J_{0}^{-n}\mathbb{Z}_{{}}^{N-D}\cong J_{R}%
^{-n}\mathbb{Z}_{{}}^{D} \label{eqPalle.pound3}%
\end{equation}
with the isomorphism induced by $\rho_{n}$ of (\ref{eqPalle.pound1}). It is an
isomorphism, for if $J_{R}^{-n}l=0$ then $l=0$ since $J_{R}$ is nonsingular.
So then%
\[
J^{-n}%
\begin{pmatrix}
k\\
0
\end{pmatrix}
=%
\begin{pmatrix}
J_{0}^{-n}k\\
0
\end{pmatrix}
\]
by (\ref{eqPalle.pound2}). This proves (\ref{eqPalle.pound3}).

By (\ref{eqPediferient.35}), we get
\begin{multline*}
\qquad K_{0}\left(  \mathfrak{A}_{L}\right)  \diagup\ker\left(  \tau\right)
\cong(\text{ the inductive limit}\\
\text{constructed from }J^{-n}\mathbb{Z}^{N}\diagup J_{0}^{-n}\mathbb{Z}_{{}%
}^{N-D}),\qquad
\end{multline*}
and so%
\begin{align*}
K_{0}\left(  \mathfrak{A}_{L}\right)  \diagup\ker\left(  \tau\right)   &
\cong\operatorname*{ind}_{n}J_{R}^{-n}\mathbb{Z}_{{}}^{D}\\
&  =\bigcup_{n\geq0}J_{R}^{-n}\mathbb{Z}_{{}}^{D}%
\end{align*}
by (\ref{eqPalle.pound3}). To see this, we must also check that the defining
homomorphism (\ref{eqPalle.pound1}) does indeed pass to the respective
inductive limit groups $\bigcup_{n\geq0}J^{-n}\mathbb{Z}^{N}$ and
$\bigcup_{n\geq0}J_{R}^{-n}\mathbb{Z}_{{}}^{D}$. But note that%
\[
J^{-n}%
\begin{pmatrix}
k\\
l
\end{pmatrix}
=J^{-\left(  n+1\right)  }%
\begin{pmatrix}
J_{0}k+Ql\\
J_{R}l
\end{pmatrix}
\]
and the right-hand side is mapped into%
\[
J_{R}^{-\left(  n+1\right)  }J_{R}^{{}}l=J_{R}^{-n}l
\]
under $\rho_{n+1}$ from (\ref{eqPalle.pound1}). So we have the commutative
diagram%
\[%
\begin{array}
[c]{ccccc}%
J^{-n}\mathbb{Z}^{N} &  & \makebox[0pt]{$\hooklongrightarrow$} &  &
J^{-\left(  n+1\right)  }\mathbb{Z}^{N}\\
& \underset{\rho_{n}\phantom{\rho_{n}}}{\searrow} &  & \underset{\phantom
{\rho_{n+1}}\rho_{n+1}}{\swarrow} & \\
&  & J_{R}^{-n}\mathbb{Z}_{{}}^{D} &  &
\end{array}
\]
of homomorphisms of abelian groups. As a result, there is an induced
homomorphism of the respective inductive limit groups%
\[
K_{0}\left(  \mathfrak{A}_{J}\right)  \overset{\rho}{\longrightarrow}%
\bigcup_{n\geq0}J_{R}^{-n}\mathbb{Z}_{{}}^{D},
\]
where $J=J_{L}$ for short. The formula (\ref{eqPalle.pound2}) shows that%
\[
\ker\left(  \rho\right)  \cong\bigcup_{n\geq0}J_{0}^{-n}\mathbb{Z}_{{}}%
^{N-D}\cong\ker\left(  \tau\right)
\]
where we used (\ref{eqCyc.plus4}) in the last step. Hence, by the homomorphism
theorem, we have%
\[
K_{0}\left(  \mathfrak{A}_{J}\right)  \diagup\ker\left(  \tau\right)
\overset{\tilde{\rho}}{\cong}\bigcup_{n\geq0}J_{R}^{-n}\mathbb{Z}_{{}}^{D},
\]
which is the assertion (\ref{eqPalle.pound}).

Let $L$ and $L^{\prime}$ be as in Theorem \ref{ThmSub.2} with associated
matrices $J=J_{L}$ and $J^{\prime}=J_{L^{\prime}}$, and suppose the $C^{\ast}%
$-algebras $\mathfrak
{A}_{L}$ and $\mathfrak{A}_{L^{\prime}}$ are isomorphic. The corresponding
order isomorphism
\[
\theta\colon K_{0}\left(  \mathfrak{A}_{L}\right)  \longrightarrow
K_{0}\left(  \mathfrak{A}_{L^{\prime}}\right)
\]
therefore induces isomorphisms
\[
\theta_{\text{(restriction)}}\colon\ker\left(  \tau\right)  \longrightarrow
\ker\left(  \tau^{\prime}\right)
\]
and%
\[
\theta_{\text{(quotient)}}\colon K_{0}\left(  \mathfrak{A}_{L}\right)
\diagup\ker\left(  \tau\right)  \longrightarrow K_{0}\left(  \mathfrak
{A}_{L^{\prime}}\right)  \diagup\ker\left(  \tau^{\prime}\right)  .
\]
It follows further that $\theta_{\text{(quotient)}}$ then induces an
\emph{isomorphism}%
\[
\bigcup_{n\geq0}J_{R}^{-n}\mathbb{Z}_{{}}^{D}\cong\bigcup_{n\geq0}\left(
J_{R}^{\prime}\right)  ^{-n}\mathbb{Z}_{{}}^{D}.
\]
This makes sense since we have already concluded that $N=N^{\prime}$ and
$D=D^{\prime}$. (Recall that $N-D=\operatorname*{rank}\left(  \ker\left(
\tau\right)  \right)  $.)

Now the argument after (\ref{eqPediferient.39}) applies to $J_{R}$ and
$J_{R}^{\prime}$, the same way as we used it to get identity of the sets of
primes for $\left|  \det J_{0}\right|  $ and $\left|  \det J_{0}^{\prime
}\right|  $. Using finally%
\[
R_{D}=\left|  \det J_{R}\right|  ,\qquad R_{D}^{\prime}=\left|  \det
J_{R}^{\prime}\right|  ,
\]
we conclude that%
\[
\operatorname*{Prim}\left(  R_{D}\right)  =\operatorname*{Prim}\left(
R_{D}^{\prime}\right)  ,
\]
which is the claim. The final statement of Theorem \ref{CorCyc.9} is clear
from the proof in the special case that $J$ has the form (\ref{eqCyc.2}). The
only thing that separates the general case from the special one is the special
form of $J_{0}$, $Q$, $J_{R}$ in (\ref{eqCyc.22}), and thus the formulae
$\left|  \det\left(  J\right)  \right|  =m_{N}$, $\left|  \det\left(
J_{0}\right)  \right|  =\left|  Q_{N-D}\right|  $ and $\left|  \det
J_{R}\right|  =\left|  R_{D}\right|  $.
\end{proof}

\begin{corollary}
\label{CorPediferient.nada}Assume that $J$ satisfies the hypotheses of Theorem
\textup{\ref{ThmCyc.7}} and let $\lambda$ be the Perron--Frobenius eigenvalue
of $J$. Define $J_{R}$ as in the proof of Proposition
\textup{\ref{ThmCyc.plus1}.}

\begin{enumerate}
\item \label{CorPediferient.nada(1)}There is a natural isomorphism between the
two \textup{(}unordered\/\textup{)} groups $\mathbb{Z}\left[  \tfrac
{1}{\lambda}\right]  $ and $\operatorname*{ind}\left(  J_{R}\right)
=\bigcup_{n\geq0}J_{R}^{-n}\mathbb{Z}^{D}$.

\item \label{CorPediferient.nada(2)}If $\alpha_{D}:=\left(  1,1/\lambda
,\dots,1/\lambda^{D-1}\right)  $, so that $\alpha_{D}J_{R}=\lambda\alpha_{D}$,
the isomorphism is determined by $\ip{\alpha_{D}}{\,\cdot\,}$ as follows:%
\[%
\begin{array}
[c]{ccccccc}%
0 & \longrightarrow & \operatorname*{ind}\left(  J_{R}\right)  &
\longrightarrow & \mathbb{Z}\left[  1/\lambda\right]  & \longrightarrow & 0\\
&  & \hookdownarrow &  & \hookdownarrow &  & \\
&  & \mathbb{Q}^{D} & \longrightarrow & \mathbb{R} &  & \\
&  & \makebox[0pt]{\hss$\displaystyle\cup$\hss}\makebox[0pt]
{\hss\rule[-0.15pt]{0.225pt}{6pt}\hss} &  & \makebox[0pt]{\hss$\displaystyle
\cup$\hss}\makebox[0pt]
{\hss\rule[-0.15pt]{0.225pt}{6pt}\hss} &  & \\
&  & v & \longmapsto & \ip{\alpha_{D}}{v} &  &
\end{array}
\]
where $\ip{\,\cdot\,}{\,\cdot\,}$ is the usual inner product in $\mathbb{R}%
^{D}$.
\end{enumerate}
\end{corollary}

\begin{proof}
The result follows from the equality $\tau\left(  K_{0}\left(  \mathfrak
{A}_{J}\right)  \right)  =\mathbb{Z}\left[  \tfrac{1}{\lambda}\right]  $ in
(\ref{eqRemNewBis.22}), and the natural isomorphism%
\[
\tau\left(  K_{0}\left(  \mathfrak{A}_{J}\right)  \right)  \cong K_{0}\left(
\mathfrak{A}_{J}\right)  \diagup\ker\left(  \tau\right)  \cong\bigcup_{n\geq
0}J_{R}^{-n}\mathbb{Z}_{{}}^{D}%
\]
coming from (\ref{eqPediferient.35}) in Theorem \ref{CorCyc.9} above. See also
(\ref{eqPediferientOla.54})--(\ref{eqPediferientOla.57}) below.
\end{proof}

\begin{corollary}
\label{CorPediferientOla.15}If $J$, $J^{\prime}$ are two matrices of the form
\textup{(\ref{eqCyc.2}),} one of them has a rational Perron--Frobenius
eigenvalue, and they define isomorphic AF-algebras, then both the
Perron--Frobenius eigenvalues $\lambda$, $\lambda^{\prime}$ are integral, and
$R_{D}=\lambda$ and $R_{D}^{\prime}=\lambda^{\prime}$, so%
\[
\operatorname*{Prim}\left(  R_{D}\right)  =\operatorname*{Prim}\left(
\lambda\right)  =\operatorname*{Prim}\left(  \lambda^{\prime}\right)
=\operatorname*{Prim}\left(  R_{D}^{\prime}\right)  .
\]
\end{corollary}

\begin{proof}
If for example $\lambda$ is rational then $\lambda$ is integral since it is a
solution of a monic polynomial. If $J$ and $J^{\prime}$ define isomorphic
AF-algebras, then it follows from \cite[Theorem 10]{BJKR00} that
$\mathbb{Q}\left[  \lambda^{\prime}\right]  =\mathbb{Q}\left[  \lambda\right]
=\mathbb{Q}$, i.e., $\lambda^{\prime}$ is rational and thus integral, and
$\lambda$, $\lambda^{\prime}$ are products of the same primes,
$\operatorname*{Prim}\left(  \lambda\right)  =\operatorname*{Prim}\left(
\lambda^{\prime}\right)  $. But in this case $p_{a}\left(  x\right)
=1-\lambda x$, and by (\ref{eqrp}) $r_{m}\left(  x\right)  =\lambda
x/x=\lambda=R_{D}=R$.
\end{proof}

\begin{remark}
\label{RemPediferient.20}If $g\in K_{0}\left(  \mathfrak{A}_{L}\right)
\subset\mathbb{Q}^{N}$ has coordinates $g^{\mathcal{A}}=\left(  k_{0}%
,k_{1},\dots,k_{N-1}\right)  $ relative to the old basis $\mathcal{A}=\left(
1,x,\dots,x^{N-1}\right)  $ and coordinates $g^{\mathcal{B}}=\left(
l_{0},l_{1},\dots,l_{N-1}\right)  $ relative to the new basis in Theorem
\textup{\ref{CorCyc.9},} then $g$ corresponds to the polynomial%
\begin{equation}
p_{g}\left(  x\right)  =\sum_{i=0}^{N-1}k_{i}x^{i}=\sum_{i=0}^{N-D-1}%
l_{i}x^{i}p_{a}\left(  x\right)  +\sum_{i=N-D}^{N-1}l_{i}x^{i-N+D}.
\label{eqPediferientOla.50}%
\end{equation}
If%
\begin{equation}
p_{a}\left(  x\right)  =\sum_{m=0}^{D}a_{m}x^{m}, \label{eqPediferientOla.51}%
\end{equation}
we hence compute%
\begin{equation}
p_{g}\left(  x\right)  =\sum_{j=0}^{N-1}\left(  \sum_{i=\left(  j-D\right)
\vee0}^{j\wedge\left(  N-D-1\right)  }a_{j-i}l_{i}\right)  x^{j}+\sum
_{j=0}^{D-1}l_{j+N-D}x^{j}, \label{eqPediferientOla.52}%
\end{equation}
and hence the transformation matrix between the new and the old coordinate
system is%
\begin{multline}
I_{\mathcal{B}}^{\mathcal{A}}=\vspace{0.25in}%
\vrule height0pt width0pt depth18pt%
\label{eqPediferientOla.53}\\
\quad\settowidth{\asubwidth}{$a_{D-1}$}\addtolength{\asubwidth}{-0.1\asubwidth
}\left(
\begin{tabular}
[c]{cccccccccccc}%
\llap{$\scriptstyle0$\kern1.5em}$\vphantom{\vdots}\makebox[0.5\asubwidth]{{}%
}\makebox{\raisebox{4ex}[0pt][0pt]{\makebox[0pt]{\hss$\scriptstyle0$\hss}}%
}\makebox[0pt]{\hss$a_{0}$\hss}\makebox[0.5\asubwidth]{{}}$ & $0$ & $0$ &  &
$\cdots$ &  & $0$ & $\makebox[\asubwidth]{$0$}\makebox{\raisebox
{4ex}[0pt][0pt]{\llap{$\scriptstyle N-D-1$}}}$ & $1\makebox{\raisebox
{4ex}[0pt][0pt]{\makebox[0pt]{\hss$\scriptstyle N-D$\hss}}}$ & $0$ & $\cdots$%
& $\makebox{\raisebox{4ex}[0pt][0pt]{\makebox[0pt]{\hss$\scriptstyle N-1$\hss
}}}0$\rlap{\kern1.5em$\scriptstyle0$}\\
$\vphantom{\vdots}a_{1}$ & $a_{0}$ & $0$ &  & $\cdots$ &  & $0$ & $0$ & $0$ &
$1$ &  & $0$\\
$\vdots$ &  &  & \llap{\raisebox{-9.33pt}[0pt][0pt]{$\ddots$\kern5.83pt}} &  &
&  & $\vdots$ & $\vdots$ &  & $\ddots$ & $\vdots$\\
\llap{$\scriptstyle D-1$\kern1.5em}$\vphantom{\vdots}\makebox[\asubwidth
]{$a_{D-1}$}$ & $\makebox[\asubwidth]{$a_{D-2}$}$ &  &  & $0$ &  &  &  & $0$ &
$0$ & $\cdots$ & $1$\rlap{\kern1.5em$\scriptstyle D-1$}\\
\llap{$\scriptstyle D$\kern1.5em}$\vphantom{\vdots}\makebox[\asubwidth
]{$a_{D}$}$ & $\makebox[\asubwidth]{$a_{D-1}$}$ & $\makebox[\asubwidth
]{$a_{D-1}$}$ & $\cdots$ & $a_{0}$ & $0$ &  &  & $0$ & \rlap{\kern6pt$\cdots$}%
&  & $0$\rlap{\kern1.5em$\scriptstyle D$}\\
$\vphantom{\vdots}0$ & $a_{D}$ &  &  & \llap{\raisebox{-9.33pt}%
[0pt][0pt]{$\ddots$\kern5.83pt}} & $\ddots$ &  & $\vdots$ & $\vdots$ &  &  &
$\vdots$\\
$\vdots$ &  & $\ddots$ &  &  &  & $a_{0}$ & $0$ &  &  &  & \\
$\vphantom{\vdots}$ &  & $0$ & $a_{D}$ & $\cdots$ & $a_{2}$ & $a_{1}$ &
$a_{0}$ &  &  &  & \\
$\vphantom{\vdots}$ &  &  & $0$ &  &  & $a_{2}$ & $a_{1}$ &  &  &  & \\
$\vdots$ &  &  &  & \rlap{\raisebox{9.33pt}[0pt][0pt]{\kern5.83pt$\ddots$}} &
&  & $\vdots$ &  &  &  & \\
$\vphantom{\vdots}0$ & $0$ &  & $\cdots$ &  & $0$ & $a_{D}$ & $\makebox
[\asubwidth]{$a_{D-1}$}$ & $\vdots$ &  &  & $\vdots$\\
\llap{$\scriptstyle N-1$\kern1.5em}$\vphantom{\vdots}\makebox[0.5\asubwidth
]{{}}\makebox{\raisebox{-3ex}[0pt][0pt]{\makebox[0pt]{\hss$\scriptstyle0$\hss
}}}\makebox[0pt]{\hss$0$\hss}\makebox[0.5\asubwidth]{{}}$ & $0$ &  & $\cdots$%
&  & $\makebox[\asubwidth
]{$0$}$ & $\makebox[\asubwidth
]{$0$}$ & $\makebox[\asubwidth]{$a_{D}$}\makebox{\raisebox{-3ex}%
[0pt][0pt]{\llap{$\scriptstyle N-D-1$}}}$ & $0\makebox{\raisebox
{-3ex}[0pt][0pt]{\makebox[0pt]{\hss$\scriptstyle N-D$\hss}}}$ & \rlap
{\kern6pt$\cdots$} &  & $\makebox{\raisebox{-3ex}[0pt][0pt]{\makebox
[0pt]{\hss$\scriptstyle N-1$\hss}}}0$\rlap{\kern1.5em$\scriptstyle N-1$}%
\end{tabular}
\right)  \quad.\\%
\vrule height0pt width0pt depth18pt%
\vspace{0.25in}%
\end{multline}
More interestingly, let us illustrate the power of the polynomial
representation in the computation of the trace functional $\tau\left(
g\right)  $ in the new representation. Recall from Lemma
\textup{\ref{CorCyc.4}} that%
\begin{equation}
\tau\left(  g\right)  =\ip{\alpha}{g^{\mathcal{A}}}=p_{g}\left(  a\right)  .
\label{eqPediferientOla.54}%
\end{equation}
If $\beta\in\left(  \mathbb{R}^{N}\right)  ^{\ast}$ is the row vector such
that%
\begin{equation}
\tau\left(  g\right)  =\ip{\beta}{g^{\mathcal{B}}}=\ip{\beta}{I_{\mathcal{A}%
}^{\mathcal{B}}g_{}^{\mathcal{A}}} \label{eqPediferientOla.55}%
\end{equation}
we have $\beta=\alpha I_{\mathcal{B}}^{\mathcal{A}}=\alpha\left(
I_{\mathcal{A}}^{\mathcal{B}}\right)  ^{-1}$. But note that%
\begin{equation}
\tau\left(  g\right)  =p_{a}\left(  a\right)  =\sum_{i=0}^{N-1}k_{i}a^{i}%
=\sum_{i=0}^{N-D-1}l_{i}a^{i}p_{a}\left(  a\right)  +\sum_{i=N-D}^{N-1}%
l_{i}a^{i-N+D}=\sum_{i=0}^{N-1}l_{i}a^{i-N+D} \label{eqPediferientOla.56}%
\end{equation}
and hence%
\begin{equation}
\beta=(\underset{0\mathstrut}{0},\underset{1\mathstrut}{0},\dots
,\underset{N-D-1\mathstrut}{0},\underset{N-D\mathstrut}{1},\underset
{N-D+1\mathstrut}{a},\dots,\underset{N-1\mathstrut}{\phantom{{}^{D-1}}a^{D-1}%
}) \label{eqPediferientOla.57}%
\end{equation}
gives the trace functional in the representation \textup{(\ref{eqpbasis}%
)--(\ref{eqCyc.22}).} This is useful in explicit computations of the dimension
group from the formulae \textup{(\ref{eqRemNewBis.19})--(\ref{eqRemNewBis.20}%
):}\label{LOSK0AL_5}%
\begin{equation}
K_{0}\left(  \mathfrak{A}_{L}\right)  =\bigcup_{n=0}^{\infty}J^{-n}%
\mathbb{Z}^{N}, \label{eqPediferientOla.58}%
\end{equation}
which in the new representation becomes%
\begin{equation}
K_{0}\left(  \mathfrak{A}_{L}\right)  =\bigcup_{n=0}^{\infty}\left(
J_{\mathcal{B}}^{\mathcal{B}}\right)  _{{}}^{-n}I_{\mathcal{A}}^{\mathcal{B}%
}\mathbb{Z}_{{}}^{N}. \label{eqPediferientOla.59}%
\end{equation}
Here $J_{\mathcal{B}}^{\mathcal{B}}$ is $J$ in the alias
\textup{(\ref{eqCyc.22})} and $I_{\mathcal{A}}^{\mathcal{B}}=\left(
I_{\mathcal{B}}^{\mathcal{A}}\right)  ^{-1}$ where $I_{\mathcal{B}%
}^{\mathcal{A}}$ is given in \textup{(\ref{eqPediferientOla.53}).} Note that%
\begin{equation}
\left|  \det I_{\mathcal{B}}^{\mathcal{A}}\right|  =\left|  a_{D}\right|
^{N-D} \label{eqPediferientOla.60}%
\end{equation}
so $I_{\mathcal{A}}^{\mathcal{B}}\mathbb{Z}_{{}}^{N}$ is a lattice containing
$\mathbb{Z}^{N}$ as a proper sublattice if $\left|  a_{D}\right|  \neq1$.
\end{remark}

\chapter{\label{Ker}The invariants $K_{0}\left(  \mathfrak{A}_{L}\right)
\otimes_{\mathbb{Z}}\mathbb{Z}_{n}$ and $\left(  \ker\tau\right)
\otimes_{\mathbb{Z}}\mathbb{Z}_{n}$ for $n=2,3,4,\dots$}

In this chapter we will mainly study invariants deducible from the groups
$K_{0}\left(  \mathfrak{A}_{L}\right)  =G$ and $\ker\tau=G_{0}$\label{LOSG0_3}
alone without the order structure. Of course any invariant associated to these
groups will be an invariant of the algebra. For example viewing $G$ as a
$\mathbb{Z}$-module, the groups $G\otimes\mathbb{Z}_{p}$, for $p=2,3,\dots$,
are invariants. We will also discuss structure on $G$ coming from the
embedding $\mathbb{Z}^{N}\hookrightarrow G$ given by (\ref{eqRemNewBis.3}%
)--(\ref{eqRemNewBis.6}) and the shift automorphism defined by $G$, but since
this is extraneous structure it is not clear that it leads to invariants (the
shift itself is not an invariant by the example (\ref{eqBrunt.22}%
)--(\ref{eqBrunt.23})).

Both in the construction of $K_{0}\left(  \mathfrak{A}_{L}\right)  $, and in
that of $\ker\left(  \tau_{L}\right)  $, the following inductive limit is
involved:\label{LOSLinductivelimit_1}%
\begin{equation}
\mathcal{L}\subset J^{-1}\left(  \mathcal{L}\right)  \subset J^{-2}\left(
\mathcal{L}\right)  \subset\cdots\label{eqCyc.4}%
\end{equation}
where $\mathcal{L}$ is a lattice of the same rank as the matrix $J$. But both
$\mathcal{L}$ and $J$ change in passing from $K_{0}\left(  \mathfrak{A}%
_{L}\right)  $ to $\ker\left(  \tau_{L}\right)  $ for the inductive limit construction.

We next show that quotients of these lattices, which are obviously finite
groups, are necessarily \emph{cyclic} when $J$ is the original $J_{L}$.

\begin{proposition}
\label{ProPediferient.1}The quotient%
\begin{equation}
J^{-\left(  k+1\right)  }\left(  \mathcal{L}\right)  \diagup J^{-k}\left(
\mathcal{L}\right)  \label{eqCyc.5}%
\end{equation}
is isomorphic to the cyclic group $\mathbb{Z}_{m_{N}}=\mathbb{Z}\diagup
m_{N}\mathbb{Z}$ for each $k=0,1,\dots$.
\end{proposition}

\begin{proof}
In general, if $\Gamma$ is a lattice in $\mathbb{R}^{N}$ and if $M$ is an
$N\times N$ matrix such that $M\left(  \Gamma\right)  \subset\Gamma$, then
$\Gamma\diagup M\left(  \Gamma\right)  $ is a finite abelian group of order
$\left|  \det\left(  M\right)  \right|  $. See, e.g., \cite[Proposition 5.5,
p.\ 109]{Woj97}. Applying this to (\ref{eqCyc.5}) for each $k$, we get that%
\begin{equation}
A_{k}:=J^{-\left(  k+1\right)  }\mathcal{L}\diagup J^{-k}\mathcal{L}
\label{eqCyc.6}%
\end{equation}
has order $=\left|  \det J\right|  =m_{N}$. Note from Corollary \ref{CorCyc.8}
that%
\begin{equation}
\det\left(  J\right)  =\left(  -1\right)  ^{N-1}m_{N}. \label{eqCyc.7}%
\end{equation}

A further calculation shows that the usual matrix multiplication,
$\mathcal{L}\ni l\mapsto Jl$, induces an isomorphism of abelian groups
$A_{k}\cong A_{k+1}$ for each $k$; so, in proving cyclicity, it is enough to
deal with $k=0$ where the assertion amounts to the

\begin{claim}
\label{ClaCyc.2}There is an isomorphism%
\begin{equation}
\mathbb{Z}^{N}\diagup J\mathbb{Z}^{N}\longrightarrow\mathbb{Z}_{m_{N}}
\label{eqCyc.8}%
\end{equation}
given by%
\begin{equation}%
\begin{pmatrix}
0\\
0\\
\vdots\\
0\\
i
\end{pmatrix}
+J\mathbb{Z}^{N}\longmapsto i+m_{N}\mathbb{Z},\qquad i\in\mathbb{Z}.
\label{eqCyc.9}%
\end{equation}
\end{claim}

\begin{proof}
[Proof of Claim]Since $\mathbb{Z}^{N}\diagup J\mathbb{Z}^{N}$ has order
$m_{N}$ it is enough to show that the lattice element $v_{i}:=\left(
\begin{smallmatrix}
0\\
0\\
\raisebox{0pt}[10pt]{$\vdots$}\\
0\\
i
\end{smallmatrix}
\right)  =ie_{N}$ is in $J\mathbb{Z}^{N}$ if and only if the number $i$ is
divisible by $m_{N}$. Hence, we must show that, if $i\in\mathbb{Z}$, then the
equation $v_{i}=Jk$ is solvable in $\mathbb{Z}^{N}$ if and only if $m_{N}|i$.
But $k=J^{-1}\left(  v_{i}\right)  =iJ^{-1}\left(  e_{N}\right)  $ is a
solution in $\mathbb{R}^{N}$; in fact, the explicit formula is given by
(\ref{eqSubNew.28}) or (\ref{eqThetaNyet}) as follows:%
\begin{equation}
\left\{
\begin{aligned}
k_{1} &=\frac{i}{m_{N}},\\
k_{2} &=-m_{1}\frac{i}{m_{N}},\\
k_{3} &=-m_{2}\frac{i}{m_{N}},\\
\vdots&  \qquad\vdots\\
k_{N} &=-m_{N-1}\frac{i}{m_{N}}.
\end{aligned}%
\right.  \label{eqCyc.10}%
\end{equation}
This proves that $k\in\mathbb{Z}^{N}$ if and only if $m_{N}|i$ as claimed.
\end{proof}

As mentioned, the claim proves Proposition \ref{ProPediferient.1}.
\end{proof}

This can be used to give a unique representation of elements $g\in G$.

\begin{corollary}
\label{CorCyc.4half}Let $G$ be the inductive limit group formed from
\textup{(\ref{eqCyc.4}).}

\begin{enumerate}
\item \label{CorCyc.4half(1)}In terms of the elements $v_{i}=ie_{N}$,
$i=1,\dots,m_{N}$, introduced in the proof of Claim \textup{\ref{ClaCyc.2},}
the following unique representation for points $g$ in $G$ is valid:%
\begin{equation}
g=l+J^{-1}v_{i_{1}}+J^{-2}v_{i_{2}}+\cdots\label{eqCyc.star}%
\end{equation}
where $l\in\mathcal{L}$, $i_{1},i_{2},\ldots\in\left\{  0,1,\dots
,m_{N}-1\right\}  $, and the sum is finite.

\item \label{CorCyc.4half(2)}The $\mathbb{Z}^{N}$-coordinates of $g$ in
\textup{(\ref{eqCyc.star})} are elements of $\mathbb{Z}\left[  \frac{1}{m_{N}%
}\right]  $.

\item \label{CorCyc.4half(4)}If $g$ is represented as in
\textup{(\ref{eqCyc.star})} and $l=\left(  l_{1},\dots,l_{N}\right)
\in\mathbb{Z}^{N}$, and the trace $\tau$ is given by the Frobenius eigenvector
$\alpha$ in \textup{(\ref{eqCyc.15}) as in (\ref{eqRemNewBis.20}),} then%
\begin{equation}
\tau\left(  g\right)  =\sum_{j=1}^{N}l_{j}a^{j-1}+\sum_{k\geq1}i_{k}a^{N+k-1}.
\label{eqPediferient.pound}%
\end{equation}
\end{enumerate}
\end{corollary}

\begin{proof}
Follows directly from Proposition \ref{ProPediferient.1}, Lemma \ref{LemCyc.3}%
, and Lemma \ref{CorCyc.4}. See in particular (\ref{eqCyc.10}) and
(\ref{eqRemNewBis.3})--(\ref{eqRemNewBis.6}).
\end{proof}

\begin{remark}
\label{RemPediferient.pound}Note that the right-hand side of
\textup{(\ref{eqPediferient.pound}) is related to the }$\beta$-expansion
\textup{(}$a=1/\beta$, or $\beta=\lambda$\textup{)} of the number $\tau\left(
g\right)  $ \textup{(}\cite{Par60} or \cite{FrSo92}\textup{). But the
expansion here is finite.}
\end{remark}

\begin{corollary}
\label{CorKer.A}If $G$ is described in the polynomial representation
\textup{(\ref{eqRemFeb.X})} as%
\begin{equation}
G\cong\mathbb{Z}\left[  x\right]  \diagup p_{L}\left(  x\right)
\mathbb{Z}\left[  x\right]  , \label{eqKer.a}%
\end{equation}
where\label{LOSpLx_3}%
\begin{equation}
p_{L}\left(  x\right)  =\sum_{j=1}^{N}m_{j}x^{j}-1, \label{eqKer.b}%
\end{equation}
then the element $g$ in \textup{(\ref{eqCyc.star})} is given in this
representation by the polynomial%
\begin{equation}
p\left(  x\right)  =\sum_{j=1}^{N}l_{j}x^{j-1}+\sum_{k\geq1}i_{k}x^{N+k-1}.
\label{eqKer.c}%
\end{equation}
This representation is unique within the constraint $0\leq i_{k}\leq m_{N}-1$.
\end{corollary}

\begin{proof}
Immediate from (\ref{eqRemFeb.X})--(\ref{eqRemFeb.uther}) and Corollary
\ref{CorCyc.4half}.
\end{proof}

We will consider an analogue of Corollary \ref{CorCyc.4half} for $\ker\left(
\tau\right)  $ later, in (\ref{eqKer.f}) and Corollary \ref{CorKer.C}.

\begin{remark}
\label{RemKer.B}If $q\left(  x\right)  $ is some representative polynomial in
$\mathbb{Z}\left[  x\right]  $ for $g\in G=\mathbb{Z}\left[  x\right]
\diagup\left(  p_{L}\left(  x\right)  \right)  $, one may obtain $p\left(
x\right)  $ as follows: Let $nx^{M}$ be the leading term in $q\left(
x\right)  $. If $M\leq N-1$, there is nothing to do. If $M\geq N$, add an
integer multiple of $p_{L}\left(  x\right)  x^{M-N}$ to $q\left(  x\right)  $
such that the leading coefficient is contained in $\left\{  0,1,\dots
,m_{N}-1\right\}  $. Then do the same thing for the second leading term in the
new polynomial, etc. We are soon going to adapt this procedure to the case
where $\mathbb{Z}$ is replaced by $\mathbb{Z}_{p}=\mathbb{Z}\diagup
p\mathbb{Z}$.
\end{remark}

So far we have mainly considered ordered groups in this chapter and Chapter
\ref{Pediferient}. The order is not essential for most of the results,
however. Let us first consider Theorem \ref{ThmCyc.7}.

Let $J$ be a nonsingular $N\times N$ matrix over $\mathbb{Z}$, and set
$\operatorname*{ind}\left(  J\right)  :=\bigcup_{n\geq0}J^{-n}\mathbb{Z}^{N}$.
We identify it concretely as a subgroup of $\mathbb{Q}^{N}$ (actually of
$\mathbb{Z}\left[  1/\left|  \det J\right|  \right]  ^{N}$), by the natural
inclusion mapping%
\[
\operatorname*{ind}\left(  J\right)  \hooklongrightarrow\mathbb{Q}^{N}.
\]

\begin{corollary}
\label{CorPediferient.ylem}Let $J$ be a nonsingular $N\times N$ matrix over
$\mathbb{Z}$, and suppose that there is some $D$ such that $J$ has the
triangular representation\label{LOSJL_4}
\begin{equation}
\vphantom{\left( \begin{tabular} [c]{c|c}$J_{0_{\mathstrut}}$ & $V\rlap
{$\quad\}\scriptstyle N-1$}$\\ \hline$\underset{N-D_{}^{\mathstrut}%
}{0J_{D_{\mathstrut}}^{\mathstrut}}$ & $\underset{D_{}^{\mathstrut}%
}{J_{D_{\mathstrut}}^{\mathstrut}}$\end{tabular} \right) }\left(
\begin{tabular}
[c]{c|c}%
$J_{0_{\mathstrut}}$ & $V\rlap{$\qquad\scriptstyle N-D$}$\\\hline
$\vphantom{J_{D_{\mathstrut}}^{\mathstrut}}\smash{\underset{N-D_{}%
^{\mathstrut}}{0\vphantom{J_{D_{\mathstrut}}^{\mathstrut}}}}$ & $\vphantom
{J_{D_{\mathstrut}}^{\mathstrut}}\smash{\underset{D_{}^{\mathstrut}%
}{J_{D_{\mathstrut}}^{\mathstrut}}}\rlap{$\qquad\scriptstyle D$}$%
\end{tabular}
\right)  \label{eqPediferiendPalle.39}%
\end{equation}
where the entry block matrices are also over $\mathbb{Z}$, and their sizes are
as indicated. Then there is a natural short exact sequence in the category of
abelian groups%
\[%
\begin{array}
[c]{ccccccccc}%
0 & \longrightarrow & \operatorname*{ind}\left(  J_{0}\right)  &
\longrightarrow & \operatorname*{ind}\left(  J\right)  & \longrightarrow &
\operatorname*{ind}\left(  J_{D}\right)  & \longrightarrow & 0\\
&  & \hookdownarrow &  & \hookdownarrow &  & \hookdownarrow &  & \\
0 & \longrightarrow & \mathbb{Q}^{N-D} & \longrightarrow & \mathbb{Q}^{N} &
\longrightarrow & \mathbb{Q}^{D} & \longrightarrow & 0\\
&  & \makebox[0pt]{\hss$\displaystyle\cup$\hss}\makebox[0pt]
{\hss\rule[-0.15pt]{0.225pt}{6pt}\hss} &  & \makebox[0pt]{\hss\rule
[-0.15pt]{0.225pt}{6pt}\kern2pt\rule
[-0.15pt]{0.225pt}{6pt}\hss} &  & \makebox[0pt]{\hss$\displaystyle\cup$\hss
}\makebox[0pt]
{\hss\rule[-0.15pt]{0.225pt}{6pt}\hss} &  & \\
&  & u & \llap{$\overset{\iota}{\mapsto}$}%
\begin{pmatrix}
u\\
0
\end{pmatrix}
\rlap{$\in$} &
\begin{pmatrix}
\mathbb{Q}^{N-D}\\
\mathbb{Q}^{D}%
\end{pmatrix}
& \llap{$\ni$}%
\begin{pmatrix}
u\\
v
\end{pmatrix}
\rlap{$\overset{\pi}{\mapsto}$} & v &  &
\end{array}
\]
where the morphisms in the first row are restrictions of those in the second
one, as follows: If $k\in\mathbb{Z}^{N-D}$, and $l\in\mathbb{Z}^{D}$, then%
\[
\iota\left(  J_{0}^{-n}k\right)  =J_{{}}^{-n}%
\begin{pmatrix}
k\\
0
\end{pmatrix}
\]
and%
\[
\pi\left(  J_{{}}^{-n}%
\begin{pmatrix}
k\\
l
\end{pmatrix}
\right)  =J_{D}^{-n}l.
\]
\end{corollary}

\begin{proof}
The details are contained in the last part of the proof of Theorem
\ref{CorCyc.9}.
\end{proof}

We now specialize to the case where $J$ has a form similar to
(\ref{eqRemNewBis.40}),\label{LOSJ0_1}%
\begin{equation}
J_{0}=%
\begin{pmatrix}
\vphantom{\vdots}q_{1} & 1 & 0 & \cdots & 0 & 0\\
\vphantom{\vdots}q_{2} & 0 & 1 & \cdots & 0 & 0\\
\vdots &  & \ddots & \ddots & \vdots & \vdots\\
\vphantom{\vdots}q_{M-2} & 0 &  & \ddots & 1 & 0\\
\vphantom{\vdots}q_{M-1} & 0 & 0 &  & 0 & 1\\
\vphantom{\vdots}q_{M} & 0 & 0 & \cdots & 0 & 0
\end{pmatrix}
, \label{eqKer.d}%
\end{equation}
but now we merely assume that $q_{1},\dots,q_{M}$ are (not necessarily
positive) integers and that $J_{0}$ is nonsingular, i.e., $q_{M}\neq0$. Again
one verifies that $q_{M}$ times the characteristic polynomial of $J_{0}^{-1}$
is\label{LOSp0x_1}%
\begin{equation}
p_{0}\left(  x\right)  =\sum_{j=1}^{M}q_{j}x^{j}-1, \label{eqKer.e}%
\end{equation}
and one verifies as in (\ref{eqRemFeb.X})--(\ref{eqRemFeb.uther}) that
$G_{0}=\operatorname*{ind}J_{0}$\label{LOSG0_4} identifies with the additive
group\label{LOSkertau_5}%
\begin{equation}
G_{0}\cong\mathbb{Z}\left[  x\right]  \diagup\left(  p_{0}\left(  x\right)
\right)  \label{eqKer.f}%
\end{equation}
in such a way that application of $J_{0}^{-1}$ corresponds to multiplication
by $x$.

\begin{corollary}
\label{CorKer.C}Adopt the notation and assumptions in the preceding paragraph,
in particular\label{LOSG0_5}%
\begin{equation}
G_{0}=\operatorname*{ind}\left(  J_{0}\right)  =\bigcup_{n=0}^{\infty}%
J_{0}^{-n}\mathbb{Z}^{M}. \label{eqKer.g}%
\end{equation}
Then the results \textup{(\ref{CorCyc.4half(1)}), (\ref{CorCyc.4half(2)})} in
Corollary \textup{\ref{CorCyc.4half}} and \textup{(\ref{eqKer.c})} in
Corollary \textup{\ref{CorKer.A}} remain valid, i.e.:

\begin{enumerate}
\item  In terms of the elements $v_{i}=ie_{M}$, $i=1,\dots,q_{M}$, the
following unique representations of elements $g$ of $G_{0}$ are valid:%
\begin{equation}
g=l+J_{0}^{-1}v_{i_{1}}+J_{0}^{-2}v_{i_{2}}+\cdots\label{eqKer.h}%
\end{equation}
where $l\in\mathbb{Z}^{M}$, $i_{1},i_{2},\ldots\in\left\{  0,1,\dots
,q_{M}-1\right\}  $, and the sum is finite.

\item  The polynomial representative of $g$ in \textup{(\ref{eqKer.h})} is%
\begin{equation}
\sum_{j=1}^{M}l_{j}x^{j-1}+\sum_{k\geq1}i_{k}x^{M+k-1}, \label{eqKer.i}%
\end{equation}
and this form of the representative \textup{(}i.e., with $0\leq i_{k}<q_{M}%
$\textup{)} is unique.
\end{enumerate}
\end{corollary}

\begin{proof}
As the proofs of Corollary \ref{CorCyc.4half} and Corollary \ref{CorKer.A}.
\end{proof}

In the following we will consider derived groups of the form $G_{0}%
\otimes_{\mathbb{Z}}C$ where $C$ is an abelian group. Recall from
\cite{CaEi56} that $G_{0}\otimes_{\mathbb{Z}}C$ is the free abelian group
generated by $g\otimes C$, with $g\in G_{0}$, $c\in C$, modulo the relations
$\left(  g_{1}+g_{2}\right)  \otimes C=g_{1}\otimes C+g_{2}\otimes C$,
$g\otimes\left(  C_{1}+C_{2}\right)  =g\otimes C_{1}+g\otimes C_{2}$. (This
and all the other remarks also apply to $G$ as an unordered abelian group with
the obvious modifications.) We will be interested in the case $C=\mathbb{Z}%
_{n}$, where $n\in\left\{  2,3,\dots\right\}  $.

Since $G_{0}=\mathbb{Z}\left[  x\right]  \diagup p_{0}\left(  x\right)
\mathbb{Z}\left[  x\right]  $ we have a short exact sequence\label{LOSG0_6}%
\begin{equation}
0\longrightarrow p_{0}\left(  x\right)  \mathbb{Z}\left[  x\right]
\longrightarrow\mathbb{Z}\left[  x\right]  \longrightarrow G_{0}%
\longrightarrow0. \label{eqKer.j}%
\end{equation}
But by \cite[Proposition II.4.5]{CaEi56} the functor $\cdot\otimes
_{\mathbb{Z}}C$ is right exact for any abelian group $C$, so in particular,%
\begin{equation}
p_{0}\left(  x\right)  \mathbb{Z}\left[  x\right]  \otimes_{\mathbb{Z}%
}\mathbb{Z}_{n}\longrightarrow\mathbb{Z}\left[  x\right]  \otimes
\mathbb{Z}_{n}=\mathbb{Z}_{n}\left[  x\right]  \longrightarrow G_{0}%
\otimes_{\mathbb{Z}}\mathbb{Z}_{n}\longrightarrow0 \label{eqKer.k}%
\end{equation}
is exact for $n=2,3,4,\dots$. Thus $G_{0}\otimes_{\mathbb{Z}}\mathbb{Z}_{n}$
is isomorphic to $\mathbb{Z}_{n}\left[  x\right]  $ modulo the image of
$p_{0}\left(  x\right)  \mathbb{Z}\left[  x\right]  \otimes_{\mathbb{Z}%
}\mathbb{Z}_{n}$ in $\mathbb{Z}_{n}\left[  x\right]  $, and this image is
easily seen to be $p_{0}^{\left(  n\right)  }\left(  x\right)  \mathbb{Z}%
_{n}\left[  x\right]  $, where $p_{0}^{\left(  n\right)  }\left(  x\right)  $
is the polynomial $p_{0}\left(  x\right)  $ with the coefficients reduced
modulo $n$. (This is because the map $m\rightarrow m\mod{n}$ is a \emph{ring}
morphism $\mathbb{Z}\rightarrow\mathbb{Z}_{n}$.) Thus\label{LOSG0_7}%
\begin{equation}
G_{0}\otimes_{\mathbb{Z}}\mathbb{Z}_{n}\cong\mathbb{Z}_{n}\left[  x\right]
\diagup\left(  p_{0}^{\left(  n\right)  }\left(  x\right)  \mathbb{Z}%
_{n}\left[  x\right]  \right)  \label{eqKer.l}%
\end{equation}

\begin{corollary}
\label{CorKer.D}Adopt the notation and assumptions in Corollary
\textup{\ref{CorKer.C},} and let $n\in\left\{  2,3,\dots\right\}  $. Let%
\begin{equation}
\operatorname{div}=\gcd\left\{  n,q_{M}\right\}  , \label{eqKer.m}%
\end{equation}
where $q_{M}$ is the leading coefficient in $p_{0}\left(  x\right)  $
\textup{(}see \textup{(\ref{eqKer.e})).} Then any%
\begin{equation}
g\in G_{0}\otimes\mathbb{Z}_{n}\cong\mathbb{Z}_{n}\left[  x\right]  \diagup
p_{0}^{\left(  n\right)  }\left(  x\right)  \mathbb{Z}_{n}\left[  x\right]
\label{eqKer.n}%
\end{equation}
has a unique polynomial representative of the form%
\begin{equation}
\sum_{j=1}^{M}l_{j}x^{j-1}+\sum_{k\geq1}i_{k}x^{M+k-1}, \label{eqKer.o}%
\end{equation}
where $0\leq l_{j}<n$, $0\leq i_{k}<\operatorname{div}=\gcd\left\{
n,q_{M}\right\}  $, and the right-hand sum is finite. In particular, if
$\gcd\left\{  n,q_{M}\right\}  =1$, then%
\begin{equation}
G_{0}\otimes\mathbb{Z}_{n}\cong\mathbb{Z}_{n}^{M}. \label{eqKer.p}%
\end{equation}
\end{corollary}

\begin{proof}
If $q\left(  x\right)  $ is a polynomial in $\mathbb{Z}_{n}\left[  x\right]
\diagup(p_{0}^{\left(  n\right)  }\left(  x\right)  )$, we may assume that all
the coefficients of $q$ are in $\left\{  0,1,\dots,n-1\right\}  $ by reducing
modulo $n$. Let $mx^{N}$ be the leading term in $q\left(  x\right)  $. If
$N\leq M-1$, there is nothing to do. If $N\geq M$ add an integer multiple of
$p_{0}\left(  x\right)  x^{N-M}$ to $q\left(  x\right)  \mod{n}$ such that the
leading coefficient is contained in $\left\{  0,1,\dots,\gcd\left(
n,q_{M}\right)  -1\right\}  $. Then do the same thing for the second leading
term in the new polynomial, etc. It is clear that this procedure determines
the coefficients $i_{k}$ uniquely.

In the special case that $\gcd\left(  n,q_{M}\right)  =1$, the expansion
(\ref{eqKer.o}) reduces to
\[
\sum_{j=1}^{M}l_{j}x^{j-1}%
\]
and hence%
\[
G_{0}\otimes\mathbb{Z}_{n}\cong\mathbb{Z}_{n}^{M}%
\]
in that case.
\end{proof}

\begin{remark}
\label{RemKer.E}We will see later, in Chapter \textup{\ref{APP.EXA},} that
Corollary \textup{\ref{CorKer.D}} gives an efficient method of distinguishing
cases which are not distinguished by the invariants in Chapters
\textup{\ref{Brunt}} and \textup{\ref{Pediferient}.}
\end{remark}

\chapter{\label{Str}\label{LOSK0AL_6}\label{LOSkertau_6}Associated structure
of the groups $K_{0}\left(  \mathfrak{A}_{L}\right)  $ and $\ker\tau$}

In this chapter we will study associated structure of the groups $K_{0}\left(
\mathfrak{A}_{L}\right)  $ and $\ker\tau$ which is related to the action of
$J$, to the embeddings $\mathbb{Z}^{N}\subseteq K_{0}\left(  \mathfrak{A}%
_{L}\right)  $ and $\mathbb{Z}^{N-D}\subseteq\ker\tau$, and to invariant
subgroup structure of $J$. It is not clear that these additional structures
define invariants \emph{per se,} but we will see in Chapter \ref{APP.EXA} that
they can be used to establish a quite effective machine to determine
non-isomorphism when the basic invariants from Chapter \ref{Pediferient} are
the same.

\begin{remark}
\label{Remzjzjz}Let $J$ be a nonsingular matrix with nonnegative integer
entries, and suppose, for some $k\in\mathbb{N}$, that $J^{k}$ has only
positive entries. We saw near \textup{(\ref{eqRemNewBis.1}%
)--(\ref{eqRemNewBis.6})} that then $G_{J}$ may be obtained as the inductive
limit%
\begin{equation}
\mathbb{Z}^{N}\hooklongrightarrow J^{-1}\mathbb{Z}^{N}\hooklongrightarrow
J^{-2}\mathbb{Z}^{N}\hooklongrightarrow\cdots, \label{eqzjzjz}%
\end{equation}
and $N$ is the rank of $G_{J}$. Moreover, \textup{(\ref{eqzjzjz})} defines an
embedding of $\mathbb{Z}^{N}$ as a subgroup of $G_{J}$, and we can consider
the quotient group $F\left(  J\right)  :=G_{J}\diagup\mathbb{Z}^{N}$.
\textup{(}It is not clear that the group $F\left(  J\right)  $ is an
isomorphism invariant.\textup{)} Using Theorem \textup{\ref{ThmCyc.7} we can
similarly show} that $\ker\left(  \tau\right)  $ has an analogous
representation. Its rank is $M=N-D$, and it is obtained as an inductive limit
\begin{equation}
\mathbb{Z}^{M}\hooklongrightarrow J_{0}^{-1}\mathbb{Z}^{M}\hooklongrightarrow
J_{0}^{-2}\mathbb{Z}^{M}\hooklongrightarrow\cdots, \label{eqzj0zj0z}%
\end{equation}
where $J_{0}$ is the upper left-hand submatrix of \textup{(\ref{eqCyc.22}%
):}\label{LOSJ0_2}
\begin{equation}
J_{0}=%
\begin{pmatrix}
\vphantom{\vdots}Q_{1} & 1 & \cdots & 0 & 0\\
\vphantom{\vdots}Q_{2} & 0 & \ddots & 0 & 0\\
\vdots &  & \ddots & \ddots & \vdots\\
\vphantom{\vdots}Q_{M-1} & 0 &  & 0 & 1\\
\vphantom{\vdots}Q_{M} & 0 & \cdots & 0 & 0
\end{pmatrix}
,\qquad Q_{i}\in\mathbb{Z},\;Q_{M}\neq0. \label{eqj0qqqq}%
\end{equation}
The number $m_{N}$ cannot be derived from the groups $F\left(  L\right)
:=G_{L}\diagup\mathbb{Z}^{N}$ and $F_{\tau}\left(  L\right)  :=\ker\left(
\tau\right)  \diagup\mathbb{Z}^{M}$ by the example below, where $N$,
resp.\ $M$, is still the rank of $G_{L}$, resp.\ $\ker\left(  \tau\right)  $.

Three important properties we establish in Lemma \textup{\ref{LemSubNew.14}}
below, which relate $m_{N}$ and the group $F\left(  J\right)  $, are the
following \textup{(}see details below\/\textup{):}

\begin{enumerate}
\item \label{RemSubDisclaimers(1)}$F\left(  J\right)  $ has elements of
minimal order $m_{N}$;

\item \label{RemSubDisclaimers(2)}every element of $F\left(  J\right)  $ has a
finite order which is a divisor of a power of $m_{N}$; and

\item \label{RemSubDisclaimers(3)}$m_{N}F\left(  J\right)  =F\left(  J\right)
$.
\end{enumerate}
\end{remark}

\begin{lemma}
\label{LemSubNew.14}Let the $\left(  L_{i}\right)  _{i=1}^{d}$ system be as in
Theorem \textup{\ref{ThmSub.2}.} Let $J$ be the incidence matrix
\textup{(\ref{eqCyc.2}),} and assume $\left|  \det J\right|  =m_{N}>1$.
Consider the group $F\left(  L\right)  =K_{0}\left(  \mathfrak{A}_{L}\right)
\diagup\mathbb{Z}^{N}$, where we use the concrete realization
\textup{(\ref{eqRemNewBis.6})} or \textup{(\ref{eqzjzjz})} of $K_{0}\left(
\mathfrak{A}_{L}\right)  $.

Let $e_{1},\dots,e_{N}$ be the standard basis for $\mathbb{Z}^{N}$, and
define\label{LOSgi}%
\begin{equation}
g_{i}=J^{-i}e_{N} \label{eqStr.A}%
\end{equation}
for $i=1,2,3,\dots$, and%
\begin{equation}
g_{-i}=e_{N-i} \label{eqStr.B}%
\end{equation}
for $i=0,1,\dots,N-1$. Then the elements $g_{i}$ generate $K_{0}\left(
\mathfrak{A}\right)  $ as an abelian group, and satisfy the relations%
\begin{equation}
m_{N}g_{N+i}=g_{i}-m_{1}g_{i+1}-\dots-m_{N-1}g_{N+i-1} \label{eqStr.C}%
\end{equation}
for $i=-\left(  N-1\right)  ,-\left(  N-1\right)  +1,\dots$. Moreover,
$K_{0}\left(  \mathfrak{A}\right)  $ can be characterized as a group as the
abelian group generated by elements $g_{-N+1},g_{-N+2},\dots$ satisfying these
relations, and the order in $K_{0}\left(  \mathfrak{A}\right)  $ is given by%
\begin{equation}
\sum_{i\geq1-N}c_{i}g_{i}>0\iff\sum_{i\geq1-N}c_{i}\lambda^{-i}>0,
\label{eqStr.D}%
\end{equation}
where the sums are finite and $\lambda$ is the Perron--Frobenius eigenvalue of
$J$.

Correspondingly, if we put\label{LOSxi_1}%
\begin{equation}
x_{i}=g_{i}\mod{\mathbb{Z}^{N}}, \label{eqStr.E}%
\end{equation}
the $x_{i}$ satisfy the relations%
\begin{equation}
x_{i}=0 \label{eqStr.F}%
\end{equation}
for $i=1-N,2-N,\dots,0$ and%
\begin{equation}
m_{N}x_{N+i}=x_{i}-m_{1}x_{i+1}-\dots-m_{N-1}x_{N+i-1} \label{eqStr.G}%
\end{equation}
for $i=1-N,2-N,\dots$, and $F\left(  L\right)  $ can be characterized as the
abelian group generated by these relations.
\end{lemma}

\begin{proof}
Let $g_{i}=J^{-i}e_{N}\in G\left(  L\right)  $, $x_{i}=g_{i}\mod%
{\mathbb{Z}^{N}}$, and $m=m_{N}$. Then%
\begin{align*}
mg_{1}  &  =e_{1}-m_{1}e_{2}-\dots-m_{N-1}e_{N},\\
mg_{2}  &  =J^{-1}e_{1}-m_{1}J^{-1}e_{2}-\dots-m_{N-1}J^{-1}e_{N}\\
&  =e_{2}-m_{1}e_{3}-\dots-m_{N-2}e_{N}-m_{N-1}g_{1},\\%
\intertext{and}%
mg_{3}  &  =e_{3}-m_{1}e_{4}-\dots-m_{N-2}g_{1}-m_{N-1}g_{2},\\
\vdots &  \qquad\vdots\\
mg_{N+i}  &  =g_{i}-m_{1}g_{i+1}-\dots-m_{N-1}g_{N+i-1},
\end{align*}
and in $F\left(  L\right)  =G\left(  L\right)  \diagup\mathbb{Z}^{N}$:%
\begin{align*}
mx_{1}  &  =0,\\
mx_{2}  &  =-m_{N-1}x_{1},\\
mx_{3}  &  =-m_{N-2}x_{1}-m_{N-1}x_{2},\\
\vdots &  \qquad\vdots\\
mx_{N}  &  =-m_{1}x_{1}-\dots-m_{N-1}x_{N-1},\\
\vdots &  \qquad\vdots\\
mx_{N+i}  &  =x_{i}-m_{1}x_{i+1}-\dots-m_{N-1}x_{N+i-1}.
\end{align*}

This proves the relations (\ref{eqStr.C}), (\ref{eqStr.F}) and (\ref{eqStr.G}%
). To prove that the relations (\ref{eqStr.C}) actually characterize
$K_{0}\left(  \mathfrak{A}_{L}\right)  $ we use the polynomial representation
(\ref{eqRemFeb.X})--(\ref{eqRemFeb.uther}). There $K_{0}\left(  \mathfrak
{A}_{L}\right)  $ is characterized as the additive group $\mathbb{Z}\left[
x\right]  $ modulo the linear combinations of the elements%
\[
x^{n}p_{L}\left(  x\right)
\]
for $n=0,1,2$, where $p_{L}\left(  x\right)  $ is given by
(\ref{eqRemNewBis.41}) as%
\[
p_{L}\left(  x\right)  =\sum_{j=1}^{N}m_{j}x^{j}-1.
\]
Thus $K_{0}\left(  \mathfrak{A}_{L}\right)  $ is characterized as the abelian
group generated by elements $1,x,x^{2},\dots$ with the relations%
\[
m_{N}x^{N+i}=x^{i}-m_{1}x^{i+1}-\dots-m_{N-1}x^{N+i-1}%
\]
for $i=0,1,2,\dots$. But then the abelian group defined by the relations
(\ref{eqStr.C}) above is isomorphic to this polynomial group through the map%
\[
g_{i}\longmapsto x^{i+N-1}%
\]
for $i=1-N,2-N,\dots$. This proves that the abelian group defined by
(\ref{eqStr.C}) is isomorphic to $K_{0}\left(  \mathfrak{A}_{L}\right)  $, and
furthermore, an isomorphism between the groups is given by%
\[
\sum_{i\geq1-N}c_{i}g_{i}\longmapsto\sum_{i\geq1-N}c_{i}x^{i+N-1}.
\]
Using (\ref{eqRemFeb.Y}), we thus see that (\ref{eqStr.D}) is valid.

Since $\mathbb{Z}^{N}\subseteq K_{0}\left(  \mathfrak{A}\right)  $ identifies
with the free abelian group generated by\linebreak \ $g_{1-N},g_{2-N}%
,\dots,g_{0}$ in the above picture, the remaining statement about $F\left(
L\right)  $ is immediate.
\end{proof}

\begin{remark}
\label{RemPediferient.nada}For the example $J=\left(
\begin{smallmatrix}
1 & 1 & 0\\
0 & 0 & 1\\
4 & 0 & 0
\end{smallmatrix}
\right)  $, the relations for the $x_{i}$ take the form
\begin{align*}
4x_{1}  &  =0,\\
4x_{2}  &  =0,\\
4x_{3}  &  =-x_{1},\\
\vdots &  \qquad\vdots\\
4x_{i}  &  =x_{i-3}-x_{i-2},\qquad i=4,5,\dots.
\end{align*}
This example has Perron--Frobenius eigenvalue $\lambda=2$. Thus we see that
the group $F\left(  L\right)  $ for this example is isomorphic to $\left(
\mathbb{Z}\left[  \frac{1}{2}\right]  \diagup\mathbb{Z}\right)  ^{2}$ by the
isomorphism%
\begin{align*}
x_{1}  &  \longrightarrow\left(  \frac{1}{4},0\right)  ,\\
x_{2}  &  \longrightarrow\left(  0,\frac{1}{4}\right)  ,\\
x_{3}  &  \longrightarrow\left(  -\frac{1}{4^{2\mathstrut}},0\right)  ,\\
x_{4}  &  \longrightarrow\left(  \frac{1}{4^{2\mathstrut}},-\frac{1}%
{4}\right)  ,\\
x_{5}  &  \longrightarrow\left(  \frac{1}{4^{3\mathstrut}},\frac
{1}{4^{2\mathstrut}}\right)  ,\\
x_{6}  &  \longrightarrow\left(  -\frac{2}{4^{3\mathstrut}},\frac
{1}{4^{2\mathstrut}}\right)  ,\text{\qquad etc.}%
\end{align*}
We will later, in Proposition \textup{\ref{prop10.22}} and Remark
\textup{\ref{remark10.20},} use these relations to prove the useful scaling
property
\[
\left\{  g\in K_{0}\left(  \mathfrak{A}_{J}\right)  \mid4^{i}g\in
\mathbb{Z}^{3}\right\}  =J^{-2i}\mathbb{Z}^{3}%
\]
for $i=1,2,\dots$; see, e.g., Corollary \textup{\ref{corollary10.19}} for the
use of the scaling. We know from \textup{(\ref{eqRemNewBis.31})} that
$K_{0}\left(  \mathfrak{A}_{J}\right)  $ is an extension of $\mathbb{Z}\left[
\frac{1}{2}\right]  $ by $\ker\tau$, and we will show in Example
\textup{\ref{Exa1}} by using relations analogous to the above for $\ker\tau$
that $\ker\tau$ is an extension of $\mathbb{Z}\left[  \frac{1}{2}\right]  $ by
$\mathbb{Z}$.

For the example $J=\left(
\begin{smallmatrix}
4 & 1\\
4 & 0
\end{smallmatrix}
\right)  $, $\lambda=2\cdot\left(  1+\sqrt{2}\right)  $, we have $N=2$. Here
$\operatorname*{order}\left(  x_{1}\right)  =\operatorname*{order}\left(
x_{2}\right)  =4$, while%
\[%
\begin{cases}
4x_{2i+1}=x_{2i-1}  \\
4x_{2i+2}=3x_{2i-1}+x_{2i}
\end{cases}%
\]
for $i\in\mathbb{N}$; that is, the transition matrix in this example is
$\left(
\begin{smallmatrix}
1 & 0\\
3 & 1
\end{smallmatrix}
\right)  $. It follows from Proposition \textup{\ref{prop10.22}} and Remark
\textup{\ref{remark10.20}} that also this example has the scaling property
$\left\{  g\in K_{0}\left(  \mathfrak{A}_{J}\right)  \mid4^{i}g\in
\mathbb{Z}^{2}\right\}  =J^{-2i}\mathbb{Z}^{2}$.
\end{remark}

\begin{remark}
\label{RemPediferient.ylem}The spectra of the respective matrices of the
decomposition $J=\left(
\begin{tabular}
[c]{c|c}%
$J_{0_{\mathstrut}}$ & $Q$\\\hline
$0$ & $J_{R}$%
\end{tabular}
\right)  $ in Theorem \textup{\ref{CorCyc.9}} and Corollary
\textup{\ref{CorPediferient.nada} }may be summarized by the factorization%
\begin{equation}
\det\left(  t\openone_{N}-J\right)  =\det\left(  t\openone_{N-D}-J_{0}\right)
\det\left(  t\openone_{D}-J_{R}\right)  . \label{eqPediferientPalle.38}%
\end{equation}
Since all three matrices $J$, $J_{0}$, and $J_{R}$ have the form
\textup{(\ref{eqCyc.2}),} the coefficients in the respective characteristic
polynomials are just the numbers from the first columns in the three matrices.
It is also clear from Theorem \textup{\ref{ThmCyc.7}} that the
Perron--Frobenius eigenvalue $\lambda$ is in the spectrum of $J_{R}$, and so
the points $\sigma$ in the spectrum of $J_{0}$ satisfy $\left|  \sigma\right|
<\lambda$.
\end{remark}

\begin{corollary}
\label{CorCyc.plus2}Let $\left(  J,\mathcal{L}\right)  $, $\mathcal{L}%
=\mathbb{Z}^{N}$,\label{LOSLinductivelimit_2} be as described in Proposition
\textup{\ref{ThmCyc.plus1}.} Then there is a finite decomposition series of
lattices $\mathcal{L}_{1},\mathcal{L}_{2},\dots,\mathcal{L}_{p}$ such that

\begin{enumerate}
\item \label{CorCyc.plus2(1)}the characteristic polynomial of $J|_{\mathcal{L}%
_{p}}$ is irreducible;

\item \label{CorCyc.plus2(2)}$N>\operatorname*{rank}\mathcal{L}_{1}%
>\operatorname*{rank}\mathcal{L}_{2}>\dots>\operatorname*{rank}\mathcal{L}%
_{p}$ \textup{(}if more than one term\textup{)} with each $\mathcal{L}_{i}%
$\label{LOSLk_2} invariant under $J$;

\item \label{CorCyc.plus2(3)}when the corresponding inductive limit groups
$G_{i}$ are formed from each $\mathcal{L}_{i}$, they satisfy $G_{p}%
\hookrightarrow\dots\hookrightarrow G_{2}\hookrightarrow G_{1}\hookrightarrow
K_{0}\left(  \mathfrak{A}_{L}\right)  $ and $G_{1}=\ker\left(  \tau\right)  $.
The step from $\mathcal{L}_{i}$ to $\mathcal{L}_{i+1}$ is that of Theorems
\textup{\ref{ThmCyc.7}} and \textup{\ref{CorCyc.9}.} The first column in
$J_{i}=J|_{\mathcal{L}_{i}}$ defines an element of $\mathbb{Z}\left[
x\right]  $ by \textup{(\ref{eqCyc.12}).} If this polynomial is irreducible,
then the algorithm stops. If not, it has a real root $a_{i}$, and we use the
corresponding minimal polynomial $p_{a_{i}}\left(  t\right)  \in
\mathbb{Z}\left[  t\right]  $ in passing to the next step $i+1$ of the
algorithm as done above in the proof of Theorems \textup{\ref{ThmCyc.7}} and
\textup{\ref{CorCyc.9}.} The corresponding absolute determinants $\left|  \det
J_{i}\right|  $ and polynomials $p_{i}\left(  t\right)  $ form successions of divisors.
\end{enumerate}
\end{corollary}

\begin{proof}
The proof is similar to the proof of Proposition \ref{ThmCyc.plus1}. We use
the fact that if $J=\left(
\begin{tabular}
[c]{c|c}%
$J_{0_{\mathstrut}}$ & $Q$\\\hline
$0$ & $J_{R}$%
\end{tabular}
\right)  $ represents a step in the algorithm, and if $p_{J}\left(  t\right)
$, $p_{J_{0}}\left(  t\right)  $ and $p_{J_{R}}\left(  t\right)  $ are the
corresponding characteristic polynomials, then $p_{J}\left(  t\right)
=p_{J_{0}}\left(  t\right)  \cdot p_{J_{R}}\left(  t\right)  $.

As described in (\ref{CorCyc.plus2(3)}), the argument is by recursion: Suppose%
\[
J\cong\left(
\begin{tabular}
[c]{c|c}%
$J_{1_{\mathstrut}}$ & $V_{1_{\mathstrut}}$\\\hline
$0$ & $K_{1}$%
\end{tabular}
\right)
\]
is a triangular representation as in Corollary \ref{CorPediferient.ylem}. Then
formula (\ref{eqPediferientPalle.38}) yields divisibility for the respective
characteristic polynomials%
\begin{equation}
\operatorname*{ch}\nolimits_{J}\left(  t\right)  =\operatorname*{ch}%
\nolimits_{J_{1}}\left(  t\right)  \operatorname*{ch}\nolimits_{K_{1}}\left(
t\right)  . \label{eqchchch}%
\end{equation}
If this first reduction decreases the rank, then (\ref{eqchchch}) shows that
$\operatorname*{ch}\nolimits_{J}\left(  t\right)  $ could not be irreducible.
At the first step in the reduction, Theorem \ref{CorCyc.9} and Corollary
\ref{CorPediferient.nada} show that the Perron--Frobenius eigenvalue $\lambda$
is a root of $\operatorname*{ch}\nolimits_{K_{1}}\left(  t\right)  $. We must
show that, if $\operatorname*{ch}\nolimits_{J_{1}}\left(  t\right)  $ factors
nontrivially, i.e., $\operatorname*{ch}\nolimits_{J_{1}}\left(  t\right)
=q\left(  t\right)  p\left(  t\right)  $, with $q\left(  t\right)  ,p\left(
t\right)  \in\mathbb{Z}\left[  t\right]  $, and say $p\left(  t\right)  $
irreducible, then the process may continue. Since the matrices $J_{1}$ and
$K_{1}$ have the same form as $J$ at the outset, we would get%
\begin{equation}
J_{1}\cong\left(
\begin{tabular}
[c]{c|c}%
$J_{2_{\mathstrut}}$ & $V_{2_{\mathstrut}}$\\\hline
$0$ & $K_{2}$%
\end{tabular}
\right)  , \label{eqJ1J2}%
\end{equation}
again with the properties from the proof of Theorem \ref{CorCyc.9} and
Corollary \ref{CorPediferient.ylem}. Let $f\left(  t\right)
=\operatorname*{ch}\nolimits_{J_{1}}\left(  t\right)  $. Then $J_{1}$ may be
represented, via (\ref{eqpp}), as multiplication by $t$ on $\mathbb{Z}\left[
t\right]  \diagup\left(  f\left(  t\right)  \right)  $. Let $W$ denote the
following induced linear mapping (quotient by ideals):%
\[
\mathbb{Z}\left[  t\right]  \diagup\left(  q\left(  t\right)  \right)
\overset{W}{\hooklongrightarrow}\mathbb{Z}\left[  t\right]  \diagup\left(
f\left(  t\right)  \right)  ,
\]
$W\left(  h\left(  t\right)  +q\right)  :=p\left(  t\right)  h\left(
t\right)  +\left(  f\right)  $, for $h\left(  t\right)  \in\mathbb{Z}\left[
t\right]  $. It is well defined and injective due to the assumptions made on
$f\left(  t\right)  $. Since $J_{1}$ is represented as multiplication by $t$
in $\mathbb{Z}\left[  t\right]  \diagup\left(  f\left(  t\right)  \right)  $,
the range of $W$ is then a nontrivial invariant subspace (over $\mathbb{Z}$)
for $J_{1}$, and we arrive at the triangular form (\ref{eqJ1J2}). The argument
from the proof of Theorem \ref{CorCyc.9} shows that the entries of
(\ref{eqJ1J2}) must have the same standard form as described in the previous
step. Hence the process may continue until at some step, $p$ say,
$\operatorname*{ch}\nolimits_{J_{p}}\left(  t\right)  $ is irreducible.
\end{proof}

\chapter{\label{Ext}\label{LOSExt}The invariant $\operatorname*{Ext}\left(
\tau\left(  K_{0}\left(  \mathfrak{A}_{L}\right)  \right)  ,\ker\tau\right)  $}

\setcounter{figure}{\value{figurelink}} In this chapter and the next we study
the set $\mathcal{L}\left(  \lambda\right)  $ of matrices $J_{m}$ of the form
(\ref{eqCyc.2}) such that $\lambda^{N}-m_{1}\lambda^{N-1}-\dots-m_{N-1}%
\lambda-m_{N}=0$. For the case when $\lambda\in\mathbb{Z}_{+}$, we will show
in Theorem \ref{theorem10.10}, Corollary \ref{corollary10.11}, Corollary
\ref{corollary10.12}, Theorem \ref{theorem10.15}, Corollary
\ref{corollary10.17}, and Proposition \ref{prop10.18} that $\tau\left(
v\right)  $ can be used to show non-isomorphism where $\tau$ is the normalized
trace, and $v$ is the right Perron--Frobenius eigenvector, i.e., $Jv=\lambda
v$, $v_{1}=1$. See (\ref{eqSubNew.pound}) for the explicit form of $v$.

There are examples $J$, $J^{\prime}$ such that all the three
$\operatorname*{Prim}$-invariants agree on $J$ and $J^{\prime}$ while the
$C^{\ast}$-algebras $\mathfrak{A}_{J}$ and $\mathfrak{A}_{J^{\prime}}$ are
non-isomorphic. Take, for example, $N=3$, $D=1$, $\lambda=\lambda^{\prime}=2$,
and
\[
\begin{picture}(102,38)
\put(74,1){\includegraphics
[bb=0 0 15 13,height=13bp,width=15bp]{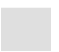}}
\put(0,0){\makebox(102,38)[l]{$\displaystyle
m=\begin{pmatrix}
m_{1}  \\
m_{2}  \\
m_{3}
\end{pmatrix}
=\left( \mkern8mu
\begin{matrix}
1  \\
0  \\
4
\end{matrix}
\mkern8mu \right) ,$}}
\end{picture}\raisebox{16.4pt}{\quad and\quad}\begin{picture}(105,38)
\put(77,1){\includegraphics
[bb=0 0 15 13,height=13bp,width=15bp]{remExtsh.eps}}
\put(0,0){\makebox(105,38)[l]{$\displaystyle
m^{\prime}=\begin{pmatrix}
m_{1}^{\prime}  \\
m_{2}^{\prime}  \\
m_{3}^{\prime}
\end{pmatrix}
=\left( \mkern8mu
\begin{matrix}
0  \\
2  \\
4
\end{matrix}
\mkern8mu \right) .$}}
\end{picture}
\]
(For more examples, see also Chapter \ref{APP.EXA} and Table
\ref{TableCycpoundprime} in Chapter \ref{Sesquilabial}.) Then the respective
triangular forms are%
\[
\begin{picture}(96,38)
\put(67,-0.9){\includegraphics
[bb=0 0 15 13,height=13bp,width=15bp]{remExtsh.eps}}
\put(67,14.1){\includegraphics
[bb=0 0 15 13,height=13bp,width=15bp]{remExtsh.eps}}
\put(0,0){\makebox(96,38)[l]{$\displaystyle
J\cong\left(
\begin{tabular}
[c]{cc|c}%
$-1$ & $1$ & $1$\\
$-2$ & $0$ & $2$\\\hline
$0$ & $0$ & $2$%
\end{tabular}
\mkern4mu \right) ,$}}
\end{picture}\qquad\begin{picture}(98,38)
\put(70,-0.9){\includegraphics
[bb=0 0 15 13,height=13bp,width=15bp]{remExtsh.eps}}
\put(70,14.1){\includegraphics
[bb=0 0 15 13,height=13bp,width=15bp]{remExtsh.eps}}
\put(0,0){\makebox(98,38)[l]{$\displaystyle
J^{\prime}\cong\left(
\begin{tabular}
[c]{cc|c}%
$-2$ & $1$ & $2$\\
$-2$ & $0$ & $2$\\\hline
$0$ & $0$ & $2$%
\end{tabular}
\mkern4mu \right) .$}}
\end{picture}
\]
and therefore%
\[
Q_{2}^{{}}=Q_{2}^{\prime}=R_{1}^{{}}=R_{1}^{\prime}=2.
\]
In the next chapters, we identify additional quantities which can be used to
distinguish $\mathfrak{A}_{J}$ and $\mathfrak{A}_{J^{\prime}}$. If $v$ denotes
the right Perron--Frobenius eigenvector, then one of these quantities is
$\tau\left(  v\right)  $. The actual non-isomorphism of the two specimens
above can, however, be established by using (\ref{eqKer.p}) in Corollary
\ref{CorKer.D}; see the $N=3$ case in Chapter \ref{APP.EXA}.

We mentioned in Chapter \ref{Rem} that the dimension group $D\left(
\mathfrak{A}_{L}\right)  $,\label{LOSDAL_5} that is the group $K_{0}\left(
\mathfrak{A}_{L}\right)  $\label{LOSK0AL_7} with the Riesz order and the
element $\left[  \openone\right]  _{0}$, is a complete isomorphism invariant
by the general theory. Objects that can be derived from $D\left(  \mathfrak
{A}_{L}\right)  $, like $\tau\left(  K_{0}\left(  \mathfrak{A}_{L}\right)
\right)  $, $\ker\tau$,\label{LOSkertau_7} $\operatorname*{Ext}$,
$N=\operatorname*{rank}\left(  K_{0}\left(  \mathfrak{A}_{L}\right)  \right)
$, and the sets of prime factors of $m_{N}$, $R_{D}$, and $Q_{N-D}$, are
secondary invariants. In this chapter and the next we shall treat the
invariant in $\operatorname*{Ext}\left(  \tau\left(  K_{0}\left(  \mathfrak
{A}_{L}\right)  \right)  ,\ker\tau\right)  $ defined by $K_{0}\left(
\mathfrak{A}_{L}\right)  $.

Aside from the two groups $\ker\left(  \tau\right)  $ and $\tau\left(
K_{0}\left(  \mathfrak{A}_{L}\right)  \right)  $ themselves, $D\left(
\mathfrak{A}_{L}\right)  $ determines the intrinsic exact sequence:%
\begin{equation}
0\longrightarrow\ker\left(  \tau\right)  \overset{\iota}{\hooklongrightarrow
}G\overset{\tau}{\longrightarrow}\tau\left(  G\right)  \longrightarrow0,
\label{eqExt.1}%
\end{equation}
where we use the shorthand notation $G=G_{L}:=K_{0}\left(  \mathfrak{A}%
_{L}\right)  $.\label{LOSK0AL_8} Hence the complete invariant $D\left(
\mathfrak{A}_{L}\right)  $ for isomorphism of the AF-$C^{\ast}$-algebras
$\mathfrak{A}_{L}$ includes (\ref{eqExt.1}), characterized as an element of
$\operatorname*{Ext}\left(  \tau\left(  G\right)  ,\ker\left(  \tau\right)
\right)  $. We shall need a few facts from homology about the
$\operatorname*{Ext}$-groups\label{LOSExtgroups}, and we refer to
\cite{MacL67} for background material: if $A$ and $C$ are abelian groups, an
element of $\operatorname*{Ext}\left(  C,A\right)  $ is an equivalence class
of short exact sequences of abelian groups%
\begin{equation}
0\longrightarrow A\overset{\iota}{\hooklongrightarrow}E\overset{\tau
}{\longrightarrow}C\longrightarrow0. \label{eqExt.2}%
\end{equation}
It is conventional to use the same letter $E$ also to denote this exact
sequence and $E$ or $\left[  E\right]  $ to denote the equivalence class. Two
elements $E$ and $E^{\prime}$ are said to be equivalent in
$\operatorname*{Ext}\left(  C,A\right)  $ if there is an isomorphism
$\theta\colon E\rightarrow E^{\prime}$ of abelian groups such that
\begin{equation}%
\begin{array}
[c]{ccccccccc}%
&  &  &  & E &  &  &  & \\
&  &  & \overset{\iota\phantom{\iota}}{\nearrow} &  & \overset{\phantom{\tau
}\tau}{\searrow} &  &  & \\
0 & \longrightarrow &  A &  & \smash{\vcenter{\hrule width0.275pt height32pt}%
\makebox[0pt] {\hss\raisebox{-12pt}{$\downarrow$}\hss}}\rlap{$\scriptstyle
\;\theta$} &  & C & \longrightarrow & 0\\
&  &  & \underset{\iota^{\prime}\phantom{\iota^{\prime}}}{\searrow} &  &
\underset{\phantom{\tau^{\prime}}\tau^{\prime}}{\nearrow} &  &  & \\
&  &  &  & E^{\prime} &  &  &  &
\end{array}
\label{eqExt.3}%
\end{equation}
commutes, or more globally if%
\begin{equation}%
\begin{array}
[c]{ccccccccc}%
0 & \longrightarrow &  A & \longrightarrow &  E & \longrightarrow &  C &
\longrightarrow & 0\\
&  & \llap{$\scriptstyle\alpha\;$}\downarrow &  & \llap{$\scriptstyle
\theta\;$}\downarrow &  & \downarrow\rlap{$\scriptstyle\;\gamma$} &  & \\
0 & \longrightarrow &  A^{\prime} & \longrightarrow &  E^{\prime} &
\longrightarrow &  C^{\prime} & \longrightarrow & 0
\end{array}
\label{eqExt.4}%
\end{equation}
commutes, where $\alpha$, $\gamma$ are isomorphisms of abelian groups. Note if
we have $\theta\in\operatorname*{Hom}\left(  E,E^{\prime}\right)  $, and if
$\alpha$ and $\gamma$ are isomorphisms, then $\theta$ will be an isomorphism
by the Short Five Lemma; see \cite{MacL67}. With a standard addition
$E+E^{\prime}$, $\operatorname*{Ext}\left(  C,A\right)  $ itself acquires the
structure of an abelian group. $E^{\prime\prime}=E+E^{\prime}$ is defined by%
\[
E^{\prime\prime}=\left\{  \left(  x,y\right)  \in E\oplus E^{\prime}\mid
\tau\left(  x\right)  =\tau^{\prime}\left(  y\right)  \right\}  \diagup
\left\{  \left(  \iota\left(  a\right)  ,-\iota^{\prime}\left(  a\right)
\right)  \mid a\in A\right\}
\]
with%
\[
\iota^{\prime\prime}\colon A\longrightarrow E^{\prime\prime}\colon
a\longmapsto\left[  \left(  \iota\left(  a\right)  ,0\right)  \right]
\]
and
\[
\tau^{\prime\prime}\colon E^{\prime\prime}\longrightarrow C\colon\left[
\left(  x,y\right)  \right]  \longrightarrow\tau\left(  x\right)  .
\]
(In these considerations, $\tau$, $\tau^{\prime}$, and $\tau^{\prime\prime}$
are only viewed as maps of abelian semigroups.) This makes
$\operatorname*{Ext}\left(  C,A\right)  $ into an abelian semigroup, with
identity element the trivial extension%
\begin{align*}
E_{0}  &  =A\oplus C,\\
\iota_{0}\colon A  &  \longrightarrow E_{0}\colon a\longrightarrow\left(
a,0\right)  ,\\
\tau_{0}\colon E_{0}  &  \longrightarrow C\colon\left(  a,b\right)
\longrightarrow b.
\end{align*}
Any element has an inverse given by%
\begin{align*}
E^{\prime}  &  =E,\\
\iota^{\prime}  &  =-\iota,\\
\tau^{\prime}  &  =\tau,
\end{align*}
and this makes $\operatorname*{Ext}\left(  C,A\right)  $ into an abelian group.

We say that (\ref{eqExt.2}) \emph{splits} if there is a $\psi\in
\operatorname*{Hom}\left(  C,E\right)  $ such that $\tau\circ\psi
=\operatorname*{id}$. This is equivalent to%
\[
E\cong A\oplus C\text{\qquad(direct sum of abelian groups),}%
\]
with trivial maps $\iota$, $\tau$, and then the corresponding $E$ is the zero
element of the abelian group $\operatorname*{Ext}\left(  C,A\right)  $. Note
that $\operatorname*{Ext}\left(  \mathbb{Z}_{2},\mathbb{Z}\right)
\cong\mathbb{Z}_{2}$. The corresponding two group elements, $0$, resp., $1$,
are (the equivalence classes of) (\ref{eqExt.5}) and (\ref{eqExt.6}):%
\begin{equation}
0\longrightarrow2\mathbb{Z}\hooklongrightarrow\left(  2\mathbb{Z}\right)
\oplus\mathbb{Z}_{2}\overset{0+\operatorname*{id}}{\longrightarrow}%
\mathbb{Z}_{2}\longrightarrow0 \label{eqExt.5}%
\end{equation}
and%
\begin{equation}
0\longrightarrow2\mathbb{Z}\hooklongrightarrow\mathbb{Z}\overset{\text{proj.}%
}{\longrightarrow}\mathbb{Z}_{2}\longrightarrow0, \label{eqExt.6}%
\end{equation}
where the second, (\ref{eqExt.6}), is non-split. More generally,%
\begin{equation}
\operatorname*{Ext}\left(  \mathbb{Z}_{m},A\right)  \cong A\diagup mA;
\label{eqExt.7}%
\end{equation}
see \cite{MacL67}. A refinement of (\ref{eqExt.7}), also due to Mac Lane et
al., is the characterization of
\begin{equation}
\operatorname*{Ext}\left(  \mathbb{Z}\left[  \tfrac{1}{k}\right]
,\mathbb{Z}^{l}\right)  \label{eqExt.8}%
\end{equation}
as a solenoid, depending on $k,l\in\mathbb{N}$, $k>1$. In particular,
(\ref{eqExt.8}) is overcountable. The description of our $\mathfrak
{A}_{L}$'s associated with $\left|  \det J_{L}\right|  =k$ in the special case
that $\ker\left(  \tau\right)  \cong\mathbb{Z}^{l}$, must be given in terms of
(\ref{eqExt.8}).

In the general case, we have $J_{L}$ of the form\label{LOSJ_4}
\begin{equation}
J_{L}=%
\begin{pmatrix}
\vphantom{\vdots}m_{1} & 1 & 0 & \cdots & 0 & 0\\
\vphantom{\vdots}m_{2} & 0 & 1 & \cdots & 0 & 0\\
\vphantom{\vdots}m_{3} & 0 & 0 & \ddots & 0 & 0\\
\vdots & \vdots &  & \ddots & \ddots & \vdots\\
\vphantom{\vdots}m_{N-1} & 0 & 0 &  & 0 & 1\\
\vphantom{\vdots}m_{N} & 0 & 0 & \cdots & 0 & 0
\end{pmatrix}
, \label{eqExt.9}%
\end{equation}
where $m_{N}=\left(  -1\right)  ^{N-1}\det J_{L}$, and the characteristic
polynomial $p_{L}\left(  \lambda\right)  $ is
\begin{equation}
p_{L}\left(  \lambda\right)  =\det\left(  \lambda-J_{L}\right)  =\lambda
^{N}-m_{1}\lambda^{N-1}-\dots-m_{N-1}\lambda-m_{N}. \label{eqExt.10}%
\end{equation}
Here $N$ is the rank of $G_{L}$. Then $\mathbb{Z}^{N}$ is embedded in $G_{L}$,
and we can introduce the quotient group $G_{L}\diagup\mathbb{Z}^{N}$. Using
this, we show that $G_{L}\diagup\mathbb{Z}^{N}$ is a specific extension of the
inductive limit group $\mathcal{C}_{m_{N}}$ defined by
\[
\mathbb{Z}_{m_{N}}\overset{m_{N}}{\longrightarrow}\mathbb{Z}_{m_{N}^{2}%
}\overset{m_{N}}{\longrightarrow}\mathbb{Z}_{m_{N}^{3}}\longrightarrow\cdots.
\]
Let $\tau$ be the normalized trace on $G_{L}$. Then there is a short exact
sequence
\[
0\longrightarrow\ker\left(  \tau\right)  \longrightarrow G_{L}\overset{\tau
}{\longrightarrow}\tau\left(  G_{L}\right)  \longrightarrow0.
\]
We further show that, if $M$ is the rank of $\ker\left(  \tau\right)  $, then
$\ker\left(  \tau\right)  \diagup\mathbb{Z}^{M}$ is an extension of a second
inductive limit group $\mathcal{C}_{k}$ formed from finite cyclic groups:%

\[
\mathbb{Z}_{k}\overset{k}{\longrightarrow}\mathbb{Z}_{k^{2}}\overset
{k}{\longrightarrow}\mathbb{Z}_{k^{3}}\overset{k}{\longrightarrow}\cdots,
\]
where $k$ divides $m_{N}$. It will follow in particular from the construction
that every element of $G_{L}\diagup\mathbb{Z}^{N}$ has a (finite) order which
divides a power of $m_{N}$; and, similarly, that every element of $\ker\left(
\tau\right)  \diagup\mathbb{Z}^{M}$ has an order which is a divisor of $k^{i}$
for some $i$ (depending on the element).

\begin{figure}[ptb]
\begin{picture}(360,528)
\put(0,0){\includegraphics
[bb=267 120 453 720,clip,width=167.4bp,height=540bp]{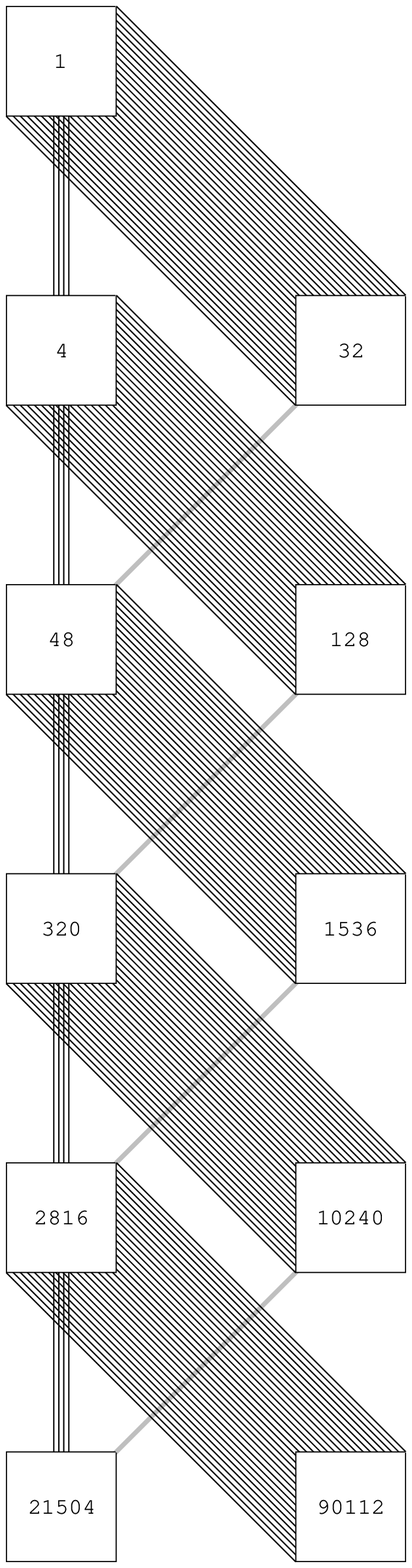}}
\put(192.6,0){\includegraphics
[bb=267 120 453 720,clip,width=167.4bp,height=540bp]{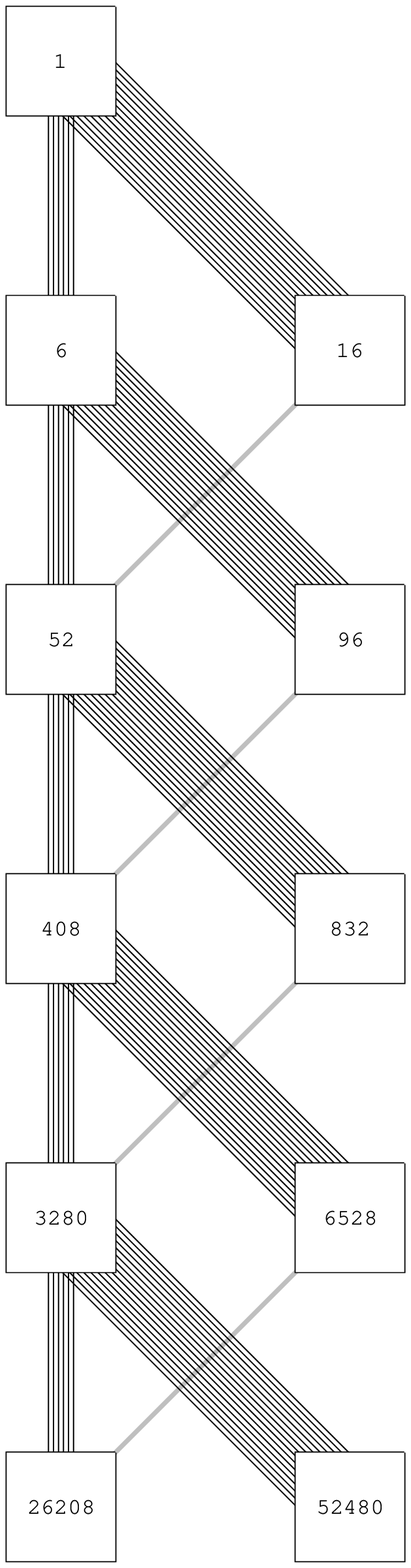}}
\end{picture}\caption{$L=\{1,1,1,1,\frogbrace{32}{2}\}$, first column
$=(4\;32)^{\mathrm{t}}$ (left); $L=\{1,1,1,1,1,1,\frogbrace{16}{2}\}$, first
column $=(6\;16)^{\mathrm{t}}$ (right). See (\ref{eqExtNew.pound}). These
diagrams represent isomorphic algebras.}%
\label{BratDiagsn2}%
\end{figure}

Note that, as a consequence of (\ref{eqExt.10}), the vector $\left(
m_{1},\dots,m_{N}\right)  $ is a similarity invariant for $J_{L}$, i.e., two
nonsingular matrices $J_{L}$ and $J_{L^{\prime}}$ of the form (\ref{eqExt.9})
are similar if and only if they are equal. But similarity of two $J_{L}$'s is
a condition which \emph{a priori} is much more restrictive than isomorphism of
the corresponding pair of $C^{\ast}$-algebras, $\mathfrak{A}_{L}$ and
$\mathfrak{A}_{L^{\prime}}$. In \cite{BJKR00} this is discussed in detail, and
we show for example that the matrices we discussed in (\ref{eqjjp}),%
\begin{equation}
J_{L}=%
\begin{pmatrix}
4 & 1\\
32 & 0
\end{pmatrix}
\text{\qquad and\qquad}J_{L^{\prime}}=%
\begin{pmatrix}
6 & 1\\
16 & 0
\end{pmatrix}
, \label{eqExtNew.pound}%
\end{equation}
define isomorphic $C^{\ast}$-algebras. (See Figure \ref{BratDiagsn2}, below.)
Other examples are in Example \ref{ExaRemNewBis.3}, in (\ref{eqBrunt.22}), in
Chapters \ref{App.Cla} and \ref{APP.EXA}, and in \cite{BJKR00}.

Let $L$ be such that $\tau\left(  G_{L}\right)  =\mathbb{Z}\left[  \frac{1}%
{k}\right]  $ for some $k$. Then $\mathfrak{A}_{L}$ corresponds to a nonzero
element of $\operatorname*{Ext}\left(  \mathbb{Z}\left[  \frac{1}{k}\right]
,\ker\left(  \tau\right)  \right)  $, if and only if%
\begin{equation}
0\longrightarrow\ker\left(  \tau\right)  \longrightarrow G_{L}\overset{\tau
}{\longrightarrow}\mathbb{Z}\left[  \tfrac{1}{k}\right]  \longrightarrow0
\label{eqExt.pound}%
\end{equation}
is non-split. Let $J_{L}$ be given as usual (see (\ref{eqExt.9})), and let
$G_{L}$ be the inductive limit from%
\begin{equation}
\mathbb{Z}^{N}\subset J_{L}^{-1}\left(  \mathbb{Z}^{N}\right)  \subset
J_{L}^{-2}\left(  \mathbb{Z}^{N}\right)  \subset\cdots.
\label{eqExt.poundpound}%
\end{equation}

For the pair (\ref{eqExtNew.pound}) we show in \cite{BJKR00}, or, more
generically, in Proposition \ref{ProApp.Cla2}, that (\ref{eqExt.pound}) is the
following exact sequence:%
\begin{equation}
0\longrightarrow\mathbb{Z}\left[  \tfrac{1}{2}\right]  \xrightarrow
{z\mapsto\left( z,-8z\right) }\mathbb{Z}\left[  \tfrac{1}{2}\right]
\times\mathbb{Z}\left[  \tfrac{1}{2}\right]  \xrightarrow{\tau=\ip{\,\cdot
\,}{\left( 1,1/8\right) }}\mathbb{Z}\left[  \tfrac{1}{2}\right]
\longrightarrow0. \label{eqExtNew.poundpound}%
\end{equation}
This is the zero element of $\operatorname*{Ext}\left(  \mathbb{Z}\left[
\frac{1}{2}\right]  ,\mathbb{Z}\left[  \frac{1}{2}\right]  \right)  $: an
injection $\psi\colon\mathbb{Z}\left[  \frac{1}{2}\right]  \rightarrow
\mathbb{Z}\left[  \frac{1}{2}\right]  \times\mathbb{Z}\left[  \frac{1}%
{2}\right]  $ may be defined by%
\begin{equation}
\psi\left(  u\right)  :=\left(  \frac{7}{8}u,u\right)  ,\qquad u\in
\mathbb{Z}\left[  \tfrac{1}{2}\right]  . \label{eqExtNew.lift}%
\end{equation}
Then clearly $\tau\left(  \psi\left(  u\right)  \right)  =u$, so $\psi$
defines a section, and (\ref{eqExtNew.poundpound}) splits.

The inductive limit $G_{L}$ from \textup{(\ref{eqExt.poundpound}),} and
$\tau\left(  G_{L}\right)  =\mathbb{Z}\left[  \frac{1}{k}\right]  $, $k>1$, in
general define a nontrivial element of $\operatorname*{Ext}$, i.e.,
(\ref{eqExt.pound}) is non-split in general. It is split if and only if there
is an element $g\in G_{L}$ such that $\left(  k^{-i}\right)  g\in G_{L}$,
$\forall\,i\in\mathbb{N}$, and $\tau\left(  g\right)  =1$.

\setcounter{figurelink}{\value{figure}}

\chapter{\label{Sesquilabial}Scaling and non-isomorphism}

\setcounter{figure}{\value{figurelink}} In this chapter we introduce a number
$\tau\left(  v\right)  $, and prove in Theorem \ref{theorem10.10} and
Corollary \ref{corollary10.19} that it can be used to establish
non-isomorphism for classes of algebras where the basic invariants in Theorem
\ref{CorCyc.9} are the same..

Let $\lambda\in\mathbb{R}_{+}$, and let
\[
\mathcal{L}\left(  \lambda\right)  =\{J\mid J\text{ is of the form
(\ref{eqCycNew.1}) with Perron--Frobenius eigenvalue }\lambda\}.
\]
In particular, the standard matrix $J_{m}$ is in $\mathcal{L}\left(
\lambda\right)  $ if and only if
\begin{equation}
\lambda^{N}-m_{1}\lambda^{N-1}-m_{2}\lambda^{N-2}-\dots-m_{N-1}\lambda
-m_{N}=0. \label{eqExt.lambdam}%
\end{equation}
The admissible numbers $\lambda$ must therefore be algebraic. These algebraic
integers $\lambda$ may be specified further; see, e.g., \cite{Han92},
\cite{Ren57}, and \cite{Par60} for more details on this point.

We are not restricting the size $N\times N$ of the matrices $J$ in
$\mathcal{L}\left(  \lambda\right)  $.

Our main result, Theorem \ref{theorem10.10}, in this chapter, is that
$\tau\left(  v\right)  =\ip{\alpha}{v}$, introduced in (\ref{eqajla}%
)--(\ref{eqjvlv}), can be used to show non-isomorphism of a class of cases in
$\mathcal{L}\left(  \lambda\right)  $ when $\lambda\in\mathbb{Z}_{+}$,
$\lambda>1$.

We consider matrices $J=J_{m}=J_{L}$ having the form\label{LOSJ_5}%
\begin{equation}
J_{m}=%
\begin{pmatrix}
\vphantom{\vdots}m_{1} & 1 & 0 & \cdots & 0 & 0\\
\vphantom{\vdots}m_{2} & 0 & 1 & \cdots & 0 & 0\\
\vphantom{\vdots}m_{3} & 0 & 0 & \ddots & 0 & 0\\
\vdots & \vdots &  & \ddots & \ddots & \vdots\\
\vphantom{\vdots}m_{N-1} & 0 & 0 &  & 0 & 1\\
\vphantom{\vdots}m_{N} & 0 & 0 & \cdots & 0 & 0
\end{pmatrix}
, \label{eqCycNew.1}%
\end{equation}
with $m_{i}\in\mathbb{Z}$, $m_{i}\geq0$, $m_{N}>0$, and satisfying the further
requirement that for some $k\in\mathbb{N}$, $J^{k}$ has only positive entries
(equivalently, $\gcd\left\{  i\mid m_{i}\neq0\right\}  =1$). Non-unimodularity
means $m_{N}>1$.

Recall from (\ref{eqRemNewBis.17}) that the vector $\alpha=\alpha_{\lambda
}=\left(  1,\lambda^{-1},\dots,\lambda^{-\left(  N-1\right)  }\right)
$\label{LOSalphaeigenvector_3}
satisfies\label{LOSPerronFrobeniuslefteigenvector_5}%
\begin{equation}
\alpha J=\lambda\alpha, \label{eqajla}%
\end{equation}
and also there is a unique $v\in\mathbb{Z}\left[  \lambda\right]  ^{N}$ such
that\label{LOSPerronFrobeniusrighteigenvector_3}%
\begin{equation}
Jv=\lambda v\text{\qquad and\qquad}v_{1}=1. \label{eqjvlv}%
\end{equation}
An explicit expression for $v$ is given in (\ref{eqSubNew.pound}). When $J$ is
given, let $\mathfrak{A}_{J}$ be the corresponding AF-algebra.

Our list of invariants, so far, cannot separate the AF-isomorphism classes
corresponding to the following pair of examples (Examples A and B) where:%
\begin{equation}%
\begin{pmatrix}
m_{1}\\
m_{2}\\
m_{3}%
\end{pmatrix}
_{A}=%
\begin{pmatrix}
0\\
3\\
2
\end{pmatrix}
\text{ and }%
\begin{pmatrix}
m_{1}\\
m_{2}\\
m_{3}%
\end{pmatrix}
_{B}=%
\begin{pmatrix}
0\\
15\\
4
\end{pmatrix}
,\qquad\lambda_{A}=2\text{ and }\lambda_{B}=4. \label{eqExaABspecs}%
\end{equation}%
\begin{equation}%
\begin{tabular}
[c]{c}%
The stabilized Bratteli diagrams:\\
\raisebox{57bp}{Example A:}$\begin{picture}(180,114)(0,-6) \put
(0,0){\includegraphics
[bb=135 68 297 170,clip,height=102bp,width=162bp]{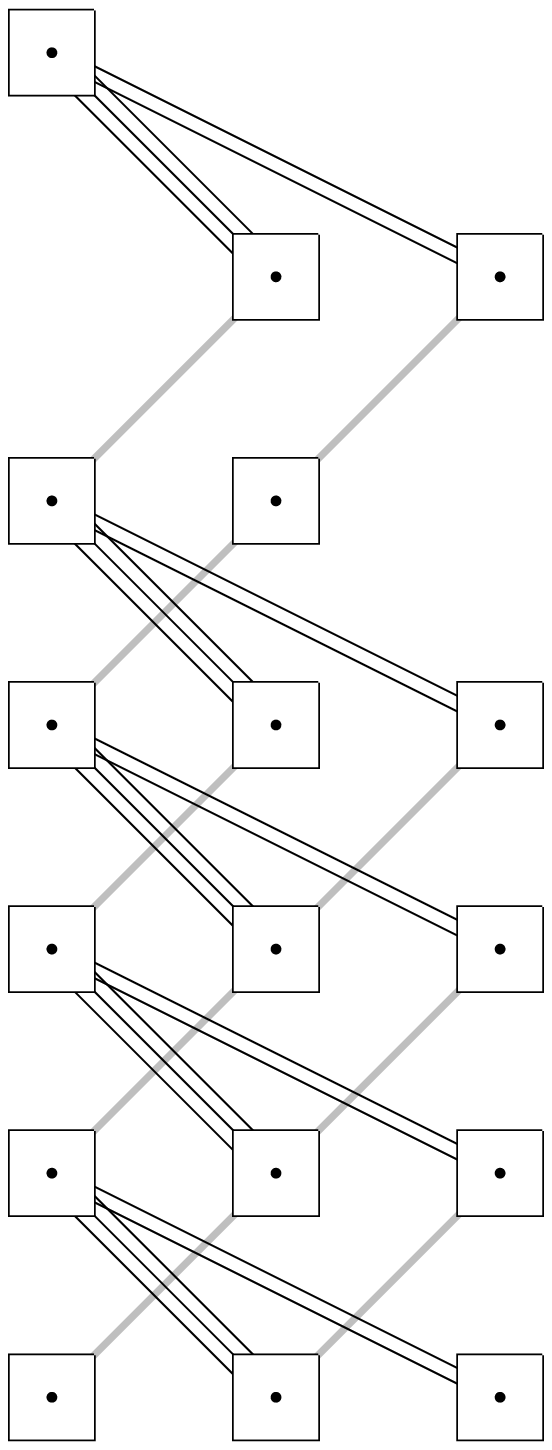}} \end{picture}%
$\\
\raisebox{57bp}{Example B:}$\begin{picture}(180,114)(0,-6) \put
(0,0){\includegraphics
[bb=135 68 297 170,clip,height=102bp,width=162bp]{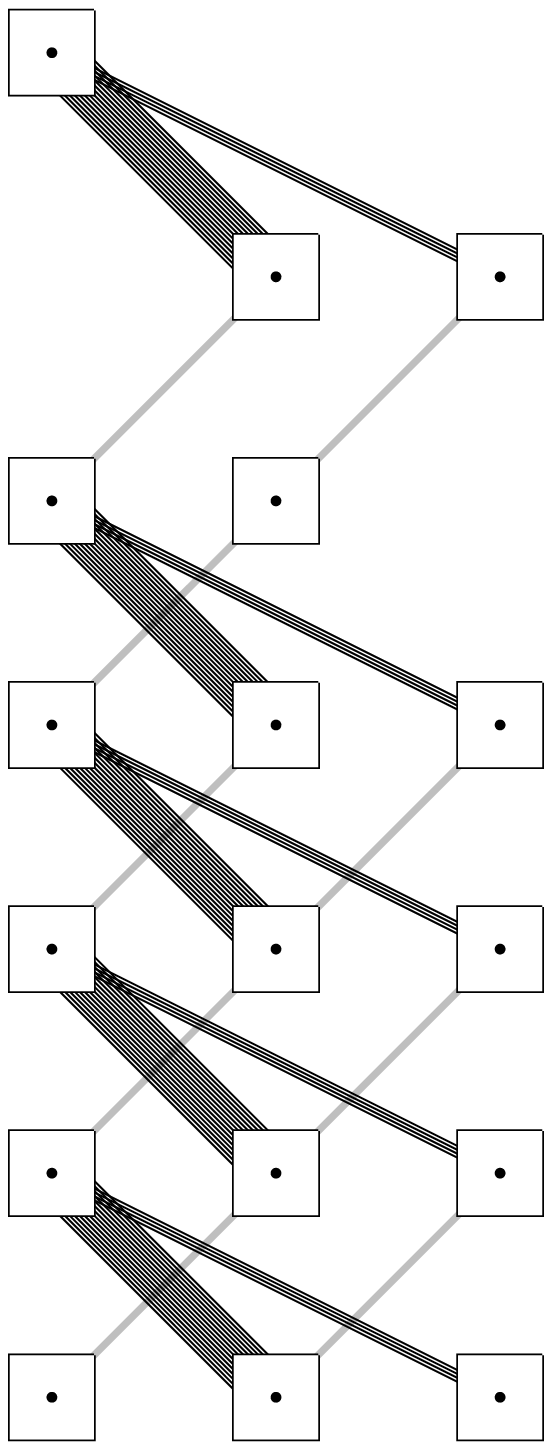}} \end{picture}$%
\end{tabular}
\label{eqTableBratDiag}%
\end{equation}

\begin{table}[ptb]
\caption{Parameters for Examples A and B in (\ref{eqExaABspecs}%
)--(\ref{eqTableBratDiag}).}%
\label{TableCycpoundprime}%
\begin{tabular}
[c]{ccccccc}%
& $\lambda$ & $m_{N}$ & $R_{D}$ & $\ker\left(  \tau\right)  $ & $F=K_{0}%
\diagup\mathbb{Z}^{3}$ & $\ip{\alpha}{v}\vphantom{\left( \begin{smallmatrix}%
0\\ 0\\ 0 \end{smallmatrix} \right) }$\\
Example A: & $2$ & $2$ & $2$ & $\mathbb{Z}^{2}$ & $\mathbb{Z}\left[  \frac
{1}{2}\right]  \diagup\mathbb{Z}$ & $\frac{9}{4}\vphantom{\left
( \begin{smallmatrix}0\\ 0\\ 0 \end{smallmatrix}\right) }$\\
Example B: & $4$ & $4$ & $4$ & $\mathbb{Z}^{2}$ & $\mathbb{Z}\left[  \frac
{1}{2}\right]  \diagup\mathbb{Z}$ & $\frac{33}{16}\vphantom{\left
( \begin{smallmatrix}0\\ 0\\ 0 \end{smallmatrix} \right) }$%
\end{tabular}
\end{table}

For both these examples, the basic invariants in Chapter \ref{Pediferient}
have the values $N=3$, $D=1$, $\operatorname*{Prim}\left(  m_{N}\right)
=\left\{  2\right\}  $, $\operatorname*{Prim}\left(  R_{D}\right)  =\left\{
2\right\}  $, $\operatorname*{Prim}\left(  Q_{N-D}\right)  =\varnothing$.
These invariants do not directly separate the isomorphism classes of the
examples. However, since one Perron--Frobenius eigenvalue is a power of the
other, we will show that $\tau\left(  v\right)  $ can be used to check
non-isomorphism of the two AF-algebras $\mathfrak{A}_{A}$ and $\mathfrak
{A}_{B}$. This is shown for this specific example in Observation
\ref{ObsCyc.3}, and more generally in Theorem \ref{theorem10.10}.

It is easy to check that both examples have $\ker\left(  \tau\right)
\cong\mathbb{Z}^{2}$, and $\tau\left(  K_{0}\right)  \cong\mathbb{Z}\left[
\frac{1}{2}\right]  $. Strictly speaking, $\tau\left(  K_{0}\right)  $ is
$\mathbb{Z}\left[  \frac{1}{2}\right]  $ for Example A, and $\mathbb{Z}\left[
\frac{1}{4}\right]  $ for Example B; but $\mathbb{Z}\left[  \frac{1}%
{2}\right]  =\mathbb{Z}\left[  \frac{1}{4}\right]  $ with the natural
isomorphism specified by%
\[
\frac{1}{4^{i}}\longmapsto\frac{1}{2^{2i}}.
\]
Hence, both of the examples are characterized as elements of
$\operatorname*{Ext}\left(  \mathbb{Z}\left[  \frac{1}{2}\right]
,\mathbb{Z}^{2}\right)  $, in the usual manner. Let $G_{A}$ and $G_{B}$ be the
respective $K_{0}$-groups. The rank of each group is clearly $N=3$.

The next Observation illustrates the previous remarks about $m_{N}$. Let
$F_{A}:=G_{A}\diagup\mathbb{Z}^{3}$,\label{LOSFL_1} and similarly for $F_{B}%
$.\renewcommand
{\theenumi}{\alph{enumi}}

\begin{observation}
\label{ObsCyc.1}$G_{A}\diagup\mathbb{Z}^{3}\cong G_{B}\diagup\mathbb{Z}%
^{3}\cong\mathbb{Z}\left[  \frac{1}{2}\right]  \diagup\mathbb{Z}$.
\end{observation}

\begin{proof}
The respective quotient groups have the following generators:%
\begin{align*}
F_{A}\colon x_{i}  &  =J_{A}^{-i}e_{3}\mod{\mathbb{Z}^{3}}\\%
\intertext{and}%
F_{B}\colon y_{i}  &  =J_{B}^{-i}e_{3}\mod{\mathbb{Z}^{3}},
\end{align*}
and a use of Lemma \ref{LemSubNew.14} yields:

\begin{enumerate}
\item \label{ObsCyc.1proof(1)}$2x_{1}=0$, $2x_{i+1}=x_{i}$, $i\in\mathbb{N}$, and

\item \label{ObsCyc.1proof(2)}$4y_{1}=0$, $4y_{i+1}=y_{i}$, $i\in\mathbb{N}.$
\end{enumerate}

\noindent Hence $G_{A}\diagup\mathbb{Z}^{3}\cong\mathbb{Z}\left[  \frac{1}%
{2}\right]  \diagup\mathbb{Z}$, and $G_{B}\diagup\mathbb{Z}^{3}\cong
\mathbb{Z}\left[  \frac{1}{4}\right]  \diagup\mathbb{Z}=\mathbb{Z}\left[
\frac{1}{2}\right]  \diagup\mathbb{Z}$.
\end{proof}

The crucial property derived from (\ref{ObsCyc.1proof(1)}%
)--(\ref{ObsCyc.1proof(2)}) above is not really Observation \ref{ObsCyc.1},
but that scaling by a power of $2$ (in Example A) and $4$ (in B) determines a
\emph{filtration} of $G_{A}=\bigcup_{n\geq0}J_{A}^{-n}\mathbb{Z}^{3}$.
Specifically, let $\pi_{A}\colon G_{A}\rightarrow F_{A}=G_{A}\diagup
\mathbb{Z}^{3}$ be the quotient mapping. If elements of $F_{A}$ are
represented as $\left(  i_{1},i_{2},\dots,i_{n-1},1\right)  $, $i_{j}%
\in\mathbb{Z}_{2}$, we write%
\begin{multline*}
g_{n}=x+i_{1}J_{A}^{-1}e_{3}+i_{2}J_{A}^{-2}e_{3}+\dots\\
+i_{n-1}J_{A}^{-n-1}e_{3}+1\cdot J_{A}^{-n}e_{3}\text{\qquad for }%
x\in\mathbb{Z}^{3},\;i_{j}\in\left\{  0,1\right\}  .
\end{multline*}
Recall from Corollary \ref{CorCyc.4half} that this representation is unique.
Then%
\[
\pi_{A}\left(  g_{n}\right)  =\left(  i_{1}^{\prime},i_{2}^{\prime}%
,\dots,i_{n-1}^{\prime},1\right)  ,\qquad g_{n}\in G_{n}\left(  A\right)
=J_{A}^{-n}\mathbb{Z}^{3},
\]
and we note that%
\[
\pi_{A}\left(  2g_{n}\right)  =(\underset{n-1\text{ places}}{\underbrace
{i_{2}^{\prime},i_{3}^{\prime},\dots,i_{n-1}^{\prime},1})}=\pi_{A}\left(
J_{A}g_{n}\right)  .
\]
A similar remark applies to $F_{B}=G_{B}\diagup\mathbb{Z}^{3}$, but there the
scaling is by a power of $4$. We have proved the following:%
\[
\text{An element }g\text{ of }G_{A}\text{ is in }G_{n}\left(  A\right)  \text{
if and only if }2^{n}g\in\mathbb{Z}^{3}.
\]
Similarly, when $g\in G_{B}$, then $g\in G_{n}\left(  B\right)  =J_{B}%
^{-n}\mathbb{Z}^{3}$ if and only if $4^{n}g\in\mathbb{Z}^{3}$. We shall need
this in the proof of Observation \ref{ObsCyc.3} below. (See Figure
\ref{FigCycpoundpound}.)

In summary, both examples have $N=3$, $D=1$, and $\tau\left(  K_{0}\left(
\mathfrak{A}\right)  \right)  =\mathbb{Z}\left[  \frac{1}{2}\right]  $, and
the other data are as in Table \ref{TableCycpoundprime} (above).

The proof of non-isomorphism for A and B uses the fact that $%
\ip{\alpha_{A}}{v_{A}}%
$ and $%
\ip{\alpha_{B}}{v_{B}}%
$ have different prime factors than $2$ for their numerators.

\begin{observation}
\label{ObsCyc.3}Examples A and B correspond to non-isomorphic AF-algebras
$\mathfrak{A}_{A}$ and $\mathfrak{A}_{B}$.
\end{observation}

The argument proves that there is not even a nontrivial homomorphism
$\theta\colon G_{B}\rightarrow G_{A}$ which makes the diagrams in
(\ref{eqCycpound}) and Figure \ref{FigCycpoundpound} (below) commutative.%

\begin{gather}
\text{The nonexistent isomorphism }\theta\colon G_{B}\rightarrow G_{A}%
\text{:}\nonumber\\%
\begin{array}
[c]{ccccccccc}%
0 & \longrightarrow & \ker\left(  \tau_{B}\right)  & \overset{i_{B}%
}{\hooklongrightarrow} & G_{B} & \overset{\tau_{B}}{\longrightarrow} &
\tau_{B}\left(  G_{B}\right)  & \longrightarrow & 0\\
&  & \vrule depth6.5pt width0.275pt height11.5pt\rlap{$\scriptstyle\;\cong$} &
&  &  & \vrule depth5pt width0.275pt height11.5pt\rlap{$\scriptstyle\;\cong$}%
&  & \\
&  & \mathbb{Z}^{2} &  & \llap{$\scriptstyle
\ncong\;$}\smash{\vcenter{\hrule width0.275pt height46pt}\makebox
[0pt] {\hss\raisebox{-19pt}{$\downarrow$}\hss}}\rlap{$\scriptstyle
\;\theta$} &  & \mathbb{Z}\left[  \frac{1}{2}\right]  &  & \\
&  & \vrule depth6.5pt width0.275pt height11.5pt\rlap{$\scriptstyle\;\cong$} &
&  &  & \vrule depth6.5pt width0.275pt height9.5pt\rlap{$\scriptstyle\;\cong$}%
&  & \\
0 & \longrightarrow & \ker\left(  \tau_{A}\right)  & \overset{i_{A}%
}{\hooklongrightarrow} & G_{A} & \overset{\tau_{A}}{\longrightarrow} &
\tau_{A}\left(  G_{A}\right)  & \longrightarrow & 0
\end{array}
\label{eqCycpound}%
\end{gather}

\begin{figure}[ptb]
\newcounter{frogtag} \setlength{\unitlength}{25.14417bp} \begin{picture}%
(10.9369,16.2492)(-2.8832,-14.8746)
\put(-2.8832,-14.8746){\includegraphics
[bb=66 0 366 432,height=396bp,width=275bp]{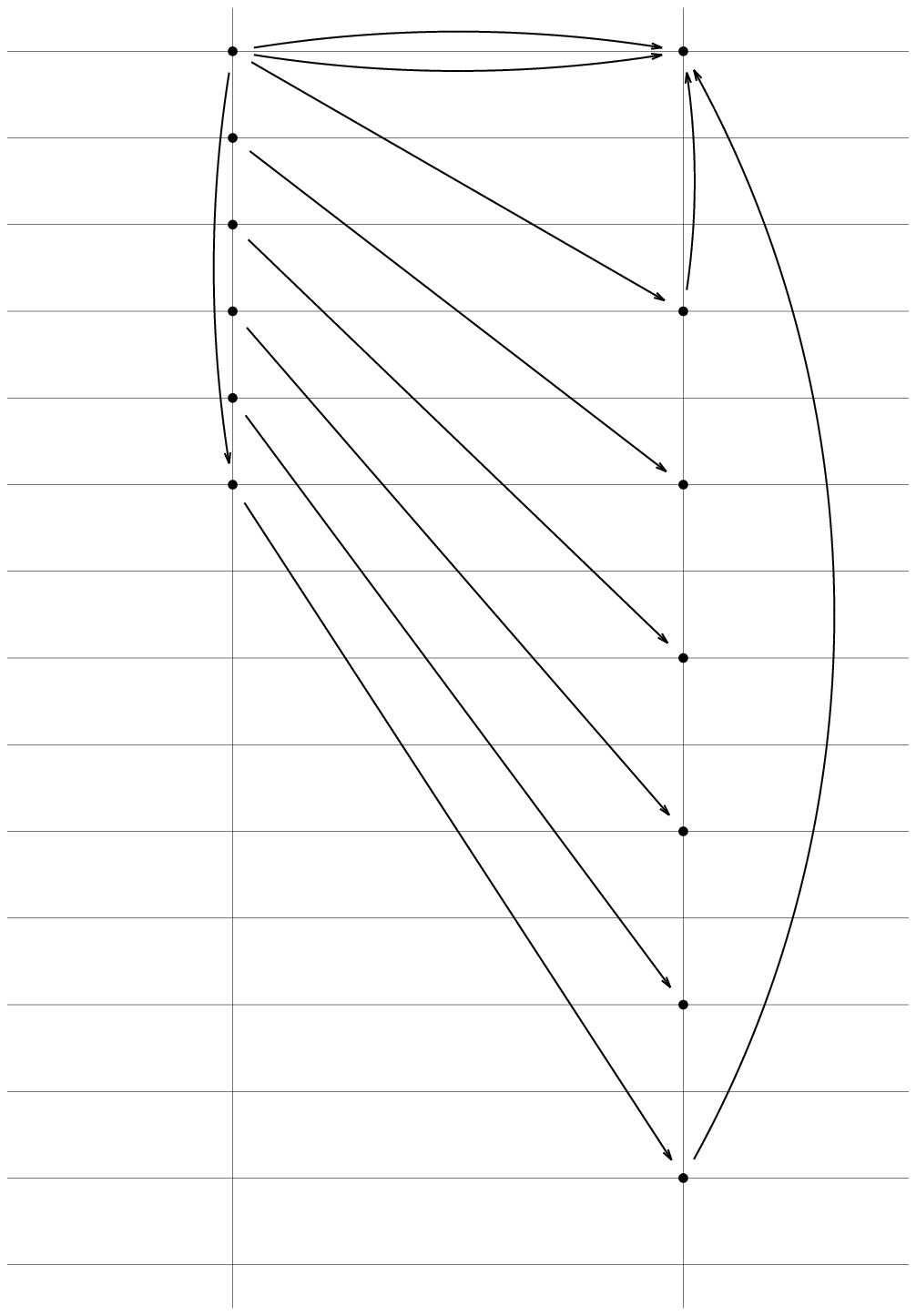}}
\put(0,0.573){\makebox(0,0)[b]{$G_{n}\left( B\right) $}}
\put(5.19615,0.573){\makebox(0,0)[b]{$G_{n}\left( A\right) $}}
\put(-0.3457,-2.5){\makebox(0,0)[r]{$J_{B}^{-n}$}}
\put(7.0109,-6.5){\makebox(0,0)[l]{$J_{A}^{k+2n}$}}
\put(5.4004,-1.5){\makebox(0,0)[l]{$J_{A}^{k}$}}
\put(2.59808,-0.3){\makebox(0,0)[tl]{$\psi_{n}$}}
\put(2.59808,0.3){\makebox(0,0)[bl]{$\psi_{0}$}}
\multiput(2.59808,-1.427)(0,-1.5){6}{\makebox(0,0)[bl]{$\theta$}}
\multiput(-2.59808,0.073)(0,-1){15}{\makebox(0,0)[bl]
{$n=\arabic{frogtag}$\addtocounter{frogtag}{1}}}
\end{picture}
\caption{The (nonexistent) matrices $\psi_{n}$ in the proof of Observation
\ref{ObsCyc.3} (example with $k=3$).}%
\label{FigCycpoundpound}%
\end{figure}

\begin{proof}
It is enough to show that $G_{A}$ and $G_{B}$ represent different elements of
\linebreak $\operatorname*{Ext}\left(  \mathbb{Z}\left[  \frac{1}{2}\right]
,\mathbb{Z}^{2}\right)  $. This can be done by recursion, and use of the
relations (\ref{ObsCyc.1proof(1)})--(\ref{ObsCyc.1proof(2)}). Alternatively,
it can be checked directly by the argument from the proof of Theorem
\ref{theorem10.10} below that the two $\operatorname*{Ext}$-elements $G_{A}$
and $G_{B}$ are different in $\operatorname*{Ext}\left(  \mathbb{Z}\left[
\frac{1}{2}\right]  ,\mathbb{Z}^{2}\right)  $. Both arguments are essentially
based on the $%
\ip{\alpha}{v}%
$-number, even though $\lambda_{A}\neq\lambda_{B}$. In the present case,
$\lambda_{A}^{2}=4=\lambda_{B}^{{}}$, which is good enough. The present
argument is essentially a ``baby'' version of the argument in the rest of this chapter.

We sketch the details. It is a proof by contradiction. Suppose $\theta$ were
an isomorphism of the ordered $K_{0}$-groups, say $\theta\colon G_{B}%
\rightarrow G_{A}$, which made them the same element of $\operatorname*{Ext}%
\left(  \mathbb{Z}\left[  \frac{1}{2}\right]  ,\mathbb{Z}^{2}\right)  $. Since
$G_{A}=\bigcup_{n\geq0}J_{A}^{-n}\mathbb{Z}^{3}$, there is a $k$ such that
$\theta\left(  \mathbb{Z}^{3}\right)  \subset J_{A}^{-k}\mathbb{Z}^{3}$. We
then claim that $\theta\left(  J_{B}^{-n}\mathbb{Z}^{3}\right)  \subset
J_{A}^{-\left(  k+2n\right)  }\mathbb{Z}^{3}$ for all $n$. The argument for
this is based on properties (\ref{ObsCyc.1proof(1)})--(\ref{ObsCyc.1proof(2)})
for the generators: Let $x\in\mathbb{Z}^{3}$. To verify that $\theta\left(
J_{B}^{-n}x\right)  \in J_{A}^{-k-2n}\mathbb{Z}^{3}$, we must check that
$2^{2n+k}\theta\left(  J_{B}^{-n}x\right)  \in\mathbb{Z}^{3}$. This holds
since $2^{2n+k}\theta\left(  J_{B}^{-n}x\right)  =2^{k}\theta\left(
4^{n}J_{B}^{-n}x\right)  \in2^{k}\theta\left(  \mathbb{Z}^{3}\right)
\subset2^{k}J_{A}^{-k}\mathbb{Z}^{3}\subset\mathbb{Z}^{3}$. Hence
$2^{2n+k}\theta\left(  J_{B}^{-n}x\right)  \in\mathbb{Z}^{3}$, and therefore
$\theta\left(  J_{B}^{-n}x\right)  \in J_{A}^{-\left(  k+2n\right)
}\mathbb{Z}^{3}$ as claimed.

These maps may be represented with matrices $\psi_{n}\in M_{3}\left(
\mathbb{Z}\right)  $ as follows:
\[
\theta\left(  J_{B}^{-n}x\right)  =J_{A}^{-\left(  k+2n\right)  }\psi
_{n}\left(  x\right)  ,\qquad x\in\mathbb{Z}^{3},\;n=0,1,2,\dots,
\]
with the consistency conditions%
\[
\psi_{n}=J_{A}^{2n}\psi_{0}J_{B}^{-n},\qquad n\in\mathbb{N}.
\]
Thus $\psi_{0}=J_{A}^{k}\theta|_{\mathbb{Z}^{3}}$, and $\psi_{n}=J_{A}%
^{k+2n}\theta J_{B}^{-n}|_{\mathbb{Z}^{3}}$ (see Figure \ref{FigCycpoundpound}%
). This defines the sequence $\psi_{n}$ as a sequence of linear endomorphisms
of $\mathbb{Z}^{3}$, and so each $\psi_{n}$ is represented by a matrix in
$M_{3}\left(  \mathbb{Z}\right)  $. That turns out to be very restrictive. It
is not satisfied for $\psi_{0}=I_{3}$. In fact, even $J_{A}^{2}J_{B}^{-1}$ has
a non-integral entry. (The matrices $\psi_{n}$ play the role of the
intertwiners $A_{n}$ in the diagram (\ref{eqbeforeRemNewTer.42}), with
$A_{1}=\psi_{0}$, $A_{2}=\psi_{1}$, etc., but in the reasoning here positivity
does not play a role.)

Let $v_{A}=\left(
\begin{smallmatrix}
1\\
2\\
1
\end{smallmatrix}
\right)  $ and $v_{B}=\left(
\begin{smallmatrix}
1\\
4\\
1
\end{smallmatrix}
\right)  $ be the normalized Perron--Frobenius eigenvectors:%
\[
J_{A}v_{A}=2v_{A},\qquad J_{B}v_{B}=4v_{B}.
\]
Hence%
\[
\psi_{n}\left(  v_{B}\right)  =J_{A}^{2n}\psi_{0}J_{B}^{-n}v_{B}=4^{-n}%
J_{A}^{2n}\psi_{0}\left(  v_{B}\right)  \underset{n\rightarrow\infty
}{\longrightarrow}\frac{\tau_{A}\left(  \psi_{0}v_{B}\right)  }{\tau
_{A}\left(  v_{A}\right)  }v_{A}%
\]
by the Perron--Frobenius theorem. But%
\[
\tau_{A}\left(  \psi_{0}v_{B}\right)  =2^{k}\tau_{B}\left(  v_{B}\right)
\]
since $\theta$ preserves the normalized trace; see (\ref{eqCycpound}).

Hence, by taking $n$ large, $\psi_{n}\left(  v_{B}\right)  $ ($\in
\mathbb{Z}^{3}$) will be arbitrarily close to $2^{k}\frac{\tau_{B}\left(
v_{B}\right)  }{\tau_{A}\left(  v_{A}\right)  }v_{A}=2^{k}\frac{11}%
{12}\left(
\begin{smallmatrix}
1\\
2\\
1
\end{smallmatrix}
\right)  $, where we used the numbers from the last column in Table
\ref{TableCycpoundprime}. If $\left(  Q_{ij}^{\left(  n\right)  }\right)
_{i,j=1}^{3}$ denote the matrix entries of $\psi_{n}$, then $Q_{ij}^{\left(
n\right)  }\in\mathbb{Z}$, and $\psi_{n}\left(  v_{B}\right)  _{1}%
=Q_{11}^{\left(  n\right)  }+4Q_{12}^{\left(  n\right)  }+Q_{13}^{\left(
n\right)  }$. But they are integers, so there is an $n_{0}\in\mathbb{N}$ such
that $\mathbb{Z}\ni\psi_{n}\left(  v_{B}\right)  _{1}=\frac{11\cdot2^{k}}{12}$
for all $n\geq n_{0}$. Since $3$ does not divide $11$, this is a
contradiction. We have proved that the $\operatorname*{Ext}$-elements are
different as claimed.
\end{proof}

\begin{remark}
\label{RemCyc.4}We will get back to the idea of defining $\theta$ by a matrix
in $\operatorname*{GL}\left(  N,\mathbb{R}\right)  $ in a more systematic way
later, in Proposition \textup{\ref{prop10.7}.}
\end{remark}

Let $J$ be a matrix specified as in (\ref{eqCycNew.1}) and let $G=G_{J}%
$\label{LOSGJ} be the corresponding inductive limit group (see Chapter
\ref{Rem}). Recall $m_{N}=\left|  \det J\right|  $. Setting $G_{i}%
:=J^{-i}\mathbb{Z}^{N}$,\label{LOSGi} we shall use the homomorphisms,
$g\mapsto m_{N}^{k}g$, in localizing the scaling $G_{i}$, $i=1,2,\dots$. It is
immediate from the proof of Claim \ref{ClaCyc.2} that $m_{N}J^{-1}$ has
integral entries, i.e., is in $M_{N}\left(  \mathbb{Z}\right)  $. Since
$m_{N}^{i}J_{{}}^{-i}=\left(  m_{N}J^{-1}\right)  ^{i}$, this is also true for
the iterations. We have proved the following implication:%
\begin{equation}
g\in G_{i}\Longrightarrow g\in G\text{\qquad and\qquad}m_{N}^{i}g\in
\mathbb{Z}^{N}. \label{eqSesquilabial.pound}%
\end{equation}

It will be useful that a scaled version of the converse also is true.

\begin{proposition}
\label{prop10.5}\renewcommand{\theenumi}{\roman{enumi}}Let $J$ be given as
specified in \textup{(\ref{eqCycNew.1}).} Then there exists a $p\in\mathbb{N}$
such that the following implication \textup{(\ref{prop10.5(1)})} $\Rightarrow$
\textup{(\ref{prop10.5(2)})} holds for $g\in G_{J}$:
\end{proposition}

\begin{enumerate}
\item \label{prop10.5(1)}$m_{N}g\in\mathbb{Z}^{N}$

$\quad\Downarrow$

\item \label{prop10.5(2)}$g\in J^{-p}(\mathbb{Z}^{N})$
\end{enumerate}

\noindent Furthermore, for the same $p\in\mathbb{N}$ we have the following
implications for $g\in G_{J}$:

\begin{enumerate}
\setcounter{enumi}{2}

\item \label{prop10.5(3)}$m_{N}^{i}g\in\mathbb{Z}^{N}$

$\quad\Downarrow$

\item \label{prop10.5(4)}$g\in J^{-ip}(\mathbb{Z}^{N})$ for $i=1,2,3,\ldots$.
\end{enumerate}

\begin{proof}
If $m_{N}g\in\mathbb{Z}^{N}$, then $g\in\frac{1}{m_{N}}\mathbb{Z}^{N}$, and
hence $g$ has the form
\[
g=\frac{(k_{1},\ldots,k_{N})}{m_{N}}+m
\]
where $k_{i}\in\{0,1,\ldots,m_{N}-1\}$ for $i=1,\ldots,N$, and $m\in
\mathbb{Z}^{N}$. But $G_{J}$ can contain at most $m_{N}^{N}$ elements of the
form $(k_{1},\ldots,k_{N})\left/  m_{N}\right.  $ where $k_{i}\in
\{0,1,\ldots,m_{N}-1\}$, and since this number is finite and
\[
G_{J}=\bigcup_{n}J^{-n}\mathbb{Z}^{N}%
\]
is an increasing union, it follows that there is a $p\in\mathbb{N}$ such that
all these elements are contained in $J^{-p}\mathbb{Z}^{N}$. But as
$\mathbb{Z}^{N}\subseteq J^{-p}\mathbb{Z}^{N}$, the implication
(\ref{prop10.5(1)}) $\Rightarrow$ (\ref{prop10.5(2)}) follows.

Next, choose $p\in\mathbb{N}$ such that (\ref{prop10.5(1)}) $\Rightarrow$
(\ref{prop10.5(2)}) holds. We prove by induction with respect to $i$ that
(\ref{prop10.5(3)}) $\Rightarrow$ (\ref{prop10.5(4)}) holds. $i=1$ is
(\ref{prop10.5(1)}) $\Rightarrow$ (\ref{prop10.5(2)}), so assume
(\ref{prop10.5(3)}) $\Rightarrow$ (\ref{prop10.5(4)}) holds for $i-1$, and
assume that $g\in G_{J}$ and
\[
m_{N}^{i}g=m_{N}^{i-1}m_{N}g\in\mathbb{Z}^{N}\text{.}%
\]
By the induction hypothesis, we then have
\[
m_{N}g\in J^{-p(i-1)}\mathbb{Z}^{N}\text{.}%
\]
But applying $J^{p(i-1)}$ to both sides, we have
\[
m_{N}J^{p(i-1)}g\in\mathbb{Z}^{N}\text{.}%
\]
Applying the case $i=1$, one obtains
\[
J^{p(i-1)}g\in J^{-p}\mathbb{Z}^{N}\text{,}%
\]
and applying $J^{-p(i-1)}$ to both sides
\[
g\in J^{-pi}\mathbb{Z}^{N}%
\]
and this proves (\ref{prop10.5(3)}) $\Rightarrow$ (\ref{prop10.5(4)}).
\end{proof}

Note that the implication (\ref{prop10.5(1)}) $\Rightarrow$ (\ref{prop10.5(2)}%
) holds if there is a $p\in\mathbb{N}$, and an $E\in M_{N}\left(
\mathbb{Z}\right)  $, such that%
\begin{equation}
J^{p}=m_{N}E. \label{eqSesquilabial.5}%
\end{equation}
If (\ref{prop10.5(3)}) $\Rightarrow$ (\ref{prop10.5(4)}) holds for \emph{all}
$g\in\mathbb{Q}^{N}$, then $J^{ip}\left(  \mathbb{Z}^{N}\right)  \subset
m_{N}^{i}\mathbb{Z}_{{}}^{N}$ for all $i$, so \textup{(\ref{eqSesquilabial.5}%
)} is valid. Thus \textup{(\ref{eqSesquilabial.5})} is stronger than
(\ref{prop10.5(3)}) $\Rightarrow$ (\ref{prop10.5(4)})\textup{.}

\begin{definition}
\label{definition10.6}Suppose $J$ is a matrix of the form
\textup{(\ref{eqCycNew.1}).} We say that $J$ has \emph{scaling degree} $\leq
p$, and write $\deg(J)\leq p$\label{LOSdeg_2} if there exists an $n_{0}%
\in\mathbb{N}$ such that, for $g\in G_{J}$,
\[
m_{N}^{i}g\in\mathbb{Z}^{N}\Longrightarrow g\in J^{-\lfloor ip\rfloor-n_{0}%
}\mathbb{Z}^{N}\text{.}%
\]
More generally, if $m$ is a positive integer containing exactly the same prime
factors as $m_{N}$, we say that $J$ has $m$\emph{-scaling degree} $\leq
p=p^{(m)}$ and write $m$-$\deg(J)\leq p$\label{LOSmdeg_2} if there exists an
$n_{0}\in\mathbb{N}$ such that, for $g\in G_{J}$%
\[
m^{i}g\in\mathbb{Z}^{N}\Longrightarrow g\in J^{-\lfloor ip\rfloor-n_{0}%
}\mathbb{Z}^{N}\text{.}%
\]
\end{definition}

Note that as $m$ contains the same prime factors as $m_{N}$, $J$ has a finite
$m$-scaling degree if and only if it has a finite scaling degree, and the last
is true by Proposition \ref{prop10.5}.

We note \emph{en passant} that the above remark implies the following
corollary, which is true whether the Perron--Frobenius eigenvalue $\lambda$ is
rational or not. (But if $\lambda$ is rational under the conditions in the
corollary, it follows from the characteristic equation that
$\operatorname*{Prim}(\lambda)\subseteq\operatorname*{Prim}\left(
m_{N}\right)  $).

\begin{corollary}
\label{corollary10.A}Let $J$ be a matrix of the form \textup{(\ref{eqCycNew.1}%
),} and assume that each $m_{i}$ is either $0$ or contains all the prime
factors of $m_{N}$. It follows that
\[
G=\mathbb{Z}\left[  \tfrac{1}{m_{N}}\right]  ^{N}%
\]
when $G$ is identified as a subgroup of $\mathbb{Q}^{N}$ by
\textup{(\ref{eqRemNewBis.6}).}
\end{corollary}

\begin{proof}
Since $J^{-1}=\frac{1}{m_{N}}E$ where $E$ is a matrix with integer
coefficients, it is clear from (\ref{eqRemNewBis.6}) that
\[
G\subseteq\mathbb{Z}\left[  \tfrac{1}{m_{N}}\right]  ^{N}\text{.}%
\]
But it follows from (\ref{eqSesquilabial.5}) that
\[
\frac{1}{m_{N}^{n}}\mathbb{Z}^{N}\subseteq J^{-pn}\mathbb{Z}^{N}\subseteq G
\]
for $n=1,2,\ldots,$ and hence the converse inclusion is valid.
\end{proof}

It follows from the inductive limit construction for $G$, i.e., $G=\bigcup
_{i\geq0}G_{i}$, $G_{i}\subset G_{i+1}$, that if the implication
(\ref{prop10.5(3)}) $\Rightarrow$ (\ref{prop10.5(4)}) holds for some
$p\in\mathbb{N}$, then it also holds for $p+1$, and so the scaling degree is
well defined.

While the two groups from (\ref{eqExtNew.pound}), $K_{0}\left(  \mathfrak
{A}_{J_{L}}\right)  $ and $K_{0}\left(  \mathfrak{A}_{L^{\prime}}\right)  $,
agree, the relationship between the corresponding scale of subgroups is more
subtle. Using (\ref{eqSesquilabial.5}) we can establish the following subgroup
inclusions:%
\begin{equation}
\text{(a) }G_{i}^{{}}\subset G_{6i}^{\prime}\text{\qquad and\qquad(b) }%
G_{i}^{\prime}\subset G_{3i}^{{}}, \label{eqSesquilabial.star}%
\end{equation}
where%
\[
G_{k}^{{}}:=J_{L}^{-k}\mathbb{Z}_{{}}^{2}\text{\qquad and\qquad}G_{k}^{\prime
}:=J_{L^{\prime}}^{-k}\mathbb{Z}_{{}}^{2}.
\]

To prove this, let $J=J_{L}$ and $K=J_{L^{\prime}}$. We proved after the
statement of (\ref{eqSesquilabial.5}) that there is an $E\in M_{2}\left(
\mathbb{Z}\right)  $ such that $J^{3}=16E$. Since $16K^{-1}\in M_{2}\left(
\mathbb{Z}\right)  $, we conclude that%
\[
J^{3}K^{-1}=E\cdot\left(  16K^{-1}\right)  \in M_{2}\left(  \mathbb{Z}%
\right)
\]
and so%
\[
K^{-1}\mathbb{Z}^{2}\subset J^{-3}\mathbb{Z}^{2}.
\]
Similarly%
\[
J^{3i}K^{-i}=E^{i}\left(  16^{i}K^{-1}\right)  \in M_{2}\left(  \mathbb{Z}%
\right)
\]
since each factor is in $M_{2}\left(  \mathbb{Z}\right)  $. This yields%
\[
K^{-i}\mathbb{Z}^{2}\subset J^{-3i}\mathbb{Z}^{2},
\]
which is the assertion (\ref{eqSesquilabial.star})(b). The claim
(\ref{eqSesquilabial.star})(a) follows by the same argument applied to the
factorization $K^{6}=32F$ for some $F\in M_{2}\left(  \mathbb{Z}\right)  $,
and $32J^{-1}\in M_{2}\left(  \mathbb{Z}\right)  $.

The discussion above leads to the notion of the \emph{degree} of an
isomorphism or homomorphism $\theta\colon G$ $\rightarrow G^{\prime}$ as
follows. When talking about isomorphisms and homomorphisms, we will henceforth
always assume that $\theta\left(  \lbrack\openone]\right)  =[\openone^{\prime
}]$, i.e., in the concrete representations,
\begin{equation}
G=\bigcup_{n}J^{-n}\mathbb{Z}^{N}=\bigcup_{n}G_{n}\text{,}\qquad G^{\prime
}=\bigcup_{n}J^{\prime\,-n}\mathbb{Z}^{N}=\bigcup_{n}G_{n}^{\prime}\text{,}
\label{eq10.8}%
\end{equation}
$\theta$ maps $\left(
\begin{smallmatrix}
1\\
0\\
\raisebox{0pt}[10pt]{$\vdots$}\\
0\\
0
\end{smallmatrix}
\right)  $ into $\left(
\begin{smallmatrix}
1\\
0\\
\raisebox{0pt}[10pt]{$\vdots$}\\
0\\
0
\end{smallmatrix}
\right)  $. We may assume $N=N^{\prime}$ in the discussion since $N$ is an
isomorphism invariant. But as $\theta$ is an order isomorphism, this means
that $\tau^{\prime}\circ\theta=\tau$ for the associated normalized traces.
Recall from (\ref{eqRemNewBis.20}) that
\begin{equation}
\tau^{\prime}=\left\langle 1,a^{\prime},a^{\prime\,2},\ldots,a^{\prime
\,N-1}\right|  ,\qquad\tau=\left\langle 1,a,\ldots,a^{N-1}\right|
\label{eq10.9}%
\end{equation}
where $a^{\prime}=1/\lambda^{\prime},a=1/\lambda$ where $\lambda$ is the
Perron--Frobenius eigenvalue of $J^{\prime},J$, respectively.

The following proposition is a globalization of Corollary \ref{CorRemNewBis.1}.

\begin{proposition}
\label{prop10.7}A map
\begin{equation}
\theta\colon G=\bigcup_{n}J^{-n}\mathbb{Z}^{N}\rightarrow G^{\prime}%
=\bigcup_{n}J^{\prime-n}\mathbb{Z}^{N} \label{eq10.10}%
\end{equation}
is an isomorphism between the ordered groups $(G,G_{+})$ and $(G^{\prime
},G_{+}^{\prime})$ $($mapping $[\openone]$ into $[\openone^{\prime}])$ if and
only if there exists a matrix $\Lambda\in\operatorname*{GL}(N,\mathbb{R}%
)$\label{LOSLambda_1} and a sequence $(n_{i})$ in $\mathbb{N}$ with the
following properties:
\begin{equation}%
\begin{array}
[c]{ll}%
1. & \theta(g)=\Lambda g\qquad\text{for }g\in G\subseteq\mathbb{R}^{N}%
\text{;}\\
2. & \alpha^{\prime}\Lambda=\alpha\text{\qquad where }\alpha=(1,a,\ldots
,a^{N-1})\text{, etc.;}\\
3. & J^{\prime\,n_{i}}\Lambda J^{-i}\in M_{N}(\mathbb{Z})\qquad\text{for
}i=1,2,\ldots\text{;}\\
4. & J^{n_{i}}\Lambda^{-1}J^{\prime\,-i}\in M_{N}(\mathbb{Z})\qquad\text{for
}i=1,2,\ldots\text{;}\\
5. & \Lambda\left(
\begin{smallmatrix}
1\\
0\\
\raisebox{0pt}[10pt]{$\vdots$}\\
0
\end{smallmatrix}
\right)  =\left(
\begin{smallmatrix}
1\\
0\\
\raisebox{0pt}[10pt]{$\vdots$}\\
0
\end{smallmatrix}
\right)  \text{.}%
\end{array}
\label{eq10.11}%
\end{equation}
Furthermore, $\theta$ is a homomorphism from $(G,G_{+})$ into $(G^{\prime
},G_{+}^{\prime})$ $($mapping $[\openone]$ into $[\openone^{\prime}])$ if and
only if there exists a matrix $\Lambda\in M_{N}(\mathbb{R})$ with the
properties \textup{(1), (2), (3)} and \textup{(5).} In both cases we actually
have
\[
\Lambda\in M_{N}\left(  \frac{1}{m_{N}^{\prime\,n_{1}}}\mathbb{Z}\right)
\]
\end{proposition}

\begin{proof}
Assume first that $\theta$ is given, and define the matrix $\Lambda$
by\label{LOSLambda_2}%
\begin{equation}
\Lambda=\left(  \theta\left(
\begin{array}
[c]{c}%
1\\
0\\
0\\
\vdots\\
0
\end{array}
\right)  ,\theta\left(
\begin{array}
[c]{c}%
0\\
1\\
0\\
\vdots\\
0
\end{array}
\right)  ,\dots,\theta\left(
\begin{array}
[c]{c}%
0\\
0\\
0\\
\vdots\\
1
\end{array}
\right)  \right)  \text{.} \label{eq10.12}%
\end{equation}
If $g\in G$, then $g\in G_{i}=J^{-i}\mathbb{Z}^{N}$ for some $i$, and hence
$g\in m_{N}^{-i}\mathbb{Z}^{N}$, i.e., $m_{N}^{i}g\in\mathbb{Z}^{N}$. But
then
\begin{equation}
m_{N}^{i}\theta(g)=\theta(m_{N}^{i}g)=\Lambda(m_{N}^{i}g)=m_{N}^{i}\Lambda g
\label{eq10.13}%
\end{equation}
and as $G$ has no torsion,
\begin{equation}
\theta(g)=\Lambda g \label{eq10.14}%
\end{equation}
which shows (1). But then (2) follows from $\tau^{\prime}\circ\theta=\tau$.
Furthermore, as $G_{i}=J^{-i}\mathbb{Z}^{N}$ is finitely generated, it follows
that $\theta(G_{i})\subseteq G_{n_{i}}^{\prime}$ for some $n_{i}$, i.e.,
\begin{equation}
\Lambda J^{-i}\mathbb{Z}_{N}\subseteq J^{-n_{i}}\mathbb{Z}_{N} \label{eq10.15}%
\end{equation}
which shows (3). The property (4) follows likewise from $\theta^{-1}%
(G_{i}^{\prime})\subseteq G_{m_{i}}$ for some $m_{i}$. Property (5) follows
from
\begin{equation}
\theta([\openone])=[\openone^{\prime}]\text{.} \label{eq10.16}%
\end{equation}
Conversely, if $\Lambda\in\operatorname*{GL}(N,\mathbb{R})$ is given with the
properties (2)--(4), one deduces that $\theta$ defined by (1) has the
properties:
\begin{equation}%
\begin{aligned}
\tau^{\prime}\circ\theta&=\tau, \\
\theta(G_{i})&\subseteq G_{n_{i}}^{\prime} , \\
\theta^{-1}(G_{i}^{^{\prime}})&\subseteq G_{m_{i}} , \\
\theta([\openone])&=[\openone^{\prime}] ,
\end{aligned}%
\label{eq10.17}%
\end{equation}
so $\theta$ is an order-automorphism from the property (\ref{eqRemNewBis.34}),
i.e.,%
\begin{equation}
G_{+}=\{g\in G\mid\tau(g)>0\}\cup\{0\}\text{.} \label{eq10.18}%
\end{equation}
The last statements in Proposition \ref{prop10.7} are straightforward from
$\left|  \det(J^{\prime\,n_{i}})\right|  =m_{N}^{\prime\,m_{i}}$.
\end{proof}

Our next aim is to show that the constants $n_{i}$, $m_{i}$ in Proposition
\ref{prop10.7} can be chosen such that they increase linearly with $i$. First,
a definition:

\begin{definition}
\label{definition10.8}Adopt the notation \textup{(\ref{eq10.8}),} and let
$\theta$ be a homomorphism from $G$ into $G^{\prime}$. We say that the degree
of $\theta\leq s,$ and write $\deg(\theta)\leq s$,\label{LOSdeg_3} if there
exists an $n_{0}\in\mathbb{N}$ such that
\begin{equation}
\theta(G_{i})\subseteq G_{n_{0}+\lfloor si\rfloor}^{\prime} \label{eq10.19}%
\end{equation}
for all $i\in\mathbb{N}$. \textup{(}Here $\lfloor si\rfloor$ is the largest
integer $\leq si$\textup{).}
\end{definition}

We next show that \textit{any} homomorphism $\theta\colon G\rightarrow
G^{\prime}$ (mapping $[\openone]$ into $[\openone^{\prime}]$) has a finite
degree which can be computed in concrete examples.

\begin{proposition}
\label{prop10.9}Let $G=G_{J}$, $G^{\prime}=G_{J^{\prime}}$ and let $\theta$ be
a morphism from $G$ into $G^{\prime}$ mapping $[\openone]$ into $[\openone
^{\prime}]$. Assume that $N=N^{\prime},\operatorname*{Prim}(m_{N}%
)=\operatorname*{Prim}(m_{N}^{\prime})$ and let $m=\operatorname{lcm}%
(m_{N},m_{N}^{\prime})$. Then\label{LOSdeg_4}\label{LOSmdeg_3}
\begin{equation}
\deg\theta\leq m\text{-}\deg(J^{\prime}). \label{eq10.20}%
\end{equation}
\textup{(}The last statement means that if $m$-$\deg(J^{\prime})\leq s$, then
$\deg\theta\leq s$.\textup{)}
\end{proposition}

\begin{proof}
If
\begin{equation}
\theta\colon\bigcup_{n}J^{-n}\mathbb{Z}^{N}\longrightarrow\bigcup_{n}%
J^{\prime\,-n}\mathbb{Z}^{N} \label{eq10.21}%
\end{equation}
is a morphism, there is a $k_{0}$ such that
\begin{equation}
\theta(\mathbb{Z}^{N})\subseteq J^{\prime\,-k_{0}}\mathbb{Z}^{N}=G_{k_{0}%
}^{\prime}\text{.} \label{eq10.22}%
\end{equation}
Then
\begin{align}
m^{k_{0}+i}\theta(G_{i})  &  =%
\genfrac{(}{)}{}{}{m}{m_{N}}%
^{i}m^{k_{0}}\theta(m_{N}^{i}G_{i})\label{eq10.23}\\
&  \subseteq%
\genfrac{(}{)}{}{}{m}{m_{N}}%
^{i}m^{k_{0}}\theta(\mathbb{Z}^{N})\nonumber\\
&  =%
\genfrac{(}{)}{}{}{m}{m_{N}}%
^{i}%
\genfrac{(}{)}{}{}{m}{m_{N}^{\prime}}%
^{k_{0}}m_{N}^{\prime\,k_{0}}\theta(\mathbb{Z}^{N})\nonumber\\
&  \subseteq%
\genfrac{(}{)}{}{}{m}{m_{N}}%
^{i}%
\genfrac{(}{)}{}{}{m}{m_{N}^{\prime}}%
^{k_{0}}m_{N}^{\prime\,k_{0}}G_{k_{0}}^{\prime}\nonumber\\
&  \subseteq%
\genfrac{(}{)}{}{}{m}{m_{N}}%
^{i}%
\genfrac{(}{)}{}{}{m}{m_{N}^{\prime}}%
^{k_{0}}\mathbb{Z}^{N}\nonumber\\
&  \subseteq\mathbb{Z}^{N}\text{.}\nonumber
\end{align}
Thus, if $m$-$\deg(J^{\prime})\leq p$, then
\begin{equation}
\theta(G_{i})\subseteq G_{\lfloor(i+k_{0})p\rfloor+n_{0}}^{\prime}
\label{eq10.24}%
\end{equation}
for some $n_{0}$ and all $i$, and hence $\deg\theta\leq p$.
\end{proof}

Note now that if
\begin{equation}
m\text{-}\deg(J^{\prime})\leq\frac{p}{q} \label{eq10.25}%
\end{equation}
where $p,q\in\mathbb{N}$, then the conclusion in Proposition \ref{prop10.9}
says that there is an $n_{0}$ such that
\begin{equation}
J^{\prime\,np+n_{0}}\Lambda J^{-nq}\in M_{N}(\mathbb{Z}) \label{eq10.26}%
\end{equation}
for $n=1,2,\ldots,$ where $\Lambda$ is the matrix in Proposition
\ref{prop10.7}. This implies the main result on the $\operatorname*{Ext}%
$-invariant in this chapter, which, together with Remark
\ref{RemSesquilabialNew.11}, is surprisingly effective in establishing
non-isomorphism when all the elementary invariants in Chapter
\ref{Pediferient} are the same.

\begin{theorem}
\label{theorem10.10}Let $J$, $J^{\prime}$ be matrices of the form
\textup{(\ref{eqCycNew.1})} with $N=N^{\prime}$ and $\operatorname*{Prim}%
m_{N}=\operatorname*{Prim}m_{N}^{\prime}$. Assume that the associated
Perron--Frobenius eigenvalues $\lambda$, $\lambda^{\prime}$ are rational
\textup{(}and thus integers\textup{),} and let $m=\operatorname{lcm}%
(m_{N},m_{N}^{\prime})$. Assume there exist rational numbers $p$, $q$ with
\begin{equation}
m\text{-}\deg(J^{\prime})\leq\frac{p}{q} \label{eq10.27}%
\end{equation}
and such that
\begin{equation}
\lambda^{\prime\,p}=\lambda^{q}\text{.} \label{eq10.28}%
\end{equation}
Let $\alpha^{\prime}$, $\alpha$, resp.$\ v^{\prime}$, $v$, be the left,
respectively right, Perron--Frobenius eigenvectors of $J$, $J^{\prime}$ given
by \textup{(\ref{eqRemNewBis.17}), (\ref{eqSubNew.pound}),} respectively.

If there exists a unital isomorphism $\mathfrak{A}_{J}\rightarrow\mathfrak
{A}_{J^{\prime}}$, then
\begin{equation}
\frac{\left\langle \alpha\mid v\right\rangle }{\left\langle \alpha^{\prime
}\mid v^{\prime}\right\rangle }\in\mathbb{Z}\left[  \tfrac{1}{\lambda^{\prime
}}\right]  \text{.} \label{eq10.29}%
\end{equation}
\end{theorem}

\begin{remark}
\label{RemSesquilabialNew.11}A rather effective ``workhorse'' to show
non-isomorphism in cases where all the basic invariants in Theorem
\textup{\ref{CorCyc.9}} are the same and the Perron--Frobenius eigenvalues are
integers, is to use Theorem \textup{\ref{theorem10.10}} together with the fact
that $m_{N}^{\prime}$-$\deg(J^{\prime})=1$ if $m_{N-1}^{\prime}$ is nonzero
and $m_{N}^{\prime},m_{N-1}^{\prime}$ are mutually prime. See Proposition
\textup{\ref{prop10.22}} and Remark \textup{\ref{remark10.20}} below. More
generally, $m_{N}^{\prime}$-$\deg(J^{\prime})=n$ if $m_{N-k}^{\prime}%
=0\mod{m_{N}}$ for $k=1,\ldots,n-1$ and $m_{N-n}^{\prime}$ and $m_{N}^{\prime
}$ are mutually prime.

Note for example that this theorem covers the two matrices:
\[
J_{A}=\left(
\begin{array}
[c]{ccc}%
0 & 1 & 0\\
3 & 0 & 1\\
2 & 0 & 0
\end{array}
\right)  \qquad J_{B}=\left(
\begin{array}
[c]{ccc}%
0 & 1 & 0\\
15 & 0 & 1\\
2 & 0 & 0
\end{array}
\right)
\]
considered in Observation \textup{\ref{ObsCyc.3}}. We have $\lambda
_{A}=2,\lambda_{B}=4$, and since $3$ and $15$ are mutually prime with $2$, we
have $2\deg J_{A}=1=4$-$\deg J_{B}$. \textup{(}This also follows from Lemma
\textup{\ref{lemma10.13}).} Using the computation in Observation
\textup{\ref{ObsCyc.3},}
\[
\left\langle \alpha_{A}\mid v_{A}\right\rangle =\frac{9}{4},\qquad\left\langle
a_{B}\mid v_{B}\right\rangle =\frac{13}{4}\text{.}%
\]
Since $1$-$\deg J_{A}=1$, we have $2$-$\deg J_{A}=2$ by Proposition
\textup{\ref{prop10.5}}, and applying Theorem \textup{\ref{theorem10.10}} with
$m=\operatorname{lcm}(2,4)=4$, $J^{\prime}=J_{A}$, $J=J_{B}$, $p=2$, $q=1$, we
see that \textup{(\ref{eq10.29})} takes the form
\[
\frac{13}{9}\in\mathbb{Z}\left[  \tfrac{1}{2}\right]
\]
which is clearly false. Hence, there is no unital morphism $\mathfrak{A}%
_{B}\rightarrow\mathfrak{A}_{A}$. Of course, the proof of Theorem
\textup{\ref{theorem10.10}} is just a more general version of the proof of
Observation \textup{\ref{ObsCyc.3}.}
\end{remark}

\begin{proof}
By (\ref{eq10.27}) in the form (\ref{eq10.26}), it follows that there exists
an $n_{0}$ such that
\begin{equation}
J^{\prime\,np+n_{0}}\Lambda J^{-nq}\in M_{N}(\mathbb{Z}) \label{eq10.30}%
\end{equation}
for $n=1,2,\ldots$, when $\Lambda$ is the matrix associated to the
homomorphism
\[
\theta\colon G_{\Gamma}\longrightarrow G_{\Gamma^{\prime}}%
\]
by Proposition \ref{prop10.7}. We have
\begin{equation}
J^{\prime\,np+n_{0}}\Lambda J^{-nq}v=\lambda^{-nq}J^{\prime\,np+n_{0}}\Lambda
v \label{eq10.31}%
\end{equation}
Since $\lambda^{-nq}\lambda^{\prime\,np}=1$ for all $n$, it follows from
Perron--Frobenius theory that there exists a constant $c\in\mathbb{R}$ such
that
\begin{equation}
\lim_{n\rightarrow\infty}\lambda^{-nq}J^{\prime\,np+n_{0}}\Lambda
v=cv^{\prime}\text{.} \label{eq10.32}%
\end{equation}
But the first component of $v^{\prime}$ is $1$ by (\ref{eqSubNew.pound}) and
\[
J^{\prime\,np+n_{0}}\Lambda J^{-nq}\in M_{N}(\mathbb{Z})
\]
for all $n$, and since \emph{all} components of $v$ are integer by
(\ref{eqSubNew.pound}) and $\lambda\in\mathbb{N}$, it follows that $c$ in an
integer. But it follows from (\ref{eq10.32}) that
\begin{align*}
c\left\langle \alpha^{\prime}\mid v^{\prime}\right\rangle  &  =\lim
_{n\rightarrow\infty}\lambda^{-nq}\left\langle \alpha^{\prime}\mid
J^{\prime\,np+n_{0}}\Lambda v\right\rangle \\
&  =\lim_{n\rightarrow\infty}\lambda^{-nq}\lambda^{\prime\,np+n_{0}%
}\left\langle \alpha^{\prime}\mid\Lambda v\right\rangle \\
&  =\lambda^{\prime\,n_{0}}\left\langle \alpha\mid v\right\rangle
\end{align*}
where we used (\ref{eq10.11})(2). It follows that
\[
\frac{\left\langle \alpha\mid v\right\rangle }{\left\langle \alpha^{\prime
}\mid v^{\prime}\right\rangle }=\frac{c}{\lambda^{\prime\,n_{0}}}%
\]
so
\[
\frac{\left\langle \alpha\mid v\right\rangle }{\left\langle \alpha^{\prime
}\mid v^{\prime}\right\rangle }\in\mathbb{Z}\left[  \tfrac{1}{\lambda^{\prime
}}\right]  =\mathbb{Z}\left[  \tfrac{1}{\lambda}\right]  .\settowidth
{\qedskip}{$\displaystyle
\frac{\left\langle\alpha\mid v\right\rangle}{\left\langle\alpha^{\prime}\mid
v^{\prime}\right\rangle}\in\mathbb{Z}\left[ \textstyle\frac{1}{\lambda
^{\prime}}\right] =\mathbb{Z}\left[ \textstyle\frac{1}{\lambda}\right
] .$}\addtolength{\qedskip}{-\textwidth}\rlap{\hbox to-0.5\qedskip{\hfill\qed
}}%
\]
\renewcommand{\qed}{}
\end{proof}

We will often apply this theorem in the following special form:

\begin{corollary}
\label{corollary10.11}Let $J$, $J^{\prime}$ be matrices of the form
\textup{(\ref{eqCycNew.1})} with $N=N^{\prime}$ and $m_{N}\mid m_{N}^{\prime}%
$, and same Perron--Frobenius eigenvalue $\lambda=\lambda^{\prime}$. If
$\deg(J^{\prime})\leq1$ and there exists a unital morphism $\mathfrak{A}%
_{J}\rightarrow\mathfrak{A}_{J^{\prime}}$, then
\begin{equation}
\frac{\left\langle \alpha\mid v\right\rangle }{\left\langle \alpha^{\prime
}\mid v^{\prime}\right\rangle }\in\mathbb{Z}\left[  \tfrac{1}{\lambda}\right]
=\mathbb{Z}\left[  \tfrac{1}{\lambda^{\prime}}\right]  \text{.}
\label{eq10.33}%
\end{equation}
\end{corollary}

\begin{proof}
In this case $m=m_{N}^{\prime}$.
\end{proof}

The following even more special corollary will be useful in Chapter
\ref{APP.EXA}.

\begin{corollary}
\label{corollary10.12}Let $J$, $J^{\prime}$ be matrices of the form
\textup{(\ref{eqCycNew.1})} with $N=N^{\prime}$ and $m_{N}=m_{N^{\prime}%
}=\lambda=\lambda^{\prime}$. If there exists a unital morphism $\mathfrak
{A}_{J}\rightarrow\mathfrak{A}_{J^{\prime}}$, then
\begin{equation}
\frac{\left\langle \alpha\mid v\right\rangle }{\left\langle \alpha^{\prime
}\mid v^{\prime}\right\rangle }\in\mathbb{Z}\left[  \tfrac{1}{\lambda}\right]
\text{.} \label{eq10.34}%
\end{equation}
\end{corollary}

\begin{proof}
First apply Lemma \ref{lemma10.13} below to deduce $m_{N}$-$\deg J=1$ and
$m_{k}^{\prime}$-$\deg J^{\prime}=1$. But then we may apply Theorem
\ref{theorem10.10} with $p=q=1$ and $m=m_{N}=m_{N}^{\prime}$.
\end{proof}

\begin{lemma}
\label{lemma10.13}Let $J$ be a matrix of the form \textup{(\ref{eqCycNew.1})}
with $m_{N}=\lambda$. Then
\begin{equation}
m_{N}\text{-}\deg(J)=\deg(J)=1\text{.} \label{eq10.35}%
\end{equation}
\end{lemma}

\begin{remark}
\label{RemSesquilabialNew.15}If $\lambda=m_{N}$, it follows that
$R_{D}=\lambda=\pm m_{N}$ and $Q_{N-D}=\pm1$ in \textup{(\ref{eqCyc.22}),}
hence $\left|  \det J_{0}\right|  =1$ where $J_{0}$ is given by
\textup{(\ref{eqj0qqqq})} and $\ker\tau\cong\mathbb{Z}^{N-1}$ by Theorem
\textup{\ref{CorCyc.9}.} Conversely, if $\ker\tau\cong\mathbb{Z}^{N-1}$ we
must have $\left|  \det J_{0}\right|  =1$, hence $R_{D}=\pm\lambda=\pm m_{N}$,
i.e., $\lambda=m_{N}$.

Conclusion:%
\[
\lambda=m_{N}\Longleftrightarrow\ker\tau\cong\mathbb{Z}^{N-1}\text{.}%
\]
The theory for the $\lambda=m_{N}$ case will be developed in much more detail
in Chapter \textup{\ref{ClmN}.}
\end{remark}

\begin{proof}
If $\lambda=m_{N}$, the matrix $J_{\mathcal{B}}^{\mathcal{B}}$ in
(\ref{eqCyc.22}) takes the form (now $D=1$):
\begin{equation}
J_{\mathcal{B}}^{\mathcal{B}}=\left(
\begin{array}
[c]{ccccccc}%
Q_{1_{\mathstrut}} & 1 & 0 & \cdots & 0 & 0 & Q_{1}\vphantom{\vdots}\\
Q_{2_{\mathstrut}} & 0 & 1 & \cdots & 0 & 0 & Q_{2}\vphantom{\vdots}\\
Q_{3} & 0 & 0 & \ddots & 0 & 0 & Q_{3}\vphantom{\vdots}\\
\vdots & \vdots &  & \ddots & \ddots & \vdots & \vdots\vphantom{\vdots}\\
Q_{N-2_{\mathstrut}} & 0 & 0 &  & 0 & 1 & Q_{N-2}\vphantom{\vdots}\\
-1 & 0 & 0 & \cdots & 0 & 0 & -1\vphantom{\vdots}\\
0 & 0 & 0 & \cdots & 0 & 0 & m_{N}\vphantom{\vdots}%
\end{array}
\right)  \label{eq10.36}%
\end{equation}
where
\begin{equation}%
\begin{aligned}
Q_{1} & =  m_{1_{\mathstrut}}-\lambda\\
Q_{2} & =  m_{2}+\lambda m_{1_{\mathstrut}}-\lambda^{2^{\mathstrut}}\\
Q_{3} & =  m_{3}+\lambda m_{2}+\lambda^{2^{\mathstrut}}m_{1_{\mathstrut
}}-\lambda^{3}\\
\vdots&  \qquad\vdots\\
Q_{N-1} & =  m_{N-1}+\lambda m_{N-2}+\cdots+\lambda^{N-2^{\mathstrut}%
}m_{1_{\mathstrut}}-\lambda^{N-1}
=  -1^{\mathstrut}\end{aligned}%
\label{eq10.37}%
\end{equation}
Note in passing that
\begin{equation}
v=\left(
\begin{array}
[c]{c}%
1\\
-Q_{1}\\
-Q_{2}\\
\vdots\\
-Q_{N-1}%
\end{array}
\right)  \label{eq10.38}%
\end{equation}
is the right Perron--Frobenius eigenvector of $J$ by (\ref{eqSubNew.pound}).

Using (\ref{eqCyc.plus6}) one computes that
\begin{equation}
J_{\mathcal{B}}^{\mathcal{B\,}-1}=\left(
\begin{tabular}
[c]{cccccc|c}%
$0$ & $0$ & $\cdots$ & $0$ & $0$ & $-1$ & $-\frac{1}{m_{N}}\vphantom{\vdots}%
$\\
$1$ & $0$ &  & $0$ & $0$ & $Q_{1}$ & $0\vphantom{\vdots}$\\
$0$ & $1$ & $\ddots$ &  & $0$ & $Q_{2}$ & $0\vphantom{\vdots}$\\
$\vdots$ & $\vdots$ & $\ddots$ & $\ddots$ &  & $\vdots$ & $\vdots
\vphantom{\vdots}$\\
$0$ & $0$ & $\cdots$ & $1$ & $0$ & $Q_{N-3}$ & $0\vphantom{\vdots}$\\
$0$ & $0$ & $\cdots$ & $0$ & $1$ & $Q_{N-2}$ & $0\vphantom{\vdots}$\\\hline
$0$ & $0$ & $\cdots$ & $0$ & $0$ & $0$ & $\frac{1^{\mathstrut}}{m_{N}%
}\vphantom{\vdots}$%
\end{tabular}
\right)  \label{eq10.39}%
\end{equation}

Iterating, one computes that
\begin{equation}
J_{\mathcal{B}}^{\mathcal{B\,}-k}=\left(
\begin{tabular}
[c]{llllll|c}%
&  &  &  &  &  & $p_{1}^{(k)}%
\genfrac{(}{)}{}{}{1}{m_{N}}%
$\\
&  & \rlap{\raisebox{-6pt}[0pt][0pt]{$\phantom{\cdots}J_{0}^{-k}$}} &  &  &  &
$\vdots$\\
&  &  &  &  &  & $\vdots$\\
&  &  &  &  &  & $p_{N-1}^{(k)}%
\genfrac{(}{)}{}{}{1}{m_{N}}%
_{\mathstrut}$\\\hline
$0$ & $0$ & $\cdots$ & $\cdots$ & $0$ & $0$ & $\frac{1^{\mathstrut}}%
{m_{N}^{k^{\mathstrut}}}$%
\end{tabular}
\right)  \label{eq10.40}%
\end{equation}
where $p_{l}^{(k)}(x)\in\mathbb{Z}[x]$ and the degree of the polynomial
$p_{l}^{(k)}$ is:
\[
\deg p_{l}^{(k)}=(k-l+1)\vee0
\]
for $l=1,\ldots,N-1$. Using this and the transformation matrix
\begin{multline}
I_{\mathcal{A}}^{\mathcal{B}}=(I_{\mathcal{B}}^{\mathcal{A}})^{-1}%
\label{eq10.41}\\
=\left(
\begin{array}
[c]{ccccccc}%
0 & -\frac{1}{m_{N}} & -\frac{1}{m_{N}^{2^{\mathstrut}}} & \cdots & -\frac
{1}{m_{N}^{N-3^{\mathstrut}}} & -\frac{1}{m_{N}^{N-2^{\mathstrut}}} &
-\frac{1}{m_{N}^{N-1^{\mathstrut}}}\vphantom{\vdots}\\
0 & 0 & -\frac{1}{m_{N}} & \phantom{-\frac{1}{m_{N}^{N-4^{\mathstrut}}}} &
-\frac{1}{m_{N}^{N-4^{\mathstrut}}} & -\frac{1}{m_{N}^{N-3^{\mathstrut}}} &
-\frac{1}{m_{N}^{N-2^{\mathstrut}}}\vphantom{\vdots}\\
\vdots &  & \ddots & \ddots &  & \vdots & \vdots\\
\vdots &  &  & \ddots & \ddots & \vdots & \vdots\\
0 & 0 & 0 &  & 0 & -\frac{1}{m_{N}} & -\frac{1}{m_{N}^{2^{\mathstrut}}%
}\vphantom{\vdots}\\
0 & 0 & 0 & \cdots & 0 & 0 & -\frac{1}{m_{N}}\vphantom{\vdots}\\
1 & \frac{1}{m_{N}} & \frac{1}{m_{N}^{2^{\mathstrut}}} & \cdots & \frac
{1}{m_{N}^{N-3^{\mathstrut}}} & \frac{1}{m_{N}^{N-2^{\mathstrut}}} & \frac
{1}{m_{N}^{N-1^{\mathstrut}}}\vphantom{\vdots}%
\end{array}
\right)
\end{multline}
and the definition of $m_{N}$-$\deg$ together with (\ref{eqPediferientOla.59})
one has to show that for any $l,k\in\mathbb{N}$ with $l\leq k$ and any
$n\in\mathbb{Z}^{N}$ that
\begin{equation}
\left\{
\begin{array}
[c]{c}%
m_{N_{\mathstrut}}^{l}J_{\mathcal{B}}^{\mathcal{B\,}-k}I_{\mathcal{A}%
}^{\mathcal{B}}n\in I_{\mathcal{A}}^{\mathcal{B}}\mathbb{Z}^{N}\\
\Downarrow_{\mathstrut}^{\mathstrut}\\
J_{\mathcal{B}}^{\mathcal{B\,}-k}I_{\mathcal{A}}^{\mathcal{B}}n\in
J_{\mathcal{B}_{\mathstrut}}^{\mathcal{B\,}-l}I_{\mathcal{A}}^{\mathcal{B}%
}\mathbb{Z}^{N}\\
\Updownarrow_{\mathstrut}^{\mathstrut}\\
J_{\mathcal{B}}^{\mathcal{B\,}l-k}I_{\mathcal{A}}^{\mathcal{B}}n\in
I_{\mathcal{A}}^{\mathcal{B}}\mathbb{Z}^{N}\text{.}%
\end{array}
\right.  \label{eq10.42}%
\end{equation}
(The last equivalence is trivial.) This can be done by brute force, looking at
highest order terms in $\frac{1}{m_{N}}$. We do omit the painful details,
however, since the result can also be proved by another method described below.
\end{proof}

The alternative way of proving Lemma \ref{lemma10.13} is based on:

\begin{proposition}
\label{prop10.14}Let $J$ be a matrix of the form \textup{(\ref{eqCycNew.1})}
with $m_{N}=\lambda=m$. Use the concrete realization
\textup{(\ref{eqRemNewBis.6})} or \textup{(\ref{eqzjzjz})} of $G=K_{0}%
(\mathfrak{A}_{J})$, and define $F=G/\mathbb{Z}^{N}$. Then the generators
$x_{1}$, $x_{2}$, $x_{3}$ of $F$ defined in Lemma \textup{\ref{LemSubNew.14}}
satisfy
\begin{equation}
mx_{1}=0 \label{eq10.43}%
\end{equation}
and
\begin{equation}
mx_{i}=x_{i-1},\qquad i=2,3,\ldots\text{.} \label{eq10.44}%
\end{equation}
Thus,\label{LOSFL_2}
\begin{equation}
F\cong\mathbb{Z}\left[  \tfrac{1}{m}\right]  \diagup\mathbb{Z} \label{eq10.45}%
\end{equation}
the isomorphism being given by\label{LOSxi_2}
\begin{equation}
x_{i}\longrightarrow\frac{1}{m^{i}}\text{.} \label{eq10.46}%
\end{equation}
\end{proposition}

\begin{proof}
{}From the relations (\ref{eq10.37}), it follows that
\begin{align}
m_{1}  &  =Q_{1}+m\label{eq.10.47}\\
m_{2}  &  =-mQ_{1}+Q_{2}\nonumber\\
m_{3}  &  =-mQ_{2}+Q_{3}\nonumber\\
\vdots &  \qquad\vdots\nonumber\\
m_{N-1}  &  =-mQ_{N-2}+Q_{N-1}\nonumber\\
m_{N}=m  &  =-mQ_{N-1}\Rightarrow Q_{N-1}=-1\text{.}\nonumber
\end{align}
Inserting these relations in the relations \textup{(\ref{eqStr.G})} for
$x_{i}$ in Lemma \ref{LemSubNew.14} gives
\begin{align}
mx_{i}=x_{i-N}  &  -(Q_{1}+m)x_{i-N+1}\label{eq10.48}\\
&  -(-mQ_{1}+Q_{2})x_{i-N+2}\nonumber\\
&  -(-mQ_{2}+Q_{3})x_{i-N+3}\nonumber\\
&  -\cdots\nonumber\\
&  -(-mQ_{N-2}+Q_{N-1})x_{i-1}\text{.}\nonumber
\end{align}
We know already from \textup{(\ref{eqStr.F})} that $x_{i}=0$ for
$i=1-N,2-N,\ldots,0$. Assume by induction that $x_{j-1}=mx_{j}$ for all $j<i$.
It follows from (\ref{eq10.48}) that
\begin{align*}
mx_{i}=x_{i-N}  &  -Q_{1}x_{i-N+1}-x_{i-N}\\
&  -Q_{2}x_{i-N-2}+Q_{1}x_{i-N+1}\\
&  -Q_{3}x_{i-N+3}+Q_{2}x_{i-N+2}\\
&  -\cdots\\
&  -Q_{N-1}x_{i-1}+Q_{N-2}x_{i}\\
=  &  -(-1)x_{i-1}=x_{i-1}\text{.}%
\end{align*}
This shows (\ref{eq10.44}), and the remaining statements in Proposition
\ref{prop10.14} are obvious.
\end{proof}

\begin{proof}
[Alternative Proof of Lemma \textup{\ref{lemma10.13}}]Use the notation of
Proposition \ref{prop10.14} and define
\begin{equation}
G_{i}=J^{-i}\mathbb{Z}^{N}\text{,\qquad}F_{i}=G_{i}\diagup\mathbb{Z}%
^{N}\text{.} \label{eq10.50}%
\end{equation}
It follows from Proposition \ref{prop10.14} that
\begin{equation}
F_{i}\cong\mathbb{Z}_{N^{i}}=\mathbb{Z}\diagup N^{i}\mathbb{Z}\text{.}
\label{eq10.51}%
\end{equation}
The conclusion in Lemma \ref{lemma10.13} is
\begin{equation}
G_{i}=\{g\in G\mid m^{i}g\in\mathbb{Z}^{n}\}\text{.} \label{eq10.52}%
\end{equation}
But
\begin{equation}
m^{i}g\in\mathbb{Z}^{N}\Longleftrightarrow m^{i}(g+\mathbb{Z}^{n})=0\text{,}
\label{eq10.53}%
\end{equation}
so this conclusion is equivalent to
\begin{equation}
F_{i}=\{h\in F\mid m^{i}F=0\}\text{.} \label{eq10.54}%
\end{equation}
But the last statement is obvious from (\ref{eq10.46}) and (\ref{eq10.51}).
This proves Lemma \ref{lemma10.13}.
\end{proof}

We will now prove a theorem somewhat close in spirit to Theorem
\ref{theorem10.10}. If $G$ is a torsion free abelian group, and
$n=2,3,4,\ldots$, we define\label{LOSDlambdaG_1}
\begin{align}
D_{n}(G)  &  =\bigcap_{k=1}^{\infty}n^{k}G\label{eq10.55}\\
&  =%
\begin{array}
[t]{l}%
\text{the set of elements}\\
\text{of }G\text{ which are }\\
\text{infinitely divisible by }n\text{.}%
\end{array}
\nonumber
\end{align}

$D_{n}(G)$, as well as its rank, is clearly an isomorphism invariant of $G$,
and any homomorphism $G\rightarrow G^{\prime}$ will map $D_{n}(G)$ into
$D_{n}(G^{\prime})$. $D_{n}(G)$ only depends on $G$ and the prime factors of
$n$, and $D_{n}(G)$ is in a natural way a $\mathbb{Z}\left[  \frac{1}%
{n}\right]  $-module. Note that even if the rank of $D_{n}(G)$ is $1$,
$D_{n}(G)$ is not necessarily isomorphic to $\mathbb{Z}\left[  \frac{1}%
{n}\right]  $, as seen from the example $G=\mathbb{Z}\left[  \frac{1}%
{6}\right]  $ and $n=2$. But in the rank $1$ case, $D_{n}(G)$ is isomorphic
with a subgroup of the additive group $\mathbb{R}$ containing $\mathbb{Z}%
[\frac{1}{n}]$.

In the special case that $G=K_{0}(\mathfrak{A}_{J})$, where $J$ is a matrix of
the form (\ref{eqCycNew.1}), and the Perron--Frobenius eigenvalue $\lambda$ of
$J$ is rational, and thus an integer, we note that the right Perron--Frobenius
eigenvector $v$, normalized as in (\ref{eqSubNew.pound}), is contained in
$D_{\lambda}(G)$. In fact, since $J^{-1}v=\lambda^{-1}v$ we have
\begin{equation}
D_{\lambda}(G)\supseteq\mathbb{Z}\left[  \tfrac{1}{\lambda}\right]  v\text{.}
\label{eq10.56}%
\end{equation}

If, furthermore, $\operatorname*{rank}$ $(D_{\lambda}(G))=1$, there exists a
subgroup $D_{\lambda}^{\mathbb{Q}}(G)\subseteq\mathbb{Q}$ such that
\begin{equation}
D_{\lambda}(G)=D_{\lambda}^{\mathbb{Q}}(G)v \label{eq10.57}%
\end{equation}
and this identity defines an isomorphism between $D_{\lambda}(G)$ and
$D_{\lambda}^{\mathbb{Q}}(G)$.

Let $D_{\lambda}^{\mathbb{Q}U}(G)$ be the set of multiplicatively invertible
elements of $D_{\lambda}^{\mathbb{Q}}(G)$, so if for example $D_{\lambda
}^{\mathbb{Q}U}(G)=\mathbb{Z}\left[  \frac{1}{n}\right]  $, this is the set of
numbers of the form $p_{1}^{n_{1}}p_{2}^{n_{2}}\ldots p_{m}^{n_{m}}$ where
$p_{1}\cdots p_{m}$ are the prime factors in $n$ and $n_{1},n_{2},\ldots
,n_{m}\in\mathbb{Z}$.

\begin{theorem}
\label{theorem10.15}Let $J$, $J^{\prime}$ be matrices of the form
\textup{(\ref{eqCycNew.1})} with rational Perron--Frobenius eigenvalues
$\lambda$, $\lambda^{\prime}$ and $\operatorname*{Prim}(\lambda
)=\operatorname*{Prim}(\lambda^{\prime})=\{p_{1},\ldots,p_{m}\}$%
.\label{LOSPrim_3} Assume\label{LOSrank_1}
\[
\operatorname*{rank}(D_{\lambda^{\prime}}(K_{0}(\mathfrak{A}_{J^{\prime}%
})))=1
\]
and let $\alpha$, $\alpha^{\prime}$, resp.$\ v$, $v^{\prime}$, be the left,
respectively right, Perron--Frobenius eigenvectors of $J$, $J^{\prime}$ given
by \textup{(\ref{eqRemNewBis.17}), (\ref{eqSubNew.pound}),} respectively. If
there exists a unital morphism $\mathfrak{A}_{J}\rightarrow\mathfrak
{A}_{J^{\prime}}$, then
\begin{equation}
\frac{\left\langle \alpha\mid v\right\rangle }{\left\langle \alpha^{\prime
}\mid v^{\prime}\right\rangle }\in D_{\lambda^{\prime}}^{\mathbb{Q}}%
(G^{\prime}) \label{eq10.59}%
\end{equation}
In particular, if
\begin{equation}
D_{\lambda^{\prime}}(K_{0}(\mathfrak{A}_{J^{\prime}}))=\mathbb{Z}\left[
\tfrac{1}{\lambda^{\prime}}\right]  v^{\prime} \label{eq10.60}%
\end{equation}
we conclude
\begin{equation}
\frac{\left\langle \alpha\mid v\right\rangle }{\left\langle \alpha^{\prime
}\mid v^{\prime}\right\rangle }\in\mathbb{Z}\left[  \tfrac{1}{\lambda}\right]
=\mathbb{Z}\left[  \tfrac{1}{\lambda^{\prime}}\right]  \text{.}
\label{eq10.61}%
\end{equation}
If in addition the unital homomorphism $\mathfrak{A}_{J}\rightarrow
\mathfrak{A}_{J^{\prime}}$ is an isomorphism, and
\begin{equation}
D_{\lambda}(K_{0}(\mathfrak{A}_{J}))=\mathbb{Z}\left[  \tfrac{1}{\lambda
}\right]  v\text{,} \label{eq10.62}%
\end{equation}
then there exist integers $n_{1},\ldots,n_{m}$ such that
\begin{equation}
\frac{\left\langle \alpha\mid v\right\rangle }{\left\langle \alpha^{\prime
}\mid v^{\prime}\right\rangle }=p_{1}^{n_{1}}p_{2}^{n_{2}}\cdots p_{m}^{n_{m}%
}\text{.} \label{eq10.63}%
\end{equation}
\end{theorem}

\begin{remark}
It does not suffice instead of \textup{(\ref{eq10.60})} to assume that
$D_{\lambda^{\prime}}(K_{0}(\mathfrak{A}_{J^{\prime}}))$ is isomorphic to
$\mathbb{Z}\left[  \frac{1}{\lambda^{\prime}}\right]  $, because, say,
$\mathbb{Z}\left[  \frac{1}{\lambda^{\prime}}\right]  \diagup p$ is isomorphic
to $\mathbb{Z}\left[  \frac{1}{\lambda^{\prime}}\right]  $ if $p$ is a prime
not in the set $\{p_{1},\ldots,p_{m}\}$. For the same reason, even though the
existence of an isomorphism $\mathfrak{A}_{J}\rightarrow\mathfrak
{A}_{J^{\prime}}$ implies that $D_{\lambda}(K_{0}(\mathfrak{A}_{J}%
))\approx\mathbb{Z}\left[  \frac{1}{\lambda^{\prime}}\right]  $ one cannot
omit condition \textup{(\ref{eq10.62})} to obtain \textup{(\ref{eq10.63}).}
See the proof for explanation of this.
\end{remark}

\begin{proof}
Let $\psi\colon K_{0}(\mathfrak{A}_{J})\rightarrow K_{0}(\mathfrak
{A}_{J^{\prime}})$ be the $K$-theory morphism defined by the morphism
$\mathfrak{A}_{J}\rightarrow\mathfrak{A}_{J^{\prime}}$. Since
$\operatorname*{Prim}(\lambda)=\operatorname*{Prim}(\lambda^{\prime})$ we
have
\[
D_{\lambda}(K_{0}(\mathfrak{A}_{J}))=D_{\lambda^{\prime}}(K_{0}(\mathfrak
{A}_{J}))
\]
and since
\[
D_{\lambda^{\prime}}(K_{0}(\mathfrak{A}_{J^{\prime}}))=D_{\lambda^{\prime}%
}^{\mathbb{Q}}(G^{\prime})v^{\prime}%
\]
and thus
\[
\psi(v)\in\psi(D_{\lambda}(K_{0}(\mathfrak{A}_{J}))\subseteq D_{\lambda
^{\prime}}^{\mathbb{Q}}(G^{\prime})v^{\prime}%
\]
there is a $\xi\in D_{\lambda^{\prime}}^{\mathbb{Q}}(G^{\prime})$ such that
\[
\psi(v)=\xi v^{\prime}\text{.}%
\]
Apply $\left\langle \alpha^{\prime}\mid\right.  $ to both sides
\[
\left\langle \alpha^{\prime}\mid\psi(v)\right\rangle =\xi\left\langle
\alpha^{\prime}\mid v^{\prime}\right\rangle \text{.}%
\]
But since the morphism is assumed to be unital, we have $\left\langle
\alpha^{\prime}\mid\psi\right.  =\left\langle \alpha\mid\right.  $ by
uniqueness of the trace, and hence
\[
\left\langle \alpha\mid v\right\rangle =\xi\left\langle \alpha^{\prime}\mid
v^{\prime}\right\rangle \text{.}%
\]
This proves (\ref{eq10.59}). Since (\ref{eq10.60}) means $D_{\lambda^{\prime}%
}^{\mathbb{Q}}(K_{0}(\mathfrak{A}_{J^{\prime}}))=\mathbb{Z}\left[  \frac
{1}{\lambda^{\prime}}\right]  $, (\ref{eq10.61}) follows. Finally, if the
homomorphism $\mathfrak{A}_{J}\rightarrow\mathfrak{A}_{J^{\prime}}$ is an
isomorphism and (\ref{eq10.62}) holds, it follows by reverting the proof that
$\frac{\left\langle \alpha^{\prime}\mid v^{\prime}\right\rangle }{\left\langle
\alpha\mid v\right\rangle }\in\mathbb{Z}\left[  \frac{1}{\lambda}\right]  $,
and thus $\frac{\left\langle \alpha\mid v\right\rangle }{\left\langle
\alpha^{\prime}\mid v^{\prime}\right\rangle }$ has a multiplicative inverse in
$\mathbb{Z}\left[  \frac{1}{\lambda}\right]  $. But multiplicative invertible
elements of $\mathbb{Z}\left[  \frac{1}{\lambda}\right]  $ have the form on
the right-hand side of (\ref{eq10.63}).
\end{proof}

There is one interesting circumstance where (\ref{eq10.61}) or (\ref{eq10.63})
is automatically satisfied.

\begin{lemma}
\label{lemma10.16}Let $J$ be a matrix of the form \textup{(\ref{eqCycNew.1})}
with rational \textup{(}and thus integer\/\textup{)} Perron--Frobenius
eigenvalue $\lambda$. Assume that
\begin{equation}
\operatorname*{Prim}(m_{N})=\operatorname*{Prim}(\lambda) \label{eq10.64}%
\end{equation}
and that
\begin{equation}
\operatorname*{rank}(D_{\lambda}(G))=1\text{.} \label{eq10.65}%
\end{equation}
It follows that
\begin{equation}
D_{\lambda}(G)=\mathbb{Z}\left[  \tfrac{1}{\lambda}\right]  v \label{eq10.66}%
\end{equation}
where $v$ is the right Perron--Frobenius eigenvector given by
\textup{(\ref{eqSubNew.pound}).}
\end{lemma}

\begin{proof}
Since $J^{-n}v=\lambda^{-n}v$, it is clear that
\begin{equation}
\mathbb{Z}\left[  \tfrac{1}{\lambda}\right]  v\subseteq D_{\lambda}(G)\text{.}
\label{eq10.67}%
\end{equation}
Conversely, if $g\in D_{\lambda}(G)\subseteq G$, then $m_{N}^{i}g\in
\mathbb{Z}^{N}$ for some $N$ since
\[
m_{N}=(-1)^{N+1}\det(J).
\]
But as $\operatorname*{Prim}(\lambda)=\operatorname*{Prim}(m_{N})$, it follows
that $\lambda^{j}g\in\mathbb{Z}^{N}$ for some $j$. But $D_{\lambda}(G)$ has
$\operatorname*{rank}1$ and $v\in D_{\lambda}(G)$, so $D_{\lambda}%
(G)\subseteq\mathbb{Q}v$. Thus
\[
\lambda^{j}g\in\mathbb{Q}v\cap\mathbb{Z}^{N}\text{.}%
\]
But as $v\in\mathbb{Z}^{N}$ and the first component of $v$ is $1$, it follows
that $\mathbb{Q}v\cap\mathbb{Z}^{N}=\mathbb{Z}v$. Thus $\lambda^{j}g=nv$ for
an $n\in\mathbb{Z}$, so
\[
g=\frac{n}{\lambda^{j}}v\in\mathbb{Z}\left[  \tfrac{1}{\lambda}\right]
v\text{.}%
\]
This together with (\ref{eq10.67}) proves the lemma.
\end{proof}

\begin{corollary}
\label{corollary10.17}Let $J$, $J^{\prime}$ be matrices of the form
\textup{(\ref{eqCycNew.1})} with integer Perron--Frobenius eigenvalues
$\lambda$, $\lambda^{\prime}$, and let $\alpha$, $\alpha^{\prime}$,
resp.$\ v$, $v^{\prime}$, be the left, resp.\ right, Perron--Frobenius
eigenvectors of $J$, $J^{\prime}$ given by \textup{(\ref{eqRemNewBis.17}),
(\ref{eqSubNew.pound}),} respectively. Assume that\label{LOSPrim_4}
\begin{equation}
\operatorname*{Prim}(\lambda)=\operatorname*{Prim}(\lambda^{\prime
})=\operatorname*{Prim}(m_{N})=\operatorname*{Prim}(m_{N^{\prime}}^{\prime
})=\{p_{1},\ldots,p_{m}\} \label{eq10.68}%
\end{equation}
and that\label{LOSrank_2}
\begin{equation}
\operatorname*{rank}(D_{\lambda^{\prime}}(K_{0}(\mathfrak{A}_{J^{\prime}%
})))=1\text{.} \label{eq10.69}%
\end{equation}
If there exists a unital morphism $\mathfrak{A}_{J}\rightarrow\mathfrak
{A}_{J^{\prime}}$, then
\begin{equation}
\frac{\left\langle \alpha\mid v\right\rangle }{\left\langle \alpha^{\prime
}\mid v^{\prime}\right\rangle }\in\mathbb{Z}\left[  \frac{1}{\lambda}\right]
\label{eq10.70}%
\end{equation}
and if this morphism is an isomorphism, there exist integers $n_{1}%
,\ldots,n_{m}$ such that
\begin{equation}
\frac{\left\langle \alpha\mid v\right\rangle }{\left\langle \alpha^{\prime
}\mid v^{\prime}\right\rangle }=p_{1}^{n_{1}}p_{2}^{n_{2}}\cdots p_{m}^{n_{m}%
}\text{.} \label{eq.10.71}%
\end{equation}
\end{corollary}

\begin{proof}
If there exists an isomorphism, then
\begin{align}
D_{\lambda}(K_{0}(\mathfrak{A}_{J}))  &  =D_{\lambda^{\prime}}(K_{0}%
(\mathfrak{A}_{J}))\label{eq10.72}\\
&  \cong D_{\lambda^{\prime}}(K_{0}(\mathfrak{A}_{J^{\prime}}))\nonumber
\end{align}
and hence it follows from (\ref{eq10.69}) that $\operatorname*{rank}%
D_{\lambda}(K_{0}(\mathfrak{A}_{J}))$. The rest is straightforward from Lemma
\ref{lemma10.16} and Theorem \ref{theorem10.15}.
\end{proof}

In applying Corollary \ref{corollary10.17}, the most difficult condition to
verify is of course (\ref{eq10.69}). To this end, the following criterion is
often useful.

\begin{proposition}
\label{prop10.18}Let $J$ be a matrix of the form \textup{(\ref{eqCycNew.1})}
with integer Perron--Frobenius eigenvalue $\lambda$, and let $v$ be the
corresponding right eigenvector given by \textup{(\ref{eqSubNew.pound}).} Put
\[
G_{i}=J^{-i}\mathbb{Z}^{N}\text{ and }G=\bigcup_{i}G_{i}%
\]
as usual. Assume that there is an $n\in\mathbb{N}$ such that
\begin{equation}
\{g\in G\mid\lambda^{ni}g\in\mathbb{Z}^{N}\}\subseteq G_{ni} \label{eq10.73}%
\end{equation}
for all $i\in\mathbb{N}$, and assume that
\begin{equation}
\lambda^{n}J^{-n}\text{ has integer entries.} \label{eq10.74}%
\end{equation}
It follows that
\begin{equation}
D_{\lambda}(G)=\mathbb{Z}\left[  \tfrac{1}{\lambda}\right]  v\text{,}
\label{eq10.75}%
\end{equation}
so in particular \textup{(\ref{eq10.69})} is fulfilled.
\end{proposition}

\begin{proof}
It suffices to show that
\begin{equation}
D_{\lambda}(G)\cap\mathbb{Z}^{N}=\mathbb{Z}v\text{.} \label{eq10.76}%
\end{equation}
For this, let $w\in D_{\lambda}(G)\cap\mathbb{Z}^{N}$. Then $w=\lambda
^{ni}g_{i}$ for a $g_{i}\in G$ for all $i\in\mathbb{N}$. Using (\ref{eq10.73})
we have
\[
g_{i}\in G_{ni}=J^{-ni}\mathbb{Z}^{N}%
\]
so
\[
w=\lambda^{ni}g_{i}\in(\lambda^{n}J^{-n})^{i}\mathbb{Z}^{N}\text{.}%
\]
Thus
\[
w\in\bigcap_{i\geq0}(\lambda^{n}J^{-n})^{i}\mathbb{Z}^{N}\text{.}%
\]
But $\lambda^{n}$, being the Perron--Frobenius eigenvalue of the primitive
matrix $J^{n}$, is strictly larger in absolute value than any other
eigenvalue, and since $\lambda^{n}J^{-n}$ is a matrix with integer matrix
elements, it follows that
\[
\bigcap_{i\geq0}(\lambda^{n}J^{-n})^{i}\mathbb{Z}^{N}\subseteq\mathbb{R}%
v\cap\mathbb{Z}^{N}\text{.}%
\]
But since the first component of $v$ is 1 by (\ref{eqSubNew.pound}), it
follows that
\[
\mathbb{R}v\cap\mathbb{Z}^{N}=\mathbb{Z}v\cap\mathbb{Z}^{N}%
\]
so
\[
w\in\mathbb{Z}v\cap\mathbb{Z}^{N}%
\]
and this proves (\ref{eq10.76}) and thereby (\ref{eq10.75})$.$
\end{proof}

\begin{corollary}
\label{corollary10.19}Let $J$, $J^{\prime}$ be matrices of the form
\textup{(\ref{eqCycNew.1})} with integer Perron--Frobenius eigenvalues
$\lambda$, $\lambda^{\prime}$ and let $\alpha$, $\alpha^{\prime}$, resp.$\ v$,
$v^{\prime}$, be the left, resp.\ right, Perron--Frobenius eigenvectors of
$J$, $J^{\prime}$ given by \textup{(\ref{eqRemNewBis.17}),
(\ref{eqSubNew.pound}),} respectively. Assume that\label{LOSPrim_5}
\begin{equation}
\operatorname*{Prim}(\lambda)=\operatorname*{Prim}(\lambda^{\prime
})=\operatorname*{Prim}(m_{N})=\operatorname*{Prim}(m_{N}^{\prime}%
)=\{p_{1},\ldots,p_{m}\} \label{eq10.76bis}%
\end{equation}
and there is an $n\in\mathbb{N}$ such that
\begin{equation}
\{g\in G^{\prime}\mid\lambda^{ni}g\in\mathbb{Z}^{N}\}\subseteq G_{ni}^{\prime}
\label{eq10.77}%
\end{equation}
for all $i\in\mathbb{N}$, where
\[
G_{ni}^{\prime}=J^{\prime\,-ni}\mathbb{Z}^{N}\text{,}%
\]
and
\begin{equation}
\lambda^{\prime\,n}J^{\prime\,-n}\text{\qquad has integer entries.}
\label{eq10.78}%
\end{equation}

If there exists a unital morphism $\mathfrak{A}_{J}\rightarrow\mathfrak
{A}_{J^{\prime}}$, then
\begin{equation}
\frac{\left\langle \alpha\mid v\right\rangle }{\left\langle \alpha^{\prime
}\mid v^{\prime}\right\rangle }\in\mathbb{Z}\left[  \tfrac{1}{\lambda}\right]
\label{eq10.80}%
\end{equation}
and if this morphism is an isomorphism, there exist integers $n_{1}%
,\ldots,n_{m}$ such that
\begin{equation}
\frac{\left\langle \alpha\mid v\right\rangle }{\left\langle \alpha^{\prime
}\mid v^{\prime}\right\rangle }=p_{1}^{n_{1}}p_{2}^{n_{2}}\cdots p_{m}^{n_{m}%
}\text{.} \label{eq10.81}%
\end{equation}
\end{corollary}

\begin{proof}
This follows from Corollary \ref{corollary10.17} and Proposition
\ref{prop10.18}.
\end{proof}

\begin{remark}
\label{remark10.20}Note that condition \textup{(\ref{eq10.78})} is equivalent
to $\lambda^{\prime\,in}J^{\prime\,-in}$ having integer entries for
$i=1,2,\ldots$ and hence this condition is equivalent to
\[
G_{ni}^{\prime}\subseteq\{g\in G^{\prime}\mid\lambda^{ni}g\in\mathbb{Z}^{N}\}
\]
for $i=1,2,\ldots$. Thus, conditions \textup{(\ref{eq10.77})} and
\textup{(\ref{eq10.78})} taken together are equivalent to the single
condition
\begin{equation}
\{g\in G^{\prime}\mid\lambda^{ni}g\in\mathbb{Z}^{N}\}=G_{ni}^{\prime}
\label{eq10.82}%
\end{equation}
for $i=1,2,\ldots$.

Now one proves as in Proposition \textup{\ref{prop10.5}} that
\textup{(\ref{eq10.77})} is equivalent to the same condition for $i=1$, and
hence \textup{(\ref{eq10.82})} is equivalent to the single condition
\begin{equation}
\{g\in G^{\prime}\mid\lambda^{n}g\in\mathbb{Z}^{N}\}=G_{n}^{\prime}\text{.}
\label{eq10.83}%
\end{equation}

By the same token, \textup{(\ref{eq10.73})} $\wedge$ \textup{(\ref{eq10.74})}
is equivalent to \textup{(\ref{eq10.83}).}
\end{remark}

\begin{scholium}
\label{scholium10.21}We saw in Lemma \textup{\ref{lemma10.13}} that if
$\lambda=m_{N}$, condition \textup{(\ref{eq10.83})} is automatically fulfilled
with $n=1$, and this was used in the proof of Corollary
\textup{\ref{corollary10.12}.} Thus, Corollary \textup{\ref{corollary10.12}}
may be viewed as a special case of Corollary \textup{\ref{corollary10.19}}.
But in order to verify the hypotheses \textup{(\ref{eq10.83}),} we need an
efficient algorithm. One such algorithm is one that in particularly ``pure''
form occurs in the proof of Proposition \textup{\ref{prop10.14}.} So in terms
of $G=K_{0}(\mathfrak{A}_{L})$,\label{LOSFL_3}
\begin{align*}
F  &  =G\diagup\mathbb{Z}^{N}\\
F_{i}  &  =G_{i}\diagup\mathbb{Z}^{N}=J^{-i}\mathbb{Z}^{N}\diagup
\mathbb{Z}^{N}%
\end{align*}
we want to establish \textup{(\ref{eq10.83}),} i.e.,
\begin{equation}
\{h\in F\mid\lambda^{n}h=0\}=F_{n}\text{.} \label{eq10.84}%
\end{equation}

In the special case $\lambda^{n}=m_{N}=m$, we may proceed like this: If
$x_{1}$, $x_{2}$\label{LOSxi_3} are the generators of $F$ in Lemma
\textup{\ref{LemSubNew.14}, (\ref{eqStr.E})} and \textup{(\ref{eqStr.F}),}
then
\begin{equation}%
\begin{array}
[c]{l}%
mx_{1}=0\\
mx_{2}=-m_{N-1}x_{1}\\
mx_{3}=-m_{N-2}x_{1}-m_{N-1}x_{2}\\
\quad\vdots\qquad\vdots\\
mx_{k}=x_{k-N}-m_{1}x_{k-N+1}-\cdots-m_{N-1}x_{k-1}%
\end{array}
\label{eq10.85}%
\end{equation}
\textup{(}where we use $x_{l}=0$ for $l\leq0$\textup{).} As seen for example
in Corollary \textup{\ref{CorKer.C},} any element $h\in F$ can be uniquely
represented as
\begin{equation}
h=t_{1}x_{1}+t_{2}x_{2}+\cdots+t_{k}x_{k} \label{eq10.86}%
\end{equation}
for some $k$, where $t_{i}\in\{0,1,\ldots,m-1\}$. It follows from
\textup{(\ref{eq10.85})} that $mh=0$ is equivalent to
\begin{equation}
-t_{2}m_{N-1}x_{1}+\cdots+t_{k}(x_{k-N}-m_{1}x_{k-N+1}-\cdots-m_{N-1}%
x_{k-1})=0\text{.} \label{eq10.87}%
\end{equation}
The leading term here is $-t_{k}m_{N-1}x_{k-1}$. If now
\[
t_{k}m_{N-1}\neq0\mod{m}
\]
for $t_{k}=1,2,\ldots,m-1$, it follows that $mh=0\Rightarrow t_{k}=0$. Thus
$t_{k}x_{k}$ cannot be the leading term in $h$. Continuing in this manner, one
loops off $t_{k-1}$, $t_{k-2}$, etc., and one ends up showing
$mh=0\Leftrightarrow h\in F_{1}$. If now $m_{N-1}=0\mod{m}$, one gets from the
outset that $mx_{2}=0$, and hence there are no restrictions on $t_{2}$, and
since now $m_{N-1}x_{k-1}$ can be expanded in $x_{k-2},x_{k-3}$, etc., one may
try to go one step further and show that $t_{k}=0$, etc., and then $h\in
F_{2}$ if $tm_{_{N-2}}\neq0\mod{m}$ for $t\in\{1,\ldots,m-1\}$, etc. In
general, it may of course happen for example that $tm_{N-2}\neq0\mod{m}$ for
some $t\in\{1,\ldots,m-1\}$ and $tm_{N-2}=0\mod{m}$ for some other $t$, and
then the computation of $\{h\in F\mid mh=0\}$ becomes much more complicated.
Let us single out the simple case.
\end{scholium}

\begin{proposition}
\label{prop10.22}Let $J$ be a matrix of the form \textup{(\ref{eqCycNew.1})}
and put $m=m_{N}$.\label{LOSmN_1} Choose $n\in\{1,\ldots,N\}$ such that
\begin{equation}
m_{N-k}=0\mod{m_{N}} \label{eq10.88}%
\end{equation}
for $k=0,\ldots,n-1$, and assume that
\begin{equation}
\gcd(m_{N-n},m_{N})=1\text{.} \label{eq10.89}%
\end{equation}
It follows that
\begin{equation}
\{h\in F\mid mh=0\}=F_{n}\text{.} \label{eq10.90}%
\end{equation}
If $n=N,$ condition \textup{(\ref{eq10.89})} may be omitted.
\end{proposition}

\begin{proof}
Once one notes that condition (\ref{eq10.89}),
\[
\gcd(m_{N-n},m_{N})=1\text{,}%
\]
means that
\[
tm_{N-n}\neq0\mod{m_{N}}%
\]
for $t=1,2,\ldots,m_{N-1}$, this is clear from the discussion preceding the proposition.
\end{proof}

\begin{example}
\label{example10.23}Consider
\[
J=\left(
\begin{array}
[c]{cccc}%
0 & 1 & 0 & 0\\
1 & 0 & 1 & 0\\
2 & 0 & 0 & 1\\
8 & 0 & 0 & 0
\end{array}
\right)  \text{,}\qquad J^{\prime}=\left(
\begin{array}
[c]{cccc}%
1 & 1 & 0 & 0\\
1 & 0 & 1 & 0\\
0 & 0 & 0 & 1\\
4 & 0 & 0 & 0
\end{array}
\right)  \text{.}%
\]
Here $\lambda=\lambda^{\prime}=2$, so
\[
\operatorname*{Prim}(\lambda)=\operatorname*{Prim}(\lambda^{\prime
})=\operatorname*{Prim}(m_{N})=\operatorname*{Prim}(m_{N}^{\prime
})=\{2\}\text{.}%
\]
Note that
\[
4J^{\prime\,-2}=\left(
\begin{array}
[c]{cccc}%
0 & 0 & 1 & 0\\
0 & 0 & -1 & 1\\
4 & 0 & -1 & -1\\
0 & 4 & 0 & -1
\end{array}
\right)
\]
so \textup{(\ref{eq10.74})} is satisfied with $n=2$. Further note that
\[
\left\langle \alpha\mid v\right\rangle =\left(  1,\frac{1}{2},\frac{1}%
{4},\frac{1}{8}\right)  \left(
\begin{array}
[c]{c}%
1\\
2\\
3\\
4
\end{array}
\right)  =\frac{13}{4}%
\]
and
\[
\left\langle a^{\prime}\mid v^{\prime}\right\rangle =\left(  1,\frac{1}%
{2},\frac{1}{4},\frac{1}{8}\right)  \left(
\begin{array}
[c]{c}%
1\\
1\\
1\\
2
\end{array}
\right)  =2
\]
and hence
\[
\frac{\left\langle \alpha\mid v\right\rangle }{\left\langle \alpha^{\prime
}\mid v^{\prime}\right\rangle }=\frac{13}{2}\in\mathbb{Z}\left[  \frac{1}%
{2}\right]  \text{,}%
\]
but this number does not have a multiplicative inverse in $\mathbb{Z}\left[
\frac{1}{2}\right]  $. Thus, if we can show \textup{(\ref{eq10.77}),}
\[
\{g\in G^{\prime}\mid\lambda^{2i}g\in\mathbb{Z}^{N}\}\subseteq G_{2i}^{\prime
}\text{,}%
\]
it follows from Corollary \textup{\ref{corollary10.19}} that $\mathfrak{A}%
_{J}$ and $\mathfrak{A}_{J^{\prime}}$ are non-isomorphic. But Proposition
\textup{\ref{prop10.22}} applied to $J^{\prime}$ shows that
\[
\{h\in F\mid4h=0\}=F_{2}\text{.}%
\]
Thus \textup{(\ref{eq10.83})} holds with $n=2$ and in particular
\textup{(\ref{eq10.77})} is valid. This shows the non-isomorphism of the two
dimension groups.
\end{example}

In addition to the criteria of non-isomorphism given by Theorem
\ref{theorem10.10}, Corollary \ref{corollary10.11}, Corollary
\ref{corollary10.12}, Theorem \ref{theorem10.15}, Corollary
\ref{corollary10.17} and Proposition \ref{prop10.18}, it is frequently
possible to decide non-isomorphism by another route, namely by establishing
that the exact sequence (\ref{eqRemNewBis.31}):\label{LOSkertau_8}%
\label{LOSK0AL_9}
\[
0\longrightarrow\ker\tau\hooklongrightarrow K_{0}(\mathfrak{A}_{J}%
)\overset{\tau}{\longrightarrow}Z\left[  \tfrac{1}{\lambda}\right]
\longrightarrow0
\]
splits for one specimen but not for another. With $K_{0}(\mathfrak{A}_{J})$
realized as (\ref{eqRemNewBis.19}), there is a simple criterion for this.

\begin{lemma}
\label{lemma10.24}Let $J$ be a matrix of the form \textup{(\ref{eqCycNew.1})}
with integer Perron--Frobenius eigenvalue $\lambda$ with
\[
\operatorname*{Prim}(\lambda)=\{p_{1},\ldots,p_{m}\}\text{.}%
\]
The following conditions \textup{(i)} and \textup{(ii)} are equivalent.

\begin{enumerate}
\item [(i)]The exact sequence
\[
0\longrightarrow\ker\tau\hooklongrightarrow K_{0}(\mathfrak{A}_{J}%
)\overset{\tau}{\longrightarrow}\mathbb{Z}\left[  \tfrac{1}{\lambda}\right]
\longrightarrow0
\]
splits.

\item[(ii)] There is a $w\in D_{\lambda}(G)$ such that $\tau(w)$ has a
multiplicative inverse in $Z\left[  \frac{1}{\lambda}\right]  $, i.e.,
$\tau(w)=p_{1}^{n_{1}}p_{2}^{n_{2}}\cdots p_{m}^{n_{m}}$ for suitable
$n_{1},\ldots,n_{m}\in\mathbb{Z}$.
\end{enumerate}
\end{lemma}

\begin{proof}
(i)$\Rightarrow$(ii): If the sequence splits, let
\[
\psi\colon\mathbb{Z}\left[  \tfrac{1}{\lambda}\right]  \longrightarrow
K_{0}(\mathfrak{A}_{J})
\]
be a section, and put $w=\psi(1)$. Then
\[
w\in\psi\left(  \mathbb{Z}\left[  \tfrac{1}{\lambda}\right]  \right)
\subseteq D_{\lambda}(G)\text{ and }\tau(w)=1\text{.}%
\]

(ii)$\Rightarrow$(i): If (ii) holds, we can define a section $\psi$ by
\[
\psi(\tau(w))=w
\]
and since $\mathbb{Z}\left[  \frac{1}{\lambda}\right]  \tau(w)=\mathbb{Z}%
\left[  \frac{1}{\lambda}\right]  $, and $w\in D_{\lambda}(G)$, $\psi$ extends
uniquely to $\mathbb{Z}\left[  \frac{1}{\lambda}\right]  $.
\end{proof}

\begin{corollary}
\label{corollary10.25}Let $J$ be a matrix of the form
\textup{(\ref{eqCycNew.1})} with integer Perron--Frobenius eigenvalue
$\lambda$ with $\operatorname*{Prim}(\lambda)=\left\{  p_{1},\ldots
,p_{m}\right\}  $ and left and right Frobenius eigenvectors $\alpha,v$ given
by \textup{(\ref{eqRemNewBis.17}), (\ref{eqSubNew.pound}),} respectively. The
exact sequence \textup{(\ref{eqRemNewBis.31}):}
\begin{equation}
0\longrightarrow\ker\tau\hooklongrightarrow K_{0}(\mathfrak{A}_{J}%
)\overset{\tau}{\longrightarrow}\mathbb{Z}\left[  \tfrac{1}{\lambda}\right]
\rightarrow0 \label{eq10.91}%
\end{equation}
splits if
\begin{equation}
\left\langle \alpha\mid v\right\rangle =p_{1}^{n_{1}}p_{2}^{n_{2}}\cdots
p_{m}^{n_{m}} \label{eq10.92}%
\end{equation}
for suitable $n_{1},n_{2},\ldots,n_{m}$ in $\mathbb{Z}$.

Conversely, if \textup{(\ref{eq10.91})} splits and $\operatorname*{rank}%
(D_{\lambda}(G))=1$ and $\operatorname*{Prim}(m_{N})=\operatorname*{Prim}%
(\lambda)$, then $\left\langle \alpha\mid v\right\rangle $ has the form
\textup{(\ref{eq10.92}).}
\end{corollary}

\begin{proof}
Since always $v\in D_{\lambda}(G)$, the first statement follows from Lemma
\ref{lemma10.24}.

The last statement follows from Lemma \ref{lemma10.16}: If
$\operatorname*{rank}(D_{\lambda}(G))=1$ and $\operatorname*{Prim}%
(m_{N})=\operatorname*{Prim}(\lambda)$, then $D_{\lambda}(G)=Z\left[  \frac
{1}{\lambda}\right]  v$ and if $\psi$ is a section for (\ref{eq10.91}), then
$\psi(1)\in D_{\lambda}(G)=Z\left[  \frac{1}{\lambda}\right]  v$. Thus, there
is a $t\in Z\left[  \frac{1}{\lambda}\right]  $ such that $\psi(1)=tv$. But
then
\[
1=\left\langle \alpha\mid\psi(1)\right\rangle =\left\langle \alpha\mid
tv\right\rangle =t\left\langle \alpha\mid v\right\rangle
\]
and it follows that $\left\langle \alpha\mid v\right\rangle $ has the
multiplicative inverse $t$ in $Z\left[  \frac{1}{\lambda}\right]  $, so
$\left\langle \alpha\mid v\right\rangle $ has the form (\ref{eq10.92}).
\end{proof}

\begin{example}
\label{ExaCyc.end1bis}An instance where the isomorphism question for two
$C^{\ast}$-algebras $\mathfrak{A}_{J}$ and $\mathfrak{A}_{J^{\prime}}$ is not
immediate is when%
\[
J=%
\begin{pmatrix}
1 & 1 & 0\\
1 & 0 & 1\\
2 & 0 & 0
\end{pmatrix}
\text{\qquad and\qquad}J^{\prime}=%
\begin{pmatrix}
0 & 1 & 0\\
3 & 0 & 1\\
2 & 0 & 0
\end{pmatrix}
.
\]
Both matrices are regular, and both have determinant $2$ and Perron--Frobenius
eigenvalue $2$. The respective right eigenvectors are
\[
v_{J}=%
\begin{pmatrix}
1\\
1\\
1
\end{pmatrix}
\text{\qquad and\qquad}v_{J^{\prime}}=%
\begin{pmatrix}
1\\
2\\
1
\end{pmatrix}
.
\]
Since
\[
\frac{\tau\left(  v_{J}\right)  }{\tau\left(  v_{J^{\prime}}\right)  }%
=\frac{7}{9},
\]
it follows from Corollary \textup{\ref{corollary10.12}} that the two $C^{\ast
}$-algebras $\mathfrak{A}_{J}$ and $\mathfrak{A}_{J^{\prime}}$ are
non-isomorphic. The non-isomorphism here is perhaps a bit surprising since the
two matrices $J^{6}$ and $J^{\prime\,6}$ have the same spectrum.

Later, in Theorem \textup{\ref{ThmClmN.2},} we will show that $\left\{
\operatorname*{Prim}\lambda,\lambda^{2}\tau\left(  v\right)  \right\}  $ is a
\emph{complete} isomorphism invariant for $3\times3$ matrices with
$\lambda=m_{3}$. For the two matrices above one computes that $\lambda^{2}%
\tau\left(  v\right)  $ is $7$, $9$, respectively, confirming that they are
non-isomorphic. See also Theorem \textup{\ref{ThmClmN.3}.}
\end{example}

\setcounter{figurelink}{\value{figure}} \setcounter{tablelink}{\value{table}}

\chapter{\label{SubG0}\label{LOSG0_8}Subgroups of $G_{0}=\bigcup_{n=0}%
^{\infty}J_{0}^{-n}\mathcal{L}$}

Before applying our general theory of isomorphism of stationary AF-algebras to
more specific examples in Chapters \ref{App.Cla}--\ref{App.Fur}, we will
mention one more example of how to decide nontriviality of extensions which is
sometimes useful. In many examples we compute that $G=K_{0}\left(
\mathfrak{A}_{L}\right)  $ or $G_{0}=\ker\tau$ or some other group is an
extension of the form%
\begin{equation}
0\longrightarrow H\longrightarrow G_{0}\longrightarrow\mathbb{Z}\left[
\tfrac{1}{d}\right]  \longrightarrow0. \label{eqSubG0.1}%
\end{equation}

We first state a proposition which is a variation of a result due to David
Handelman \cite[Proposition 10.1]{BrJo99a}.

\begin{proposition}
\label{ProSubG0.0}Let $\mathbb{E}$ be an $N\times N$ matrix with integer
entries and assume $\det\left(  \mathbb{E}\right)  \neq0$. Let $f\left(
x\right)  =\det\left(  x\openone-\mathbb{E}\right)  $ be the characteristic
polynomial of $\mathbb{E}$ and let $f\left(  x\right)  =f_{1}\left(  x\right)
f_{2}\left(  x\right)  \cdots f_{n}\left(  x\right)  $ be the decomposition of
$f$ into irreducible polynomials in $\mathbb{Z}\left[  x\right]  $. Define
\[
q\left(  x\right)  =\prod_{\substack{i\\\left|  f_{i}\left(  0\right)
\right|  =1}}f_{i}\left(  x\right)  ,\qquad p\left(  x\right)  =\prod
_{\substack{i\\\left|  f_{i}\left(  0\right)  \right|  \neq1}}f_{i}\left(
x\right)  .
\]
Then%
\[
p\left(  \mathbb{E}\right)  \mathbb{Z}^{N}\subseteq\left\{  m\in\mathbb{Z}%
^{N}\mid q\left(  \mathbb{E}\right)  m=0\right\}  =\bigcap_{k}\mathbb{E}%
^{k}\mathbb{Z}^{N}.
\]
\end{proposition}

\begin{proof}
The left inclusion follows from $q\left(  \mathbb{E}\right)  p\left(
\mathbb{E}\right)  =f\left(  \mathbb{E}\right)  =0$. Next note that
$W=\left\{  m\in\mathbb{Z}^{N}\mid q\left(  \mathbb{E}\right)  m=0\right\}  $
is an $\mathbb{E}$-invariant sublattice of $\mathbb{Z}^{N}$. Note that if
$q\left(  x\right)  =\sum_{i=0}^{k}q_{i}x^{i}$, then $q_{0}=\pm1$, so we may
assume $q_{0}=-1$, and hence $m\in W$ if and only if%
\[
m=\sum_{i=1}^{k}q_{i}\mathbb{E}^{i}m.
\]
But then by iteration%
\[
m=\left(  \sum_{i=1}^{k}q_{i}\mathbb{E}^{i}\right)  ^{l}m
\]
for $l=1,2,\dots$, and expanding those polynomials we see that
\[
m\in\bigcap_{l}\mathbb{E}^{l}\mathbb{Z}^{N}.
\]
Thus%
\[
W\subseteq V\equiv\bigcap_{i}\mathbb{E}^{i}\mathbb{Z}^{N}.
\]
But $V$ is also an $\mathbb{E}$-invariant sublattice of $\mathbb{Z}^{N}$, and
thus a free abelian group, and the restriction of $\mathbb{E}$ to $V$ is
clearly surjective. Since $\mathbb{E}$ is injective, it follows that
$\mathbb{E}|_{V}$ is invertible and thus $\left|  \det\left(  \mathbb{E}%
|_{V}\right)  \right|  =1$. But the characteristic polynomial of
$\mathbb{E}|_{V}$ is a factor of the characteristic polynomial of $\mathbb{E}%
$, and since the constant term of the former polynomial is $\pm1$, it is a
factor of $q\left(  x\right)  $. It follows that $q\left(  \mathbb{E}%
|_{V}\right)  =0$, which means $V\subseteq W$.
\end{proof}

In order that an extension such as (\ref{eqSubG0.1}) shall be trivial, it is
necessary that $G_{0}$ contain $\mathbb{Z}\left[  \frac{1}{d}\right]  $ as a
subgroup. In order to decide that, the following more local proposition is
sometimes useful. The condition on $\mathcal{L}$ means that $\mathcal{L}%
\cong\mathbb{Z}^{M}$ and that $\mathcal{L}$ spans $\mathbb{Q}^{M}$ over
$\mathbb{Q}$.

\begin{proposition}
\label{ProSubG0.1}\renewcommand{\theenumi}{\roman{enumi}}Let $J_{0}$ be a
nonsingular $M\times M$ matrix with integer matrix elements, and let
$\mathcal{L}$ be a free abelian subgroup of rank $M$ of $\mathbb{Q}^{M}$.
Consider the additive subgroup%
\begin{equation}
G_{0}=\bigcup_{n=0}^{\infty}J_{0}^{-n}\mathcal{L} \label{eqSubG0.2}%
\end{equation}
of $\mathbb{Q}^{M}$. Let $d\in\mathbb{N}$ be a number such that%
\begin{equation}
E=dJ_{0}^{-1} \label{eqSubG0.3}%
\end{equation}
is a matrix with integer matrix elements. Assume that

\begin{enumerate}
\item \label{ProSubG0.1(1)}there is a prime factor $f$ of the monic
polynomial
\[
\det\left(  \lambda\openone-E\right)  =\frac{\left(  -\lambda\right)  ^{M}%
}{\det J_{0}}\det\left(  \frac{d}{\lambda}\openone-J_{0}\right)
\]
such that $\left|  f\left(  0\right)  \right|  =1$.
\end{enumerate}

\noindent It follows that

\begin{enumerate}
\addtocounter{enumi}{1}

\item \label{ProSubG0.1(2)}$G_{0}$ contains a subgroup isomorphic to
$\mathbb{Z}\left[  \frac{1}{d}\right]  $
\end{enumerate}

\noindent but \textup{(\ref{ProSubG0.1(2)})} does not imply
\textup{(\ref{ProSubG0.1(1)})} in general.
\end{proposition}

\begin{remark}
\label{RemSubG0.2}Let us exhibit a partial example showing that
\textup{(\ref{ProSubG0.1(2)})} does not imply \textup{(\ref{ProSubG0.1(1)}).}
Let%
\[
J_{0}=%
\begin{pmatrix}
0 & 1 & 0 & \cdots & 0 & 0\\
0 & 0 & 1 & \cdots & 0 & 0\\
0 & 0 & 0 & \ddots & 0 & 0\\
\vdots & \vdots &  & \ddots & \ddots & \vdots\\
0 & 0 & 0 &  & 0 & 1\\
d & 0 & 0 & \cdots & 0 & 0
\end{pmatrix}
\]
and $\mathcal{L}=\mathbb{Z}^{M}$. Then $J_{0}^{M}=d\openone$, so
$G_{0}=\bigcup_{k}J_{0}^{-Mk}\mathbb{Z}^{M}=\mathbb{Z}\left[  \frac{1}%
{d}\right]  ^{M}$ and hence \textup{(\ref{ProSubG0.1(2)})} holds. But since
\[
\bigcap_{k}E^{k}\mathbb{Z}^{M}=\bigcap_{k}E^{Mk}\mathbb{Z}^{M}=\bigcap
_{k}d^{\left(  M-1\right)  k}\mathbb{Z}^{M}=0
\]
it follows from the equivalent form \textup{(\ref{ProSubG0.1(3)})} of
\textup{(\ref{ProSubG0.1(1)})} in the proof below that
\textup{(\ref{ProSubG0.1(1)})} does not hold.
\end{remark}

\begin{proof}
\renewcommand{\theenumi}{\roman{enumi}}We know from Proposition
\ref{ProSubG0.0} that (\ref{ProSubG0.1(1)}) is equivalent to

\begin{enumerate}
\addtocounter{enumi}{2}

\item \label{ProSubG0.1(3)}$\bigcap_{k=1}^{\infty}E^{k}\mathbb{Z}^{M}\neq0$.
\end{enumerate}

\noindent We now argue that (\ref{ProSubG0.1(3)}) is equivalent to

\begin{enumerate}
\addtocounter{enumi}{3}

\item \label{ProSubG0.1(4)}$\bigcap_{k=1}^{\infty}E^{k}\mathcal{L}\neq0$.
\end{enumerate}

\noindent But since both $\mathcal{L}$ and $\mathbb{Z}^{M}$ are the free
abelian groups generated by $M$ points in $\mathbb{Q}^{M}$, and there is an
element of $\operatorname*{GL}\left(  \mathbb{Q},M\right)  $ transforming
these $M$-tuples into each other, it follows that there is a natural number
$n$ such that
\begin{equation}
\mathcal{L}\subseteq\frac{1}{n}\mathbb{Z}^{M}\text{\quad and\quad}%
\mathbb{Z}^{M}\subseteq\frac{1}{n}\mathcal{L}. \label{eqSubG0.4}%
\end{equation}
Thus the equivalence of (\ref{ProSubG0.1(3)}) and (\ref{ProSubG0.1(4)})
follows by linearity. Next, put%
\begin{equation}
H=\bigcap_{k\geq0}d^{k}G_{0}. \label{eqSubG0.5}%
\end{equation}
Clearly, $H$ is a subgroup of $G_{0}$ containing $0$. We now show%
\begin{equation}
\bigcap_{k}^{\infty}E^{k}\mathcal{L}\subseteq H. \label{eqSubG0.6}%
\end{equation}
But this follows from%
\[
\bigcap_{k}^{\infty}E^{k}\mathcal{L}=\bigcap_{k}^{\infty}d^{k}J_{0}%
^{-k}\mathcal{L}\subseteq\bigcap_{k}^{\infty}d^{k}G_{0}=H.
\]
But since \textup{(\ref{ProSubG0.1(1)}) }$\Leftrightarrow$
(\ref{ProSubG0.1(4)}), it follows from \textup{(\ref{ProSubG0.1(1)}) and
(\ref{eqSubG0.6}) }that $H\neq0$. But if $g\in H$, then $d^{-k}g\in G_{0}$ for
all $k$ by (\ref{eqSubG0.5}). Thus $G_{0}$ contains a subgroup isomorphic to
$\mathbb{Z}\left[  \frac{1}{d}\right]  $.
\end{proof}

\chapter{\label{App.Cla}Classification of the AF-algebras $\mathfrak{A}_{L}$
with $\operatorname*{rank}\left(  K_{0}\left(  \mathfrak{A}_{L}\right)
\right)  =2$}

\setcounter{figure}{\value{figurelink}} \setcounter{table}{\value{tablelink}}
Let us consider matrices of the form (\ref{eqRemNewBis.40}) with $N=2$,%
\begin{equation}
J=%
\begin{pmatrix}
m_{1} & 1\\
m_{2} & 0
\end{pmatrix}
, \label{eqApp.Cla1}%
\end{equation}
where $m_{1},m_{2}\in\mathbb{N}$. We divide the discussion into two
cases.\setcounter{case}{0}

\begin{case}
\label{CasApp.Cla1}The Perron--Frobenius eigenvalue $\lambda$ is rational, and
thus $\lambda\in\mathbb{N}$.
\end{case}

In this case one computes that $J$ has the form%
\begin{equation}
J=%
\begin{pmatrix}
\lambda-k & 1\\
k\lambda & 0
\end{pmatrix}
,\qquad k=1,\dots,\lambda-1, \label{eqApp.Cla2}%
\end{equation}
and the spectrum is
\begin{equation}
\operatorname*{spec}\left(  J\right)  =\left\{  -k,\lambda\right\}  \text{.}
\label{eqApp.Cla3}%
\end{equation}
Referring to Theorem \textup{\ref{ThmCyc.7},} we have $D=1$, $N=2$, and the
triangular form \textup{(\ref{eqCyc.22})} is \textup{(}with $p_{a}\left(
x\right)  =\lambda x-1$\textup{):}%
\begin{equation}%
\begin{pmatrix}
-k & k\\
0 & \lambda
\end{pmatrix}
. \label{eqApp.Cla4}%
\end{equation}
Hence the invariants of Theorem \textup{\ref{CorCyc.9}} are:

\begin{enumerate}
\item \label{App.Cla(1)}$N=2$,

\item \label{App.Cla(2)}$\operatorname*{Prim}\left(  k\lambda\right)  $,\label{LOSPrim_6}

\item \label{App.Cla(3)}$\operatorname*{Prim}\left(  k\right)  $,

\item[\textup{(\ref{App.Cla(3)})}$^{\prime}$] $\operatorname*{Prim}\left(
\lambda\right)  $,

\item \label{App.Cla(4)}$D=1$.
\end{enumerate}

Furthermore, we will argue below that%
\begin{equation}
\ker\left(  \tau\right)  \cong\mathbb{Z}\left[  \tfrac{1}{k}\right]
,\qquad\tau\left(  K_{0}\left(  \mathfrak{A}\right)  \right)  =\mathbb{Z}%
\left[  \tfrac{1}{\lambda}\right]  , \label{eqApp.Cla6}%
\end{equation}
so $K_{0}\left(  \mathfrak{A}\right)  $ is an extension of $\mathbb{Z}\left[
\tfrac{1}{\lambda}\right]  $ by $\mathbb{Z}\left[  \tfrac{1}{k}\right]  $:%
\begin{equation}
0\longrightarrow\mathbb{Z}\left[  \tfrac{1}{k}\right]  \longrightarrow
K_{0}\left(  \mathfrak{A}\right)  \longrightarrow\mathbb{Z}\left[  \tfrac
{1}{\lambda}\right]  \longrightarrow0. \label{eqApp.Cla7}%
\end{equation}
To see this, one computes%
\begin{align}%
\begin{pmatrix}
\lambda-k & 1\\
k\lambda & 0
\end{pmatrix}
^{-n}  &  =\frac{1}{\lambda+k}%
\begin{pmatrix}
1 & 1\\
k & -\lambda
\end{pmatrix}%
\begin{pmatrix}
\lambda^{-n} & 0\\
0 & \left(  -k\right)  ^{-n}%
\end{pmatrix}%
\begin{pmatrix}
\lambda & 1\\
k & -1
\end{pmatrix}
\label{eqApp.Cla8}\\
&  =\frac{1}{\lambda+k}%
\begin{pmatrix}
\lambda^{-n+1}-\left(  -k\right)  ^{-n+1} & \lambda^{-n}-\left(  -k\right)
^{-n}\\
k\lambda^{-n+1}+\lambda\left(  -k\right)  ^{-n+1} & k\lambda^{-n}%
+\lambda\left(  -k\right)  ^{-n}%
\end{pmatrix}
\nonumber\\
&  =\frac{1}{\lambda+k}%
\begin{pmatrix}
\lambda^{-n+1}-\left(  -k\right)  ^{-n+1} & \lambda^{-n}-\left(  -k\right)
^{-n}\\
k\lambda\left(  \lambda^{-n}-\left(  -k\right)  ^{-n}\right)  & k\lambda
\left(  \lambda^{-n-1}-\left(  -k\right)  ^{-n-1}\right)
\end{pmatrix}
.\nonumber
\end{align}
Hence, using \textup{(\ref{eqApp.Cla8}),} one computes%
\begin{align}
g  &  =%
\begin{pmatrix}
\lambda-k & 1\\
k\lambda & 0
\end{pmatrix}
^{-n}%
\begin{pmatrix}
n_{1}\\
n_{2}%
\end{pmatrix}
\label{eqApp.Cla9}\\
&  =\frac{1}{\lambda+k}\left[  \left(  \lambda n_{1}+n_{2}\right)
\lambda^{-n}%
\begin{pmatrix}
1\\
k
\end{pmatrix}
+\left(  kn_{1}-n_{2}\right)  \left(  -k\right)  ^{-n}%
\begin{pmatrix}
1\\
-\lambda
\end{pmatrix}
\right]  ,\nonumber
\end{align}
and thus, using \textup{(\ref{eqRemNewBis.17}),}%
\begin{equation}
\tau\left(  g\right)  =\left(  1,\frac{1}{\lambda}\right)  g=\left(  \lambda
n_{1}+n_{2}\right)  \lambda^{-n-1}. \label{eqApp.Cla10}%
\end{equation}
This confirms \textup{(\ref{eqRemNewBis.22}):} $\tau\left(  K_{0}\left(
\mathfrak{A}\right)  \right)  =\mathbb{Z}\left(  \frac{1}{\lambda}\right)  $,
and we see that $g\in\ker\left(  \tau\right)  $ if and only if $\lambda
n_{1}+n_{2}=0$, i.e., $n_{2}=-\lambda n_{1}$, so $g\in\ker\left(  \tau\right)
$ if and only if%
\begin{equation}
g=\frac{1}{\lambda+k}\left(  kn_{1}+\lambda n_{1}\right)  \left(  -k\right)
^{-n}%
\begin{pmatrix}
1\\
-\lambda
\end{pmatrix}
=n_{1}\left(  -k\right)  ^{-n}%
\begin{pmatrix}
1\\
-\lambda
\end{pmatrix}
\label{eqApp.Cla11}%
\end{equation}
for an $n_{1}\in\mathbb{Z}$, $n\in\mathbb{N}$, which confirms $\ker\left(
\tau\right)  \cong\mathbb{Z}\left[  \frac{1}{k}\right]  $, so the sequence
\textup{(\ref{eqApp.Cla7})} is well defined and exact.

Now, using \textup{(\ref{eqApp.Cla9})} we can prove

\begin{proposition}
\label{ProApp.Cla1}If $G=K_{0}\left(  \mathfrak{A}\right)  $ is realized
concretely in $\mathbb{Q}^{2}$ as above we have%
\begin{equation}
\mathbb{Z}\left[  \tfrac{1}{\lambda}\right]
\begin{pmatrix}
1\\
k
\end{pmatrix}
+\mathbb{Z}\left[  \tfrac{1}{k}\right]
\begin{pmatrix}
1\\
-\lambda
\end{pmatrix}
\subset G\subset\frac{1}{\left(  \lambda+k,\lambda\right)  }\mathbb{Z}\left[
\tfrac{1}{\lambda}\right]
\begin{pmatrix}
1\\
k
\end{pmatrix}
+\frac{1}{\left(  \lambda+k,k\right)  }\mathbb{Z}\left[  \tfrac{1}{k}\right]
\begin{pmatrix}
1\\
-\lambda
\end{pmatrix}
, \label{eqApp.Cla12}%
\end{equation}
where%
\begin{equation}
\left(  n_{1},n_{2}\right)  =\frac{n_{1}}{\gcd\left(  n_{1}^{{}},n_{2}%
^{\infty}\right)  } \label{eqApp.Cla13}%
\end{equation}
for $n_{1},n_{2}\in\mathbb{N}$, where $\gcd\left(  n_{1}^{{}},n_{2}^{\infty
}\right)  $ is the \textup{(}unique\/\textup{)} greatest common divisor of
$n_{1}$ and $n_{2}^{m}$ for large $m$. Furthermore an element $a%
\begin{pmatrix}
1\\
k
\end{pmatrix}
+b%
\begin{pmatrix}
1\\
-\lambda
\end{pmatrix}
$ of $G$ is \textup{(}nonzero\/\textup{)} positive if and only if%
\begin{equation}
a>0. \label{eqApp.Cla14}%
\end{equation}

In particular, the following conditions are equivalent:

\begin{enumerate}
\item \label{ProApp.Cla1(1)}$G=\mathbb{Z}\left[  \frac{1}{\lambda}\right]
\begin{pmatrix}
1\\
k
\end{pmatrix}
+\mathbb{Z}\left[  \frac{1}{k}\right]
\begin{pmatrix}
1\\
-\lambda
\end{pmatrix}
$,

and

\item \label{ProApp.Cla1(2)}$\operatorname*{Prim}\left(  \lambda+k\right)
\subset\operatorname*{Prim}\left(  \lambda\right)  \cap\operatorname*{Prim}%
\left(  k\right)  $.
\end{enumerate}

$G$ can also in general be characterized as the set of linear combinations of
elements of $\mathbb{Z}^{2}$ and the elements%
\begin{equation}
\frac{1}{\lambda+k}%
\begin{pmatrix}
\lambda^{-n}-\left(  -k\right)  ^{-n}\\
k\lambda^{-n}+\left(  -k\right)  ^{-n}\lambda
\end{pmatrix}
\label{eqApp.Cla15}%
\end{equation}
with integer coefficients, $n=1,2,\dots$.

In the special case
\begin{equation}
\operatorname*{Prim}\left(  \lambda\right)  =\operatorname*{Prim}\left(
k\right)  , \label{eqApp.Cla16}%
\end{equation}
put
\begin{equation}
\lambda_{0}=\prod\left\{  p\mid p\in\operatorname*{Prim}\left(  \lambda
\right)  =\operatorname*{Prim}\left(  k\right)  \right\}  .
\label{eqApp.Cla17}%
\end{equation}
Then%
\begin{equation}
G=%
\begin{pmatrix}
\mathbb{Z}\left[  \frac{1}{\lambda_{0}}\right]  _{\mathstrut}\\
\mathbb{Z}\left[  \frac{1}{\lambda_{0}}\right]  ^{\mathstrut}%
\end{pmatrix}
\cong\mathbb{Z}\left[  \tfrac{1}{\lambda_{0}}\right]  ^{2}.
\label{eqApp.Cla18}%
\end{equation}
\end{proposition}

\begin{remark}
\label{RemApp.ClaNew.1}Note that the condition \textup{(\ref{ProApp.Cla1(2)})}
in Proposition \textup{\ref{ProApp.Cla1}} is equivalent to each of the conditions

\begin{enumerate}
\addtocounter{enumi}{2}

\item \label{ProApp.Cla1(3)}$\operatorname*{Prim}\left(  \lambda+k\right)
\subset\operatorname*{Prim}\left(  \lambda\right)  $,

\item \label{ProApp.Cla1(4)}$\operatorname*{Prim}\left(  \lambda+k\right)
\subset\operatorname*{Prim}\left(  k\right)  $.
\end{enumerate}

\noindent Clearly \textup{(\ref{ProApp.Cla1(2)})} $\Rightarrow$
\textup{(\ref{ProApp.Cla1(3)})} and \textup{(\ref{ProApp.Cla1(2)})}
$\Rightarrow$ \textup{(\ref{ProApp.Cla1(4)}).} For
\textup{(\ref{ProApp.Cla1(3)})} $\Rightarrow$ \textup{(\ref{ProApp.Cla1(2)}),}
use $k=\left(  \lambda+k\right)  -\lambda$, etc.
\end{remark}

\begin{proof}
Setting $n_{2}=-\lambda n_{1}$ in \textup{(\ref{eqApp.Cla9})} we obtain%
\begin{equation}
g=\frac{1}{\lambda+k}\left(  kn_{1}+\lambda n_{1}\right)  \left(  -k\right)
^{-n}%
\begin{pmatrix}
1\\
-\lambda
\end{pmatrix}
=n_{1}\left(  -k\right)  ^{-n}%
\begin{pmatrix}
1\\
-\lambda
\end{pmatrix}
\label{eqApp.Cla14bis}%
\end{equation}
and hence%
\begin{equation}
\mathbb{Z}\left[  \tfrac{1}{k}\right]
\begin{pmatrix}
1\\
-\lambda
\end{pmatrix}
\subset G.. \label{eqApp.Cla15bis}%
\end{equation}
Next setting $n_{2}=kn_{1}$ in \textup{(\ref{eqApp.Cla9})} we obtain
similarly
\begin{equation}
g=\frac{1}{\lambda+k}\left(  \lambda n_{1}+kn_{1}\right)  \lambda^{-n}%
\begin{pmatrix}
1\\
k
\end{pmatrix}
=n_{1}\lambda^{-n}%
\begin{pmatrix}
1\\
k
\end{pmatrix}
, \label{eqApp.Cla16bis}%
\end{equation}
so%
\begin{equation}
\mathbb{Z}\left[  \tfrac{1}{\lambda}\right]
\begin{pmatrix}
1\\
k
\end{pmatrix}
\subset G. \label{eqApp.Cla17bis}%
\end{equation}
We conclude that%
\begin{equation}
\mathbb{Z}\left[  \tfrac{1}{\lambda}\right]
\begin{pmatrix}
1\\
k
\end{pmatrix}
+\mathbb{Z}\left[  \tfrac{1}{k}\right]
\begin{pmatrix}
1\\
-\lambda
\end{pmatrix}
\subset G. \label{eqApp.Cla18bis}%
\end{equation}
Conversely, if $g\in G$ it follows from \textup{(\ref{eqApp.Cla9})} that $g$
has the form%
\begin{equation}
g=a%
\begin{pmatrix}
1\\
k
\end{pmatrix}
+b%
\begin{pmatrix}
1\\
-\lambda
\end{pmatrix}
, \label{eqApp.Cla19}%
\end{equation}
where the pair $a,b$ has the form%
\begin{equation}
a=\frac{1}{\lambda+k}\left(  \lambda n_{1}+n_{2}\right)  \lambda^{-n},\qquad
b=\frac{1}{\lambda+k}\left(  kn_{1}-n_{2}\right)  \left(  -k\right)  ^{-n}
\label{eqApp.Cla20}%
\end{equation}
for suitable $n_{1},n_{2}\in\mathbb{Z}$, $n\in\mathbb{Z}$. But writing%
\[
a=\frac{\lambda^{m}}{\lambda+k}\left(  \lambda n_{1}+n_{2}\right)
\lambda^{-\left(  n+m\right)  }%
\]
and choosing $m$ large enough, it follows that%
\[
a\in\frac{1}{\left(  \lambda+k,\lambda\right)  }\mathbb{Z}\left[  \tfrac
{1}{\lambda}\right]  ,
\]
and using the same reasoning on $b$, the inclusion
\begin{equation}
G\subset\frac{1}{\left(  \lambda+k,\lambda\right)  }\mathbb{Z}\left[
\tfrac{1}{\lambda}\right]
\begin{pmatrix}
1\\
k
\end{pmatrix}
+\frac{1}{\left(  \lambda+k,\lambda\right)  }\mathbb{Z}\left[  \tfrac{1}%
{k}\right]
\begin{pmatrix}
1\\
-\lambda
\end{pmatrix}
\label{eqApp.Cla21}%
\end{equation}
follows. Applying \textup{(\ref{eqRemNewBis.17})} and
\textup{(\ref{eqRemNewBis.34})} to $\alpha=\left(  1,1/\lambda\right)  $, it
follows that $g=a\left(
\begin{smallmatrix}
1\\
k
\end{smallmatrix}
\right)  +b\left(
\begin{smallmatrix}
1\\
-\lambda
\end{smallmatrix}
\right)  $ is positive if and only if $a>0$.

We have $\left(  \lambda+k,k\right)  =1$ if and only if $\operatorname*{Prim}%
\left(  \lambda+k\right)  \subset\operatorname*{Prim}\left(  k\right)  $, and
hence \textup{(\ref{ProApp.Cla1(2)})} $\Rightarrow$
\textup{(\ref{ProApp.Cla1(1)}).} But if \textup{(\ref{ProApp.Cla1(2)})} is not
fulfilled, then $\left(  \lambda+k,k\right)  >1$ or $\left(  \lambda
+k,\lambda\right)  >1$, and choosing $n_{1}=0$ \textup{(}or $n_{2}%
=0$\textup{)} in \textup{(\ref{eqApp.Cla9})} we see that $G$ contains elements
that are not in $\mathbb{Z}\left[  \tfrac{1}{\lambda}\right]  \left(
\begin{smallmatrix}
1\\
k
\end{smallmatrix}
\right)  +\mathbb{Z}\left[  \tfrac{1}{k}\right]  \left(
\begin{smallmatrix}
1\\
-\lambda
\end{smallmatrix}
\right)  $. Thus \textup{(\ref{ProApp.Cla1(1)})} $\Rightarrow$
\textup{(\ref{ProApp.Cla1(2)}).}

Next, define%
\begin{equation}
g\left(  n_{1},n_{2},n\right)  =\frac{1}{\lambda+k}\left[  \left(  \lambda
n_{1}+n_{2}\right)  \lambda^{-n}%
\begin{pmatrix}
1\\
k
\end{pmatrix}
+\left(  kn_{1}-n_{2}\right)  \left(  -k\right)  ^{-n}%
\begin{pmatrix}
1\\
-\lambda
\end{pmatrix}
\right]  \label{eqApp.Cla22}%
\end{equation}
for $n\in\mathbb{N}$, $n_{1},n_{2}\in\mathbb{Z}$. Then
\[
g\left(  n_{1},n_{2},n\right)  =n_{1}g\left(  1,0,n\right)  +n_{2}g\left(
0,1,n\right)
\]
and%
\[
g\left(  1,0,n+1\right)  =g\left(  0,1,n\right)  ,
\]
so $G$ is spanned over $\mathbb{Z}$ by $\mathbb{Z}^{2}$ and the elements%
\begin{align}
g\left(  0,1,n\right)   &  =\frac{1}{\lambda+k}\left(  \lambda^{-n}%
\begin{pmatrix}
1\\
k
\end{pmatrix}
-\left(  -k\right)  ^{-n}%
\begin{pmatrix}
1\\
-\lambda
\end{pmatrix}
\right) \label{eqApp.Cla23}\\
&  =\frac{1}{\lambda+k}%
\begin{pmatrix}
\lambda^{-n}-\left(  -k\right)  ^{-n}\\
k\lambda^{-n}+\left(  -k\right)  ^{-n}\lambda
\end{pmatrix}
.\nonumber
\end{align}
It remains to prove the last statement in the proposition. So assume
$\operatorname*{Prim}\left(  \lambda\right)  =\operatorname*{Prim}\left(
k\right)  $ and define $\lambda_{0}$ by \textup{(\ref{eqApp.Cla17})} as the
product of the primes in this set. It follows that the matrix elements in the
left column of%
\[
J=%
\begin{pmatrix}
\lambda-k & 1\\
k\lambda & 0
\end{pmatrix}
\]
all are divisible by $\lambda_{0}$ and hence all matrix elements of $J^{2n}$
are divisible by $\lambda_{0}^{n}$, i.e.,%
\[
\lambda_{0}^{n}\mathbb{Z}_{{}}^{2}\supset J^{2n}\mathbb{Z}_{{}}^{2},
\]
and hence, applying $\lambda_{0}^{-n}J_{{}}^{-2n}$ to both sides,%
\[
J^{-2n}\mathbb{Z}^{2}\supset\lambda_{0}^{-n}\mathbb{Z}_{{}}^{2}.
\]
It follows that%
\begin{equation}
G\supset\mathbb{Z}\left[  \tfrac{1}{\lambda_{0}}\right]  \mathbb{Z}^{2}=%
\begin{pmatrix}
\mathbb{Z}\left[  \frac{1}{\lambda_{0}}\right]  _{\mathstrut}\\
\mathbb{Z}\left[  \frac{1}{\lambda_{0}}\right]  ^{\mathstrut}%
\end{pmatrix}
. \label{eqApp.Cla24}%
\end{equation}
Conversely, $G$ is spanned over $\mathbb{Z}$ by $\mathbb{Z}^{2}$ and the
elements \textup{(\ref{eqApp.Cla15})} for $n=1,2,\dots$. But%
\begin{align}
\frac{1}{\lambda+k}%
\begin{pmatrix}
\lambda^{-n}-\left(  -k\right)  ^{-n}\\
k\lambda^{-n}+\left(  -k\right)  ^{-n}\lambda
\end{pmatrix}
&  =\frac{1}{\lambda^{n}k^{n}}\frac{1}{\lambda+k}%
\begin{pmatrix}
k^{n}-\left(  -\lambda\right)  ^{n}\\
k^{n+1}+\left(  -1\right)  ^{n}\lambda^{n+1}%
\end{pmatrix}
\label{eqApp.Cla25}\\
&  =\frac{\left(  -1\right)  ^{n+1}}{\lambda^{n}k^{n}}\frac{1}{\lambda+k}%
\begin{pmatrix}
\lambda^{n}-\left(  -k\right)  ^{n}\\
-\left(  \lambda^{n+1}-\left(  -k\right)  ^{n+1}\right)
\end{pmatrix}
.\nonumber
\end{align}
But since $\lambda=-k\pmod{\lambda+k}$, we have $\lambda^{n}=\left(
-k\right)  ^{n}\pmod{\lambda+k}$, so the vector $\dfrac{1}{\lambda+k}%
\begin{pmatrix}
\lambda^{n}-\left(  -k\right)  ^{n}\\
-\left(  \lambda^{n+1}-\left(  -k\right)  ^{n+1}\right)
\end{pmatrix}
$ has integral components. It follows that the elements
\textup{(\ref{eqApp.Cla15})} are contained in $\mathbb{Z}\left[  \frac
{1}{\lambda_{0}}\right]  ^{2}$. Thus%
\begin{equation}
G\subset\mathbb{Z}\left[  \tfrac{1}{\lambda_{0}}\right]  . \label{eqApp.Cla26}%
\end{equation}
Now, \textup{(\ref{eqApp.Cla24})} and \textup{(\ref{eqApp.Cla26})} finally
establish \textup{(\ref{eqApp.Cla18})}. \textup{(}The last argument was also
used in \cite[Remark after Proposition 5]{BJKR00}.\textup{)}
\end{proof}

In general, the sequence \textup{(\ref{eqApp.Cla7})} does not split, i.e.,
there does not exist a well defined homomorphism $\psi\colon\mathbb{Z}\left[
\frac{1}{\lambda}\right]  \rightarrow K_{0}\left(  \mathfrak{A}\right)  $ with
$\tau\circ\psi=\operatorname*{id}$. \textup{(}Well defined means for example
$\psi\left(  m\lambda^{-n}\right)  =\psi\left(  \left(  m\lambda\right)
\lambda^{-n-1}\right)  $.\textup{)} In general, the class of $K_{0}\left(
\mathfrak{A}\right)  $ in $\operatorname*{Ext}\left(  \mathbb{Z}\left[
\tfrac{1}{\lambda}\right]  ,\mathbb{Z}\left[  \tfrac{1}{k}\right]  \right)  $
depends on properties of the prime decompositions of $\lambda$ and $k$, and
seems to have to be treated on a case-by-case basis. There are, however, two
special cases that behave nicely. The first of these is the last case in
Proposition \textup{\ref{ProApp.Cla1},}%
\begin{equation}
\operatorname*{Prim}\left(  \lambda\right)  =\operatorname*{Prim}\left(
k\right)  . \label{eqApp.Cla27}%
\end{equation}
Then%
\begin{equation}
K_{0}\left(  G\right)  =%
\begin{pmatrix}
\mathbb{Z}\left[  \frac{1}{k}\right]  _{\mathstrut}\\
\mathbb{Z}\left[  \frac{1}{k}\right]  ^{\mathstrut}%
\end{pmatrix}
=%
\begin{pmatrix}
\mathbb{Z}\left[  \frac{1}{\lambda}\right]  _{\mathstrut}\\
\mathbb{Z}\left[  \frac{1}{\lambda}\right]  ^{\mathstrut}%
\end{pmatrix}
\label{eqApp.Cla28}%
\end{equation}
with trace functional $\left(  1,1/\lambda\right)  $. Since the dimension
group is a complete invariant, the following Proposition ensues.

\begin{proposition}
\label{ProApp.Cla2}Let $J_{1}$, $J_{2}$ be $2\times2$ matrices of the form
\textup{(\ref{eqApp.Cla1})--(\ref{eqApp.Cla2})} and let the subindices $1$,
$2$ refer to $J_{1}$, $J_{2}$ respectively.

\begin{enumerate}
\item \label{ProApp.Cla2(1)}If $J_{1}$, $J_{2}$ define isomorphic $C^{\ast}%
$-algebras, then $\operatorname*{Prim}\left(  \lambda_{1}\right)
=\operatorname*{Prim}\left(  \lambda_{2}\right)  $ and $\operatorname*{Prim}%
\left(  k_{1}\right)  =\operatorname*{Prim}\left(  k_{2}\right)  $.

\item \label{ProApp.Cla2(2)}If $\operatorname*{Prim}\left(  k_{1}\right)
=\operatorname*{Prim}\left(  k_{2}\right)  =\operatorname*{Prim}\left(
\lambda_{2}\right)  =\operatorname*{Prim}\left(  \lambda_{1}\right)  $, then
$J_{1}$, $J_{2}$ define isomorphic algebras.
\end{enumerate}

\noindent In the latter case, the dimension group is%
\begin{equation}
G\cong\mathbb{Z}\left[  \left\{  \tfrac{1}{k}\mid k\in\operatorname*{Prim}%
\left(  \lambda_{1}\right)  \right\}  \right]  ^{2} \label{eqApp.Cla29}%
\end{equation}
with positivity determined as follows: $g=\left(  g_{1},g_{2}\right)  \in G$
is positive if and only if $g_{1}+\lambda_{0}g_{2}>0$, where $\lambda
_{0}=\prod_{p\in\operatorname*{Prim}\left(  \lambda_{1}\right)  }p$.
\end{proposition}

\begin{proof}
As already remarked, \textup{(\ref{ProApp.Cla2(1)})} is a special case of
Theorem \textup{\ref{ThmCyc.7} and \cite[Theorem 10]{BJKR00}.} As for
\textup{(\ref{ProApp.Cla2(2)}),} Proposition \textup{\ref{ProApp.Cla1}} shows
that%
\begin{equation}
G_{1}=G_{2}=\mathbb{Z}\left[  \tfrac{1}{\lambda_{0}}\right]  ^{2}
\label{eqApp.Cla30}%
\end{equation}
as unordered groups, with positive cones determined by%
\begin{equation}
\left(  G_{i}\right)  _{+}=\left\{  g=\left(  g_{1},g_{2}\right)  \in
G_{i}\mid g_{1}+\frac{1}{\lambda_{i}}g_{2}>0\right\}  . \label{eqApp.Cla31}%
\end{equation}
But then the map%
\begin{equation}
\left(  g_{1},g_{2}\right)  \longmapsto\left(  g_{1},\frac{\lambda_{i}%
}{\lambda_{0}}g_{2}\right)  \label{eqApp.Cla32}%
\end{equation}
defines an isomorphism of ordered groups $G\rightarrow G_{i}$ for $i=1,2$.
Thus $\left(  G_{i},\left(  G_{i}\right)  _{+}\right)  $, $i=1,2$, are both
isomorphic to $\left(  G,G_{+}\right)  $, and Proposition
\textup{\ref{ProApp.Cla2}} follows.
\end{proof}

The second special case that behaves nicely is when condition
(\ref{ProApp.Cla1(2)}) in Proposition \ref{ProApp.Cla1} is fulfilled, i.e.,
$\operatorname*{Prim}\left(  \lambda+k\right)  \subset\operatorname*{Prim}%
\left(  \lambda\right)  \cap\operatorname*{Prim}\left(  k\right)  $. Let us
first mention a simple algorithm to construct \emph{all} pairs $\left(
\lambda,k\right)  $ of positive integers with $1\leq k\leq\lambda-1$
satisfying these properties: One first picks such a pair $\left(
\lambda^{\prime},k^{\prime}\right)  $ with $\gcd\left(  \lambda^{\prime
},k^{\prime}\right)  =1$, then lets $\mu$ be the product of all the prime
factors of $\lambda^{\prime}+k^{\prime}$, and then the pair%
\[
\lambda=n\mu\lambda^{\prime},\qquad k=n\mu k^{\prime},
\]
where $n$ is an arbitrary positive integer, will have the property
(\ref{ProApp.Cla1(2)}). One obtains all pairs $\left(  \lambda,k\right)  $
having the property (\ref{ProApp.Cla1(2)}) in this way, since given one such
pair one may divide by $\gcd\left(  \lambda,k\right)  $ to obtain $\left(
\lambda^{\prime},k^{\prime}\right)  $, and then get back to $\left(
\lambda,k\right)  $ by the process above.

As a simple example of the procedure above, take $\lambda^{\prime}=2$ and
$k^{\prime}=1$. This gives $\mu=3$, so all pairs
\[
\lambda=n\cdot6,\qquad k=n\cdot3
\]
are examples.

If (\ref{ProApp.Cla1(2)}) is fulfilled, there is an exact
sequence\label{LOSZ1lambda_1}%
\[
0\longrightarrow\mathbb{Z}\left[  \tfrac{1}{k}\right]  \hooklongrightarrow
G\overset{\tau}{\longrightarrow}\mathbb{Z}\left[  \tfrac{1}{\lambda}\right]
\longrightarrow0
\]
which splits. This is because $G$ has the form (\ref{ProApp.Cla1(1)}) in
Proposition \ref{ProApp.Cla1}, and one verifies directly that the map%
\[
\psi\colon\mathbb{Z}\left[  \tfrac{1}{\lambda}\right]  \longrightarrow G\colon
g\longrightarrow\frac{\lambda g}{\lambda+k}%
\begin{pmatrix}
1\\
k
\end{pmatrix}
\]
is a section (it is well defined since $\operatorname*{Prim}\left(
\lambda+k\right)  \subseteq\operatorname*{Prim}\left(  \lambda\right)  $).
Hence we get a different criterion from that in Proposition \ref{ProApp.Cla2}:

\begin{proposition}
\label{ProApp.ClaNew3}Let $J_{1}$, $J_{2}$ be $2\times2$ matrices of the form
\textup{(\ref{eqApp.Cla1})--(\ref{eqApp.Cla2})} and let the subindices $1$,
$2$ refer to $J_{1}$, $J_{2}$ respectively.

If $\operatorname*{Prim}\left(  k_{1}\right)  =\operatorname*{Prim}\left(
k_{2}\right)  $ and $\operatorname*{Prim}\left(  \lambda_{1}\right)
=\operatorname*{Prim}\left(  \lambda_{2}\right)  $, and%
\[
\operatorname*{Prim}\left(  \lambda_{i}+k_{i}\right)  \subseteq
\operatorname*{Prim}\left(  \lambda_{i}\right)  \cap\operatorname*{Prim}%
\left(  k_{i}\right)
\]
for $i=1,2$ \textup{(}see Remark \textup{\ref{RemApp.ClaNew.1}} after
Proposition \textup{\ref{ProApp.Cla1}),} then $J_{1}$, $J_{2}$ define
isomorphic algebras. Furthermore in this case the dimension group is%
\[
G\cong\mathbb{Z}\left[  \tfrac{1}{\lambda_{1}}\right]  \oplus\mathbb{Z}\left[
\tfrac{1}{k_{1}}\right]
\]
with positivity determined by positivity of the first coordinate.
\end{proposition}

\begin{proof}
This follows from the discussion before the Proposition.
\end{proof}

\begin{example}
\label{ExaApp.ClaNew.4}If%
\[
J_{1}=%
\begin{pmatrix}
6n_{1} & 1\\
3n_{1} & 0
\end{pmatrix}
\qquad J_{2}=%
\begin{pmatrix}
6n_{2} & 1\\
3n_{2} & 0
\end{pmatrix}
\]
where the positive integers $n_{1}$, $n_{2}$ contain the same prime factors,
then the corresponding AF-algebras are isomorphic. The case $n_{1}=2$,
$n_{2}=4$ is illustrated in Figure \textup{\ref{BratDiagsn2n4},} below.
\end{example}

\begin{figure}[ptb]
\begin{picture}(360,528)
\put(0,0){\includegraphics
[bb=267 120 453 720,clip,width=167.4bp,height=540bp]{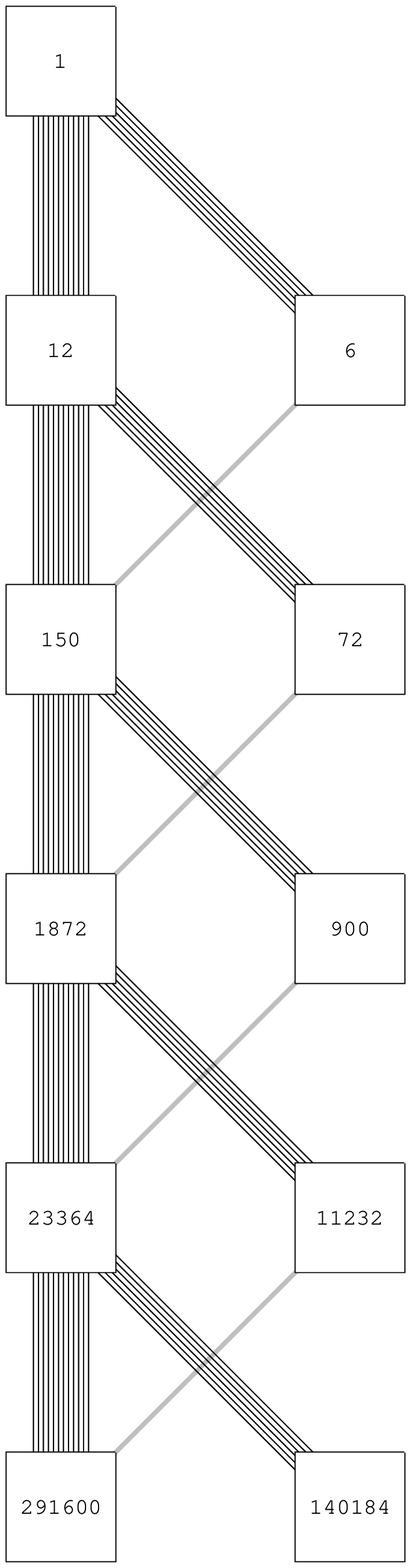}}
\put(192.6,0){\includegraphics
[bb=267 120 453 720,clip,width=167.4bp,height=540bp]{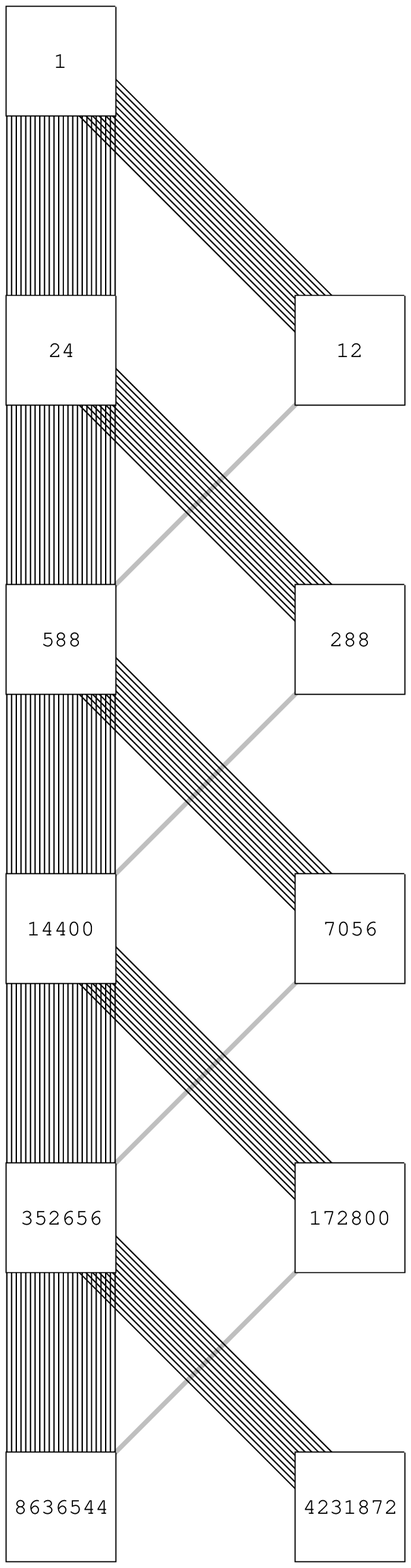}}
\end{picture}\caption{$L=\{\frogbrace{12}{1},2,2,2,2,2,2\}$, first column
$=(12\;6)^{\mathrm{t}}$ (left); $L=\{\frogbrace{24}{1},\frogbrace{12}{2}\}$,
first column $=(24\;12)^{\mathrm{t}}$ (right). See Example
\ref{ExaApp.ClaNew.4}. These diagrams represent isomorphic algebras.}%
\label{BratDiagsn2n4}%
\end{figure}

\begin{remark}
\label{RemApp.ClaNew.5}If%
\[
J=%
\begin{pmatrix}
\lambda-k & 1\\
\lambda k & 0
\end{pmatrix}
\]
then the right Perron--Frobenius eigenvector $v$ in
\textup{(\ref{eqSubNew.pound})} is $v=\left(
\begin{smallmatrix}
1\\
k
\end{smallmatrix}
\right)  $ and hence%
\[%
\ip{\alpha}{v}%
=1+\frac{k}{\lambda}=\frac{\lambda+k}{\lambda}.
\]
But if $\operatorname*{Prim}\left(  \lambda+k\right)  \subseteq
\operatorname*{Prim}\lambda$, this number has a multiplicative inverse in
$\mathbb{Z}\left[  \frac{1}{\lambda}\right]  $. Hence the split property used
in the proof of Proposition \textup{\ref{ProApp.ClaNew3}} can also be deduced
from Corollary \textup{\ref{corollary10.25}.}
\end{remark}

\begin{table}[ptb]
\caption{$\operatorname*{Prim}$ invariants for various $\mathfrak{A}_{L}$
algebras\firkelbow with $\operatorname*{rank}\left(  K_{0}\left(  \mathfrak
{A}_{L}\right)  \right)  =2$.}%
\label{TableApp.Cla}%
\begin{tabular}
[c]{cccc}%
Matrix (\ref{eqApp.Cla2}): $\left(
\begin{smallmatrix}
\lambda-k & 1\\
k\cdot\lambda & 0
\end{smallmatrix}
\right)  $ & Block form (\ref{eqApp.Cla4}): $\left(
\begin{smallmatrix}
-k & k\\
0 & \lambda
\end{smallmatrix}
\right)  $ & $\operatorname*{Prim}\left(  k\right)  $ & $\operatorname*{Prim}%
\left(  \lambda\right)  $\\
{\underline{$\lambda=8$}\hfill} &  &  & \\
\vphantom{$\begin{pmatrix}0\\0\\0\end{pmatrix}$}$%
\begin{pmatrix}
7 & 1\\
8 & 0
\end{pmatrix}
$ & $%
\begin{pmatrix}
-1 & 1\\
0 & 8
\end{pmatrix}
$ & $\varnothing$ & $\left\{  2\right\}  $\\
\vphantom{$\begin{pmatrix}0\\0\\0\end{pmatrix}$}$%
\begin{pmatrix}
6 & 1\\
16 & 0
\end{pmatrix}
$ & $%
\begin{pmatrix}
-2 & 2\\
0 & 8
\end{pmatrix}
$ & $\left\{  2\right\}  $ & $\left\{  2\right\}  $\\
\vphantom{$\begin{pmatrix}0\\0\\0\end{pmatrix}$}$%
\begin{pmatrix}
5 & 1\\
24 & 0
\end{pmatrix}
$ & $%
\begin{pmatrix}
-3 & 3\\
0 & 8
\end{pmatrix}
$ & $\left\{  3\right\}  $ & $\left\{  2\right\}  $\\
\vphantom{$\begin{pmatrix}0\\0\\0\end{pmatrix}$}$%
\begin{pmatrix}
4 & 1\\
32 & 0
\end{pmatrix}
$ & $%
\begin{pmatrix}
-4 & 4\\
0 & 8
\end{pmatrix}
$ & $\left\{  2\right\}  $ & $\left\{  2\right\}  $\\
\vphantom{$\begin{pmatrix}0\\0\\0\end{pmatrix}$}$%
\begin{pmatrix}
3 & 1\\
40 & 0
\end{pmatrix}
$ & $%
\begin{pmatrix}
-5 & 5\\
0 & 8
\end{pmatrix}
$ & $\left\{  5\right\}  $ & $\left\{  2\right\}  $\\
\vphantom{$\begin{pmatrix}0\\0\\0\end{pmatrix}$}$%
\begin{pmatrix}
2 & 1\\
48 & 0
\end{pmatrix}
$ & $%
\begin{pmatrix}
-6 & 6\\
0 & 8
\end{pmatrix}
$ & $\left\{  2,3\right\}  $ & $\left\{  2\right\}  $\\
\vphantom{$\begin{pmatrix}0\\0\\0\end{pmatrix}$}$%
\begin{pmatrix}
1 & 1\\
56 & 0
\end{pmatrix}
$ & $%
\begin{pmatrix}
-7 & 7\\
0 & 8
\end{pmatrix}
$ & $\left\{  7\right\}  $ & $\left\{  2\right\}  $\\
{\underline{$\lambda=13$}\hfill} &  &  & \\
\vphantom{$\begin{pmatrix}0\\0\\0\end{pmatrix}$}$%
\begin{pmatrix}
7 & 1\\
6\cdot13 & 0
\end{pmatrix}
$ & $%
\begin{pmatrix}
-6 & 6\\
0 & 13
\end{pmatrix}
$ & $\left\{  2,3\right\}  $ & $\left\{  13\right\}  $\\
\vphantom{$\begin{pmatrix}0\\0\\0\end{pmatrix}$}$%
\begin{pmatrix}
1 & 1\\
12\cdot13 & 0
\end{pmatrix}
$ & $%
\begin{pmatrix}
-12 & 12\\
0 & 13
\end{pmatrix}
$ & $\left\{  2,3\right\}  $ & $\left\{  13\right\}  $%
\end{tabular}
\noindent\end{table}

Let us focus on an example of the use of Propositions
\textup{\ref{ProApp.Cla1}} and \textup{\ref{ProApp.Cla2}:} Consider the list
of matrices in Table \textup{\ref{TableApp.Cla},} below. It follows from
Proposition \textup{\ref{ProApp.Cla2}(\ref{ProApp.Cla2(1)})} that the only
candidates for nontrivial pairs defining isomorphic AF-algebras from this list
are
\begin{equation}%
\begin{pmatrix}
6 & 1\\
16 & 0
\end{pmatrix}
,\qquad%
\begin{pmatrix}
4 & 1\\
32 & 0
\end{pmatrix}
\label{eqApp.Cla34}%
\end{equation}
and
\begin{equation}
J=%
\begin{pmatrix}
7 & 1\\
78 & 0
\end{pmatrix}
,\qquad J^{\prime}=%
\begin{pmatrix}
1 & 1\\
156 & 0
\end{pmatrix}
. \label{eqApp.Cla35}%
\end{equation}
The first of these pairs actually defines isomorphic algebras by Proposition
\textup{\ref{ProApp.Cla2}(\ref{ProApp.Cla2(2)}).} \textup{(}This was already
proved in \cite[Proposition 5]{BJKR00}.\textup{)} For the latter of these pairs
the special criteria of Proposition \textup{\ref{ProApp.Cla1}} cannot be
employed. But as $78=2\cdot3\cdot13$ and $156=2\cdot2\cdot3\cdot13$ we have
$\gcd\left(  7,78\right)  =1$ and $\gcd\left(  1,156\right)  =1$, and it
follows from Proposition \ref{prop10.22} that $\deg J^{\prime}=156$-$\deg
J^{\prime}=1=\deg J=78$-$\deg J$. But $m=\operatorname{lcm}\left(
78,156\right)  =156$ and hence we may apply Theorem \ref{theorem10.10} to the
pair $J$, $J^{\prime}$. Using the formula in Remark \ref{RemApp.ClaNew.5} we
see that
\[%
\ip{\alpha}{v}%
=\frac{19}{13},\qquad%
\ip{\alpha^{\prime}}{v^{\prime}}%
=\frac{25}{13},
\]
and hence%
\[
\frac{%
\ip{\alpha}{v}%
}{%
\ip{\alpha^{\prime}}{v^{\prime}}%
}=\frac{19}{25}\notin\mathbb{Z}\left[  \tfrac{1}{13}\right]  .
\]
It follows from Theorem \ref{theorem10.10} that $\mathfrak{A}_{J}$ and
$\mathfrak{A}_{J^{\prime}}$ are non-isomorphic.

Finally, note that the set of $2\times2$ matrices of the form
(\ref{eqCycNew.1}), or (\ref{eqApp.Cla2}), with $\lambda=m_{2}$, i.e., the
matrices
\[%
\begin{pmatrix}
\lambda-1 & 1\\
\lambda & 0
\end{pmatrix}
\]
for $\lambda=2,3,4,\dots$ all give rise to non-isomorphic algebras. This is
proved in Section \ref{ClmN.2}.

\begin{case}
\label{CasApp.Cla2}The Perron--Frobenius eigenvalue $\lambda$ is irrational,
and hence in a quadratic extension of $\mathbb{Z}$, since $\lambda$ satisfies
a monic quadratic equation.
\end{case}

In this case, the exact sequence \textup{(\ref{eqRemNewBis.31})} is%
\[
0\longrightarrow0\longrightarrow K_{0}\left(  \mathfrak{A}\right)
\longrightarrow\mathbb{Z}\left[  \tfrac{1}{\lambda}\right]  \longrightarrow0,
\]
so $K_{0}\left(  \mathfrak{A}\right)  $ equals $\mathbb{Z}\left[  \frac
{1}{\lambda}\right]  $, in the sense of ordered groups. But since $\left(
K_{0}\left(  \mathfrak{A}\right)  ,K_{0}\left(  \mathfrak{A}\right)
_{+}\right)  $ is a complete invariant, and the \emph{ordered} group
$\mathbb{Z}\left[  \lambda\right]  $ determines $1/\lambda$ uniquely when
$1/\lambda$ is in a quadratic extension of $\mathbb{Q}$, it follows that the
irrational number%
\[
\frac{1}{\lambda}=\frac{\sqrt{m_{1}^{2}+4m_{2}^{{}}}-m_{1}^{{}}}{2m_{2}^{{}}}%
\]
is a complete invariant. But since the equation
\[
1-m_{1}x-m_{2}x^{2}=0
\]
for $1/\lambda$ is irreducible, and $1/\lambda$ is characterized as the
positive solution of this equation it follows that $\lambda$ determines
$m_{1}$, $m_{2}$ in this case. Conclusion:

\begin{proposition}
\label{ProApp.Cla3}If%
\begin{equation}
J=%
\begin{pmatrix}
m_{1} & 1\\
m_{2} & 0
\end{pmatrix}
,\qquad J^{\prime}=%
\begin{pmatrix}
m_{1}^{\prime} & 1\\
m_{2}^{\prime} & 0
\end{pmatrix}
\label{eqApp.Cla36}%
\end{equation}
are matrices with $m_{1}^{{}},m_{2}^{{}},m_{1}^{\prime},m_{2}^{\prime}%
\in\mathbb{N}$ and at least one of the numbers $\sqrt{m_{1}^{2}+4m_{2}^{{}}}$
or $\sqrt{m_{1}^{\prime\,2}+4m_{2}^{\prime}}$ is irrational, then the
AF-algebra $\mathfrak{A}_{J}$ $is$ isomorphic to $\mathfrak{A}_{J^{\prime}}$
if and only if
\begin{equation}
m_{1}^{{}}=m_{1}^{\prime}\text{\qquad and\qquad}m_{2}^{{}}=m_{2}^{\prime}.
\label{eqApp.Cla37}%
\end{equation}
\end{proposition}

Finally, recall that if $m_{1}=1$, $m_{2}=m\in\mathbb{N}$, then $\mathfrak
{A}_{J}$ is the Pimsner--Voiculescu algebra associated with the continued
fraction%
\begin{equation}
\lambda=\frac{m_{1}^{{}}+\sqrt{m_{1}^{2}+4m_{2}^{{}}}}{2}=m_{1}+\frac{m_{2}%
}{\lambda}=m_{1}+\cfrac{m_{2}}{m_{1}+\cfrac{m_{2}}{m_{1}+\cfrac{m_{2}}%
{m_{1}+\cdots}}}; \label{eqApp.Cla38}%
\end{equation}
see \cite{PiVo80}. See Figure \textup{\ref{BratDiag12}} for the special case
$m_{1}=m_{2}=1$. \setcounter{figurelink}{\value{figure}} \setcounter
{tablelink}{\value{table}}

\chapter{\label{App.Inv}Linear algebra of $J$}

\setcounter{table}{\value{tablelink}} We have introduced several parameters
related to $C^{\ast}$-isomorphism invariants on the AF-algebras $\mathfrak
{A}_{L}$, $L=\left(  L_{1},\dots,L_{d}\right)  $ such that $\gcd\left(
L_{i}\right)  =1$. Those parameters compute out as shown in Table
\ref{TableInvariants} for the examples in Figures \ref{BratDiag11}%
--\ref{BratDiag235} and \ref{BratDiagsd6}. We have included the value of $d$
from $\mathcal{O}_{d}$, and the Perron--Frobenius eigenvalue $\lambda$, in the
table (we have shown in Example \ref{ExaRemNewBis.3} that $d$ is not an
invariant in general). The actual invariants are $\mathbb{Z}\left[  \frac
{1}{\lambda}\right]  $, $D$, $N$, $\operatorname*{Prim}\left(  m_{N}%
/R_{D}\right)  $, $\operatorname*{Prim}\left(  R_{D}\right)  $, while
$\tau\left(  v\right)  $ is a restricted invariant (see Theorems
\ref{CorCyc.9} and \ref{theorem10.10}).

\renewcommand{\arraystretch}{1.5}

\begin{table}[ptb]
\caption{Some parameters related to invariants\firkelbow for various
$\mathfrak{A}_{L}$ algebras.}%
\label{TableInvariants}%
\begin{tabular}
[c]{rlclccccl}%
\begin{minipage}[t]{16pt}Fig.\newline\makebox[16pt][r]{\#} \end{minipage} &
Equation & $d$ & $\lambda$ & $D$ & $N$ & $m_{N}$ & $R_{D}$ & $\tau\left(
v\right)  $\\
\ref{BratDiag11} & $2x=1$ & $2$ & $\lambda=2$ & $1$ & $1$ & $2$ & $2$ & $1$\\
\ref{BratDiag12} & $x+x^{2}=1$ & $2$ & $\lambda=\frac{1+\sqrt{5}}{2}$ & $2$ &
$2$ & $1$ & $1$ & $\frac{5-\sqrt{5}}{2}=\frac{\sqrt{5}}{\lambda}$\\
\ref{BratDiag4458} & $2x^{4}+x^{5}+x^{8}=1$ & $4$ & \begin{minipage}%
[t]{60pt}$\lambda=a_{3}^{-1}$,\newline\hspace*{0.5em}\rlap{$a_{3}%
=0.7549{}^{\ast}$}\end{minipage} & $3$ & $8$ & $1$ & $1$ & \begin{minipage}%
[t]{60pt}$4+a_{3}^{5}+4a_{3}^{8}$\newline\hspace*{0.5em}$\approx
4.6669 $ \end{minipage}\\
\ref{BratDiag23} & $x^{2}+x^{3}=1$ & $2$ & $\lambda=a_{3}^{-1}{}^{\ast}$ & $3$%
& $3$ & $1$ & $1$ & \begin{minipage}[t]{60pt}$2+a_{3}^{3}$\newline
\hspace*{0.5em}$\approx2.4302$ \end{minipage}\\
\ref{BratDiag235} & $x^{2}+x^{3}+x^{5}=1$ & $3$ & \begin{minipage}%
[t]{60pt}$\lambda=a^{-1}$,\newline\hspace*{0.5em}$a=0.6997$ \end{minipage} &
$5$ & $5$ & $1$ & $1$ & \begin{minipage}[t]{60pt}$2+a^{3}+3a^{5}%
$\newline\hspace*{0.5em}$\approx2.8459$ \end{minipage}\\
& $x+x^{5}=1$ & $2$ & $\lambda=a_{3}^{-1}{}^{\ast}$ & $3$ & $5$ & $1$ & $1$ &
\begin{minipage}[t]{60pt}$1+4a_{3}^{5}$\newline\hspace*{0.5em}$\approx
1.9805$ \end{minipage}\\
& $x+4x^{3}=1$ & $5$ & $\lambda=2$ & $1$ & $3$ & $4$ & $2$ & $2$\\
& $3x^{2}+2x^{3}=1$ & $5$ & $\lambda=2$ & $1$ & $3$ & $2$ & $2$ & $\frac{9}%
{4}$\\
\ref{BratDiagsd6} & $x+3x^{3}+2x^{4}=1$ & $6$ & $\lambda=2$ & $1$ & $4$ & $2$%
& $2$\phantom{\rlap{\raisebox{-6pt}[0pt][0pt]{$\Bigr\rbrace$ \parbox
{75pt}{(non-iso\-mor\-phic $\mathfrak{A}_{L}$'s)\hfill}}}} & $\frac{17}{8}$\\
\ref{BratDiagsd6} & $3x^{2}+x^{3}+2x^{4}=1$ & $6$ & $\lambda=2$ & $1$ & $4$ &
$2$ & $2$ & $\frac{19}{8}$\\
& $x+2x^{3}+4x^{4}=1$ & $7$ & $\lambda=2$ & $1$ & $4$ & $4$ & $2$%
\phantom{\rlap{\raisebox
{-6pt}[0pt][0pt]{$\Bigr\rbrace$ \parbox{75pt}%
{(still non-iso\-mor\-phic; see \S\ref{Rem})\hfill}}}} & $\frac{9}{4}$\\
& $3x^{2}+4x^{4}=1$ & $7$ & $\lambda=2$ & $1$ & $4$ & $4$ & $2$ & $\frac{5}%
{2}$%
\end{tabular}
\noindent\parbox{\textwidth}{\small${}^*\,$Note:
$x^{2}+x^{3}-1$ is a factor of both
$2x^{4}+x^{5}+x^{8}-1$ and
$x+x^{5}-1$; see lines 3, 4, 6.}\end{table}

\renewcommand{\arraystretch}{1}

\begin{remark}
\label{RemSubNew.11}It is to be stressed that the parameters are computed for
examples $\left(  L_{1},\dots,L_{d}\right)  $, subject to the restriction that
the greatest common divisor is one, i.e., $\gcd\left(  L_{i}\right)  =1$. It
is immediate from Chapter \ref{KMS} that the pair of AF-algebras $\mathfrak
{A}_{L}$ and $\mathfrak
{A}_{sL}$ computed from the two, $\left(  L_{1},\dots,L_{d}\right)  $ and
$\left(  sL_{1},\dots,sL_{d}\right)  $, are isomorphic. But, in a sense, the
two parameters $N$ and $\varphi_{\tau}=\lambda^{-1}p_{L}^{\prime}\left(
\frac{1}{\lambda}\right)  $ scale by $s$. Since $m_{N}=\left|  \det J\right|
$, $m_{N}$ does not. To see this, note that the scaling $L\mapsto sL$,
$s\in\mathbb{N}$, corresponds to a simple inflation of $J$, as illustrated by
\textup{(}$s=2$\textup{):}%
\[
\begin{picture}(194,76)
\put(88,0){\includegraphics
[bb=0 11 98 87,height=76bp,width=98bp]{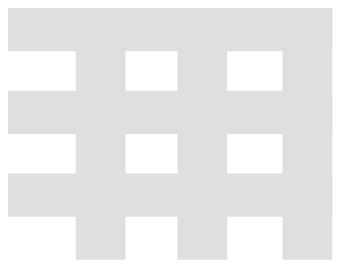}}
\put(0,0){\makebox(194,76)[l]{$\displaystyle
\begin{pmatrix}
m_{1} & 1 & 0\\
m_{2} & 0 & 1\\
m_{3} & 0 & 0
\end{pmatrix}
\longmapsto\left(
\begin{matrix}
0 & 1 & 0 & 0 & 0 & 0\\
m_{1} & 0 & 1 & 0 & 0 & 0\\
0 & 0 & 0 & 1 & 0 & 0\\
m_{2} & 0 & 0 & 0 & 1 & 0\\
0 & 0 & 0 & 0 & 0 & 1\\
m_{3} & 0 & 0 & 0 & 0 & 0
\end{matrix}
\mkern8mu \right) ,$}}
\end{picture}
\]
and these matrices define isomorphic stable AF-algebras. \textup{(}For the
$6\times6$ matrix the maps in the inductive sequence are not
injective.\textup{)} In the last example, to get back from the inflated matrix
to $J$, simply delete the rows numbered $1$, $3$, and $5$ \textup{(}%
shaded\textup{),} and the columns numbered $2$, $4$, and $6$.

Hence, to compute parameters for a general divisible $\left(  L_{i}\right)
$-system, first pass to $\left(  L_{i}^{\prime}\right)  $ where $L_{i}%
^{\prime}=\frac{1}{\gcd\left(  L_{j}\right)  }L_{i}$, and then use the
prescribed formulas \textup{(}for the parameters\/\textup{)} on the $\left(
L_{i}^{\prime}\right)  $-system.
\end{remark}

We next consider the following observation regarding the parameter
$\tau\left(  v\right)  $.

\begin{proposition}
[Scaling Property for the Parameter]\label{ProSubNew.12}Let $\left(
L_{1},\dots,L_{d}\right)  $ be given, and let $s\in\mathbb{Q}_{+}$. Let
$\tau\left(  v\right)  $ and $\tau\left(  v_{s}\right)  $ be the respective
numbers for $L$ and $sL$, as follows: let $v\in\mathbb{R}^{N}$ satisfy
$Jv=\lambda v$, $v_{1}=1$, where $\lambda$ is the Perron--Frobenius eigenvalue
for $\alpha J=\lambda\alpha$. Then%
\begin{equation}
\tau\left(  v\right)  =\lambda^{-1}p_{L}^{\prime}\left(  \frac{1}{\lambda
}\right)  \label{eqSub.star}%
\end{equation}
where\label{LOSpLx_4}%
\begin{equation}
p_{L}\left(  x\right)  =\sum_{i=1}^{d}x^{L_{i}}-1, \label{eqSub.starstar}%
\end{equation}
and the corresponding number for $sL$ is%
\begin{equation}
\tau\left(  v_{s}\right)  =s\tau\left(  v\right)  . \label{eqSub.starstarstar}%
\end{equation}
Suppose $\deg\left(  L\right)  >1$. Let the other roots of $p_{L}\left(
x\right)  $ be $\left\{  a_{i}\right\}  _{i=1}^{L_{d}-1}$. Then, by the
assumptions, $\left|  a_{i}^{-1}\right|  <\lambda$, and%
\[
p_{L}^{\prime}\left(  e^{-\beta_{L}}\right)  =\prod_{i=1}^{L_{d}-1}\left(
e^{-\beta_{L}}-a_{i}\right)  .
\]
\end{proposition}

\begin{proof}
Writing $J$ in the form\label{LOSJ_6}%
\begin{equation}
J=%
\begin{pmatrix}
\vphantom{\vdots}m_{1} & 1 & 0 & \cdots & 0 & 0\\
\vphantom{\vdots}m_{2} & 0 & 1 & \cdots & 0 & 0\\
\vdots &  & \ddots & \ddots & \vdots & \vdots\\
\vphantom{\vdots}m_{N-2} & 0 & 0 & \ddots & 1 & 0\\
\vphantom{\vdots}m_{N-1} & 0 & 0 &  & 0 & 1\\
\vphantom{\vdots}m_{N} & 0 & 0 & \cdots & 0 & 0
\end{pmatrix}
, \label{eqSub.pound.st}%
\end{equation}
note that then (\ref{eqSub.starstar}) becomes\label{LOSpLx_5}%
\[
p_{L}\left(  x\right)  =\sum_{j=1}^{N}m_{j}x^{j}-1,
\]
and the eigenvector $v$ may be computed by directly solving $Jv=\lambda
v$:\label{LOSPerronFrobeniusrighteigenvector_4}%
\begin{equation}
v=%
\begin{pmatrix}
\vphantom{\vdots}1\\
\vphantom{\vdots}\lambda-m_{1}\\
\vphantom{\vdots}\lambda^{2}-m_{1}\lambda-m_{2}\\
\vphantom{\vdots}\lambda^{3}-m_{1}\lambda^{2}-m_{2}\lambda-m_{3}\\
\vdots\\
\vphantom{\vdots}\lambda^{N-1}-m_{1}\lambda^{N-2}-\dots-m_{N-2}\lambda-m_{N-1}%
\end{pmatrix}
\in\mathbb{R}^{N}. \label{eqSubNew.pound}%
\end{equation}
Since $\tau\left(  x\right)  =%
\ip{\alpha}{x}%
$, with $\alpha=\left(  1,a,a^{2},\dots,a^{N-1}\right)  $, and $a=\frac
{1}{\lambda}$, we get
\begin{align*}
\tau\left(  v\right)   &  =%
\ip{\alpha}{v}%
\\
&  =N-\left(  N-1\right)  m_{1}a-\left(  N-2\right)  m_{2}a^{2}-\dots
-2m_{N-2}a^{N-2}-m_{N-1}a^{N-1}\\
&  =Nm_{N}a^{N}+m_{1}a+2m_{2}a^{2}+\dots+\left(  N-1\right)  m_{N-1}a^{N-1}\\
&  =ap_{L}^{\prime}\left(  a\right)  ,
\end{align*}
where we used the fact that
\[
m_{1}a+m_{2}a^{2}+\dots+m_{N}a^{N}=1.
\]
We claimed in (\ref{eqSub.star})--(\ref{eqSub.starstar}) that $\tau\left(
v\right)  =ap_{L}^{\prime}\left(  a\right)  $, and this now follows. The
scaling property (\ref{eqSub.starstarstar}) is immediate from this.
\end{proof}

\begin{corollary}
\label{CorSubNew.3}Let $J$, $J^{\prime}$ be two matrices specified as in
\textup{(\ref{eqSub.pound.st}),} and suppose they have the same value for the
rank $N$ and the same Perron--Frobenius eigenvalue $\lambda$. Let $v$,
resp.\ $v^{\prime}$, be the right Perron--Frobenius eigenvectors
\textup{(}with $v_{1}^{{}}=v_{1}^{\prime}=1$\textup{).} Then%
\[
v=v^{\prime}\Longleftrightarrow J=J^{\prime}.
\]
\end{corollary}

\begin{proof}
If $v=v^{\prime}$, then recursion in (\ref{eqSubNew.pound}) yields $m_{i}^{{}%
}=m_{i}^{\prime}$ for $i=1,2,\dots,N-1$. Since $N$ and $\lambda$ take the same
values on $J$ and $J^{\prime}$, the \emph{identity} (in each case),%
\[
\lambda^{N-1}-m_{1}\lambda^{N-2}-\dots-m_{N-2}\lambda-m_{N-1}=\frac{m_{N}%
}{\lambda},
\]
shows that then also $m_{N}^{{}}=m_{N}^{\prime}$, and therefore, by
(\ref{eqSub.pound.st}), $J=J^{\prime}$. The converse is clear.
\end{proof}

\setcounter{tablelink}{\value{table}}

\chapter{\label{App.Lat}Lattice points}

\setcounter{figure}{\value{figurelink}} Let $\left(  L_{i}\right)  _{i=1}^{d}$
be a standard system, and let $\mathfrak{A}_{L}$ be the corresponding
AF-algebra. We saw that the trace is unique and determined by the value of the
$L_{i}$'s. It is clear that when $\left(  L_{i}\right)  $ is given ($L_{i}>0$
say), there is a unique $\beta$ such that $\sum_{i=1}^{d}e^{-\beta L_{i}}=1$.
This means that $x_{\beta}:=e^{-\beta}$ is a root of\label{LOSpLx_6}%
\[
p_{L}\left(  x\right)  =x^{L_{1}}+\dots+x^{L_{d}}-1.
\]
But with the restrictions $L_{i}\in\mathbb{N}$, $L_{1}\leq L_{2}<\dots\leq
L_{d}$, it follows from Example \ref{ExaRemNewBis.3}, and, later, Chapter
\ref{APP.EXA}, that the $L_{i}$'s are \emph{not} determined by $\beta$. We
have already seen examples illustrating that, up to the obvious permutations,
there is for fixed $0<x<1$ and $d$, a multiplicity of lattice points on the
variety $\left(  L_{i}\right)  \subset\mathbb{R}^{d}$, $L_{i}>0$ with
Perron--Frobenius eigenvalue $1/x$. The pair of lattice points $\left(
2,3\right)  $, $\left(  1,5\right)  $ in Figure \ref{Lfromx} are on the same
curve. For $d=2$, we know of no other pair of distinct lattice points over the
$45^{\circ}$ line lying on the same curve.

\begin{example}
\label{ExaInv.6}Consider the AF-algebra of $x^{2}+x^{3}=1$ in Figure
\ref{BratDiag23}. There $e^{-\beta}=a\approx0.7549$ is the positive root and
$p_{\beta}\left(  x\right)  :=1-x^{2}-x^{3}$ is the minimal polynomial.
$K_{0}$ of this example is therefore given by \textup{(\ref{eqRemNew.pound}),}
and $\ker\left(  \tau\right)  =0$. But there is also an example for $N=5$,
$J_{5}=\left(
\begin{smallmatrix}
1 & 1 & 0 & 0 & 0\\
0 & 0 & 1 & 0 & 0\\
0 & 0 & 0 & 1 & 0\\
0 & 0 & 0 & 0 & 1\\
1 & 0 & 0 & 0 & 0
\end{smallmatrix}
\right)  $ with the same $\beta$ \textup{(}and therefore root $a=e^{-\beta}%
$\textup{).} Now for this related example, there are infinitesimal elements,
i.e., $\ker\left(  \tau_{5}\right)  \neq0$ \textup{(}in fact $\ker\left(
\tau_{5}\right)  \cong\mathbb{Z}^{2}$\textup{),} and hence the corresponding
two AF-algebras are non-isomorphic. Specifically, $\ker\left(  \tau
_{5}\right)  $ may be computed from \textup{(\ref{eqRem.star})} where the
restriction matrix $J_{0}$ is $\left(
\begin{smallmatrix}
1 & 1\\
-1 & 0
\end{smallmatrix}
\right)  $. Since $\det J_{0}=1$, it follows that $\ker\left(  \tau
_{5}\right)  \cong\mathbb{Z}^{2}$, as claimed.

The examples $x^{5}+x-1$ and $x^{3}+x^{2}-1$ show that we \emph{must} at least
add $N=N^{\prime}$ as a condition for isomorphism, because these two have the
same $\beta$. \textup{(}In fact $x^{5}+x-1=\left(  x^{3}+x^{2}-1\right)
\left(  x^{2}-x+1\right)  $.\textup{)} The triangular form $J_{5}=\left(
\begin{tabular}
[c]{c|c}%
$J_{0_{\mathstrut}}$ & $Q$\\\hline
$0$ & $J_{R}$%
\end{tabular}
\right)  $ from Theorem \textup{\ref{ThmCyc.7}} which corresponds to this
factorization is%
\[
J\colon\left(
\begin{tabular}
[c]{cc|ccc}%
\llap{\raisebox{-7pt}[0pt][0pt]{$J_{0}\rightarrow$\kern1.5em}}$-1$ & $1$ & $1$%
& $0$ & $0$\rlap{\raisebox{-7pt}[0pt][0pt]{\kern1.5em$\leftarrow Q$}}\\
$1$ & $0$ & $-1$ & $0$ & $0$\\\hline
$0$ & $0$ & $0$ & $1$ & $0$\\
$0$ & $0$ & $1$ & $0$ & $1$\rlap{\raisebox{-7pt}[0pt][0pt]{\kern
1.5em$\leftarrow J_{R}$}}\\
$0$ & $0$ & $1$ & $0$ & $0$%
\end{tabular}
\right)
\]
and the three $\operatorname*{Prim}$-invariants from Theorem
\textup{\ref{CorCyc.9}} are all $\varnothing$.
\end{example}

However, in Chapter \ref{APP.EXA}, for $d\geq5$, we will give other examples
of multiple lattice points. (We also gave such examples in Example
\ref{ExaRemNewBis.3}.) In those examples, for each of the cases $d=5$ and
$d=6$, there are three such multiple points on the $e^{-\beta}=\frac{1}{2}$
variety. There are others for different values of $\beta$, but none with
$d=2$. The value of $d$ is different for the two in the above pair; and if it
is further assumed that $l=L_{1}+L_{2}$ be the same, convexity and symmetry
show that there cannot be double points.

\begin{figure}[ptb]
\setlength{\unitlength}{5.45pt} \begin{picture}(59,57)(-5,-3)
\put(-1.9,-2.3){\includegraphics
[bb=1 0 287 288,height=288bp,width=286bp]{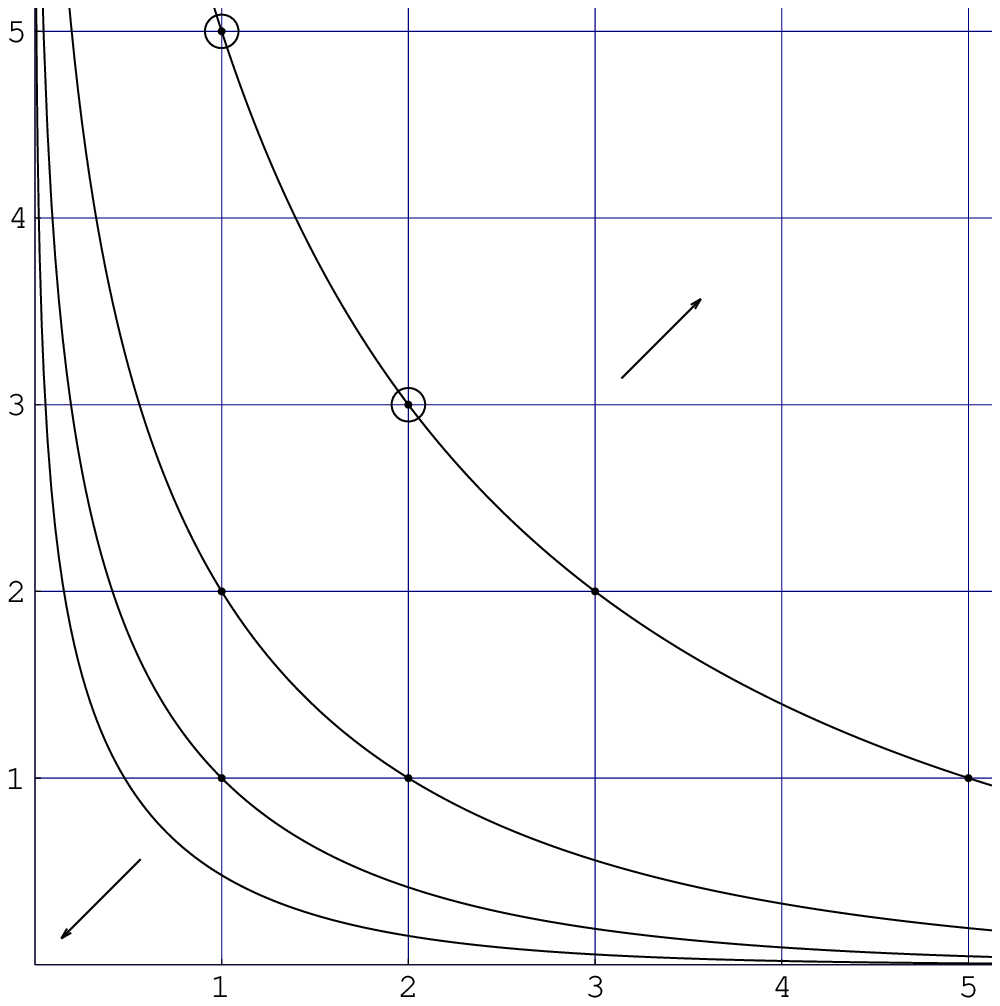}}
\put(52,0){\makebox(0,0)[l]{$L_{1}$}}
\put(0,52){\makebox(0,0)[b]{$L_{2}$}}
\put(40,14){\makebox(0,0)[bl]{$x\approx0.7549$}}
\put(19.9,10.1){\makebox(0,0)[bl]{$x=\frac{\sqrt{5}-1}{2}\approx0.6180$}}
\put(14,6.4){\makebox(0,0)[bl]{$x=\frac{1}{2}$}}
\put(3.2,13.4){\makebox(0,0)[tr]{$2^{-\sqrt{2}}=x$}}
\put(3,3){\makebox(0,0)[tl]{$x\searrow0$}}
\put(33,33){\makebox(0,0)[tl]{$x\nearrow1$}}
\end{picture}
\caption{Examples from lattice points: $\frac{\mbox{\protect\raisebox
{0pt}[10pt][2pt]{$\scriptstyle
dL_{2}$}}}{dL_{1}}=-x^{L_{1}-L_{2}}$.}%
\label{Lfromx}%
\end{figure}

For $d=2$, the picture is as shown in Figure \ref{Lfromx}. There is some
non-uniqueness as follows: if%
\[
m_{1}x^{L_{1}}+m_{2}x^{L_{2}}+\dots+m_{k}x^{L_{k}}=1
\]
where $L_{i},m_{i}\in\mathbb{N}$, then, if $Q\left(  x\right)  $ is any
polynomial of the form \linebreak $Q\left(  x\right)  =1+\sum_{m=1}^{k}%
n_{m}x^{L_{m}}$, where $n_{m}$ are positive integers, $0\leq n_{m}\leq L_{m}$,
then\ $Q\left(  x\right)  \left(  m_{1}x^{L_{1}}+\dots+m_{k}x^{L_{k}%
}-1\right)  =0$, which gives another polynomial of the form above. Even the
added condition%
\[
\sum_{i}m_{i}=d
\]
does not imply uniqueness, by Example \ref{ExaRemNewBis.3} and Chapter
\ref{APP.EXA}. \setcounter{figurelink}{\value{figure}}

\chapter{\label{APP.EXA}Complete classification in the cases $\lambda=2$,
$N=2,3,4$}

\setcounter{figure}{\value{figurelink}} \setcounter{table}{\value{tablelink}}
The examples when the Perron--Frobenius eigenvalue $\lambda=2$ entail some of
the essential features of the dimension groups associated with the
corresponding AF-algebras $\mathfrak{A}_{L}$.

The construction in the examples below is a special case of the following: Let
$p\left(  x\right)  \in\mathbb{Z}\left[  x\right]  $ be given, and assume it
is irreducible. Let $N\in\mathbb{N}$, and let $\mathcal{F}_{N}\left(
p\right)  $ be the set of $N\times N$ matrices over $\mathbb{Z}$ of the
form\label{LOSJ_7}
\begin{equation}
J=%
\begin{pmatrix}
\vphantom{\vdots}m_{1} & 1 & \cdots & 0 & 0\\
\vphantom{\vdots}m_{2} & 0 & \ddots & 0 & 0\\
\vdots &  & \ddots & \ddots & \vdots\\
\vphantom{\vdots}m_{N-1} & 0 &  & 0 & 1\\
\vphantom{\vdots}m_{N} & 0 & \cdots & 0 & 0
\end{pmatrix}
\label{eqResv.1}%
\end{equation}
with $m_{N}\geq1$, and $m_{i}\geq0$, such that $p\left(  x\right)  $ divides
$p_{J}\left(  x\right)  =\sum_{i=1}^{N}m_{i}x^{i}-1$. We saw that, for
$p\left(  x\right)  =2x-1$, $x_{N}=\#\mathcal{F}_{N}\left(  2x-1\right)  $ is
finite for all $N=2,3,\dots$. An analogue of this holds true in general. If
$p\left(  x\right)  =x^{2}+x^{3}-1$ \textup{(}see Example
\textup{\ref{ExaRemNewBis.3}),} then $\mathcal{F}_{4}\left(  p\right)
=\varnothing$, while $x_{5}\left(  p\right)  =\#\mathcal{F}_{5}\left(
p\right)  =2$. The two elements of $\mathcal{F}_{5}\left(  p\right)  $ are%
\[%
\begin{pmatrix}
m_{1}\\
m_{2}\\
m_{3}\\
m_{4}\\
m_{5}%
\end{pmatrix}
\colon%
\begin{pmatrix}
1\\
0\\
0\\
0\\
1
\end{pmatrix}
,\;%
\begin{pmatrix}
0\\
0\\
1\\
1\\
1
\end{pmatrix}
\;,
\]
corresponding to isomorphic $\mathfrak{A}_{L}$'s. This approach is in general
most useful if $p\left(  x\right)  $ has the form $\sum_{i=1}^{D}n_{i}x^{i}-1$
where the $n_{i}$ are nonnegative integers, $n_{D}\neq0$ and $\gcd\left\{
i\mid n_{i}\neq0\right\}  =1$, because then we can say at the outset that the
Perron--Frobenius eigenvalue of \textup{(\ref{eqResv.1})} is $1/a$, where $a$
is the unique positive root of $p\left(  x\right)  $.

For each $N=2,3,4,\dots$, there is only a finite number $x_{N}$ of
possibilities for the matrix $J$. Since the matrix $J$ is of the form
\textup{(\ref{eqResv.1})} they are described by the numbers $m_{i}$ of the
first column. They are given by the following algorithm: If $Q_{1}%
,\dots,Q_{N-1}\in\mathbb{Z}$, and $q\left(  x\right)  =1+Q_{N-1}x+\dots
+Q_{1}x^{N-1}$, then $p_{L}\left(  x\right)  =\left(  2x-1\right)  q\left(
x\right)  $ has the form\label{LOSpLx_7}%
\[
p_{L}\left(  x\right)  =-1+m_{1}x+\dots+m_{N}x^{N}%
\]
with $m_{i}\geq0$, $m_{N}>0$ if and only if $Q_{1}>0$ and%
\begin{equation}
\left\{
\begin{aligned}
\vphantom{\vdots}Q_{1} &\leq2Q_{2},\\
\vphantom{\vdots}Q_{2} &\leq2Q_{3},\\
\vdots& \qquad\vdots\\
\vphantom{\vdots}Q_{N-2} &\leq2Q_{N-1},\\
\vphantom{\vdots}Q_{N-1} &\leq2.
\end{aligned}%
\right.  \label{eqInvPalle.pound}%
\end{equation}
This is proved by simple algebra. The numbers $x_{N}$ are $x_{2}=2$, $x_{3}%
=6$, $x_{4}=26$, $\dots$. In a slightly more condensed form, the conditions
are
\[
0<Q_{1}\leq2\cdot Q_{2}\leq4\cdot Q_{3}\leq8\cdot Q_{4}\leq\dots\leq
2^{N-2}\cdot Q_{N-1}\leq2^{N-1}.
\]
It follows that the specific cases may be summarized for $N=2$, $3$, and $4$
\textup{(}lexicographic order from \textup{(\ref{eqInvPalle.pound})):}%
\begin{align}%
\begin{pmatrix}
m_{1}\\
m_{2}%
\end{pmatrix}
&  \colon%
\begin{pmatrix}
1\\
2
\end{pmatrix}
,%
\begin{pmatrix}
0\\
4
\end{pmatrix}
^{\ast},\label{eqInvPalle.a}\\%
\begin{pmatrix}
m_{1}\\
m_{2}\\
m_{3}%
\end{pmatrix}
&  \colon%
\begin{pmatrix}
1\\
1\\
2
\end{pmatrix}
,%
\begin{pmatrix}
1\\
0\\
4
\end{pmatrix}
,%
\begin{pmatrix}
0\\
3\\
2
\end{pmatrix}
,%
\begin{pmatrix}
0\\
2\\
4
\end{pmatrix}
,%
\begin{pmatrix}
0\\
1\\
6
\end{pmatrix}
,%
\begin{pmatrix}
0\\
0\\
8
\end{pmatrix}
^{\ast},\label{eqInvPalle.b}\\%
\intertext{and}%
\begin{pmatrix}
m_{1}\\
m_{2}\\
m_{3}\\
m_{4}%
\end{pmatrix}
&  \colon%
\begin{pmatrix}
1\\
1\\
1\\
2
\end{pmatrix}
,%
\begin{pmatrix}
1\\
1\\
0\\
4
\end{pmatrix}
,%
\begin{pmatrix}
1\\
0\\
3\\
2
\end{pmatrix}
,%
\begin{pmatrix}
1\\
0\\
2\\
4
\end{pmatrix}
,%
\begin{pmatrix}
1\\
0\\
1\\
6
\end{pmatrix}
,%
\begin{pmatrix}
1\\
0\\
0\\
8
\end{pmatrix}
,%
\begin{pmatrix}
0\\
3\\
1\\
2
\end{pmatrix}
,%
\begin{pmatrix}
0\\
3\\
0\\
4
\end{pmatrix}
^{\llap{$\scriptstyle\ast$}},%
\begin{pmatrix}
0\\
2\\
3\\
2
\end{pmatrix}
,\label{eqInvPalle.c}\\
&  \phantom{{}\colon{}}%
\begin{pmatrix}
0\\
2\\
2\\
4
\end{pmatrix}
,%
\begin{pmatrix}
0\\
2\\
1\\
6
\end{pmatrix}
,%
\begin{pmatrix}
0\\
2\\
0\\
8
\end{pmatrix}
^{\llap{$\scriptstyle\ast$}},%
\begin{pmatrix}
0\\
1\\
5\\
2
\end{pmatrix}
,%
\begin{pmatrix}
0\\
1\\
4\\
4
\end{pmatrix}
,%
\begin{pmatrix}
0\\
1\\
3\\
6
\end{pmatrix}
,%
\begin{pmatrix}
0\\
1\\
2\\
8
\end{pmatrix}
,%
\begin{pmatrix}
0\\
1\\
1\\
10
\end{pmatrix}
,%
\begin{pmatrix}
0\\
1\\
0\\
12
\end{pmatrix}
^{\llap{$\scriptstyle\ast$}},\nonumber\\
&  \phantom{{}\colon{}}%
\begin{pmatrix}
0\\
0\\
7\\
2
\end{pmatrix}
,%
\begin{pmatrix}
0\\
0\\
6\\
4
\end{pmatrix}
,%
\begin{pmatrix}
0\\
0\\
5\\
6
\end{pmatrix}
,%
\begin{pmatrix}
0\\
0\\
4\\
8
\end{pmatrix}
,%
\begin{pmatrix}
0\\
0\\
3\\
10
\end{pmatrix}
,%
\begin{pmatrix}
0\\
0\\
2\\
12
\end{pmatrix}
,%
\begin{pmatrix}
0\\
0\\
1\\
14
\end{pmatrix}
,%
\begin{pmatrix}
0\\
0\\
0\\
16
\end{pmatrix}
^{\llap{$\scriptstyle\ast$}}.\nonumber
\end{align}
But inspection reveals that, of the two $N=2$ cases
\textup{(\ref{eqInvPalle.a})}, only the first one has $L_{i}$-values with
greatest common divisor equal to one. For the six $N=3$ cases
\textup{(\ref{eqInvPalle.b})}, all but the last of them have this property.
Finally, for the $N=4$ examples \textup{(\ref{eqInvPalle.c})}, the property
holds for all but the $8$th, $12$th, $18$th and the last one. Note that the
$\mathfrak{A}_{L}$'s associated to the last vector in each of the three
sequences are all isomorphic, and isomorphic to the algebra defined by the
constant $1\times1$ incidence matrix $%
\begin{pmatrix}
2
\end{pmatrix}
$, i.e., the Glimm algebra of type $2^{\infty}$ illustrated in Figure
\textup{\ref{BratDiag11}. The algebras }corresponding to the third and seventh
vectors in the $N=4$ sequence \textup{(\ref{eqInvPalle.c})} are illustrated in
Figure \textup{\ref{BratDiagsd6}.}

In order to distinguish the isomorphism classes of the remaining specimens, we
will use the invariants developed in Chapters \textup{\ref{Brunt}%
--\ref{SubG0}.} Since $\lambda=2$ in these cases, we always have $\tau\left(
K_{0}\left(  \mathfrak{A}_{L}\right)  \right)  =\mathbb{Z}\left[  \frac{1}%
{2}\right]  $ and $D=1$. Thus the invariants in Chapter
\textup{\ref{Pediferient}} reduce to $N$, $\operatorname*{Prim}\left(
m_{N}\right)  $, and $\operatorname*{Prim}\left(  Q_{N-1}\right)
=\operatorname*{Prim}\left(  m_{N}/2\right)  $. It follows from Remark
\textup{\ref{RemSesquilabialNew.15}} that in these cases $m_{N}=2$ if and only
if $\ker\tau\cong\mathbb{Z}^{N-1}$, and then Corollary
\textup{\ref{corollary10.12}} can be used to distinguish some cases. Here
$\alpha=\left(  1,\frac{1}{2},\frac{1}{4},\dots,\frac{1}{2^{N-1\mathstrut}%
}\right)  $ by \textup{(\ref{eqRemNewBis.17})} and $v$ is given by
\textup{(\ref{eqSubNew.pound})} with $\lambda=2$. (The case $\lambda=m_{N}$
will be studied in detail in Chapter \ref{ClmN}.) We will also use some
secondary invariants derived from the group $G_{0}=\ker\tau$,\label{LOSG0_9}
since of course any group invariant derived from $G_{0}$ is an invariant for
$\mathfrak{A}_{L}$. Since $G_{0}$ is a natural $\mathbb{Z}$-module, the tensor
product group $G_{0}\otimes\mathbb{Z}_{p}$ \textup{(}tensor products as
$\mathbb{Z}$-modules\/\textup{)} is a secondary invariant. For example
$\mathbb{Z}_{q}\otimes\mathbb{Z}_{p}\cong\mathbb{Z}_{\gcd\left(  q,p\right)
}$, $\mathbb{Z\otimes Z}_{p}\cong\mathbb{Z}_{p}$ and $\mathbb{Z}\left[
\frac{1}{q}\right]  \otimes\mathbb{Z}_{p}\cong\mathbb{Z}_{\left(  p,q\right)
}$ where $\left(  p,q\right)  =p/\gcd\left(  p,q^{\infty}\right)  $ is defined
by \textup{(\ref{eqApp.Cla13}).} For our specimens, it will be useful to use
$p=2$, and using \textup{(\ref{eqRemFeb.a}):} $G_{0}\cong\mathbb{Z}\left[
x\right]  \diagup\left(  p_{0}\left(  x\right)  \right)  $, we obtain
$G_{0}\otimes\mathbb{Z}_{2}\cong\mathbb{Z}_{2}\left(  x\right)  \diagup\left(
p_{0,2}\left(  x\right)  \right)  $ where $p_{0,2}\left(  x\right)  $ is the
polynomial obtained from $p_{0}\left(  x\right)  $ by reducing the
coefficients modulo $2$. Thus $G_{0}\otimes\mathbb{Z}_{2}$ is the direct sum
of a finite number of copies of $\mathbb{Z}_{2}$, and this finite number is an
invariant. See Corollary \ref{CorKer.D}.

Then to work on the list \textup{(\ref{eqInvPalle.a})--(\ref{eqInvPalle.c}%
).\smallskip}

\noindent\textbf{Rank 2} ($N=2$): By \textup{(\ref{eqInvPalle.a})} there is
only one specimen, $%
\begin{pmatrix}
1\\
2
\end{pmatrix}
$.\textup{\smallskip}

\noindent\textbf{Rank 3} ($N=3$): By \textup{(\ref{eqInvPalle.b}),} there are
five specimens, which can first be classified as follows:\label{LOSPrim_7}%
\[%
\begin{tabular}
[c]{cll}%
\begin{tabular}
[c]{l}%
Group\\
number
\end{tabular}
& $%
\begin{array}
[c]{l}%
\operatorname*{Prim}\left(  m_{3}\right) \\
\operatorname*{Prim}\left(  m_{3}/2\right)
\end{array}
$ & \\\hline
1 & $%
\begin{array}
[c]{l}%
\left\{  2\right\} \\
\varnothing
\end{array}
$ & $\vphantom{\begin{pmatrix}0\\0\\0\\0\end{pmatrix}}%
\begin{pmatrix}
1\\
1\\
2
\end{pmatrix}
^{\text{(a)}},%
\begin{pmatrix}
0\\
3\\
2
\end{pmatrix}
^{\text{(b)}},$\\\hline
2 & $%
\begin{array}
[c]{l}%
\left\{  2\right\} \\
\left\{  2\right\}
\end{array}
$ & $\vphantom{\begin{pmatrix}0\\0\\0\\0\end{pmatrix}}%
\begin{pmatrix}
1\\
0\\
4
\end{pmatrix}
^{\text{(a)}},%
\begin{pmatrix}
0\\
2\\
4
\end{pmatrix}
^{\text{(b)}},$\\\hline
3 & $%
\begin{array}
[c]{l}%
\left\{  2,3\right\} \\
\left\{  3\right\}
\end{array}
$ & $\vphantom{\begin{pmatrix}0\\0\\0\\0\end{pmatrix}}%
\begin{pmatrix}
0\\
1\\
6
\end{pmatrix}
.$\\\hline
\end{tabular}
\]

For \textbf{Group 1}, we may use Corollary \textup{\ref{corollary10.12}.} But
this was done already in Example \textup{\ref{ExaCyc.end1bis}} with the result
that these two specimens define non-isomorphic algebras. (Specimen (b) will be
considered further in Example \ref{Exa2}.)

For \textbf{Group 2}, we compute, using \textup{(\ref{eqKer.p})} in Corollary
\textup{\ref{CorKer.D},}%
\[%
\begin{tabular}
[t]{rccc}%
$\text{Group 2}$ & Specimen & $p_{0}\left(  x\right)  $ & $G_{0}%
\otimes\mathbb{Z}_{2}$\\\hline
(a) & $\vphantom{\begin{pmatrix}0\\0\\0\\0\end{pmatrix}}%
\begin{pmatrix}
1\\
0\\
4
\end{pmatrix}
$ & $2x^{2}+x+1$ & $\mathbb{Z}_{2}$\\\hline
(b) & $\vphantom{\begin{pmatrix}0\\0\\0\\0\end{pmatrix}}%
\begin{pmatrix}
0\\
2\\
4
\end{pmatrix}
$ & $2x^{2}+2x+1$ & $0$\\\hline
\end{tabular}
\]
and hence $G_{0}$ is non-isomorphic for the two examples. (Specimen (a) here
has already been studied in Remark \ref{RemPediferient.nada}, and will again
be considered in Example \ref{Exa1}.)

We conclude that all the $5$ specimens in the $N=3$ case are mutually
non-isomorphic.\textup{\smallskip}

\noindent\textbf{Rank 4} ($N=4$): Here we have $22$ specimens which can be
divided into $6$ groups according to the invariants $\operatorname*{Prim}%
\left(  m_{4}\right)  $, $\operatorname*{Prim}\left(  m_{4}/2\right)  $:%
\[%
\begin{tabular}
[c]{lll}%
$%
\begin{tabular}
[c]{l}%
Group\\
number
\end{tabular}
$ & $%
\begin{array}
[c]{l}%
\operatorname*{Prim}\left(  m_{4}\right) \\
\operatorname*{Prim}\left(  m_{4}/2\right)
\end{array}
$ & \\\hline
\multicolumn{1}{c}{1} & $%
\begin{array}
[c]{l}%
\left\{  2\right\} \\
\varnothing
\end{array}
$ & $\vphantom{\begin{pmatrix}0\\0\\0\\0\\0\end{pmatrix}}%
\begin{pmatrix}
1\\
1\\
1\\
2
\end{pmatrix}
^{\!\!\!(1)}%
\begin{pmatrix}
1\\
0\\
3\\
2
\end{pmatrix}
^{\!\!\!(2)}%
\begin{pmatrix}
0\\
3\\
1\\
2
\end{pmatrix}
^{\!\!\!(3)}%
\begin{pmatrix}
0\\
2\\
3\\
2
\end{pmatrix}
^{\!\!\!(4)}%
\begin{pmatrix}
0\\
1\\
5\\
2
\end{pmatrix}
^{\!\!\!(5)}%
\begin{pmatrix}
0\\
0\\
7\\
2
\end{pmatrix}
^{\!\!\!(6)}$\\\hline
\multicolumn{1}{c}{2} & $%
\begin{array}
[c]{l}%
\left\{  2\right\} \\
\left\{  2\right\}
\end{array}
$ & $\vphantom{\begin{pmatrix}0\\0\\0\\0\\0\end{pmatrix}}%
\begin{pmatrix}
1\\
1\\
0\\
4
\end{pmatrix}
^{\!\!\!\text{(a)}}%
\begin{pmatrix}
1\\
0\\
2\\
4
\end{pmatrix}
^{\!\!\!\text{(b)}}%
\begin{pmatrix}
1\\
0\\
0\\
8
\end{pmatrix}
^{\!\!\!\text{(c)}}%
\begin{pmatrix}
0\\
2\\
2\\
4
\end{pmatrix}
^{\!\!\!\text{(d)}}%
\begin{pmatrix}
0\\
1\\
4\\
4
\end{pmatrix}
^{\!\!\!\text{(e)}}%
\begin{pmatrix}
0\\
1\\
2\\
8
\end{pmatrix}
^{\!\!\!\text{(f)}}%
\begin{pmatrix}
0\\
0\\
6\\
4
\end{pmatrix}
^{\!\!\!\text{(g)}}%
\begin{pmatrix}
0\\
0\\
4\\
8
\end{pmatrix}
^{\!\!\!\text{(h)}}$\\\hline
\multicolumn{1}{c}{3} & $%
\begin{array}
[c]{l}%
\left\{  2,3\right\} \\
\left\{  3\right\}
\end{array}
$ & $\vphantom{\begin{pmatrix}0\\0\\0\\0\\0\end{pmatrix}}%
\begin{pmatrix}
1\\
0\\
1\\
6
\end{pmatrix}
^{\!\!\!\text{(a)}}%
\begin{pmatrix}
0\\
2\\
1\\
6
\end{pmatrix}
^{\!\!\!\text{(b)}}%
\begin{pmatrix}
0\\
1\\
3\\
6
\end{pmatrix}
^{\!\!\!\text{(c)}}%
\begin{pmatrix}
0\\
0\\
5\\
6
\end{pmatrix}
^{\!\!\!\text{(d)}}$\\\hline
\multicolumn{1}{c}{4} & $%
\begin{array}
[c]{l}%
\left\{  2,3\right\} \\
\left\{  2,3\right\}
\end{array}
$ & $\vphantom{\begin{pmatrix}0\\0\\0\\0\\0\end{pmatrix}}%
\begin{pmatrix}
0\\
0\\
2\\
12
\end{pmatrix}
$\\\hline
\multicolumn{1}{c}{5} & $%
\begin{array}
[c]{l}%
\left\{  2,5\right\} \\
\left\{  5\right\}
\end{array}
$ & $\vphantom{\begin{pmatrix}0\\0\\0\\0\\0\end{pmatrix}}%
\begin{pmatrix}
0\\
1\\
1\\
10
\end{pmatrix}
^{\!\!\!\text{(a)}}%
\begin{pmatrix}
0\\
0\\
3\\
10
\end{pmatrix}
^{\!\!\!\text{(b)}}$\\\hline
\multicolumn{1}{c}{6} & $%
\begin{array}
[c]{l}%
\left\{  2,7\right\} \\
\left\{  7\right\}
\end{array}
$ & $\vphantom{\begin{pmatrix}0\\0\\0\\0\\0\end{pmatrix}}%
\begin{pmatrix}
0\\
0\\
1\\
14
\end{pmatrix}
$\\\hline
\end{tabular}
\]

For \textbf{Group 1}, we may apply Corollary \textup{\ref{corollary10.12}.}
Using \textup{(\ref{eqRemNewBis.17})} with $e^{-\beta}=\frac{1}{2}$ and
\textup{(\ref{eqSubNew.pound})} together with $p_{L}\left(  2\right)  =0$, we
compute for two general matrices of the form \textup{(\ref{eqSub.pound.st})}
with $N=4$ and $\lambda=2$ that%
\[%
\ip{\alpha}{v}%
=\frac{1}{8}\left(  4\left(  m_{1}+m_{2}\right)  +3m_{3}+2m_{4}\right)  .
\]
For the $6$ specimens in group 1 this leads to%
\[%
\ip{\alpha}{v_{i}}%
=\frac{13+2i}{8}%
\]
for $i=1,2,\dots,6$. It follows that%
\[
\frac{%
\ip{\alpha}{v_{i}}%
}{%
\ip{\alpha}{v_{j}}%
}\notin\mathbb{Z}\left[  \tfrac{1}{2}\right]
\]
whenever $i\neq j$, $i,j=1,\dots,6$. We conclude from Corollary
\textup{\ref{corollary10.12}} that all these $6$ specimens define mutually
non-isomorphic algebras $\mathfrak{A}_{L}$, although each defines an exact
sequence of the form%
\[
0\longrightarrow\mathbb{Z}^{3}\longrightarrow G\longrightarrow\mathbb{Z}%
\left[  \tfrac{1}{2}\right]  \longrightarrow0.
\]
In \textbf{Group 2} there are $8$ specimens. Let us compute the polynomial
$p_{0}\left(  x\right)  =p_{L}\left(  x\right)  /\left(  2x-1\right)  $ for
these, and use the result to compute $G_{0}\otimes\mathbb{Z}_{2}$, where
$G_{0}=\ker\tau=\mathbb{Z}\left[  x\right]  \diagup\left(  p_{0}\left(
x\right)  \right)  $,\label{LOSG0_10} using \textup{(\ref{eqKer.p})} in
Corollary \textup{\ref{CorKer.D}. The result is exhibited in Table
\ref{TableR4G2Specimens} on the next page.}

\begin{table}[ptb]
\caption{The specimens in Group $2$ for Rank $4$.}%
\label{TableR4G2Specimens}
\begin{tabular}
[c]{rccc}%
Group 2 & Specimen & $p_{0}\left(  x\right)  $\label{LOSp0x_2} & $G_{0}\otimes_{Z}Z_{2}%
$\label{LOSG0_11}\\\hline
(a) & $\vphantom{\renewcommand{\arraystretch}{1.15}\begin{pmatrix}%
0\\0\\0\\0\\0\end{pmatrix}\renewcommand{\arraystretch}{1}}%
\begin{pmatrix}
1\\
1\\
0\\
4
\end{pmatrix}
$ & $1+x+x^{2}+2x^{3}$ & $Z_{2}^{2}$\\\hline
(b) & $\vphantom{\renewcommand{\arraystretch}{1.15}\begin{pmatrix}%
0\\0\\0\\0\\0\end{pmatrix}\renewcommand{\arraystretch}{1}}%
\begin{pmatrix}
1\\
0\\
2\\
4
\end{pmatrix}
$ & $%
\begin{array}
[c]{l}%
1+x+2x^{2}+2x^{3}\\
=\left(  1+x\right)  \left(  1+2x^{3}\right)
\end{array}
$ & $Z_{2}$\\\hline
(c) & $\vphantom{\renewcommand{\arraystretch}{1.15}\begin{pmatrix}%
0\\0\\0\\0\\0\end{pmatrix}\renewcommand{\arraystretch}{1}}%
\begin{pmatrix}
1\\
0\\
0\\
8
\end{pmatrix}
$ & $1+x+2x^{2}+4x^{3}$ & $Z_{2}$\\\hline
(d) & $\vphantom{\renewcommand{\arraystretch}{1.15}\begin{pmatrix}%
0\\0\\0\\0\\0\end{pmatrix}\renewcommand{\arraystretch}{1}}%
\begin{pmatrix}
0\\
2\\
2\\
4
\end{pmatrix}
$ & $1+2x+2x^{2}+2x^{3}$ & $0$\\\hline
(e) & $\vphantom{\renewcommand{\arraystretch}{1.15}\begin{pmatrix}%
0\\0\\0\\0\\0\end{pmatrix}\renewcommand{\arraystretch}{1}}%
\begin{pmatrix}
0\\
1\\
4\\
4
\end{pmatrix}
$ & $%
\begin{array}
[c]{l}%
1+2x+3x^{2}+2x^{3}\\
=\left(  1+x\right)  \left(  1+x+2x^{2}\right)
\end{array}
$ & $Z_{2}^{2}$\\\hline
(f) & $\vphantom{\renewcommand{\arraystretch}{1.15}\begin{pmatrix}%
0\\0\\0\\0\\0\end{pmatrix}\renewcommand{\arraystretch}{1}}%
\begin{pmatrix}
0\\
1\\
2\\
8
\end{pmatrix}
$ & $1+2x+3x^{2}+4x^{3}$ & $Z_{2}^{2}$\\\hline
(g) & $\vphantom{\renewcommand{\arraystretch}{1.15}\begin{pmatrix}%
0\\0\\0\\0\\0\end{pmatrix}\renewcommand{\arraystretch}{1}}%
\begin{pmatrix}
0\\
0\\
6\\
4
\end{pmatrix}
$ & $1+2x+4x^{2}+2x^{3}$ & $0$\\\hline
(h) & $\vphantom{\renewcommand{\arraystretch}{1.15}\begin{pmatrix}%
0\\0\\0\\0\\0\end{pmatrix}\renewcommand{\arraystretch}{1}}%
\begin{pmatrix}
0\\
0\\
4\\
8
\end{pmatrix}
$ & $1+2x+4x^{2}+4x^{3}$ & $0$\\\hline
\end{tabular}
\end{table}

It follows from \textup{Table \ref{TableR4G2Specimens} }that we can group the
$8$ specimens into $3$ subgroups with no isomorphism between the different
subgroups:\textup{\smallskip}

\noindent\textbf{Subgroup 1}%
\[
\text{(d):}%
\begin{pmatrix}
0\\
2\\
2\\
4
\end{pmatrix}
,\text{\quad(g):}%
\begin{pmatrix}
0\\
0\\
6\\
4
\end{pmatrix}
,\text{\quad(h):}%
\begin{pmatrix}
0\\
0\\
4\\
8
\end{pmatrix}
.
\]
Here $\ker\tau\otimes\mathbb{Z}_{2}=0$, so $\ker\tau$ is a torsion-free
abelian group of rank $3$ such that all the elements are divisible by $2$, and
also $\ker\left(  \tau\right)  \subseteq\mathbb{Z}\left[  \frac{1}{2}\right]
^{3}$. It follows that%
\[
\ker\tau=\mathbb{Z}\left[  \tfrac{1}{2}\right]  ^{3}.
\]
Thus $G$ is an extension%
\[
0\longrightarrow\mathbb{Z}\left[  \tfrac{1}{2}\right]  ^{3}\longrightarrow
G\longrightarrow\mathbb{Z}\left[  \tfrac{1}{2}\right]  \longrightarrow0
\]
in all three cases. But $\operatorname*{Ext}\left(  \mathbb{Z}\left[  \frac
{1}{2}\right]  ,\mathbb{Z}\left[  \frac{1}{2}\right]  ^{3}\right)  =0$ by
\cite[Proposition VI.2.1]{CaEi56}. \textup{(}If one assumes \emph{a priori}
that $G$ is divisible by $2$ this is trivial, but in the general case one
proceeds as follows: It is clear that $\operatorname*{Ext}\left(
\mathbb{Z}\left[  \tfrac{1}{2}\right]  ,\mathbb{Z}\left[  \frac{1}{2}\right]
^{3}\right)  \cong\operatorname*{Ext}\left(  \mathbb{Z}\left[  \tfrac{1}%
{2}\right]  ,\mathbb{Z}\left[  \tfrac{1}{2}\right]  \right)  ^{3}$
\textup{(}i.e., three copies\/\textup{).} Assume that%
\begin{equation}
0\longrightarrow\mathbb{Z}\left[  \tfrac{1}{2}\right]  \longrightarrow
M\longrightarrow\mathbb{Z}\left[  \tfrac{1}{2}\right]  \longrightarrow0
\label{eqApp.ExaAprstar}%
\end{equation}
is an exact sequence of $\mathbb{Z}$-modules. Since $\mathbb{Z}\left[
\frac{1}{2}\right]  =\mathbb{Z}_{\left\{  2\right\}  }$ \textup{(}localized in
$\left\{  2\right\}  $\textup{),} $\mathbb{Z}\left[  \frac{1}{2}\right]  $
will be $\mathbb{Z}$-flat. Take the tensor product of
\textup{(\ref{eqApp.ExaAprstar})} with $\mathbb{Z}\left[  \frac{1}{2}\right]
$ over $\mathbb{Z}$ to obtain%
\begin{equation}
0\longrightarrow\mathbb{Z}\left[  \tfrac{1}{2}\right]  \longrightarrow
M\otimes_{\mathbb{Z}}\mathbb{Z}\left[  \tfrac{1}{2}\right]  \longrightarrow
\mathbb{Z}\left[  \tfrac{1}{2}\right]  \longrightarrow0,
\label{eqApp.ExaAprstarstar}%
\end{equation}
which has to be isomorphic to \textup{(\ref{eqApp.ExaAprstar}).} But
\textup{(\ref{eqApp.ExaAprstarstar})} splits by the initial remark.\textup{)}
Thus the three vectors in subgroup 1 define isomorphic algebras.

This can also be seen much more directly as follows: It follows directly from
Corollary \textup{\ref{corollary10.A}} that $G=\mathbb{Z}\left[  \frac{1}%
{2}\right]  ^{4}$ for these three specimens (d,g,h) and the Perron--Frobenius
eigenvalue $\lambda=2$ in all three cases, so it follows from
\textup{(\ref{eqRemNewBis.17})} and \textup{(\ref{eqRemNewBis.34})} that the
three specimens are isomorphic.\textup{\smallskip}

\noindent\textbf{Subgroup 2}%
\[
\text{(b):}%
\begin{pmatrix}
1\\
0\\
2\\
4
\end{pmatrix}
,\text{\quad(c):}%
\begin{pmatrix}
1\\
0\\
0\\
8
\end{pmatrix}
.
\]
In specimen (b), $p_{0}\left(  x\right)  =\left(  1+x\right)  \left(
1+2x^{2}\right)  $ so $\ker\tau$ is given by an extension%
\[
0\longrightarrow\mathbb{Z}\left[  \tfrac{1}{2}\right]  ^{2}\longrightarrow
\ker\tau\longrightarrow\mathbb{Z}\longrightarrow0
\]
\textup{(}the right morphism is evaluation of the polynomial at $-1$, where
$p_{0}\left(  -1\right)  =0$\textup{).} By $\mathbb{Z}$ being a free
$\mathbb{Z}$-module, this extension automatically splits. In specimen (c),
$p_{0}\left(  x\right)  =1+x+2x^{2}+4x^{3}$ is irreducible, and we have an
exact sequence%
\[
0\longrightarrow\mathbb{Z}\longrightarrow\ker\tau\longrightarrow
\mathbb{Z}\left[  \tfrac{1}{2}\right]  ^{2}\longrightarrow0.
\]
We will show that these specimens are non-isomorphic by using Corollary
\textup{\ref{corollary10.19}.} Let $J$ correspond to specimen (b) and
$J^{\prime}$ to specimen (c). By Proposition \textup{\ref{prop10.22}} we have%
\[
\left\{  h\in F^{\prime}\mid8h=0\right\}  =F_{3}^{\prime}%
\]
and hence \textup{(\ref{eq10.83})} is fulfilled:%
\[
\left\{  g\in G^{\prime}\mid2^{3}g\in\mathbb{Z}^{N}\right\}  =G_{3}^{\prime}.
\]
One now computes $%
\ip{\alpha}{v}%
=\frac{9}{4}$ and $%
\ip{\alpha^{\prime}}{v^{\prime}}%
=\frac{5}{2}$, so%
\[
\frac{%
\ip{\alpha}{v}%
}{%
\ip{\alpha^{\prime}}{v^{\prime}}%
}=\frac{9}{10}\notin\mathbb{Z}\left[  \tfrac{1}{2}\right]  .
\]
It follows from Corollary \textup{\ref{corollary10.19}} that there does not
exist a unital morphism $\mathfrak{A}_{J}\rightarrow\mathfrak{A}_{J^{\prime}}%
$, and in particular specimens (b) and (c) are
non-isomorphic.\textup{\smallskip}

\noindent\textbf{Subgroup 3}%
\[
\text{(a):}%
\begin{pmatrix}
1\\
1\\
0\\
4
\end{pmatrix}
,\text{\quad(e):}%
\begin{pmatrix}
0\\
1\\
4\\
4
\end{pmatrix}
,\text{\quad(f):}%
\begin{pmatrix}
0\\
1\\
2\\
8
\end{pmatrix}
.
\]
Specimen (e) has the reducible polynomial
\[
p_{0}\left(  x\right)  =2x^{3}+3x^{2}+2x+1=\left(  x+1\right)  \left(
2x^{2}+x+1\right)
\]
so there is an exact diagram%
\begin{equation}%
\begin{array}
[c]{ccc}%
& 0 & \\
& \downarrow & \\
0\longrightarrow\mathbb{Z}\longrightarrow\!\! & E & \!\!\longrightarrow
\mathbb{Z}\left[  \frac{1}{2}\right]  \longrightarrow0.\\
& \downarrow & \\
& \makebox[0pt]{\hss$\ker\tau$\hss} & \\
& \downarrow & \\
& \mathbb{Z} & \\
& \downarrow & \\
& 0 &
\end{array}
\label{eqApp.ExaJun.8}%
\end{equation}
The horizontal sequence is described in detail in the end of Example
\textup{\ref{Exa1}.} The vertical sequence necessarily splits since
$\mathbb{Z}$ is free, i.e.,\label{LOSG0_12}%
\[
G_{0}=\ker\tau\cong\mathbb{Z}\oplus E.
\]
Since $p_{0}\left(  1\right)  =0$, evaluation at $-1$ gives a homomorphism
$\ker\left(  \tau\right)  \rightarrow\mathbb{Z}$, which is the lower vertical
map in the diagram. Specimens (a) and (f) have irreducible $p_{0}%
$-polynomials, so there are no homomorphisms $G\rightarrow\mathbb{Z}$, and
therefore these are non-isomorphic to specimen (e). But specimens (a) and (f)
are mutually non-isomorphic by Example \textup{\ref{example10.23}.} Hence all
three specimens (a,e,f) are mutually non-isomorphic. Note also that for (e),
it follows from Proposition \textup{\ref{prop10.22}} that $\left\{  h\in
F^{\text{(e)}}\mid4h=0\right\}  =F_{2}^{\text{(e)}}$ and hence%
\[
\left\{  g\in F^{\text{(e)}}\mid4g=0\right\}  =G_{2}^{\text{(e)}}.
\]
For this specimen we have%
\[%
\ip{\alpha^{\text{(e)}}}{v^{\text{(e)}}}%
=%
\begin{pmatrix}
1 & \frac{1}{2} & \frac{1}{4} & \frac{1}{8}%
\end{pmatrix}%
\begin{pmatrix}
1\\
2\\
3\\
2
\end{pmatrix}
=3
\]
and since%
\[%
\ip{\alpha^{\text{(a)}}}{v^{\text{(a)}}}%
=2\text{\quad and\quad}%
\ip{\alpha^{\text{(f)}}}{v^{\text{(f)}}}%
=\frac{13}{4},
\]
it follows also directly from Corollary \textup{\ref{corollary10.19}} and
Remark \textup{\ref{remark10.20}} that these specimens all are non-isomorphic.

In \textbf{Group 3}, there are $4$ specimens. We compute the polynomial
$p_{0}\left(  x\right)  =p_{L}\left(  x\right)  /\left(  2x-1\right)  $ for
these, and use the result to compute $G_{0}\otimes\mathbb{Z}_{3}$, where
$G_{0}=\ker\tau$ =$\mathbb{Z}\left[  x\right]  \diagup\left(  p_{0}\left(
x\right)  \mathbb{Z}\left[  x\right]  \right)  $, using \textup{(\ref{eqKer.p}%
)} in Corollary \textup{\ref{CorKer.D}.}\label{LOSp0x_3}%
\[%
\begin{tabular}
[t]{rccc}%
$\text{Group 3}$ & Specimen & $p_{0}\left(  x\right)  $ &
$G_{0}\otimes\mathbb{Z}_{3}$\label{LOSG0_13}\\\hline
(a) & $\vphantom{\begin{pmatrix}0\\0\\0\\0\\0\end{pmatrix}}%
\begin{pmatrix}
1\\
0\\
1\\
6
\end{pmatrix}
$ & $1+x+2x^{2}+3x^{3}$ & $\mathbb{Z}_{3}^{2}$\\\hline
(b) & $\vphantom{\begin{pmatrix}0\\0\\0\\0\\0\end{pmatrix}}%
\begin{pmatrix}
0\\
2\\
1\\
6
\end{pmatrix}
$ & $1+2x+2x^{2}+3x^{3}$ & $\mathbb{Z}_{3}^{2}$\\\hline
(c) & $\vphantom{\begin{pmatrix}0\\0\\0\\0\\0\end{pmatrix}}%
\begin{pmatrix}
0\\
1\\
3\\
6
\end{pmatrix}
$ & $1+2x+3x^{2}+3x^{3}$ & $\mathbb{Z}_{3}$\\\hline
(d) & $\vphantom{\begin{pmatrix}0\\0\\0\\0\\0\end{pmatrix}}%
\begin{pmatrix}
0\\
0\\
5\\
6
\end{pmatrix}
$ & $%
\begin{array}
[c]{l}%
1+2x+4x^{2}+3x^{3}\\
=\left(  1+x\right)  \left(  1+x+3x^{2}\right)
\end{array}
$ & $\mathbb{Z}_{3}^{2}$\\\hline
\end{tabular}
\]
Thus we see immediately that specimen (c) is non-isomorphic to the three
others. Also, we see that specimen (d) permits a homomorphism from $\ker\tau$
into $\mathbb{Z}$, but not the two others, so it remains to consider the pair
(a,b). But for both these specimens we have $m_{4}=6$ and $m_{3}=1$, so
applying Proposition \textup{\ref{prop10.22}} we have%
\[
\left\{  h\in F\mid6h=0\right\}  =F_{1}%
\]
for both. But this means that $6$-$\deg J=\deg J=1$ for both. But $%
\ip{\alpha^{\text{(a)}}}{v^{\text{(a)}}}%
=\frac{19}{8}$ and $%
\ip{\alpha^{\text{(b)}}}{v^{\text{(b)}}}%
=\frac{23}{8}$, and hence Theorem \textup{\ref{theorem10.10}} implies that
there cannot be a unital homomorphism from either of these specimens into the
other. It follows that the four specimens in this group are mutually non-isomorphic.\smallskip

\textbf{Groups 4 and 6} have only one specimen each, so the basic invariants
$\operatorname*{Prim}\left(  m_{4}\right)  $\label{LOSPrim_8} and
$\operatorname*{Prim}\left(  m_{4}/2\right)  $ suffice to separate this group
from the others.\smallskip

In \textbf{Group 5} there are two specimens,%
\[
\text{(a):}%
\begin{pmatrix}
0\\
1\\
1\\
10
\end{pmatrix}
,\text{\quad(b):}%
\begin{pmatrix}
0\\
0\\
3\\
10
\end{pmatrix}
.
\]
For both of these $m_{4}=10$ and $m_{3}\neq0$ and $m_{3}$ is mutually prime
with $m_{4}$. Thus Proposition \textup{\ref{prop10.22}} implies that%
\[
10\text{-}\deg J=\deg J=1
\]
for both these specimens. But in this case $%
\ip{\alpha^{\text{(a)}}}{v^{\text{(a)}}}%
=\frac{27}{8}$ and $%
\ip{\alpha^{\text{(b)}}}{v^{\text{(b)}}}%
=\frac{29}{8}$, and it follows from Theorem \textup{\ref{theorem10.10}} that
there are no unital homomorphisms from one of these two into the other.\smallskip

\begin{figure}[ptb]
\begin{picture}(360,438)
\put(0,0){\includegraphics
[bb=253 47 419 672,clip,width=116.2bp,height=437.5bp]{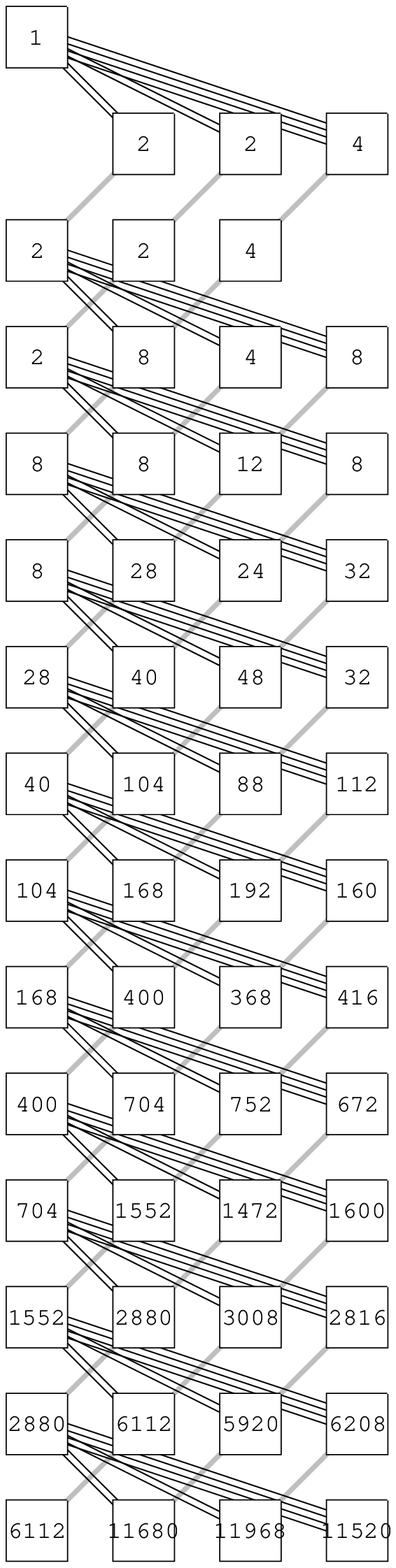}}
\put(121.9,0){\includegraphics
[bb=253 47 419 672,clip,width=116.2bp,height=437.5bp]{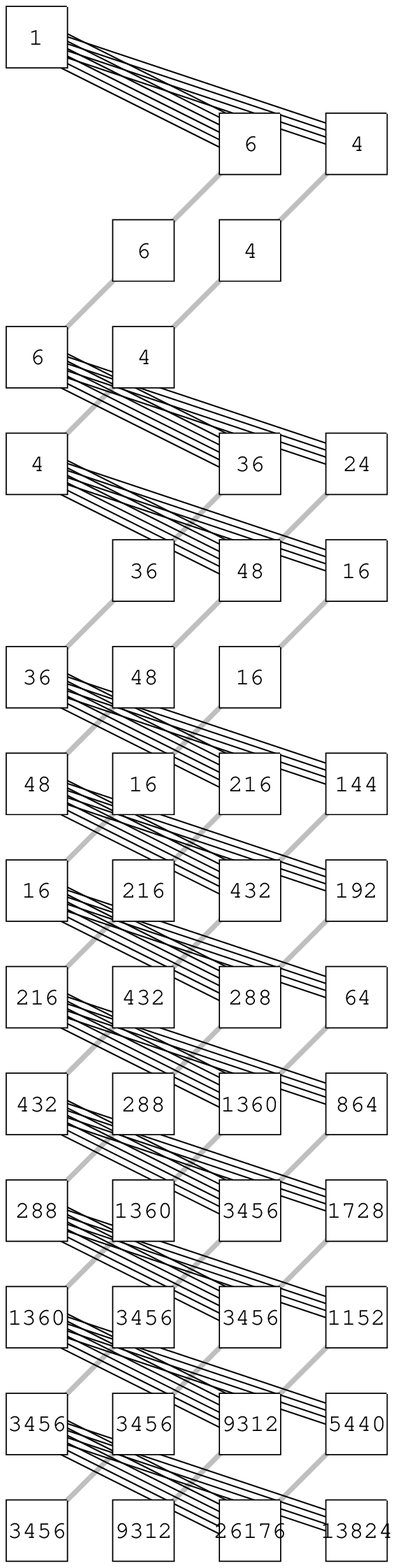}}
\put(243.8,0){\includegraphics
[bb=253 47 419 672,clip,width=116.2bp,height=437.5bp]{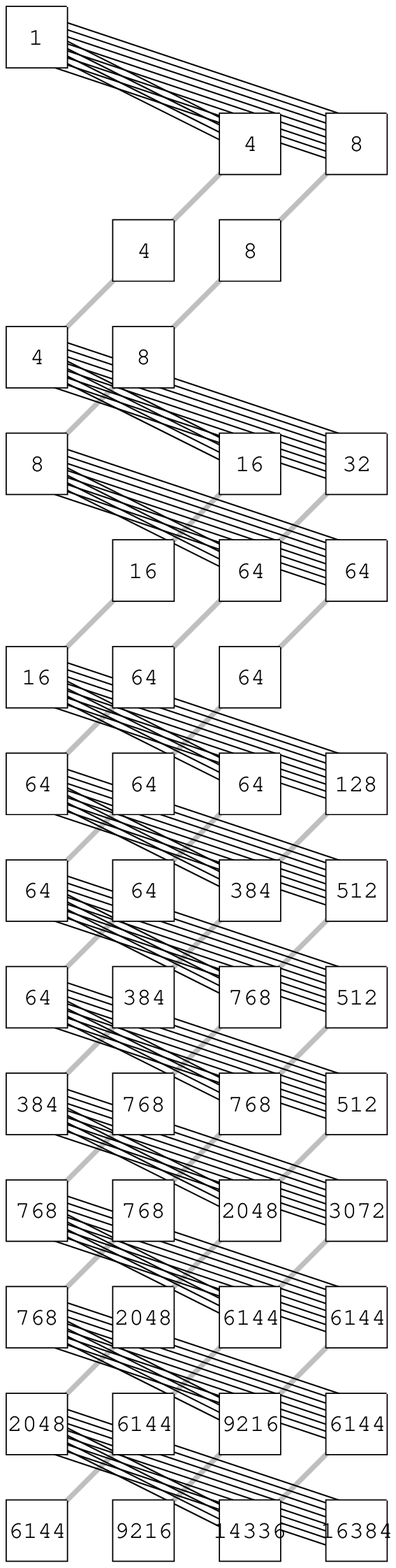}}
\end{picture}
\caption{$L=\{2,2,3,3,4,4,4,4\}$; first column $=(0\;2\;2\;4)^{\mathrm{t}}$
(left); $L=\{3,3,3,3,3,3,4,4,4,4\}$; first column $=(0\;0\;6\;4)^{\mathrm{t}}$
(center); $L=\{3,3,3,3,4,4,4,4,4,4,4,4\}$; first column
$=(0\;0\;4\;8)^{\mathrm{t}}$ (right): The three isomorphic algebras in the
final summary of Chapter \ref{APP.EXA} (see (\ref{eqApp.ExaJun.9})). }%
\label{BratDiagsN4g2s1}%
\end{figure}

In summary, all the $22$ permitted specimens in \textup{(\ref{eqInvPalle.c})}
are mutually non-isomorphic except for one group of $3$ isomorphic specimens,
namely the ones in Subgroup 1 of Group 2:%
\begin{equation}%
\begin{pmatrix}
0\\
2\\
2\\
4
\end{pmatrix}
,%
\begin{pmatrix}
0\\
0\\
6\\
4
\end{pmatrix}
,%
\begin{pmatrix}
0\\
0\\
4\\
8
\end{pmatrix}
. \label{eqApp.ExaJun.9}%
\end{equation}
These specimens are illustrated in Figure \ref{BratDiagsN4g2s1}.

\setcounter{figurelink}{\value{figure}} \setcounter{tablelink}{\value{table}}

\chapter{\label{ClmN}Complete classification in the case $\lambda=m_{N}$}

\setcounter{figure}{\value{figurelink}} We will now continue the study of the
case $\lambda=m_{N}$ begun in Corollary \ref{corollary10.12}, Lemma
\ref{lemma10.13} and Proposition \ref{prop10.14}. We will introduce a new
invariant $I\left(  J\right)  $ in (\ref{eqClmN.L}) and Corollary
\ref{CorClmN.E1} below. In the case $N=1$, the invariant always has the value
$1$ (Section \ref{ClmN.1}), but in the case $N=2$, the invariant separates all
specimens in this class, so they are all non-isomorphic (Section
\ref{ClmN.2}). More interestingly, in the case $N=3$, the pair $\left(
\operatorname*{Prim}\lambda,I\left(  J\right)  \right)  $ turns out to be a
\emph{complete} invariant (Theorem \ref{ThmClmN.2}), and this can be used to
exhibit nontrivial pairs of $3\times3$ matrices in this class giving
isomorphic algebras. It is curious that for $N=3$, the equality $I\left(
J\right)  =I\left(  J^{\prime}\right)  $ forces $\lambda$ to be unequal to
$\lambda^{\prime}$ (unless $J=J^{\prime}$; Proposition \ref{ProClmN.1}). For
$N=4$ we do exhibit an isomorphic pair with $\lambda=\lambda^{\prime}$, and we
present a proof of K.H. Kim and F. Roush that $\left(  N,\operatorname*{Prim}%
\lambda,I\left(  J\right)  \right)  $ is a complete invariant in general for
the class $\lambda=m_{N}$ (Theorem \ref{ThmClmNJun.18}). In this class, it
also turns out that $\left(  N,\operatorname*{Prim}\lambda,I\left(  J\right)
\right)  $ is a complete invariant for \emph{stable} isomorphism (see
Corollary \ref{CorClmNJun.21}).

\begin{lemma}
\label{LemClmN.A}Let $J$ be a matrix of the form \textup{(\ref{eqCycNew.1}),}
and assume that the Perron--Frobenius eigenvalue $\lambda$ of $J$ is equal to
$m_{N}$.\label{LOSmN_2} It follows that\label{LOSG}\label{LOSZ1lambda_2}%
\begin{equation}
G:=\bigcup_{n}J^{-n}\mathbb{Z}^{N}=\mathbb{Z}^{N}+\mathbb{Z}\left[  \tfrac
{1}{\lambda}\right]  \mathbf{v}, \label{eqClmN.A}%
\end{equation}
where $\mathbf{v}$ is the right Perron--Frobenius eigenvector given by
\textup{(\ref{eqSubNew.pound}).}
\end{lemma}

\begin{proof}
Clearly $\mathbb{Z}^{N}\subseteq G$ and $\mathbf{v}\in\mathbb{Z}^{N}$, so
$J^{-n}\mathbf{v}=\lambda^{-n}\mathbf{v}\in G$, and hence%
\begin{equation}
\mathbb{Z}^{N}+\mathbb{Z}\left[  \tfrac{1}{\lambda}\right]  \mathbf{v}%
\subseteq G. \label{eqClmN.B}%
\end{equation}
Since $G$ is the smallest $J^{-1}$-invariant subgroup of $\mathbb{Q}^{N}$, in
order to show the converse inclusion it suffices to show that%
\begin{equation}
J^{-1}\left(  \mathbb{Z}^{N}+\mathbb{Z}\left[  \tfrac{1}{\lambda}\right]
\mathbf{v}\right)  \subseteq\mathbb{Z}^{N}+\mathbb{Z}\left[  \tfrac{1}%
{\lambda}\right]  \mathbf{v}. \label{eqClmN.C}%
\end{equation}
But as $J^{-1}\mathbf{v}=\frac{1}{\lambda}\mathbf{v}$, it suffices to show
that%
\begin{equation}
J^{-1}\mathbb{Z}^{N}\subseteq\mathbb{Z}^{N}+\mathbb{Z}\left[  \tfrac
{1}{\lambda}\right]  \mathbf{v.} \label{eqClmN.D}%
\end{equation}
Then it suffices to show that the right column vector in $J^{-1}$ in
(\ref{eqThetaNyet}) is in $\mathbb{Z}^{N}+\mathbb{Z}\left[  \frac{1}{\lambda
}\right]  \mathbf{v}$. But using $m_{N}=\lambda$ and
(\textup{\ref{eqSubNew.pound}) we see that this column vector has the form
}$\lambda^{-1}\mathbf{v}+\mathbf{m}$, where $\mathbf{m}\in\mathbb{Z}^{N}$.
This proves (\ref{eqClmN.D}) and thus Lemma \ref{LemClmN.A} is proved.
\end{proof}

\begin{lemma}
\label{LemClmN.B}Let $J$ be a matrix of the form \textup{(\ref{eqCycNew.1}),}
and assume that $\lambda=m_{N}$. Let $\mathbf{v}$ be the right
Perron--Frobenius eigenvector of $J$ normalized as in
\textup{(\ref{eqSubNew.pound}).} It follows that $\mathbf{v}$ has the
form\label{LOSPerronFrobeniusrighteigenvector_5}%
\begin{equation}
\mathbf{v}=%
\begin{pmatrix}
1\\
v_{2}\\
\vdots\\
v_{N-1}\\
1
\end{pmatrix}
\label{eqClmN.E}%
\end{equation}
where%
\begin{equation}
v_{i+1}=\lambda v_{i}-m_{i} \label{eqClmN.F}%
\end{equation}
for $i=1,2,\dots,N-1$.
\end{lemma}

\begin{proof}
This is an immediate consequence of (\ref{eq10.37}) and (\ref{eq10.38}) in the
proof of Lemma \ref{lemma10.13}.
\end{proof}

\begin{lemma}
\label{LemClmN.C}Let $J$, $\lambda=m_{N}$, $\mathbf{v}$ be as in Lemma
\textup{\ref{LemClmN.B}} and define the $\left(  N-1\right)  \times N$ matrix
$M_{\mathbf{v}}$ by%
\begin{equation}
M_{\mathbf{v}}=%
\begin{pmatrix}
v_{2} & -1 & 0 & \cdots & 0 & 0\\
v_{3} & 0 & -1 &  & 0 & 0\\
\vdots & \vdots &  & \ddots &  & \vdots\\
v_{N-1} & 0 & 0 &  & -1 & 0\\
1 & 0 & 0 & \cdots & 0 & -1
\end{pmatrix}
. \label{eqClmN.G}%
\end{equation}
Then%
\begin{equation}
G:=\bigcup_{n}J^{-n}\mathbb{Z}^{N}=\left\{  x\in\mathbb{Z}\left[  \tfrac
{1}{\lambda}\right]  ^{N}\bigm|M_{\mathbf{v}}x\in\mathbb{Z}^{N-1}\right\}  .
\label{eqClmN.H}%
\end{equation}
\end{lemma}

\begin{proof}
Let $\mathbf{v}^{\perp}=\left\{  \mathbf{y}\in\mathbb{Z}^{N}\mid%
\ip{\mathbf{y}}{\mathbf{v}}%
=0\right\}  $. We will also prove that%
\begin{equation}
G=\left\{  \mathbf{x}\in\mathbb{Z}\left[  \tfrac{1}{\lambda}\right]  ^{N}%
\bigm|\forall\,\mathbf{y}\in\mathbf{v}^{\perp}\colon%
\ip{\mathbf{y}}{\mathbf{x}}%
\in\mathbb{Z}\right\}  \label{eqClmN.I}%
\end{equation}
by establishing the following relations between the right-hand sets of
(\ref{eqClmN.A}), (\ref{eqClmN.H}) and (\ref{eqClmN.I}):%
\[
\text{(\ref{eqClmN.A})}_{\mathrm{R}}\subseteq\text{(\ref{eqClmN.I}%
)}_{\mathrm{R}}\subseteq\text{(\ref{eqClmN.H})}_{\mathrm{R}}\subseteq
\text{(\ref{eqClmN.A})}_{\mathrm{R}}.
\]
The first inclusion to the left is immediate, and since the vectors
\[
\left(  v_{2},-1,0,\dots,0\right)  ,\left(  v_{3},0,-1,0,\dots,0\right)
,\dots,\left(  1,0,\dots,0,-1\right)
\]
are all in $\mathbf{v}^{\perp}$ by (\ref{eqClmN.E}), the second inclusion
follows. But if $\mathbf{x}\in$(\ref{eqClmN.H})$_{\mathrm{R}}$, then%
\begin{align*}
x_{1}  &  \in\mathbb{Z}\left[  \tfrac{1}{\lambda}\right] \\
x_{2}  &  =v_{2}x_{1}+m_{2}\\
x_{2}  &  =v_{3}x_{1}+m_{3}\\
\vdots &  \qquad\vdots\\
x_{N-1}  &  =v_{N-1}x_{1}+m_{N-1}\\
x_{N}  &  =x_{1}+m_{N}=v_{N}x_{1}+m_{N}%
\end{align*}
where $m_{2},\dots,m_{N}\in\mathbb{Z}$, and hence $\mathbf{x}=\mathbf{m}%
+x_{1}\mathbf{v}\in$(\ref{eqClmN.A})$_{\mathrm{R}}$.
\end{proof}

\begin{lemma}
\label{LemClmN.D}Let $J$, $J^{\prime}$ be as in Lemma \textup{\ref{LemClmN.B}%
,} and assume that there is an isomorphism $\theta\colon G\rightarrow
G^{\prime}$ \textup{(}and thus $N=N^{\prime}$, $\operatorname*{Prim}\left(
\lambda\right)  =\operatorname*{Prim}\left(  \lambda^{\prime}\right)
$\textup{).} Let $\Lambda\in\operatorname*{GL}\left(  N,\mathbb{Z}\left[
\frac{1}{\lambda}\right]  \right)  $ be the matrix in Proposition
\textup{\ref{prop10.7}} implementing the isomorphism. It follows
that\label{LOSLambda_3}%
\begin{equation}
\Lambda\mathbf{v}=\xi\mathbf{v}^{\prime} \label{eqClmN.J}%
\end{equation}
where $\xi$\label{LOSxiGreek} is an element of $\mathbb{Z}\left[  \frac
{1}{\lambda}\right]  $ with multiplicative inverse \textup{(}i.e., $\xi$ is a
product of powers of the primes in $\operatorname*{Prim}\left(  \lambda
\right)  $\textup{).}
\end{lemma}

\begin{proof}
By Lemma \ref{LemClmN.A} we have%
\begin{align*}
G  &  =\mathbb{Z}^{N}+\mathbb{Z}\left[  \tfrac{1}{\lambda}\right]
\mathbf{v},\\
G^{\prime}  &  =\mathbb{Z}^{N}+\mathbb{Z}\left[  \tfrac{1}{\lambda}\right]
\mathbf{v}^{\prime},
\end{align*}
so\label{LOSDlambdaG_2}%
\begin{align*}
D_{\lambda}\left(  G\right)   &  =\mathbb{Z}\left[  \tfrac{1}{\lambda}\right]
\mathbf{v},\\%
\intertext{and}%
D_{\lambda}\left(  G^{\prime}\right)   &  =\mathbb{Z}\left[  \tfrac{1}%
{\lambda}\right]  \mathbf{v}^{\prime}%
\end{align*}
(see (\ref{eq10.55}) for the definition of $D_{\lambda}$). But $\theta$ must
map $D_{\lambda}\left(  G\right)  $ onto $D_{\lambda}\left(  G^{\prime
}\right)  $, from which the assertion follows.
\end{proof}

Note that Lemma \ref{LemClmN.D} immediately gives a strengthening of Corollary
\ref{corollary10.12}. But we will do better: see Lemma \ref{LemClmN.H}.

\begin{corollary}
\label{CorClmN.E}Let $J$, $J^{\prime}$ be matrices of the form
\textup{(\ref{eqCycNew.1})} with $m_{N}=\lambda$ and $m_{N^{\prime}}%
=\lambda^{\prime}$. If there is a unital isomorphism $\mathfrak{A}%
_{J}\rightarrow\mathfrak{A}_{J^{\prime}}$, then $N=N^{\prime}$,
$\operatorname*{Prim}\left(  \lambda\right)  =\operatorname*{Prim}\left(
\lambda^{\prime}\right)  =\left\{  p_{1},\dots,p_{n}\right\}  $, and there are
integers $m_{1},\dots,m_{n}$ such that%
\begin{equation}
\frac{%
\ip{\alpha}{\mathbf{v}}%
}{%
\ip{\alpha^{\prime}}{\mathbf{v}^{\prime}}%
}=p_{1}^{m_{1}}p_{2}^{m_{2}}\cdots p_{n}^{m_{n}}. \label{eqClmN.K}%
\end{equation}
\end{corollary}

\begin{proof}
We have $N=N^{\prime}$ and $\operatorname*{Prim}\left(  m_{N}\right)
=\operatorname*{Prim}\left(  m_{N^{\prime}}^{\prime}\right)  $ by Theorem
\ref{CorCyc.9}. From (\ref{eqClmN.J}) and Proposition \ref{prop10.7}, it
follows that%
\[%
\ip{\alpha}{\mathbf{v}}%
=%
\ip{\alpha^{\prime}}{\Lambda\mathbf{v}}%
=%
\ip{\alpha^{\prime}}{\xi\mathbf{v}^{\prime}}%
=\xi%
\ip{\alpha^{\prime}}{\mathbf{v}^{\prime}}%
.
\]
\end{proof}

We will now give a useful alternative form of Corollary \ref{CorClmN.E} by
means of the number\label{LOSIJ_2}%
\begin{equation}
I\left(  J\right)  =\sum_{i=1}^{N}v_{i}\lambda^{N-i}=\lambda^{N-1}%
+v_{2}\lambda^{N-2}+\dots+v_{N-1}\lambda+1, \label{eqClmN.L}%
\end{equation}
where $v_{1}=1,v_{2},\dots,v_{N-1},v_{N}=1$ are the components of the right
Perron--Frobenius eigenvector in Lemma \ref{LemClmN.B}. The next corollary
says that $I\left(  J\right)  $ is an invariant in the context $\lambda=m_{N}$.

\begin{corollary}
\label{CorClmN.E1}Let $J$, $J^{\prime}$ be matrices of the form
\textup{(\ref{eqCycNew.1})} with $m_{N}=\lambda$ and $m_{N^{\prime}}^{\prime
}=\lambda^{\prime}$. If $\mathfrak{A}_{J}$ is isomorphic to $\mathfrak
{A}_{J^{\prime}}$ then%
\begin{equation}
I\left(  J\right)  =I\left(  J^{\prime}\right)  . \label{eqClmN.M}%
\end{equation}
\end{corollary}

\begin{proof}
Using Lemma \ref{CorCyc.4} and (\ref{eqClmN.E}) we have%
\begin{align*}%
\ip{\alpha}{\mathbf{v}}%
&  =1+v_{2}\frac{1}{\lambda}+\dots+v_{N-1}\frac{1}{\lambda^{N-2\mathstrut}%
}+\frac{1}{\lambda^{N-1\mathstrut}}\\
&  =\frac{1}{\lambda^{N-1\mathstrut}}\left(  \lambda^{N-1}+v_{2}\lambda
^{N-2}+\dots+v_{N-1}\lambda+1\right)  ,
\end{align*}
and a similar expression is valid for $%
\ip{\alpha^{\prime}}{\mathbf{v}^{\prime}}%
$. Combining this with (\ref{eqClmN.K}) we find an element $\xi^{\prime}%
\in\mathbb{Z}\left[  \frac{1}{\lambda}\right]  $ with a multiplicative inverse
such that%
\[
\xi^{\prime}\left(  \lambda^{N-1}+v_{2}\lambda^{N-2}+\dots+v_{N-1}%
\lambda+1\right)  =\left(  \lambda^{\prime\,N-1}+v_{2}^{\prime\,N-2}%
+\dots+v_{N-1}^{\prime}\lambda^{\prime}+1\right)  .
\]
Now, we may find two disjoint subsets of $\left\{  p_{1},\dots,p_{n}\right\}
$, say $P_{+}$ and $P_{-}$, such that%
\[
\xi^{\prime}=\prod_{p\in P_{+}}p^{n\left(  p\right)  }\prod_{q\in P_{-}%
}q^{-n\left(  q\right)  },
\]
where $n\left(  p\right)  \in\mathbb{N}$, $n\left(  q\right)  \in\mathbb{N}$.
The relation above may be written%
\begin{multline*}
\prod_{p\in P_{+}}p^{n\left(  p\right)  }\left(  \lambda^{N-1}+v_{2}%
\lambda^{N-2}+\dots+v_{N-1}\lambda+1\right) \\
=\prod_{q\in P_{-}}q^{n\left(  q\right)  }\left(  \lambda^{\prime\,N-1}%
+v_{2}^{\prime}\lambda^{\prime\,N-2}+\dots+v_{N-1}^{\prime}\lambda^{\prime
}+1\right)  .
\end{multline*}
Since the primes $p_{1},\dots,p_{n}$ are all distinct, and all of them are
factors of both $\lambda$ and $\lambda^{\prime}$, and thus none of them are
factors of the polynomials $\left(  \cdots+1\right)  $ above, it follows that
$P_{+}=P_{-}=\varnothing$. Thus $\xi^{\prime}=1$. But this means%
\[
\lambda^{N-1}%
\ip{\alpha}{\mathbf{v}}%
=\lambda^{\prime\,N-1}%
\ip{\alpha^{\prime}}{\mathbf{v}^{\prime}}%
\]
and%
\begin{align}
I\left(  J\right)   &  =\lambda^{N-1}+v_{2}\lambda^{N-2}+\dots+v_{N-1}%
\lambda+1\label{eqClmN.N}\\
&  =\lambda^{\prime\,N-1}+v_{2}^{\prime}\lambda^{\prime\,N-2}+\dots
+v_{N-1}^{\prime}\lambda^{\prime}+1\nonumber\\
&  =I\left(  J^{\prime}\right)  .%
\settowidth{\qedskip}{$\displaystyle I\left
( J\right) =I\left( J^{\prime}\right) .$}\settowidth{\qedadjust}%
{$\displaystyle I\left( J\right) =\lambda^{\prime\,N-1}+v_{2}^{\prime}%
\lambda^{\prime\,N-2}+\dots+v_{N-1}^{\prime}\lambda^{\prime}+1$}%
\addtolength{\qedadjust}{\textwidth}\addtolength{\qedskip}{-0.5\qedadjust
}\rlap{\makebox[-\qedskip][r]{\qedsymbol}}%
\nonumber
\end{align}%
\renewcommand{\qed}{}%
\end{proof}

Under some circumstances, Corollary \ref{CorClmN.E1} can be used to give more
amenable conditions for isomorphism.

\begin{corollary}
\label{LemClmN.F}Let $J$, $J^{\prime}$ be matrices of the form
\textup{(\ref{eqCycNew.1})} with $m_{N}=\lambda$ and $m_{N^{\prime}}^{\prime
}=\lambda^{\prime}$.

If there is a unital isomorphism $\mathfrak{A}_{J}\rightarrow\mathfrak
{A}_{J^{\prime}}$ and $\lambda=\lambda^{\prime}$, then
\begin{equation}%
\ip{\alpha}{\mathbf{v}}%
=%
\ip{\alpha^{\prime}}{\mathbf{v}^{\prime}}%
. \label{eqClmN.M1}%
\end{equation}

If there is a unital isomorphism $\mathfrak{A}_{J}\rightarrow\mathfrak
{A}_{J^{\prime}}$ and $\lambda^{\prime}$ is an integer multiple of $\lambda$,
then%
\begin{equation}
\lambda=\lambda^{\prime}. \label{eqClmN.M2}%
\end{equation}
\end{corollary}

\begin{proof}
The first statement follows from the formula\label{LOStauv_2}%
\[%
\ip{\alpha}{\mathbf{v}}%
=\lambda^{-\left(  N-1\right)  }I\left(  J\right)
\]
in the beginning of the proof of Corollary \ref{CorClmN.E1}, as well as from
the corollary itself, and the fact that $N$ is an isomorphism invariant
(Theorem \ref{CorCyc.9}).

For the second statement, note that (\ref{eqClmN.M}) implies%
\begin{equation}
v_{2}\lambda^{N-2}+\dots+v_{N-1}\lambda=\lambda^{\prime\,N-1}-\lambda
^{N-1}+\left(  v_{2}^{\prime}\lambda^{\prime\,N-2}+\dots+v_{N-1}^{\prime
}\lambda^{\prime}\right)  , \label{eqClmN.O}%
\end{equation}
and the expression in parentheses is positive. We have assumed that
$\lambda^{\prime}$ is an integer multiple of $\lambda$, and if this multiple
is $>1$ we will show the contradiction%
\begin{equation}
v_{2}\lambda^{N-2}+\dots+v_{N-1}\lambda<\lambda^{\prime\,N-1}-\lambda^{N-1}.
\label{eqClmN.P}%
\end{equation}
This will prove the lemma. Since%
\begin{align}
v_{2}\lambda^{N-2}+\dots+v_{N-1}\lambda &  \leq\lambda\lambda^{N-2}%
+\lambda^{2}\lambda^{N-3}+\dots+\lambda^{N-2}\lambda\label{eqClmN.Q}\\
&  =\left(  N-2\right)  \lambda^{N-1}\nonumber
\end{align}
by (\ref{eqSubNew.pound}), (\ref{eqClmN.P}) will follow if we can show that%
\[
\left(  N-1\right)  \lambda^{N-1}\leq\lambda^{\prime\,N-1}-\lambda^{N-1}%
\]
or%
\[
\left(  N-1\right)  <\left(  \lambda^{\prime}/\lambda\right)  ^{N-1}.
\]
But as $\lambda^{\prime}$ is an integer multiple of $\lambda$, this says%
\[
\left(  N-1\right)  <2^{N-1},
\]
which is obvious when $N>1$.
\end{proof}

\begin{lemma}
\label{LemClmN.G}Adopt the assumptions and notation in Lemma
\textup{\ref{LemClmN.A}.} Then $G$ has a direct sum decomposition%
\begin{equation}
G=\mathbb{Z}^{N-1}\oplus\mathbb{Z}\left[  \tfrac{1}{\lambda}\right]
\mathbf{v}, \label{eqClmN.R}%
\end{equation}
where $\mathbb{Z}^{N-1}$ identifies with the elements of $\mathbb{Z}^{N}$ with
zero first coordinate. If $\alpha=\left(  1,\frac{1}{\lambda},\dots,\frac
{1}{\lambda^{N-1\mathstrut}}\right)  $\label{LOSalphaeigenvector_4} and
$\beta=\left(  \frac{1}{\lambda},\frac{1}{\lambda^{2\mathstrut}},\dots
,\frac{1}{\lambda^{N-1\mathstrut}}\right)  $,\label{LOSbeta_1} an element
$\mathbf{x}\oplus\xi\mathbf{v}$ of $G$ is positive if and only if $%
\ip{\beta}{\mathbf{x}}%
+\xi%
\ip{\alpha}{\mathbf{v}}%
>0$.
\end{lemma}

\begin{proof}
Put $H=\mathbb{Z}^{N-1}\oplus\mathbb{Z}\left[  \frac{1}{\lambda}\right]
\mathbf{v}$. Since $v_{1}=1$, this sum is really direct, and it follows from
(\ref{eqClmN.A}) that $H\subseteq G$. Conversely, if $\mathbf{y}+\xi
\mathbf{v}\in G$, define $\xi^{\prime}=\xi+y_{1}\in\mathbb{Z}\left[  \frac
{1}{\lambda}\right]  $ and write $\mathbf{y}+\xi\mathbf{v}=\left(
\mathbf{y}-y_{1}\mathbf{v}\right)  +\xi^{\prime}\mathbf{v}$. But
$\mathbf{y}-y_{1}\mathbf{v}\in\mathbb{Z}^{N-1}$ since $v_{1}=1$, and hence
$\mathbf{y}+\xi\mathbf{v}=\left(  \mathbf{y}-y_{1}\mathbf{v}\right)
+\xi^{\prime}\mathbf{v}\in\mathbb{Z}^{N-1}+\mathbb{Z}\left[  \frac{1}{\lambda
}\right]  \mathbf{v}$. Thus $G\subseteq H$ and $G=H$. Since $\mathbf{x}\in G$
is positive if and only if $%
\ip{\alpha}{\mathbf{x}}%
>0$ the last statement is clear.
\end{proof}

\begin{lemma}
\label{LemClmN.H}Let $J$, $J^{\prime}$ be matrices of the form
\textup{(\ref{eqCycNew.1})} with $\lambda=m_{N}$ and $\lambda^{\prime
}=m_{N^{\prime}}^{\prime}$, $\operatorname*{Prim}\left(  \lambda\right)
=\operatorname*{Prim}\left(  \lambda^{\prime}\right)  $, $N=N^{\prime}$, and
$I\left(  J\right)  =I\left(  J^{\prime}\right)  $, so that%
\begin{equation}%
\begin{aligned} G
&=\mathbb{Z}^{N}+\mathbb{Z}\left[ \tfrac{1}{\lambda}\right] \mathbf{v}
=\mathbb{Z}^{N-1}\oplus\mathbb{Z}\left[ \tfrac{1}{\lambda}\right] \mathbf
{v}, \\
G^{\prime} &=\mathbb{Z}^{N}+\mathbb{Z}\left[ \tfrac
{1}{\lambda}\right] \mathbf{v}^{\prime}
=\mathbb{Z}^{N-1}\oplus\mathbb{Z}
\left[ \tfrac{1}{\lambda}\right] \mathbf{v}^{\prime}
\end{aligned}%
\label{eqClmN.S}%
\end{equation}
by Lemma \textup{\ref{LemClmN.A}} and Lemma \textup{\ref{LemClmN.G}.} Then any
unital order isomorphism $\theta\colon G\rightarrow G^{\prime}$ has the form%
\begin{equation}
\theta\left(  \mathbf{x},\xi\mathbf{v}\right)  =\left(  A\mathbf{x},\left(
\eta\left(  \mathbf{x}\right)  +\xi\lambda^{\prime\,N-1}\lambda^{-\left(
N-1\right)  }\right)  \mathbf{v}^{\prime}\right)  \label{eqClmN.T}%
\end{equation}
relative to the right decompositions, where%
\begin{align}
A  &  \in\operatorname*{GL}\left(  N-1,\mathbb{Z}\right)  ,\label{eqClmN.U}\\
\eta &  \in\operatorname*{Hom}\left(  \mathbb{Z}^{N-1},\mathbb{Z}\left[
\tfrac{1}{\lambda}\right]  \right)  ,\label{eqClmN.V}\\%
\ip{\beta}{\mathbf{x}}%
&  =%
\ip{\beta^{\prime}}{A\mathbf{x}}%
+\eta\left(  \mathbf{x}\right)
\ip{\alpha^{\prime}}{\mathbf{v}^{\prime}}%
,\label{eqClmN.W}\\
A\mathbf{\bar{v}}  &  =\mathbf{\bar{v}}^{\prime}\text{\quad and\quad}%
\eta\left(  \mathbf{\bar{v}}\right)  =\lambda^{\prime\,N-1}\lambda^{-\left(
N-1\right)  }-1 \label{eqClmN.X}%
\end{align}
where%
\begin{equation}
\mathbf{\bar{v}}=%
\begin{pmatrix}
v_{2}\\
\vdots\\
v_{N-1}\\
1
\end{pmatrix}
\text{\quad and\quad}\mathbf{\bar{v}}^{\mathbf{\prime}}=%
\begin{pmatrix}
v_{2}^{\prime}\\
\vdots\\
v_{N-1}^{\prime}\\
1
\end{pmatrix}
. \label{eqClmN.Y}%
\end{equation}
Conversely, if $\left(  A,\eta\right)  $ satisfies \textup{(\ref{eqClmN.U}%
)--(\ref{eqClmN.X}),} then \textup{(\ref{eqClmN.T})} defines a unital order
isomorphism $G\rightarrow G^{\prime}$.
\end{lemma}

\begin{proof}
Let $\Lambda\in\operatorname*{GL}\left(  N,\mathbb{Z}\left[  \frac{1}{\lambda
}\right]  \right)  $ be the matrix in Proposition \ref{prop10.7} implementing
the isomorphism. By Lemma \ref{LemClmN.D}, $\Lambda\mathbf{v}=\xi
\mathbf{v}^{\prime}$, and by the proof of Corollary \ref{CorClmN.E},%
\[
\xi=\frac{%
\ip{\alpha}{\mathbf{v}}%
}{%
\ip{\alpha^{\prime}}{\mathbf{v}^{\prime}}%
}=\frac{\lambda^{-\left(  N-1\right)  }I\left(  J\right)  }{\lambda
^{\prime\,-\left(  N-1\right)  }I\left(  J^{\prime}\right)  }=\lambda
^{\prime\,N-1}\lambda^{-\left(  N-1\right)  },
\]
and thus
\begin{equation}
\Lambda\mathbf{v}=\lambda^{\prime\,N-1}\lambda^{-\left(  N-1\right)
}\mathbf{v}^{\prime}. \label{eqClmN.Z}%
\end{equation}
This shows that%
\[
\theta\left(  \xi\mathbf{v}\right)  =\xi\lambda^{\prime\,N-1}\lambda^{-\left(
N-1\right)  }\mathbf{v}^{\prime}%
\]
for all $\xi\in\mathbb{Z}\left[  \frac{1}{\lambda}\right]  $. Furthermore we
must have%
\begin{equation}
\theta|_{\mathbb{Z}^{N-1}}=A\oplus\eta\left(  \,\cdot\,\right)  \mathbf{v}%
^{\prime}, \label{eqClmN.uther}%
\end{equation}
where $A\in M_{N-1}\left(  \mathbb{Z}\right)  $ and $\eta\in
\operatorname*{Hom}\left(  \mathbb{Z}^{N-1},\mathbb{Z}\left[  \frac{1}%
{\lambda}\right]  \right)  $. Now, (\ref{eqClmN.T}) follows from
(\ref{eqClmN.Z}) and (\ref{eqClmN.uther}). For $\theta$ to be onto, $A$ must
be surjective and hence $A\in\operatorname*{GL}\left(  N,\mathbb{Z}\right)  $.
The condition (\ref{eq10.11})(2) in Proposition \ref{prop10.7} is equivalent
to%
\[%
\ip{\beta}{\mathbf{x}}%
+\xi%
\ip{\alpha}{\mathbf{v}}%
=%
\ip{\beta^{\prime}}{A\mathbf{x}}%
+%
\ip{\alpha^{\prime}}{\eta\left( \mathbf{x}\right) \mathbf{v}^{\prime}}%
+\xi\lambda^{\prime\,N-1}\lambda^{-\left(  N-1\right)  }%
\ip{\alpha^{\prime}}{\mathbf{v}^{\prime}}%
\]
for all $\mathbf{x}\oplus\xi\mathbf{v}\in\mathbb{Z}^{N-1}\oplus\mathbb{Z}%
\left[  \frac{1}{\lambda}\right]  \mathbf{v}$ and since $\lambda^{N-1}%
\ip{\alpha}{\mathbf{v}}%
=I\left(  J\right)  =I\left(  J^{\prime}\right)  =\lambda^{\prime\,N-1}%
\ip{\alpha^{\prime}}{\mathbf{v}^{\prime}}%
$, this is (\ref{eqClmN.W}). Finally noting that%
\[%
\begin{pmatrix}
1\\
0\\
\vdots\\
0
\end{pmatrix}
=\mathbf{v}-\mathbf{v}^{\prime},
\]
(\ref{eqClmN.X}) is a transcription of (\ref{eq10.11})(5):%
\begin{align*}
\Lambda\mathbf{v}-\Lambda\mathbf{\bar{v}}  &  =\lambda^{\prime\,N-1}%
\lambda^{-\left(  N-1\right)  }\mathbf{v}^{\prime}-\left(  A\mathbf{\bar{v}%
},\eta\left(  \mathbf{\bar{v}}\right)  \mathbf{v}^{\prime}\right) \\
&  =\left(  -A\mathbf{\bar{v}},\left(  -\eta\left(  \mathbf{\bar{v}}\right)
+\lambda^{\prime\,N-1}\lambda^{-\left(  N-1\right)  }\mathbf{v}^{\prime
}\right)  \right)  .
\end{align*}
But since $\Lambda\left(
\begin{smallmatrix}
1\\
0\\
\raisebox{0pt}[10pt]{$\vdots$}\\
0
\end{smallmatrix}
\right)  =\left(
\begin{smallmatrix}
1\\
0\\
\raisebox{0pt}[10pt]{$\vdots$}\\
0
\end{smallmatrix}
\right)  $, this is equal to
\[
\mathbf{v}^{\prime}-\mathbf{\bar{v}}^{\prime}=\left(  -\mathbf{\bar{v}%
}^{\prime},\mathbf{v}^{\prime}\right)  ,
\]
which is equivalent to (\ref{eqClmN.X}).

For the converse statement, one has to verify that if $\theta$ is defined by
(\ref{eqClmN.T}), then $\theta$ satisfies the conditions in Proposition
\ref{prop10.7}, but this follows by the same computations as above. (Note that
as $\theta\left(  G\right)  =G^{\prime}$, the conditions (3) and (4) in
Proposition \ref{prop10.7} are automatic. To show $\theta\left(  G\right)
=G^{\prime}$, note first that $\theta\left(  \mathbb{Z}\left[  \frac
{1}{\lambda}\right]  \mathbf{v}\right)  =\mathbb{Z}\left[  \frac{1}{\lambda
}\right]  \mathbf{v}^{\prime}$, and next, since $A$ is onto, there is for any
$\mathbf{m}\in\mathbb{Z}^{N-1}$ an $\mathbf{n}\in\mathbb{Z}^{N-1}$ with
$A\mathbf{n}=\mathbf{m}$, but then $\theta\left(  \mathbf{n}\oplus\left(
-\eta\left(  \mathbf{m}\right)  \mathbf{v}\right)  \right)  =A\mathbf{v}%
\oplus0=\mathbf{m}\oplus0$, thus $\theta$ is surjective. It is clearly injective.)
\end{proof}

Note that Lemma \ref{LemClmN.H} hints at a method of constructing unital order
isomorphisms $\left(  A,\eta\right)  $. First find an $A\in\operatorname*{GL}%
\left(  N-1,\mathbb{Z}\right)  $ satisfying $A\mathbf{\bar{v}}=\mathbf{\bar
{v}}^{\prime}$, and then solve (\ref{eqClmN.W}) for $\eta\left(
\mathbf{x}\right)  $. However, one then has to check (\ref{eqClmN.V}) and the
remaining condition in (\ref{eqClmN.X}), and these conditions are very
restrictive. This is illustrated by the following lemma.

\begin{lemma}
\label{LemClmN.I}Adopt the notation and general assumptions in Lemma
\textup{\ref{LemClmN.H}.} If $\left(  A,\eta\right)  $ is a solution of the
conditions \textup{(\ref{eqClmN.U})--(\ref{eqClmN.X}),} then%
\begin{equation}
\eta\left(  \mathbb{Z}^{N-1}\right)  \subseteq\left(  \frac{\gcd\left(
\lambda,\lambda^{\prime}\right)  }{\lambda}\right)  ^{N-1}\mathbb{Z},
\label{eqClmN.okseskaft}%
\end{equation}
and thus one may without loss of generality replace \textup{(\ref{eqClmN.V})}
by%
\begin{equation}
\eta\in\operatorname*{Hom}\left(  \mathbb{Z}^{N-1},\left(  \frac{\gcd\left(
\lambda,\lambda^{\prime}\right)  }{\lambda}\right)  ^{N-1}\mathbb{Z}\right)  .
\label{eqClmN.angstrom}%
\end{equation}
In particular, if $\lambda=\lambda^{\prime}$, then%
\begin{equation}
\eta\in\operatorname*{Hom}\left(  \mathbb{Z}^{N-1},\mathbb{Z}\right)  .
\label{eqClmN.1}%
\end{equation}
\end{lemma}

\begin{proof}
Note first that by the beginning of the proof of Lemma \ref{LemClmN.F}, $%
\ip{\alpha^{\prime}}{\mathbf{v}^{\prime}}%
$ has the form\label{LOSalphaeigenvector_5}\label{LOSbeta_2}%
\begin{equation}%
\ip{\alpha^{\prime}}{\mathbf{v}^{\prime}}%
=\frac{\left(  a\lambda+1\right)  }{\lambda^{\prime\,N-1\mathstrut}}%
=\frac{I\left(  J^{\prime}\right)  }{\lambda^{\prime\,N-1\mathstrut}},
\label{eqClmN.2}%
\end{equation}
where $a$ is a positive integer. But it follows from (\ref{eqClmN.W}) that%
\begin{align*}
\eta\left(  \mathbf{x}\right)   &  =\frac{\left(
\ip{\beta}{\mathbf{x}}%
-%
\ip{\beta^{\prime}}{A\mathbf{x}}%
\right)  }{%
\ip{\alpha^{\prime}}{\mathbf{v}^{\prime}}%
}\\
&  =\frac{\lambda^{\prime\,N-1}\left(
\ip{\beta}{\mathbf{x}}%
-%
\ip{\beta^{\prime}}{A\mathbf{x}}%
\right)  }{I\left(  J^{\prime}\right)  }\\
&  \subseteq\frac{\lambda^{\prime\,N-1}\left(  \frac{1}{\lambda^{N-1\mathstrut
}}\mathbb{Z}+\frac{1}{\lambda^{\prime\,N-1\mathstrut}}\mathbb{Z}\right)
}{I\left(  J^{\prime}\right)  }\\
&  =\frac{\left(  \left(  \frac{\lambda^{\prime}}{\lambda}\right)
^{N-1}\mathbb{Z}+\mathbb{Z}\right)  }{I\left(  J^{\prime}\right)  }\\
&  =\frac{\gcd\left(  \lambda,\lambda^{\prime}\right)  ^{N-1}}{\lambda
^{N-1\mathstrut}I\left(  J\right)  }\mathbb{Z}.
\end{align*}
But on the other hand%
\[
\eta\left(  \mathbf{x}\right)  \in\mathbb{Z}\left[  \tfrac{1}{\lambda
}\right]
\]
by (\ref{eqClmN.V}), and since $I\left(  J\right)  =\left(  a\lambda^{\prime
}+1\right)  $ is mutually prime with $\lambda$, (\ref{eqClmN.okseskaft})
follows. The remaining statements in Lemma \ref{LemClmN.I} are obvious.
\end{proof}

Note that among the solutions of the relations in Lemma \ref{LemClmN.H} one
can always look for one of them with $\eta=0$. Because of (\ref{eqClmN.X})
such solutions only exist if $\lambda=\lambda^{\prime}$, and then $A$ must
satisfy%
\begin{equation}%
\begin{cases}
A\in\operatorname*{GL}\left( N-1,\mathbb{Z}\right)   \\
\diracb{\beta}A=\diracb{\beta}  \\
A\mathbf{\bar{v}}=\mathbf{\bar{v}}^{\prime}.
\end{cases}
\label{eqClmN.3}%
\end{equation}
Thus if $\lambda=\lambda^{\prime}$, the existence of an $A$ with the
properties in (\ref{eqClmN.3}) is a sufficient condition for isomorphism. We
next formulate a rather complicated condition which is both sufficient and
necessary for general $\lambda$, $\lambda^{\prime}$.

\begin{lemma}
\label{LemClmN.J}Adopt the notation and general assumptions in Lemma
\textup{\ref{LemClmN.H}.} Then a necessary and sufficient condition that $G$
and $G^{\prime}$ are unital order isomorphic is that there exists an integer
$\left(  N-1\right)  \times\left(  N-1\right)  $ matrix $A=\left[
a_{ij}\right]  _{i,j=1}^{N}$ with the properties%
\begin{align}
A  &  \in\operatorname*{GL}\left(  N-1,\mathbb{Z}\right)  ,\label{eqClmN.JU}\\
\lambda^{\prime\,N-1}\lambda^{1-i}-\sum_{j=2}^{N}a_{ji}\lambda^{\prime\,N-j}
&  \in I\left(  J\right)  \mathbb{Z}\left[  \tfrac{1}{\lambda}\right]
\text{\qquad for }i=2,3,\dots,N,\label{eqClmN.JW}\\
A\mathbf{\bar{v}}  &  =\mathbf{\bar{v}}^{\prime}. \label{eqClmN.JX}%
\end{align}
Thus a necessary condition for isomorphism is that
\begin{equation}
\lambda^{\prime\,N-1}\lambda^{1-N}\in I\left(  J\right)  \mathbb{Z}\left[
\tfrac{1}{\lambda}\right]  +\mathbb{Z}. \label{eqClmN.JY}%
\end{equation}
\end{lemma}

\begin{remark}
\label{RemClmN.J1}Note that since all the terms on the left side of
\textup{(\ref{eqClmN.JW}) are integers except }$\lambda^{1-i}$, and $I\left(
J\right)  $ is mutually prime with $\lambda$ and $\lambda^{\prime}$, and
$\mathbb{Z}\left[  \frac{1}{\lambda}\right]  $ is closed under division by
$\lambda$, the condition \textup{(\ref{eqClmN.JW})} can be formulated in the
following more user-friendly way:%
\begin{equation}
\lambda^{i-1}\sum_{j=2}^{N}a_{ji}\lambda^{\prime\,N-j}\in\lambda^{\prime
\,N-1}+I\left(  J\right)  \mathbb{Z} \label{eqClmN.JW1}%
\end{equation}
for $i=2,3,\dots,N$.
\end{remark}

\begin{proof}
We know that the conditions (\ref{eqClmN.U})--(\ref{eqClmN.X}) in Lemma
\ref{LemClmN.H} are necessary and sufficient. But with $A$ given, one may use
(\ref{eqClmN.W}) to define $\eta$,%
\begin{align*}
\eta\left(  \mathbf{x}\right)   &  =\frac{%
\ip{\beta}{\mathbf{x}}%
-%
\ip{\beta^{\prime}}{A\mathbf{x}}%
}{%
\ip{\alpha^{\prime}}{\mathbf{v}^{\prime}}%
}\\
&  =\frac{\lambda^{\prime\,N-1}\left(
\ip{\beta}{\mathbf{x}}%
-%
\ip{\beta^{\prime}}{A\mathbf{x}}%
\right)  }{I\left(  J^{\prime}\right)  }%
\end{align*}
and since $I\left(  J^{\prime}\right)  =I\left(  J\right)  $ is relatively
prime to both $\lambda$ and $\lambda^{\prime}$, we see that (\ref{eqClmN.JW})
is necessary and sufficient for (\ref{eqClmN.V}) by putting $\mathbf{x}%
=\mathbf{e}_{2},\mathbf{e}_{3},\dots,\mathbf{e}_{N}$. Finally if $\eta$ is
defined as above we have%
\begin{align*}
\eta\left(  \mathbf{\bar{v}}\right)   &  =\frac{\lambda^{\prime\,N-1}%
\ip{\beta}{\mathbf{\bar{v}}}%
-\lambda^{\prime\,N-1}%
\ip{\beta^{\prime}}{\mathbf{\bar{v}}^{\prime}}%
}{I\left(  J\right)  }\\
&  =\frac{\lambda^{\prime\,N-1}\left(
\ip{\alpha}{\mathbf{v}}%
-1\right)  -\lambda^{\prime\,N-1}\left(
\ip{\alpha^{\prime}}{\mathbf{v}^{\prime}}%
-1\right)  }{I\left(  J\right)  }\\
&  =\frac{\lambda^{\prime\,N-1}\lambda^{-\left(  N-1\right)  }I\left(
J\right)  -I\left(  J\right)  }{I\left(  J\right)  }\\
&  =\lambda^{\prime\,N-1}\lambda^{-\left(  N-1\right)  }-1,
\end{align*}
so the last condition in (\ref{eqClmN.X}) is fulfilled.

Finally, (\ref{eqClmN.JY}) follows by putting $i=N$ in (\ref{eqClmN.JW}).
\end{proof}

We will now apply this theory to more specific examples.

\section{\label{ClmN.1}The case $N=1$}

Here one has $\lambda=m_{N}$ automatically, and the corresponding $C^{\ast}%
$-algebra $\mathfrak{A}_{L}$ is the UHF-algebra of Glimm type $\lambda
^{\infty}$ \cite{Gli60}. Thus the algebras corresponding to $\lambda$ and
$\lambda^{\prime}$ are isomorphic if and only if Prim$\left(  \lambda\right)
=\operatorname*{Prim}\left(  \lambda^{\prime}\right)  $. (Note that $I\left(
J\right)  =1$ in all these cases, so this invariant does not separate
isomorphism classes.)

\section{\label{ClmN.2}The case $N=2$}

Here it follows from (\ref{eqApp.Cla2}) that the possible $J$'s with $N=2$ and
$\lambda=m_{N}$ are%
\begin{equation}
J=%
\begin{pmatrix}
\lambda-1 & 1\\
\lambda & 0
\end{pmatrix}
\label{eqClmN.4}%
\end{equation}
for $\lambda=2,3,\dots$. By Lemma \ref{LemClmN.B}, $\mathbf{v}=\left(
\begin{smallmatrix}
1\\
1
\end{smallmatrix}
\right)  $, and by (\ref{eqClmN.L}),%
\begin{equation}
I\left(  J\right)  =\lambda+1. \label{eqClmN.5}%
\end{equation}
It follows from Corollary \ref{CorClmN.E1} that \emph{all} the algebras
corresponding to (\ref{eqClmN.4}) for $\lambda=2,3,4,5,\dots$ are pairwise non-isomorphic.

\section{\label{ClmN.3}The case $N=3$}

Using Lemma \ref{LemClmN.B} one observes that if $\lambda\in\left\{
2,3,\dots\right\}  $ is given and $\lambda=m_{3}$, then $J$ has the form%
\begin{equation}
J=%
\begin{pmatrix}
\lambda-v_{2} & 1 & 0\\
\lambda v_{2}-1 & 0 & 1\\
\lambda & 0 & 0
\end{pmatrix}
\text{\quad with\quad}\mathbf{v}=%
\begin{pmatrix}
1\\
v_{2}\\
1
\end{pmatrix}
\label{eqClmN.6}%
\end{equation}
where $1\leq v_{2}\leq\lambda$. Using (\ref{eqClmN.L}) one computes%
\[
I\left(  J\right)  =\lambda^{2}+v_{2}\lambda+1.
\]
Hence an immediate corollary of Corollary \ref{CorClmN.E1} is:

\begin{proposition}
\label{ProClmN.1}If $J$, $J^{\prime}$ are matrices of the form
\textup{(\ref{eqClmN.6})} with $\lambda=\lambda^{\prime}$ and $N=3$, then
$\mathfrak{A}_{J}$ and $\mathfrak{A}_{J^{\prime}}$ are isomorphic if and only
if $J=J^{\prime}$. The same is true if one $\lambda$ is an integer multiple of
the other.
\end{proposition}

\begin{proof}
The last statement follows from Lemma \ref{LemClmN.F}, (\ref{eqClmN.M2}).
\end{proof}

Let us now look at the example%
\[
\lambda=48=2^{4}\cdot3,\qquad\lambda^{\prime}=54=2\cdot3^{3}.
\]
Then one computes that $I\left(  J\right)  =I\left(  J^{\prime}\right)  $ for
exactly the following four pairs:%
\[%
\begin{tabular}
[c]{r|rrrr}%
$v_{2}$ & $15$ & $24$ & $33$ & $42$\\\hline
$v_{2}^{\prime}$ & $2$ & $10$ & $18$ & $26$%
\end{tabular}
\]

In general the solutions of (\ref{eqClmN.JW1}) together with (\ref{eqClmN.JU})
and (\ref{eqClmN.JX}) in Lemma \ref{LemClmN.J} may be found as follows.

Let $A=%
\begin{pmatrix}
a & b\\
c & d
\end{pmatrix}
^{\mathstrut}=%
\begin{pmatrix}
a_{22} & a_{23}\\
a_{32} & a_{33}%
\end{pmatrix}
$. Then (\ref{eqClmN.JX}) is fulfilled if and only if%
\begin{align*}
b  &  =v_{2}^{\prime}-av_{2},\\
d  &  =1-cv_{2},
\end{align*}
and then (\ref{eqClmN.JU}) holds if and only if%
\[
\det A=a-cv_{2}^{\prime}=\pm1,
\]
and it follows that the conjunction of (\ref{eqClmN.JU}) and (\ref{eqClmN.JX})
is valid if and only if%
\begin{equation}
A=%
\begin{pmatrix}
cv_{2}^{\prime}\pm1 & v_{2}^{\prime}-cv_{2}v_{2}^{\prime}\mp v_{2}\\
c & 1-cv_{2}%
\end{pmatrix}
\label{eqClmN.3A}%
\end{equation}
for an integer $c$. To determine this integer, we use (\ref{eqClmN.JW1}) to
deduce (for $i=2,3$):%
\begin{align}
c\lambda\left(  v_{2}^{\prime}\lambda^{\prime}+1\right)   &  \in\left(
\mp\lambda+\lambda^{\prime}\right)  \lambda^{\prime}+I\left(  J\right)
\mathbb{Z},\label{eqClmN.3B}\\
cv_{2}\lambda^{2}\left(  v_{2}^{\prime}\lambda^{\prime}+1\right)   &
\in\left(  v_{2}^{\prime}\mp v_{2}\right)  \lambda^{2}\lambda^{\prime}%
+\lambda^{2}-\lambda^{\prime\,2}+I\left(  J\right)  \mathbb{Z},
\label{eqClmN.3C}%
\end{align}
where%
\begin{equation}
I\left(  J\right)  =\lambda^{2}+v_{2}\lambda+1=\lambda^{\prime\,2}%
+v_{2}^{\prime}\lambda^{\prime}+1=I\left(  J^{\prime}\right)  .
\label{eqClmN.3D}%
\end{equation}
Now, let us note that (\ref{eqClmN.3C}) actually follows from (\ref{eqClmN.3B}%
), with the same $c$. This is because (\ref{eqClmN.3B}) implies that
\begin{align}
c\lambda\left(  v_{2}\lambda\right)  \left(  v_{2}^{\prime}\lambda^{\prime
}+1\right)   &  \subseteq v_{2}\lambda\left(  \mp\lambda+\lambda^{\prime
}\right)  \lambda^{\prime}+I\left(  J\right)  v_{2}\lambda\mathbb{Z}%
\label{eqClmN.3E}\\
&  \subseteq v_{2}\lambda\left(  \mp\lambda+\lambda^{\prime}\right)
\lambda^{\prime}+I\left(  J\right)  \mathbb{Z}.\nonumber
\end{align}
But note that%
\begin{align*}
D  &  :=v_{2}\lambda\left(  \mp\lambda+\lambda^{\prime}\right)  \lambda
^{\prime}-\left(  \left(  v_{2}^{\prime}\mp v_{2}\right)  \lambda^{2}%
\lambda^{\prime}+\lambda^{2}-\lambda^{\prime\,2}\right) \\
&  =v_{2}\lambda\lambda^{\prime\,2}-v_{2}^{\prime}\lambda^{2}\lambda^{\prime
}-\lambda^{2}+\lambda^{\prime\,2}.
\end{align*}
As $I\left(  J\right)  =\lambda^{2}+v_{2}\lambda+1=\lambda^{\prime\,2}%
+v_{2}^{\prime}\lambda^{\prime}+1=I\left(  J^{\prime}\right)  $ we have that%
\[
\lambda^{2}-\lambda^{\prime\,2}=v_{2}^{\prime}\lambda^{\prime}-v_{2}\lambda,
\]
and hence%
\begin{align*}
D  &  =v_{2}\lambda\lambda^{\prime\,2}-v_{2}^{\prime}\lambda^{2}%
\lambda^{\prime}-v_{2}^{\prime}\lambda^{\prime}+v_{2}\lambda\\
&  =v_{2}\lambda\left(  \lambda^{\prime\,2}+1\right)  -v_{2}^{\prime}%
\lambda^{\prime}\left(  \lambda^{2}+I\right) \\
&  =v_{2}\lambda I\left(  J^{\prime}\right)  -v_{2}\lambda v_{2}^{\prime
}\lambda^{\prime}-v_{2}^{\prime}\lambda^{\prime}I\left(  J\right)
+v_{2}^{\prime}\lambda^{\prime}v_{2}\lambda\\
&  \in I\left(  J\right)  \mathbb{Z}.
\end{align*}
Thus (\ref{eqClmN.3C}) follows from (\ref{eqClmN.3B}), i.e., the two
conditions (\ref{eqClmN.3B}) and (\ref{eqClmN.3C}) are equivalent to
(\ref{eqClmN.3B}) alone. But now observe that%
\[
v_{2}^{\prime}\lambda^{\prime}+1=I\left(  J^{\prime}\right)  -\lambda
^{\prime\,2},
\]
and hence (\ref{eqClmN.3B}) is equivalent to%
\[
-c\lambda\lambda^{\prime\,2}\in\left(  \mp\lambda+\lambda^{\prime}\right)
\lambda^{\prime}+I\left(  J\right)  \mathbb{Z}.
\]
Since $\lambda\lambda^{\prime\,2}$ is relatively prime to $I\left(  J\right)
$, it follows from the Euclidean algorithm within $\mathbb{Z}$ that this
relation always has a solution $c$ in $\mathbb{Z}$! Thus we have proved:

\begin{theorem}
\label{ThmClmN.2}If $J$, $J^{\prime}$ are matrices of the form
\textup{(\ref{eqClmN.6}),} the following conditions are equivalent.

\begin{enumerate}
\item \label{ThmClmN.2(1)}$\mathfrak{A}_{J}$ and $\mathfrak{A}_{J^{\prime}}$
are isomorphic.

\item \label{ThmClmN.2(2)}$\operatorname*{Prim}\left(  \lambda\right)
=\operatorname*{Prim}\left(  \lambda^{\prime}\right)  $ and $I\left(
J\right)  =\lambda^{2}+v_{2}\lambda+1=\lambda^{\prime\,2}+v_{2}^{\prime
}\lambda^{\prime}+1=I\left(  J^{\prime}\right)  $.
\end{enumerate}
\end{theorem}

\begin{proof}
The necessity of the conditions was established in Theorem \ref{CorCyc.9} and
Corollary \ref{CorClmN.E1}. The sufficiency follows from Lemma \ref{LemClmN.J}
as argued before the theorem.
\end{proof}

\begin{remark}
\label{RemClmN.3}It follows from the argument that the $A\in\operatorname*{GL}%
\left(  2,\mathbb{Z}\right)  $ and the $\eta\in\operatorname*{Hom}\left(
\mathbb{Z}^{N-1},\mathbb{Z}\left[  \frac{1}{\lambda}\right]  \right)  $
corresponding to the isomorphism can and must be taken to be%
\[
A=%
\begin{pmatrix}
cv_{2}^{\prime}\pm1 & v_{2}^{\prime}-cv_{2}v_{2}^{\prime}\mp v_{2}\\
c & 1-cv_{2}%
\end{pmatrix}
,
\]
where $c$ is any solution of
\[
c\lambda\lambda^{\prime\,2}\in\left(  \pm\lambda-\lambda^{\prime}\right)
\lambda^{\prime}+I\left(  J\right)  \mathbb{Z}%
\]
and $\eta$ is then defined by \textup{(\ref{eqClmN.W}).}
\end{remark}

In particular the four pairs mentioned after Proposition \ref{ProClmN.1}
define isomorphic algebras. For example, the algebra defined by%
\[
J=%
\begin{pmatrix}
33 & 1 & 0\\
719 & 0 & 1\\
48 & 0 & 0
\end{pmatrix}
\text{\qquad with }\lambda=48\text{ and }\mathbf{v}=%
\begin{pmatrix}
1\\
15\\
1
\end{pmatrix}
\]
is isomorphic to%
\[
J^{\prime}=%
\begin{pmatrix}
52 & 1 & 0\\
107 & 0 & 1\\
54 & 0 & 0
\end{pmatrix}
\text{\qquad with }\lambda^{\prime}=54\text{ and }\mathbf{v}^{\prime}=%
\begin{pmatrix}
1\\
2\\
1
\end{pmatrix}
.
\]
Here $I\left(  J\right)  =I\left(  J^{\prime}\right)  =3025$, and the
Euclidean algorithm gives%
\[
c=-313308=1292\pmod{3025}%
\]
and thus one $A$ that can be used is%
\[
A=%
\begin{pmatrix}
2585 & -38773\\
1292 & -19379
\end{pmatrix}
.
\]

\section{\label{ClmN.4}The case $N\geq4$}

In this section we prove that Theorem \ref{ThmClmN.2} remains valid also in
the general $\lambda=m_{N}$ case, i.e., $\left\{  N,\operatorname*{Prim}%
\lambda,I\left(  J\right)  \right\}  $ is a complete invariant. If $A$ is a
matrix in $M_{N-1}\left(  \mathbb{Z}\right)  $ of the form%
\begin{equation}
A=%
\begin{pmatrix}
\mathbf{a}_{2} & \mathbf{a}_{3} & \cdots & \mathbf{a}_{N}%
\end{pmatrix}
\label{eqClmN.Abis}%
\end{equation}
in terms of column vectors $\mathbf{a}_{2},\dots,\mathbf{a}_{N}$, then an easy
computation shows that $A\mathbf{v}=\mathbf{v}^{\prime}$ as in
(\ref{eqClmN.JX}) if and only if%
\begin{equation}
\mathbf{a}_{N}=\mathbf{v}^{\prime}-v_{2}\mathbf{a}_{2}-v_{3}\mathbf{a}%
_{3}-\dots-v_{N-1}\mathbf{a}_{N-1}. \label{eqClmN.Bbis}%
\end{equation}
Thus condition (\ref{eqClmN.JU}) says that%
\begin{equation}
\det%
\begin{pmatrix}
\mathbf{a}_{2} & \mathbf{a}_{3} & \cdots & \mathbf{a}_{N-1} & \mathbf{v}%
^{\prime}%
\end{pmatrix}
=\pm1, \label{eqClmN.Cbis}%
\end{equation}
and what remains is to choose $\mathbf{a}_{2},\mathbf{a}_{3},\dots
,\mathbf{a}_{N-1}$ such that (\ref{eqClmN.JW1}) is satisfied. For example, try%
\begin{equation}
A=%
\begin{pmatrix}
a_{2} & 0 & \cdots & 0 & v_{2}^{\prime}-a_{2}v_{2}\\
0 & a_{3} &  & 0 & v_{3}^{\prime}-a_{3}v_{3}\\
\vdots &  & \ddots &  & \vdots\\
0 & 0 &  & a_{N-1} & v_{N-1}^{\prime}-a_{N-1}v_{N-1}\\
c_{2} & c_{3} & \cdots &  c_{N-1} & 1-c_{2}v_{2}-c_{3}v_{3}-\dots
-c_{N-1}v_{N-1}%
\end{pmatrix}
. \label{eqClmN.Dbis}%
\end{equation}
Fixing $c_{2},c_{3},\dots,c_{N-1}$, one may now choose $a_{N-1},a_{N-2}%
,\dots,a_{2}$ successively such that the determinant of the lower right
$k\times k$ matrix is $1$ for $k=2,3,\dots,N-1$. This leads to the recursion
relations%
\begin{equation}%
\begin{aligned}
a_{N-1}  & =1+c_{N-1}v_{N-1}^{\prime},\\
a_{N-2}  & =1+c_{N-2}a_{N-1}v_{N-2}^{\prime},\\
a_{N-3}  & =1+c_{N-3}a_{N-2}a_{N-1}v_{N-3}^{\prime},\\
\vdots& \qquad\vdots\\
a_{2}  & =1+c_{2}a_{3}a_{4}\cdots a_{N-1}v_{2}^{\prime}.
\end{aligned}%
\label{eqClmN.Ebis}%
\end{equation}
Inserting this in (\ref{eqClmN.JW1}) gives a way of determining suitable
values of $c_{2},\dots,c_{N-1}$ as in the $N=3$ case, but we have by now
already made several arbitrary choices which we did not need to do in the
$N=3$ case, and it is not quite clear that this approach leads to the goal.
There is, however, one case where it does, namely if $\lambda=\lambda^{\prime
}$. Then one may simply choose $c_{2}=c_{3}=\dots=c_{N-1}=0$, and so
$a_{2}=a_{3}=\dots=a_{N-1}=1$ and%
\begin{equation}
A=%
\begin{pmatrix}
1 & 0 & \cdots & 0 & v_{2}^{\prime}-v_{2}\\
0 & 1 &  & 0 & v_{3}^{\prime}-v_{3}\\
\vdots &  & \ddots &  & \vdots\\
0 & 0 &  & 1 & v_{N-1}^{\prime}-v_{N-1}\\
0 & 0 & \cdots & 0 & 1
\end{pmatrix}
. \label{eqClmN.Fbis}%
\end{equation}
Since $\lambda=\lambda^{\prime}$, (\ref{eqClmN.JW1}) is automatically
satisfied, and we have

\begin{theorem}
\label{ThmClmN.3}Let $J$, $J^{\prime}$ be matrices of the form
\textup{(\ref{eqCycNew.1}),} and assume that $\lambda=m_{N}=\lambda^{\prime
}=m_{N}^{\prime}$. Then the following three conditions are equivalent:

\begin{enumerate}
\item \label{ThmClmN.3(1)}$\mathfrak{A}_{J}$ and $\mathfrak{A}_{J^{\prime}}$
are isomorphic;

\item \label{ThmClmN.3(2)}$N=N^{\prime}$ and $I\left(  J\right)  =I\left(
J^{\prime}\right)  $;
\end{enumerate}

\noindent and

\begin{enumerate}
\addtocounter{enumi}{2}

\item \label{ThmClmN.3(3)}$N=N^{\prime}$ and $%
\ip{\alpha}{\mathbf{v}}%
=%
\ip{\alpha^{\prime}}{\mathbf{v}^{\prime}}%
$.
\end{enumerate}
\end{theorem}

\begin{proof}
This follows from the remarks before the theorem, Lemma \ref{LemClmN.J},
Corollary \ref{CorClmN.E1}, and the formula%
\begin{equation}%
\ip{\alpha}{\mathbf{v}}%
=I\left(  J\right)  /\lambda^{N-1}.%
\settowidth{\qedskip}{$\displaystyle
\ip{\alpha}{\mathbf{v}}=I\left( J\right) /\lambda^{N-1}.$}\addtolength
{\qedskip}{-\textwidth}\rlap{\makebox[-0.5\qedskip][r]{\qedsymbol}}%
\label{eqClmN.Gbis}%
\end{equation}%
\renewcommand{\qed}{}%
\end{proof}

\begin{figure}[ptb]
\begin{picture}(360,514)
\put(0,0){\includegraphics
[bb=243 62 461 704,clip,width=174.4bp,height=513.6bp]{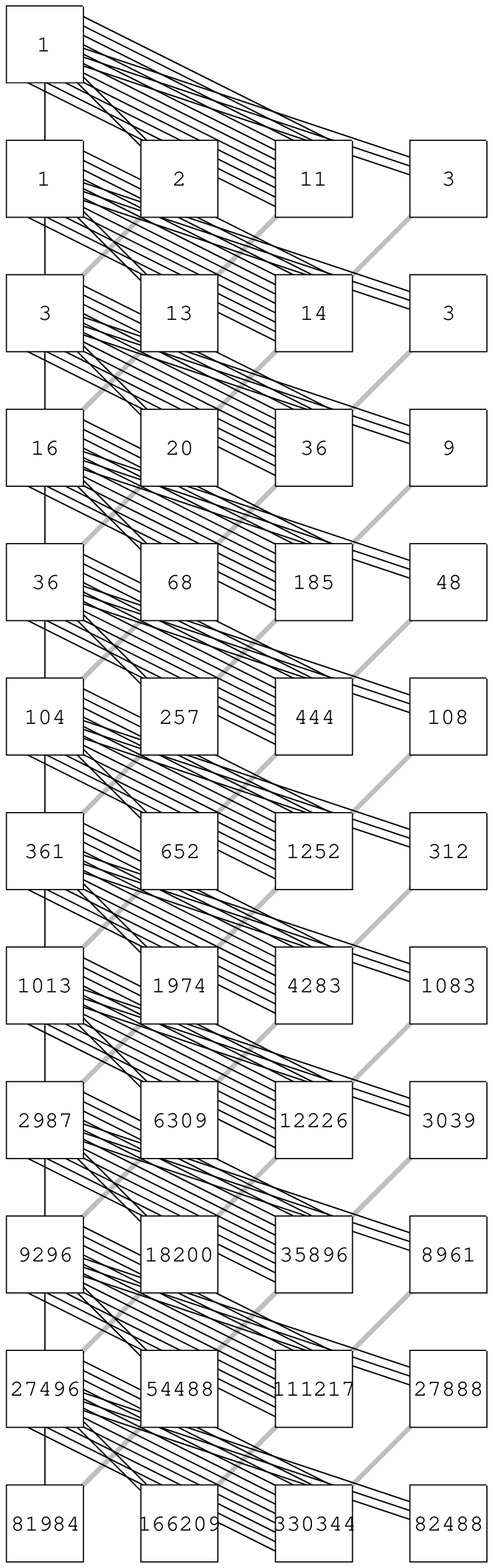}}
\put(185.6,0){\includegraphics
[bb=243 62 461 704,clip,width=174.4bp,height=513.6bp]{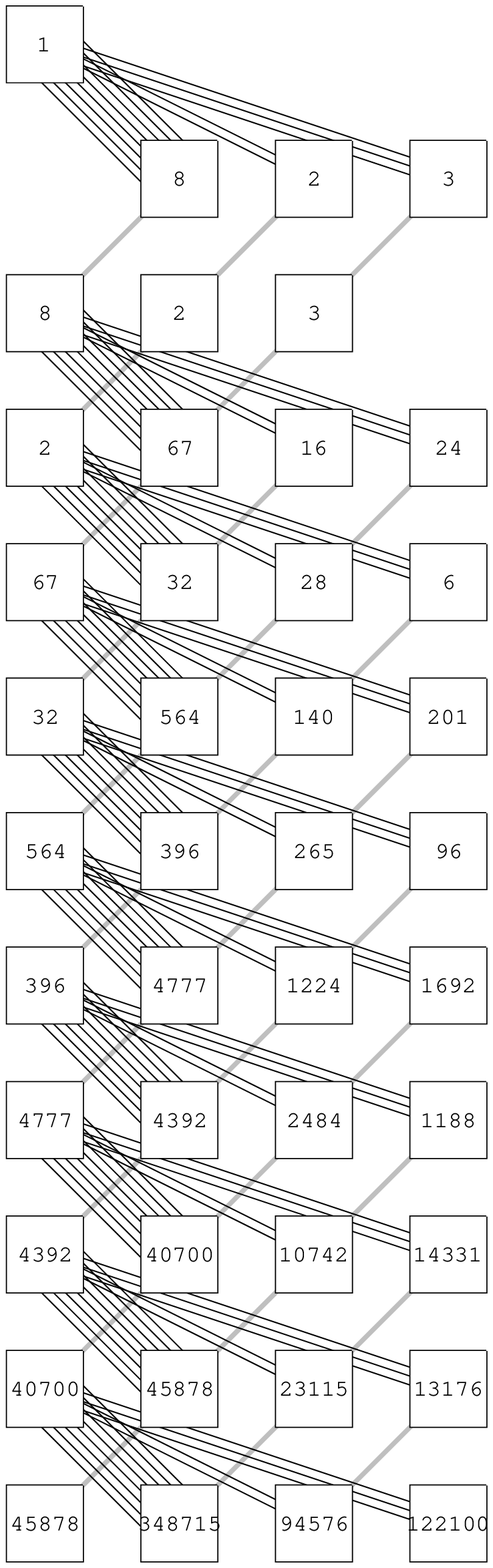}}
\end{picture}\caption{$L=\{1,2,2,3,3,3,3,3,3,3,3,3,3,3,4,4,4\}$, first column
$=(1\;2\;11\;3)^{\mathrm{t}}$ (left); $L=\{2,2,2,2,2,2,2,2,3,3,4,4,4\}$, first
column $=(0\;8\;2\;3)^{\mathrm{t}}$ (right). See (\ref{eqClmNJun.55}%
)--(\ref{eqClmNJun.61}). These diagrams represent isomorphic AF-algebras.}%
\label{BratDiagN4samlam}%
\end{figure}

\begin{example}
\label{ExaClmN}The present example shows that Proposition
\textup{\ref{ProClmN.1}} does not extend to the case $N=4$, i.e., isomorphism
does not imply $\lambda\neq\lambda^{\prime}$, and it also shows that Theorem
\textup{\ref{ThmClmN.3}} is not merely concerned with the empty set.

The example is%
\begin{equation}
J=%
\begin{pmatrix}
1 & 1 & 0 & 0\\
2 & 0 & 1 & 0\\
11 & 0 & 0 & 1\\
3 & 0 & 0 & 0
\end{pmatrix}
,\qquad J^{\prime}=%
\begin{pmatrix}
0 & 1 & 0 & 0\\
8 & 0 & 1 & 0\\
2 & 0 & 0 & 1\\
3 & 0 & 0 & 0
\end{pmatrix}
. \label{eqClmNJun.55}%
\end{equation}
Here $\lambda=\lambda^{\prime}=3$, and%
\begin{equation}
\mathbf{v}=%
\begin{pmatrix}
1\\
2\\
4\\
1
\end{pmatrix}
,\qquad\mathbf{v}^{\prime}=%
\begin{pmatrix}
1\\
3\\
1\\
1
\end{pmatrix}
. \label{eqClmNJun.56}%
\end{equation}
See Figure \textup{\ref{BratDiagN4samlam}.} Using the unique decomposition%
\begin{equation}%
\begin{pmatrix}
1\\
0\\
0\\
0
\end{pmatrix}
=\frac{\mathbf{v}}{%
\ip{\alpha}{\mathbf{v}}%
}+\mathbf{w}, \label{eqClmNJun.57}%
\end{equation}
where $%
\ip{\alpha}{\mathbf{w}}%
=0$, one sees that the $n$'th row vector in the left diagram in Figure
\textup{\ref{BratDiagN4samlam}} behaves asymptotically like%
\begin{equation}
3^{n}\frac{\mathbf{v}^{\mathrm{t}}}{%
\ip{\alpha}{\mathbf{v}}%
}=\frac{27}{58}\cdot3^{n}\cdot\left(  1,2,4,1\right)  , \label{eqClmNJun.58}%
\end{equation}
and the $n$'th row vector in the right diagram in Figure
\textup{\ref{BratDiagN4samlam}} behaves asymptotically like%
\begin{equation}
\frac{27}{58}\cdot3^{n}\cdot\left(  1,3,1,1\right)  . \label{eqClmNJun.59}%
\end{equation}
\textup{(}The latter matrix has an eigenvalue $-2.769\dots$, which is negative
and close to $3$ in absolute value, therefore the slow and oscillatory
convergence to the asymptotic behaviour indicated by the figure.\textup{)} One
can now check that the conditions in Lemma \textup{\ref{LemClmN.H}} are
fulfilled with%
\begin{equation}
A=%
\begin{pmatrix}
1 & 0 & 1\\
0 & 1 & -3\\
0 & 0 & 1
\end{pmatrix}
,\qquad\eta=0. \label{eqClmNJun.60}%
\end{equation}
In fact, since $\lambda=\lambda^{\prime}=3$ and $\eta=0$ it suffices to check
\textup{(\ref{eqClmN.3}),} and that is straightforward. Thus the two
AF-algebras $\mathfrak{A}_{J}$ and $\mathfrak{A}_{J^{\prime}}$ are isomorphic,
showing that Proposition \textup{(\ref{ProClmN.1})} does not extend to $N=4$.
Computing the matrix $\Lambda$ in Proposition \textup{\ref{prop10.7}} for this
example, one finds\label{LOSLambda_4}%
\begin{equation}
\Lambda=%
\begin{pmatrix}
1 & 0 & 0 & 0\\
0 & 1 & 0 & 1\\
0 & 0 & 1 & -3\\
0 & 0 & 0 & 1
\end{pmatrix}
. \label{eqClmNJun.61}%
\end{equation}
\end{example}

We are now ready to state and give the proof of the main result in this chapter.

\begin{theorem}
\label{ThmClmNJun.18}\textup{(K.H. Kim and F. Roush)} Let $J$, $J^{\prime}$ be
matrices of the form \textup{(\ref{eqCycNew.1})} satisfying the standard
requirements there, and assume that $\lambda=m_{N}$ and $\lambda^{\prime
}=m_{N^{\prime}}^{\prime}$. Then the following conditions are equivalent.

\begin{enumerate}
\item \label{ThmClmNJun.18(1)}$\mathfrak{A}_{J}$ and $\mathfrak{A}_{J^{\prime
}}$ are isomorphic.

\item \label{ThmClmNJun.18(2)}$N=N^{\prime}$, $\operatorname*{Prim}\left(
\lambda\right)  =\operatorname*{Prim}\left(  \lambda^{\prime}\right)  $, and
$I\left(  J\right)  =I\left(  J^{\prime}\right)  $.
\end{enumerate}
\end{theorem}

We will establish later, in Corollary \ref{CorClmNJun.21}, that these
conditions are also equivalent to stable isomorphism of $\mathfrak{A}_{J}$ and
$\mathfrak{A}_{J^{\prime}}$.

In order to prove this theorem, we will need some elementary facts about%
\begin{equation}
\mathcal{S}_{n}=\left\{  \mathbf{x}\in\mathbb{Z}^{n}\setminus\left\{
0\right\}  \mid\gcd\mathbf{x}=1\right\}  , \label{eqClmNJun.62}%
\end{equation}
where $n$ is a natural number, and $\gcd\mathbf{x}$ is the greatest common
divisor of the components of $\mathbf{x}$. If $\mathbf{x}\in\mathbb{Z}^{n}$ we
write $%
\diracb{\mathbf{x}}%
$ if we think about $\mathbf{x}$ as a row vector and $%
\dirack{\mathbf{x}}%
$ if we consider $\mathbf{x}$ as the column vector which is the transpose of $%
\diracb{\mathbf{x}}%
$. Thus, by the Euclidean algorithm,%
\begin{equation}
\mathcal{S}_{n}=\left\{
\dirack{\mathbf{x}}%
\in\mathbb{Z}^{n}\mid\exists\,%
\diracb{\mathbf{t}}%
\in\mathbb{Z}^{n}\ni\joinrel\relbar%
\ip{\mathbf{t}}{\mathbf{x}}%
=1\right\}  . \label{eqClmNJun.63}%
\end{equation}
Note that $\operatorname*{GL}\left(  n,\mathbb{Z}\right)  $ (i.e., the group
of matrices in $M_{n}\left(  \mathbb{Z}\right)  $ with determinant $\pm1$)
acts on $%
\dirack{\mathcal{S}_{n}}%
$ by multiplication from the left. The reason is that if $\mathbf{x}%
\in\mathcal{S}_{n}$, there is a $t\in\mathbb{Z}^{n}$ with $%
\ip{\mathbf{t}}{\mathbf{x}}%
=1$, and hence $%
\diracb{\mathbf{t}}%
A^{-1}A%
\dirack{\mathbf{x}}%
=1$, so $A\mathbf{x}\in\mathcal{S}_{n}$. We next argue that%
\begin{equation}
\text{The action of }\operatorname*{GL}\left(  n,\mathbb{Z}\right)  \text{ on
}\mathcal{S}_{n}\text{ is transitive.} \label{eqClmNJun.64}%
\end{equation}

\begin{proof}
By \cite[Theorem II.1]{New72}, there exists for any $%
\dirack{\alpha}%
\in%
\dirack{\mathcal{S}_{n}}%
$ a matrix $A_{\alpha}\in\operatorname*{GL}\left(  n,\mathbb{Z}\right)  $ such
that the first column of $A_{\alpha}$ is $%
\dirack{\alpha}%
$: this means $%
\dirack{\alpha}%
=A_{\alpha}%
\dirack{\mathbf{e}_{1}}%
$. But if $%
\dirack{\beta}%
\in%
\dirack{\mathcal{S}_{n}}%
$, this means $A_{\alpha}^{{}}A_{\beta}^{-1}%
\dirack{\beta}%
=%
\dirack{\alpha}%
$ and transitivity follows.
\end{proof}

Now, let us prove

\begin{lemma}
\label{LemClmNJun.19}If $\alpha^{\left(  1\right)  }$, $\alpha^{\left(
2\right)  }$, $\mathbf{v}_{1}$, $\mathbf{v}_{2}\in\mathcal{S}_{n}$, the
following conditions \textup{(\ref{eqClmNJun.65})} and
\textup{(\ref{eqClmNJun.66})} are equivalent if $n\geq3$.%
\begin{equation}
\text{There is an }A\in\operatorname*{GL}\left(  n,\mathbb{Z}\right)  \text{
such that }%
\diracb{\smash{\alpha^{\left( 1\right) }}}%
A=%
\diracb{\smash{\alpha^{\left( 2\right) }}}%
\text{ and }A%
\dirack{\mathbf{v}_{1}}%
=%
\dirack{\mathbf{v}_{2}}%
.\!\!\!\!\! \label{eqClmNJun.65}%
\end{equation}%
\begin{equation}%
\ip{\smash{\alpha^{\left( 1\right) }}}{\mathbf{v}_{2}}%
=%
\ip{\smash{\alpha^{\left( 2\right) }}}{\mathbf{v}_{1}}%
. \label{eqClmNJun.66}%
\end{equation}
\textup{(}The implication \textup{(\ref{eqClmNJun.65})} $\Rightarrow$
\textup{(\ref{eqClmNJun.66})} is true for all $n$, but the converse
implication may fail for $n=2$.\textup{)}
\end{lemma}

\begin{proof}
(Due to K.H. Kim and F. Roush.) If (\ref{eqClmNJun.65}) holds, then%
\[%
\ip{\smash{\alpha^{\left( 1\right) }}}{\mathbf{v}_{2}}%
=%
\diracb{\smash{\alpha^{\left( 1\right) }}}%
A%
\dirack{\mathbf{v}_{1}}%
=%
\ip{\smash{\alpha^{\left( 2\right) }}}{\mathbf{v}_{1}}%
.
\]
If, conversely, (\ref{eqClmNJun.66}) holds, first choose matrices
$U,V\in\operatorname*{GL}\left(  n,\mathbb{Z}\right)  $ with%
\[%
\diracb{\smash{\alpha^{\left( 1\right) }}}%
U=%
\diracb{\mathbf{e}_{1}}%
\text{,\qquad}V%
\dirack{\mathbf{v}_{1}}%
=%
\dirack{\mathbf{e}_{1}}%
.
\]
This is possible by (\ref{eqClmNJun.64}). It follows from Lemma
\ref{LemClmNJun.20}, below, that there exists a matrix $B\in\operatorname*{GL}%
\left(  n,\mathbb{Z}\right)  $ such that the first column in $B$ is $U^{-1}%
\dirack{\mathbf{v}_{2}}%
$ and the first row in $B$ is $%
\diracb{\smash{\alpha^{\left( 2\right) }}}%
V^{-1}$. For this, we note that $U^{-1}%
\dirack{\mathbf{v}_{2}}%
\in%
\dirack{\mathcal{S}_{n}}%
$ and $%
\diracb{\smash{\alpha^{\left( 2\right) }}}%
V^{-1}\in%
\diracb{\mathcal{S}_{n}}%
$, and the first component of $U^{-1}%
\dirack{\mathbf{v}_{2}}%
$ is $%
\diracb{\mathbf{e}_{1}}%
U^{-1}%
\dirack{\mathbf{v}_{2}}%
=%
\diracb{\smash{\alpha^{\left( 1\right) }}}%
UU^{-1}%
\dirack{\mathbf{v}_{2}}%
=%
\ip{\smash{\alpha^{\left( 1\right) }}}{\mathbf{v}_{2}}%
=%
\ip{\smash{\alpha^{\left( 2\right) }}}{\mathbf{v}_{1}}%
=%
\diracb{\smash{\alpha^{\left( 2\right) }}}%
V^{-1}V%
\dirack{\mathbf{v}_{1}}%
=%
\diracb{\smash{\alpha^{\left( 2\right) }}}%
V^{-1}%
\dirack{\mathbf{e}_{1}}%
=$ the first component of $%
\diracb{\smash{\alpha^{\left( 2\right) }}}%
V^{-1}$. Now put%
\[
A=UBV.
\]
Then%
\[%
\diracb{\smash{\alpha^{\left( 1\right) }}}%
A=%
\diracb{\smash{\alpha^{\left( 1\right) }}}%
UBV=%
\diracb{\mathbf{e}_{1}}%
BV=%
\diracb{\smash{\alpha^{\left( 2\right) }}}%
V^{-1}V=%
\diracb{\smash{\alpha^{\left( 2\right) }}}%
,
\]
and%
\[
A%
\dirack{\mathbf{v}_{1}}%
=UBV%
\dirack{\mathbf{v}_{1}}%
=UB%
\dirack{\mathbf{e}_{1}}%
=UU^{-1}%
\dirack{\mathbf{v}_{2}}%
=%
\dirack{\mathbf{v}_{2}}%
.%
\settowidth{\qedskip}{$\displaystyle A\dirack
{\mathbf{v}_{1}}=UBV\dirack{\mathbf{v}_{1}}=UB\dirack{\mathbf{e}_{1}}%
=UU^{-1}\dirack{\mathbf{v}_{2}}=\dirack{\mathbf{v}_{2}}.$}\addtolength
{\qedskip}{-\textwidth}\rlap{\makebox[-0.5\qedskip][r]{\qedsymbol}}%
\]%
\renewcommand{\qed}{}%
\end{proof}

Thus, we have to prove

\begin{lemma}
\label{LemClmNJun.20}Let $\alpha$, $\beta$ be vectors in $\mathcal{S}_{n}$
with $\alpha_{1}=\beta_{1}$, and assume that $n\geq3$. Then there exists a
matrix $V\in\operatorname*{GL}\left(  n,\mathbb{Z}\right)  $ such that the
first column in $B$ is $%
\dirack{\alpha}%
$ and the first row is $%
\diracb{\beta}%
$.
\end{lemma}

\begin{proof}
(Due to K.H. Kim and F. Roush.) We will use the fact that row, respectively
column, operations on a matrix can be effectuated by multiplying from the
left, respectively right, by matrices in $\operatorname*{GL}\left(
n,\mathbb{Z}\right)  $. For example, interchanging the first two rows in $A$
corresponds to multiplying $A$ from the left by $\left(
\begin{smallmatrix}
0 & 1\\
1 & 0
\end{smallmatrix}
\right)  \oplus\openone_{n-2}$, and adding $\mu$ times row $2$ to row $1$
corresponds to left-multiplying by $\left(
\begin{smallmatrix}
1 & \mu\\
0 & 1
\end{smallmatrix}
\right)  \oplus\openone_{n-2}$. The corresponding column operations follow by
taking the transpose. If now $A\in\operatorname*{GL}\left(  n,\mathbb{Z}%
\right)  $ is a matrix of the form%
\[
A=\left(
\begin{tabular}
[c]{c|c}%
$\vphantom{\vdots}\alpha_{1}$ &
\begin{tabular}
[c]{ccc}%
$\beta_{2}$ & $\cdots$ & $\beta_{n}$%
\end{tabular}
\\\hline%
\begin{tabular}
[c]{c}%
$\vphantom{\vdots}\alpha_{2}$\\
$\vdots$\\
$\vphantom{\vdots}\alpha_{n}$%
\end{tabular}
&
\begin{tabular}
[c]{ccc}%
$\ast$ &  & $\ast$\\
&  & \\
$\ast$ &  & $\ast$%
\end{tabular}
\end{tabular}
\right)  ,
\]
let $\gamma_{1}=\gcd\left(  \beta_{2},\dots,\beta_{n}\right)  $, and choose
$\rho_{2},\rho_{3},\dots,\rho_{d}$ such that $\rho_{2}\beta_{2}+\dots+\rho
_{n}\beta_{n}=\gamma_{1}$. Let $U_{1}$ be a matrix in $\operatorname*{GL}%
\left(  n-1,\mathbb{Z}\right)  $ with first column $\left(
\begin{smallmatrix}
\rho_{2}\\
\vdots\\
\rho_{n}%
\end{smallmatrix}
\right)  $ (it exists by \cite[Theorem II.1]{New72}). Now multiply by $A$ from
the left to obtain%
\[
A\left(
\begin{tabular}
[c]{c|ccc}%
$\vphantom{\vdots}1$ & $0$ & $\cdots$ & $0$\\\hline
$\vphantom{\vdots}0$ &  &  & \\
$\vdots$ &  & $U_{1}$ & \\
$\vphantom{\vdots}0$ &  &  &
\end{tabular}
\right)  =A\left(
\begin{tabular}
[c]{c|ccc}%
$\vphantom{\vdots}1$ & $0$ & $\cdots$ & $0$\\\hline
$\vphantom{\vdots}0$ & $\rho_{2}$ &  & $\ast$\\
$\vdots$ & $\vdots$ &  & \\
$\vphantom{\vdots}0$ & $\rho_{n}$ &  & $\ast$%
\end{tabular}
\right)  =\left(
\begin{tabular}
[c]{c|cccc}%
$\vphantom{\vdots}\alpha_{1}$ & $\gamma_{1}$ & $\gamma_{2}$ & $\cdots$ &
$\gamma_{n-1}$\\\hline
$\vphantom{\vdots}\alpha_{2}$ & \raisebox{-6pt}[0pt][0pt]{\rlap{$\mkern
18mu\ast$}} &  & \raisebox{-6pt}[0pt][0pt]{\rlap{$\mkern18mu\ast$}} & \\
$\vphantom{\vdots}\alpha_{3}$ &  &  &  & \\
$\vdots$ & \raisebox{-6pt}[0pt][0pt]{\rlap{$\mkern18mu\ast$}} &  &
\raisebox{-6pt}[0pt][0pt]{\rlap{$\mkern18mu\ast$}} & \\
$\vphantom{\vdots}\alpha_{n}$ & $\hphantom{\cdots}$ & $\hphantom{\cdots}$ &  &
$\hphantom{\cdots}$%
\end{tabular}
\right)  ,
\]
where the remaining elements $\gamma_{2},\dots,\gamma_{n-1}$ on the first row
are linear combinations of $\beta_{2},\dots,\beta_{n}$, and thus multiples of
$\gamma_{1}$. By subtracting integer multiples of the second column from the
remaining columns, one finally finds a matrix $U\in\operatorname*{GL}\left(
n,\mathbb{Z}\right)  $ such that%
\[
AU=\left(
\begin{tabular}
[c]{c|cccc}%
$\vphantom{\vdots}\alpha_{1}$ & $\gamma_{1}$ & $0$ & $\cdots$ & $0$\\\hline
$\vphantom{\vdots}\alpha_{2}$ & \raisebox{-6pt}[0pt][0pt]{\rlap{$\mkern
18mu\ast$}} &  & \raisebox{-6pt}[0pt][0pt]{\rlap{$\mkern18mu\ast$}} & \\
$\vphantom{\vdots}\alpha_{3}$ &  &  &  & \\
$\vdots$ & \raisebox{-6pt}[0pt][0pt]{\rlap{$\mkern18mu\ast$}} &  &
\raisebox{-6pt}[0pt][0pt]{\rlap{$\mkern18mu\ast$}} & \\
$\vphantom{\vdots}\alpha_{n}$ & $\hphantom{\cdots}$ & $\hphantom{\cdots}$ &  &
$\hphantom{\cdots}$%
\end{tabular}
\right)  .
\]
Putting $\gamma_{2}=\gcd\left(  \alpha_{2},\dots,\alpha_{n}\right)  $ and
transposing all these operations, one finds a $V\in\operatorname*{GL}\left(
n,\mathbb{Z}\right)  $ such that%
\[
VAU=\left(
\begin{tabular}
[c]{c|cccc}%
$\vphantom{\vdots}\alpha_{1}$ & $\gamma_{1}$ & $0$ & $\cdots$ & $0$\\\hline
$\vphantom{\vdots}\gamma_{2}$ & \raisebox{-6pt}[0pt][0pt]{\rlap{$\mkern
18mu\ast$}} &  & \raisebox{-6pt}[0pt][0pt]{\rlap{$\mkern18mu\ast$}} & \\
$\vphantom{\vdots}0$ &  &  &  & \\
$\vdots$ & \raisebox{-6pt}[0pt][0pt]{\rlap{$\mkern18mu\ast$}} &  &
\raisebox{-6pt}[0pt][0pt]{\rlap{$\mkern18mu\ast$}} & \\
$\vphantom{\vdots}0$ & $\hphantom{\cdots}$ & $\hphantom{\cdots}$ &  &
$\hphantom{\cdots}$%
\end{tabular}
\right)  ,
\]
where $\gcd\left(  \gamma_{2},\alpha_{1}\right)  =1$. Thus, if we can prove
Lemma \ref{LemClmNJun.20} for this kind of matrices, the general Lemma
\ref{LemClmNJun.20} follows by multiplying from the left and right by the
inverses $V^{-1}$, $U^{-1}$. This reduces the proof of Lemma
\ref{LemClmNJun.20} to the case%
\[
\alpha_{3}=\alpha_{4}=\dots=\alpha_{n}=0=\beta_{3}=\beta_{4}=\dots=\beta_{n},
\]
and $\gcd\left(  \alpha_{1},\alpha_{2}\right)  =\gcd\left(  \beta_{1}%
,\beta_{2}\right)  =1$ where still $\alpha_{1}=\beta_{1}$. But to this end one
can use the matrix%
\[%
\begin{pmatrix}
\alpha_{1} & \beta_{2} & 0 & 0 & \cdots & 0\\
\alpha_{2} & 1 & x & 0 & \cdots & 0\\
0 & 1 & y & 0 & \cdots & 0\\
0 & 0 & 0 & 1 &  & 0\\
\vdots & \vdots & \vdots &  & \ddots & \vdots\\
0 & 0 & 0 & 0 & \cdots & 1
\end{pmatrix}
.
\]
The determinant is $\alpha_{1}\left(  y-x\right)  -\alpha_{2}\beta_{2}y$, but
as $\gcd\left(  \alpha_{1},\alpha_{2}\beta_{2}\right)  =1$, this can be made
equal to $1$ by choosing the integers $\left(  y-x\right)  $ and $\left(
-y\right)  $ by the Euclidean algorithm. This ends the proof of Lemma
\ref{LemClmNJun.20}, and thus of Lemma \ref{LemClmNJun.19}. Note that if $n=2$
the proof above does not work: One must choose an integer $x$ such that%
\[%
\begin{vmatrix}
\alpha_{1} & \beta_{2}\\
\alpha_{2} & x
\end{vmatrix}
=\alpha_{1}x-\alpha_{2}\beta_{2}=\pm1
\]
when $\alpha_{1}$, $\alpha_{2}$, $\beta_{2}$ are given with $\gcd\left(
\alpha_{1},\alpha_{2}\right)  =1=\gcd\left(  \alpha_{1},\beta_{2}\right)  $,
and this is clearly impossible in general.
\end{proof}

\begin{proof}
[Proof of Theorem \textup{\ref{ThmClmNJun.18}}](Due to K.H. Kim and F. Roush.)
If $J$, $J^{\prime}$ are matrices of the form (\ref{eqCycNew.1}) with
$\lambda=m_{N}$ and $\lambda^{\prime}=m_{N^{\prime}}^{\prime}$, then
$N=N^{\prime}$ and $\operatorname*{Prim}\left(  \lambda\right)
=\operatorname*{Prim}\left(  \lambda^{\prime}\right)  $ by Theorem
\ref{CorCyc.9}, and then $I\left(  J\right)  =I\left(  J^{\prime}\right)  $ by
Corollary \ref{CorClmN.E1}. This proves (\ref{ThmClmNJun.18(1)}) $\Rightarrow$
(\ref{ThmClmNJun.18(2)}). The converse implication follows in the cases
$N=1,2,3$ by Theorem \ref{ThmClmN.2} and the discussion around (\ref{eqClmN.5}%
), so we may assume $N\geq4$ from now on. Assuming (\ref{ThmClmNJun.18(2)}) it
follows from Lemma \ref{LemClmN.J} and Remark \ref{RemClmN.J1} that
$\mathfrak{A}_{J}$ and $\mathfrak{A}_{J^{\prime}}$ are isomorphic if and only
if there exists a matrix $A\in\operatorname*{GL}\left(  N-1,\mathbb{Z}\right)
$ such that%
\begin{equation}%
\diracb{\beta^{\prime}}%
A=%
\diracb{\beta}%
\mod{I\left( J\right) \mathbb{Z}^{N-1}} \label{eqClmNJun.67}%
\end{equation}
where we now define%
\begin{equation}
\beta^{\prime}=\lambda^{N-1}\left(  \lambda^{\prime\,N-2},\lambda
^{\prime\,N-3},\dots,\lambda^{\prime},1\right)  \label{eqClmNJun.68}%
\end{equation}
and%
\begin{equation}
\beta=\lambda^{\prime\,N-1}\left(  \lambda^{N-2},\lambda^{N-3},\dots
,\lambda,1\right)  \label{eqClmNJun.69}%
\end{equation}
and%
\begin{equation}
A%
\dirack{\mathbf{\bar{v}}}%
=%
\dirack{\mathbf{\bar{v}}^{\prime}}%
. \label{eqClmNJun.70}%
\end{equation}
Now, one checks that%
\begin{equation}%
\ip{\beta^{\prime}}{\mathbf{\bar{v}}^{\prime}}%
=\lambda_{{}}^{N-1}\left(  \lambda_{{}}^{\prime\,N-2}c_{2}^{\prime}%
+\dots+\lambda_{{}}^{\prime}c_{N-1}^{\prime}+1\right)  =\lambda^{N-1}I\left(
J\right)  -\lambda^{N-1}\lambda^{\prime\,N-1} \label{eqClmNJun.71}%
\end{equation}
and%
\begin{equation}%
\ip{\beta}{\mathbf{\bar{v}}}%
=\lambda^{\prime\,N-1}I\left(  J\right)  -\lambda^{\prime\,N-1}\lambda^{N-1},
\label{eqClmNJun.72}%
\end{equation}
and thus%
\begin{equation}%
\ip{\beta^{\prime}}{\mathbf{\bar{v}}^{\prime}}%
=%
\ip{\beta}{\mathbf{\bar{v}}}%
+\left(  \lambda^{N-1}-\lambda^{\prime\,N-1}\right)  I\left(  J\right)  =%
\ip{\beta}{\mathbf{\bar{v}}}%
\mod{I\left( J\right) }. \label{eqClmNJun.73}%
\end{equation}
Now we cannot apply Lemma \ref{LemClmNJun.19} directly on $\alpha^{\left(
1\right)  }=\beta^{\prime}$, $\alpha^{\left(  2\right)  }=\beta$,
$\mathbf{v}_{1}=\mathbf{\bar{v}}$, $\mathbf{v}_{2}=\mathbf{\bar{v}}^{\prime}$,
for two reasons: we do not have $\beta,\beta^{\prime}\in\mathcal{S}_{n}$, and
the condition (\ref{eqClmNJun.66}) is only fulfilled modulo $I\left(
J\right)  $. But let us remedy this by modifying $\beta$, $\beta^{\prime}$ as
follows: First add an integer multiple of $I\left(  J\right)  \left(
0,0,\dots,1\right)  $ to $\beta^{\prime}$ to obtain a new $\beta^{\prime}$,
called $\beta^{\left(  1\right)  \prime}$, such that $%
\ip{\smash{\beta^{\left( 1\right) \prime}}}{\mathbf{\bar{v}}^{\prime}}%
=%
\ip{\beta}{\mathbf{\bar{v}}}%
$. This is possible since the last component of $\mathbf{\bar{v}}^{\prime}$ is
$1$ by (\ref{eqClmN.Y}). Now modify the new $\beta^{\left(  1\right)  \prime}$
to $\beta^{\left(  2\right)  \prime}$ by adding integer multiples of the
vector $I\left(  J\right)  \left(  0,0,\dots,0,-1,v_{N-1}^{\prime}\right)  $
to $\beta^{\left(  1\right)  \prime}$ until the second-to-last component
contains none of the prime factors in $\operatorname*{Prim}\left(
\lambda\right)  =\operatorname*{Prim}\left(  \lambda^{\prime}\right)  $. This
is possible since $I\left(  J\right)  $ is relatively prime to $\lambda$ and
$\lambda^{\prime}$. Then $%
\ip{\smash{\beta^{\left( 2\right) \prime}}}{\mathbf{\bar{v}}^{\prime}}%
=%
\ip{\smash{\beta^{\left( 1\right) \prime}}}{\mathbf{\bar{v}}^{\prime}}%
$ since $\left(  0,0,\dots,0,-1,v_{N-1}^{\prime}\right)  $ is orthogonal to
$\mathbf{\bar{v}}^{\prime}$ and $\beta^{\left(  2\right)  \prime}%
\in\mathcal{S}_{n}$. Finally, modify $\beta$ to $\beta^{\left(  2\right)  }$
by adding multiples of $I\left(  J\right)  \left(  0,0,\dots,0,-1,v_{N-1}%
\right)  $ until the second-to-last component is relatively prime to the first
$N-3$ components. Then $%
\ip{\smash{\beta^{\left( 2\right) }}}{\mathbf{\bar{v}}}%
=%
\ip{\beta}{\mathbf{\bar{v}}}%
$ and hence%
\[%
\ip{\smash{\beta^{\left( 2\right) \prime}}}{\mathbf{\bar{v}}^{\prime}}%
=%
\ip{\beta}{\mathbf{\bar{v}}}%
.
\]
But since $N-1\geq3$ we may now apply Lemma \ref{LemClmNJun.19} to find an
$A\in\operatorname*{GL}\left(  N-1,\mathbb{Z}\right)  $ such that%
\[%
\diracb{\smash{\beta^{\left( 2\right) \prime}}}%
A=%
\diracb{\smash{\beta^{\left( 2\right) }}}%
\text{\quad and\quad}A%
\dirack{\mathbf{\bar{v}}}%
=%
\dirack{\mathbf{\bar{v}}^{\prime}}%
.
\]
But since%
\[%
\diracb{\smash{\beta^{\left( 2\right) \prime}}}%
=%
\diracb{\beta^{\prime}}%
\mod{I\left( J\right) \mathbb{Z}^{N-1}}%
\]
and%
\[%
\diracb{\smash{\beta^{\left( 2\right) }}}%
=%
\diracb{\beta}%
\mod{I\left( J\right) \mathbb{Z}^{N-1}},
\]
it follows that%
\[%
\diracb{\beta^{\prime}}%
A=%
\diracb{\beta}%
\mod{I\left( J\right) \mathbb{Z}^{N-1}}.
\]
Thus (\ref{eqClmNJun.67}) and (\ref{eqClmNJun.70}) are fulfilled and Theorem
\ref{ThmClmNJun.18} is proved.
\end{proof}

Let us end this chapter by mentioning that the equivalent conditions
(\ref{ThmClmNJun.18(1)}) and (\ref{ThmClmNJun.18(2)}) again are equivalent to
the condition that $\mathfrak{A}_{J}$ and $\mathfrak{A}_{J^{\prime}}$ are
stably isomorphic, i.e., to that $\mathfrak{A}_{J}\otimes\mathcal{K}\left(
\mathcal{H}\right)  $ is isomorphic to $\mathfrak{A}_{J^{\prime}}%
\otimes\mathcal{K}\left(  \mathcal{H}\right)  $, where $\mathcal{K}\left(
\mathcal{H}\right)  $ is the $C^{\ast}$-algebra of compact operators on a
separable Hilbert space $\mathcal{H}$. This is due to the very special
position of $\left[  \openone\right]  $ in $K_{0}\left(  \mathfrak{A}%
_{J}\right)  $, and this property has no general analogue: For example, it is
impossible to find an automorphism of $\mathbb{Z}\left[  \frac{1}{2}\right]  $
mapping $1$ into $3$.

\begin{corollary}
\label{CorClmNJun.21}Let $J$, $J^{\prime}$ be matrices of the form
\textup{(\ref{eqCycNew.1})} satisfying the standard requirement there, and
assume that $\lambda=m_{N}$ and $\lambda^{\prime}=m_{N^{\prime}}^{\prime}$.
Then the following three conditions are equivalent.

\begin{enumerate}
\item \label{CorClmNJun.21(1)}The triples $\left(  K_{0}\left(  \mathfrak
{A}_{J}\right)  ,K_{0}\left(  \mathfrak{A}_{J}\right)  _{+},\left[
\openone\right]  \right)  $ and $\left(  K_{0}\left(  \mathfrak{A}_{J^{\prime
}}\right)  ,K_{0}\left(  \mathfrak{A}_{J^{\prime}}\right)  _{+},\left[
\openone\right]  \right)  $ are isomorphic, i.e., the dimension groups are
isomorphic by an isomorphism mapping $\left[  \openone\right]  $ into $\left[
\openone\right]  $.

\item \label{CorClmNJun.21(2)}The dimension groups $\left(  K_{0}\left(
\mathfrak{A}_{J}\right)  ,K_{0}\left(  \mathfrak{A}_{J}\right)  _{+}\right)  $
and $\left(  K_{0}\left(  \mathfrak{A}_{J^{\prime}}\right)  ,K_{0}\left(
\mathfrak{A}_{J^{\prime}}\right)  _{+}\right)  $ are isomorphic.

\item \label{CorClmNJun.21(3)}$N=N^{\prime}$, $\operatorname*{Prim}\left(
\lambda\right)  =\operatorname*{Prim}\left(  \lambda^{\prime}\right)  $, and
$I\left(  J\right)  =I\left(  J^{\prime}\right)  $.
\end{enumerate}
\end{corollary}

\begin{proof}
The equivalence (\ref{CorClmNJun.21(1)}) $\Leftrightarrow$
(\ref{CorClmNJun.21(3)}) is Theorem \ref{ThmClmNJun.18}, and
(\ref{CorClmNJun.21(1)}) $\Rightarrow$ (\ref{CorClmNJun.21(2)}) is trivial.
Thus it suffices to show that (\ref{CorClmNJun.21(2)}) $\Rightarrow$
(\ref{CorClmNJun.21(3)}), so assume (\ref{CorClmNJun.21(2)}). Then one
establishes $N=N^{\prime}$ and $\operatorname*{Prim}\left(  \lambda\right)
=\operatorname*{Prim}\left(  \lambda^{\prime}\right)  $ exactly as in Theorem
\ref{CorCyc.9}, and it remains to establish $I\left(  J\right)  =I\left(
J^{\prime}\right)  $. For this one notes that Proposition \ref{prop10.7}
remains true in the context of nonunital isomorphism with the exception that
the condition $5$ is just removed, and condition $2$ is replaced by%
\[%
\diracb{\alpha^{\prime}}%
\Lambda=\mu%
\diracb{\alpha}%
,
\]
where $\mu$ is a positive scalar. But since $\mu$ induces an automorphism on
the range $\mathbb{Z}\left[  \frac{1}{\lambda}\right]  =\mathbb{Z}\left[
\frac{1}{\lambda^{\prime}}\right]  $ of the trace, it follows that $\mu$ is an
invertible element of the ring $\mathbb{Z}\left[  \frac{1}{\lambda}\right]  $.
Now the Lemmas \ref{LemClmN.A}--\ref{LemClmN.C} do not involve $\left[
\openone\right]  $ and are still valid, and then Lemma \ref{LemClmN.D} is
valid with the same proof. Now the equation in the proof of Corollary
\ref{CorClmN.E} becomes%
\[%
\ip{\alpha}{\mathbf{v}}%
=\mu^{-1}%
\ip{\alpha^{\prime}}{\Lambda\mathbf{v}}%
=\mu^{-1}%
\ip{\alpha^{\prime}}{\xi\mathbf{v}^{\prime}}%
=\mu^{-1}\xi%
\ip{\alpha^{\prime}}{\mathbf{v}^{\prime}}%
,
\]
so Corollary \ref{CorClmN.E} is still valid. The proof that $I\left(
J\right)  =I\left(  J^{\prime}\right)  $ is now exactly as in the proof of
Corollary \ref{CorClmN.E1}.
\end{proof}

\setcounter{figurelink}{\value{figure}}

\chapter{\label{App.Fur}Further comments on two examples from Chapter
\ref{APP.EXA}}

We now consider two sub-examples with $N=3$, $d=5$. Although the groups
$K_{0}\left(  \mathfrak{A}_{L}\right)  $\label{K0frakALter} and $K_{0}\left(
\mathfrak
{A}_{L^{\prime}}\right)  $ have the same rank, we will show directly in
Examples \ref{Exa1} and \ref{Exa2} that they are non-isomorphic. The details
also serve to illustrate what goes into the computation of some particular
inductive limit which is not immediately transparent.

The two algebras corresponding to $x+4x^{3}=1$ and $3x^{2}+2x^{3}=1$ are the
two with stabilized diagrams%
\[
\begin{picture}(360,114)(0,-6)
\put(0,0){\includegraphics
[bb=135 68 297 170,clip,height=102bp,width=162bp]{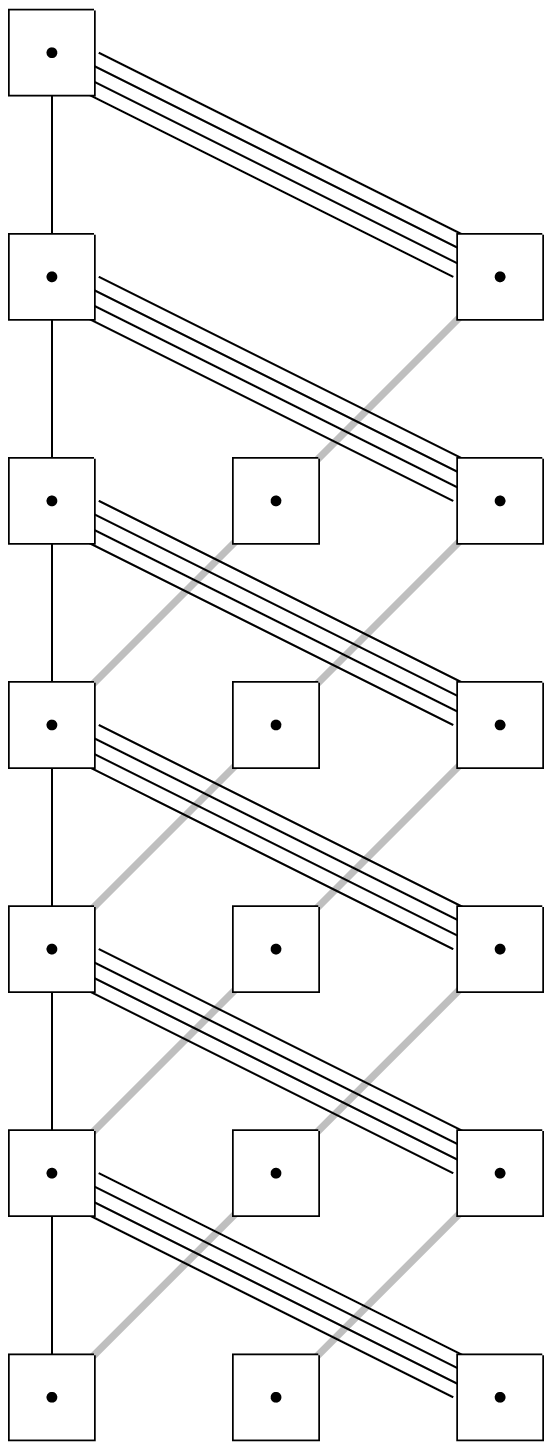}}
\put(198,0){\includegraphics
[bb=135 68 297 170,clip,height=102bp,width=162bp]{bb22233.eps}}
\end{picture},
\]
i.e., with incidence matrices%
\[
J_{1}=%
\begin{pmatrix}
1 & 1 & 0\\
0 & 0 & 1\\
4 & 0 & 0
\end{pmatrix}
,\qquad J_{2}=%
\begin{pmatrix}
0 & 1 & 0\\
3 & 0 & 1\\
2 & 0 & 0
\end{pmatrix}
.
\]

The two examples above both have $\beta=\ln2$, and $\operatorname*{rank}%
\left(  K_{0}\left(  \mathfrak{A}_{L}\right)  \right)  =3$ (and $d=5$), but
the two dimension groups are actually non-isomorphic since the
$\operatorname*{Prim}\left(  Q\right)  $ invariants are different; see the
$N=3$ case in Chapter \ref{APP.EXA}. We will, however, establish this the hard
way in Examples \ref{Exa1} and \ref{Exa2} by showing that the respective
$\ker\left(  \tau\right)  $-groups in the two examples are non-isomorphic.

\begin{example}
\label{Exa1}This is an elaboration on Remark \textup{\ref{RemPediferient.nada}%
.}
\[
J_{1}=%
\begin{pmatrix}
1 & 1 & 0\\
0 & 0 & 1\\
4 & 0 & 0
\end{pmatrix}
\]
The Frobenius eigenvalue is $2$ \textup{(}see \textup{(\ref{eqRemNew.1})) and
the corresponding }normalized left eigenvector $\alpha$ is%
\[
\alpha=\left(  1,\tfrac{1}{2},\tfrac{1}{4}\right)
\]
by \textup{(\ref{eqRemNewBis.17}).} $J_{1}$ has eigenvalues $2$ and
$\lambda_{\pm}=\dfrac{-1\pm\sqrt{-7}}{2}$, $\left|  \lambda_{\pm}\right|  <2$;
$J_{1}$ leaves the subspace orthogonal to $\alpha=\left(  1,\frac{1}{2}%
,\frac{1}{4}\right)  $ invariant. This is spanned by the vectors%
\[
e_{1}=%
\begin{pmatrix}
1\\
-2\\
0
\end{pmatrix}
,\qquad e_{2}=%
\begin{pmatrix}
0\\
1\\
-2
\end{pmatrix}
,
\]
and one computes that the matrix of the restriction of $J$ to this subspace is%
\[
J_{0}=%
\begin{pmatrix}
-1 & 1\\
-2 & 0
\end{pmatrix}
.
\]
We may compute the iterated inverses $J_{0}^{-n}$ by straightforward spectral
theory, and if we put%
\[
\mu=\left(  1-i\sqrt{7}\right)  /4
\]
the result is%
\[
J_{0}^{-n}=%
\begin{pmatrix}
-1 & 1\\
-2 & 0
\end{pmatrix}
^{-n}=\frac{1}{\left(  -2\right)  ^{n\mathstrut}}\frac{1}{\mu-\bar{\mu}}%
\begin{pmatrix}
\frac{1}{\mu^{n-1\mathstrut}}-\frac{1}{\bar{\mu}^{n-1\mathstrut}} &
-\frac{\bar{\mu}}{\mu^{n-1\mathstrut}}+\frac{\mu}{\bar{\mu}^{n-1\mathstrut}}\\
\frac{1}{\mu^{n\mathstrut}}-\frac{1}{\bar{\mu}^{n\mathstrut}} & -\frac
{\bar{\mu}}{\mu^{n\mathstrut}}+\frac{\mu}{\bar{\mu}^{n\mathstrut}}%
\end{pmatrix}
.
\]
The eigenvector of $J_{1}$ with eigenvalue $2$ is
\[
e_{0}=%
\begin{pmatrix}
1\\
1\\
2
\end{pmatrix}
.
\]
Thus, if%
\[
\Delta=%
\begin{pmatrix}
e_{0} & e_{1} & e_{2}%
\end{pmatrix}
=%
\begin{pmatrix}
1 & 1 & 0\\
1 & -2 & 1\\
2 & 0 & -2
\end{pmatrix}
,
\]
then
\[
\Delta^{-1}J_{1}\Delta=%
\begin{pmatrix}
2 & 0 & 0\\
0 & -1 & 1\\
0 & -2 & 0
\end{pmatrix}
\]
and%
\[
\Delta_{{}}^{-1}J_{1}^{-n}\Delta=\left(  \Delta^{-1}J_{1}\Delta\right)
^{-n}=\left(
\begin{tabular}
[c]{c|c}%
$2^{-n}$ & $%
\begin{array}
[c]{cc}%
0 & 0
\end{array}
$\\\hline
$%
\begin{array}
[c]{c}%
0\\
0
\end{array}
$ & $J_{0}^{-n}$%
\end{tabular}
\right)  .
\]
Thus%
\[
J_{1}^{-n}=\Delta\left(
\begin{tabular}
[c]{c|c}%
$2^{-n}$ & $%
\begin{array}
[c]{cc}%
0 & 0
\end{array}
$\\\hline
$%
\begin{array}
[c]{c}%
0\\
0
\end{array}
$ & $J_{0}^{-n}$%
\end{tabular}
\right)  \Delta^{-1}.
\]
We have%
\[
\Delta^{-1}=\frac{1}{8}%
\begin{pmatrix}
4 & 2 & 1\\
4 & -2 & -1\\
4 & 2 & -3
\end{pmatrix}
.
\]
Thus the dimension group, as a subgroup of $\mathbb{R}^{3}$ ($N=3$), is%
\[
\bigcup_{n=0}^{\infty}\frac{1}{8}%
\begin{pmatrix}
1 & 1 & 0\\
1 & -2 & 1\\
2 & 0 & -2
\end{pmatrix}
\left(
\begin{tabular}
[c]{c|c}%
$2^{-n}$ & $%
\begin{array}
[c]{cc}%
0 & 0
\end{array}
$\\\hline
$%
\begin{array}
[c]{c}%
0\\
0
\end{array}
$ & $J_{0}^{-n}$%
\end{tabular}
\right)
\begin{pmatrix}
4 & 2 & 1\\
4 & -2 & -1\\
4 & 2 & -3
\end{pmatrix}
\mathbb{Z}^{3}%
\]
(where the union is increasing), i.e., $\bigcup_{n}M_{n}\mathbb{Z}^{3}$ when
$M_{n}$ is the product matrix. The range of the trace on $K_{0}$ is $\alpha$
applied to this set (the range of the trace on $\mathfrak{A}_{L}$ is gotten by
intersecting with $\left[  0,1\right]  $). We have%
\begin{align*}
{}  &  \alpha\frac{1}{8}%
\begin{pmatrix}
1 & 1 & 0\\
1 & -2 & 1\\
2 & 0 & -2
\end{pmatrix}%
\begin{pmatrix}
2^{-n} & 0 & 0\\
0 & \ast & \ast\\
0 & \ast & \ast
\end{pmatrix}%
\begin{pmatrix}
4 & 2 & 1\\
4 & -2 & -1\\
4 & 2 & -3
\end{pmatrix}
\\
{}  &  \qquad=\frac{1}{4}%
\begin{pmatrix}
1 & 0 & 0
\end{pmatrix}%
\begin{pmatrix}
2^{-n} & 0 & 0\\
0 & \ast & \ast\\
0 & \ast & \ast
\end{pmatrix}%
\begin{pmatrix}
4 & 2 & 1\\
4 & -2 & -1\\
4 & 2 & -3
\end{pmatrix}
\\
{}  &  \qquad=%
\begin{pmatrix}
2^{-n-2} & 0 & 0
\end{pmatrix}%
\begin{pmatrix}
4 & 2 & 1\\
\ast & \ast & \ast\\
\ast & \ast & \ast
\end{pmatrix}
\\
{}  &  \qquad=%
\begin{pmatrix}
2^{-n} & 2^{-n+1} & 2^{-n-2}%
\end{pmatrix}
,
\end{align*}
i.e., $\alpha M_{n}=%
\begin{pmatrix}
2^{-n} & 2^{-n-1} & 2^{-n-2}%
\end{pmatrix}
$, and applying this to $\mathbb{Z}^{3}$ we reconfirm that the range of the
trace is the set of all dyadic rationals. The range of the trace on the $m$'th
term in
\[
\mathbb{Z}^{3}\longrightarrow\mathbb{Z}^{3}\longrightarrow\cdots
\]
is%
\begin{align*}
2^{-m+1}%
\ip{\alpha}{\mathbb{Z}^{3}}%
&  =2^{-m+1}\left(  1,\tfrac{1}{2},\tfrac{1}{4}\right)  \mathbb{Z}^{3}\\
&  =2^{-m+1}\mathbb{Z}.
\end{align*}
Now, the infinitesimal elements of the $m$'th $\mathbb{Z}^{3}$ are the
elements of the kernel of $\alpha J_{1}^{-m+1}=2^{-m+1}\alpha$, i.e., the
elements of the kernel of $\alpha$, that is elements of $\mathbb{Z}^{3}$ of
the form%
\[%
\begin{pmatrix}
n_{1}\\
n_{2}\\
-4n_{1}-2n_{2}%
\end{pmatrix}
=n_{1}e_{1}+\left(  2n_{1}+n_{2}\right)  e_{2},
\]
where $n_{1},n_{2}\in\mathbb{Z}$. Thus%
\[
J_{1}^{-n}%
\begin{pmatrix}
n_{1}\\
n_{2}\\
-4n_{1}-2n_{2}%
\end{pmatrix}
=\Delta\left(
\begin{tabular}
[c]{c|c}%
$2^{-n}$ & $%
\begin{array}
[c]{cc}%
0 & 0
\end{array}
$\\\hline
$%
\begin{array}
[c]{c}%
0\\
0
\end{array}
$ & $J_{0}^{-n}$%
\end{tabular}
\right)
\begin{pmatrix}
0\\
n_{1}\\
2n_{1}+n_{2}%
\end{pmatrix}
.
\]
Thus the dimension group of $\mathfrak{A}_{L}$ is an
extension\label{LOSG0_14}
\[
0\hooklongrightarrow G_{0}\hooklongrightarrow K_{0}\left(  \mathfrak{A}%
_{L}\right)  \overset{\tau}{\longrightarrow}\mathbb{Z}\left[  \tfrac{1}%
{2}\right]  \longrightarrow0,
\]
where $G_{0}=\ker\tau$\label{LOSkertau_9} is the group of infinitesimal
elements. It is rather complicated to describe $G_{0}$ in the matrix formalism
above, so let us instead use the algebraic description in
\textup{(\ref{eqRemFeb.a}),} that is,%
\[
G_{0}\cong\mathbb{Z}\left[  x\right]  \diagup p_{0}\left(  x\right)
\mathbb{Z}\left[  x\right]  ,
\]
where%
\[
p_{0}\left(  x\right)  =\frac{p_{L}\left(  x\right)  }{p_{a}\left(  x\right)
}=\frac{x+4x^{3}-1}{2x-1}=2x^{2}+x+1.
\]
Now, embed $\mathbb{Z}$ in $G_{0}$ as the group $H$ generated by $1\pmod
{p_{0}\left( x\right) }$ \textup{(}note that $n1\neq0\pmod{p_{0}\left
( x\right) }$ for all $n\in\mathbb{Z}\diagdown\left\{  0\right\}  $, so this
is really an embedding\textup{).} We argue that%
\[
G_{0}\diagup H\cong\mathbb{Z}\left[  \tfrac{1}{2}\right]  .
\]
We have
\begin{align*}
G_{0}\diagup H  &  =\left(  \mathbb{Z}\left[  x\right]  \diagup p_{0}\left(
x\right)  \mathbb{Z}\left[  x\right]  \right)  \diagup\mathbb{Z}1\\
&  =\left(  \mathbb{Z}\left[  x\right]  \diagup\left(  2x^{2}+x+1\right)
\mathbb{Z}\left[  x\right]  \right)  \diagup\mathbb{Z}1\\
&  =\mathbb{Z}\left[  x\right]  \diagup\left(  \left(  2x^{2}+x+1\right)
\mathbb{Z}\left[  x\right]  +1\mathbb{Z}\right)
\end{align*}
as abelian groups. Let $\left(  p\left(  x\right)  \right)  $ denote the
residue class of the polynomial $p$ in $G_{0}\diagup H$.

Let%
\[
u_{0}=\left(  x\right)  =\left(  -2x^{2}\right)  =2u_{1},
\]
where $u_{1}=\left(  -x^{2}\right)  $. Since%
\[
\left(  2x^{n+1}+x^{n}+x^{n-1}\right)  =0,
\]
we obtain%
\[
\left(  x^{n}\right)  =-\left(  x^{n-1}\right)  +2\left(  -x^{n+1}\right)
\]
for $n=1,2,\dots$. It follows by induction that the elements $\left(
x^{n}\right)  $ are divisible by $2$ for $n=1,2,\dots$ \textup{(}this is not
true for $n=0$ by a use of Lemma \textup{\ref{LemSubNew.14}.)} Furthermore, we
can find monic polynomials $p_{n}$ of degree $n+1$ such that the sequence
$u_{n}=\left(  p_{n}\left(  x\right)  \right)  $ has the property%
\[
u_{n+1}=2u_{n}%
\]
for all $n$. But then $\{1,p_{0}\left(  x\right)  =x,p_{1}\left(  x\right)
=-x^{2},\dots,p_{2k-1}\left(  x\right)  =-x^{2}+x^{3}-\dots-x^{2k}%
,p_{2k}\left(  x\right)  =x^{2}-x^{3}+\dots-x^{2k+1},\cdots\}$ span all of
$\mathbb{Z}\left[  x\right]  $, so $u_{0},u_{1},\dots$ span all of
$G_{0}\diagup H$. It follows that%
\[
G_{0}\diagup H\cong\mathbb{Z}\left[  \tfrac{1}{2}\right]  .
\]
Thus $G_{0}$ is an extension of the form%
\[
0\longrightarrow\mathbb{Z}\longrightarrow G_{0}\longrightarrow\mathbb{Z}%
\left[  \tfrac{1}{2}\right]  \longrightarrow0.
\]
In conclusion we have the exact diagram\label{K0frakALtetra}%
\[%
\begin{array}
[c]{ccc}%
& 0 & \\
& \downarrow & \\
& \mathbb{Z} & \\
& \downarrow & \\
0\longrightarrow\!\! & G_{0} & \!\!\longrightarrow K_{0}\left(  \mathfrak
{A}_{L}\right)  \overset{\tau}{\longrightarrow}\mathbb{Z}\left[  \frac{1}%
{2}\right]  \longrightarrow0.\\
& \downarrow & \\
& \makebox[0pt]{\hss$\mathbb{Z}\left[ \frac{1}{2}\right] $\hss} & \\
& \downarrow & \\
& 0 &
\end{array}
\]
\end{example}

\begin{example}
\label{Exa2}
\[
J_{2}=%
\begin{pmatrix}
0 & 1 & 0\\
3 & 0 & 1\\
2 & 0 & 0
\end{pmatrix}
\]
Again the Frobenius eigenvalue is $2$ and the normalized solution of%
\[
\alpha J_{2}=2\alpha
\]
is%
\[
\alpha=\left(  1,\tfrac{1}{2},\tfrac{1}{4}\right)  .
\]
$J_{2}$ has eigenvalues $2$, $-1$, $-1$. With $e_{1}$, $e_{2}$ as before,
$J_{2}$ leaves the subspace spanned by $e_{1}$ and $e_{2}$ invariant, and the
matrix of the restriction is%
\[%
\begin{pmatrix}
-2 & 1\\
-1 & 0
\end{pmatrix}
.
\]
The determinant is $1$, and one computes that
\[
e_{1}+e_{2}=%
\begin{pmatrix}
1\\
-1\\
-2
\end{pmatrix}
\]
is the unique eigenvector with eigenvalue $-1$. Using%
\[
f_{1}=e_{1}+e_{2}=%
\begin{pmatrix}
1\\
-1\\
-2
\end{pmatrix}
\]
and%
\[
f_{2}=e_{2}=%
\begin{pmatrix}
0\\
1\\
-2
\end{pmatrix}
\]
as basis instead, one computes that the matrix is%
\[%
\begin{pmatrix}
-1 & 1\\
0 & -1
\end{pmatrix}
\]
and hence%
\[%
\begin{pmatrix}
-1 & 1\\
0 & -1
\end{pmatrix}
^{-n}=\left(  -1\right)  ^{n}%
\begin{pmatrix}
1 & n\\
0 & 1
\end{pmatrix}
\]
for $n=0,1,2,\dots$. The right eigenvector of $J_{2}$ with eigenvalue $2$ is%
\[
f_{0}=%
\begin{pmatrix}
1\\
2\\
1
\end{pmatrix}
.
\]
Thus, if%
\[
\Delta=%
\begin{pmatrix}
f_{0} & f_{1} & f_{2}%
\end{pmatrix}
=%
\begin{pmatrix}
1 & 1 & 0\\
2 & -1 & 1\\
1 & -2 & -2
\end{pmatrix}
,
\]
then
\[
\Delta^{-1}J_{2}\Delta=%
\begin{pmatrix}
2 & 0 & 0\\
0 & -1 & 1\\
0 & 0 & -1
\end{pmatrix}
\]
and%
\[
\Delta_{{}}^{-1}J_{2}^{-n}\Delta=\left(  \Delta^{-1}J_{2}\Delta\right)  ^{-n}=%
\begin{pmatrix}
2^{-n} & 0 & 0\\
0 & \left(  -1\right)  ^{n} & \left(  -1\right)  ^{n}n\\
0 & 0 & \left(  -1\right)  ^{n}%
\end{pmatrix}
.
\]
Thus%
\[
J_{2}^{-n}=\Delta%
\begin{pmatrix}
2^{-n} & 0 & 0\\
0 & \left(  -1\right)  ^{n} & \left(  -1\right)  ^{n}n\\
0 & 0 & \left(  -1\right)  ^{n}%
\end{pmatrix}
\Delta^{-1}.
\]
We have%
\[
\Delta^{-1}=\frac{1}{9}%
\begin{pmatrix}
4 & 2 & 1\\
5 & -2 & -1\\
-3 & 3 & -3
\end{pmatrix}
.
\]
Thus the dimension group as a subgroup of $\mathbb{R}^{3}$ is%
\[
\bigcup_{n=0}^{\infty}\frac{1}{9}%
\begin{pmatrix}
1 & 1 & 0\\
2 & -1 & 1\\
1 & -2 & -2
\end{pmatrix}%
\begin{pmatrix}
2^{-n} & 0 & 0\\
0 & \left(  -1\right)  ^{n} & \left(  -1\right)  ^{n}n\\
0 & 0 & \left(  -1\right)  ^{n}%
\end{pmatrix}%
\begin{pmatrix}
4 & 2 & 1\\
5 & -2 & -1\\
-3 & 3 & -3
\end{pmatrix}
\mathbb{Z}^{3}.
\]
The range of the trace is $\alpha=\left(  1,\frac{1}{2},\frac{1}{4}\right)  $
applied to this set, which is%
\[
\bigcup_{n=0}^{\infty}\frac{1}{4}%
\begin{pmatrix}
1 & 0 & 0
\end{pmatrix}%
\begin{pmatrix}
2^{-n} & 0 & 0\\
\ast & \ast & \ast\\
\ast & \ast & \ast
\end{pmatrix}%
\begin{pmatrix}
4 & 2 & 1\\
\ast & \ast & \ast\\
\ast & \ast & \ast
\end{pmatrix}
\mathbb{Z}^{3}=\bigcup_{n=0}^{\infty}%
\begin{pmatrix}
2^{-n} & 2^{-n+1} & 2^{-n-2}%
\end{pmatrix}
\mathbb{Z}^{3},
\]
which is not unexpectedly the set of dyadic rationals. Since the determinant
of the matrix $\left(
\begin{smallmatrix}
-1 & 1\\
0 & -1
\end{smallmatrix}
\right)  $ (or $\left(
\begin{smallmatrix}
-2 & 1\\
-1 & 0
\end{smallmatrix}
\right)  $) is $1$, this matrix defines a $1$---$1$ map on the infinitesimal
elements, so the group of infinitesimal elements is isomorphic to
$\mathbb{Z}^{2}$. Thus $K_{0}\left(  \mathfrak{A}_{L}\right)  $ is an
extension\label{LOSkertau_10}
\[
0\hooklongrightarrow\mathbb{Z}^{2}\hooklongrightarrow K_{0}\left(
\mathfrak{A}_{L}\right)  \overset{\tau}{\longrightarrow}\mathbb{Z}\left[
\tfrac{1}{2}\right]  \longrightarrow0.
\]
\end{example}

We see that the dimension groups of the two examples \ref{Exa1} and \ref{Exa2}
are non-isomorphic, so the algebras are non-isomorphic.

\begin{acknowledgements}
The authors are grateful to A.~Kishimoto, K.-H. Kim, M. Laca, A. Laudal, and
F.~Roush for very helpful discussions, to B. Treadway and C. Wehr for
excellent typesetting, and to B. Treadway for graphics production and
programming, and for useful suggestions. We are especially thankful to Kim and
Roush for permitting us to publish their proof of Theorem \ref{ThmClmNJun.18}
via Lemmas \ref{LemClmNJun.19} and \ref{LemClmNJun.20}, which they provided
after seeing a preprint version of this memoir.
The referee kindly offered suggestions as to how the subject
could better be presented to a wider audience.
\end{acknowledgements}

\backmatter
\bibliographystyle{amsplain}
\bibliography{jorgen}

\renewcommand{\frogleg}{\ } \renewcommand{\frogelbow}{\\}%

\listoffigures

\listoftables

\chapter*{List of Terms and Symbols}

\LOSitem{\textbf{Term or Symbol}}{\textbf{Usage}}{\textbf{Page}}\bigskip

\LOSitem
{
AF-algebra
}{
\
}{
\pageref{LOSAFalgebra_0},
\pageref{LOSAFalgebra_00},
\pageref{LOSAFalgebra_1},
\pageref{LOSAFalgebra_2},
\pageref{LOSAFalgebra_3}
}

\LOSitem
{
Bratteli diagram
}{
\
}{
\pageref{LOSBrattelidiagrams_1},
\pageref{LOSBrattelidiagrams_2},
\pageref{LOSBrattelidiagrams_3},
\pageref{LOSBrattelidiagrams_4},
\pageref{LOSBrattelidiagrams_5}
}

\LOSitem
{
dimension group
}{
\
}{
\pageref{LOSdimensiongroup_0},
\pageref{LOSdimensiongroup_1},
\pageref{LOSdimensiongroup_2}
}

\LOSitem
{
matrix units
}{
\
}{
\pageref{LOSmatrixunits}
}

\LOSitem
{
Perron--Frobenius theorem
}{
\
}{
\pageref{LOSPerronFrobeniustheorem}
}

\LOSitem
{
$\#_{j}\left(  \alpha\right)  =\#\left\{  \alpha_{i}\mid\alpha_{i}%
=j\right\}  $
}{
count
}{
\pageref{LOShashjalpha}
}

\LOSitem
{
$\left[  \openone\right]  $
}{
unit in $K_{0}\left(  \mathfrak{A}_{L}\right)  $
}{
\pageref{LOSoneonesquare}
}

\LOSitem
{
$a:=\lambda^{-1}=e^{-\beta}$
}{
inverse Perron--Frobenius eigenvalue
}{
\pageref{LOSa}
}

\LOSitem
{
$\mathfrak{A}_{L}=\mathfrak{A}_{J}=\mathcal{O}_{d}^{\sigma^{L}}\left(
0\right)  $
}{
the AF-algebra with symbol list $L=\left\{  L_{1},\dots,L_{d}\right\}  $ and
matrix $J$
}{
\pageref{LOSAL_1},
\pageref{LOSAL_2},
\pageref{LOSAL_3}
}

\LOSitem
{
$A_{n}:=\left\{  \left(  \alpha\right)  \mid\alpha\in L^{-1}\left(  n\right)
\cup E_{n}\right\}  $
}{
family of projections
}{
\pageref{LOSAn}
}

\LOSitem
{
$\mathfrak{A}_{n}=\bigoplus_{k=0}^{L_{d}-1}\mathfrak{A}^{\left(  n,k\right
)  }$
}{
direct sum decomposition of $\mathfrak{A}_{n}$
}{
\pageref{LOSAnm}
}

\LOSitem
{
$D\left(  \mathfrak{A}_{L}\right)  =\left(  K_{0}\left(  \mathfrak{A}%
_{L}\right)  ,K_{0}\left(  \mathfrak{A}_{L}\right)  _{+},\left[  \openone
\right]  \right)  $
}{
dimension group
}{
\pageref{LOSDAL_1},
\pageref{LOSDAL_2},
\pageref{LOSDAL_3},
\pageref{LOSDAL_4},
\pageref{LOSDAL_5}
}

\LOSitem
{
$\mathcal{D}_{d}\cong C\left(  \prod_{1}^{\infty}\mathbb{Z}_{d}\right)
=C^{\ast}\left(  s_{\alpha}^{{}}s_{\alpha}^{\ast}\bigm|\alpha\in
\vphantom{\coprod}\smash{\coprod_{1}^{\infty}}\mathbb{Z}_{d}\right)
\vphantom{\coprod_{1}^{\infty}}$
}{
abelian subalgebra of diagonal elements
}{
\pageref{LOSDd_1},
\pageref{LOSDd_2},
\pageref{LOSDd_3}
}

\LOSitem
{
$\deg$
}{
scaling degree
}{
\pageref{LOSdeg_1},
\pageref{LOSdeg_2},
\pageref{LOSdeg_3},
\pageref{LOSdeg_4}
}

\LOSitem
{
$\det\left(  t\openone-J\right)  =t^{N}-m_{1}t^{N-1}-m_{2}t^{N-2}%
-\dots-m_{N-1}t-m_{N}$
}{
characteristic polynomial of $J$
}{
\pageref{LOSdett1J}
}

\LOSitem
{
$D_{n}(G) =\bigcap_{i}n^{i}G$
}{
the elements of $G$ infinitely divisible by $n$
}{
\pageref{LOSDlambdaG_1},
\pageref{LOSDlambdaG_2}
}

\LOSitem
{
$e_{j}:=x^{j}p_{a}\left(  x\right)  $
}{
basis for $\left\{  \alpha\right\}  ^{\perp}\cap\mathbb{Z}^{N}$
}{
\pageref{LOSej}
}

\LOSitem
{
$E_{n}=\{  \gamma\mid\gamma=\left(  \alpha i\right)  $ where
$L\left(  \alpha\right)  <n$ and $L\left(  \alpha\right)  +L_{i}>n\}  $
}{
complementary set of projections
}{
\pageref{LOSEn}
}

\LOSitem
{
$E_{n}\left(  m\right)  =\left\{  \gamma\in E_{n}\mid L\left(  \gamma\right)
=n+m\right\}  $
}{
indexed complementary set of projections
}{
\pageref{LOSEnm}
}

\LOSitem
{
$\operatorname*{Ext}\left(  \tau\left(  K_{0}\left(  \mathfrak{A}_{L}\right)
\right)  ,\ker\tau\right)  $
}{
the $\operatorname*{Ext}$-functor
}{
\pageref{LOSExt}
}

\LOSitem
{
$\operatorname*{Ext}$-groups
}{
\
}{
\pageref{LOSExtgroups}
}

\LOSitem
{
$\mathbb{F}_{d}$
}{
free group on $d$ generators
}{
\pageref{LOSFd}
}

\LOSitem
{
$F_{L}:=G_{L}\diagup\mathbb{Z}^{N}$
}{
a torsion group quotient
}{
\pageref{LOSFL_1},
\pageref{LOSFL_2},
\pageref{LOSFL_3}
}

\LOSitem
{
$f_{m}\left(  x\right)  :=m_{1}+m_{2}x+\dots+m_{N}x^{N-1}=q_{m}\left(
x\right)  p_{a}\left(  x\right)  +r_{m}\left(  x\right)  $
}{
characteristic polynomial
}{
\pageref{LOSfm_1},
\pageref{LOSfm_2},
\pageref{LOSfm_3}
}

\LOSitem
{
$G=K_{0}\left(  \mathfrak{A}_{L}\right)  $
}{
$K_{0}$ group
}{
\pageref{LOSG}
}

\LOSitem
{
$G_{0}=\ker\tau\cong\mathbb{Z}\left[  x\right]  \diagup\left(  p_{0}\left(
x\right)  \right)  $
}{
kernel of trace $\tau=$ infinitesimal elements in $G$
}{
\pageref{LOSG0_1},
\pageref{LOSG0_2},
\pageref{LOSG0_3},
\pageref{LOSG0_4},
\pageref{LOSG0_5},
\pageref{LOSG0_6},
\pageref{LOSG0_7},
\pageref{LOSG0_8},
\pageref{LOSG0_9},
\pageref{LOSG0_10},
\pageref{LOSG0_11},
\pageref{LOSG0_12},
\pageref{LOSG0_13},
\pageref{LOSG0_14}
}

\LOSitem
{
$g_{i}=\left\{
\begin{array}
[c]{ll}J^{-i}e_{N} & \text{if }i=1,2,\dots\\
e_{i+N} & \text{if }i=1-N,\dots,-1,0
\end{array}
\right.  $
}{
\
}{
\pageref{LOSgi}
}

\LOSitem
{
$G_{i}:=J^{-i}\mathbb{Z}^{N}$
}{
indexed subgroups in $G_{J}$
}{
\pageref{LOSGi}
}

\LOSitem
{
$G_{J}=\bigcup_{n}J^{-n}\mathbb{Z}^{N}$
}{
inductive limit group
}{
\pageref{LOSGJ}
}

\LOSitem
{
$\left(  G\left(  L\right)  ,\left(  \sigma_{L}\right)  _{\ast}\right)  $
}{
shift dynamical system
}{
\pageref{LOSGLsigmaLstar}
}

\LOSitem
{
$\operatorname*{grade}\left(  s_{\alpha}^{{}}s_{\gamma}^{\ast}\right)
=L\left(  \alpha\right)  $ if $L\left(  \alpha\right)  =L\left(
\gamma\right)  $
}{
grade function on monomials
}{
\pageref{LOSgrade}
}

\LOSitem
{
$\mathcal{H}=\int_{\Omega}^{\oplus}\mathcal{H}\left(  x\right)  \,d\mu\left(
x\right)  $
}{
direct integral of Hilbert spaces
}{
\pageref{LOSH_1},
\pageref{LOSH_2},
\pageref{LOSH_3},
\pageref{LOSH_4}
}

\LOSitem
{
$\mathcal{H}_{-}$
}{
closed linear span in $\ell^{2}\left(  \mathbb{F}_{d}\right)  $ of the vectors
$\left\{  \lambda\left(  s^{-1}\right)  \xi_{e}\mid s\in\mathbb{S}_{d}%
\right\}  $
}{
\pageref{LOSHminus_1},
\pageref{LOSHminus_2}
}

\LOSitem
{
$\mathcal{H}_{\Omega_{0}}$
}{
cyclic subspace of the representation of $\mathcal{O}_{d}$ induced from the
state $\omega_{\left(  p\right)  }\left(  s_{\alpha}^{{}}s_{\gamma}^{\ast
}\right)  =p_{\alpha} \delta_{\alpha\gamma}$
}{
\pageref{LOSHOmega_1},
\pageref{LOSHOmega_2},
\pageref{LOSHOmega_3}
}

\LOSitem
{
$\mathcal{H}_{\omega_{p}}$
}{
Hilbert space of the state $\omega_{p}$
}{
\pageref{LOSHomegap}
}

\LOSitem
{
$I\left(  J\right)  =\sum_{i=1}^{N}v_{i}\lambda^{N-i}=\lambda^{N-1}%
\TeXButton{ip}{\ip{\alpha}{v}}$
}{
invariant
}{
\pageref{LOSIJ_1},
\pageref{LOSIJ_2}
}

\LOSitem
{
$J=J_{m}=J_{L}=
\begin{pmatrix}
\vphantom{\vdots}m_{1} & 1 & 0 & \cdots& 0 & 0\\
\vphantom{\vdots}m_{2} & 0 & 1 & \cdots& 0 & 0\\
\vdots&  & \ddots& \ddots& \vdots& \vdots\\
\vphantom{\vdots}m_{N-2} & 0 &  & \ddots& 1 & 0\\
\vphantom{\vdots}m_{N-1} & 0 & 0 &  & 0 & 1\\
\vphantom{\vdots}m_{N} & 0 & 0 & \cdots& 0 & 0
\end{pmatrix}
$
}{
matrix with integer entries $m_{i}$ with $m_{N}\neq0$ and $\gcd\left\{
i\mid m_{i}\neq0\right\}  =1$
}{
\pageref{LOSJ_1},
\pageref{LOSJ_2},
\pageref{LOSJ_3},
\pageref{LOSJ_4},
\pageref{LOSJ_5},
\pageref{LOSJ_6},
\pageref{LOSJ_7}
}

\LOSitem
{
$J_{0}=\begin{pmatrix}
\vphantom{\vdots}q_{1} & 1 & 0 & \cdots& 0 & 0\\
\vphantom{\vdots}q_{2} & 0 & 1 & \cdots& 0 & 0\\
\vdots&  & \ddots& \ddots& \vdots& \vdots\\
\vphantom{\vdots}q_{M-2} & 0 &  & \ddots& 1 & 0\\
\vphantom{\vdots}q_{M-1} & 0 & 0 &  & 0 & 1\\
\vphantom{\vdots}q_{M} & 0 & 0 & \cdots& 0 & 0
\end{pmatrix}
$
}{
matrix of form similar to $J$
}{
\pageref{LOSJ0_1},
\pageref{LOSJ0_2}
}

\LOSitem
{
$J_{L}=\left(
\begin{tabular}
[c]{c|c}$J_{0_{\mathstrut}}$ & $Q$\\\hline
$0$ & $J_{D}$\end{tabular}
\right)  $
}{
triangular representation of matrix $J$
}{
\pageref{LOSJL_1},
\pageref{LOSJL_2},
\pageref{LOSJL_3},
\pageref{LOSJL_4}
}

\LOSitem
{
$J^{-1}=\left(
\begin{tabular}
[c]{c|c}$J_{0_{\mathstrut}}^{-1}$ & $-J_{0_{\mathstrut}}^{-1}QJ_{R_{\mathstrut
}}^{-1}$\\\hline
$0$ & $J_{R}^{-1^{\mathstrut}}$\end{tabular}
\right)  $
}{
triangular representation of $J^{-1}$
}{
\pageref{LOSJLinverse}
}

\LOSitem
{
$\TeXButton{ip}{\ip{k}{\alpha}}=\sum_{i=1}^{N}k_{i}a^{i-1}=f_{k}%
\left(  a\right)  $
}{
inner product
}{
\pageref{LOSkalpha}
}

\LOSitem
{
$K_{0}\left(  \mathfrak{A}_{L}\right)  $
}{
the $K_{0}$ group of $\mathfrak{A}_{L}$ (see $G$)
}{
\pageref{LOSK0AL_1},
\pageref{LOSK0AL_2},
\pageref{LOSK0AL_3},
\pageref{LOSK0AL_4},
\pageref{LOSK0AL_5},
\pageref{LOSK0AL_6},
\pageref{LOSK0AL_7},
\pageref{LOSK0AL_8},
\pageref{LOSK0AL_9}
}

\LOSitem
{
$K_{0}\left(  \mathfrak{A}_{L}\right)  _{+}=\left\{  0\right\}  \cup\left\{
v\in K_{0}\left(  \mathfrak{A}_{L}\right)  \mid\TeXButton{ip}{\ip{\alpha}{v}%
}>0\right\}  $
}{
positive elements of $K_{0}\left(  \mathfrak{A}\right)  $
}{
\pageref{LOSK0ALplus_1},
\pageref{LOSK0ALplus_2}
}

\LOSitem
{
$\ker\tau$
}{
kernel of $\tau=$ infinitesimal elements of $G$
}{
\pageref{LOSkertau_1},
\pageref{LOSkertau_2},
\pageref{LOSkertau_3},
\pageref{LOSkertau_4},
\pageref{LOSkertau_5},
\pageref{LOSkertau_6},
\pageref{LOSkertau_7},
\pageref{LOSkertau_8},
\pageref{LOSkertau_9},
\pageref{LOSkertau_10}
}

\LOSitem
{
$L=\left\{  L_{1},\dots,L_{d}\right\}  $
}{
list of symbols
}{
\pageref{LOSL}
}

\LOSitem
{
$L\left(  \alpha\right)  =\sum_{j=1}^{d}\#_{j}\left(  \alpha\right)
L_{j}=\sum_{m=1}^{k}L_{\alpha_{m}}$
}{
weight function
}{
\pageref{LOSLalpha_1},
\pageref{LOSLalpha_2}
}

\LOSitem
{
$\mathcal{L}$
}{
inductive limit lattice
}{
\pageref{LOSLinductivelimit_1},
\pageref{LOSLinductivelimit_2}
}

\LOSitem
{
$\mathcal{L}_{k}$
}{
lattice with index
}{
\pageref{LOSLk_1},
\pageref{LOSLk_2}
}

\LOSitem
{
$L^{-1}\left(  n\right)  =\left\{  \alpha\mid L\left(  \alpha\right)
=n\right\}  $
}{
projections with specified $L$-number
}{
\pageref{LOSLinversen}
}

\LOSitem
{
$M_{d}=M_{d}\left(  \mathbb{C}\right)  $
}{
algebra of $d\times d$ complex matrices
}{
\pageref{LOSMd}
}

\LOSitem
{
$m$-$\deg$
}{
$m$-scaling degree
}{
\pageref{LOSmdeg_1},
\pageref{LOSmdeg_2},
\pageref{LOSmdeg_3}
}

\LOSitem
{
$M_{f}=\int_{\Omega}^{\oplus}f\left(  x\right)  \openone_{\mathcal{H}\left(
x\right)  }\,d\mu\left(  x\right)  $
}{
representation of $f$ in $\mathcal{H}$
}{
\pageref{LOSMf}
}

\LOSitem
{
$m_{N}=\left|  \det J\right|  $
}{
\
}{
\pageref{LOSmN_1},
\pageref{LOSmN_2}
}

\LOSitem
{
$\mathcal{O}_{d}$
}{
the Cuntz algebra
}{
\pageref{LOSOd_1},
\pageref{LOSOd_2},
\pageref{LOSOd_3},
\pageref{LOSOd_4},
\pageref{LOSOd_5},
\pageref{LOSOd_6}
}

\LOSitem
{
$\mathcal{O}_{d}^{\sigma^{L}}\left(  n\right)  =\mathcal{O}_{d}^{\sigma
}\left(  n\right)  =\left\{  x\in\mathcal{O}_{d}\mid\sigma_{t}^{\left(
L\right)  }\left(  x\right)  =e^{int}x\right\}  $
}{
spectral subspaces
}{
\pageref{LOSOdsigmaLn}
}

\LOSitem
{
$\mathcal{O}_{d}^{\mathbb{T}}=\operatorname*{UHF}_{d}$
}{
fixed-point algebra of the gauge action of $\mathbb{T}$
}{
\pageref{LOSOdT}
}

\LOSitem
{
$\mathcal{O}_{d}^{\mathbb{T}^{d}}=\operatorname*{GICAR}_{d}$
}{
$\mathbb{T}^{d}$-gauge-invariant elements of $\mathcal{O}_{d}$
}{
\pageref{LOSOdTd}
}

\LOSitem
{
$\mathcal{P}=\left\{  \sum_{k=1}^{d}n_{k}L_{k}\bigm|n_{k}\in\mathbb{N}%
\cup\left\{  0\right\}  \right\}  $
}{
semigroup generated by $L_{1},\dots,L_{d}$
}{
\pageref{LOSP}
}

\LOSitem
{
$p_{0}\left(  x\right)  =\sum_{j=1}^{M}q_{j}x^{j}-1$
}{
$q_{M}$ times the characteristic polynomial of $J_{0}^{-1}$
}{
\pageref{LOSp0x_1},
\pageref{LOSp0x_2},
\pageref{LOSp0x_3}
}

\LOSitem
{
$p_{a}\left(  x\right)  =x^{D}p_{\lambda}\left(  \frac{1}{x}\right)  $
}{
minimal polynomial of inverse of Perron--Frobenius eigenvalue $a$
}{
\pageref{LOSpa}
}

\LOSitem
{
$p^{\alpha}=p_{\alpha_{1}}p_{\alpha_{2}}\cdots p_{\alpha_{k}}$
}{
monomials
}{
\pageref{LOSpalpha}
}

\LOSitem
{
$p_{j}=e^{-\beta L_{j}}$
}{
probability weights
}{
\pageref{LOSpj_1},
\pageref{LOSpj_2},
\pageref{LOSpj_3}
}

\LOSitem
{
$p_{\lambda}\left(  x\right)  =x^{D}p_{a}\left(  \frac{1}{x}\right)  $,
$a=\frac{1}{\lambda}$
}{
minimal polynomial of Perron--Frobenius eigenvalue $\lambda$
}{
\pageref{LOSplambda}
}

\LOSitem
{
$p_{L}\left(  x\right)  =\sum_{j=1}^{N}m_{j}x^{j}-1=p_{0}\left(  x\right)
p_{a}\left(  x\right)  $
}{
$\left|  \det J\right|  $ times the characteristic polynomial of $J^{-1}$
}{
\pageref{LOSpLx_1},
\pageref{LOSpLx_2},
\pageref{LOSpLx_3},
\pageref{LOSpLx_4},
\pageref{LOSpLx_5},
\pageref{LOSpLx_6},
\pageref{LOSpLx_7}
}

\LOSitem
{
$\operatorname*{Prim}\left(  n\right)  $
}{
the set of prime factors of $n$
}{
\pageref{LOSPrim_1},
\pageref{LOSPrim_2},
\pageref{LOSPrim_3},
\pageref{LOSPrim_4},
\pageref{LOSPrim_5},
\pageref{LOSPrim_6},
\pageref{LOSPrim_7},
\pageref{LOSPrim_8}
}

\LOSitem
{
$\operatorname*{Prim}\left(  m_{N}\right)  $
}{
the set of prime factors of $m_{N}=\left|  \det J\right|  $
}{
\pageref{LOSPrimmN}
}

\LOSitem
{
$\operatorname*{Prim}\left(  Q_{N-D}\right)  $
}{
the set of prime factors of $Q_{N-D}$
}{
\pageref{LOSPrimQND}
}

\LOSitem
{
$\operatorname*{Prim}\left(  R_{D}\right)  $
}{
the set of prime factors of $R_{D}=\left|  \det J_{R}\right|  $
}{
\pageref{LOSPrimRD}
}

\LOSitem
{
$q_{m}\left(  x\right)  =\sum_{k=1}^{N-D}Q_{k}x^{k-1}$
}{
special polynomial
}{
\pageref{LOSqm}
}

\LOSitem
{
$\operatorname*{rank}$
}{
rank of an abelian group
}{
\pageref{LOSrank_1},
\pageref{LOSrank_2}
}

\LOSitem
{
$r_{m}\left(  x\right)  =\sum_{k=1}^{D}R_{k}x^{k-1}$
}{
residue polynomial
}{
\pageref{LOSrm_1},
\pageref{LOSrm_2}
}

\LOSitem
{
$s_{i}$
}{
generators of $\mathcal{O}_{d}$
}{
\pageref{LOSsi}
}

\LOSitem
{
$s_{\alpha_{{}}}=s_{\alpha_{1}}s_{\alpha_{2}}\cdots s_{\alpha_{n}}$
}{
monomial in $\mathcal{O}_{d}$
}{
\pageref{LOSsalpha_1},
\pageref{LOSsalpha_2}
}

\LOSitem
{
$\mathbb{S}_{d}$
}{
free semigroup
}{
\pageref{LOSSd}
}

\LOSitem
{
$S_{i}$
}{
operators representing $s_{i}$
}{
\pageref{LOSSi_1},
\pageref{LOSSi_2}
}

\LOSitem
{
$S_{i}\sp\ast$
}{
adjoint operators of $S_{i}$
}{
\pageref{LOSSistar}
}

\LOSitem
{
$S_{\alpha_{{}}}=S_{\alpha_{1}}S_{\alpha_{2}}\cdots S_{\alpha_{n}}$
}{
monomial of operators
}{
\pageref{LOSSalpha_1},
\pageref{LOSSalpha_2}
}

\LOSitem
{
$\mathbb{T}^{d}$
}{
the $d$-torus
}{
\pageref{LOSTd}
}

\LOSitem
{
$U=\sum_{i=1}^{d}S_{i}^{{}}T_{i}^{\ast}=\int_{\Omega}^{\oplus}U\left(
x\right)  \,d\mu\left(  x\right)  $
}{
intertwining unitary having decomposition
}{
\pageref{LOSU_1},
\pageref{LOSU_2}
}

\LOSitem
{
$\operatorname*{UHF}\nolimits_{d}\cong\bigotimes_{1}^{\infty}M_{d}$
}{
uniformly hyperfinite $C^{\ast}$-algebra
}{
\pageref{LOSUHFd}
}

\LOSitem
{
$v=
\begin{pmatrix}
\vphantom{\vdots}1\\
\vphantom{\vdots}\lambda-m_{1}\\
\vphantom{\vdots}\lambda^{2}-m_{1}\lambda-m_{2}\\
\vphantom{\vdots}\lambda^{3}-m_{1}\lambda^{2}-m_{2}\lambda-m_{3}\\
\vdots\\
\vphantom{\vdots}\lambda^{N-1}-\dots-m_{N-2}\lambda-m_{N-1}\end{pmatrix}
$
}{
Perron--Frobenius right eigenvector: $Jv=\lambda v$
}{
\pageref{LOSPerronFrobeniusrighteigenvector_1},
\pageref{LOSPerronFrobeniusrighteigenvector_2},
\pageref{LOSPerronFrobeniusrighteigenvector_3},
\pageref{LOSPerronFrobeniusrighteigenvector_4},
\pageref{LOSPerronFrobeniusrighteigenvector_5}
}

\LOSitem
{
$\mathcal{V}_{N}:=\left\{  f\left(  x\right)  \in\mathbb{Z}\left[  x\right]
\mid\deg f\leq N-1\right\}  \cong\mathbb{Z}^{N}$
}{
representation of special elements of $G$
}{
\pageref{LOSVN_1},
\pageref{LOSVN_2}
}

\LOSitem
{
$x_{i}=g_{i}\mod{\mathbb{Z}^{N}}$
}{
generators of the group $F_{J}:=G_{J}\diagup\mathbb{Z}^{N}$
}{
\pageref{LOSxi_1},
\pageref{LOSxi_2},
\pageref{LOSxi_3}
}

\LOSitem
{
$\mathbb{Z}\left[  \tfrac{1}{\lambda}\right]  $
}{
ring generated by $1$ and $\tfrac{1}{\lambda}$
}{
\pageref{LOSZ1lambda_1},
\pageref{LOSZ1lambda_2}
}

\LOSitem
{
$\mathbb{Z}_{d}=\left\{  1,\dots,d\right\}  $
}{
set of $d$ elements
}{
\pageref{LOSZd}
}

\LOSitem
{
$\mathbb{Z}^{N}$
}{
rank-$N$ integer lattice
}{
\pageref{LOSVN_2}
}

\LOSitem
{
$\alpha=\left(  \alpha_{1}\alpha_{2}\dots\alpha_{n}\right)  $
}{
multi-index
}{
\pageref{LOSalpha}
}

\LOSitem
{
$\alpha=\left(  1,a,a^{2},\dots,a^{N-1}\right)  =\left(  1,e^{-\beta
},e^{-2\beta},\dots,e^{-\left(  N-1\right)  \beta}\right)  $
}{
Perron--Frobenius left eigenvector: $\TeXButton{bra}{\diracb{\alpha}}%
J=\lambda\TeXButton{bra}{\diracb{\alpha}}$
}{
\pageref{LOSPerronFrobeniuslefteigenvector_1},
\pageref{LOSPerronFrobeniuslefteigenvector_2},
\pageref{LOSalphaeigenvector_1},
\pageref{LOSPerronFrobeniuslefteigenvector_3},
\pageref{LOSalphaeigenvector_2},
\pageref{LOSalphaeigenvector_3},
\pageref{LOSalphaeigenvector_4},
\pageref{LOSalphaeigenvector_5},
}

\LOSitem
{
$\beta=\left(  a,a^{2},\dots,a^{N-1}\right)  =\left(  \frac{1}{\lambda}%
,\frac{1}{\lambda^{2\mathstrut}},\dots,\frac{1}{\lambda^{N-1\mathstrut}}%
\right)  $
}{
foreshortened vector
}{
\pageref{LOSbeta_1},
\pageref{LOSbeta_2}
}

\LOSitem
{
$\lambda$
}{
Perron--Frobenius eigenvalue of $J$
}{
\pageref{LOSPerronFrobeniuseigenvalue_1},
\pageref{LOSPerronFrobeniuseigenvalue_2},
\pageref{LOSPerronFrobeniuseigenvalue_3},
\pageref{LOSPerronFrobeniuseigenvalue_4},
\pageref{LOSPerronFrobeniuseigenvalue_5}
}

\LOSitem
{
$\Lambda=\left(  \theta\left(
\begin{smallmatrix}
1\\
0\\
0\\
\vdots\\
0
\end{smallmatrix}
\right)  ,\theta\left(
\begin{smallmatrix}
0\\
1\\
0\\
\vdots\\
0
\end{smallmatrix}
\right)  ,\dots,\theta\left(
\begin{smallmatrix}
0\\
0\\
0\\
\vdots\\
1
\end{smallmatrix}
\right)  \right)  $
}{
matrix implementing isomorphism $\theta$ between dimension groups: $\theta
(g)=\Lambda g$
}{
\pageref{LOSLambda_1},
\pageref{LOSLambda_2},
\pageref{LOSLambda_3},
\pageref{LOSLambda_4}
}

\LOSitem
{
$\xi\colon\Lambda\mathbf{v}=\xi\mathbf{v}^{\prime}$
}{
element of $\mathbb{Z}\left[  \frac{1}{\lambda}\right
]  $ with multiplicative inverse
}{
\pageref{LOSxiGreek}
}

\LOSitem
{
$\sigma$
}{
right shift on $\Omega$
}{
\pageref{LOSsigma}
}

\LOSitem
{
$\sigma_{i}$
}{
sections of $\sigma$
}{
\pageref{LOSsigmai}
}

\LOSitem
{
$\sigma_{t}^{L}\left(  s_{j}\right)  =\exp\left(  itL_{j}\right)  s_{j}
=\sigma\left(  e^{itL_{1}},\dots,e^{itL_{d}%
}\right)  \left(  s_{j}\right)  $
}{
one-parameter group $\sigma^{L}$ of $\ast$-au\-to\-mor\-phisms of
$\mathcal{O}_{d}$
}{
\pageref{LOSsigmatL_1},
\pageref{LOSsigmatL_2},
\pageref{LOSsigmatL_3}
}

\LOSitem
{
$\left(  \sigma^{L},\beta\right)  $-KMS state
}{
also $\beta$-KMS state
}{
\pageref{LOSsigmaLbetaKMSstate_1},
\pageref{LOSsigmaLbetaKMSstate_2},
\pageref{LOSsigmaLbetaKMSstate_3}
}

\LOSitem
{
$\tau$
}{
trace on $K_{0}\left(  \mathfrak{A}_{L}\right)  $
}{
\pageref{LOStau}
}

\LOSitem
{
$\tau\left(  g\right)  =\TeXButton{ip}{\ip{\alpha}{g}}$
}{
trace on $K_{0}\left(  \mathfrak{A}_{L}\right)  $ or on $\mathfrak{A}_{L}$
}{
\pageref{LOStaug}
}

\LOSitem
{
$\tau\left(  v\right)  =\TeXButton{ip}{\ip{\alpha}{v}}$
}{
inner product of left and right Perron--Frobenius vectors
}{
\pageref{LOStauv_1},
\pageref{LOStauv_2}
}

\LOSitem
{
$\omega=\omega_{p}$
}{
state given by $p$ (see $p_{j}$)
}{
\pageref{LOSomega}
}

\LOSitem
{
$\omega\left(  s_{\alpha}^{{}}s_{\gamma}^{\ast}\right)  =\delta_{\alpha\gamma
}e^{-\beta\sum_{k=1}^{\left|  \alpha\right|  }L_{\alpha_{k}}}$
}{
KMS state
}{
\pageref{LOSomegasalphasgamma}
}

\LOSitem
{
$\Omega=\prod_{1}^{\infty}\mathbb{Z}_{d}$
}{
product space
}{
\pageref{LOSOmega_1},
\pageref{LOSOmega_2}
}

\LOSitem
{
$\Omega_{i}=\sigma_{i}\left(  \Omega\right)  $
}{
partition of the product space $\Omega$ into clopen sets
}{
\pageref{LOSOmegai}
}

\LOSitem
{
$\Omega_{0}$
}{
cyclic vector
}{
\pageref{LOSOmega0_1},
\pageref{LOSOmega0_2}
}
\end{document}